%% file: the100.tex
\documentclass[12pt]{book}

\usepackage{a4wide}
\usepackage{graphicx}
\usepackage{color}
\usepackage{pictex}
\usepackage{amsthm}
\usepackage{amsmath}
\usepackage{amssymb}
\usepackage{mathrsfs}

\newtheorem{lemma}{Lemma}[section]
\newtheorem{theorem}{Theorem}[section]
\newtheorem{proposition}{Proposition}[section]
\newtheorem{remark}{Remark}[section]
\newtheorem{definition}{Definition}[section]
\newtheorem{conjecture}{Conjecture}[section]
\newtheorem{corollary}{Corollary}[section]

\begin{document}

\title{Quantum Superalgebras at Roots of Unity and Topological Invariants of
Three-manifolds}

\author{ {Sacha Carl Blumen} \\ 
{B.A. (Hons), A.Mus.A.} \\
{  } \\
{  } \\
{  } \\
{ {A thesis submitted in fulfilment of the requirements}} \\
{ {for the degree of Doctor of Philosophy}} \\
{  } \\
{  } \\
{  } \\
{  } \\
{  } \\
{  } \\
{  } \\
{  } \\
{\small {School of Mathematics and Statistics}} \\ 
{\small {University of Sydney, NSW, 2006, Australia}} \\
{\small {November 17th, 2005}}  }

\date{ }

\maketitle

\newpage

\pagestyle{myheadings}

\input{chaporiginality100.tex}

\addcontentsline{toc}{chapter}{Statement of Originality}

\input{chapacknowledgements100.tex}

\addcontentsline{toc}{chapter}{Acknowledgements}

\input{chapabstract1100.tex}

\addcontentsline{toc}{chapter}{Abstract}

\tableofcontents

\input{chapIntro1100.tex}

\input{chapter1100.tex}

\input{chapter2100.tex}

\input{chap2A100.tex}

\input{chapter4100.tex}
\appendix

\input{chapterappendixA100.tex}

\input{chapterappendixB100.tex}

\input{chapterappendixD100.tex}

\input{chapterappendixC100.tex}

\end{document}

%% file: chaporiginality100.tex
\pagenumbering{roman} \setcounter{page}{1} 
\chapter*{Statement of Originality\markboth{Statement of Originality}{Statement of Originality}} 

To the best of my knowledge and belief, this thesis is original and my own work
except as acknowledged in the text.  
This thesis has not been submitted, in whole or in part, for any other degree at this or
any other university.

\vspace{30mm}

\noindent
\underline{\hspace{50mm}}

\vspace{3mm}

\noindent
Sacha Carl Blumen

%% file: chapacknowledgements100.tex
\chapter*{Acknowledgements\markboth{Acknowledgements}{Acknowledgements}}

I think that it is almost always true that a Doctoral thesis is the combined result of the
candidate's work and the support of the people around them, and my thesis is no exception.
Many people have offered a huge amount of support during my studies, and there are a number
whom I particularly wish to thank.
Firstly, I would like to thank Troy Roderick, who has been a constant source of support,
especially when I felt that I might never complete the thesis.  
He helped me through a lot of times when things were very tough, 
and also proffered much fun in all the other parts of life that you can forget about 
while writing a thesis.
I would also like to thank Br\"{o}nte Metcalf-Roderick.

Thanks to my parents, Colleen and Carl Blumen, for their constant interest in my
studies, from undergraduate days to now.

Thanks to my friends for their support and understanding, 
especially when it sounded as if I were singing a continual refrain of ``just another 6 months!''.  

Thanks to my friends and colleagues in the academic world for just being there and being a
funny sort of `family'.
In particular, thanks to people in the Department of Mathematics, 
University of Queensland, where I started my Ph.D.: 
in no particular order (and to just mention a few), 
thanks to James Wood, Andrew Scott, Maithili Mehta, Nick Cavenagh, 
Andrew Blinco, Michiru Takizawa, and Professors Tony Bracken and Mark Gould.
Thanks to the Department of Mathematics for giving me a start in mathematics and
mathematical physics, something I had wanted to do for many years!

Thanks to friends and colleagues in 
the School of Mathematics and Statistics, the University of Sydney, where I completed this thesis:
fellow Ph.D. students,
my office-mates Mark Hopkins and Annelies Tjetjep for putting up with my tastes in
music and the way I appropriated the office, and
everyone else for their friendly faces and willingness to lend an ear. 
Thanks to the administrative staff for their help and friendliness, 
especially Sonia Morr and Janet Doyle.

Thanks to my Ph.D. supervisor, Dr. Rui Bin Zhang, 
who suggested the topic of this thesis 
and assisted with the multitudes of mathematical problems that appeared, and for also
offering a great deal of insightful criticism of my work.
Thanks to the School of Mathematics and Statistics for providing a space in 
which I was able to complete this thesis.

Finally, thanks to everyone else who helped me.  
Many people gave a lot to help me during my studies, 
and I hope that this work stands your support in good stead!

%% file: chapabstract1100.tex
\chapter*{Abstract\markboth{Abstract}{Abstract}} 

The general method of Reshetikhin and Turaev is followed to develop 
topological invariants of closed, connected, orientable $3$-manifolds 
from a new class of algebras called pseudo-modular Hopf algebras.
Pseudo-modular Hopf algebras are a class of $\mathbb{Z}_{2}$-graded 
ribbon Hopf algebras that generalise the concept of a modular Hopf 
algebra.

The quantum superalgebra $U_{q}(osp(1|2n))$ over $\mathbb{C}$ is 
considered with $q$ a primitive $N^{th}$ root of unity for all 
integers $N \geq 3$.  For such a $q$, a certain left ideal ${\cal{I}}$ 
of $U_{q}(osp(1|2n))$ is also a two-sided Hopf ideal, and the quotient 
algebra $U_{q}^{(N)}(osp(1|2n)) = U_{q}(osp(1|2n)) / 
{\cal{I}}$ is a $\mathbb{Z}_{2}$-graded ribbon Hopf algebra.

For all $n$ and all $N \geq 3$, a finite collection of 
finite dimensional representations of $U_{q}^{(N)}(osp(1|2n))$ is
defined.
Each such representation of $U_{q}^{(N)}(osp(1|2n))$ is labelled by an 
integral dominant weight belonging to the truncated dominant Weyl chamber.
Properties of these representations are considered:  the quantum 
superdimension of each representation is calculated, each 
representation is shown to be self-dual, and more importantly, 
the decomposition of the tensor product
of an arbitrary number of such representations is obtained for even $N$.
 
It is proved that the quotient algebra $U_{q}^{(N)}(osp(1|2n))$, 
together with the set of finite dimensional representations discussed above, form
a pseudo-modular Hopf algebra when $N \geq 6$ is twice an odd number.

Using this pseudo-modular Hopf algebra, we construct a
topological invariant of $3$-manifolds. 
This invariant is shown to be different to the topological invariants of 
$3$-manifolds arising from quantum $so(2n+1)$ at roots of unity.

%% file: chapIntro1100.tex
\begin{chapter}{Introduction}
\pagestyle{myheadings}
\label{chapterIntro1}
\markboth{\text{Chapter \ref{chapterIntro1}. Introduction}}
{\text{ }}
\pagenumbering{arabic} \setcounter{page}{1}

A striking development in mathematics in recent decades has been use of
ideas in mathematical physics to study the topology of low dimensional manifolds.
In 1983, Donaldson used self-dual Yang-Mills theory to
study $4$-manifolds, leading to the development of   
new topological invariants of $4$-manifolds \cite{donaldson83,donaldson87,donaldson90}.
In 1988, Witten introduced the concept of a topological quantum field theory (TQFT) \cite{witten88},
which provided a unifying principle for Donaldson's $4$-manifold invariants and Floer's 
invariants of oriented integral homology $3$-spheres \cite{floer88}. 
In 1989, Witten showed a connection between 
the Jones polynomial of links and quantum Chern-Simons theory that 
{\emph{``gave a (formal) $3$-dimensional interpretation of the Jones polynomial''}} \cite{witten,o}.
The quantum Chern-Simons theory has also been a rich source of topological invariants of
$3$-manifolds.
Its study has led to much progress in the area of knot theory and $3$-manifold theory.  
In particular, the influential work of Reshetikhin and Turaev \cite{rt} gave a mathematically rigorous
construction of the $3$-manifold invariants using quantum algebras.

This provides the first of the two general subject areas 
encompassing the work presented in this thesis: the study of 
topological invariants of  $3$-manifolds.  
The other general subject area is quite distinct from the first 
and can be easily studied without reference to the first.
The genesis of this thesis is in 
some very interesting work by Reshetikhin amongst others who showed unexpected connections between
this second area and the area of
topological invariants of links and $3$-manifolds, and it is this connection
that has essentially led to the theme of our research.
This second general subject area
is the study of quantum algebras and quantum superalgebras and their representations.

Each of these two subject areas is very interesting and offers many unanswered questions, 
some of which we discuss below. 
In this thesis we obtain some new results in both of these areas and
then weave them together yielding further results.

In this thesis, we firstly 
study the quantum superalgebra $U_{q}(osp(1|2n))$ and its representations at roots of unity
in Chapters \ref{chap2:titlelabel} and \ref{chap2A:titlelabel}
and then use the results obtained to study topological invariants of $3$-manifolds
in Chapter \ref{chaptitle:topologicalinvariants}.

Quantum algebras have been the subject of much research since Drinfel'd \cite{d1} 
and Jimbo \cite{j1,j2} introduced quantum groups to the mathematical world in approximately 1985.
Quantum superalgebras were introduced in \cite{brackgouldzhang,y} and other papers.
Quantum algebras and quantum superalgebras related to classical Lie algebras and simple basic Lie
superalgebras, respectively, are $q$-deformations of the 
universal enveloping algebras of the relevant Lie algebra or Lie superalgebra.
The representation theory of quantum algebras and quantum superalgebras 
has attracted a lot of attention, 
not only for the intrinsic interest in the representations of these new algebras, 
but also for the possibility that it could solve already existing problems, 
and in addition, potential applications.
The potential for applications was quite strong
as the quantum groups introduced by Drinfel'd admitted 
an element  satisfying the Yang-Baxter equation, an important equation 
 in statistical mechanics.  
The existence and explication of these elements,  
called universal $R$-matrices, has been one of the goals in the study
of quantum algebras and quantum superalgebras.

The representation theories of quantum algebras and quantum superalgebras are at present 
not completely known, but aspects of them are somewhat known.  
In some ways the known representation theory depends dramatically on the nature of $q$:
if $q$ is non-zero and not a root of unity,  the representation theory of quantum
algebras and quantum superalgebras is much better known than if $q$ is a root of unity.
In one sense, this may arise from the fact that 
the centres of the quantum algebras and quantum superalgebras are much
larger when $q$ is a root of unity.

The representation theory of quantum algebras at roots of unity has attracted much research
(eg \cite{a,andersenjantzensoergel,ap}).
However, the representation theory of
 quantum superalgebras at roots of unity is not nearly as well known.
 For example, the representation theory of $U_{q}(osp(1|2n))$ at roots of unity was barely studied
 before.
 The first part of this thesis studies the structure and representations of $U_{q}(osp(1|2n))$ at roots
 of unity.

In particular, we fix $q = \exp{(2 \pi i/N)}$ where $N \geq 3$ is an integer.
In this case, the quantum superalgebra $U_{q}(osp(1|2n))$ has a much larger centre than that
at generic $q$.
A two-sided ideal ${\cal{I}}$ is generated by certain central elements, 
which is also a Hopf ideal.
The quotient quantum superalgebra
 $U_{q}^{(N)}(osp(1|2n)) = U_{q}(osp(1|2n))/{\cal{I}}$ 
is again a $\mathbb{Z}_{2}$-graded Hopf algebra, and,
importantly, admits a universal $R$-matrix.
It transpires that $U_{q}^{(N)}(osp(1|2n))$ is a $\mathbb{Z}_{2}$-graded ribbon Hopf algebra.

We construct a set of $U_{q}^{(N)}(osp(1|2n))$-modules, each of which is characterised by an element
from a truncated Weyl chamber $\Lambda_{N}^{+}$.
We prove that each of these modules is self-dual,
and, more importantly, their tensor products decompose in a very nice way.
It is the tensor product decomposition theorems which we need for the cosntruction of
$3$-manifold invariants.

The second major theme in this thesis is the construction of topological invariants of closed,
connected, orientable $3$-manifolds.  
We review the construction of the topological invariants from modular Hopf algebras introduced by
Reshetikhin and Turaev \cite{rt}, and show that topological invariants, similar to those introduced in
\cite{rt}, can be constructed from a more general class of algebras which we call pseudo-modular Hopf
algebras.  The conditions for an algebra to be modular are quite prescriptive, and 
we define pseudo-modular Hopf algebras to be much more general, 
so that  topological invariants can be constructed following our work for as many ribbon Hopf algebras
as possible.

After these two independent strands of study, we tie them together by showing that
the quotient quantum superalgebra $U_{q}^{(N)}(osp(1|2n))$, 
when $q = \exp{(2 \pi i/N)}$ and $N \geq 6$ satisfies
$N \equiv 2 \pmod{4}$, together with a set of finite dimensional representations, is  a
pseudo-modular Hopf algebra, and thus yields 
topological invariants of closed, connected, orientable $3$-manifolds.

We now turn to a discussion of $3$-manifold invariants, proper.

Despite the relative age of the programme of studying $3$-manifolds 
(it being initiated in the 1880s \cite{gordon}), there are still  
many avenues of research, as evidenced by the recent claim of the
proof of the Poincar\'{e} Conjecture \cite{morgan}.  
In addition, it is still unknown whether the problem of classifying 
$3$-manifolds into homeomorphism equivalence classes, one of the major goals of
the research programme, is solvable \cite{o}.
This contrasts with what is known for simply connected $n$-manifolds 
for all integers $n \geq 2$ except $n=3$.
Here the classification problem has been solved one way or the other:
it has been solved 
for $n=2$ and $n \geq 5$, and is known to be algorithmically unsolvable for $n=4$, as
{\emph{``the set of fundamental groups of $4$-manifolds is sufficiently large to be algorithmically
unsolvable as a word problem in group theory''}} \cite{o}.
The discovery of new topological invariants of $3$-manifolds should shed light on 
the classification problem;
this is the general goal towards which part of the work presented in this thesis is directed.

While it is not known whether the classification problem for $3$-manifolds is solvable, 
Lins has discussed
a potential conjecture for completely classifying connected $3$-manifolds \cite{lins}.  
This potential conjecture uses a
combinatorial algorithm for classifying $3$-manifolds using the so-called $3$-gems, or
{\emph{$3$-dimensional graph-encoded manifolds}} \cite{kl,lins}.
Using $3$-gems, Kauffman and Lins calculated topological invariants 
for a large number of $3$-manifolds  \cite{kl} and 
while the algorithm works on the elements of a proper subset of all 
connected $3$-manifolds where the calculations are computationally tractable, 
it is not known whether the algorithm works in general.

Even if one does not have a complete topological invariant of 
connected $3$-manifolds, a collection of 
topological invariants may altogether be a complete topological invariant.
In any event, different topological invariants may have different properties, 
so it is useful to have as many topological invariants as possible  
when studying individual $3$-manifolds and the classification problem.

This provides the raison d'etre for part of this thesis:
we will construct a topological invariant of closed, connected, orientable $3$-manifolds. 
Our invariants are one of the class of  {\emph{quantum invariants}}  \cite{o},
which are a large collection of topological invariants of
knots and $3$-manifolds that have been discovered since
Jones discovered a new polynomial of links in 1985 \cite{jones1}
(for examples of quantum invariants of tangles and knots see 
\cite{homfly,blm,hoste,lm,t88,pt,da,resh,reshtur,kau,z0,zgb,lg,z1,gould93,gtb,t,glz,z5,ddw}).  
The Jones polynomial was the first link invariant discovered since Alexander
discovered, in 1928,  the Alexander-Conway polynomial 
\cite{alex} (see also \cite{con}).  
The topological invariants of $3$-manifolds generated by TQFTs are also quantum invariants.
 
The topological invariants of $3$-manifolds we construct in this thesis
were inspired by the construction of Reshetikhin and Turaev in 1991 \cite{rt}.  
Their construction was formulated using a {\emph{modular Hopf algebra}}, 
which is a ribbon Hopf algebra together with a finite collection of finite
dimensional irreducible representations with non-zero quantum dimension, satisfying certain
criteria including a subtle tensor product rule.  
These criteria encapsulate the necessary conditions for the construction of the
$3$-manifold invariants following  \cite{rt}, although, as we shall discuss
below, it is possible to construct $3$-manifold invariants using ribbon Hopf algebras satisfying
 less restrictive criteria than those satisfied by modular Hopf algebras.

The essential method for constructing the topological invariants presented in \cite{rt} is as follows.  
One uses a result of Lickorish \cite{l1} to present each closed, connected, orientable $3$-manifold 
$M_{L}$ as the result of performing surgery on the $3$-sphere $S^{3}$ along a framed link 
$L \subset S^{3}$.
From the work of Kirby \cite{k} and Fenn and Rourke \cite{fr}, 
two such $3$-manifolds, $M_{L_{1}}$ and $M_{L_{2}}$, are homeomorphic if and only if the links 
$L_{1}$ and $L_{2}$ are equivalent with respect to the so-called
Kirby moves, which are a set of equivalence relations on planar projections of a link.  
This enables one to turn the study of homeomorphism classes of closed,
connected, orientable $3$-manifolds into the study of equivalence classes of framed links in 
$S^{3}$ where the equivalence is generated by the Kirby moves.
Then, upon taking any isotopy invariant of framed links that is invariant under the Kirby moves,
one obtains a topological invariant of closed, connected, orientable $3$-manifolds.
Reshetikhin and Turaev put this scheme into effect in \cite{rt} 
by using an isotopy invariant of framed
links they had constructed in 1990 from quantum algebras  \cite{reshtur}.
Interestingly, the 
isotopy invariant that they used from \cite{reshtur} 
can be seen as a generalisation of the Jones polynomial.
By taking a linear sum of isotopy invariants of framed links that remained unchanged 
if the Kirby moves were applied to the framed link, 
Reshetikhin and Turaev obtained a topological invariant of 
closed, connected, orientable $3$-manifolds.

Since Reshetikhin and Turaev's work in 1991, 
many topological invariants of $3$-manifolds have been constructed
from a variety of different viewpoints, for example,  
they have been constructed following the general
idea of \cite{rt} by taking linear sums of invariants of framed links that are unchanged under
the Kirby moves (eg see \cite{l2,l91(2),walk,bhmv,l92,tw,z3}), 
by taking averages over all possible cablings of invariants of framed links \cite{w10},
by taking a state-sum over an arbitrary triangulation of the $3$-manifold using quantum $6j$-symbols
\cite{tv} and by using a Heegard decomposition of a $3$-manifold \cite{kohno}.

The topological invariants of $3$-manifolds constructed since 1991 
have been much studied: in particular the invariants of Reshetikhin and
Turaev have been greatly studied due to the comparative ease of their calculation
(eg see \cite{km,l93,kerler}).
In addition to their status of new topological invariants,
this study yielded applications to classical $3$-manifold topology -- for example 
the Reshetikhin-Turaev invariants were used to develop obstructions to
embedding one 3-manifold in another \cite{fkb}.

The construction of $3$-manifold invariants in \cite{rt} is general, 
in that it provides a way to construct the invariants 
from any modular Hopf algebra.
It was shown in \cite{rt} that a quotient of $U_{q}(sl_{2})$ at $4k^{th}$ primitive roots of unity
was modular.
This naturally led researchers to ask whether other quantum
algebras were also modular. 
But the modularity criteria involves representation-theoretic conditions, 
and the relevant representation theory of quotients of other quantum algebras 
at primitive roots of unity was not known.

However, Turaev and Wenzl managed to somewhat circumvent this problem in 1993 by 
constructing topological invariants of closed, connected, orientable $3$-manifolds
from  {\emph{quasimodular Hopf algebras}}, which are more general than
modular Hopf algebras \cite{tw}. 
The essential difference between modular and quasimodular Hopf algebras is that the 
condition that the modules be irreducible is relaxed for quasimodular Hopf algebras.
The quantum algebras related to the families of complex Lie algebras
$A_{n}, B_{n}, C_{n}$ and $D_{n}$ are quasimodular at even primitive roots of unity, 
and the authors constructed $3$-manifold invariants using these algebras \cite{tw}.

In 1994 Turaev reformulated modular Hopf algebras 
using the language of categories, and showed that
modular categories, which are tensor categories with additional structures, 
could be used to construct $3$-manifold invariants \cite{t}.
Informally, the key features of a modular category are a finite set of objects satisfying a certain
tensor product property and the invertibility of a certain matrix
commonly called the $S$-matrix.
This was followed by the study of pre-modular categories by Brugui\`{e}res \cite{brug00,brug00(2)}, 
which are a class of
categories more general than modular categories, in which the $S$-matrix is not necessarily invertible.  
Invariants of links and tangles and sometimes $3$-manifolds
can be constructed using a pre-modular category \cite{bb01}. 
Although we have not explored the possibility, it appears plausible that 
the pseudo-modular Hopf algebras we define in this thesis may yield examples of
pre-modular categories.

In 1995, the construction of Reshetikhin-Turaev $3$-manifold invariants 
was extended to modular Hopf algebras derived from 
$\mathbb{Z}_{2}$-graded ribbon Hopf algebras by Zhang \cite{z3}, 
and the first such topological invariants
were constructed from quotients of the quantum superalgebras 
$U_{q}(osp(1|2))$ \cite{z2}, and $U_{q}(gl(2|1))$ \cite{z3} at odd roots of unity.
This construction slightly differed from that of \cite{rt}, 
in that it was not shown that the algebras were modular,  
but the author found a way to avoid this requirement 
by using known features of the representation theory of these algebras.

In contrast to the construction of topological invariants of $3$-manifolds using 
quotients of quantum algebras at roots of unity and their representations, 
the construction of topological invariants of $3$-manifolds using 
quotients of quantum superalgebras at roots of unity and their representations,
has, unfortunately, been quite limited, this limitation having arisen from the limited
knowledge of the representation theory of quantum superalgebras and their quotients at roots of unity.
This was the motivation for our earlier-mentioned study 
of  representations of a quotient of
$U_{q}(osp(1|2n))$ at roots of unity in this thesis.

As mentioned previously,
we then show that topological invariants of closed, connected, orientable $3$-manifolds can be
constructed from a class of $\mathbb{Z}_{2}$-graded ribbon Hopf algebras 
(together with a finite set of their representations) satisfying more general conditions than
those satisfied by a modular or quasimodular Hopf algebras.  
We call these algebras {\emph{pseudo-modular Hopf algebras}}.
The key difference between a pseudo-modular and a modular Hopf algebra is that the $S$-matrix of a
pseudo-modular Hopf algebra is not necessarily invertible, while it is invertible for
modular Hopf algebras.

After defining pseudo-modular Hopf algebras and showing that
$3$-manifold invariants can be constructed from pseudo-modular Hopf algebras, 
we show that the quotient algebra
$U_{q}^{(N)}(osp(1|2n))$, where $q = \exp{(2 \pi i/N)}$ and $N \geq 6$ satisfies 
$N \equiv 2 \pmod{4}$, together with a set of its representations, 
is pseudo-modular and thus yields $3$-manifold invariants.

A preliminary announcement of the work presented in this thesis was given in \cite{blumen}.

\begin{section}{Summary of new results}

We now briefly summarise the new results of each chapter, and refer  to the body of the
thesis for the precise statement of each result mentioned here.

\begin{subsubsection}{Chapter \ref{chap2:titlelabel}}

Let $U_{\mathbb{C}[[h]]}(osp(1|2n))$  denote the Drinfel'd version of the quantum superalgebra
over the ring $\mathbb{C}[[h]]$, for an indeterminate $h$.
Write the universal $R$-matrix of $U_{\mathbb{C}[[h]]}(osp(1|2n))$ as 
$R = \widetilde{K} \widetilde{R}$ (eg see \cite{kt}), where 
$\widetilde{R}$ is an infinite sum of root vectors in $U_{\mathbb{C}[[h]]}(osp(1|2n))$. 
One of the problems in dealing with Jimbo's quantum superalgebra
$U_{q}(osp(1|2n))$ over $\mathbb{C}$ is that $U_{q}(osp(1|2n))$ 
is not complete and thus does not admit a universal $R$-matrix.

In Section \ref{subsec:RRRR} we 
define $R$-matrices for finite dimensional representations of $U_{q}(osp(1|2n))$ 
at generic $q$, as is done in quantum algebras over $\mathbb{C}$ \cite{cp,ks}.  
We do this by defining a completion  $\overline{U}_{q}^{+}(osp(1|2n))$ of
$U_{q}(osp(1|2n))$ which contains the factor $\widetilde{R}$ of the universal $R$-matrix.
There is the problem of dealing with the factor $\widetilde{K}$ of the 
universal $R$-matrix of $U_{\mathbb{C}[[h]]}(osp(1|2n))$, as it is not clear which element 
 of  $\overline{U}_{q}^{+}(osp(1|2n))$, if any, corresponds to $\widetilde{K}$.  
We define an element, ${\cal{E}}_{\mu} \in U_{q}(osp(1|2n))$ 
for each integral dominant weight $\mu$ such that ${\cal{E}}_{\mu}$ has the same action
on $V_{\mu} \otimes V_{\lambda}$, for any integral dominant $\lambda$,
as would $\widetilde{K}$  if it existed.
Then  ${\cal{E}}_{\mu} \widetilde{R}$ is an `$R$-matrix' for $V_{\mu} \otimes V_{\lambda}$ in that
it is an $R$-matrix for representations.

In Section \ref{sec:'central'element} we define some useful elements of 
$\overline{U}_{q}^{+}(osp(1|2n))$.
These elements act in finite dimensional irreducible representations of
$U_{q}(osp(1|2n))$ in a similar way as the elements $u$ and $v$ 
of a $\mathbb{Z}_{2}$-graded ribbon Hopf algebra.
Let $V_{\lambda}$ be the finite dimensional irreducible $U_{q}(osp(1|2n))$-module with
integral dominant highest weight $\lambda$.
For each integral dominant $\lambda$, define the elements
$u_{\lambda}, v_{\lambda} \in \overline{U}_{q}^{+}(osp(1|2n))$.
Even though $v_{\lambda}$ does not necessarily
commute with each element of $U_{q}(osp(1|2n))$, it
acts on  $V_{\lambda}$ 
as  multiplication by the scalar
$q^{-(\lambda + 2\rho, \lambda)}$.

In Section \ref{sec:spectam(a)} we state and prove the spectral decomposition of the element
$$\check{\cal{R}}_{V,V} \in End_{U_{q}(osp(1|2n))}(V \otimes V)$$ 
defined by
$\check{\cal{R}}_{V,V} (x \otimes y) = P \circ R_{V,V}(x \otimes y)$, 
where $P$ is the graded permutation
operator and $R_{V,V}$ is the $R$-matrix for the tensor product of
$U_{q}(osp(1|2n))$-representations $V \otimes V$.

In Section \ref{eq:theXfactor(a)} we show that there exists a representation of the 
Birman-Wenzl-Murakami algebra $\mathscr{BW}_{f}(-q^{2n},q)$ in 
$End_{U_{q}(osp(1|2n))}(V^{\otimes f})$.

In Section \ref{subsec:projectontome} we use the  elements
$v_{\lambda}$ of $\overline{U}_{q}^{+}(osp(1|2n))$
to define a set of mutually orthogonal idempotents 
in $End_{U_{q}(osp(1|2n))}(V^{\otimes f})$ that
project down from $V^{\otimes f}$ (which is completely reducible) 
onto its irreducible $U_{q}(osp(1|2n))$-submodules.

In Section \ref{subsec:rhapsodyinred} we define matrix units in 
$End_{U_{q}(osp(1|2n))}(V^{\otimes f})$
which we obtain from the projections in Section \ref{subsec:projectontome} 
and the matrix units in a semisimple quotient of 
$\mathscr{BW}_{f}(-q^{2n},q)$.
The intertwiner matrix units in $End_{U_{q}(osp(1|2n))}(V^{\otimes f})$ 
are $U_{q}(osp(1|2n))$-linear maps between isomorphic irreducible $U_{q}(osp(1|2n))$-submodules of
$V^{\otimes f}$.

\end{subsubsection}

\begin{subsubsection}{Chapter \ref{chap2A:titlelabel}}

In this chapter we define the quotient algebra $U_{q}^{(N)}(osp(1|2n))$
and its relevant representations. 
The new results in this chapter are as follows.
Fix $q = \exp{(2 \pi i/N)}$ where $N \geq 3$ is an integer.
In Theorem \ref{th:cooljazz(200)} we prove that a certain left ideal 
${\cal{I}} \subset U_{q}(osp(1|2n))$ is 
a two-sided Hopf ideal,  thus  
$$U_{q}^{(N)}(osp(1|2n)) = U_{q}(osp(1|2n)) / {\cal{I}},$$
is a $\mathbb{Z}_{2}$-graded Hopf algebra.
This quotient algebra has apeared in the literature but only for the case $n=1$ where 
$N \geq 3$ is odd \cite{z2}. 
The proof of Theorem \ref{th:cooljazz(200)}  involves many intricate and involved calculations and we
leave these calculations to Appendix \ref{chap:appendixC}.

In Proposition \ref{prop:R-matrichereicome} we state
the universal $R$-matrix of $U_{q}^{(N)}(osp(1|2n))$, which was originally given in \cite{z1}.
In Theorem \ref{th:ribbonalgebraribbonalgebra} we prove that $U_{q}^{(N)}(osp(1|2n))$ is a
$\mathbb{Z}_{2}$-graded ribbon Hopf algebra.

We then construct representations of $U_{q}^{(N)}(osp(1|2n))$ that we will use in constructing the
topological invariant.  
Lemma \ref{lem:manydimroot} discusses the fundamental irreducible
$U_{q}^{(N)}(osp(1|2n))$-module $V$;  all the $U_{q}^{(N)}(osp(1|2n))$-modules that we later define
are submodules of tensor powers of $V$.
In Subsection \ref{subsec:capitallambda} we introduce the set $\overline{\Lambda_{N}^{+}}$ of integral
weights in the truncated fundamental Weyl alcove.  
In Subsection \ref{subsec:7zark7} we define a $U_{q}^{(N)}(osp(1|2n))$-module $V_{\lambda}$ associated
with each integral dominant $\lambda \in \overline{\Lambda_{N}^{+}}$.
In Lemma \ref{lem:kilcutiecutie} we give the quantum superdimension of each of these 
$U_{q}^{(N)}(osp(1|2n))$-modules.
The quantum superdimension of $V_{\lambda}$ has the same expression as does the quantum superdimension of 
the irreducible $U_{q}(osp(1|2n))$-module with highest weight $\lambda$
but we specialise $q$ to the appropriate root of unity.
We do not know whether the $U_{q}^{(N)}(osp(1|2n))$-modules
defined in Subsection \ref{subsec:7zark7} are irreducible or even indecomposable, but we
conjecture that they are all irreducible.

In Proposition \ref{prop:moo20} we prove that each $U_{q}^{(N)}(osp(1|2n))$-module
defined in Subsection \ref{subsec:7zark7} is self-dual.

In Section \ref{sec:fedupwiththetensors} we prove tensor product theorems for the
$U_{q}^{(N)}(osp(1|2n))$-modules defined in Subsection \ref{subsec:7zark7} at even $N$.  
The proofs are based on the proofs of similar tensor product theorems for representations of quantum
algebras at even roots of unity \cite[Thm. 5.2.2]{tw}.
The new tensor product theorems are reminiscent of similar tensor product theorems for modular and
quasimodular Hopf algebras, and are perhaps analogues of them.
\end{subsubsection}

\begin{subsubsection}{Chapter \ref{chaptitle:topologicalinvariants}}

In this chapter we construct a topological invariant of closed,
connected, orientable $3$-manifolds.
In Subsection \ref{subsec:allanbridges2} we define 
{\emph{pseudo-modular Hopf algebras}} which generalise
modular Hopf algebras. In Section \ref{sec:rushingwater} we prove that a topological
invariant of closed, connected, orientable $3$-manifolds
can be constructed from a pseudo-modular Hopf algebra and its representations.

The essential differences between a pseudo-modular Hopf algebra and a modular Hopf algebra are
three-fold: (i) the irreducibility conditions on the modules are relaxed for a pseudo-modular Hopf
algebra, (ii) the condition that a certain matrix equation has a unique set of solutions 
is relaxed to a condition that the matrix equation has {\emph{at least one}} set of solutions, and
(iii) a certain sum must be non-zero.
In modular Hopf algebras the result corresponding to (iii) is automatically true,  
which is a consequence of the fact that the matrix equation has a unique set of solutions 
(see \cite{tw}). 
The condition in (ii) for modular Hopf algebras that the matrix equation has a unique set of solutions 
boils down to a condition that the $S$-matrix be invertible, and the importance of this is that the
solutions to the matrix equation are exactly the weights in the linear sum of isotopy invariants of
framed links giving rise to the $3$-manifold invariant.
For pseudo-modular Hopf algebras we merely require that there exists at least one set of (non-zero) 
solutions to the matrix equation.

Theorem \ref{theorem:bigbigbig} contains the core result of this thesis:
we prove that $U_{q}^{(N)}(osp(1|2n))$, together with a set of its representations
$\left\{ V_{\lambda} \big| \ \lambda \in \Lambda_{N}^{+} \right\}$
(here $\Lambda_{N}^{+} \subset \overline{\Lambda_{N}^{+}}$),
and a collection of other data,
is a pseudo-modular Hopf algebra when $N \geq 6$ satisfies $N \equiv 2 \pmod{4}$.  
Furthermore, in Theorem \ref{lem:dalektron4} we prove that 
$U_{q}^{(N)}(osp(1|2n))$ 
is {\emph{not}} a pseudo-modular Hopf algebra when $N \geq 4$ satisfies
$N \equiv 0 \pmod{4}$, and that $3$-manifold
invariants cannot be constructed from $U_{q}^{(N)}(osp(1|2n))$ when $N \equiv 0 \pmod{4}$.

A relevant question to ask is whether our topological invariants are complete 
or even new invariants.  
These are very difficult questions to answer, and it is more tractable  to 
ask whether our topological invariants are the same as existing  invariants.
The set $\Lambda_{N}^{+}$ for $U_{q}^{(N)}(osp(1|2n))$ when $N \geq 6$ satisfies
$N \equiv 2 \pmod{4}$ is identical to the corresponding set for $U_{q}^{(N/2)}(so(2n+1))$, and 
in Section \ref{chapter4sectionlabel(comparing_the_)} we compare the
topological invariants constructed using $U_{q}^{(N)}(osp(1|2n))$ to the invariants
constructed using $U_{q}^{(N/2)}(so(2n+1))$.  
We show that the invariants from these algebras are {\emph{not}} the same.

\end{subsubsection}

\begin{subsubsection}{Appendices \ref{appen:gausssss}--\ref{chap:appendixC}}

In Appendix \ref{chap:appendixB} we prove two generalisations of the $q$-binomial theorem
that we use in Appendix \ref{chap:appendixC}.
In Appendix \ref{chap:appendixC} we prove that the left ideal 
${\cal{I}} \subset U_{q}(osp(1|2n))$  at roots of unity is a two-sided Hopf ideal.

\end{subsubsection}

\end{section}

\end{chapter}

%% file: chapter1100.tex
\begin{chapter}{Notation and background}
\label{chap:groovytuesday(1)}
\pagestyle{myheadings}
\markboth{\text{Chapter \ref{chap:groovytuesday(1)}. Notation and background}}
{\text{ }}

In this chapter we introduce some notation that we use
in this thesis and we also give relevent background material.

The plan of this chapter is as follows.
In Section \ref{sec:preliminaries} we introduce some of the notations and conventions 
that we use in this thesis.
In Section \ref{sec:definandbackgroundmaterial} we introduce the elementary
algebraic structures that we will refer to and often use in this thesis.
In Section \ref{chap1:quastriangaulrhopfalgebras} we define  
$\mathbb{Z}_{2}$-graded quasitriangular Hopf algebras, 
which are a class of $\mathbb{Z}_{2}$-graded Hopf algebras
admitting a certain element called the universal $R$-matrix.
These algebras are intrinsically mathematically interesting and are also important in physics
(eg in statistical mechanics \cite{kulskl}), 
as the universal $R$-matrix furnishes a
solution to the Yang-Baxter equation in each representation of the algebra.

We then define $\mathbb{Z}_{2}$-graded ribbon Hopf algebras.
This last class of algebras is useful in constructing knot invariants (eg see \cite{o}) and is
very important in our construction of 
topological invariants of $3$-manifolds in Chapter \ref{chaptitle:topologicalinvariants}.

\begin{section}{Notation}
\label{sec:preliminaries}
\markright{\text{Notation}}

In this section we introduce the notation and conventions that we will use in this thesis:
\begin{itemize}
\item[]  $\mathbb{Z}$ the set of integers,
\item[]  $\mathbb{Z}_{+}$ the set of non-negative integers,
\item[]  $\mathbb{Z}_{N}=\mathbb{Z}/N\mathbb{Z}$,
\item[]  $\mathbb{N}$ the set of positive integers,
\item[]  $\mathbb{C}$ the field of complex numbers,
\item[]  $\mathbb{R}$ the field of real numbers.
\end{itemize}

\noindent
For each $N \in \mathbb{Z}$, set
$$N' = \left\{ 
\begin{array}{ll}
N, & \mbox{if $N$ is odd}, \\
N/2, & \mbox{if $N$ is even},
\end{array}
\right.  \hspace{15mm} 
\overline{N} = \left\{ 
\begin{array}{ll}
2N, & \mbox{if $N$ is odd}, \\
N, & \mbox{if } N \equiv 0 \pmod{4}, \\
N/2, & \mbox{if } N \equiv 2 \pmod{4}.
\end{array}
\right. $$ 

\noindent
For each $q = \rho e^{i \theta} \in \mathbb{C}$ where $0 \leq \theta < 2\pi$ and $\rho >0$,
set  $q^{1/2} = +\sqrt{\rho}e^{i \theta/2}$.

\noindent
For each $q \in \mathbb{C}$ not equal to $0$ or $1$, define 
$$[n]^{q} 
= \frac{1-q^{n}}{1-q}, 
\hspace{10mm} [n]^{q}! = [n]^{q}[n-1]^{q}\cdots [1]^{q}, \hspace{5mm} n \in \mathbb{N},
\hspace{15mm} [0]^{q}! = 1,$$
and the Gaussian binomial
\begin{equation}
\label{eq:guasbinomial}
\left[ \begin{array}{c}
n \\
i
\end{array}\right]^{q} = \frac{\left[n\right]^{q}!}{\left[i\right]^{q}! \left[n-i\right]^{q}!},  \hspace{10mm} i \leq n.
\end{equation}

\noindent
We also define
$$(n)_{q} 
= \frac{1-(-q)^{n}}{1+q}, \hspace{10mm}
(n)_{q}! = (n)_{q} (n-1)_{q} \cdots (1)_{q}, \hspace{5mm} n \in \mathbb{N}, 
\hspace{15mm} (0)_{q}! = 1,$$
and the pseudo-Gaussian binomial
\begin{equation}
\label{eq:pseduogausbinom}
\left( \begin{array}{c}
 n \\
 i
\end{array} \right)_{q} = \frac{(n)_{q}!}{(i)_{q}!(n-i)_{q}!}, \hspace{10mm} i \leq n.
\end{equation}
Note that $(n)_{q} = [n]^{-q}$ 
but we introduce the two notations as this is convenient when $q$ is a root of unity.
We also define 
$$[n]_{q}= q^{(1-n)/2}(n)_{q} =\frac{q^{-n/2}-(-1)^{n}q^{n/2}}{q^{-1/2}+q^{1/2}}, 
\hspace{10mm} n \in \mathbb{N},$$
and
$[n]_{q}! = [n]_{q}[n-1]_{q} \cdots [1]_{q}, \  [0]_{q}! = 1$ and
$
\left[ \begin{array}{c}
n \\
i
\end{array}
\right]_{q} = \displaystyle{\frac{[n]_{q}!}{[i]_{q}![n-i]_{q}!}}$, $i \leq n$.

\noindent
For each non-empty set $A$, we define the delta function
$\delta: A \times A \rightarrow \mathbb{C}$ by
$$\delta_{x,y} = \left\{ \begin{array}{ll}
1, & \mbox{if } x=y, \\
0, & \mbox{if } x \neq y,
\end{array} \right. \hspace{5mm}  \mbox{for all } x, y \in A.$$

\end{section}

\begin{section}{Background algebraic structures}
\label{sec:definandbackgroundmaterial}
\markright{\text{Background algebraic structures}}

In this section we introduce the background
algebraic structures that we refer to and often use in this thesis \cite{sch,gzb}.
Throughout this section we denote the two elements of $\mathbb{Z}_{2}$ by $0$ and $1$.

\begin{definition}
Let $V$ be a vector space over $\mathbb{C}$.  A $\mathbb{Z}_{2}$-gradation of the vector space
$V$ is a family $(V_{\gamma})_{\gamma \in \mathbb{Z}_{2}}$ of subspaces of $V$ such that
$$V = \bigoplus_{\gamma \in \mathbb{Z}_{2}}V_{\gamma}.$$
The vector space $V$ is said to be $\mathbb{Z}_{2}$-graded if it is equipped with a
$\mathbb{Z}_{2}$-gradation.
\end{definition}

An element of $V$ is called 
{\emph{homogeneous of degree $\gamma$}}, $\gamma \in \mathbb{Z}_{2}$, if it
is an element of $V_{\gamma}$.  The elements of $V_{0}$ (resp. $V_{1}$) are called even
(resp. odd).  We define a mapping $[ \cdot ]: V_{\gamma} \rightarrow \{0, 1\}$ for each $\gamma \in
\mathbb{Z}_{2}$ by setting $[y] = \gamma$ for each $y \in V_{\gamma}$.

Each element $y \in V$ has a unique decomposition of the form
$$y=\sum_{\gamma \in \mathbb{Z}_{2}} y_{\gamma}, \hspace{10mm} y_{\gamma} \in V_{\gamma}.$$
The element $y_{\gamma}$ is called the {\emph{homogeneous component}} of $y$ of degree $\gamma$.

A subspace $U$ of $V$ is called $\mathbb{Z}_{2}$-graded if it contains the homogeneous components of
all of its elements, that is, if
$$U = \bigoplus_{\gamma \in \mathbb{Z}_{2}} U_{\gamma}, \hspace{10mm}  U_{\gamma} = (U \cap V_{\gamma}).$$

\begin{definition}
Let $V = \bigoplus_{\gamma \in \mathbb{Z}_{2}} V_{\gamma}$ and 
$W = \bigoplus_{\delta \in \mathbb{Z}_{2}} W_{\delta}$ be two $\mathbb{Z}_{2}$-graded vector spaces 
over $\mathbb{C}$.   
A $\mathbb{C}$-linear mapping
$$\psi: V \rightarrow W,$$
is said to be homogeneous of degree $\gamma \in \mathbb{Z}_{2}$, if
$$\psi(V_{\alpha}) \subseteq W_{\alpha + \gamma}, \hspace{5mm} \mbox{ for all } \alpha \in
\mathbb{Z}_{2}.$$

The mapping $\psi$ is called a homomorphism of the $\mathbb{Z}_{2}$-graded vector space $V$ into the
$\mathbb{Z}_{2}$-graded vector space $W$ if $\psi$ is homogeneous of degree $0$.  
\end{definition}

Let $V = \bigoplus_{\gamma \in \mathbb{Z}_{2}} V_{\gamma}$ and 
$W = \bigoplus_{\delta \in \mathbb{Z}_{2}} W_{\delta}$ be 
two $\mathbb{Z}_{2}$-graded vector spaces over $\mathbb{C}$.
The vector space $\mathrm{Hom}_{\mathbb{C}}(V, W)$ of $\mathbb{C}$-linear maps from $V$ to $W$ 
admits a $\mathbb{Z}_{2}$-gradation 
$$\mathrm{Hom}_{\mathbb{C}}(V, W) = 
\bigoplus_{\alpha \in \mathbb{Z}_{2}} \mathrm{Hom}_{\mathbb{C}}(V, W)_{\alpha},$$
where $\mathrm{Hom}_{\mathbb{C}}(V, W)_{\alpha}=
\big\{ \psi \in \mathrm{Hom}_{\mathbb{C}}(V, W) | 
 \ \psi \mbox{ is homogeneous of degree } \alpha \in \mathbb{Z}_{2} \big\}$.
 We denote $End_{\mathbb{C}}(V) = \mathrm{Hom}_{\mathbb{C}}(V, V)$.

If $W=\mathbb{C}$ is regarded as a $\mathbb{Z}_{2}$-graded vector space with zero odd subspace,  
$\mathrm{Hom}_{\mathbb{C}}(V, \mathbb{C})$ 
is called the dual space of $V$ and is denoted by $V^{*}$, which is clearly $\mathbb{Z}_{2}$-graded.

\begin{definition}
Let $V$ be a $\mathbb{Z}_{2}$-graded vector space over $\mathbb{C}$. 
Define the map $\gamma \in End_{\mathbb{C}}(V)$  by
$$\gamma(v) = (-1)^{[v]}(v),$$ 
for all homogeneous  $v \in V$.
We define the supertrace by 
$$\mathrm{str}: End_{\mathbb{C}}(V) \rightarrow \mathbb{C}, \hspace{10mm}
\mathrm{str}(x) = \mathrm{tr}(\gamma x ),$$
where $\mathrm{tr}: End_{\mathbb{C}}(V) \rightarrow \mathbb{C}$ is the standard trace functional.
\end{definition}

We write $V \otimes W$ to denote $V \otimes_{\mathbb{C}} W$
where each of $V$ and $W$ are vector spaces over $\mathbb{C}$, unless otherwise specified.

\begin{remark}
Let $V$ and $W$ be two $\mathbb{Z}_{2}$-graded vector spaces over $\mathbb{C}$, then
$V \otimes W$ has a natural $\mathbb{Z}_{2}$-gradation defined by
$$(V \otimes W)_{\gamma} = \bigoplus_{(\alpha + \beta) \equiv \gamma \pmod{2}} V_{\alpha} \otimes W_{\beta},
\hspace{10mm} \gamma \in \mathbb{Z}_{2}.$$

Let $V$ and $W$ be two $\mathbb{Z}_{2}$-graded vector spaces over $\mathbb{C}$.
The graded permutation operator $P: V \otimes W \rightarrow W \otimes V$
is defined for all homogeneous $v \in V$ and $w \in W$ by
$$P(v \otimes w) = (-1)^{[v][w]}w \otimes v,$$
and extended to all inhomogeneous elements by linearity.
\end{remark}

\begin{definition}
An associative $\mathbb{Z}_{2}$-graded algebra $A$ over $\mathbb{C}$ 
is a $\mathbb{Z}_{2}$-graded vector
space equipped with graded vector space homomorphisms $m: A \otimes A \rightarrow A$ and 
$u: \mathbb{C} \rightarrow A$ satisfying the following properties:
\begin{equation}
\label{eq:associativem}
m \circ (m \otimes \mathrm{id}) = m \circ (\mathrm{id} \otimes m),
\end{equation}
where $\mathrm{id}: A \rightarrow A$ is the identity homomorphism and
we regard both sides of (\ref{eq:associativem}) as maps $A \otimes A \otimes A \rightarrow A$, and
\begin{equation}
\label{eq:unitunitunit}
m \circ (u \otimes \mathrm{id}) = m \circ (\mathrm{id} \otimes u),
\end{equation}
where $A \rightarrow \mathbb{C} \otimes A$ and 
$A \rightarrow A \otimes \mathbb{C}$ are the natural maps and
we regard both sides of (\ref{eq:unitunitunit}) as maps $A \rightarrow A$. 
The maps $m$ and $u$ are called the multiplication and unit of $A$, respectively.
\end{definition}

A $\mathbb{Z}_{2}$-graded algebra $A$ is called commutative if 
$$m \circ (v \otimes w) = m \circ P(v \otimes w), \hspace{10mm} \forall v, w \in A.$$

In this thesis we assume that all $\mathbb{Z}_{2}$-graded algebras are associative algebras.

\begin{definition}
Let $A$ be a $\mathbb{Z}_{2}$-graded algebra with unit and 
let $V$ be a left $A$-module.  
The $A$-module $V$ is said to be $\mathbb{Z}_{2}$-graded
if the underlying vector space is $\mathbb{Z}_{2}$-graded and if
$$A_{\alpha}V_{\beta} \subseteq V_{\alpha + \beta}, 
\hspace{10mm} A_{\alpha} \in A, \hspace{10mm} \mbox{ for all } \alpha, \beta \in \mathbb{Z}_{2}.$$
\end{definition}
Right $\mathbb{Z}_{2}$-graded $A$-modules are similarly defined.
In this thesis we assume that all modules are finite dimensional left modules unless otherwise stated.

A homomorphism of $\mathbb{Z}_{2}$-graded $A$-modules is by definition a homomorphism of the
$A$-modules as well as of the underlying $\mathbb{Z}_{2}$-graded vector spaces.
Each such homomorphism is $A$-linear and homogeneous of degree $0$.
Let $V = V_{0} \oplus V_{1}$ be a $\mathbb{Z}_{2}$-graded $A$-module.
If $V' = V_{0}' \oplus V_{1}'$ is a $\mathbb{Z}_{2}$-graded $A$-module
with $V_{0}' = V_{1}$ and $V_{1}' = V_{0}$, then we regard $V$ and $V'$ as 
{\emph{not}} being isomorphic as $A$-modules.

Note that if $A$ and $B$ are two $\mathbb{Z}_{2}$-graded algebras, then 
$A \otimes B$ is a $\mathbb{Z}_{2}$-graded  algebra with the multiplication 
$$m_{A \otimes B}=(m_{A} \otimes m_{B}) \circ (\mathrm{id} \otimes P \otimes \mathrm{id}),$$
and the unit
$$u_{A \otimes B} = u_{A} \otimes u_{B}.$$
Let $V$ (resp. $W$) be a $\mathbb{Z}_{2}$-graded $A$-module (resp. $B$-module), then 
$V \otimes W$ has a unique structure as a $\mathbb{Z}_{2}$-graded $(A \otimes B)$-module given by
$$(a \otimes b)(v \otimes w) = (-1)^{[b][v]} av \otimes bw, 
\hspace{10mm} a \in A, \hspace{5mm} b \in B, \hspace{5mm} v \in V, \hspace{5mm} w \in W.$$

\begin{definition}
A $\mathbb{Z}_{2}$-graded co-algebra $C$ is a 
$\mathbb{Z}_{2}$-graded vector space $V$ together with 
$\mathbb{Z}_{2}$-graded vector space homomorphisms 
$\Delta: C \rightarrow C \otimes C$ and $\epsilon: C \rightarrow \mathbb{C}$,
where $\Delta$ satisfies
\begin{equation}
\label{eq:coassociativity}
(\Delta \otimes \mathrm{id}) \circ \Delta = (\mathrm{id} \otimes \Delta) \circ \Delta,
\end{equation}
with $\mathrm{id}: C \rightarrow C$ being the identity homomorphism 
(we regard both sides of (\ref{eq:coassociativity}) as maps $C \rightarrow C \otimes C \otimes C$) 
and $\epsilon$ and $\Delta$ satisfy
\begin{equation}
\label{eq:someking}
(\epsilon \otimes \mathrm{id}) \circ \Delta = \mathrm{id} = (\mathrm{id} \otimes \epsilon)
\circ \Delta,
\end{equation}
where we regard each side in (\ref{eq:someking}) as a map $C \rightarrow C$.
The maps $\Delta$ and $\epsilon$ are called the co-multiplication and
co-unit of $C$, respectively.  
We refer to Eq. (\ref{eq:coassociativity}) by saying that 
the co-multiplication $\Delta: C \rightarrow C \otimes C$ is co-associative.
\end{definition}
For a co-multiplication $\Delta: C \rightarrow C \otimes C$
we define the opposite co-multiplication $\Delta' : C \rightarrow C \otimes C$ by 
$\Delta' = P \circ \Delta$.
A $\mathbb{Z}_{2}$-graded co-algebra $C$ is called co-commutative if 
$\Delta(x) = \Delta'(x)$ for all $x \in C$.
A $\mathbb{C}$-linear subspace $I$ of a $\mathbb{Z}_{2}$-graded co-algebra $C$ is called a
two-sided co-ideal if
\begin{equation}
\Delta(x) \subseteq I \otimes A + A \otimes I,  
\hspace{2mm} \mbox{ and } \hspace{2mm} \epsilon(x)=0, \hspace{10mm} \forall x \in I.
\end{equation}

We inductively define the $n^{th}$ co-multiplication 
$\Delta^{(n)}: {C}^{\otimes n} \rightarrow {C}^{\otimes (n+1)}$ by
$$
\Delta^{(n)} = (\Delta \otimes \mathrm{id}^{\otimes (n-1)}) \circ \Delta^{(n-1)}, 
$$
where we fix $\Delta^{(1)}=\Delta$. 
Similarly, we inductively define the $n^{th}$ opposite co-multiplication
$\Delta'^{(n)}: {C}^{\otimes n} \rightarrow {C}^{\otimes (n+1)}$ by
$$\Delta'^{(n)}= (\Delta' \otimes \mathrm{id}^{\otimes (n-1)}) \circ \Delta'^{(n-1)},$$
where we fix $\Delta'^{(1)}=\Delta'$.

We may use {\emph{Sweedler's Sigma notation}} in writing down the co-multiplication of an element of a 
$\mathbb{Z}_{2}$-graded co-algebra $C$ \cite{sweed}.  In this notation we write
$$\Delta(x) = \sum_{(x)} x_{(1)} \otimes x_{(2)}, \hspace{10mm}
\Delta^{(2)}(x) = \sum_{(x)}  x_{(1)} \otimes x_{(2)} \otimes x_{(3)},$$
$$\Delta^{(n)}(x)=\sum_{(x)} x_{(1)} \otimes x_{(2)} \otimes \cdots \otimes x_{(n+1)}, 
\hspace{5mm} n \geq 2, \hspace{10mm} \forall x \in C.$$
A consequence of the co-associativity of the co-multiplication is that
$$\Delta^{(n)}(x) = 
(\mathrm{id}^{\otimes m} \otimes \Delta \otimes \mathrm{id}^{\otimes (n-m-1)} ) 
\circ \Delta^{(n-1)}(x), \ \ \mbox{for each }
m=0, 1, \ldots, n-1,$$
for each $x \in C$.

Observe that if $A$ and $B$ are $\mathbb{Z}_{2}$-graded co-algebras with co-multiplications 
$\Delta_{A}$, $\Delta_{B}$,  
and co-units $\epsilon_{A}$, $\epsilon_{B}$, respectively,
then $A \otimes B$ is a $\mathbb{Z}_{2}$-graded co-algebra 
with the co-multiplication 
$$\Delta_{A \otimes B} = (\mathrm{id} \otimes P \otimes \mathrm{id}) \circ
(\Delta_{A} \otimes \Delta_{B}),$$
and the co-unit
$$\epsilon_{A \otimes B} = \epsilon_{A} \otimes \epsilon_{B}.$$

\begin{definition}
Let $A$ and $B$ be $\mathbb{Z}_{2}$-graded co-algebras over $\mathbb{C}$.
A homomorphism $f: A \rightarrow B$ of the co-algebras is 
defined to be a homomorphism of the underlying $\mathbb{Z}_{2}$-graded vector spaces 
satisfying
\begin{equation}
\label{eq:bettybay}
(f \otimes f) \circ \Delta_{A} = \Delta_{B} \circ f,
\end{equation}
and 
\begin{equation}
\label{eq:pottspt}
\epsilon_{B} \circ f = \epsilon_{A},
\end{equation}
where both sides of (\ref{eq:bettybay}) are maps $A \rightarrow B \otimes B$, and 
both sides of (\ref{eq:pottspt}) are maps $A \rightarrow \mathbb{C}$.
\end{definition}

\begin{definition}
Let $A$ be a $\mathbb{Z}_{2}$-graded algebra with multiplication $m$ and unit $u$, 
and also a $\mathbb{Z}_{2}$-graded co-algebra with co-multiplication $\Delta$ and co-unit $\epsilon$. 
Then $A$ is called a $\mathbb{Z}_{2}$-graded bi-algebra if one of the following
equivalent conditions is satisfied:
\begin{itemize}
\item[(i)] $\Delta$ and $\epsilon$ are $\mathbb{Z}_{2}$-graded algebra homomorphisms,
\item[(ii)] $m$ and $u$ are $\mathbb{Z}_{2}$-graded co-algebra homomorphisms.
\end{itemize}
If, furthermore, there exists a  
$\mathbb{Z}_{2}$-graded vector space homomorphism $S: A \rightarrow A$ satisfying
\begin{equation}
\label{eq:sSSS}
m \circ(\mathrm{id} \otimes S) \circ \Delta = u \circ \epsilon = m \circ (S \otimes \mathrm{id})
\circ \Delta,
\end{equation}
then $A$ is called a $\mathbb{Z}_{2}$-graded Hopf algebra.
Such a $\mathbb{Z}_{2}$-graded vector space homomorphism is called the antipode of $A$.  
\end{definition}
The antipode $S$ of $A$ has the following properties:
\begin{itemize}
\item[(i)]    $S \circ m = m \circ P \circ (S \otimes S)$,
\item[(ii)]   $S \circ u = u$,
\item[(iii)]  $\epsilon \circ S = \epsilon$,
\item[(iv)]   $P \circ (S \otimes S) \circ \Delta = \Delta \circ S$,
\item[(v)]    if $H$ is commutative or co-commutative, $S^{2} = \mathrm{id}$.
\end{itemize}

\begin{definition}
\label{chapter1:def:biggorilla(1)}
Let $A$ be a $\mathbb{Z}_{2}$-graded Hopf algebra over $\mathbb{C}$ with antipode $S$.  
Given the $A$-module $V$, we define the dual $A$-module $V^{*}$ as follows.
Let $V^{*}=\mathrm{Hom}_{\mathbb{C}}(V, \mathbb{C})$ and
let $\langle \cdot, \cdot \rangle : V^{*} \times V \rightarrow \mathbb{C}$ 
denote the dual space pairing.
We define the action of $A$ on $V^{*}$ by
$$\langle a v^{*}, w\rangle = 
(-1)^{[a][v^{*}]} \langle v^{*}, S(a)w \rangle, \hspace{5mm} 
\mbox{for all } \hspace{2mm} a \in A, \hspace{2mm} v^{*} \in V^{*}, \hspace{2mm} w \in V.$$
\end{definition}

\end{section}

\begin{section}{$\mathbb{Z}_{2}$-graded quasitriangular Hopf algebras}
\markright{\text{$\mathbb{Z}_{2}$-graded quasitriangular Hopf algebras}}
\label{chap1:quastriangaulrhopfalgebras}

\begin{definition}
We call a $\mathbb{Z}_{2}$-graded Hopf algebra $A$ a
$\mathbb{Z}_{2}$-graded quasitriangular Hopf algebra if there exists an invertible
even element $R \in A \otimes A$ satisfying the equations
\begin{eqnarray}
R \Delta(x) & = & \Delta'(x)R, \hspace{10mm} \forall x \in A, \label{eq:asdf;lkj} \\
(\Delta \otimes \mathrm{id}) R & = & R_{13} R_{23}, \label{eq:ohgygod} \\
(\mathrm{id} \otimes \Delta)R & = & R_{13} R_{12}, \label{eq:deddy}
\end{eqnarray}
where $\mathrm{id}$ is the identity map on $A$.
Write $R = \sum_{t} \alpha_{t} \otimes \beta_{t}$, then
$R_{12} = \sum_{t} \alpha_{t} \otimes \beta_{t} \otimes 1$,
$R_{13} = \sum_{t} \alpha_{t} \otimes 1 \otimes \beta_{t}$ and
$R_{23} = \sum_{t} 1 \otimes \alpha_{t} \otimes \beta_{t}$.
The element $R$ is called the universal $R$-matrix of $A$.
\end{definition}

\begin{lemma}
Let $A$ be a $\mathbb{Z}_{2}$-graded quasitriangular Hopf algebra with universal $R$-matrix
$R \in A \otimes A$.  The universal $R$-matrix satisfies the following equations
\begin{equation}
\label{eq:gnatghat}
R_{12}R_{13}R_{23} = R_{23}R_{13}R_{12},
\end{equation}
\begin{equation}
\label{eq:quasiz2gradedequationrmatrix}
(S \otimes S)R=R, \hspace{10mm} (S \otimes 1)R=R^{-1},  \hspace{10mm} (1 \otimes S)R^{-1}=R,
\end{equation}
\begin{equation}
(\epsilon \otimes \mathrm{id})R=(\mathrm{id} \otimes \epsilon)R=1.
\end{equation}
Equation (\ref{eq:gnatghat}) is called the graded Quantum Yang-Baxter equation.
\end{lemma}
\begin{proof}
The proofs are standard (eg see \cite{ks}), for instance the proof of Eq. (\ref{eq:gnatghat}) is
$$
R_{12}R_{13}R_{23} = R_{12} \big( ( \Delta \otimes \mathrm{id}) R  \big)
                   = \big((\Delta' \otimes \mathrm{id})R\big)R_{12}  
                   = R_{23}R_{13}R_{12}.
$$
The proofs of the other equations similarly follow 
those of the corresponding equations for ungraded quasitriangular Hopf algebras \cite{ks}.
\end{proof}

Let $A$ be a $\mathbb{Z}_{2}$-graded quasitriangular Hopf algebra over 
$\mathbb{C}$ with universal $R$-matrix $R$.  
Let $V$ and $W$ be finite dimensional $A$-modules and let $\pi_{V}$, $\pi_{W}$ be the 
representations of $A$ afforded by $V$ and $W$, respectively. 
Fix $R_{VW} = (\pi_{V} \otimes \pi_{W})R$ and let 
$\check{R}_{VW} \in \mathrm{Hom}_{\mathbb{C}}(V \otimes W, W \otimes V)$  be the map defined by
$$\check{R}_{VW}(v \otimes w) = P \circ \big(R_{VW}(v \otimes w)\big), 
\hspace{5mm} \mbox{for all } v \in V, \ w \in W.$$
The actions of $R_{VW}$ and $\check{R}_{VW}$ on $v \otimes w \in V \otimes W$ are respectively
$$R_{VW} (v \otimes w) = \sum_{t} \pi_{V}(\alpha_{t})v \otimes \pi_{W}(\beta_{t})w  (-1)^{[\beta_{t}][v]},$$
$$
\check{R}_{VW} (v \otimes w) =
\sum_{t} \pi_{W}(\beta_{t})w \otimes \pi_{V}(\alpha_{t})v  (-1)^{[\alpha_{t}] + [w]([v] + [\alpha_{t}])}.$$

\begin{lemma}
\label{lem:113}
Let $A$ be a $\mathbb{Z}_{2}$-graded quasitriangular Hopf algebra with universal $R$-matrix $R$.
Let $V$ be a finite dimensional $A$-module and let $\pi$ 
be the representation of $A$ afforded by $V$.
For each integer $t \geq 2$, define
$$\check{R}_{i} = \mathrm{id}^{\otimes (i-1)} \otimes \check{R}_{VV} \otimes
\mathrm{id}^{\otimes (t-(i+1))} \in End_{\mathbb{C}}(V^{\otimes t}),$$
$$P_{i, i+1} = \mathrm{id}^{\otimes (i-1)} \otimes P \otimes
\mathrm{id}^{\otimes (t-(i+1))} \in End_{\mathbb{C}}(V^{\otimes t}) ,$$
for each $1 \leq i \leq t-1$, 
where  $\mathrm{id} = \mathrm{id}_{V}$, then
\begin{itemize}
\item[(i)]     $\check{R}_{i}$ is an element of the centraliser algebra 
               ${\cal{L}}_{t} = End_{A}(V^{\otimes t})$, for each $1 \leq i \leq t-1$,
\item[(ii)]    $\check{R}_{i} \check{R}_{j} = \check{R}_{j} \check{R}_{i}$, 
               \hspace{29mm} for all $|i-j| > 1$,
\item[(iii)]   $\check{R}_{i} \check{R}_{i+1} \check{R}_{i} = \check{R}_{i+1} \check{R}_{i}
                \check{R}_{i+1}$, \hspace{10mm} for all $1 \leq i \leq t-2$.
\end{itemize}
\end{lemma}
\begin{proof}
\begin{itemize}
\item[(i)] 
The proof is standard.  For an arbitrary $x \in A$, we have
\begin{eqnarray*}
\check{R}_{VV} \big( (\pi \otimes \pi) \Delta(x) \big)
& = & P \big((\pi \otimes \pi) R \Delta(x) \big) \\
& = & P \big((\pi \otimes \pi) \Delta'(x) R\big) \\
& = & \big((\pi \otimes \pi) \Delta(x)\big) \big(P (\pi \otimes \pi) R \big) \\
& = & \big((\pi \otimes \pi) \Delta(x)\big) \check{R}_{VV}.
\end{eqnarray*}
\item[(ii)]  This is obvious.
\item[(iii)] 
We consider the case $i=1$, the other cases are similar.  We write $R_{i j}$ to mean 
$(\pi_{V} \otimes \pi_{V} \otimes \pi_{V}) R_{i j}$.  Now
\begin{eqnarray*}
\check{R}_{1} \check{R}_{2} \check{R}_{1} & = &
P_{12}  R_{12}  P_{23}  R_{23}  P_{12}  R_{12}  \\
& = & P_{12}  P_{23}  R_{13}  R_{23} P_{12} R_{12} \\
& = & P_{12}  P_{23}  P_{12}  R_{23}R_{13}R_{12}, 
\end{eqnarray*}
where we have used the easily verified results:  
$R_{12} P_{23} = P_{23} R_{13}$ and $R_{13}R_{23}P_{12} = P_{12}R_{23}R_{13}$.
A similar calculation shows that
$$
\check{R}_{2}\check{R}_{1}\check{R}_{2} 
= P_{23}  P_{12}  P_{23}  R_{12}R_{13}R_{23}, 
$$
The equality $P_{12} P_{23} P_{12}=P_{23} P_{12} P_{23}$ 
then yields what we seek to prove.
\end{itemize}
\end{proof}

\noindent
From \cite{lr} we obtain the following lemma.
\begin{lemma}
\label{lem:fugue}
Let $A$ be a $\mathbb{Z}_{2}$-graded quasitriangular Hopf algebra with universal $R$-matrix 
$R = \sum_{t} a_{t} \otimes b_{t}$, and
let $R^{T} = \sum_{t} b_{t} \otimes a_{t} (-1)^{[a_{t}][b_{t}]}$.
Let $V$ and $W$ be finite dimensional irreducible $A$-modules and let $\pi_{V}$, $\pi_{W}$ be the
representations of $A$ afforded by $V$ and $W$, respectively.  Then
\begin{itemize}
\item[(i)]  $(\pi_{V} \otimes \pi_{W})\big(R^{T}R\big) \in End_{A}(V \otimes W)$, and
\item[(ii)] for each integer $t \geq 3$, the element
$$
\big(\pi_{V}^{\otimes t} \otimes \pi_{V}\big)(\Delta^{(t-1)} \otimes \mathrm{id})\big(R^{T}R\big) = 
\check{R}_{t} \check{R}_{t-1} \cdots \check{R}_{1}
\check{R}_{1} \cdots \check{R}_{t-1} \check{R}_{t} \in End_{\mathbb{C}}(V^{\otimes (t+1)}),
$$
is an element of $End_{A}(V^{\otimes (t+1)})$.
\end{itemize}
\end{lemma}
\begin{proof}
(i)  Let $x \in A$ be arbitrary.
Applying $P$ to the equation $R \Delta(x)  = \Delta'(x) R$ gives
$R^{T} \Delta'(x) = \Delta(x) R^{T}$.
Hence $\Delta(x) R^{T}R = R^{T} \Delta'(x)R = R^{T}R \Delta(x)$.

(ii) 
We introduce some notation we will use later on.  
For each integer $t \geq 3$, let $i, j \in \mathbb{N}$ satisfy $1 \leq i < j \leq t$, and fix
$$R_{ij} = \sum_{t} \mathrm{id}^{\otimes (i-1)} \otimes a_{t} \otimes \mathrm{id}^{(j-1-i)} \otimes b_{t} \otimes
\mathrm{id}^{(t-j)} \in A^{\otimes t},$$
$$R_{ji} = \sum_{t} \mathrm{id}^{\otimes (i-1)} \otimes b_{t} \otimes \mathrm{id}^{(j-1-i)} \otimes a_{t} \otimes
\mathrm{id}^{(t-j)} (-1)^{[a_{t}][b_{t}]} \in A^{\otimes t}.$$
We inductively prove that 
$(\Delta^{(t-2)} \otimes \mathrm{id}) R_{12} = R_{1t}R_{2t} \cdots R_{(t-1)t}$ for each $t \geq 3$.
The induction hypothesis is true for $t=3$ from (\ref{eq:ohgygod}) 
and assume that it is true for some $t \geq 3$, then
\begin{eqnarray*}
(\Delta^{(t-1)} \otimes \mathrm{id})R & = & (\Delta \otimes \mathrm{id}^{\otimes (t-1)})
(\Delta^{(t-2)} \otimes \mathrm{id})R \\
& = & (\Delta \otimes \mathrm{id}^{\otimes (t-1)}) R_{1t}R_{2t} \cdots R_{(t-1)t} \\
& = & R_{1(t+1)}R_{2(t+1)} \cdots R_{t(t+1)}.
\end{eqnarray*}
In a similar way, we can prove that
$$
(\Delta^{(t-1)} \otimes \mathrm{id})R^{T} 
= R_{(t+1)t} R_{(t+1)(t-1)} \cdots R_{(t+1)1}.
$$

Let $\widetilde{P} = P_{12} P_{23} \cdots P_{t(t+1)}$.  Then the above formulae lead to
$$\widetilde{P} \big(\pi^{\otimes (t+1)} ( \Delta^{(t-1)} \otimes \mathrm{id})R \big) =  
                \check{R}_{1} \check{R}_{2} \cdots \check{R}_{t},$$
$$\pi^{\otimes (t+1)} \big( ( \Delta^{(t-1)} \otimes \mathrm{id}) R^{T} \big) \widetilde{P}^{-1} = 
		\check{R}_{t} \check{R}_{t-1} \cdots \check{R}_{1}.$$
Combining these equations together completes the proof.

\end{proof}

\begin{lemma}
\label{lem:theuoperator}
Let $A$ be a $\mathbb{Z}_{2}$-graded quasitriangular Hopf algebra with universal $R$-matrix 
$R=\sum_{t} a_{t} \otimes b_{t}$, and let $u \in A$ be defined by 
$u = m \circ P \circ (\mathrm{id} \otimes S) R$, that is
$$u = \sum_{t}(-1)^{\left[a_{t}\right]}S(b_{t})a_{t}.$$ 
Then $u$ is invertible, and
$$\epsilon(u)=1, \hspace{10mm} \Delta(u)=\left(u \otimes u\right)\left(R^{T}R\right)^{-1}.$$
Furthermore, $$S^{2}(x)=uxu^{-1}, \hspace{10mm}  \mbox{for all } x \in A.$$
\end{lemma}
\begin{proof}  
The proofs are standard: we follow the proofs of the corresponding equations 
for ungraded quasitriangular Hopf algebras \cite{d,ks}.
From Lemma \ref{lemma:uuuuuuuuuu}, $S^{2}(x)u=ux$ for all $x \in A$.
Let $\widetilde{u} = \sum_{s} S^{-1}(d_{s})c_{s}(-1)^{[d_{s}]}$ 
where $R^{-1}=\sum_{s}c_{s} \otimes d_{s}$; 
we will show that $\widetilde{u}$ is a two-sided inverse of $u$.  Firstly
\begin{eqnarray*}
u \widetilde{u} & = & \sum_{s} u S^{-1}(d_{s})c_{s}(-1)^{[d_{s}]}  \\
& = & \sum_{s}S(d_{s})u c_{s}(-1)^{[d_{s}]}  \\
& = & \sum_{s,t} S(b_{t} d_{s})a_{t} c_{s} (-1)^{[d_{s}]+[b_{t}]+[d_{s}][b_{t}]}. 
\end{eqnarray*} 
Applying the map
$m \circ P \circ (1 \otimes S)$ to
$R R^{-1} = \sum_{s,t} a_{t}c_{s} \otimes b_{t} d_{s} (-1)^{[b_{t}][c_{s}]}$ gives
$$\sum_{s,t} S(b_{t} d_{s}) a_{t} c_{s} (-1)^{[b_{t}] + [d_{s}] + [a_{t}][d_{s}]}=1.$$
Hence $u \widetilde{u} = 1$.

In exactly the same way we can also show that $\widetilde{u}$ is a left inverse of $u$,
thus $\widetilde{u}$ is a two-sided inverse of $u$ and we write $\widetilde{u} = u^{-1}$.
We now show that $\epsilon(u)=1$: firstly
$$\epsilon(u) = \sum_{t} \epsilon\big(S(b_{t})\big) \epsilon(a_{t}) (-1)^{[a_{t}][b_{t}]}
 = \sum_{t} \epsilon(b_{t}) \epsilon(a_{t}) = \sum_{t} \epsilon(a_{t}) \epsilon(b_{t}),$$
as $\epsilon \circ S = \epsilon$ and $\epsilon(x)=0$ if $[x]=1$.
Applying $\epsilon$ to the equation $(\epsilon \otimes \mathrm{id})R=1$  gives 
$$\epsilon(u)= \sum_{t} \epsilon(a_{t}) \epsilon(b_{t})=1.$$ 

The proof of the equation $\Delta(u)=(u \otimes u)\left(R^{T}R\right)^{-1}$  
is given explicitly in \cite{zg} and will not be repeated.

\end{proof}

\begin{lemma}
\label{lemma:uuuuuuuuuu}
Let $A$ be a $\mathbb{Z}_{2}$-graded quasitriangular Hopf algebra with universal $R$-matrix
$R=\sum_{t} a_{t} \otimes b_{t}$, then $S^{2}(x)u=ux$ for all $x \in A$.
\end{lemma}
\begin{proof}
The proof is well known.  But as this is important to us later, we give a proof here.
For each $x \in A$,
$$\sum_{(x)} R \Delta(x_{(1)}) \otimes x_{(2)}=\sum_{(x)}\Delta'(x_{(1)})R \otimes x_{(2)},$$
which we rewrite explicitly as
\begin{eqnarray*}
\lefteqn{
\sum_{t, (x)} a_{t} x_{(1)} \otimes b_{t}x_{(2)} \otimes x_{(3)} (-1)^{[b_{t}][x_{(1)}]} }  \\
& & \hspace{10mm} = \sum_{t,(x)} x_{(2)}a_{t} \otimes x_{(1)}b_{t} \otimes 
     x_{(3)} (-1)^{[x_{(1)}][x_{(2)}]+[x_{(1)}][a_{t}]}.
\end{eqnarray*}
Using this we obtain
\begin{eqnarray}
\lefteqn{
\sum_{t,(x)} x_{(3)} \otimes b_{t}x_{(2)} \otimes a_{t}x_{(1)}
(-1)^{[x_{(3)}]([x_{(1)}]+[x_{(2)}])+([a_{t}]+[x_{(1)}])([b_{t}]+[x_{(2)}])+[b_{t}][x_{(1)}] } } \nonumber \\
& & = \sum_{t,(x)} x_{(3)} \otimes x_{(1)}b_{t} \otimes x_{(2)}a_{t}
(-1)^{[x_{(3)}]([x_{(1)}] + [x_{(2)}]) + ([a_{t}]+[x_{(2)}])([b_{t}]+[x_{(1)}])} \nonumber \\
& & \hspace{15mm} \times (-1)^{[x_{(1)}][x_{(2)}] + [a_{t}][x_{(1)}]}. \label{eq:Suequalssu}
\end{eqnarray}
Note that $R$ is even, thus $\left([a_{t}]+[b_{t}]\right) \equiv 0 \pmod{2}$ for all $t$.  
In addition, $\left([x_{(1)}]+[x_{(2)}]+[x_{(3)}]\right) \equiv [x] \pmod{2}$.  
Using these facts we can simplify the sign factors in the equation considerably.  
Applying the map $(\mathrm{id} \otimes m) \circ (\mathrm{id} \otimes S \otimes \mathrm{id})$ 
to both sides of (\ref{eq:Suequalssu}) gives
\begin{eqnarray}
\lefteqn{
\sum_{(x)} x_{(3)} \otimes S(x_{(2)})ux_{(1)} (-1)^{[x_{(3)}]([x]+[x_{(3)}])+[x_{(1)}][x_{(2)}]}} \nonumber \\
& & \hspace{12mm} = \sum_{t,(x)} x_{(3)} \otimes S(b_{t})S(x_{(1)})x_{(2)}a_{t} (-1)^{([a_{t}]+[x_{(3)}])([x]+[x_{(3)}])+[a_{t}]} \nonumber \\
& & \hspace{12mm} = \sum_{t,(x)} x_{(2)} \otimes S(b_{t})a_{t} \epsilon(x_{(1)}) (-1)^{([a_{t}] + [x_{(2)}])([x]+[x_{(2)}])+[a_{t}]} \nonumber \\
& & \hspace{12mm} = x \otimes u. \label{eq:sadsong}
\end{eqnarray}
Applying the map $m \circ (S^{2} \otimes \mathrm{id})$ to (\ref{eq:sadsong}) gives
\begin{eqnarray*}
S^{2}(x)u & = & \sum_{(x)} S(x_{(2)}S(x_{(3)})) u x_{(1)} 
(-1)^{[x_{(1)}]([x_{(3)}]+[x_{(2)}])} \\
& = & \sum_{(x)} S(x_{(2)}S(x_{(3)})) u x_{(1)} (-1)^{ [x_{(1)}]([x]+[x_{(1)}])} \\
& = & ux.
\end{eqnarray*}
\end{proof}

\noindent
We can readily prove the following lemma using induction.
\begin{lemma}
\label{lemam:voperator}
Let $A$ be a $\mathbb{Z}_{2}$-graded quasitriangular Hopf algebra with 
universal $R$-matrix $R$ and the element $u = m \circ P \circ ( \mathrm{id} \otimes S)R$, then
$$\Delta^{(i)}(u) = u^{\otimes (i+1)}\big( (R^{T}R)^{-1} \otimes \mathrm{id}^{\otimes (i-1)} \big) \prod_{j=1}^{i-1}
\big((\Delta^{(j)} \otimes \mathrm{id}^{\otimes (i-j)}) \left(R^{T}R\right)^{-1}\big),
\hspace{5mm} \mbox{for all } i \geq 1,$$
where we fix the product to be equal to $1 \otimes 1$ if $i=1$.
\end{lemma}

\begin{definition}
A $\mathbb{Z}_{2}$-graded quasitriangular Hopf algebra $A$ is called a
$\mathbb{Z}_{2}$-graded ribbon Hopf algebra if it is equipped with an 
invertible central even element $v \in A$ satisfying
\begin{equation}
v^{2} = u S(u), \hspace{10mm} S(v)=v, \hspace{10mm} \epsilon(v)=1,
\end{equation}
\begin{equation}
\Delta(v) = (v \otimes v)\left(R^{T}R\right)^{-1}.
\end{equation}
\end{definition}

\begin{definition}
Let $A$ be a $\mathbb{Z}_{2}$-graded ribbon Hopf algebra over $\mathbb{C}$.
Let $V$ be an $A$-module, $\pi$ the representation of $A$ afforded by $V$, and 
$f$ an element of $End_{\mathbb{C}}(V)$.
The quantum supertrace of $f$ is defined to be
$$str_{q}(f) = \mathrm{str}\left(\pi(v^{-1} u) \circ f \right).$$

The quantum superdimension of $V$ is defined to be
the quantum supertrace of the identity endomorphism on $V$:
$$sdim_{q}(V) = \mathrm{str}\left(\pi(v^{-1} u)\right).$$
\end{definition}

\end{section}

\end{chapter}

%% file: chapter2100.tex
\begin{chapter}{Quantum  $osp(1|2n)$ at generic $q$}
\pagestyle{myheadings}
\label{chap2:titlelabel}
\markboth{\text{Chapter \ref{chap2:titlelabel}.
Quantum  $osp(1|2n)$ at generic $q$}}
{\text{ }}

In this chapter we introduce the quantum superalgebra $U_{q}(osp(1|2n))$ over $\mathbb{C}$ and discuss
its finite dimensional irreducible representations.  
We define $R$-matrices for representations of $U_{q}(osp(1|2n))$, and construct projections from
tensor powers of the fundamental irreducible $U_{q}(osp(1|2n))$-module $V$ onto irreducible
$U_{q}(osp(1|2n))$-submodules of $V^{\otimes t}$.
By using matrix units in the Birman-Wenzl-Murakami algebra $\mathscr{BW}_{t}(-q^{2n},q)$,
we construct $U_{q}(osp(1|2n))$-linear
maps between isomorphic irreducible $U_{q}(osp(1|2n))$-submodules of $V^{\otimes t}$.

The structure of this chapter is as follows.
In Section \ref{sec:algebraschapter2} we introduce  
the quantum superalgebra $U_{q}(osp(1|2n))$.
In Section \ref{subsec:100} we discuss finite dimensional irreducible
$U_{q}(osp(1|2n))$-modules, including the fundamental irreducible module $V$, and show that
$V^{\otimes t}$ is completely reducible.
In Section \ref{subsec:RRRR} we introduce $R$-matrices for representations of
$U_{q}(osp(1|2n))$.
In Section \ref{sec:'central'element} we investigate the properties of two useful elements of a completion
$\overline{U}^{+}_{q}(osp(1|2n))$ of $U_{q}(osp(1|2n))$.
In Section \ref{sec:spectam(a)} we determine the spectral decomposition of $\check{\cal{R}}_{V,V}$.
In Section \ref{eq:theXfactor(a)} we show that there is
a representation of  $\mathscr{BW}_{t}(-q^{2n},q)$ in an algebra 
generated by the $\check{\cal{R}}_{V,V}$-matrices 
acting on the $i^{th}$ and $(i+1)^{st}$ tensor powers of
$V^{\otimes t}$ for $i=1, \ldots, t-1$.
In Section \ref{eq:caseydonovan(a)} we recall Bratteli diagrams and path algebras.
In Section \ref{subsec:projectontome} we construct 
projections from $V^{\otimes t}$
onto irreducible $U_{q}(osp(1|2n))$-submodules of $V^{\otimes t}$.
In Section \ref{subsec:rhapsodyinred} we construct matrix units in
$End_{U_{q}(osp(1|2n))}(V^{\otimes t})$ from  matrix units in $\mathscr{BW}_{t}(-q^{2n},q)$
and prove that the algebra $End_{U_{q}(osp(1|2n))}(V^{\otimes t})$ 
is generated by the $\check{\cal{R}}_{V,V}$-matrices.

\begin{section}{The quantum superalgebra $U_{q}(osp(1|2n))$}
\label{sec:algebraschapter2}
\markright{\text{The quantum superalgebra $U_{q}(osp(1|2n))$}}

In this section we introduce the algebra that is at the core of this thesis:
the quantum superalgebra $U_{q}(osp(1|2n))$ \cite{z1,kt,y}.

Let us begin by describing the root system of $osp(1|2n)$.
Let $H^{*}$ be a vector space over $\mathbb{C}$ with a basis
$\left\{ \epsilon_{i} | \ 1 \leq i \leq n \right\}$ and let 
\begin{equation}
\label{eq:starstar(1)}
( \cdot, \cdot): H^{*} \times H^{*} \rightarrow \mathbb{C},
\end{equation}
be a non-degenerate bilinear form defined by
$(\epsilon_{i}, \epsilon_{j})=\delta_{i,j}$.  

The set of simple roots of $osp(1|2n)$ is $\left\{ \alpha_{i} | \ 1 \leq i \leq n \right\}$ where
\begin{eqnarray*}
\alpha_{i} & = & \left\{
\begin{array}{ll}
\epsilon_{i}-\epsilon_{i+1}, &  i=1, \ldots, n-1,   \\
\epsilon_{n},                &  i=n,
\end{array} \right.
\end{eqnarray*}
which forms another basis of $H^{*}$.

Let $\Phi^{+}$ denote the set of the positive roots of $osp(1|2n)$, we have
$$\Phi^{+}=\left\{
\epsilon_{i} \pm \epsilon_{j}, \epsilon_{k}, 2\epsilon_{k}| \ \ 1 \leq i < j \leq n, \ \ 1 \leq k \leq n \right\}.$$
We further define the subsets $\Phi^{+}_{0}, \Phi^{+}_{1} \subset \Phi^{+}$ and
$\overline{\Phi}_{0}^{+} \subset \Phi^{+}_{0}$ by
$$\begin{array}{ll}
\Phi^{+}_{0} = \{\epsilon_{i} \pm \epsilon_{j}, 2\epsilon_{k}| \ \ 1 \leq i < j \leq n, \ \ 1 \leq k \leq n \}, &
\hspace{10mm} \Phi^{+}_{1} = \{ \epsilon_{k} | \ \ 1 \leq k \leq n\}, \\
\overline{\Phi}_{0}^{+} = \{\alpha \in \Phi^{+}_{0} | \ \ \alpha/2 \notin \Phi^{+}_{1}\}.  &
\end{array}$$
We call $\Phi^{+}_{0}$ (resp. $\Phi^{+}_{1}$) the set of 
{\emph{positive even roots}} (resp. {\emph{positive odd roots}}).
The set of negative roots of $osp(1|2n)$ is
$\Phi^{-} = -\Phi^{+}$, and $\Phi = \Phi^{+} \cup \Phi^{-}$ is the set of all the roots.

We denote by $2\rho \in H^{*}$ the graded sum of the positive roots of $osp(1|2n)$:
$$2 \rho = \sum_{\alpha \in \Phi^{+}_{0}} \alpha - \sum_{\beta \in \Phi^{+}_{1}} \beta.$$
Explicitly, $2\rho=\sum_{i=1}^{n}(2n-2i+1) \epsilon_{i}$. 
The element $2\rho$ satisfies
$\left(2 \rho,\alpha_{i}\right)=(\alpha_{i},\alpha_{i})$ for each $1 \leq i \leq n$, 
and will play an important role in this thesis.

Let $A=\left(a_{ij}\right)_{i,j=1}^{n}$ be the Cartan matrix of $osp(1|2n)$.
The components of $A$ are defined by
$a_{ij}=2\left(\alpha_{i}, \alpha_{j}\right) / \left(\alpha_{i}, \alpha_{i}\right)$, and we have explicitly
$$A=\left( \begin{array}{rrrrrr}
2 & -1 & 0 & \cdots & 0 & 0 \\
-1 & 2 & -1 & \cdots & 0 & 0 \\
0 & -1 & 2 & \cdots & 0 & 0 \\
\vdots & \vdots & \vdots & \ddots & \vdots & \vdots \\
0 & 0 & 0 & \cdots & 2 & -1 \\
0 & 0 & 0 & \cdots & -2 & 2
\end{array} \right).$$

The Lie superalgebra $\mathfrak{g}=osp(1|2n)$ over $\mathbb{C}$ 
can be defined in terms of a Serre presentation
with generators $\{E_{i}, F_{i}, H_{i} | \ 1 \leq i \leq n \}$
subject to the relations
$$\left[E_{i},F_{j}\right]=\delta_{ij}H_{i},  \hspace{10mm} \left[H_{i},H_{j}\right]=0, \hspace{10mm}  \forall i,j, $$
$$
\left[H_{i},E_{j}\right]=(\alpha_{i}, \alpha_{j}) E_{j}, 
\hspace{10mm} \left[H_{i},F_{j}\right]=-(\alpha_{i}, \alpha_{j})F_{j}, \hspace{10mm} \forall i,j, 
$$
\begin{equation}
\label{eq:heylittlepiggy(1)}
(\mbox{ad} \ E_{i})^{1-a_{ij}}E_{j}=0, \hspace{10mm}  (\mbox{ad} \ F_{i})^{1-a_{ij}}F_{j}=0, 
\hspace{10mm}  i \neq j,
\end{equation}
where $[ \cdot, \cdot ]$ represents the $\mathbb{Z}_{2}$-graded Lie bracket and 
$(\mbox{ad} \ a)b= [a,b]$.
The $\mathbb{Z}_{2}$-grading of the generators is
$$[E_{i}] = [F_{i}] = [H_{j}] = 0, \hspace{10mm} [E_{n}] = [F_{n}] = 1, 
\hspace{10mm}  1 \leq i \leq n-1, \hspace{10mm} 1 \leq j \leq n. $$
\noindent

The universal enveloping algebra $U(\mathfrak{g})$ of $\mathfrak{g}$
 is a unital associative 
$\mathbb{Z}_{2}$-graded algebra,
which may be considered as being generated by 
$\{E_{i}, F_{i}, H_{i} | \ 1 \leq i \leq n \}$ subject to relations that are formally the same as
(\ref{eq:heylittlepiggy(1)})
but with the bracket $[ \cdot, \cdot]$ interpreted as a
$\mathbb{Z}_{2}$-graded commutator
$[ \cdot, \cdot ]: U(\mathfrak{g}) \times U(\mathfrak{g}) \rightarrow U(\mathfrak{g})$ defined by
\begin{equation}
\label{eq:starstar(2)}
[X, Y] = XY -(-1)^{[X][Y]} YX.
\end{equation}
If each of two elements $X,Y \in U(\mathfrak{g})$ has a grading, 
then the grading of $XY \in U(\mathfrak{g})$ is defined by
$$[XY]= \big( [X] + [Y] \big) \bmod{2}.$$
The graded commutator $[X, Y]$ of any two homogeneous elements of $U(\mathfrak{g})$ is defined by
(\ref{eq:starstar(2)}) and is extended to inhomogeneous elements of $U(\mathfrak{g})$ by linearity.
Note that $U(\mathfrak{g}) \otimes U(\mathfrak{g})$ has a natural 
$\mathbb{Z}_{2}$-graded associative algebra structure, with the grading defined for
homogeneous $X, Y \in U(\mathfrak{g})$ by
$$[X \otimes Y] = \big( [X] + [Y] \big) \bmod{2}.$$

The universal enveloping algebra $U(\mathfrak{g})$ has the structures of a $\mathbb{Z}_{2}$-graded Hopf algebra.
The co-multiplication $\Delta$, co-unit $\epsilon$ and antipode $S$ are defined 
on each generator $X \in \mathfrak{g} \hookrightarrow U(\mathfrak{g})$ by
$$\Delta(X) = X \otimes 1 + 1 \otimes X, \hspace{10mm}
\epsilon(X)=0, \hspace{10mm} \epsilon(1)=1, \hspace{10mm} S(X)=-X.$$

The quantum superalgebra $U_{q}(\mathfrak{g})$ is some ``$q$-deformation'' of $U(\mathfrak{g})$.
We describe its Jimbo version here.

\begin{definition}
The quantum superalgebra $U_{q}(\mathfrak{g})$ over $\mathbb{C}$, 
in the sense of Jimbo, is an associative $\mathbb{Z}_{2}$-graded unital
algebra  generated by the elements
$\{e_{i},f_{i},K_{i},K_{i}^{-1}| \ 1 \leq i \leq n\}$ subject to the relations
$$\left[e_{i},f_{j}\right]= \delta_{ij}\frac{K_{i}-K_{i}^{-1}}{q-q^{-1}},$$
$$K_{i}e_{j}K_{i}^{-1}= q^{\left(\alpha_{i}, \alpha_{j}\right)}e_{j}, \hspace{10mm}
K_{i}f_{j}K_{i}^{-1}=q^{-\left(\alpha_{i}, \alpha_{j}\right)}f_{j},$$
$$\big[K_{i}^{\pm 1},K_{j}^{\pm 1}\big]=\big[K_{i}^{\pm 1},K_{j}^{\mp 1}\big]=0, 
\hspace{10mm} K_{i}^{\pm 1}K_{i}^{\mp 1}=1,$$
\begin{equation}
\label{equat:two}
(ad_{q}e_{i})^{1-a_{ij}}e_{j}=0, \hspace{10mm} (ad_{q}f_{i})^{1-a_{ij}}f_{j}=0, \hspace{10mm} \forall i,j,
\end{equation}
where $0 \neq q \in \mathbb{C}$ and $q^{2} \neq 1$, and the adjoint actions are defined by
$$(ad_{q}e_{i})X  = e_{i}X -(-1)^{[e_{i}][X]} K_{i}      X K_{i}^{-1} e_{i},$$
$$(ad_{q}f_{i})X =  f_{i}X -(-1)^{[f_{i}][X]} K_{i}^{-1} X K_{i}      f_{i},$$
for all $X \in U_{q}(\mathfrak{g})$.  
In (\ref{equat:two}), the bracket $[ \cdot, \cdot]$ is as defined in Eq. (\ref{eq:starstar(2)}).
\end{definition}

As is well known,
there exists a $\mathbb{Z}_{2}$-graded Hopf algebra structure on $U_{q}(\mathfrak{g})$ with 
the co-multiplication $\Delta$, the co-unit $\epsilon$, and the antipode $S$ defined on each generator by
$$
\Delta(e_{i}) =  e_{i} \otimes K_{i} + 1 \otimes e_{i},   \hspace{5mm}
\Delta(f_{i})  = f_{i} \otimes 1 + K_{i}^{-1} \otimes f_{i}, 
\hspace{5mm} \Delta\left(K_{i}^{\pm 1}\right) =  K_{i}^{\pm 1} \otimes K_{i}^{\pm 1}.$$
$$\epsilon(e_{i}) = \epsilon(f_{i}) =  0, \hspace{5mm} \epsilon\left(K_{i}^{\pm 1}\right)=\epsilon(1) = 1.$$
$$ S(e_{i}) = -e_{i}K^{-1}_{i},      \hspace{5mm}
S(f_{i}) = -K_{i}f_{i},           \hspace{5mm}
S(K_{i}^{\pm 1}) = K_{i}^{\mp 1}.$$

There are a number of subalgebras of $U_{q}(\mathfrak{g})$ which will be quite useful.
We define $U_{q}(\mathfrak{b}_{+})$ to be the subalgebra of $U_{q}(\mathfrak{g})$ 
generated by $\left\{ e_{i}, K^{\pm 1}_{i} | \ 1 \leq i \leq n \right\}$, and 
$U_{q}(\mathfrak{b}_{-})$ to be the subalgebra of $U_{q}(\mathfrak{g})$ 
generated by $\left\{ f_{i}, K^{\pm 1}_{i} | \ 1 \leq i \leq n \right\}$.

We note that the action of $S^{-1}$ on each generator of $U_{q}(\mathfrak{g})$ is
$$S^{-1}(e_{i}) = -K_{i}^{-1}e_{i}, \hspace{10mm} S^{-1}(f_{i}) = -f_{i}K_{i}, \hspace{10mm}
S^{-1}(K_{i}^{\pm 1}) = K_{i}^{\mp 1}.$$

Finally, for each $\beta = \sum_{i=1}^{n} m_{i} \alpha_{i}$ where $m_{i} \in \mathbb{Z}$, we define
$K_{\beta} = \prod_{i=1}^{n} \left( K_{i} \right)^{m_{i}}$.

The quantum superalgebra $U_{h}(\mathfrak{g})$, in the sense of Drinfel'd, 
is a $\mathbb{Z}_{2}$-graded algebra over the ring
$\mathbb{C}[[h]]$ for an indeterminate $h$, 
completed with respect to the $h$-adic topology \cite{kt}.
The $\mathbb{Z}_{2}$-graded algebra $U_{h}(\mathfrak{g})$ 
is generated by $\{E_{i}, F_{i}, H_{i} | \ 1 \leq i \leq n \}$ subject to the relations
$$\left[E_{i},F_{j}\right] = \delta_{ij}\frac{e^{h H_{i}} - e^{-h H_{i}}}{e^{h}-e^{-h}},   
\hspace{10mm} [H_{i}, H_{j}]=0,$$
$$[H_{i}, E_{j}]=(\alpha_{i}, \alpha_{j}) E_{j}, \hspace{10mm}
[H_{i}, F_{j}]=-(\alpha_{i}, \alpha_{j}) F_{j},$$
$$(ad_{q}E_{i})^{1-a_{ij}}E_{j}=0, \hspace{10mm} (ad_{q}F_{i})^{1-a_{ij}}F_{j}=0, \hspace{10mm} \forall i,j,$$
where the adjoint functions are defined by
$$(ad_{q}E_{i})X  = E_{i}X-(-1)^{[E_{i}][X]} e^{h H_{i}}  X e^{-h H_{i}} E_{i},$$
$$(ad_{q}F_{i})X  = F_{i}X-(-1)^{[F_{i}][X]} e^{-h H_{i}} X e^{h H_{i}}  F_{i},$$
where we fix $q = e^{h}$.

There is a $\mathbb{Z}_{2}$-graded Hopf algebra structure on $U_{h}(\mathfrak{g})$ with the
co-multiplication 
$\Delta: U_{h}(\mathfrak{g}) \rightarrow U_{h}(\mathfrak{g}) \otimes U_{h}(\mathfrak{g})$, 
the co-unit $\epsilon$, and the antipode $S$ defined on each generator by
$$
\Delta(E_{i}) =  E_{i} \otimes e^{h H_{i}} + 1 \otimes E_{i},   \hspace{5mm}
\Delta(F_{i}) = F_{i} \otimes 1 + e^{-h H_{i}} \otimes F_{i}, \hspace{5mm} 
\Delta(H_{i}) =  H_{i} \otimes 1 + 1 \otimes H_{i}.$$
$$\epsilon(E_{i}) = \epsilon(F_{i}) = \epsilon(H_{i}) = 0, \hspace{10mm} \epsilon(1)=1,$$
$$ S(E_{i}) = -E_{i}e^{-h H_{i}}, \hspace{10mm} S(F_{i}) = -e^{h H_{i}}F_{i},  \hspace{10mm} S(H_{i}) = -H_{i}.$$

The quantum superalgebra $U_{h}(\mathfrak{g})$ admits a universal $R$-matrix \cite{kt} and is a
$\mathbb{Z}_{2}$-graded ribbon Hopf algebra \cite{zg}.

\end{section}

\begin{section}{Finite dimensional irreducible $U_{q}(osp(1|2n))$-modules}
\label{subsec:100}
\markright{\text{Finite dimensional irreducible $U_{q}(osp(1|2n))$-modules}}

Throughout this section we assume that $q$ is non-zero and not a root of unity.
In this case the representation theory of $U_{q}(\mathfrak{g})$
is completely understood \cite{zsuper,zou}.

For a $\mathbb{Z}_{2}$-graded algebra $A$, we will write
$V_{\lambda}$ in this thesis to denote an $A$-module labelled by $\lambda \in I$, 
for some index set $I$, and
we write $\pi_{\lambda}$ to denote the representation of $A$ afforded by $V_{\lambda}$.

We say that an element $\lambda \in H^{*}$ is {\emph{integral}} if
$$l_{i} = \frac{2(\lambda,\alpha_{i})}{(\alpha_{i},\alpha_{i})} \in \mathbb{Z}, 
\ \ \forall i < n, \hspace{10mm} l_{n} = \frac{(\lambda,\alpha_{n})}{(\alpha_{n},\alpha_{n})} \in 
\mathbb{Z},$$
where $(\cdot, \cdot): H^{*} \times H^{*} \rightarrow \mathbb{C}$ is the
bilinear form in (\ref{eq:starstar(1)}). 
Let ${\cal{P}}$ denote the set of all integral elements of $H^{*}$.  
We say that an element $\lambda \in {\cal{P}}$ is {\emph{integral dominant}} if 
$l_{i} \in \mathbb{Z}_{+}$ for all $i$.  
We denote the set of all integral dominant elements of $H^{*}$ by ${\cal{P}}^{+}$.

We call a $U_{q}(\mathfrak{g})$-module $V$ a {\emph{highest weight module}} 
if there is a non-zero vector $v \in V$ satisfying  
\begin{itemize}
\item[(i)]   $e_{i}v = 0$,  $\hspace{5mm} 1 \leq i \leq n$,
\item[(ii)]  $K_{i}v = \omega_{i} v$, $\hspace{5mm} \omega_{i} \in \mathbb{C}$,  $\hspace{5mm} 1 \leq i \leq n$,
\item[(iii)] $V \subseteq U_{q}(\mathfrak{b}_{-})v$.
\end{itemize}
The vector $v$ is called a {\emph{highest weight vector}} of $V$.
Similarly, we call a $U_{q}(\mathfrak{g})$-module $V$ a 
{\emph{lowest weight module}} if there is a non-zero vector $w \in V$ satisfying  
\begin{itemize}
\item[(i)]   $f_{i}w = 0$,  $\hspace{5mm} 1 \leq i \leq n$,
\item[(ii)]  $K_{i}w = \omega'_{i} w$, $\hspace{5mm} \omega'_{i} \in \mathbb{C}$,  
$\hspace{5mm} 1 \leq i \leq n$,
\item[(iii)] $V \subseteq U_{q}(\mathfrak{b}_{+})w$.
\end{itemize}
In this case $w$ is called a {\emph{lowest weight vector}}.

Let $V$ be a $U_{q}(\mathfrak{g})$-module on which all the $K_{i}$ act semisimply.
For each sequence 
$\beta = (\beta_{1}, \beta_{2}, \ldots, \beta_{n})$ where $\beta_{j} \in \mathbb{C}$ for each $j$, 
define
$$V_{\beta} = \left\{x \in V | \ K_{i}x = \beta_{i}x, \ 1 \leq i \leq n \right\}.$$
If $V_{\beta} \neq 0$, we say that $\beta$ is a weight of $V$, 
and that $V_{\beta}$ is a weight space of $V$.
The nonzero elements of $V_{\beta}$ are called weight vectors.  

Zou investigated the representation theory of the quantum superalgebra 
$U_{q}(\mathfrak{g})$ over the quotient field 
$\mathbb{C}(v)$ for an indeterminate $v$ \cite{zou}, which is related to $q$ 
via $q=v^{2}$.  
Zou's results can be adapted to our setting where we take $q$ to be a transcendental number.  
Let $\sqrt{q}$ be any square root of the complex number $q$.  
Call a $U_{q}(\mathfrak{g})$-module $V$ {\emph{integrable}} if $V$ 
is a direct sum of its weight spaces and
if $e_{i}$ and  $f_{i}$ act as locally nilpotent endomorphisms of $V$ for each 
$i=1, \ldots, n$.

Let $\overline{V}(\omega)$ be a highest weight $U_{q}(\mathfrak{g})$-module 
with highest weight vector $v$ where 
$K_{i} v = \omega_{i} v$, $\omega_{i} \in \mathbb{C}$, for each $i=1, \ldots, n$. 
The $U_{q}(\mathfrak{g})$-module $\overline{V}(\omega)$ has a unique maximal proper
$U_{q}(\mathfrak{g})$-submodule $\overline{M}(\omega)$, and the quotient
$$V(\omega) = \overline{V}(\omega) / \overline{M}(\omega)$$
is an irreducible $U_{q}(\mathfrak{g})$-module with highest weight
$\omega = (\omega_{1}, \omega_{2}, \ldots, \omega_{n})$.
Theorem 3.1 of \cite{zou} can be stated in our setting as follows.
\begin{theorem}
The irreducible highest weight $U_{q}(\mathfrak{g})$-module $V(\omega)$, 
with $\omega = (\omega_{1}, \omega_{2}, \ldots, \omega_{n})$, is integrable if and only if
$$\omega_{i} = \zeta_{i} q^{m_{i}}, \hspace{5mm} 1 \leq i \leq n-1,$$
where $m_{i} \in \mathbb{Z}_{+}$, $\zeta_{i}^{2}=1$, and
$$ \omega_{n} = \left\{ \begin{array}{ll}
\pm q^{m}, & \mbox{ if } m \in \mathbb{Z}_{+}, \\
\pm\sqrt{-1} \ q^{m}, & \mbox{ if } m \in \mathbb{Z}_{+} + \frac{1}{2}.
\end{array} \right.$$
\end{theorem}
\noindent
Note that every finite dimensional integrable $U_{q}(\mathfrak{g})$-module is semisimple
\cite[Sec. 5]{zou}.

If $\omega = (q^{m_{1}}, q^{m_{2}}, \ldots, q^{m_{n}})$ with $m_{i} \in \mathbb{Z}_{+}$
for each $i$, there exists an irreducible $U(\mathfrak{g})$-module 
$V(\omega)_{\mathfrak{g}}$
 with highest weight $\lambda \in {\cal{P}}^{+}$ satisfying $(\lambda, \alpha_{i}) = m_{i}$ 
 for each $i$, and
 $V(\omega)$ and $V(\omega)_{\mathfrak{g}}$ have the same weight space decomposition \cite{zsuper}.
 In this thesis, we are mostly interested in these irreducible $U_{q}(\mathfrak{g})$-modules, and
 for $\lambda \in {\cal{P}}^{+}$ 
 we let $V_{\lambda}$ denote $V(\omega)$ and call $\lambda$ 
 the highest weight of $V_{\lambda}$.

The following lemma is from \cite[Thms. 2.2, 3.3--3.5]{zsuper}.
\begin{lemma}
\label{lem:trianglesforall}  
Let $V_{\mu}$ and $V_{\nu}$ be finite dimensional irreducible 
$U_{q}(\mathfrak{g})$-modules with highest weights $\mu$ and $\nu$ respectively, where 
$\mu, \nu \in {\cal{P}}^{+}$.
Then $V_{\mu} \otimes V_{\nu}$ can be decomposed into a direct sum of irreducible
$U_{q}(\mathfrak{g})$-modules $V_{\lambda}$ with highest weights $\lambda \in {\cal{P}}^{+}$.
\end{lemma}

Throughout this thesis we always take the grading of the highest weight vector of the 
finite dimensional irreducible
$U_{q}(\mathfrak{g})$-module $V_{\lambda}$ with integral dominant highest weight
$\lambda = \sum_{i=1}^{n} \lambda_{i} \epsilon_{i} \in {\cal{P}}^{+}$ to be
$\left\{ \begin{array}{ll}
\mbox{even}, & \mbox{if } \left( \sum_{i=1}^{n} \lambda_{i} \right) \bmod{2}=0, \\
\mbox{odd},  & \mbox{if } \left( \sum_{i=1}^{n} \lambda_{i} \right) \bmod{2}=1. \\
\end{array} \right.$

In the next lemma adapted from \cite{lz} we consider the 
{\emph{fundamental}} $U_{q}(\mathfrak{g})$-module $V$ which will play
a very important role in this thesis.  This lemma is proved by elementary calculations.
\begin{lemma}
\label{lem:fundamentaldimensional}
There exists a $(2n+1)$-dimensional irreducible $U_{q}(osp(1|2n))$-module $V=V_{\epsilon_{1}}$
with highest weight $\epsilon_{1}$.
This module admits a basis 
$\{ v_{i} | \ -n \leq i \leq n\}$ with $v_{1}$ being the highest weight vector.
The actions of the generators of $U_{q}(osp(1|2n))$ on the basis elements are 
$$f_{i} v_{i} = v_{i+1},  \hspace{5mm}
f_{n} v_{n} = v_{0},  \hspace{5mm}
f_{n} v_{0} = v_{-n},  \hspace{5mm}
f_{i} v_{-i-1} = v_{-i},  
$$
$$
e_{i} v_{i+1} = v_{i}, \hspace{5mm}
e_{n} v_{0}   = v_{n}, \hspace{5mm}
e_{n} v_{-n}  = -v_{0}, \hspace{5mm}
e_{i} v_{-i}  = v_{-i-1}, $$  
$$
K_{j}^{\pm 1} v_{k} = q^{\pm (\alpha_{j},\epsilon_{k})} v_{k},$$
where $1 \leq i< n$,  $1 \leq j \leq n$, $-n \leq k \leq n$, 
$\epsilon_{0}=0$, and $\epsilon_{-i} = - \epsilon_{i}$.
All remaining actions are zero.
\end{lemma}

Note that the highest weight vector $v_{1}$ of the fundamental 
irreducible $U_{q}(osp(1|2n))$-module $V$ always has an odd grading in this thesis.

\begin{proposition}
\label{prop:mooV}
There exists a $U_{q}(\mathfrak{g})$-invariant, non-degenerate bilinear form
$\langle \langle \ , \ \rangle \rangle : V \times V \rightarrow \mathbb{C}$.
Thus the dual $U_{q}(\mathfrak{g})$-module of $V$ is isomorphic to $V$.
\end{proposition}
\begin{proof}
Let $\{ v_{i} | \ -n \leq i \leq n \}$ be the basis of $V$ given in Lemma \ref{lem:fundamentaldimensional}.
Now define a non-degenerate $\mathbb{C}$-bilinear form
$$\langle \langle \ , \ \rangle \rangle : V \times V \rightarrow \mathbb{C},$$
by
$$
\begin{array}{rcll}
\langle \langle v_{1}, v_{-1} \rangle \rangle & = & 1, & \\
\langle \langle v_{i}, v_{-i} \rangle \rangle & = & 
-q^{-1} \langle \langle v_{i-1}, v_{-(i-1)} \rangle \rangle, &  2 \leq i \leq n, \\
\langle \langle v_{0}, v_{0} \rangle \rangle & = & q^{-1} \langle \langle v_{n}, v_{-n} \rangle \rangle, & \\
\langle \langle v_{-n}, v_{n} \rangle \rangle & = & - \langle \langle v_{0}, v_{0} \rangle \rangle, & \\
\langle \langle v_{-j}, v_{j} \rangle \rangle & = & 
-q^{-1} \langle \langle v_{-(j+1)}, v_{j+1} \rangle \rangle, &  1 \leq j \leq n-1, \\
\langle \langle v_{k}, v_{l} \rangle \rangle & = & 0, & \mbox{if } k + l \neq 0.
\end{array}
$$
A direct calculation shows that
$$\langle \langle x v_{i}, v_{j} \rangle \rangle = 
(-1)^{[x][v_{i}]}  \langle \langle v_{i}, S(x) v_{j} \rangle \rangle,  
\hspace{5mm} \forall x \in U_{q}(\mathfrak{g}), \hspace{5mm} v_{i}, v_{j} \in V,$$
thus proving the $U_{q}(\mathfrak{g})$-invariance of the bilinear form.
This form identifies $V$ with its dual module.
\end{proof}

Let us discuss in more detail the dual module of $V$.
Recall the definition of the dual $U_{q}(\mathfrak{g})$-module $V^{*}$ to $V$.
Let $\{v_{i}^{*} | \ -n \leq i \leq n\}$ be a basis of $V^{*}$ such that
$\langle v_{i}^{*}, v_{j} \rangle = \delta_{ij}$ and $[v_{i}^{*}]=[v_{i}]$ where
$\langle \ , \  \rangle: V^{*} \times V \rightarrow \mathbb{C}$ is the dual space pairing.
Now define a homogeneous bijection $T \in \mathrm{Hom}_{\mathbb{C}}(V, V^{*})$ of degree $0$ by
\begin{equation}
\label{eq:johnfaulkner(20)}
T:  v_{i}  \mapsto (-1)^{i-1} q^{-(i-1)}v_{-i}^{*}, \hspace{5mm}
       v_{0}  \mapsto (-1)^{n-1} q^{-n} v_{0}^{*}, \hspace{5mm}
       v_{-i} \mapsto (-1)^{i} q^{-(2n-i)} v_{i}^{*}, \hspace{5mm} 1 \leq i \leq n.
\end{equation}
A direct calculation shows that
this map is an element of $\mathrm{Hom}_{U_{q}(\mathfrak{g})}(V, V^{*})$ and that it satisfies
$$ \langle T(v_{i}), v_{j} \rangle = \langle \langle v_{i}, v_{j} \rangle \rangle,
 \ \hspace{10mm} \mbox{for all } v_{i}, v_{j} \in V.$$

\end{section}

\begin{section}{$R$-matrices for representations of $U_{q}(osp(1|2n))$}
\label{subsec:RRRR}
\markright{\text{$R$-matrices for representations of $U_{q}(osp(1|2n))$}}

Drinfel'd's quantum superalgebra $U_{h}(\mathfrak{g})$ 
has a universal $R$-matrix \cite{kt}.
We will show that there does not exist an element of $U_{q}(\mathfrak{g})$ over $\mathbb{C}$
that corresponds to the universal $R$-matrix of $U_{h}(\mathfrak{g})$ in an obvious way.
However, there is a completion 
$\overline{U}_{q}^{+}(\mathfrak{g})$ of $U_{q}(\mathfrak{g})$ such
that one of the multiplicative factors of the universal $R$-matrix of $U_{h}(\mathfrak{g})$
maps to an element $\widetilde{R}$ of $\overline{U}_{q}^{+}(\mathfrak{g})$.
Although $\widetilde{R}$ is not an element of $U_{q}(\mathfrak{g})$, 
only a finite number of terms of $\widetilde{R}$ act as non-zero endomorphisms on each 
tensor product of finite dimensional irreducible $U_{q}(\mathfrak{g})$-modules, and thus
the action of $\widetilde{R}$ on these tensor products is well-defined.

\begin{subsection}{The universal $R$-matrix of $U_{h}(osp(1|2n))$}

The Drinfel'd quantum algebras admit universal $R$-matrices \cite{ls,kirresh},
and the explicit expressions of these universal $R$-matrices employ infinite sums of root vectors
in the quantum algebra
corresponding to the positive roots of the associated Lie algebra \cite{cp}.  

The simple root vectors in the quantum algebra are just the generators $e_{i}$ and $f_{i}$,
and the non-simple root vectors are obtained by applying  Lusztig's automorphisms
to the generators \cite{lusztig1,cp}.  
The non-simple root vectors are  not uniquely defined in general \cite{dk,cp}.  
Different choices for the decomposition of the longest element of the Weyl group of the associated
Lie algebra into a product of reflections generated by the simple roots may
lead to non-simple root vectors that may not even be proportional to each other, and 
{\emph{a priori}} there is no   
canonical choice for the decomposition of the longest element of the Weyl group into a product of such
reflections \cite{cp}.

Khoroshkin and Tolstoy wrote down the universal $R$-matrix of $U_{h}(osp(1|2n))$  \cite{kt}
using a different method.
They employed infinite sums of root vectors in $U_{h}(\mathfrak{g})$,
but these root vectors were defined in a different way to how the root vectors in quantum algebras
were defined. 
Khoroshkin and Tolstoy's procedure is general for quantum superalgebras and we write it down for
$U_{h}(\mathfrak{g})$ here.

In $U_{h}(\mathfrak{g})$, root vectors are only defined for the elements of the
{\emph{reduced root system}} $\phi$ of $\mathfrak{g}$.
The reduced root system $\phi$ of $\mathfrak{g}$ is the set of all 
positive roots of $\mathfrak{g}$ except those roots $\alpha$ for which $\alpha/2$ is also a positive root,
and the reduced root system of $\mathfrak{g}=osp(1|2n)$ 
is just $\phi=\overline{\Phi}^{+}_{0} \cup \Phi^{+}_{1}$.

A total ordering of $\phi$ called a {\emph{normal ordering}} 
is then introduced, and the root vectors of $U_{h}(\mathfrak{g})$
are recursively defined using the normal ordering of $\phi$ 
and a map involving the $q$-bracket that we introduce below. 
We denote a normal ordering of $\phi$ by ${\cal{N}}(\phi)$.
A difference between the root vectors in quantum algebras and 
the root vectors in $U_{q}(\mathfrak{g})$ is that
the latter are defined by a map that is not necessarily an algebra automorphism.
The way the universal $R$-matrix of $U_{h}(\mathfrak{g})$ is formally written down depends on
${\cal{N}}(\phi)$.

A normal ordering of a reduced root system of $\mathfrak{g}$ 
is defined as follows \cite[Def. 3.1]{kt}.
\begin{definition}
A normal ordering ${\cal{N}}(\phi)$ of $\phi =  \overline{\Phi}^{+}_{0} \cup \Phi^{+}_{1}$
is a total order $\prec$ of the set $\phi$ such that if $\alpha \prec \beta$ and 
$\alpha + \beta \in \phi$, then $\alpha \prec \alpha + \beta \prec \beta$.
\end{definition}
 
In general, there is more than one normal ordering of $\phi$ \cite{kt}.  
For example, the reduced root system of
$osp(1|4)$ is $\phi = \{ \epsilon_{1}, \epsilon_{2}, \epsilon_{1} \pm \epsilon_{2}\}$ 
and there are two different normal orderings of $\phi$:
$$\alpha_{1} \prec \alpha_{1} + \alpha_{2} \prec \alpha_{1} + 2\alpha_{2} \prec \alpha_{2},$$
$$\alpha_{2} \prec \alpha_{1} + 2\alpha_{2} \prec \alpha_{1} + \alpha_{2} \prec \alpha_{1},$$
where we write the elements of $\phi$ as sums of the simple roots.

We now write down the universal $R$-matrix of $U_{h}(osp(1|2n))$ \cite{kt}, 
which we adapt slightly to take account of the different co-multiplication used in this thesis.
Writing $q=e^{h} \in \mathbb{C}[[h]]$,  let us firstly define
$\displaystyle{\exp_{q}{(x)} = \sum_{k=0}^{\infty} \frac{x^{k}}{\left[k\right]^{q}!} }$. 

We now construct the root vectors in $U_{h}(osp(1|2n))$.  
The easiest root vectors are the simple root vectors:
we fix $E_{\alpha_{i}}=E_{i}$, $F_{\alpha_{i}}=F_{i}$ and $H_{\alpha_{i}}=H_{i}$ 
for each simple root $\alpha_{i}$.
Now we recursively construct the non-simple root vectors.
Let $\alpha, \beta, \gamma \in \phi$ be roots such that $\gamma = \alpha + \beta$ and 
$\alpha \prec \beta$, and in addition let no other roots 
$\alpha', \beta' \in \phi$ exist which  
satisfy (i) $\alpha' + \beta' = \gamma$, 
(ii) $\alpha \prec \alpha' \prec \beta$ and (iii) $\alpha \prec \beta' \prec \beta$.  
Then, if all of the root vectors 
$E_{\alpha}, E_{\beta}, F_{\alpha}, F_{\beta} \in U_{h}(\mathfrak{g})$ 
have already been defined, we define
$$E_{\gamma} = \left[E_{\alpha},E_{\beta} \right]_{q}, \hspace{10mm} 
F_{\gamma} = \left[F_{\beta},F_{\alpha}\right]_{q^{-1}},$$
 where the $q$-bracket
 $[ \cdot, \cdot ]_{q}$
 is defined by
 $$[X_{\alpha}, X_{\beta}]_{q} = 
    X_{\alpha} X_{\beta} - (-1)^{[X_{\alpha}][X_{\beta}]} q^{(\alpha,\beta)} X_{\beta}X_{\alpha},$$
where $X$ stands for $E$ or $F$.

For each $\gamma \in \phi$, set
\begin{equation}
\label{eq:johnfaulkner(10)}
R_{\gamma} = 
\exp_{q_{\gamma}}{\big((-1)^{[\gamma]} (a_{\gamma})^{-1}(q-q^{-1}) E_{\gamma} \otimes F_{\gamma}\big)}
\in U_{h}(\mathfrak{g}) \otimes U_{h}(\mathfrak{g}),
\end{equation}
where $q_{\gamma} = (-1)^{[\gamma]} q^{-(\gamma, \gamma)}$ and
$a_{\gamma} \in \mathbb{C}[[h]]$ is defined by 
$$E_{\gamma} F_{\gamma}-(-1)^{[E_{\gamma}]}F_{\gamma}E_{\gamma}=
\frac{a_{\gamma}\left(q^{H_{\gamma}}-q^{-H_{\gamma}}\right)}{q-q^{-1}}.$$
It is important to observe that $a_{\gamma}$ is a rational function of $q$.
Now we can write down the universal $R$-matrix of $U_{h}(\mathfrak{g})$  \cite[Thm. 8.1]{kt}.
\begin{theorem}
Define $\mathbf{H}_{i} = \sum_{j=i}^{n} H_{j} \in U_{h}(\mathfrak{g})$ for each $i=1, \ldots, n$, then
\begin{equation}
\label{eq:reneeboyle(1)}
R=\exp{\left(h \sum_{i=1}^{n} \mathbf{H}_{i} \otimes \mathbf{H}_{i} \right)}
\prod_{\gamma \in \phi} R_{\gamma},
\end{equation}
is the universal $R$-matrix of $U_{h}(\mathfrak{g})$, where the product 
is ordered with respect to the same normal ordering ${\cal{N}}(\phi)$ 
 used to define the root vectors in $U_{h}(\mathfrak{g})$ so that
 $\prod_{\gamma \in \phi} R_{\gamma} = R_{\gamma_{1}} R_{\gamma_{2}} \cdots R_{\gamma_{k}}$ where
 $\gamma_{1} \prec \gamma_{2} \prec \cdots \prec \gamma_{k}$ in ${\cal{N}}(\phi)$.
 
\end{theorem}

\end{subsection}

\begin{subsection}{$R$-matrices for representations of $U_{q}(osp(1|2n))$}
\label{subsect:Rmatricesfrorepresofquantumosp}

It is unknown whether Jimbo's quantum algebras over $\mathbb{C}$ have universal $R$-matrices.
However, there exist 
{\emph{$R$-matrices for representations}} of these quantum algebras.
Let $\pi_{\lambda}$ and $\pi_{\mu}$ be finite dimensional irreducible representations 
of the quantum algebra $A$, then there is an invertible element 
${\cal{R}}_{\lambda,\mu} \in End_{\mathbb{C}}(V_{\lambda} \otimes V_{\mu})$ satisfying  
 \begin{equation}
\label{chap2:FyHendy(1)}
{\cal{R}}_{\lambda,\mu} \cdot (\pi_{\lambda} \otimes \pi_{\mu}) \big(\Delta(x)\big) = 
        (\pi_{\lambda} \otimes \pi_{\mu})  \big(\Delta'(x)\big)  \cdot {\cal{R}}_{\lambda,\mu}
	\hspace{10mm} \forall x \in A,
\end{equation}
which we note is just Eq. (\ref{eq:asdf;lkj}) in $\pi_{\lambda} \otimes \pi_{\mu}$.
We follow a similar scheme for $U_{q}(osp(1|2n))$.

For each pair of finite dimensional irreducible $U_{q}(\mathfrak{g})$-modules, we will construct
an element ${\cal{R}}_{\lambda,\mu} \in End_{\mathbb{C}}(V_{\lambda} \otimes V_{\mu})$ satisfying 
(\ref{chap2:FyHendy(1)}) for all $x \in U_{q}(\mathfrak{g})$.
We do this following the method laid out in \cite{cp,ks} for quantum algebras.

We firstly define a completion $\overline{U}^{+}_{q}(\mathfrak{g})$ of $U_{q}(\mathfrak{g})$
following \cite[Subsec. 6.3.3]{ks}.
Let $U_{q}(\mathfrak{n}_{+})$ (resp. $U_{q}(\mathfrak{n}_{-})$) be the subalgebra  of
$U_{q}(\mathfrak{g})$ generated by the elements $\{ e_{i} | \ 1 \leq i \leq n \}$
(resp. $\{ f_{i} | \ 1 \leq i \leq n \}$).  
We say that a non-zero element $x \in U_{q}(\mathfrak{g})$ has degree 
$\lambda = \sum_{i=1}^{n} m_{i}\alpha_{i}, \ m_{i} \in \mathbb{Z}$, 
if $K_{i}xK_{i}^{-1} = q^{(\lambda,\alpha_{i})}x$ for all $i=1, 2, \ldots, n$. 
We define $\overline{U}^{\pm}_{q}(\mathfrak{g})$ by
$$\overline{U}^{\pm}_{q}(\mathfrak{g}) = 
\prod_{\beta \in Q_{+}} U_{q}(\mathfrak{b}_{\pm})U_{q}^{\mp \beta}(\mathfrak{n}_{\mp}),$$
where $U_{q}^{\mp \beta}(\mathfrak{n}_{\mp})$ is defined by
$$U_{q}^{\pm \beta}(\mathfrak{n}_{\pm}) = 
\{ x \in U_{q}(\mathfrak{n}_{\pm}) |
 \ K_{i} x K_{i}^{-1} = q^{\pm (\alpha_{i},\beta)}x \}, \hspace{10mm} i=1, 2, \ldots, n,$$
 and $Q_{+}$ is defined by
 $Q_{+} = \{ \sum_{i=1}^{n} n_{i} \alpha_{i} | \ n_{i} \in \mathbb{Z}_{+} \}$.

The elements of $\overline{U}^{\pm}_{q}(\mathfrak{g})$ are sequences 
$x = (x_{\beta})_{\beta \in Q_{+}}$ where 
$x_{\beta} \in U_{q}(\mathfrak{b}_{\pm})U_{q}^{\mp \beta}(\mathfrak{n}_{\mp})$.
Let us write this sequence formally as an infinite sum $x = \sum_{\beta} x_{\beta}$.  
The results of \cite[Subsec. 6.1.5]{ks}) imply that $U_{q}(\mathfrak{g})$ can be expressed
as a direct sum
$$U_{q}(\mathfrak{g}) = \bigoplus_{\beta \in Q_{+}} U_{q}(\mathfrak{b}_{+}) U_{q}^{-\beta}(\mathfrak{n}_{-}),$$
and thus
$U_{q}(\mathfrak{g})$ can be considered as the subspace of
$\overline{U}^{+}_{q}(\mathfrak{g})$ formed by the sums $x=\sum_{\beta} x_{\beta}$ 
for which all but finitely many terms vanish.

The multiplication in $U_{q}(\mathfrak{g})$ extends to a multiplication in 
$\overline{U}^{\pm}_{q}(\mathfrak{g})$.
Let $\beta, \gamma \in Q_{+}$, let $x_{\beta} \in U_{q}^{-\beta}(\mathfrak{n}_{-})$, and 
let $y_{\gamma} \in U_{q}(\mathfrak{b}_{+})$ have degree $\gamma$.  
From the commutation relations of $U_{q}(\mathfrak{g})$, we have
$$x_{\beta} y_{\gamma} \in \bigoplus_{\delta}
U_{q}(\mathfrak{b}_{+})U_{q}^{-\delta}(\mathfrak{n}_{-}),$$
where the direct sum ranges over all $\delta \in Q_{+}$ satisfying
$|\beta| - |\gamma| \leq |\delta| \leq |\beta|$, where $|\beta| = \sum_{i=1}^{n} m_{i}$ 
for $\beta = \sum_{i=1}^{n} m_{i} \alpha_{i}$.  Thus 
$\overline{U}^{+}_{q}(\mathfrak{g})$ is an algebra, and similarly
$\overline{U}^{-}_{q}(\mathfrak{g})$ and
$\overline{U}^{\pm}_{q}(\mathfrak{g}) \overline{\otimes} 
\cdots \overline{\otimes} \overline{U}^{\pm}_{q}(\mathfrak{g})$ 
are algebras.  
The algebras $\overline{U}^{\pm}_{q}(\mathfrak{g})$ and 
$\overline{U}^{\pm}_{q}(\mathfrak{g}) \overline{\otimes} 
\cdots \overline{\otimes} \overline{U}^{\pm}_{q}(\mathfrak{g})$ (with $m$ factors) contain
$U_{q}(\mathfrak{g})$ and $U_{q}(\mathfrak{g}) \otimes \cdots \otimes U_{q}(\mathfrak{g})$ 
(with $m$ factors), respectively, as subalgebras.

We now construct an element in $\overline{U}^{+}_{q}(\mathfrak{g})$ corresponding to the element
$\prod_{\gamma \in \phi} R_{\gamma}$ in (\ref{eq:reneeboyle(1)}).
Given a normal ordering ${\cal{N}}(\phi)$ for a reduced root system $\phi$, 
we construct root vectors  $E_{\gamma}, F_{\gamma} \in U_{q}(\mathfrak{g})$ 
following the same procedure as in $U_{h}(\mathfrak{g})$ by setting
$E_{\alpha_{i}} = e_{i}$ and $F_{\alpha_{i}} = f_{i}$ and thinking of $q$ a complex number.
Then $R_{\gamma}$ is well-defined as an element of
$\overline{U}^{+}_{q}(\mathfrak{g}) \overline{\otimes} \overline{U}^{+}_{q}(\mathfrak{g})$,
and to simplify $R_{\gamma}$ we normalise the root vectors:
$$e_{\gamma} = E_{\gamma}, \hspace{10mm} f_{\gamma} = F_{\gamma}/a_{\gamma}.$$
This is well-defined as $a_{\gamma} \neq 0$ \cite[Eqs. (8.3)--(8.4)]{kt}.
Simplifying the expression for $R_{\gamma}$, we have
$$R_{\gamma} = \left\{ \begin{array}{ll}
\displaystyle{\sum_{k=0}^{\infty} 
\frac{\left(q-q^{-1}\right)^{k}\left(e_{\gamma} \otimes f_{\gamma}\right)^{k}}
{\left[k\right]^{q^{-2}}!},} & \mbox{if } [e_{\gamma}]=0,  \\
\displaystyle{\sum_{k=0}^{\infty} 
\frac{\left(q^{-1}-q\right)^{k}\left(e_{\gamma} \otimes f_{\gamma}\right)^{k}}
{\left[k\right]^{-q^{-1}}!}, } & \mbox{if } [e_{\gamma}]=1. 
\end{array} \right.$$

Define $\widetilde{R} \in 
\overline{U}^{+}_{q}(\mathfrak{g}) \overline{\otimes} \overline{U}^{+}_{q}(\mathfrak{g})$ by
$\widetilde{R} = \prod_{\gamma \in \phi} R_{\gamma}$ where the product is ordered 
using the same normal order ${\cal{N}}(\phi)$ that we used to define the root vectors in
$U_{q}(\mathfrak{g})$ and such that
$\prod_{\gamma \in \phi} R_{\gamma} = R_{\gamma_{1}} R_{\gamma_{2}} \cdots R_{\gamma_{k}}$
where ${\cal{N}}(\phi) = \gamma_{1} \prec \gamma_{2} \prec \cdots \prec \gamma_{k}$.
Clearly $\widetilde{R}$ is invertible as 
$\prod_{\gamma \in \phi} R_{\gamma} \in U_{h}(\mathfrak{g}) \otimes U_{h}(\mathfrak{g})$
is invertible and $q$ is not a root of unity.

\begin{lemma}
\label{lemlem:twenty}
Define an automorphism $\Psi$ of $U_{q}(\mathfrak{g}) \otimes U_{q}(\mathfrak{g})$ by
$$
\begin{array}{lll}
\Psi(K_{i}^{\pm 1} \otimes 1) = K_{i}^{\pm 1} \otimes 1, & & 
\Psi(1 \otimes K_{i}^{\pm 1}) = 1 \otimes K_{i}^{\pm 1}, \\
\Psi(e_{i} \otimes 1) = e_{i} \otimes K_{i}^{-1}, & &
\Psi(1 \otimes e_{i}) = K_{i}^{-1} \otimes e_{i}, \\
\Psi(f_{i} \otimes 1) = f_{i} \otimes K_{i}, & &
\Psi(1 \otimes f_{i}) = K_{i} \otimes f_{i}.
\end{array}
$$
The automorphism $\Psi$ satisfies the following relations:
\begin{itemize}
\item[(i)] $\widetilde{R} \Delta(x) = \Psi \big(\Delta'(x) \big) \cdot \widetilde{R}, 
\hspace{10mm} \mbox{for all } x \in U_{q}(\mathfrak{g})$,
\item[(ii)] $(\Delta \otimes \mathrm{id}) \widetilde{R}=\Psi_{23} (\widetilde{R}_{13}) \cdot \widetilde{R}_{23}$,
\item[(iii)] $(\mathrm{id} \otimes \Delta) \widetilde{R}=\Psi_{12} (\widetilde{R}_{13}) \cdot \widetilde{R}_{12}$,
\end{itemize} 
where $\Psi_{12} = \Psi \otimes \mathrm{id}$ and $\Psi_{23} = \mathrm{id} \otimes \Psi$.
\end{lemma}
\begin{proof}
We prove (i) for each generator of $U_{q}(\mathfrak{g})$.  
We firstly wish to prove the following equations:
\begin{equation}
\label{eq:oneoneoneoneone}
\widetilde{R} (e_{i} \otimes K_{i} + 1 \otimes e_{i}) = 
(e_{i} \otimes K_{i}^{-1} + 1 \otimes e_{i})\widetilde{R},
\end{equation}
\begin{equation}
\label{eq:twotwo}
\widetilde{R} (f_{i} \otimes 1 + K_{i}^{-1} \otimes f_{i}) =
(f_{i} \otimes 1 + K_{i} \otimes f_{i})\widetilde{R},
\end{equation}
\begin{equation}
\label{eq:threethreethree}
\widetilde{R} (K_{i}^{\pm 1} \otimes K_{i}^{\pm 1}) = 
(K_{i}^{\pm 1} \otimes K_{i}^{\pm 1}) \widetilde{R}.
\end{equation}
Now Eq. (\ref{eq:threethreethree}) is true by inspection and
Eqs. (\ref{eq:oneoneoneoneone})--(\ref{eq:twotwo}) follow from the corresponding results in 
$U_{h}(\mathfrak{g})$ \cite[Prop. 6.2]{kt}.
The proof of (i) then follows from the definition of $\Psi$ and 
the proofs of (ii) and (iii) follow similarly from \cite{kt}.
\end{proof}

We now examine the usual approach used to create $R$-matrices for representations
of a  quantum algebra $A$.
For each tensor product $W_{1} \otimes W_{2}$
of finite dimensional integrable $A$-modules, an invertible element 
${\cal{E}}_{W_{1}, W_{2}} \in End_{\mathbb{C}} (W_{1} \otimes W_{2})$ is constructed
implementing the automorphism $\Psi$, in the sense that ${\cal{E}}_{W_{1}, W_{2}}$ satisfies
$${\cal{E}}_{W_{1}, W_{2}}^{-1} \cdot (\pi_{W_{1}} \otimes \pi_{W_{2}})(x) 
\cdot {\cal{E}}_{W_{1}, W_{2}} = 
(\pi_{W_{1}} \otimes \pi_{W_{2}})\Psi(x), 
\hspace{5mm} \mbox{for all } x \in A \otimes  A.$$
This ${\cal{E}}_{W_{1}, W_{2}}$ is defined by fixing its action to be
$${\cal{E}}_{W_{1}, W_{2}} (w_{\lambda} \otimes w_{\mu}) = 
   q^{(\lambda,\mu)} (w_{\lambda} \otimes w_{\mu}),$$
on all weight vectors $ w_{\lambda} \in W_{1}$, $w_{\mu} \in W_{2}$ 
 with weights $\lambda$ and $\mu$, respectively \cite[Prop. 10.1.19]{cp}.

We could have used the same method to construct $R$-matrices for representations of
$U_{q}(osp(1|2n))$ but we have found a more useful approach.
Note that in the above we need to know
the weight space decompositions of both  $W_{1}$ and $W_{2}$
before defining ${\cal{E}}_{W_{1}, W_{2}}$, but if we have this information 
we can use it in a more interesting way: instead of defining an element of
$End_{\mathbb{C}}(W_{1} \otimes W_{2})$
we can define
an element $E_{W_{1},W_{2}} \in U_{q}(\mathfrak{g})$ with the property that
$(\pi_{W_{1}} \otimes \pi_{W_{2}}) E_{W_{1},W_{2}}={\cal{E}}_{W_{1}, W_{2}}$.
We define this $E_{W_{1},W_{2}}$ for each tensor product of 
finite dimensional irreducible $U_{q}(\mathfrak{g})$-modules
following a related idea in \cite{z1}.

For each $i=1, \ldots, n$, set $$J_{i} = K_{i}K_{i+1} \cdots K_{n}.$$
The action of $J_{i}$ on a weight vector 
$w_{\xi}$ with weight $\xi = \sum_{j=1}^{n} \xi_{j} \epsilon_{j} \in H^{*}$ 
 of a $U_{q}(\mathfrak{g})$-module is
$$J_{i} w_{\xi} = q^{\xi_{i}} w_{\xi}.$$

Consider the weight space decomposition of a 
finite dimensional irreducible $U_{q}(\mathfrak{g})$-module $V_{\mu}$ 
with integral dominant highest weight $\mu$.  
The weight of the weight vector $w_{\xi} \in V_{\mu}$ is 
$\xi = \sum_{i=1}^{n} \xi_{i} \epsilon_{i} \in \bigoplus_{i=1}^{n} \mathbb{Z} \epsilon_{i}$.  
Now define
\begin{equation}
\label{eq:chatswoodtrain}
E_{\mu} = \prod_{a=1}^{n} \sum^{s}_{b=p} (J_{a})^{b} \otimes P_{a}[b], \hspace{10mm} 
P_{a}[b] = \prod_{\stackrel{c = p}{c \neq b}}^{s} \frac{J_{a}-q^{c}}{q^{b}-q^{c}}, \hspace{5mm} c \in \mathbb{Z},
\end{equation}
where $p$ and $s$ are any integers satisfying $p \leq s$ and the following condition:
\begin{itemize}
\item 
$J_{i} w_{\xi} = q^{\xi_{i}} w_{\xi}$ for some $\xi_{i}$ satisfying $p \leq \xi_{i} \leq s$,
for each weight vector $w_{\xi} \in V_{\mu}$.
\end{itemize}
Once we have any such $p$ and $s$, we can use any $p'$ and $s'$ satisfying
$p' \leq p$ and $s' \geq s$ in (\ref{eq:chatswoodtrain}) instead of $p$ and $s$, respectively.

The element $E_{\mu}$ is well-defined and invertible in 
$U_{q}(\mathfrak{g}) \otimes U_{q}(\mathfrak{g})$, and
for all weight vectors $v_{\lambda'} \in V_{\lambda}$ and $v_{\mu'} \in V_{\mu}$ we have
\begin{equation}
\label{chapter2:opera(1)}
E_{\mu} (v_{\lambda'} \otimes v_{\mu'}) = q^{(\lambda',\mu')} \ (v_{\lambda'} \otimes v_{\mu'}),
\end{equation}
where the weights of $v_{\lambda'}$ and $v_{\mu'}$ are $\lambda'$ and $\mu'$, respectively.
The element $E_{\mu}$ is not a {\emph{universal}} element in that it does not 
satisfy (\ref{chapter2:opera(1)}) for all representations of $U_{q}(\mathfrak{g})$.
It would be very useful if one could construct such a universal element in $U_{q}(\mathfrak{g})$.

From this we obtain $R$-matrices for tensor products of finite dimensional irreducible
$U_{q}(\mathfrak{g})$-modules in the following sense.
Let $V_{\lambda}$ and $V_{\mu}$ be irreducible $U_{q}(\mathfrak{g})$-modules with  
integral dominant highest weights $\lambda$ and $\mu$, respectively.  
Then we have the following important theorem.
\begin{theorem}
\label{th:RRRmatrix}
Define $R_{\lambda,\mu} = E_{\mu} \widetilde{R} \in 
\overline{U}_{q}^{+}(\mathfrak{g}) \overline{\otimes} \overline{U}_{q}^{+}(\mathfrak{g})$
and  ${\cal{R}}_{\lambda,\mu} = (\pi_{\lambda} \otimes \pi_{\mu})R_{\lambda,\mu}$, then
\begin{equation}
\label{eq:sachanadra(1)}
{\cal{R}}_{\lambda,\mu} \cdot (\pi_{\lambda} \otimes \pi_{\mu}) \big(\Delta(x)\big) 
    = (\pi_{\lambda} \otimes \pi_{\mu})\big(\Delta'(x)\big) \cdot {\cal{R}}_{\lambda,\mu}, 
    \hspace{10mm} \forall x \in U_{q}(\mathfrak{g}).
\end{equation}
\end{theorem}
\begin{proof}
This is similar to the proof of the corresponding theorem in quantum algebras \cite[Prop. 10.1.19]{cp},
and we follow that proof.
By direct calculation we can readily show that
$$(\pi_{\lambda} \otimes \pi_{\mu}) \Psi \big(\Delta(x)\big)  = 
(\pi_{\lambda} \otimes \pi_{\mu})  \big( E_{\mu}^{-1} \cdot \Delta(x) \cdot E_{\mu} \big),
\hspace{10mm} \forall x \in U_{q}(\mathfrak{g}),$$
 then by using  
$\widetilde{R}\Delta(x) = \Psi \big(\Delta'(x)\big) \cdot \widetilde{R}$ from 
Lemma \ref{lemlem:twenty} we have
\begin{equation}
\label{chap2:equationannie(alpha)}
(\pi_{\lambda} \otimes \pi_{\mu}) \big(R_{\lambda,\mu} \Delta(x) \big) = 
(\pi_{\lambda} \otimes \pi_{\mu}) \big( \Delta'(x) R_{\lambda,\mu} \big),
\end{equation}
which is precisely Eq. (\ref{eq:sachanadra(1)}).
\end{proof}
\noindent
A similar calculation shows that
$$ (\pi_{\lambda} \otimes \pi_{\mu}) \big( R_{\mu,\lambda}^{T} \Delta'(x) \big)
 = (\pi_{\lambda} \otimes \pi_{\mu}) \big( \Delta(x) R_{\mu,\lambda}^{T} \big), \hspace{10mm}
 \forall x \in U_{q}(\mathfrak{g}).$$

We now determine some useful results involving $R_{\lambda,\mu}$.
\begin{proposition}
The element $R_{\lambda,\mu}$ has the following properties:
\begin{equation}
\label{chapt2:annieeq(1)}
(\epsilon \otimes \mathrm{id}) R_{\lambda,\mu} = (\mathrm{id} \otimes \epsilon) R_{\lambda,\mu} = 1,
\end{equation}
\begin{equation}
\label{chapt2:annieeq(2)}
(\pi_{\lambda} \otimes \pi_{\mu})\big( (S \otimes \mathrm{id}) R_{\lambda,\mu}\big) = 
(\pi_{\lambda} \otimes \pi_{\mu})R_{\lambda,\mu}^{-1}, \hspace{5mm}
  (\pi_{\lambda} \otimes \pi_{\mu})\big((\mathrm{id} \otimes S) R_{\lambda,\mu}^{-1}\big) = R_{\lambda,\mu},
  \end{equation}
  \begin{equation}
  \label{chapt2:annieeq(3)}
  (\pi_{\lambda} \otimes \pi_{\mu})\big((S \otimes S) R_{\lambda,\mu}\big) = 
  (\pi_{\lambda} \otimes \pi_{\mu})R_{\lambda,\mu}.
\end{equation}
\end{proposition}
\begin{proof}
Eq. (\ref{chapt2:annieeq(1)}) is proved by inspection.  
The proofs of (\ref{chapt2:annieeq(2)})--(\ref{chapt2:annieeq(3)}) are straightforward
and almost identical to the proofs of the corresponding equations 
in $\mathbb{Z}_{2}$-graded quasitriangular Hopf algebras. 
\end{proof}

Let $v_{\lambda}$ and $v_{\nu}$ be weight vectors of $U_{q}(\mathfrak{g})$-modules 
with weights
$\lambda$, $\nu \in \bigoplus_{i=1}^{n} \mathbb{Z} \epsilon_{i}$, respectively
and let $v_{\mu'} \in V_{\mu}$ be a weight vector with weight $\mu'$, then
it can be easily shown that
$$\big[(\Delta \otimes 1)E_{\mu}\big] (v_{\lambda} \otimes v_{\nu} \otimes v_{\mu'}) = 
q^{(\mu',\lambda + \nu)}(v_{\lambda} \otimes v_{\nu} \otimes v_{\mu'}),$$
$$\big[(1 \otimes \Delta)E_{\mu}\big] (v_{\lambda} \otimes v_{\nu} \otimes v_{\mu'}) =
q^{(\lambda, \nu + \mu')}(v_{\lambda} \otimes v_{\nu} \otimes v_{\mu'}),$$
where in  $(1 \otimes \Delta)E_{\mu}$ we change the limits $p$ and $s$ if necessary.

We now consider analogues in $U_{q}(\mathfrak{g})$ of the equations
$(\Delta \otimes 1) R = R_{13} R_{23}$ and $(1 \otimes \Delta)R = R_{13} R_{12}$
of a $\mathbb{Z}_{2}$-graded quasitriangular Hopf algebra.  
By definition, we have
$(\Delta \otimes 1) R_{\lambda,\mu} = 
\big[ (\Delta \otimes 1) E_{\mu} \big] \cdot \Psi_{23}(\widetilde{R}_{13}) \cdot \widetilde{R}_{23}$.
Let $V_{\lambda}$ and $V_{\nu}$ be finite dimensional irreducible $U_{q}(\mathfrak{g})$-modules, 
then from the properties of $E_{\mu}$ we have
\begin{eqnarray}
(\pi_{\lambda} \otimes \pi_{\nu} \otimes \pi_{\mu}) \big[ (\Delta \otimes 1) R_{\lambda,\mu} \big] 
& = &  (\pi_{\lambda} \otimes \pi_{\nu} \otimes \pi_{\mu})
        \left[  (\Delta \otimes 1) E_{\mu}  \cdot 
	(E_{\mu}^{23})^{-1} \widetilde{R}_{13} E_{\mu}^{23}\widetilde{R}_{23} \right] \nonumber \\
& = &  (\pi_{\lambda} \otimes \pi_{\nu} \otimes \pi_{\mu})
        \left[ E_{\mu}^{13} \widetilde{R}_{13} E_{\mu}^{23} \widetilde{R}_{23} \right], 
	\label{chapter2:R-matrix(99)}
\end{eqnarray}
where writing $E_{\mu} = \sum_{t} \alpha_{t} \otimes \beta_{t}$ we have
$E_{\mu}^{13} = \sum_{t} \alpha_{t} \otimes \mathrm{id} \otimes \beta_{t}$ and
$E_{\mu}^{23}= \sum_{t}  \mathrm{id} \otimes \alpha_{t} \otimes \beta_{t}$.
Note that (\ref{chapter2:R-matrix(99)}) uses  the result
$$(\pi_{\lambda} \otimes \pi_{\nu} \otimes \pi_{\mu})
        \left[ (\Delta \otimes 1) E_{\mu}  \cdot (E_{\mu}^{23})^{-1} \right] 
	= (\pi_{\lambda} \otimes \pi_{\nu} \otimes \pi_{\mu}) E_{\mu}^{13},$$
rather than an equality in $U_{q}(\mathfrak{g})^{\otimes 3}$.
Similarly, we have
$$(1 \otimes \Delta) R_{\lambda,\mu} = 
\left[(1 \otimes \Delta) E_{\mu} \right] \cdot \Psi_{12} 
(\widetilde{R}_{13}) \cdot \widetilde{R}_{12},$$
and 
\begin{equation}
\label{chapter2:R-matrix(100)}
(\pi_{\lambda} \otimes \pi_{\nu} \otimes \pi_{\mu}) \big[ (1 \otimes \Delta) R_{\lambda,\mu} \big]
= (\pi_{\lambda} \otimes \pi_{\nu} \otimes \pi_{\mu})
      \left[ E^{13}_{\mu} \widetilde{R}_{13} E_{\nu}^{12} \widetilde{R}_{12} \right]. 
      \end{equation}
Together with Theorem \ref{th:RRRmatrix}, this shows that  
$R_{\lambda,\mu}$ satisfies the defining relations (\ref{eq:asdf;lkj})--(\ref{eq:deddy}) 
of the universal $R$-matrix of a $\mathbb{Z}_{2}$-graded quasitriangular Hopf algebra 
in {\emph{each triple tensor product}} of finite dimensional
irreducible $U_{q}(\mathfrak{g})$ representations, 
if one carefully chooses the limits in the definition of $E_{\mu}$ 
(which can always be done).  
Furthermore,  
Eqs. (\ref{chap2:equationannie(alpha)}) and 
(\ref{chapter2:R-matrix(99)})--(\ref{chapter2:R-matrix(100)}) imply that
\begin{equation}
\label{chapter2:R-matrix(101)}
(\pi_{\lambda} \otimes \pi_{\lambda} \otimes \pi_{\lambda}) R_{12} R_{13} R_{23} 
 = (\pi_{\lambda} \otimes \pi_{\lambda} \otimes \pi_{\lambda})   R_{23}R_{13}R_{12} ,
\end{equation}
where we have fixed $R = R_{\lambda,\lambda}$.

For later use we note the following easily proved results.
Define an automorphism 
$\Psi_{m}: U_{q}(\mathfrak{g})^{\otimes m} \rightarrow U_{q}(\mathfrak{g})^{\otimes m}$ 
generalising the automorphism 
$\Psi: U_{q}(\mathfrak{g}) \otimes U_{q}(\mathfrak{g}) \rightarrow
U_{q}(\mathfrak{g}) \otimes U_{q}(\mathfrak{g})$ in Lemma \ref{lemlem:twenty} by
$$\begin{array}{lcl}
\Psi_{m} (1^{\otimes j} \otimes K_{i}^{\pm 1} \otimes 1^{\otimes (m-j-1)}) & = & 
   1^{\otimes j} \otimes K_{i}^{\pm 1} \otimes 1^{\otimes (m-j-1)}, \\
\Psi_{m} (1^{\otimes j} \otimes e_{i} \otimes 1^{\otimes (m-j-1)}) & = & 
   (K_{i}^{-1})^{\otimes j} \otimes e_{i} \otimes (K_{i}^{-1})^{\otimes (m-j-1)}, \\
\Psi_{m} (1^{\otimes j} \otimes f_{i} \otimes 1^{\otimes (m-j-1)}) & = & 
   (K_{i})^{\otimes j}  \otimes f_{i} \otimes (K_{i})^{\otimes (m-j-1)}, 
\end{array}$$
for each $1 \leq i \leq n$ and all $0 \leq j \leq m-1$. Then
\begin{eqnarray}
\big(\Delta \otimes \mathrm{id}^{\otimes t} \big) \Psi_{2,3,\ldots, t+1} (\widetilde{R}_{1(t+1)})
& = & \Psi_{2,3,\ldots, t+2}(\widetilde{R}_{1(t+2)}) \cdot 
   \Psi_{3,4,\ldots, t+2}(\widetilde{R}_{2(t+2)}), \label{chap2:eqmorganspurlock(1)} \\
\big( \mathrm{id}^{\otimes t} \otimes \Delta \big) \Psi_{1,2,\ldots, t} (\widetilde{R}_{1(t+1)})
& = & \Psi_{1,2,\ldots, t+1} (\widetilde{R}_{1(t+2)}) \cdot 
   \Psi_{1,2,\ldots, t} (\widetilde{R}_{1(t+1)}), \label{chap2:eqmorganspurlock(agent92)}
\end{eqnarray}
where  $\Psi_{k, \ldots, m} = \mathrm{id}^{\otimes (k-1)} \otimes \Psi_{m-k+1}$, 
$k \geq 2$, in (\ref{chap2:eqmorganspurlock(1)}) and
$\Psi_{1, \ldots, m} = \Psi_{m} \otimes \mathrm{id}$ in (\ref{chap2:eqmorganspurlock(agent92)}).
Then it may be easily shown that
\begin{eqnarray*}
\left( \pi^{\otimes t} \otimes \pi \right)  \big( \Delta^{(t-1)} \otimes \mathrm{id} \big)R 
 & = &
 \left( \pi^{\otimes t} \otimes \pi \right) R_{1(t+1)} R_{2(t+1)} \cdots R_{t(t+1)}, \\
\left( \pi \otimes \pi^{\otimes t} \right)  \big( \mathrm{id} \otimes \Delta^{(t-1)} \big)R 
 & = & 
  \left( \pi \otimes \pi^{\otimes t} \right) R_{1(t+1)} R_{1t} \cdots R_{12},
\end{eqnarray*}
where we fix $R = R_{V,V}$.

\end{subsection}

\end{section}

\begin{section}{Useful elements of $\overline{U}^{+}_{q}(\mathfrak{g})$ }
\markright{\text{Useful elements of $\overline{U}^{+}_{q}(\mathfrak{g})$ }}
\label{sec:'central'element}

Define a set of elements 
$\{u_{\lambda} \in \overline{U}^{+}_{q}(\mathfrak{g}) | \ \lambda \in {\cal{P}}^{+} \}$ by
$u_{\lambda} = \sum_{t} S(b_{\lambda_{t}}) a_{\lambda_{t}} (-1)^{[a_{\lambda_{t}}]}$
upon writing  $R_{\lambda,\lambda} = \sum_{t} a_{\lambda_{t}} \otimes b_{\lambda_{t}}$.
\begin{lemma}
\label{lemming:overtheygo}
The element $u_{\lambda}$ has the following properties:
\begin{itemize}
\item[(i)] $\epsilon(u_{\lambda})=1$,
\item[(ii)] $\pi_{\lambda}\big( S^{2}(x) u_{\lambda} \big) = \pi_{\lambda}\left( u_{\lambda}x \right)$, 
$\forall x \in U_{q}(\mathfrak{g})$.
\item[(iii)]  $\pi_{\lambda}\big(u_{\lambda} \widetilde{u}_{\lambda}\big) = \pi_{\lambda}(1)= 
\pi_{\lambda}\big(\widetilde{u}_{\lambda}u_{\lambda}\big)$, where 
$\widetilde{u}_{\lambda}$ is defined by
$\widetilde{u}_{\lambda} = \sum_{s} S^{-1}(d_{\lambda_{s}})c_{\lambda_{s}} (-1)^{[c_{\lambda_{s}}]}$,
where
$R_{\lambda,\lambda}^{-1} = \sum_{s} c_{\lambda_{s}} \otimes  d_{\lambda_{s}}$,
\item[(iv)] $(\pi_{\lambda} \otimes \pi_{\lambda})\big(\Delta(u_{\lambda})\big) = 
(\pi_{\lambda} \otimes \pi_{\lambda})
\left[(u_{\lambda} \otimes u_{\lambda})\big(R^{T}_{\lambda,\lambda}R_{\lambda,\lambda}\big)^{-1}\right]$,
\end{itemize} 
where $\pi_{\lambda}$ is the representation of $U_{q}(\mathfrak{g})$ afforded by the finite
dimensional irreducible $U_{q}(\mathfrak{g})$-module $V_{\lambda}$.
\end{lemma}
\begin{proof}
Part (i) is proved is by inspection.  To prove part (ii),
$u \in U_{h}(\mathfrak{g})$ can be written as 
$\displaystyle{
u = \sum_{t} c_{t} S(F_{t}) q^{-(\sum_{i=1}^{n} \mathbf{H}_{i}^{2} )} E_{t} (-1)^{[E_{t}]} }$, 
where the universal $R$-matrix is written as
$R = q^{\sum_{i=1}^{n} \mathbf{H}_{i} \otimes \mathbf{H}_{i}} \sum_{t} c_{t} E_{t} \otimes F_{t}$.
In $\overline{U}_{q}^{+}(\mathfrak{g})$,
$$u_{\lambda}= \sum_{t=0}^{\infty} 
c_{t} S(F_{t}) \left(\prod_{a=1}^{n} \sum_{b=p}^{s} P_{a}[b] (J_{a})^{-b} \right) E_{t} (-1)^{[F_{t}]}.$$
The proof follows by noting that
$\prod_{a=1}^{n} \sum_{b=p}^{s} P_{a}[b] (J_{a})^{-b}$ `implements' the action of 
$q^{-(\sum_{i=1}^{t} \mathbf{H}_{i}^{2} )}$ in the $U_{q}(\mathfrak{g})$ representation $\pi_{\lambda}$.
Part (iii) is proved
by formally undertaking the proof of $u \widetilde{u}=1$ in Lemma \ref{lem:theuoperator} 
for $\mathbb{Z}_{2}$-graded quasitriangular Hopf algebras but
using the expression for $R_{\lambda,\lambda}$ instead of $R$ and also using (ii) above.

The proof of part (iv) is similar to the proof of $\Delta(u) = (u \otimes u) \left(R^{T}R\right)^{-1}$
in $\mathbb{Z}_{2}$-graded quasitriangular Hopf algebras, but much of the calculation is done
in the representations of $U_{q}(\mathfrak{g})$ afforded by tensor products of
finite dimensional irreducible $U_{q}(\mathfrak{g})$-modules.  
Although the proof can be understood by following \cite[Lem. 1]{zg}, we give it here for completeness. 
Firstly 
$(\pi_{\lambda} \otimes \pi_{\lambda})R^{T}_{\lambda,\lambda} R_{\lambda,\lambda} 
\in End_{U_{q}(\mathfrak{g})}(V_{\lambda} \otimes V_{\lambda})$, thus
$$(\pi_{\lambda} \otimes \pi_{\lambda})
\Big[R_{\lambda,\lambda}^{T}R_{\lambda,\lambda} \Delta(u_{\lambda})\Big] = 
(\pi_{\lambda} \otimes \pi_{\lambda})
\sum_{t} (S \otimes S) \Delta'(b_{t}) \cdot R_{\lambda,\lambda}^{T}R_{\lambda,\lambda} \cdot \Delta(a_{t}) (-1)^{[a_{t}]}.$$
Introduce the operation
\begin{eqnarray*}
\lefteqn{
(x_{1} \otimes x_{2}) \circ (a_{1} \otimes a_{2} \otimes a_{3} \otimes a_{4}) } \\
& = & 
S(a_{3}) x_{1} a_{1} \otimes S(a_{4}) x_{2} a_{2} 
(-1)^{[x_{2}][a_{1}] + ([a_{3}]+[a_{4}])([x_{1}]+[x_{2}]+[a_{1}]+[a_{2}]) + [a_{4}]([x_{1}]+[a_{1}])},
\end{eqnarray*}
where $x_{1}, x_{2}, a_{i}$, $i=1, 2, 3, 4$, are homogeneous elements of $U_{q}(\mathfrak{g})$.  
Straightforward calculations show that 
$$(x_{1} \otimes x_{2}) \circ (bc) = \left[ (x_{1} \otimes x_{2}) \circ b \right] \circ c,
\hspace{10mm} \forall b, c \in U_{q}(\mathfrak{g})^{\otimes 4},$$
and that
$$\Delta(u_{\lambda}) R^{T}_{\lambda,\lambda} R_{\lambda,\lambda} = 
R^{21}_{\lambda,\lambda} \circ 
\left[R^{12}_{\lambda,\lambda} (\Delta \otimes \Delta') R_{\lambda,\lambda} \right].$$
\noindent
By using  Eqs. (\ref{chapter2:R-matrix(99)})--(\ref{chapter2:R-matrix(101)}) we obtain
\begin{eqnarray*}
(\pi_{\lambda} \otimes \pi_{\lambda}) 
 \Big[\Delta(u_{\lambda}) R^{T}_{\lambda,\lambda} R_{\lambda,\lambda}\Big]  & = & 
(\pi_{\lambda} \otimes \pi_{\lambda}) \left( 
R^{21}_{\lambda,\lambda} \circ \left[R^{12}_{\lambda,\lambda} R^{13}_{\lambda,\lambda} 
R^{23}_{\lambda,\lambda} R^{14}_{\lambda,\lambda} R^{24}_{\lambda,\lambda}  \right] \right), \\
& = & 
(\pi_{\lambda} \otimes \pi_{\lambda}) \left( 
R^{21}_{\lambda,\lambda} \circ \left[R^{23}_{\lambda,\lambda} R^{13}_{\lambda,\lambda} 
R^{12}_{\lambda,\lambda} R^{14}_{\lambda,\lambda} R^{24}_{\lambda,\lambda}  \right] \right).
\end{eqnarray*}
Straightforward calculations then show that  
$$(\pi_{\lambda} \otimes \pi_{\lambda}) \left( R^{21}_{\lambda,\lambda} \circ R^{23}_{\lambda,\lambda} \right)
= (\pi_{\lambda} \otimes \pi_{\lambda}) (1 \otimes 1), \hspace{10mm} 
(1 \otimes 1) \circ R^{13}_{\lambda,\lambda} = u_{\lambda} \otimes 1,
\hspace{5mm} \mbox{and}$$
$$(\pi_{\lambda} \otimes \pi_{\lambda}) 
\left[ (u_{\lambda} \otimes 1) \circ (R^{12}_{\lambda,\lambda} R^{14}_{\lambda,\lambda}) \right] = 
(\pi_{\lambda} \otimes \pi_{\lambda}) (u_{\lambda} \otimes 1),$$
and thus it follows that
$$(\pi_{\lambda} \otimes \pi_{\lambda}) 
  \Big[ \Delta(u_{\lambda}) R^{T}_{\lambda,\lambda} R_{\lambda,\lambda} \Big]
 = (\pi_{\lambda} \otimes \pi_{\lambda}) \left[(u \otimes 1) \circ R^{24}_{\lambda,\lambda} \right] = 
 (\pi_{\lambda} \otimes \pi_{\lambda})(u_{\lambda} \otimes u_{\lambda}),$$
 completing the proof.

\end{proof}

Define a set of elements 
$\{v_{\lambda} \in \overline{U}^{+}_{q}(\mathfrak{g}) | \ \lambda \in {\cal{P}}^{+} \}$ by
\begin{equation}
\label{eq:defofthevelements}
v_{\lambda} = u_{\lambda} K_{2 \rho}^{-1}.
\end{equation}
\begin{lemma}
\label{lem:galliano}
The element $v_{\lambda}$ has the following properties:
$$\epsilon(v_{\lambda})=1, \hspace{10mm} \pi_{\lambda}(v_{\lambda}x)=\pi_{\lambda}(x v_{\lambda}), 
\hspace{5mm} \forall x \in U_{q}(\mathfrak{g}),$$
\begin{equation}
\label{eq:laurasia}
(\pi_{\lambda} \otimes \pi_{\lambda})
\Delta(v_{\lambda})=(\pi_{\lambda} \otimes \pi_{\lambda})
\left[(v_{\lambda} \otimes v_{\lambda})\big(R_{\lambda,\lambda}^{T}R_{\lambda,\lambda}\big)^{-1}\right].
\end{equation}
\end{lemma}
\begin{proof}
The proofs of $\epsilon(v_{\lambda})=1$ and (\ref{eq:laurasia}) follow from the definition of 
$v_{\lambda}$.  To prove the remaining relation, note that
$S^{2}(e_{i}) = K_{i}e_{i}K_{i}^{-1} = K_{2\rho}e_{i}K_{2\rho}^{-1}$,
$S^{2}(f_{i}) = K_{i}f_{i}K_{i}^{-1} = K_{2\rho}f_{i}K_{2\rho}^{-1}$,
and $S^{2}(K^{\pm 1}_{i}) = K_{2\rho}K_{i}^{\pm 1}K_{2\rho}^{-1}$.
As $S^{2}$ is a homomorphism we have
$S^{2}(x) = K_{2\rho}x K_{2\rho}^{-1}$ for all $x \in U_{q}(\mathfrak{g})$ and then 
$$
\pi_{\lambda}\big(v_{\lambda} x v_{\lambda}^{-1}\big) 
= \pi_{\lambda} \big( u_{\lambda} K_{2\rho}^{-1} x K_{2\rho} u_{\lambda}^{-1}\big) 
= \pi_{\lambda} \big(u_{\lambda} S^{-2}(x) u^{-1}_{\lambda} \big) 
= \pi_{\lambda} \big( S^{2}(S^{-2}(x)) \big) 
= \pi_{\lambda}(x),
$$
completing the proof.
\end{proof}

Let $V_{\lambda}$ be an irreducible $U_{q}(\mathfrak{g})$-module with integral dominant highest weight
$\lambda$.
\begin{lemma}
\label{lem:markhoopkins(a)}
The element $v_{\lambda}$ acts on each vector in $V_{\lambda}$ 
as the multiplication by the scalar $q^{-(\lambda + 2\rho, \lambda)}$.
\end{lemma}  
\begin{proof}
Note that $v_{\lambda}$ is even and that 
$\pi_{\lambda}(v_{\lambda}) \in End_{U_{q}(\mathfrak{g})}(V_{\lambda})$.
Write $R_{\lambda,\lambda} = E_{\lambda} \widetilde{R}$ and
$\widetilde{R}=\sum_{t=0}^{\infty} a_{t} \otimes b_{t}$ where 
$a_{t} \in U_{q}(\mathfrak{b}_{+})$, $b_{t} \in U_{q}(\mathfrak{b}_{-})$ 
and $a_{0}=b_{0}=1$.  Then
$$\pi_{\lambda}(v_{\lambda}) = 
\pi_{\lambda}\left(\sum_{t=0}^{\infty} S(b_{t}) E a_{t} K_{2\rho}^{-1} (-1)^{[a_{t}]}\right),$$
where $E$ is an even element of $U_{q}(\mathfrak{g})$ satisfying
$E w_{\xi} = q^{-(\xi,\xi)} w_{\xi}$
for each weight vector $w_{\xi} \in V_{\lambda}$ of weight 
$\xi \in \bigoplus_{i=1}^{n} \mathbb{Z} \epsilon_{i}$.
Let $w_{\lambda}$ be a non-zero highest weight vector of $V_{\lambda}$, 
then $a_{t} w_{\lambda} = 0$ for all $t > 0$,  which yields 
$$v_{\lambda} \cdot w_{\lambda} = E K_{2\rho}^{-1} w_{\lambda} =
q^{-(\lambda+2\rho,\lambda)}w_{\lambda}.$$

As $V_{\lambda}$ is irreducible and
$\pi_{\lambda}(v_{\lambda}) \in End_{U_{q}(\mathfrak{g})}(V_{\lambda})$
(and $\pi_{\lambda}(v_{\lambda})$ is a homogeneous map of degree zero),
$v_{\lambda}$ acts on each weight vector $w$ in $V_{\lambda}$ as the multiplication by the claimed
scalar from Schur's lemma.
\end{proof} 

We may denote $q^{-(\lambda + 2\rho, \lambda)}$ by $\chi_{\lambda}(v_{\lambda})$.
Note that it may be true that
$v_{\mu}$ acts on each weight vector $w$ in $V_{\lambda}$ as the multiplication by the 
salar $q^{-(\lambda + 2\rho, \lambda)}$ even if  $\mu \neq \lambda$, and in this case 
we also write $\chi_{\lambda}(v_{\mu})$ to denote $q^{-(\lambda + 2\rho, \lambda)}$.

Following \cite{z3}, we define the quantum superdimension of the 
finite dimensional $U_{q}(\mathfrak{g})$-module $W$ to be 
$$sdim_{q}(W) = \mathrm{str} \big(\pi_{W} (K_{2\rho})\big).$$

We now prove the following lemma originally given in \cite[p. 323]{zsuper}.
\begin{lemma}
Let $V_{\lambda}$ be a finite dimensional irreducible 
$U_{q}(\mathfrak{g})$-module with integral dominant highest weight $\lambda$. 
The quantum superdimension of $V_{\lambda}$ is
\begin{equation}
\label{eq:qqqsuperdim}
sdim_{q}(V_{\lambda})
= (-1)^{[\lambda]}q^{-(\lambda,2\rho)}
\prod_{\alpha \in \overline{\Phi}^{+}_{0}}\left(
\frac{q^{2(\lambda+\rho,\alpha)}-1}{q^{2(\rho,\alpha)}-1}\right)
\prod_{\beta \in \Phi_{1}^{+}}\left(
\frac{q^{2(\lambda+\rho,\beta)}+1}{q^{2(\rho,\beta)}+1}\right), 
\end{equation}
where $[\lambda]$ is the grading of the highest weight vector of $V_{\lambda}$.
\end{lemma}
\begin{proof}

We follow the usual proof for the formula for the
quantum dimension of a finite dimensional irreducible
module over a quantum algebra.  
We have not seen the proof of Eq. (\ref{eq:qqqsuperdim}) in the literature
and so write it down explicitly here.

The weight space decomposition of the $U_{q}(\mathfrak{g})$-module $V_{\lambda}$ is the same
as for a $U(\mathfrak{g})$-module with highest weight $\lambda$, 
so we can use Kac's supercharacter formula 
(see Appendix \ref{appendixD:title}) to calculate $\mathrm{str} \big( \pi_{\lambda}(K_{2\rho}) \big)$.  
Fix $e^{h} \in \mathbb{C}$ by $e^{h}=q$ and 
define in the notations of Appendix \ref{appendixD:title}  
a homomorphism $f: H^{*} \times E' \rightarrow \mathbb{C}[[h]]$ by
$f_{\eta}(e^{\nu}) = e^{h(\eta,\nu)}$, $\eta, \nu \in H^{*}$.

By definition, $sch_{\lambda} = \sum_{\mu} (-1)^{[\mu]} m(\mu) e^{\mu}$, 
where the sum is over all weight spaces of
$V_{\lambda}$, where $[\mu]$ is the grading of the vectors in the weight space $\mu$, 
and $m(\mu)$ is the multiplicity of the weight space $\mu$.  
Applying  $f_{2\rho}$ to $sch_{\lambda}$ gives 
$$f_{2\rho}(sch_{\lambda}) = \sum_{\mu} (-1)^{[\mu]} m(\mu) e^{h(2\rho,\mu)},$$ 
which is just $sdim_{q}(V_{\lambda})$.

Applying  $f_{2\rho}$ to $sch_{\lambda}$, and 
using the variant of Weyl's denominator formula in Appendix \ref{appendixD:title}, gives  
 \begin{eqnarray}
 \lefteqn{
\prod_{\alpha \in \overline{\Phi}^{+}_{0}} (e^{h(\alpha,\rho)} - e^{-h(\alpha,\rho)}) 
\prod_{\beta \in \Phi^{+}_{1}}(e^{h(\beta,\rho)}+e^{-h(\beta,\rho)}) 
\sum_{\mu} (-1)^{[\mu]} m(\mu) e^{h(2\rho,\mu)} } \nonumber \\
&  & \hspace{40mm} =   
(-1)^{[\lambda]} \sum_{\sigma \in {\cal{W}}} \epsilon'(\sigma) e^{h(2\rho,\sigma(\lambda + \rho))}, \nonumber  \\
& & \hspace{40mm} = (-1)^{[\lambda]} f_{2\lambda + 2\rho} \sum_{\sigma \in {\cal{W}}} \epsilon'(\sigma) e^{\sigma(\rho)}.
\label{chap2:quantumsuper(1)dean}
\end{eqnarray}
Rewriting the right hand side of (\ref{chap2:quantumsuper(1)dean}) as
 $$ (-1)^{[\lambda]}
\prod_{\alpha \in \overline{\Phi}^{+}_{0}} (e^{h(\alpha,\lambda + \rho)} - e^{-h(\alpha, \lambda + \rho)}) 
\prod_{\beta \in \Phi^{+}_{1}}(e^{h(\beta,\lambda + \rho)}+e^{-h(\beta,\lambda + \rho)}),$$
yields the desired result.
\end{proof}

The quantum superdimension of the fundamental $U_{q}(osp(1|2n))$-module $V$ is
\begin{equation}
\label{eq:johnfaulkner(12)}
sdim_{q}(V) = 1 - \frac{q^{2n}-q^{-2n}}{q-q^{-1}},
\end{equation}
where the grading of the highest weight vector of $V$ is odd.

\end{section}

\begin{section}{Spectral decomposition of $\check{\cal{R}}_{V,V}$}
\label{sec:spectam(a)}
\markright{\text{Spectral decomposition of $\check{\cal{R}}_{V,V}$}}

Let $V_{\lambda}$ and $V_{\mu}$ be finite dimensional irreducible
$U_{q}(\mathfrak{g})$-modules with integral dominant highest weights $\lambda$ and $\mu$, respectively.
Let $R_{\lambda,\mu}$ be as in Theorem \ref{th:RRRmatrix} and define 
$\check{\cal{R}}_{V_{\lambda}, V_{\mu}} 
\in \mathrm{Hom}_{\mathbb{C}}(V_{\lambda} \otimes V_{\mu}, V_{\mu} \otimes V_{\lambda})$ by
\begin{equation}
\label{eq:bigjimmyboy(2)}
\check{\cal{R}}_{V_{\lambda},V_{\mu}}(v_{\lambda} \otimes v_{\mu}) = 
P \circ \big(R_{\lambda,\mu}(v_{\lambda} \otimes v_{\mu})\big),
\end{equation}
where $v_{\lambda} \in V_{\lambda}$ and $v_{\mu} \in V_{\mu}$ are weight vectors
and $P$ is the graded permutation operator.

\begin{lemma}
Let $V_{\lambda}$ be an irreducible $U_{q}(\mathfrak{g})$-module  
with integral dominant highest weight $\lambda$, then
$$\check{\cal{R}}_{V_{\lambda},V_{\lambda}} \in End_{U_{q}(\mathfrak{g})}(V_{\lambda} \otimes V_{\lambda}).$$
\end{lemma}
\begin{proof}
We follow  Lemma \ref{lem:113} (i) to prove that
$$ \check{\cal{R}}_{V_{\lambda},V_{\lambda}} \cdot 
(\pi_{\lambda} \otimes \pi_{\lambda})\big(\Delta(x)\big)
 = (\pi_{\lambda} \otimes \pi_{\lambda})\big(\Delta(x)\big) \cdot 
 \check{\cal{R}}_{V_{\lambda},V_{\lambda}},$$
then the result follows.
\end{proof}
For $n=1$, $V \otimes V$ decomposes into a direct sum of irreducible 
$U_{q}(\mathfrak{g})$-modules  (see \cite{zsuper}):
\begin{equation}
\label{eq:johnfaulkner(14)}
V \otimes V = V_{2 \epsilon_{1}} \oplus V_{\epsilon_{1}} \oplus V_{0},
\end{equation}
and for $n \geq 2$, we have
\begin{equation}
\label{eq:johnfaulkner(13)}
V \otimes V = V_{2 \epsilon_{1}} \oplus V_{\epsilon_{1} + \epsilon_{2}} \oplus V_{0}.
\end{equation}
\begin{lemma}
\label{lem:ohgremlinsintheworks}
Let $n \geq 2$ and let
$\{ P[\mu] \in End_{U_{q}(\mathfrak{g})}(V \otimes V) | 
 \ \mu= 2\epsilon_{1}, \epsilon_{1} + \epsilon_{2}, 0 \}$ be a set of even 
$U_{q}(\mathfrak{g})$-linear maps: 
$P[\mu]: V \otimes V \rightarrow V \otimes V$,
where the image of $P[\mu]$ is $V_{\mu}$ and the maps
satisfy  $\big(P[\mu]\big)^{2} = P[\mu]$ and $P[\mu] P[\nu] = \delta_{\mu \nu} P[\mu]$.
Then there is a spectral decomposition of $\check{\cal{R}}_{V,V}$:
$$\check{\cal{R}}_{V,V} = -q P[2\epsilon_{1}] + q^{-1}P[\epsilon_{1} + \epsilon_{2}] + q^{-2n} P[0].$$
\end{lemma}
\begin{proof}
As $\check{\cal{R}}_{V,V} \in End_{U_{q}(\mathfrak{g})}(V \otimes V)$, 
we can write
$$\check{\cal{R}}_{V,V} = 
\beta_{2 \epsilon_{1}} P[2 \epsilon_{1}] + \beta_{\epsilon_{1}+\epsilon_{2}}P[\epsilon_{1}+\epsilon_{2}] +
\beta_{0}P[0],$$
for some set of constants 
$\{ \beta_{\mu} \in \mathbb{C} | \ \mu = 2\epsilon_{1}, \epsilon_{1}+\epsilon_{2}, 0 \}$,
where $\beta_{\mu}$ is the scalar action of $\check{\cal{R}}_{V,V}$ on 
the irreducible $U_{q}(\mathfrak{g})$-submodule $V_{\mu} \subset V \otimes V$.  
We explicitly calculate each $\beta_{\mu}$ using $R_{V,V}$.

Let $\{v_{i} | \ -n \leq i \leq n\}$ be the basis of 
weight vectors of $V$ given in Lemma \ref{lem:fundamentaldimensional}. 
The highest weight vector of $V_{2\epsilon_{1}}$ is
$w_{2\epsilon_{1}} = v_{1} \otimes v_{1}$,
the highest weight vector of $V_{\epsilon_{1} + \epsilon_{2}}$ is
$w_{\epsilon_{1} + \epsilon_{2}} = v_{1} \otimes v_{2} - q^{-1} v_{2} \otimes v_{1}$
and the highest weight vector of the trivial module $V_{0} \subset V \otimes V$ is 
$$w_{0}=\sum_{i=-n}^{n} c_{i} v_{i} \otimes v_{-i},$$
where $\{c_{i} \in \mathbb{C} | \ -n \leq i \leq n\}$ 
is a set of non-zero constants inductively defined by
$$
\begin{array}{rclrcl}
c_{n} & = & -c_{0}, & c_{-n} & = & q^{-1}c_{0}, \\
c_{n-1} & = & -qc_{n}, & c_{-(n-1)} & = & -q^{-1}c_{-n}, \\
c_{i} & = & -q c_{i+1}, & c_{-i} & = & -q^{-1} c_{-(i+1)},
\end{array} $$
where $i=1, 2, \ldots, n-2$ and we fix $c_{0} \neq 0$.

To study the action of $\check{\cal{R}}_{V,V}$ on the highest weight vectors 
$w_{2 \epsilon_{1}}, w_{\epsilon_{1} + \epsilon_{2}}$ and $w_{0}$, 
we note that the weight space decomposition of $V$ gives
rise to the following results:
\begin{eqnarray}
& &  \pi ( f_{\epsilon_{i}})^{3} = \pi ( e_{\epsilon_{i}})^{3}=0, \hspace{2mm} \mbox{for all }  
i=1, \ldots, n, \label{eq:cyclonetracey(1)} \\
& &  \pi ( f_{\gamma})^{2} = \pi ( e_{\gamma})^{2}=0, \hspace{2mm}
\mbox{for all }  \gamma \in \phi \mbox{ where } \gamma \neq \epsilon_{i}. \label{eq:cyclonetracey(2)}
\end{eqnarray}
Using (\ref{eq:cyclonetracey(1)}) and (\ref{eq:cyclonetracey(2)}), we have
\begin{equation}
\label{eq:johnfaulkner(15)}
(\pi \otimes \pi) \widetilde{R} = (\pi \otimes \pi) \prod_{\gamma \in \phi} \mathfrak{R}_{\gamma}
\end{equation}
where the product is ordered according to the same normal ordering ${\cal{N}}(\phi)$ used to construct
the root vectors
so that given ${\cal{N}}(\phi) = \gamma_{1} \prec \gamma_{2} \prec \cdots \prec \gamma_{k}$,  we fix
$\prod_{\gamma \in \phi} \mathfrak{R}_{\gamma} = 
  \mathfrak{R}_{\gamma_{1}} \mathfrak{R}_{\gamma_{2}} \cdots \mathfrak{R}_{\gamma_{k}}$,
where
\begin{eqnarray*}
\mathfrak{R}_{\gamma} & = & \left\{ \begin{array}{lll}
\displaystyle{\sum_{k=0}^{2} \frac{(q^{-1}-q)^{k} (e_{\gamma} \otimes f_{\gamma})^{k}}{ [k]^{-q^{-1}}!}},
& & \mbox{if } \gamma = \epsilon_{i}, \\
\displaystyle{\sum_{k=0}^{1} \frac{ (q-q^{-1})^{k} (e_{\gamma} \otimes f_{\gamma})^{k}}{ [k]^{q^{-2}}!  
} }, & &  \mbox{if } \gamma \neq \epsilon_{i}.
\end{array} \right. \\
& = & \left\{ \begin{array}{lll}
\displaystyle{ 1 \otimes 1 + (q^{-1}-q)(e_{\gamma} \otimes f_{\gamma}) + 
                       \frac{(q^{-1}-q)^{2}(e_{\gamma} \otimes f_{\gamma})^{2}}{(1-q^{-1})} },
 & & \mbox{if } \gamma = \epsilon_{i}, \\
\displaystyle{ 1 \otimes 1 + (q-q^{-1}) (e_{\gamma} \otimes f_{\gamma})   }, 
  & & \mbox{if }  \gamma \neq \epsilon_{i}.
\end{array} \right.
\end{eqnarray*}
Using this, 
we have
$\check{\cal{R}}_{V,V} (w_{2 \epsilon_{1}}) = -q w_{2 \epsilon_{1}}$ and
$\check{\cal{R}}_{V,V} (w_{\epsilon_{1} + \epsilon_{2}}) = q^{-1} w_{\epsilon_{1} + \epsilon_{2}}$.
Calculating $\beta_{0}$ is more difficult: note that
$$
\check{\cal{R}}_{V,V}\left( c_{-1} v_{-1} \otimes v_{1} + \sum_{\stackrel{i=-n}{i \neq -1}}^{n} c_{i} v_{i} \otimes
v_{-i}\right) = -q^{-1} c_{-1} v_{1} \otimes v_{-1} + \sum_{\stackrel{j=-n}{j \neq -1}}^{n} c_{j}'
v_{-j} \otimes v_{j},
$$
for some set of non-zero constants 
$\left\{c_{j}' \in \mathbb{C} | \ -n \leq j \leq n, \ j \neq -1 \right\}$.
Recall that $\check{\cal{R}}_{V,V} (w_{0}) = \beta_{0} w_{0}$, 
so we obtain $\beta_{0}$ by comparing $-q^{-1}c_{-1}$ and $c_{1}$.
Now $c_{-1} = (-1)^{n-1} q^{-n} c_{0}$ and $c_{1} = (-1)^{n} q^{n-1} c_{0}$, thus $\beta_{0} = q^{-2n}$.
\end{proof}

\begin{lemma}
\label{lem:ohgremlinsintheworks(2)}
Let $n=1$ and let
$\{ P[\mu] \in End_{U_{q}(\mathfrak{g})}(V \otimes V) | \ \mu= 2\epsilon_{1}, \epsilon_{1}, 0 \}$ 
be a set of even 
$U_{q}(\mathfrak{g})$-linear maps:
$P[\mu]: V \otimes V \rightarrow V \otimes V$,
where the image of $P[\mu]$ is $V_{\mu}$ and the maps
satisfy  $\big(P[\mu]\big)^{2} = P[\mu]$ and $P[\mu] P[\nu] = \delta_{\mu \nu} P[\mu]$.
Then there is a spectral decomposition of $\check{\cal{R}}_{V,V}$:
$$\check{\cal{R}}_{V,V} = -q P[2\epsilon_{1}] + q^{-1}P[\epsilon_{1}] + q^{-2} P[0].$$
\end{lemma}
\begin{proof}
The proof is almost identical to the proof of Lemma \ref{lem:ohgremlinsintheworks} 
with just the following minor difference.
The decomposition of $V \otimes V$ is given in (\ref{eq:johnfaulkner(14)})
and the highest weight vector of $V_{\epsilon_{1}}$ is 
$w_{\epsilon_{1}} = v_{1} \otimes v_{0} + q^{-1} v_{0} \otimes v_{1}$.
To complete the proof we note that $\check{\cal{R}}_{V,V} (w_{\epsilon_{1}}) = q^{-1} w_{\epsilon_{1}}$.
\end{proof}

From this we can write down the following important result:
\begin{corollary}
\label{cor:stonemebudgie}
For each $n \geq 1$, $\check{\cal{R}}_{V,V}$ satisfies 
\begin{equation}
\label{eq:RRRRequat}
(\check{\cal{R}}_{V,V}+q)(\check{\cal{R}}_{V,V}-q^{-1})(\check{\cal{R}}_{V,V}-q^{-2n})=0.
\end{equation}
\end{corollary}
Note that 
the expression for $(\pi \otimes \pi) \widetilde{R}$ in (\ref{eq:johnfaulkner(15)})
readily allows the use of $(\pi \otimes \pi) R_{V,V}$ and $\check{\cal{R}}_{V,V}$ in calculations.

\end{section}

\begin{section}{A representation of the Birman-Wenzl-Murakami algebra}
\label{eq:theXfactor(a)}
\markright{\text{A representation of the Birman-Wenzl-Murakami algebra}}

In this section we will define a map from the 
Birman-Wenzl-Murakami algebra $\mathscr{BW}_{f}(r,q)$ \cite{bw,w2} to
a subalgebra ${\cal{C}}_{f}$ of the centraliser algebra of $V^{\otimes f}$.  
This will allow us to map matrix units in $\mathscr{BW}_{f}(r,q)$ (see \cite{rw}) 
to matrix units of ${\cal{C}}_{f}$.

Before defining the map from $\mathscr{BW}_{f}(r,q)$ to ${\cal{C}}_{f}$, 
we need to obtain some preliminary results.
\begin{definition}
Define ${\cal{C}}_{t}$ to be the algebra over $\mathbb{C}$ generated by the elements
$$\left\{\check{R}_{i}^{\pm 1} \in End_{U_{q}(\mathfrak{g})}(V^{\otimes t}) |
 \ 1 \leq i \leq t-1\right\}, \hspace{5mm} \mbox{where}$$ 
\begin{equation}
\label{eq:tom3}
\check{R}_{i} = \mathrm{id}^{\otimes (i-1)} \otimes \check{\cal{R}}_{V,V}\otimes \mathrm{id}^{\otimes (t-(i+1))}
\in End_{U_{q}(\mathfrak{g})}(V^{\otimes t}).
\end{equation}
\end{definition}
\noindent

Let us investigate ${\cal{C}}_{t}$.
Let $\{v_{i} | \ -n \leq i \leq n\}$ be the basis of weight vectors
of $V$ given in Lemma \ref{lem:fundamentaldimensional} and let
$\{v^{*}_{i} | \ -n \leq i \leq n\}$ be a basis of $V^{*}$ such that
$\langle v^{*}_{i}, v_{j}\rangle = \delta_{ij}$ and $[v^{*}_{i}]=[v_{i}]$;
we have
$$av_{i} = \sum_{j} \langle v_{j}^{*}, av_{i} \rangle v_{j}, \hspace{10mm} 
av_{i}^{*} = \sum_{j} \langle av_{i}^{*}, v_{j} \rangle v_{j}^{*}, 
\hspace{10mm} \forall a \in U_{q}(\mathfrak{g}).$$
Define $\check{e} \in End_{\mathbb{C}}(V \otimes V^{*})$ by
$$\check{e}(x \otimes y^{*}) = 
(-1)^{[y^{*}][x]}\langle y^{*}, v^{-1}u \ x\rangle \sum_{i=-n}^{n} v_{i} \otimes v_{i}^{*},$$ 
where $v$ and $u$ are the elements
$v_{\epsilon_{1}}, u_{\epsilon_{1}} \in \overline{U}_{q}^{+}(\mathfrak{g})$ respectively.

\begin{lemma}
\label{lem:marrickville}
The map $\check{e}$ satisfies
\begin{itemize}
\item[(i)]  $(\check{e})^{2}=sdim_{q}(V)\check{e}$,
\item[(ii)] $a\check{e}=\epsilon(a) \check{e},$  $\hspace{5mm} \forall a \in U_{q}(\mathfrak{g})$,
\item[(iii)] $\check{e} a=\epsilon(a) \check{e}$, $\hspace{5mm} \forall a \in U_{q}(\mathfrak{g})$,
\item[(iv)] $\check{e}_{2} \check{R}_{1} \check{e}_{2} = q^{2n} \check{e}_{2},$
where 
$$\check{e}_{2} = \mathrm{id}_{V} \otimes \check{e}: V \otimes V \otimes V^{*}
\rightarrow V \otimes V \otimes V^{*},$$
$$\check{R}_{1} = \check{\cal{R}}_{V,V} \otimes \mathrm{id}_{V^{*}}: V \otimes V \otimes V^{*}
\rightarrow V \otimes V \otimes V^{*}.$$
\end{itemize}
\end{lemma}
\begin{proof}
\begin{itemize}
\item[(i)]
\begin{eqnarray*}
(\check{e})^{2}(x \otimes y^{*}) & = & 
(-1)^{[y^{*}][x]}\langle y^{*}, v^{-1}u x \rangle
\sum_{i} (-1)^{[v_{i}]}\langle v_{i}^{*}, v^{-1}u \ v_{i} \rangle \sum_{j} v_{j} \otimes v_{j}^{*} \\
& = & sdim_{q}(V) (-1)^{[y^{*}][x]}\langle y^{*}, v^{-1}u \ x \rangle
\sum_{j} v_{j} \otimes v_{j}^{*} 
  = sdim_{q}(V) \check{e} (x \otimes y^{*}).
\end{eqnarray*}
\item[(ii)] By definition, $a \check{e} = a_{V \otimes V^{*}} \circ \check{e}$; we calculate that
\begin{eqnarray*}
a \check{e}(x \otimes y^{*}) 
& = & (-1)^{[y^{*}][x]} \langle y^{*}, v^{-1}u x \rangle
\sum_{(a), i, j, k} \langle v_{j}^{*}, a_{(1)} v_{i} \rangle
\langle v_{i}^{*}, S(a_{(2)}) v_{k} \rangle v_{j} \otimes v_{k}^{*} \\
& = & (-1)^{[y^{*}][x]} \langle y^{*}, v^{-1}u x \rangle
\sum_{(a),j,k} \langle v_{j}^{*}, a_{(1)}S(a_{(2)}) v_{k} \rangle 
v_{j} \otimes v_{k}^{*} \\
& = & \epsilon(a) (-1)^{[y^{*}][x]} \langle y^{*}, v^{-1}u x \rangle 
\sum_{k} v_{k} \otimes v_{k}^{*} 
  =  \epsilon(a) \check{e}(x \otimes y^{*}).
\end{eqnarray*}
\end{itemize}
Similar calculations prove (iii) and (iv) (see \cite{lr} for the corresponding calculations in 
ungraded quasitriangular Hopf algebras).
\end{proof}

Define $\hat{e} \in End_{\mathbb{C}}(V^{*} \otimes V)$ by
$$\hat{e}(x^{*} \otimes y) = 
\langle x^{*}, y \rangle \sum_{i=-n}^{n} (-1)^{[v_{i}]} v_{i}^{*} \otimes vu^{-1} \ v_{i},$$
where $v$ and $u$ are set to be $v_{\epsilon_{1}}$ and $u_{\epsilon_{1}}$, respectively.

\begin{lemma}
\label{lemming:sausage2}
The map $\hat{e}$ satisfies
\begin{itemize}
\item[(i)]  $(\hat{e})^{2}=sdim_{q}(V)\hat{e}$,
\item[(ii)] $a \hat{e}=\epsilon(a) \hat{e}$, $\hspace{5mm} \forall a \in U_{q}(\mathfrak{g})$,
\item[(iii)] $\hat{e} a= \epsilon(a) \hat{e}$, $\hspace{5mm} \forall a \in U_{q}(\mathfrak{g})$,
\item[(iv)] $\hat{e}_{2} \check{R}_{1}^{-1} \hat{e}_{2} = q^{-2n} \hat{e}_{2}$
where 
$$\hat{e}_{2} = \mathrm{id}_{V^{*}} \otimes \hat{e}: V^{*} \otimes V^{*} \otimes V
\rightarrow V^{*} \otimes V^{*} \otimes V.$$
\end{itemize}
\end{lemma}
\begin{proof}
The proofs of (i)--(iv) are similar to the proofs of parts (i)--(iv) of Lemma \ref{lem:marrickville}.
\end{proof}

\begin{remark}  The maps $\check{e}$ and $\hat{e}$
are $U_{q}(\mathfrak{g})$-invariant maps onto one-dimensional $U_{q}(\mathfrak{g})$-submodules 
in $V \otimes V^{*}$ and $V^{*} \otimes V$, respectively.
\end{remark}

Recall that $V \otimes V$ has the decomposition into irreducible $U_{q}(\mathfrak{g})$-modules  
given in (\ref{eq:johnfaulkner(14)})--(\ref{eq:johnfaulkner(13)})
and that there exists an even $U_{q}(\mathfrak{g})$-invariant 
map $P[0]: V \otimes V \rightarrow V \otimes V$, the image of which is $V_{0} \subset V \otimes V$,
 defined in Lemmas \ref{lem:ohgremlinsintheworks}--\ref{lem:ohgremlinsintheworks(2)}.
Recall that $V^{*} \cong V$ and define 
$$E = (\mathrm{id} \otimes T^{-1}) \circ \check{e} \circ (\mathrm{id} \otimes T)
= (T^{-1} \otimes \mathrm{id}) \circ \hat{e} \circ (T \otimes \mathrm{id}) = sdim_{q}(V)P[0],$$
where $T$ is the isomorphism $T: V \rightarrow V^{*}$ given in Eq. (\ref{eq:johnfaulkner(20)}).
Furthermore, define the elements
$$E_{i} =  \mathrm{id}^{\otimes (i-1)} \otimes E \otimes \mathrm{id}^{\otimes (t-(i+1))}
\in End_{U_{q}(\mathfrak{g})}(V^{\otimes t}), \hspace{10mm} i=1, \ldots, t-1.$$
\begin{lemma}
\label{lem:10001}
The elements $\check{R}_{i}$, $E_{i} \in End_{U_{q}(\mathfrak{g})}(V^{\otimes t})$ satisfy the 
relations
\begin{itemize}
\item[(i)] $\check{R}_{i} \check{R}_{i+1} \check{R}_{i} = 
\check{R}_{i+1} \check{R}_{i} \check{R}_{i+1}$, $ \ 1 \leq i \leq t-2$,
\item[(ii)] $\check{R}_{i} \check{R}_{j} = \check{R}_{j} \check{R}_{i}$, $ \ |i-j|>1$,
\item[(iii)] $(\check{R}_{i}+q)(\check{R}_{i}-q^{-1})(\check{R}_{i}-q^{-2n})=0$, 
$ \ 1 \leq i \leq t-1$,
\item[(iv)] $-\check{R}_{i} + \check{R}_{i}^{-1} = (q-q^{-1})(1-E_{i})$,
\item[(v)] $E_{i} \check{R}_{i-1}^{\pm 1} E_{i} = q^{\pm 2n} E_{i}$, 
\item[(vi)] $E_{i} \check{R}_{i}^{\pm 1} = \check{R}_{i}^{\pm 1}E_{i}
=q^{\mp 2n} E_{i}$, $ \ 1 \leq i \leq t-1$.
\end{itemize}
\end{lemma}
\begin{proof}
The proofs of (i) and (ii) follow from Lemma \ref{lem:113} and 
the proof of (iii) follows from Corollary \ref{cor:stonemebudgie}.
The proof of (v) follows from Lemmas \ref{lem:marrickville}--\ref{lemming:sausage2}.
The proofs of (iv) and (vi) follow from the definition of 
$E_{i}$, Eq. (\ref{eq:johnfaulkner(12)})
and the fact that $\check{R}_{1}$ acts on $V_{0} \subset V \otimes V$ as 
a scalar multiple of the identity:
$\check{R}_{1}w = q^{-2n}w$ for all $w \in V_{0}$.
\end{proof}

We recall the definition of the Birman-Wenzl-Murakami algebra $\mathscr{BW}_{f}(r, q)$ from \cite{rw}.
Let $r$ and $q$ be non-zero complex parameters and let $f \geq 2$ be an integer.  
The Birman-Wenzl-Murakami algebra
$\mathscr{BW}_{f}(r, q)$ is the algebra over $\mathbb{C}$ generated by the invertible elements
$\{ g_{i} | \ 1 \leq i \leq f-1 \}$ subject to the relations
\newline
$$
\begin{array}{ll}
g_{i} g_{j} = g_{j} g_{i}, 			 & |i-j| >1, 	\\
g_{i} g_{i+1} g_{i} = g_{i+1} g_{i} g_{i+1}, 	 & 1 \leq i \leq f-2, 	\\
e_{i} g_{i} = r^{-1} e_{i},                      & 1 \leq i \leq f-1, 	\\
e_{i}g_{i-1}^{\pm 1}e_{i} = r^{\pm 1} e_{i}, 	 & 1 \leq i \leq f-1,
\end{array} $$
\newline
where $e_{i}$ is defined by
$$(q-q^{-1})(1-e_{i}) = g_{i}-g_{i}^{-1}, \hspace{10mm}  1 \leq i \leq f-1.$$
It can be shown that each $g_{i}$ also satisfies
$$(g_{i} - r^{-1})(g_{i} + q^{-1})(g_{i} - q) = 0.$$
From Lemma \ref{lem:10001} we have the following result:
\begin{lemma}
\label{eq:peacockinblackbeansauce}
Let $q \in \mathbb{C}$ be non-zero and not a root of unity.  
Then there is an algebra homomorphism 
$\Upsilon: \mathscr{BW}_{f}(-q^{2n}, q) \rightarrow {\cal{C}}_{f} \subseteq End_{U_{q}(\mathfrak{g})}(V^{\otimes f})$  defined by
$$\Upsilon: g_{i} \mapsto -\check{R}_{i}.$$
\end{lemma}

\end{section}

\begin{section}{Bratteli diagrams and path algebras}
\label{eq:caseydonovan(a)}
\markright{\text{Bratteli diagrams and path algebras}}

\begin{subsection}{Bratteli diagrams}

To proceed further with the study of $\mathscr{BW}_{f}(r, q)$ 
we now need the notions of {\emph{Bratelli diagrams}} and 
{\emph{Path algebras related to Bratelli diagrams}}, both of which we take 
from \cite{lr}.  (The reader is also referred to \cite[Chap. 2]{ghj}).

A Bratteli diagram is an undirected graph that encodes information about
a sequence $\mathbb{C} = A_{0} \subset A_{1} \subset A_{2} \subset \cdots$ of 
inclusions of finite dimensional semisimple algebras \cite{rw}.  
In a graph-theoretic sense, the properties of a Bratteli diagram are that
\begin{itemize}
\item[(i)]  the vertices are elements of certain sets $\widetilde{A}_{t}$, $t \in \mathbb{Z}_{+}$, and
\item[(ii)] if we
let $n(a,b) \in \mathbb{Z}_{+}$ denote the number of edges between the vertices $a$ and $b$, 
then $n(a,b)=0$ for any vertices $a \in \widetilde{A}_{i}$ and $b \in \widetilde{A}_{j}$ where 
$|i-j| \neq 1$.
\end{itemize}
Assume that $\widetilde{A}_{0}$ consists of a unique vertex we denote by $\emptyset$.  
We call the elements of  $\widetilde{A}_{t}$ {\emph{shapes}} and say that $\widetilde{A}_{t}$
is the set of shapes on the $t^{th}$ level of the Bratteli diagram.  
If $\lambda \in \widetilde{A}_{t}$ is
connected to $\mu \in \widetilde{A}_{t+1}$ in the Bratteli diagram we write $\lambda \leq \mu$.  

A {\emph{multiplicity free Bratteli diagram}} is a Bratteli diagram in which any two
vertices are connected by no more than one edge.  
All Bratteli diagrams considered in this thesis are multiplicity free. 

Let $A$ be a Bratteli diagram and let $\lambda \in \widetilde{A}_{r}$ and 
$\mu \in \widetilde{A}_{t}$ for
some $0 \leq r<t$.  We define a {\emph{path from $\lambda$ to $\mu$}} to be a sequence of shapes
$$P = (s_{r}, s_{r+1}, \ldots, s_{t}),$$
where $\lambda = s_{r} \leq s_{r+1} \leq \cdots \leq s_{t-1} \leq s_{t} = \mu$ 
and each $s_{i}$ is a shape on the $i^{th}$ level of $A$.

Given a path $T = (\lambda, \ldots, \mu)$ from $\lambda$ to $\mu$ and a path
$S = (\mu, \ldots, \nu)$ from $\mu$ to $\nu$ we define the concatenation 
of $T$ and $S$ to be the path from $\lambda$ to $\nu$ defined by
$$T \circ S = (\lambda, \ldots, \mu, \ldots, \nu).$$

We define a tableau $T$ of shape $\lambda$ to be a path 
from $\emptyset \in \widetilde{A}_{0}$ to $\lambda$ and we write
$shp(T) = \lambda$. We say that the length of $T$ is $t$
if there are $t+1$ shapes in the tableau.

\end{subsection}

\begin{subsection}{Path algebras related to Bratteli diagrams}

We now define the concept of a {\emph{Path algebra}} for a Bratteli diagram $A$.
For each $t \in \mathbb{Z}_{+}$, let ${\cal{T}}^{t}$  
be the set of tableaux of length $t$ in $A$
and let $\Omega^{t} \subset {\cal{T}}^{t} \times {\cal{T}}^{t}$ 
be the set of pairs $(S,T)$ of tableaux where
$shp(S)=shp(T)$, that is $S$ and $T$ both end in the same shape.  
Let us also define an algebra $A_{t}$ over $\mathbb{C}$ generated by  
$$\{E_{ST} | \ (S, T) \in \Omega^{t}\},$$ 
where the algebra multiplication is defined by 
$$E_{ST}E_{PQ} = \delta_{TP} E_{SQ}.$$
Note that $A_{0} \cong \mathbb{C}$.  
Each element $a \in A_{t}$ can be written in the form
$$a = \sum_{(S,T) \in \Omega^{t}} a_{ST} E_{ST}, \hspace{10mm} a_{ST} \in \mathbb{C}.$$  
We refer to the collection of algebras $A_{t}$, $t \in \mathbb{Z}_{+}$, as the 
{\emph{tower of path algebras corresponding to the Bratteli diagram $A$}}.

Each of the algebras $A_{t}$ is isomorphic to a direct sum of matrix algebras.
The irreducible representations of $A_{t}$ are indexed by the elements of $\widetilde{A}_{t}$,
that is, the set of shapes on the $t^{th}$ level of $A$.
Let ${\cal{T}}^{\lambda}$ denote the set of tableaux of shape $\lambda$, then 
the cardinality $d_{\lambda}$ of
${\cal{T}}^{\lambda} \cap {\cal{T}}^{t}$ is equal to the dimension of the irreducible $A_{t}$-module
indexed by $\lambda \in \widetilde{A}_{t}$.
We record this in the formula
$$A_{t} \cong \bigoplus_{\lambda \in \widetilde{A}_{t}} M_{d_{\lambda}} (\mathbb{C}),$$
where $M_{d} (\mathbb{C})$ denotes the algebra of $d \times d$ matrices with complex entries.

We now define some useful sets.
Let ${\cal{T}}^{\mu}_{\lambda}$ be the set of paths in $A$ from the shape $\lambda$ to the shape $\mu$
and let ${\cal{T}}^{t}_{r}$ be the set of paths 
starting on the $r^{th}$ level of $A$ and going down to the $t^{th}$ level. 
Furthermore, let ${\cal{T}}^{t}_{\lambda}$ be the set of paths in $A$ from the shape 
$\lambda$ to any shape on the $t^{th}$ level of $A$.

We also define 
$\Omega^{\mu}_{\lambda} \subset {\cal{T}}^{\mu}_{\lambda} \times {\cal{T}}^{\mu}_{\lambda}$ 
to be the set of pairs $(S,T)$ of paths $S, T \in {\cal{T}}^{\mu}_{\lambda}$ and 
$\Omega^{t}_{r} \subset {\cal{T}}^{t}_{r} \times {\cal{T}}^{t}_{r}$ 
to be the set of pairs $(S, T)$ of paths where in both situations we have $shp(S)=shp(T)$.

We define the inclusion of path algebras
$A_{r} \subseteq A_{t}$ for $0 \leq r < t$ as follows: for each
pair $(P, Q) \in \Omega^{r}$ we fix $E_{PQ} \in A_{t}$ by
$$E_{PQ} = \sum_{T \in {\cal{T}}^{t}_{\lambda} \cap {\cal{T}}^{t}_{r} } E_{P \circ T, Q \circ T}, \hspace{5mm}
\mbox{where } \lambda = shp(P) = shp(Q).$$
In particular, we have $A_{s} \subseteq A_{s+1}$ for each $s \in \mathbb{Z}_{+}$.

Let $\lambda \in \widetilde{A}_{t}$ and let ${\cal{V}}_{\lambda}$ be an irreducible representation of
$A_{t}$ indexed by $\lambda$.  The restriction of ${\cal{V}}_{\lambda}$ to 
the subalgebra $A_{t-1} \subseteq A_{t}$ decomposes into irreducible representations of $A_{t-1}$
according to
$${\cal{V}}_{\lambda} \downarrow^{A_{t}}_{A_{t-1}} \cong \bigoplus_{\mu \in \lambda^{-}}
{\cal{V}}_{\mu}, \hspace{10mm} \mbox{where } \lambda^{-} = \{ \nu \in \widetilde{A}_{t-1} | \ \nu \leq \lambda \}.$$
This decomposition is multiplicity free as the Bratteli diagram $A$ is multiplicity free. 

For each $r \in \mathbb{Z}_{+}$ satisfying $r < t$, 
the {\emph{centraliser of $A_{r}$ contained in $A_{t}$}} is defined by 
$${\cal{L}}(A_{r} \subseteq A_{t}) = \{ a \in A_{t} | \ ab=ba, \ \forall b \in A_{r} \}.$$

Let now $(S,T)$ be a pair of paths each starting on the $r^{th}$ level of $A$ at the shape 
$\lambda$ and ending on the $t^{th}$ level of $A$ at the shape $\mu$.  
For each such pair we define $E_{ST} \in A_{t}$ by
$$E_{ST} = \sum_{P \in {\cal{T}}^{\lambda} \cap {\cal{T}}^{r}} E_{ P \circ S, P \circ T},$$
which we can think of as the sum of all `paths' ending in $(S,T)$. 
We then have the
following lemma, stated in \cite[Prop. $(1.4)$]{lr} and proved in \cite[Sect. 2.3]{ghj}.
\begin{lemma}
\label{lemonelem:one}
A basis of ${\cal{L}}(A_{r} \subseteq A_{t})$ is given by the elements 
$$\{ E_{ST} | \ (S, T) \in \Omega_{\lambda}^{\mu} \cap \Omega^{t}_{r}, 
 \ \lambda \in \widetilde{A}_{r}, \ \mu \in \widetilde{A}_{t}\}.$$
\end{lemma}

\noindent
From Lemma \ref{lemonelem:one} we have the following corollary, 
which is proved in \cite[Cor. $(1.5)-(1.6)$]{lr}.
\begin{corollary}
Let the collection of algebras $A_{t}$, $t \in \mathbb{Z}_{+}$, 
be a tower of path algebras corresponding to a Bratteli
diagram $A$ and suppose that $g_{i} \in A_{i+1}$, $i \in \mathbb{N}$, are elements such that
\begin{itemize}
\item[(i)] for each $t$, the elements $\{g_{i} | \ 1 \leq i \leq t-1\}$ generate $A_{t}$, and
\item[(ii)] $g_{i} g_{j} = g_{j} g_{i}$ for all $i, j$ satisfying $|i-j| > 1$,
\end{itemize}
then
$$g_{t-1} = \sum_{(P,Q) \in \Omega^{t}_{t-2}}(g_{t-1})_{PQ} E_{PQ}, 
\hspace{5mm} (g_{t-1})_{PQ} \in \mathbb{C}.$$
Furthermore, if the elements $g_{1}, g_{2}, \ldots, g_{t-1}$ satisfy the relation
$g_{i} g_{i+1} g_{i} = g_{i+1} g_{i} g_{i+1}$ for all $1 \leq i \leq t-2$, then the element
$$D_{t} = g_{t-1} g_{t-2} \cdots g_{1} g_{1} \cdots g_{t-2} g_{t-1} \in A_{t}$$ 
can be expressed as
$$D_{t} = \sum_{S \in {\cal{T}}^{t}_{t-1} } D_{SS} E_{SS}, \hspace{10mm} D_{SS} \in \mathbb{C}.$$

\end{corollary}

\end{subsection}

\begin{subsection}{Centraliser algebras}
\label{subsec:brattelidiagramsandcentraliser}

Let $U$ be a $\mathbb{Z}_{2}$-graded Hopf algebra over $\mathbb{C}$.
Let $V$ be a finite dimensional $U$-module with the property that
$V^{\otimes t}$ is completely reducible for each $t \in \mathbb{Z}_{+}$.
We will now define the concepts of a {\emph{Bratteli diagram for tensor powers of $V$}} and
the {\emph{Bratteli diagram for $V^{\otimes t}$}}.
Furthermore, we aim to show that the centraliser 
${\cal{L}}_{t}$ of $U$ in $End_{\mathbb{C}}(V^{\otimes t})$ defined by
${\cal{L}}_{t} = End_{U}(V^{\otimes t})$
is isomorphic to the path algebra $A_{t}$ of the Bratteli diagram for $V^{\otimes t}$.

In this subsection we regard all modules as being graded.
By convention $V^{\otimes 0} \cong \mathbb{C}$ and thus ${\cal{L}}_{0} = \mathbb{C}$.
If $V$ is an irreducible $U$-module then ${\cal{L}}_{1} \cong \mathbb{C}$ by Schur's lemma.
For all $0 \leq r < t$ we define the inclusion ${\cal{L}}_{r} \subseteq {\cal{L}}_{t}$
by $a \mapsto a \otimes \mathrm{id}^{\otimes (t-r)}$ for all $a \in {\cal{L}}_{r}$.
Now ${\cal{L}}_{t}$ acts on $V^{\otimes t}$ in the obvious way.  
Since $U$ and ${\cal{L}}_{t}$ commute, $V^{\otimes t}$ has
a natural ${\cal{L}}_{t} \otimes U$-module structure.

Let $\{\Lambda_{\lambda} | \ \lambda \in I \}$ be the set of 
non-isomorphic finite dimensional irreducible $U$-modules.
Then by the double centraliser theorem there exists a finite subset 
$\widetilde{\cal{L}}_{t}$ of $I$ such that
$$V^{\otimes t} \cong \bigoplus_{ \lambda \in \widetilde{\cal{L}}_{t}}
{\cal{L}}^{\lambda} \otimes \Lambda_{\lambda},$$
where each ${\cal{L}}^{\lambda}$ is an irreducible ${\cal{L}}_{t}$-module such that 
${\cal{L}}^{\lambda} \not\cong {\cal{L}}^{\mu}$ if $\lambda \neq \mu$.

We now assume that $V$ is an irreducible $U$-module and we
consider the Bratteli diagram for tensor powers of $V$.
Let $\lambda \in \widetilde{\cal{L}}_{t}$ for some $t$.  
Then we have the decomposition
\begin{equation}
\label{eq:100}
\Lambda_{\lambda} \otimes V = 
\bigoplus_{\mu \in \widetilde{\cal{L}}_{t+1}} \left(\Lambda_{\mu}\right)^{\oplus n_{\lambda}(\mu)},
\hspace{10mm} n_{\lambda}(\mu) \in \mathbb{Z}_{+},
\end{equation}
of $\Lambda_{\lambda} \otimes V$ into a direct sum of irreducible $U$-modules.  
The non-negative integer $n_{\lambda}(\mu)$ 
is the multiplicity of $\Lambda_{\mu}$ in the decomposition.
We say that the decomposition of $\Lambda_{\lambda} \otimes V$ is
{\emph{multiplicity free}} if $n_{\lambda}(\mu) \leq 1$ for all $\mu \in \widetilde{\cal{L}}_{t+1}$.   

The {\emph{branching rule for inclusion}} ${\cal{L}}_{t} \subseteq {\cal{L}}_{t+1}$ describes the decomposition
of the ${\cal{L}}_{t+1}$-module ${\cal{L}}^{\nu}$ into ${\cal{L}}_{t}$-modules
\begin{equation}
\label{eq:101}
{\cal{L}}^{\nu} = \bigoplus_{\lambda \in \widetilde{\cal{L}}_{t}}
\big( {\cal{L}}^{\lambda} \big)^{\oplus n_{\lambda}(\nu)},
\hspace{10mm} n_{\lambda}(\nu) \in \mathbb{Z}_{+}.
\end{equation}
Note that the positive integers $n_{\lambda}(\nu)$ appearing in (\ref{eq:100}) and (\ref{eq:101}) are the same \cite{lr}.

The {\emph{Bratteli diagram for tensor powers of $V$}} is defined as follows:
for each $t \in \mathbb{Z}_{+}$ fix the vertices on the $t^{th}$ level of the Bratteli diagram
to be the elements of 
$\widetilde{\cal{L}}_{t}$.  
Then a vertex $\lambda \in \widetilde{\cal{L}}_{t}$
is connected to a vertex $\mu \in \widetilde{\cal{L}}_{t+1}$ by $n_{\lambda}(\mu)$ edges.

For a fixed $t$, the {\emph{Bratteli diagram for $V^{\otimes t}$}} is an undirected graph with
vertices given by the elements of $\bigcup_{i=0}^{t} \widetilde{\cal{L}}_{i}$, 
and the edges are such that a vertex $\lambda \in \widetilde{\cal{L}}_{i}$ is connected
to a vertex $\mu \in \widetilde{\cal{L}}_{i+1}$ by $n_{\lambda}(\mu)$ edges 
for each $0 \leq i \leq t-1$.

Let $V$ be a finite dimensional irreducible $U$-module 
with the property that for every irreducible $U$-module $W$, 
the decomposition of the tensor product $W \otimes V$ is multiplicity free.  
In this case, we say that tensoring by $V$ is multiplicity free.  
We will show that the centraliser algebra
${\cal{L}}_{t} = End_{U}(V^{\otimes t})$ is isomorphic to the path algebra $A_{t}$ associated with 
the Bratteli diagram for $V^{\otimes t}$.

We construct an algebra isomorphism $A_{t} \rightarrow {\cal{L}}_{t}$ inductively.  
Assume that there is an identification of ${\cal{L}}_{t}$ with 
the path algebra $A_{t}$ for some $t \geq 0$, so that
$$V^{\otimes t} = \bigoplus_{\lambda \in \widetilde{\cal{L}}_{t}}
\left( \bigoplus_{T \in {\cal{T}}^{\lambda} \cap {\cal{T}}^{t} } E_{TT} V^{\otimes t}\right)$$
is a decomposition of $V^{\otimes t}$ into irreducible $U$-modules $\Lambda_{\lambda}$
where the $U$-submodule $E_{TT} V^{\otimes t}$ is isomorphic to
$\Lambda_{\lambda}$ when $shp(T)=\lambda$.
The map $E_{TT}$ is a $U$-invariant map from $V^{\otimes t}$ onto a $U$-submodule
isomorphic to $\Lambda_{\lambda}$.

Let $T=( \emptyset, s_{1}, s_{2}, \ldots, \lambda) \in {\cal{T}}^{\lambda}$ be a tableau
of length $t$
and let $E_{TT} V^{\otimes t} \cong \Lambda_{\lambda}$ for some $\lambda \in \widetilde{\cal{L}}_{t}$.
As tensoring by $V$ is multiplicity free, the decomposition
\begin{equation}
\label{eq:equationa}
(E_{TT}V^{\otimes t}) \otimes V = 
\bigoplus_{\stackrel{\nu \in \widetilde{\cal{L}}_{t+1}}{\lambda \leq \nu}}V_{T \circ \nu},
\end{equation}
is multiplicity free and thus unique, where $T \circ \nu$ is the tableau
$$T \circ \nu = (\emptyset, s_{1}, s_{2}, \ldots, \lambda, \nu), \hspace{5mm} \lambda \leq \nu,$$
and $V_{T \circ \nu} \cong \Lambda_{\nu}$.

The next step is to identify
$E_{T \circ \nu, T \circ \nu}$ with the unique $U$-invariant projection operator mapping 
$\left(E_{TT}V^{\otimes t}\right) \otimes V$ onto $V_{T \circ \nu}$.  
This way we identify each element $E_{SS}$ of the path algebra
$A_{t+1}$, where $S \in {\cal{T}}^{t+1}$, with an element of ${\cal{L}}_{t+1}$.  
Thus we have the decomposition
$$V^{\otimes (t+1)} = \bigoplus_{\nu \in \widetilde{\cal{L}}_{t+1}}
\left( \bigoplus_{S \in {\cal{T}}^{\nu} \cap {\cal{T}}^{t+1}} E_{SS} V^{\otimes (t+1)} \right),$$
of $V^{\otimes (t+1)}$ into irreducible $U$-modules 
$E_{SS}V^{\otimes (t+1)} = V_{S} \cong \Lambda_{\nu}$, where
$\nu \in \widetilde{\cal{L}}_{t+1}$ and $S \in {\cal{T}}^{\nu} \cap {\cal{T}}^{t+1}$.

We now identify the other elements in the basis $\{ E_{PQ} \in A_{t+1} | \ (P,Q) \in \Omega^{t+1}\}$ 
with elements of ${\cal{L}}_{t+1}$.  
For each pair of paths $(P,Q) \in \Omega^{t+1}$ we choose non-zero elements
$$E_{PQ} \in E_{PP} {\cal{L}}_{t+1}E_{QQ}, 
\hspace{10mm} E_{QP} \in E_{QQ}{\cal{L}}_{t+1}E_{PP},$$
normalised in such a way that $E_{PQ}E_{QP} = E_{PP}$ and $E_{QP}E_{PQ}=E_{QQ}$.
Thus there is an algebra isomorphism $A_{t+1} \rightarrow {\cal{L}}_{t+1}$.

We then have the following theorem.
\begin{theorem}
\label{th:onelast}
Let $V$ be a finite dimensional irreducible
$U$-module such that $V^{\otimes t}$ is completely reducible for each $t \in \mathbb{Z}_{+}$ and
such that tensoring by $V$ is multiplicity free.
Then for any $t \in \mathbb{Z}_{+}$, the centraliser algebra ${\cal{L}}_{t} = End_{U} (V^{\otimes t})$
is isomorphic to the path algebra $A_{t}$ corresponding to the Bratteli diagram for $V^{\otimes t}$.
\end{theorem}

\end{subsection}

\end{section}

\begin{section}{Projections from $V^{\otimes t}$ onto irreducible $U_{q}(\mathfrak{g})$-modules}
\label{subsec:projectontome}
\markright{\text{Projections onto finite dimensional irreducible $U_{q}(osp(1|2n))$-modules}}

In this section we define projections from $V^{\otimes t}$ 
onto irreducible $U_{q}(\mathfrak{g})$-modules 
$V_{\lambda} \subset V^{\otimes t}$, $\lambda \in {\cal{P}}^{+}$, using elements of ${\cal{C}}_{t}$.
No substantially new results appear in this section, however, we are not aware of 
this specific formulation of the projections in the literature. 
In addition, this work provides a model for the definition of the 
$U_{q}^{(N)}(\mathfrak{g})$-modules in Chapter \ref{chap2A:titlelabel}.

Let $V_{\mu}$ be a finite dimensional
irreducible $U_{q}(\mathfrak{g})$-module with highest weight $\mu \in {\cal{P}}^{+}$.
Since each weight space of $V$ is one-dimensional, $V_{\mu} \otimes V$ is multiplicity free.
\begin{definition}
\label{def:twopointthreepointten}
We define ${\cal{P}}^{+}_{\mu} \subset {\cal{P}}^{+}$ to be the set such that for each
$\lambda \in {\cal{P}}^{+}_{\mu}$, $V_{\lambda}$ appears in $V_{\mu} \otimes V$
as an irreducible submodule.
\end{definition}
Now each $\lambda \in {\cal{P}}^{+}_{\mu}$ can only have one of the following three forms: 
$\mu, \mu+\epsilon_{i}, \mu - \epsilon_{i}$ for some $i$.
Thus
$${\cal{P}}^{+}_{\mu} \subseteq {\cal{P}}^{0}_{\mu} = \{ \mu, \mu \pm \epsilon_{i} \in {\cal{P}}^{+} | 
\ 1 \leq i \leq n \}.$$

\begin{definition}
Let $V_{\mu}$ be an irreducible $U_{q}(\mathfrak{g})$-module with highest weight
$\mu \in {\cal{P}}^{+}$, then 
$V_{\mu} \otimes V = \bigoplus_{\lambda \in {\cal{P}}^{+}_{\mu}} V_{\lambda}$.
Let $\big\{p_{\mu}[\lambda] \in End_{U_{q}(\mathfrak{g})}(V_{\mu} \otimes V) 
| \ \lambda \in {\cal{P}}^{+}_{\mu} \big\}$ 
be a set of even maps
$$p_{\mu}[\lambda]: V_{\mu} \otimes V \rightarrow V_{\mu} \otimes V$$
such that
\begin{itemize}
\item[(i)] the image of $p_{\mu}[\lambda]$ is $V_{\lambda}$,
\item[(ii)]  $\big(p_{\mu}[\lambda]\big)^{2} = p_{\mu}[\lambda]$,
\item[(iii)] $p_{\mu}[\lambda] \cdot p_{\mu}[\nu] = \delta_{\lambda \nu} p_{\mu}[\lambda]$.
\end{itemize}
We call each such $p_{\mu}[\lambda]$ a projection.
\end{definition}

Recall that for each integral dominant $\lambda$, there exists an element
$v_{\lambda} \in \overline{U}_{q}^{+}(\mathfrak{g})$ defined in
 Eq. (\ref{eq:defofthevelements}) that acts on each
 vector in the finite dimensional irreducible $U_{q}(\mathfrak{g})$-module $V_{\lambda}$ as the
 multiplication by the scalar $q^{-(\lambda + 2\rho, \lambda)}$.
 
For each $\mu \in {\cal{P}}^{+}$ and each $\lambda \in {\cal{P}}^{+}_{\mu}$, define 
$p_{\mu}[\lambda] \in End_{U_{q}(\mathfrak{g})}(V_{\mu} \otimes V)$ by
\begin{equation}
\label{eq:onetwofive}
p_{\mu}[\lambda] = (\pi_{\mu} \otimes \pi) 
\left( \prod_{\stackrel{\nu \in {\cal{P}}^{+}_{\mu}}{\nu \neq \lambda}} 
\frac{\Delta(v_{\xi})- q^{-(\nu + 2\rho, \nu)}}
{ q^{-(\lambda + 2\rho, \lambda)} - q^{-(\nu + 2\rho, \nu)}    }
\right),
\end{equation}
where $v_{\xi}$ is the element $v_{\lambda} \in {\cal{P}}^{+}$ with $\lambda = \xi$, 
for some integral dominant $\xi$
which is chosen so that $v_{\xi}$ 
acts as the multiplication by the scalar $q^{-(\nu + 2\rho, \nu)}$ on each vector in the
irreducible $U_{q}(\mathfrak{g})$-module $V_{\nu}$, for each $\nu \in {\cal{P}}^{+}_{\mu}$.  
For each integral dominant  $\mu$ there always exists at least one such $\xi$.
To see this, all we need is some
$E_{\xi}$ given by Eq. (\ref{eq:chatswoodtrain}):
$\displaystyle{E_{\xi} = \prod_{a=1}^{n} \sum_{b=p}^{s} (J_{a})^{b} \otimes P_{a}[b]},$
such that the element
$\displaystyle{E = \prod_{a=1}^{n} \sum_{b=p}^{s} P_{a}[b](J_{a})^{-b}}$
acts as the multiplication by the scalar $q^{-(\zeta, \zeta)}$
on each weight vector $w_{\zeta} \in V_{\nu} \subseteq V_{\mu} \otimes V$,
where $w_{\zeta}$ has the weight $\zeta$, 
and this is true for each $\nu \in {\cal{P}}^{+}_{\mu}$.

The element $E$ has this action whenever $s$ is sufficiently large enough and the
absolute value of the negative integer $p$ is sufficiently large enough, 
and so all we need do is to choose some $\xi$ for which this is true. 

To do this, for each $\nu \in {\cal{P}}^{+}_{\mu}$ let  
$I_{\nu}$ be the set of distinct weights of the weight vectors of $V_{\nu}$,
then $(\zeta_{\nu}, \epsilon_{i}) \in \mathbb{Z}$ for each weight $\zeta_{\nu} \in I_{\nu}$ and 
each $i=1, \ldots, n$.
Let
$$m = \max{ \big\{ |(\zeta_{\nu}, \epsilon_{i})| \in \mathbb{Z}_{+} \big| \ \zeta_{\nu} \in I_{\nu},
 \nu \in {\cal{P}}^{+}_{\mu}, i=1, \ldots, n  \big\}  },$$
 then fixing $\xi = \sum_{i=1}^{n} m \epsilon_{i}$ yields elements $E_{\xi}$ and $E$ 
 with the desired properties.

Note that  $(\pi_{\mu} \otimes \pi) \Delta(v_{\xi})$ in (\ref{eq:onetwofive}) is diagonalisable as
$V_{\mu} \otimes V$ is completely reducible and
$\Delta(v_{\xi})$ acts on each irreducible $U_{q}(\mathfrak{g})$-submodule 
$V_{\nu} \subset V_{\mu} \otimes V$ as the 
multiplication by the scalar $q^{-(\nu + 2\rho, \nu)}$.

\begin{lemma}
\label{lem:lemmings}
The maps $p_{\mu}[\lambda]$ are well-defined and satisfy
\begin{itemize}
\item[(i)] $\big(p_{\mu}[\lambda]\big)^{2}= p_{\mu}[\lambda]$,
\item[(ii)] $p_{\mu}[\lambda] \cdot p_{\mu}[\nu] = \delta_{\lambda \nu} p_{\mu}[\lambda]$,
\item[(iii)] $\displaystyle{\sum_{\lambda \in {\cal{P}}^{+}_{\mu}} p_{\mu}[\lambda] = 
              \mathrm{id}_{V_{\mu} \otimes V}}$.
\end{itemize}
\end{lemma}
\begin{proof} 
If $\alpha$ and $\beta$ are the highest weights of irreducible 
$U_{q}(\mathfrak{g})$-submodules in $V_{\mu} \otimes V$, then
$(\alpha + 2\rho, \alpha) = (\beta + 2\rho, \beta)$ implies that $\alpha = \beta$.
Then $p_{\mu}[\lambda]$ is well defined as tensoring by $V$ is multiplicity free and $q$ is not a root of unity.  
The proof of (i) follows from the result
$\big(p_{\mu}[\lambda]\big)^{2}(V_{\mu} \otimes V) = p_{\mu}[\lambda] (V_{\lambda}) = V_{\lambda}$.
For (ii) the case $\lambda = \nu$ reduces to (i), and for $\lambda \neq \nu$ we have
$$
p_{\mu}[\lambda] \cdot p_{\mu}[\nu] =
(\pi_{\mu} \otimes \pi)  \left(
\prod_{\stackrel{\lambda' \in {\cal{P}}^{+}_{\mu}}{\lambda' \neq \lambda}}
\frac{\Delta(v_{\xi})-q^{-(\lambda' + 2\rho, \lambda')}}
{q^{-(\lambda+ 2\rho, \lambda)}-q^{-(\lambda'+2\rho, \lambda')}}
\prod_{\stackrel{\nu' \in {\cal{P}}^{+}_{\mu}}{\nu' \neq \nu}}
\frac{\Delta(v_{\xi})-q^{-(\nu' + 2\rho, \nu')}}
{q^{-(\nu + 2\rho, \nu)}- q^{-(\nu' + 2\rho, \nu')}} \right) = 0. $$
\begin{itemize}
\item[(iii)] $\displaystyle{\sum_{\lambda \in{\cal{P}}^{+}_{\mu}}}
p_{\mu}[\lambda] \left(V_{\mu} \otimes V \right)=
\displaystyle{\bigoplus_{\lambda \in {\cal{P}}^{+}_{\mu}} } V_{\lambda} = V_{\mu} \otimes V$.
\end{itemize}
\end{proof}
\noindent
Note that 
$(\pi_{\mu} \otimes \pi) \displaystyle{\prod_{\lambda \in {\cal{P}}^{+}_{\mu}} 
\Big( \Delta(v_{\xi})-q^{-(\lambda + 2\rho, \lambda)} \Big) } = 0$.

We introduce some notation.
Let ${\cal{T}}^{t}$ be the set of tableaux of length $t$ derived from the Bratteli diagram for
$V^{\otimes t}$.  Let
$$i^{t} = (0, s_{1}, s_{2}, \ldots, s_{t}) \in {\cal{T}}^{t}.$$
We write $\lambda_{i}^{t} = i^{t}$ where $\lambda = s_{t}$.

Let $i^{t} \in {\cal{T}}^{t}$ and $s_{j}$, $s_{j+1} \in i^{t}$. Define a map 
$$p^{t-(j+1)}_{s_{j}}[s_{j+1}]:\left(V_{s_{j}} \otimes V\right)\otimes V^{\otimes
(t-(j+1))} \rightarrow V_{s_{j+1}} \otimes V^{\otimes (t-(j+1))}$$
 by
$$p^{t-(j+1)}_{s_{j}}[s_{j+1}] = p_{s_{j}}[s_{j+1}] \otimes \mathrm{id}^{\otimes
(t-(j+1))}.$$
\begin{lemma}
\label{lem:velvetvelvet}
The map   $p^{t-(j+1)}_{s_{j}}[s_{j+1}]$ satisfies
\begin{itemize}
\item[(i)]
$\left(p_{s_{j}}^{t-(j+1)}[s_{j+1}]\right)^{2} = p_{s_{j}}^{t-(j+1)}[s_{j+1}]$,
\item[(ii)]
$p_{s_{j}}^{t-(j+1)}[s_{j+1}] \cdot p_{s_{j}}^{t-(j+1)}[r_{j+1}]
= \delta_{s_{j+1},{r_{j+1}}} p_{s_{j}}^{t-(j+1)}[s_{j+1}],$
\item[(iii)]
$\displaystyle{\sum_{s_{j+1} \in {\cal{P}}^{+}_{s_{j}}}} p_{s_{j}}^{t-(j+1)}[s_{j+1}] 
= \mathrm{id}_{V_{s_{j}} \otimes V^{\otimes (t-j)}}$.
\end{itemize}
\end{lemma}
\begin{proof}
The proofs of parts (i) and (ii) follow from Lemma \ref{lem:lemmings} (i) and (ii), 
respectively.
The proof of (iii) follows from  Lemma \ref{lem:lemmings} (iii): explicitly, we have
$$\displaystyle{\sum_{s_{j+1} \in {\cal{P}}^{+}_{s_{j}}}
p_{s_{j}}^{t-(j+1)}[s_{j+1}]\cdot \left(V_{s_{j}}\otimes V \right)
\otimes V^{\otimes (t-(j+1))}
 = V_{s_{j}} \otimes V \otimes V^{\otimes t-(j+1)}}.$$
\end{proof}

\begin{definition}
\label{def:sartor}
Let $\tilde{p}_{i}^{t}[\lambda] \in End_{U_{q}(\mathfrak{g})}(V^{\otimes t})$ be a map
$\tilde{p}_{i}^{t}[\lambda]: V^{\otimes t} \rightarrow V_{\lambda} \subset V^{\otimes t}$ 
defined by
$$\tilde{p}_{i}^{t}[\lambda] =
p^{0}_{s_{t-1}}[\lambda]p^{1}_{s_{t-2}}[s_{t-1}]\cdots
p^{t-2}_{\epsilon_{1}}[s_{2}],$$
where $\lambda_{i}^{t} \in {\cal{T}}^{t}$. 
We say that {\emph{$\tilde{p}_{i}^{t}[\lambda]$ projects from $V^{\otimes t}$ onto
$V_{\lambda}$ by the path $\lambda_{i}^{t} \in {\cal{T}}^{t}$}} and 
we call $\tilde{p}_{i}^{t}[\lambda]$ a {\emph{path projection of length $t$}}.
\end{definition}
\begin{lemma}
\label{lem:greatgnats}
The map $\tilde{p}_{i}^{t}[\lambda]$ satisfies
\begin{itemize}
\item[(i)] $\big(\tilde{p}_{i}^{t}[\lambda]\big)^{2} = \tilde{p}_{i}^{t}[\lambda]$,
\item[(ii)] $\tilde{p}_{i}^{t}[\lambda] \cdot \tilde{p}_{j}^{t}[\lambda] = 
\left\{  
\begin{array}{ll}
0,                          & \mbox{if } i^{t} \neq j^{t}, \\
\tilde{p}_{i}^{t}[\lambda], & \mbox{if } i^{t} = j^{t},
\end{array}
\right. $
\item[(iii)]  $\tilde{p}_{i}^{t}[\lambda]\cdot \tilde{p}_{j}^{t}[\mu]=0$ if $\lambda \neq \mu$.
\end{itemize}
Furthermore, the map $\displaystyle{P_{t} = \sum_{i^{t} \in {\cal{T}}^{t}} \tilde{p}_{i}^{t}[\lambda]}$
is the identity on $V^{\otimes t}$.
\end{lemma}
\begin{proof}
\noindent
\begin{itemize}
\item[(i)] This follows from Lemma \ref{lem:velvetvelvet} (i).
\item[(ii)] For $i^{t} = j^{t}$ the case reduces to (i), let $i^{t} \neq j^{t}$ where
\begin{eqnarray*}
i^{t} & = & (0, s_{1}, s_{2}, \ldots, s_{k-1}, s_{k}, r_{k+1}, \ldots, r_{t-2}, r_{t-1}, \lambda), \\
j^{t} & = & (0, s_{1}, s_{2}, \ldots, s_{k-1}, s_{k}, s_{k+1}, \ldots, s_{t-2}, s_{t-1}, \lambda).
\end{eqnarray*}
Now $i^{t}$, $j^{t} \in {\cal{T}}^{t}$ and $r_{k+1} \neq s_{k+1}$ for some $2 \leq k+1 \leq t$, then
\begin{eqnarray*}  
\tilde{p}_{i}^{t}[\lambda]\cdot \tilde{p}_{j}^{t}[\lambda] & = &
p^{0}_{r_{t-1}}[\lambda]p^{1}_{r_{t-2}}[r_{t-1}]\cdots 
p^{t-k-1}_{s_{k}}[r_{k+1}]p^{t-k}_{s_{k-1}}[s_{k}]
p^{t-k+1}_{s_{k-2}}[s_{k-1}]\cdots
p^{t-2}_{s_{1}}[s_{2}] \\
&  &  \times
p^{0}_{s_{t-1}}[\lambda]p^{1}_{s_{t-2}}[s_{t-1}]\cdots 
p^{t-k-1}_{s_{k}}[s_{k+1}]p^{t-k}_{s_{k-1}}[s_{k}]
p^{t-k+1}_{s_{k-2}}[s_{k-1}] \cdots
p^{t-2}_{s_{1}}[s_{2}]  \\
& = & p^{0}_{r_{t-1}}[\lambda]p^{0}_{s_{t-1}}[\lambda]
p^{1}_{r_{t-2}}[r_{t-1}]p^{1}_{s_{t-2}}[s_{t-1}]
\cdots p^{t-k-1}_{s_{k}}[r_{k+1}]p^{t-k-1}_{s_{k}}[s_{k+1}]  \\
& &  \times  p^{t-k}_{s_{k-1}}[s_{k}]p^{t-k}_{s_{k-1}}[s_{k}] 
p^{t-k+1}_{s_{k-2}}[s_{k-1}]p^{t-k+1}_{s_{k-2}}[s_{k-1}]\cdots
p^{t-2}_{s_{1}}[s_{2}]p^{t-2}_{s_{1}}[s_{2}]  \\
& = & 0,
\end{eqnarray*}
as $p^{t-k-1}_{s_{k}}[r_{k+1}] \cdot p^{t-k-1}_{s_{k}}[s_{k+1}]=0$.
\item[(iii)]  Assume that
\begin{eqnarray*}
i^{t} & = & (0, s_{1}, s_{2}, \ldots, s_{k-1}, s_{k}, r_{k+1}, \ldots, r_{t-2}, r_{t-1}, \lambda), \\
j^{t} & = & (0, s_{1}, s_{2}, \ldots, s_{k-1}, s_{k}, s_{k+1}, \ldots, s_{t-2}, s_{t-1}, \mu),
\end{eqnarray*}
where $r_{k+1}\neq s_{k+1}$ for some $2 \leq k+1 \leq t$.
The calculations are similar to those of (ii) and we have
$\tilde{p}^{t}_{i}[\lambda]\cdot \tilde{p}^{t}_{j}[\mu]= 0$.
\end{itemize}
The last claim
follows inductively from the result in Lemma \ref{lem:lemmings} that
$\displaystyle{\sum_{\lambda \in {\cal{P}}^{+}_{\mu}} p_{\mu}[\lambda] = \mathrm{id}_{V_{\mu} \otimes V}}$.
\end{proof}

Recall that ${\cal{C}}_{t}$ is the algebra over $\mathbb{C}$ generated by the elements
$$\left\{\check{R}_{i}^{\pm 1} \in End_{U_{q}(\mathfrak{g})}(V^{\otimes t}) |
 \ 1 \leq i \leq t-1\right\}.$$ 
\begin{proposition}
For each path $\lambda_{i}^{t} \in {\cal{T}}^{t}$, $\tilde{p}_{i}^{t}[\lambda] \in {\cal{C}}_{t}$.
\end{proposition}
\begin{proof}
We prove the proposition inductively.  
Firstly, for some appropriately chosen integral dominant weight $\xi$, 
$$(\pi \otimes \pi)\Delta(v_{\xi}) = 
(\pi \otimes \pi)\left[(v_{\xi} \otimes v_{\xi})\left(R^{T}_{\xi,\xi} R_{\xi,\xi}\right)^{-1}\right]
 = q^{-2(\epsilon_{1} + 2\rho, \epsilon_{1})} \check{\cal{R}}^{-2} \in {\cal{C}}_{2}.$$
Now assume that $\tilde{p}^{(t-1)}_{i}[\mu] \in {\cal{C}}_{(t-1)}$ where
$\tilde{p}^{(t-1)}_{i}[\mu]$ is a path projection 
$\tilde{p}^{(t-1)}_{i}[\mu]: V^{\otimes (t-1)} \rightarrow V_{\mu}$ and
$V_{\mu}$ is an irreducible $U_{q}(\mathfrak{g})$-submodule of $V^{\otimes (t-1)}$.
We will show that  $(\pi_{\mu} \otimes \pi)\Delta(v_{\zeta})$ is an element of ${\cal{C}}_{t}$ 
for some appropriately chosen $\zeta$.
Let $\zeta$ be an integral dominant weight such that the element 
$v_{\zeta} \in \overline{U}_{q}^{+}(\mathfrak{g})$ acts as the multiplication by the scalar
$q^{-(\lambda + 2\rho, \lambda)}$ on each vector in the irreducible $U_{q}(\mathfrak{g})$-submodule
$V_{\lambda} \subset V_{\mu} \otimes V$ for each $\lambda \in {\cal{P}}^{+}_{\mu}$.
Now 
\begin{eqnarray*}
(\pi_{\mu} \otimes \pi) \Delta(v_{\zeta}) 
& = & (\pi_{\mu} \otimes \pi) \left[ (v_{\zeta} \otimes v_{\zeta}) \left(R^{T}_{\zeta,\zeta}R_{\zeta,\zeta}\right)^{-1} \right] \\
& = & q^{-(\mu + 2\rho, \mu)-(\epsilon_{1} + 2\rho, \epsilon_{1})} 
      (\tilde{p}^{(t-1)}_{i}[\mu] \otimes \mathrm{id}) \left(\pi^{\otimes (t-1)} \otimes \pi \right) 
      \left(\Delta^{(t-2)} \otimes \mathrm{id}\right) 
      \left(R^{T}_{\zeta,\zeta}R_{\zeta,\zeta}\right)^{-1} \\
& = & q^{-(\mu + 2\rho, \mu)-(\epsilon_{1} + 2\rho, \epsilon_{1})}  
      (\tilde{p}^{(t-1)}_{i}[\mu] \otimes \mathrm{id}) 
      \check{R}_{t-1}^{-1} \check{R}_{t-2}^{-1} \cdots \check{R}_{1}^{-1} \check{R}_{1}^{-1} \cdots
       \check{R}_{t-2}^{-1}\check{R}_{t-1}^{-1},
\end{eqnarray*}
where we have used the identity (fixing $R = R_{\zeta,\zeta}$):
$$ \left( \pi^{\otimes (t-1)} \otimes \pi \right) \left(\Delta^{(t-2)} \otimes \mathrm{id} \right)R = 
\left(\pi^{\otimes (t-1)} \otimes \pi \right) R_{1t} R_{2t} \cdots R_{(t-1)t},$$ 
which arises from (\ref{chap2:eqmorganspurlock(1)}).

\end{proof}

\end{section}

\begin{section}{Matrix units for ${\cal{C}}_{t}$}
\label{subsec:rhapsodyinred}
\markright{\text{Matrix units for ${\cal{C}}_{t}$}}

It is clear that the Bratteli diagram for $V^{\otimes t}$ is multiplicity free,
as tensoring by the fundamental $U_{q}(\mathfrak{g})$-module $V$ is multiplicity free.
It follows then from Theorem \ref{th:onelast} that the centraliser algebra 
${\cal{L}}_{t} = End_{U_{q}(\mathfrak{g})}(V^{\otimes t})$ is isomorphic to 
the path algebra $A_{t}$ obtained from the Bratteli diagram for $V^{\otimes t}$.
Clearly, we have the inclusion ${\cal{C}}_{t} \subseteq {\cal{L}}_{t}$.  
The aim of this section is to show that ${\cal{C}}_{t}$ and ${\cal{L}}_{t}$ are in fact equal:
\begin{theorem}
\label{th:areallybigthe(199)}
The centraliser algebra ${\cal{L}}_{t} = End_{ U_{q}(\mathfrak{g}) }\big(V^{\otimes t}\big)$ 
is generated by the elements
$$\big\{ \check{R}_{i}^{\pm 1} \in End_{U_{q}(\mathfrak{g})}\big(V^{\otimes t}\big) \big| \ i=1, 2,
\ldots, t-1 \big\}.$$
\end{theorem}

To prove this theorem we firstly partition the matrix units in $A_{t}$ 
into two groups:
the {\emph{projectors}}  $\{E_{SS} \in A_{t} | \ (S,S) \in \Omega^{t}\}$ and the 
{\emph{intertwiners}}  $\{E_{ST} \in A_{t} | \ (S,T) \in \Omega^{t}, S \neq T \}$
and we use an invertible homomorphism to map matrix units in $A_{t}$ 
to matrix units in ${\cal{C}}_{t}$.

Recall that $V^{\otimes t}$ is completely reducible.
Each matrix unit in ${\cal{C}}_{t}$ corresponding to a
projector in $A_{t}$ projects down from $V^{\otimes t}$ onto an irreducible 
$U_{q}(\mathfrak{g})$-submodule $V_{\lambda} \subset V^{\otimes t}$. 
Each matrix unit in ${\cal{C}}_{t}$ corresponding to an
intertwiner in $A_{t}$ maps between isomorphic irreducible 
$U_{q}(\mathfrak{g})$-submodules of $V^{\otimes t}$.

Recall that the homomorphism $\Upsilon: g_{i} \mapsto -\check{R}_{i}$, 
given in Lemma \ref{eq:peacockinblackbeansauce}, yields a representation of
$\mathscr{BW}_{t}(-q^{2n},q)$ in ${\cal{C}}_{t}$.
In Subsection \ref{subsect:matrixbirmanwenzlstuff} we will write down the matrix units in a semisimple
quotient of $\mathscr{BW}_{t}(-q^{2n},q)$ that map via $\Upsilon$ onto the 
projectors and intertwiners in ${\cal{C}}_{t}$.
We will do this for the intertwiners, 
but we choose to define the projectors more straightforwardly using our previous work.

The projections $E_{SS}$ that project down
from $V^{\otimes t}$ onto irreducible $U_{q}(\mathfrak{g})$-submodules 
$V_{shp(S)} \subset V^{\otimes t}$ that we defined in Section \ref{subsec:projectontome}, 
are elements of ${\cal{C}}_{t}$,
and they satisfy the equations satisfied by the projector matrix units: $(E_{SS})^{2} = E_{SS}$ and
$\sum_{S \in {\cal{T}}^{t}} E_{SS} = \mathrm{id}_{V^{\otimes t}}$. 
We fix the projectors in ${\cal{C}}_{t}$ as follows:  
we map the projector $E_{SS} \in A_{t}$ to $\tilde{p}_{i}^{t}[\lambda] \in {\cal{C}}_{t}$, where
$\lambda_{i}^{t} = S \in {\cal{T}}^{t}$ is a path of length $t$:
$E_{SS} \leftrightarrow \tilde{p}_{i}^{t}[\lambda]$.

All we need to do now is to
construct the matrix units in ${\cal{C}}_{t}$ corresponding to the intertwiners in $A_{t}$.
We denote the matrix unit in ${\cal{C}}_{t}$ corresponding to 
$E_{MP} \in A_{t}$ by the same label $E_{MP}$.  
This should not cause confusion: the precise meaning of $E_{MP}$ in any given situation will be clear.

\begin{subsection}{Matrix units in $\mathscr{BW}_{t}(-q^{2n},q)$}
\label{subsect:matrixbirmanwenzlstuff}

In this subsection, we say that an algebra $B$ is semisimple if it is isomorphic to a direct sum of
matrix algebras, ie $B \cong \bigoplus_{i \in I} M_{b_{i}}(\mathbb{C})$, 
where $M_{b_{i}}(\mathbb{C})$ is the algebra of $b_{i} \times b_{i}$ matrices with complex entries.
The algebra $\mathscr{BW}_{t}(-q^{2n},q)$ is not semisimple at generic $q$ \cite[Cor. 5.6]{w2} but
Ram and Wenzl have constructed matrix units
for the semisimple Birman-Wenzl-Murakami algebra $\mathscr{BW}_{t}$ defined over $\mathbb{C}(r,q)$
(the field of rational functions in $r$ and $q$) for indeterminates $r$ and $q$ \cite{rw}.
By replacing the indeterminates $r$ and $q$ with the complex numbers $-q^{2n}$ and $q$,
respectively, we obtain matrix units in a semisimple quotient of $\mathscr{BW}_{t}(-q^{2n},q)$.
By then applying the map $\Upsilon$ to these matrix units, we obtain matrix units in
${\cal{C}}_{t}$.

Before doing this, let us introduce Young diagrams 
and discuss a relation between certain Young diagrams
and the integral dominant highest weights of irreducible $U_{q}(osp(1|2n))$-modules.
For each non-negative integer $m$, there exists 
a Young diagram for each partition of $m$.
Let $m = m_{1} + m_{2} + \cdots + m_{l}$ be a partition of $m$, where 
$m_{i} - m_{i+1} \in \mathbb{Z}_{+}$ for each
$i=1, 2, \ldots, l-1$ and $m_{l} \in \mathbb{Z}_{+}$.
The Young diagram representing this partition is a collection of $m$ boxes arranged in $l$
left-aligned rows where the $i^{th}$ row from the top  contains exactly $m_{i}$ boxes. 
If $m \geq 1$, 
let $c_{i}$, $i = 1, 2, \ldots, m_{1}$, 
be the number of boxes in the $i^{th}$ column from the left in the Young diagram,
then $c_{i} - c_{i+1} \in \mathbb{Z}_{+}$ for each $i=1, 2, \ldots, m_{1}-1$
and $c_{m_{1}} \in \{1, 2, \ldots, l\}$.

Recall that an integral dominant highest weight $\lambda$
of an irreducible $U_{q}(osp(1|2n))$-module $V_{\lambda}$ is of
the form $\lambda = \sum_{i=1}^{n} \lambda_{i} \epsilon_{i} \in {\cal{P}}^{+}$ where 
$\lambda_{i} - \lambda_{i+1} \in \mathbb{Z}_{+}$ for each $i=1, 2, \ldots, n-1$ and
$\lambda_{n} \in \mathbb{Z}_{+}$. 
We can use a Young diagram to label the highest weight of $V_{\lambda}$:
this Young diagram consists of $\sum_{i=1}^{n} \lambda_{i}$ boxes arranged in $n$ left-aligned rows,
where the $i^{th}$ row contains exactly $\lambda_{i}$ boxes.
Let $\mu$ be a Young diagram containing no more than $n$ rows of boxes and let
$\mu_{i}$ be the number of boxes in the $i^{th}$ row from the top 
for each $i=1, 2, \ldots, n$.
We can use $\mu$ to label an integral dominant highest weight of an irreducible
representation of $U_{q}(osp(1|2n))$: the integral dominant 
highest weight that $\mu$ represents is $\sum_{i=1}^{n} \mu_{i} \epsilon_{i} \in {\cal{P}}^{+}$.

\begin{subsubsection}{The algebra $\mathscr{BW}_{t}$}

The Birman-Wenzl-Murakami algebra $\mathscr{BW}_{t}$, 
with $r$ and $q$ indeterminates, is equipped with a functional 
$\mathrm{tr}: \mathscr{BW}_{t} \rightarrow \mathbb{C}(r,q)$ which satisfies, 
amongst other relations \cite[Lem. 3.4 (d)]{w2},
\begin{equation}
\label{eq:tracefunctionalnew}
\mathrm{tr}(a \chi b) = \mathrm{tr}(\chi) \mathrm{tr}(ab), 
\hspace{10mm} \forall a, b \in \mathscr{BW}_{t-1}, \hspace{5mm} \chi \in \{g_{t-1}, e_{t-1}\},
\end{equation}
where we regard each element of $\mathscr{BW}_{t-1}$ as an element of $\mathscr{BW}_{t}$
under the obvious inclusion.
The algebra $\mathscr{BW}_{t}$ is semisimple \cite[Thm. 3.5]{w2}.
To discuss its structure, we introduce the Young lattice.

The Young lattice is the following infinite graph  \cite[Sec. 1]{w2}.  
The vertices of the Young lattice are the Young diagrams;
the vertices are grouped into levels so that each Young diagram with exactly $t$ boxes
labels a vertex on the $t^{th}$ level of the Young lattice.
The edges of the Young lattice are completely determined as follows: 
a vertex $\lambda$ on the $t^{th}$ level is connected to a vertex $\mu$ on the 
$(t+1)^{st}$ level by one edge if and only if $\lambda$ and $\mu$ differ by exactly one box.
We show the Young lattice up to the $4^{th}$ level in Figure \ref{fig:young}, 
where the circle represents the Young diagram with no boxes.
We say that the level containing the Young diagram with no boxes is the $0^{th}$ level.
Note that the Young lattice is (apart from the $0^{th}$ level) 
identical to the Bratteli diagram for the 
sequence of inclusions of group algebras of the symmetric group:
 $\mathbb{C}S_{1} \subset \mathbb{C}S_{2} \subset  \mathbb{C}S_{3} \subset
\cdots $.
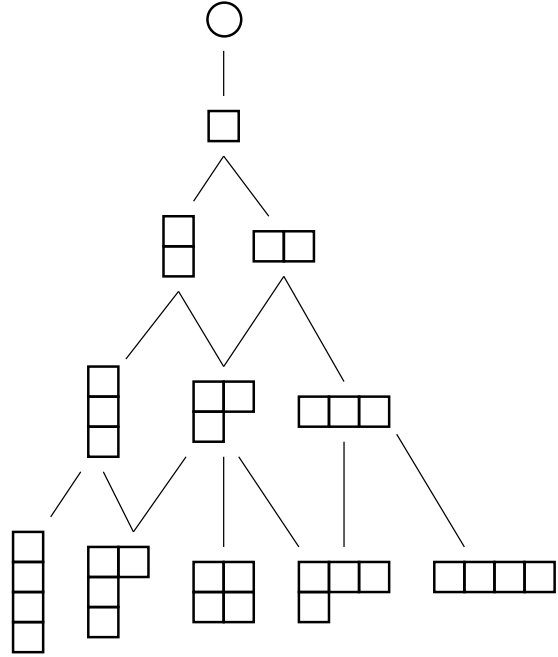
\begin{figure}[hbt]
\begin{center}
  \input{young100.pstex_t}
\caption{The Young lattice up to the $4^{th}$ level}  \label{fig:young}
\end{center}
\end{figure}

Let $\Gamma_{t}$ be the set of vertices on the $t^{th}$ level of the Young lattice, that is,
$\Gamma_{t}$ is the set of Young diagrams containing $t-2k \geq 0$ boxes, where
$k$ ranges over all of $\mathbb{Z}_{+}$.
Then $\mathscr{BW}_{t}$ is isomorphic to a direct sum of matrix algebras 
\cite[Thm. 3.5]{w2}:
$$\mathscr{BW}_{t} \cong \bigoplus_{\mu \in \Gamma_{t}} M_{b_{\mu}}(\mathbb{C}).$$
Ram and Wenzl defined matrix units for $\mathscr{BW}_{t}$  \cite{rw}.
We will write down these matrix units below.

To label the matrix units of $\mathscr{BW}_{t}$ we need to discuss the Bratteli diagram of
$\mathscr{BW}_{t}$, which is the following graph. 
The vertices of the Bratteli diagram of $\mathscr{BW}_{t}$ 
are divided into levels; for each $s=0, 1, \ldots, t$,
the elements of $\Gamma_{s}$ are the vertices on the $s^{th}$ level of
the Bratteli diagram of $\mathscr{BW}_{t}$.
The edges are as follows: a vertex
$\mu$ on the $s^{th}$ level is connected to a vertex $\lambda$ on the $(s+1)^{st}$ level if 
and only if $\mu$ and $\lambda$ differ by exactly one box.

We say that $R$ is a path of length $t$ in
the Bratteli diagram of $\mathscr{BW}_{t}$ 
if $R$ is a sequence of $t+1$ Young diagrams:
$R = ([0], [1], r_{2}, \ldots, r_{t})$ where $r_{s} \in \Gamma_{s}$ for each $s=0, 1, \ldots, t$
and where $r_{i}$ is connected to $r_{i+1}$ for each $0 \leq i \leq t-1$.
We say that $shp(R) = r_{t}$.
Let $\omega_{t}$ be the set of pairs $(R,S)$ of paths of length $t$ 
in the Bratteli diagram of $\mathscr{BW}_{t}$ satisfying $r_{t} = s_{t}$.

Ram and Wenzl defined a set of matrix units
$\{e_{ST} \in \mathscr{BW}_{t} | \ (S,T) \in \omega_{t} \}$ in \cite{rw}.
This set is a basis of $\mathscr{BW}_{t}$ and
the matrix units satisfy
$$e_{QR} e_{ST} =  \delta_{RS} e_{QT}.$$
We obtain the matrix units in 
$\mathscr{BW}_{t}(-q^{2n},q)/\mathscr{J}_{t}(-q^{2n},q)$ 
by taking a certain proper subset of the matrix units in $\mathscr{BW}_{t}$
and replacing the indeterminates $r$ and $q$ with the complex numbers
$-q^{2n}$ and $q$, respectively.

Let us fix some notation that we will use in the rest of this chapter and in 
Chapter \ref{chap2A:titlelabel}.  
Given a sequence $T$ of $t+1$ elements
$$T =(0, s_{1}, s_{2}, \ldots, s_{t-1}, s_{t}),$$ 
we fix $T'$ to be the following sequence of $t$ elements:
$$T' = (0, s_{1}, s_{2}, \ldots, s_{t-1}).$$
If $T$ is a path of length $t$, then $T'$ is the path of length $t-1$ obtained by removing the last
vertex and edge of $T$.

Before defining the matrix units of $\mathscr{BW}_{t}$ we
define some `pre-matrix units'.
Let $T$ be a path of length $t$ in the Bratteli diagram for $\mathscr{BW}_{t}$
such that $shp(T)$ has $t$ boxes.
We can identify $T$ with a standard tableau containing the numbers $1, 2, \ldots, t$
in a canonical way.
We do this by  placing the number $1$ in the top left hand box of  
$shp(T)$ and we then fill each box of $shp(T)$ with increasing numbers
according to the path $T$ \cite[Sec. 4.2]{tw}.

For each path $T$ of length $t$ in the Bratteli diagram for $\mathscr{BW}_{t}$,
 we define the number $d(T,i)$, for each $i=1, 2, \ldots, t-1$, by
\begin{equation}
\label{eq:brontethewhalecatslothdogmouse}
d(T,i) = c(i+1)-c(i)-r(i+1)+r(i),
\end{equation}
where $c(j)$ and $r(j)$ denote the column and row, respectively, 
of the box containing the number $j$ in the standard tableau corresponding to $T$.
For each $d \in \mathbb{Z} \backslash \{0\}$, we define
$$b_{d}(q) = \frac{q^{d}(1-q)}{1-q^{d}}.$$

Let $T$ be a path of length $t$ in the Bratteli diagram for $\mathscr{BW}_{t}$. 
Firstly fix $o_{[1]}=1 \in \mathscr{BW}_{t}$. 
Let $R$ be a path of length $t-1$ defined by $R = T'$ and inductively define
$$o_{T} = \prod_{S} \frac{o_{R} g_{t-1} o_{R} - b_{d(S,t-1)}(q^{2}) o_{R} }
         {b_{d(T,t-1)}(q^{2}) - b_{d(S,t-1)}(q^{2})} \in  \mathscr{BW}_{t},$$
where the product is over all paths $S$ of length $t$ where $shp(S)$ contains $t$ boxes such that
$S \neq T$ and $S' = R$. We write $o_{TT}=o_{T}$.

Let $M$ and $P$ be paths of length $t$ in the Bratteli diagram for $\mathscr{BW}_{t}$
where $(M,P) \in \omega_{t}$ and $shp(M) = shp(P)$ has exactly $t$ boxes
and $shp(M') = shp(P')$, then we define
$$o_{MP} = o_{M' P'} o_{PP}.$$

Let $M$ and $P$ be paths of length $t$ where $(M,P) \in \omega_{t}$ 
and $shp(M) = shp(P)$ has exactly $t$ boxes and
$shp(M') \neq shp(P')$, then the pre-matrix unit $o_{MP}$ is defined more intricately.
Choose paths $\overline{M}$ and $\overline{P}$ of length $t$ that satisfy
$shp(\overline{M}) = shp(M)$, $shp(\overline{P}) = shp(P)$ and the following three conditions:
\begin{itemize}
\item[(i)] $\overline{M}'' = \overline{P}''$,
\item[(ii)] $shp(\overline{M}') = shp(M')$, 
\item[(iii)] $shp(\overline{P}') = shp(P')$.
\end{itemize}
It may appear that these conditions cannot always be satisfied.  
However,  paths $\overline{M}$ and $\overline{P}$ can always be obtained
satisfying these conditions from the following construction \cite{rw}.
 
By considering $\overline{M}$ and $\overline{P}$ as standard tableaux,  
we obtain the desired paths $\overline{M}$ and $\overline{P}$ by ensuring the following
is true.
Firstly, fix $t$ to be in the same box in $\overline{M}$ (resp. $\overline{P}$) 
that $t$ is in $M$ (resp. $P$).
Then, fix $(t-1)$ to be in the same box in $\overline{M}$
(resp. $\overline{P}$) that $t$ is in $P$ (resp. $M$).
Lastly, for each $i = 1, 2, \ldots, t-2$, 
fix $i$ to be in the same box in $\overline{M}$ that it is in  $\overline{P}$.

We then define
$$o_{MP} = \frac{1-q^{2d}}{\sqrt{(1-q^{2(d+1)})(1-q^{2(d-1)})}} 
           o_{M'\overline{M}'} g_{t-1} o_{\overline{P}'P'}o_{PP},$$
where $d = d(\overline{M},t-1)$ is as given in (\ref{eq:brontethewhalecatslothdogmouse}).

This completes the definition of the `pre-matrix units'; now we
define the matrix units for $\mathscr{BW}_{t}$.  

Assume that the matrix units are known for $\mathscr{BW}_{t-1}$.
Let $M$ and $P$ be paths of length $t$
in the Bratteli diagram for $\mathscr{BW}_{t}$
where $shp(M) = \lambda = shp(P)$ 
and $\lambda$ contains strictly fewer than $t$ boxes, then we define
$$e_{MP} = \frac{Q_{\lambda}(r,q)}{\sqrt{Q_{\mu}(r,q)Q_{\widetilde{\mu}}(r,q)}}
                 e_{M'S} e_{t-1} e_{TP'},$$
where $S$ and $T$ are paths of length $t-1$ satisfying 
\begin{itemize}
\item[(i)]   $shp(S)=shp(M') = \mu$, and
\item[(ii)]  $shp(T)=shp(P') = \widetilde{\mu}$, and
\item[(iii)] $S'=T'$, and
\item[(iv)]  $shp(S')= \lambda = shp(T')$.
\end{itemize}
It may appear that these conditions cannot always be satisfied.  
However, there always exists a pair of paths $S$ and $T$ of length $t-1$
satisfying these conditions for the following reasons.
Firstly, by examining the relevant Bratteli diagrams, it is clear that 
there are no intertwiner matrix units in $\mathscr{BW}_{1}$ and $\mathscr{BW}_{2}$.  
Now for each $t \geq 3$, a shape $\lambda$ that has at most $t-2$ boxes and which labels a vertex
on the $t^{th}$ level of the Bratteli diagram for $\mathscr{BW}_{t}$ 
also labels a vertex on the $(t-2)^{nd}$ level of the Bratteli diagram.
Hence there always exists at least one path of length $t-2$ 
in the Bratteli diagram for $\mathscr{BW}_{t}$ ending at the vertex 
$shp(M)$ on the $(t-2)^{nd}$ level, as $shp(M)$ contains no more than $t-2$ boxes 
(this shows that (iii) and (iv) might be satisfied).

In the Bratteli diagram for $\mathscr{BW}_{t}$,
two vertices $\lambda$ and $\mu$ are connected by an edge 
only if their shapes differ by exactly one box.
Now the vertices $shp(M')$ and $shp(P')$ on the $(t-1)^{st}$ level
are connected to the vertex $shp(M)$ on the $t^{th}$ level by one edge each, and they
are also connected to the vertex $shp(M)$ on the $(t-2)^{nd}$ level by one edge each.
It follows, then, that by fixing $S$ and $T$ to be paths of length $t-1$ that coincide on the first
$t-2$ levels of the Bratteli diagram 
and that pass through the vertex $shp(M)$ on the $(t-2)^{nd}$ level, and also fixing
$shp(S) = shp(M')$ and $shp(T)=shp(P')$ (which is always possible), 
we obtain the desired paths $S$ and $T$.

Let $M$ and $P$ be paths of length $t$
in the Bratteli diagram for $\mathscr{BW}_{t}$, where $(M,P) \in \omega_{t}$, and where
$shp(M)$ contains $t$ boxes. 
Then we define
$$e_{MP} = (1-z_{t})o_{MP},$$
where $z_{t} = \sum_{P} e_{PP}$ with the summation going over all paths $P$ of length $t$ such that
$shp(P)$ contains fewer than $t$ boxes.

The following fact is important  \cite[Lem. 4.2]{w2}: 
let $M$ be a path of length $t$ in the Bratteli diagram for $\mathscr{BW}_{t}$
where $shp(M) = \lambda$,
then 
\begin{equation}
\label{eq:gregredinumber3}
\mathrm{tr}(e_{MM}) = Q_{\lambda}(r,q)/x^{t}, 
\end{equation}
where 
$\displaystyle{x = \frac{r-r^{-1}}{q-q^{-1}} + 1}$ and $Q_{\lambda}(r,q)$ is the polynomial
given in (\ref{eq:themagixQpolynomial}).

It is interesting to note that the quantum superdimension of the fundamental irreducible
$U_{q}(osp(1|2n))$-module $V$ is $\displaystyle{(-q^{2n}+q^{-2n})/(q-q^{-1}) + 1}$, which is just
the polynomial $x$ introduced in the preceding paragraph 
with the indeterminates $q$ and $r$ replaced with the complex numbers $q$ and
$-q^{2n}$, respectively.

\end{subsubsection}

\begin{subsubsection}{The algebra $\mathscr{BW}_{t}(r,q)$}

The algebra $\mathscr{BW}_{t}(r,q)$, with $r, q \in \mathbb{C}$, is equipped with a functional 
$\mathrm{tr}: \mathscr{BW}_{t}(r,q) \rightarrow \mathbb{C}$ which satisfies, 
amongst other relations \cite[Lem. 3.4 (d)]{w2},
\begin{equation}
\label{eq:tracefunctional}
\mathrm{tr}(a \chi b) = \mathrm{tr}(\chi) \mathrm{tr}(ab), 
\hspace{10mm} \forall a, b \in \mathscr{BW}_{t-1}(r,q), \hspace{5mm} \chi \in \{g_{t-1}, e_{t-1}\},
\end{equation}
where we regard each element of $\mathscr{BW}_{t-1}(r,q)$ as an element of $\mathscr{BW}_{t}(r,q)$
under the obvious inclusion.

Define the annihilator ideal $\mathscr{J}_{t}(r,q) \subset \mathscr{BW}_{t}(r,q)$ with respect to 
$\mathrm{tr}$ by
$$\mathscr{J}_{t}(r,q) = \left\{b \in \mathscr{BW}_{t}(r,q) | \ 
\mathrm{tr}(ab) = 0, \ \forall a \in \mathscr{BW}_{t}(r,q)\right\}.$$
If $q$ is not a root of unity and $r \neq \pm q^{k}$ for any $k \in \mathbb{Z}$,
$\mathscr{J}_{t}(r,q)=0$ and $\mathscr{BW}_{t}(r,q)$ is semisimple  \cite[Cor. 5.6]{w2}.
If $r = \pm q^{k}$ for some $k \in \mathbb{Z}$, then
$\mathscr{J}_{t}(\pm q^{k},q) \neq 0$ and
the quotient $\mathscr{BW}_{t}(\pm q^{k},q)/\mathscr{J}_{t}(\pm q^{k},q)$ is semisimple
\cite[Cor. 5.6]{w2}. Let us now fix $k=2n$ and $r=-q^{2n}$;
recall that the homomorphism $\Upsilon: g_{i} \mapsto -\check{R}_{i}$ yields a representation of
$\mathscr{BW}_{t}(-q^{2n},q)$ in ${\cal{C}}_{t}$.
The next task is to determine the structure of the quotient 
$\mathscr{BW}_{t}(-q^{2n},q)/\mathscr{J}_{t}(-q^{2n},q)$. We do this
in the following work.

We now introduce a subgraph $\Gamma(-q^{2n},q)$ of the Young lattice that we will use
to describe the structure of
$\mathscr{BW}_{t}(-q^{2n},q)/\mathscr{J}_{t}(-q^{2n},q)$.
We inductively obtain the vertices of $\Gamma(-q^{2n},q)$ as follows.
Firstly fix the Young diagram with no boxes to belong to $\Gamma(-q^{2n},q)$.
The inductive step is that if the Young diagram $\mu$ belongs to $\Gamma(-q^{2n},q)$, 
the Young diagram $\lambda$ then also belongs to $\Gamma(-q^{2n},q)$
if $\lambda$ differs from $\mu$ by exactly one box and if 
$Q_{\lambda}(-q^{2n},q) \neq 0$, where
$Q_{\lambda}(r,q)$ is given in (\ref{eq:themagixQpolynomial}).

We now give the polynomial  $Q_{\lambda}(r,q)$.
Given a Young diagram $\lambda$, let
$(i,j)$  denote the box in the $i^{th}$ row and the $j^{th}$ column of $\lambda$, 
and let $\lambda_{i}$ (resp. $\lambda_{j}'$) 
denote the number of boxes in the $i^{th}$ row (resp. $j^{th}$ column) of $\lambda$.
We introduce some notation: we may denote the Young diagram $\lambda$ by 
$\lambda = [\lambda_{1}, \lambda_{2}, \ldots, \lambda_{k}]$ where the $i^{th}$ row contains 
$\lambda_{i}$ boxes for each $i=1, 2, \ldots, k$, 
and the $l^{th}$ row contains no boxes for each $l > k$.
The polynomial  $Q_{\lambda}(r,q)$ is
\begin{eqnarray}
Q_{\lambda}(r,q) & = & 
\prod_{(j,j) \in \lambda} 
\frac{rq^{\lambda_{j} - \lambda_{j}'}-r^{-1}q^{-\lambda_{j} +\lambda_{j}'}+
q^{\lambda_{j} + \lambda_{j}'-2j+1}-q^{-\lambda_{j} -\lambda_{j}'+2j-1}}
{q^{h(j,j)}-q^{-h(j,j)}} \nonumber \\
& & \hspace{5mm} \times
\prod_{(i,j) \in \lambda, i \neq j} \frac{r q^{d(i,j)}-r^{-1}q^{-d(i,j)}}{q^{h(i,j)}-q^{-h(i,j)}},
\label{eq:themagixQpolynomial}
\end{eqnarray}
where the hooklength $h(i,j)$ is defined by
$h(i,j) = \lambda_{i} - i + \lambda_{j}' - j +1$, and where
$$d(i,j) = \left\{ \begin{array}{ll}
 \lambda_{i} + \lambda_{j}-i-j+1, & \mbox{if } i \leq j, \\
-\lambda_{i}' - \lambda_{j}' + i + j-1, & \mbox{if } i > j.
\end{array} \right.$$
More intuitively, 
the hooklength $h(i,j)$ is the number of boxes below the $(i,j)$ box in the $j^{th}$ column 
plus the number of boxes to the right of the $(i,j)$ box in the $i^{th}$ row, plus one.

Now $h(i,j) \geq 1$ for all $(i,j) \in \lambda$, so $Q_{\lambda}(-q^{2n},q)$ is well-defined for all
 $\lambda$. Also, for each $(j,j) \in \lambda$ we have
\begin{eqnarray*}
\lefteqn{ -q^{2n + \lambda_{j} - \lambda_{j}'}+q^{-2n-\lambda_{j} +\lambda_{j}'}+
q^{\lambda_{j} + \lambda_{j}'-2j+1}-q^{-\lambda_{j} -\lambda_{j}'+2j-1} } \\
& & \hspace{20mm} = (q^{-n+\lambda_{j}'-j+1/2}-q^{n-\lambda_{j}'+j-1/2})
      (q^{n+\lambda_{j}-j+1/2}+q^{-n-\lambda_{j}+j-1/2}),
\end{eqnarray*}
 so $Q_{\lambda}(-q^{2n},q) = 0$ if and only if one (or both)
   of the following conditions is satisfied:
\begin{itemize}
\item[(a)] $q^{4n + 2d(i,j)}=1$ for some $(i,j) \in \lambda$ where $i \neq j$, 
\item[(b)] $q^{2n-2\lambda_{j}'+2j-1} = 1$ or $q^{2n+2\lambda_{j}-2j+1}=-1$ for some $j$.
\end{itemize}
Now $q$ is non-zero and not a root of unity, so (b) is never satisfied for any $\lambda$, and
(a) is only satisfied if $d(i,j) = -2n$.
We now determine the circumstances for which $d(i,j) = -2n$. 
If $i > j$, we can see that
$\mathrm{min}(d(i,j)) = d(2,1) = -\lambda_{1}'-\lambda_{2}'+2$
from the constraints on the lengths of the columns of a Young diagram and
it follows that
$Q_{\lambda}(-q^{2n},q) = 0$ if $\lambda_{1}' + \lambda_{2}' = 2n+2$.
Let us call a Young diagram $\lambda$ {\emph{allowable}} if
$\lambda_{1}' + \lambda_{2}' \leq 2n+1$.

Across all the allowable Young diagrams, let us calculate $\mathrm{min}(d(i,j))$ where $i < j$.
If the first column of the allowable diagram $\lambda$ contains $2n+1$ boxes, 
ie $\lambda_{1}'=2n+1$, then all the other columns must contain no boxes
from the definition of an allowable diagram.
For such a $\lambda$, there does not exist any box $(i,j)$ in the $i^{th}$ row and the
$j^{th}$ column with $i < j$ and so there is nothing more to consider in this case. 
Now if the first column of $\lambda$ contains strictly fewer than $2n+1$ boxes, ie
$\lambda'_{1} \leq 2n$, then the following relations hold:
$i \leq 2n$, $\lambda_{i} - j \geq 0$ and $\lambda_{j} \geq 0$. 
Then $d(i,j)=\lambda_{i} + \lambda_{j}-i-j+1 \geq -2n+1$,
which means that $d(i,j) \neq -2n$ for all $i < j$.

It follows that 
$Q_{\lambda}(-q^{2n},q) = 0$ if $\lambda'_{1} + \lambda'_{2} = 2n+2$ and that
$Q_{\lambda}(-q^{2n},q) \neq 0$ for all allowable Young diagrams $\lambda$.
Consequently, the vertices of $\Gamma(-q^{2n},q)$ 
are all the allowable Young diagrams, 
that is, all the Young diagrams $\lambda$ satisfying $\lambda_{1}' + \lambda_{2}' \leq 2n+1$.

Now $\mathscr{J}_{t}(-q^{2n},q) \neq 0$ and $\mathscr{BW}_{t}(-q^{2n},q)$ is not semisimple.
However, the quotient $\mathscr{BW}_{t}(-q^{2n},q)/\mathscr{J}_{t}(-q^{2n},q)$ 
is semisimple:
$$\mathscr{BW}_{t}(-q^{2n},q)/\mathscr{J}_{t}(-q^{2n},q) \cong
 \bigoplus_{\lambda \in \Gamma(-q^{2n},q)_{t}} M_{b_{\lambda}}(\mathbb{C}),$$
where $\Gamma(-q^{2n},q)_{t}$ is the set of Young diagrams belonging to $\Gamma(-q^{2n},q)$
with $t-2k \geq 0$ boxes, where $k$ ranges over all of $\mathbb{Z}_{+}$  \cite[Cor. 5.6]{w2}.

We can obtain matrix units for $\mathscr{BW}_{t}(-q^{2n},q)/\mathscr{J}_{t}(-q^{2n},q)$
from the matrix units of  $\mathscr{BW}_{t}$. 
We replace the indeterminates $r$ and $q$ in some of the latter matrix units with the complex numbers
$-q^{2n}$ and $q$, respectively, to obtain
matrix units for $\mathscr{BW}_{t}(-q^{2n},q)/\mathscr{J}_{t}(-q^{2n},q)$.

To label the matrix units for $\mathscr{BW}_{t}(-q^{2n},q)/\mathscr{J}_{t}(-q^{2n},q)$,
we use the Bratteli diagram for $\mathscr{BW}_{t}(-q^{2n},q)/\mathscr{J}_{t}(-q^{2n},q)$,
which we define in the same way as we defined the Bratteli diagram for
$\mathscr{BW}_{t}$ but we replace $\Gamma_{s}$ with the  $\Gamma(-q^{2n},q)_{s}$ 
detailed in the next paragraph, for each
$s=0, 1, \ldots, t$.

Recall that the sets $\Gamma(-q^{2n},q)_{s}$ are given as follows.
The graph $\Gamma(-q^{2n},q)$ is a subgraph of the Young lattice and
the vertices of $\Gamma(-q^{2n},q)$  are all
the allowable Young diagrams, that is, all the Young diagrams $\lambda$ satisfying
$\lambda'_{1} + \lambda'_{2} \leq 2n+1$.
Then, for each $s=0, 1, \ldots, t$,  $\Gamma(-q^{2n},q)_{s}$ 
is the set of Young diagrams belonging
to $\Gamma(-q^{2n},q)$ that contain exactly $s-2k \geq 0$ boxes, where $k$ ranges over all of 
$\mathbb{Z}_{+}$. 

We say that $T=(0, s_{1}, s_{2}, \ldots, s_{t})$
is a path of length $t$ in the Bratteli diagram for
$\mathscr{BW}_{t}(-q^{2n},q)/\mathscr{J}_{t}(-q^{2n},q)$ 
 if $s_{i} \in \Gamma(-q^{2n},q)_{i}$ for each $i$ and if 
 $s_{j}$ is joined to $s_{j+1}$ 
for each $j=0, \ldots, t-1$.
 
Let $\omega(-q^{2n},q)_{t}$ be the set of pairs $(R,S)$ 
of paths of length $t$ in the Bratteli diagram for
$\mathscr{BW}_{t}(-q^{2n},q)/\mathscr{J}_{t}(-q^{2n},q)$ 
where $r_{t} = s_{t}$, that is $shp(R)=shp(S)$.
The matrix units 
$$\{ e_{RS} \in \mathscr{BW}_{t} | \ (R,S) \in \omega(-q^{2n},q)_{t} \}$$ 
are all well-defined and non-zero if the indeterminates
$r$ and $q$ are replaced with the complex numbers $-q^{2n}$ and $q$, respectively.
Henceforth we write $e_{RS}$ to mean the matrix unit 
$e_{RS}(-q^{2n},q) \in \mathscr{BW}_{t}(-q^{2n},q)/\mathscr{J}_{t}(-q^{2n},q)$.
It is very important to note that $\mathrm{tr}(e_{SS}) \neq 0$ 
for all $(S,S) \in \omega(-q^{2n},q)_{t}$ and that 
$e_{RS} \notin \mathscr{J}_{t}(-q^{2n},q)$ for all  $(R,S) \in \omega(-q^{2n},q)_{t}$.

\end{subsubsection}

\begin{subsubsection}{Matrix units in $\mathscr{BW}_{t}(-q^{2n},q)/\mathscr{J}_{t}(-q^{2n},q)$
and ${\cal{C}}_{t}$}

We now relate the idempotent matrix units in 
$\mathscr{BW}_{t}(-q^{2n},q)/\mathscr{J}_{t}(-q^{2n},q)$
to the projectors in $\mathcal{C}_{t}$ we defined at the start of this section.
Let $\widetilde{\mathscr{BW}}_{t}(-q^{2n},q)$ be the semisimple subalgebra of
$\mathscr{BW}_{t}(-q^{2n},q)$ spanned by 
the matrix units in $\mathscr{BW}_{t}(-q^{2n},q)/\mathscr{J}_{t}(-q^{2n},q)$, ie
$\{ e_{RS} | \ (R,S) \in \omega(-q^{2n},q)_{t} \}$.

Firstly, we will show that
$\mathscr{BW}_{t}(-q^{2n},q) = 
\widetilde{\mathscr{BW}}_{t}(-q^{2n},q) \oplus \mathscr{J}_{t}(-q^{2n},q)$.
Any $f \in \widetilde{\mathscr{BW}}_{t}(-q^{2n},q)$ can be written as
$$f = \sum_{(S,T) \in \omega(-q^{2n},q)_{t}} f_{ST} e_{ST}, \hspace{10mm} f_{ST} \in \mathbb{C},$$
where $f_{ST} \neq 0$ for at least one pair $(S,T)$ of paths.  
Fix $(A,B)$ to be such a pair, then
$$\mathrm{tr}(e_{BA}f) = \mathrm{tr}(f_{AB} e_{BA} e_{AB}) = f_{AB}\mathrm{tr}(e_{BB}) \neq 0,$$
as $\mathrm{tr}(e_{BB}) \neq 0$.  Thus any non-zero $f$ belonging to 
$\widetilde{\mathscr{BW}}_{t}(-q^{2n},q)$ does not belong to 
$\mathscr{J}_{t}(-q^{2n},q)$, yielding
$$\mathscr{BW}_{t}(-q^{2n},q) = 
\widetilde{\mathscr{BW}}_{t}(-q^{2n},q) \oplus \mathscr{J}_{t}(-q^{2n},q).$$
Then we can write $a = \widetilde{a} + a_{j}$ for each $a \in \mathscr{BW}_{t}(-q^{2n},q)$, 
  where $\widetilde{a} \in \widetilde{\mathscr{BW}}_{t}(-q^{2n},q)$ and 
  $a_{j} \in \mathscr{J}_{t}(-q^{2n},q)$.

Now define 
$$\mathscr{P}_{t} = \sum_{(S,S) \in \omega(-q^{2n},q)_{t}} e_{SS} \in 
\widetilde{\mathscr{BW}}_{t}(-q^{2n},q),$$
then $\mathscr{P}_{t} a \mathscr{P}_{t} = \widetilde{a}$, which can be seen by regarding
$\mathscr{BW}_{t}(-q^{2n},q)$ as a matrix algebra.

\end{subsubsection}

Let us now turn our attention to ${\cal{C}}_{t}$.
Define $J_{t} \subset {\cal{C}}_{t}$ to be the annihilator ideal of 
${\cal{C}}_{t}$ with respect to the quantum supertrace:
$$J_{t} = \{ b \in {\cal{C}}_{t} | \ str_{q}(ab)=0, \ \forall a \in {\cal{C}}_{t} \}.$$
Now define a map $\psi: {\cal{C}}_{t} \rightarrow \mathbb{C}$ by
$$\psi(X) = str_{q}(X)/ \big( sdim_{q}(V) \big)^{t},$$
then $\psi(X) = 0$ if and only if $str_{q}(X) = 0$, and furthermore,
\begin{equation}
\label{eq:thegreaans(1)}
\psi \big(\Upsilon(a) \big) = \mathrm{tr}(a), \hspace{10mm} \forall a \in \mathscr{BW}_{t}(-q^{2n},q),
\end{equation}
from Lemma \ref{lemma:equalityoftraces}.
Thus we can regard $J_{t}$ as the annihilator ideal of ${\cal{C}}_{t}$ with respect to $\psi$.

Now we will use Eq. (\ref{eq:thegreaans(1)}) to show that
\begin{equation}
\label{eq:shonauy(99)}
\Upsilon\big(\mathscr{J}_{t}(-q^{2n},q)\big) = J_{t}.
\end{equation}
We firstly show that $\Upsilon\big(\mathscr{J}_{t}(-q^{2n},q)\big) \subseteq J_{t}$.
Let $b$ be an arbitrary element of $\mathscr{J}_{t}(-q^{2n},q)$, then $\mathrm{tr}(ab)=0$ for all 
$a \in \mathscr{BW}_{t}(-q^{2n},q)$, and the surjectivity of 
$\Upsilon$, in addition to the fact that $\psi\big(\Upsilon(ab)\big) = \mathrm{tr}(ab)$, means that
$\Upsilon(a) \in J_{t}$.

Now let $B$ be an arbitrary element of $J_{t}$, 
then there is some $b$ belonging to $\mathscr{BW}_{t}(-q^{2n},q)$
satisfying  $B = \Upsilon(b)$, and furthermore, $b \in \mathscr{J}_{t}(-q^{2n},q)$ as
$\mathrm{tr}(ab) = \psi\big(\Upsilon(a)\Upsilon(b) \big) = 0$ for all 
$a \in \mathscr{BW}_{t}(-q^{2n},q)$. Then
$\Upsilon\big(\mathscr{J}_{t}(-q^{2n},q)\big) \supseteq J_{t}$, proving
Eq. (\ref{eq:shonauy(99)}).

The surjectivity of $\Upsilon$ implies that 
$${\cal{C}}_{t} = \Upsilon\big(\widetilde{\mathscr{BW}}_{t}(-q^{2n},q) \big) + J_{t},$$
and we will show that this sum is direct.
To see this, assume that there exists some non-zero element $F$ of ${\cal{C}}_{t}$ belonging to
$\Upsilon\big(\widetilde{\mathscr{BW}}_{t}(-q^{2n},q) \big)$ and also to $J_{t}$, then
$str_{q}(XF) = 0$ for all $X \in {\cal{C}}_{t}$.  
However, $F$ is the image of a linear combination of matrix units: 
$$F = \sum_{(S,T) \in \omega(-q^{2n},q)_{t}}f_{ST}\Upsilon(e_{ST}),
\hspace{10mm} f_{ST} \in \mathbb{C},$$ 
where $f_{ST} \neq 0$  for at least one pair $(S,T)$.  
Assume that $(A,B)$ is such a pair, then by similar reasoning as
previously, $str_{q}(\Upsilon(e_{BA})F) \neq 0$ which contradicts the assumption that $F \in J_{t}$.
Thus
$\Upsilon\big(\widetilde{\mathscr{BW}}_{t}(-q^{2n},q) \big) \cap J_{t} = 0$,
and 
\begin{equation} 
\label{eq:milindaandmilinda}
{\cal{C}}_{t} = \Upsilon\big(\widetilde{\mathscr{BW}}_{t}(-q^{2n},q) \big) \oplus J_{t}.
\end{equation}

It is clear that the image 
of each matrix unit $e_{ST} \in \widetilde{\mathscr{BW}}_{t}(-q^{2n},q)$ in ${\cal{C}}_{t}$
under the map $\Upsilon$ is again a matrix unit.
Each matrix unit $e_{SS} \in \widetilde{\mathscr{BW}}_{t}(-q^{2n},q)$ is an idempotent, thus
each $\Upsilon(e_{SS})$ is an idempotent that is also $U_{q}(\mathfrak{g})$-linear.
Now $\big(\Upsilon(e_{SS})\big) V^{\otimes t} \neq 0$ as 
$str_{q}\big(\Upsilon(e_{SS})\big) = (sdim_{q}(V))^{t}\mathrm{tr}(e_{SS}) \neq 0$, and
as $V^{\otimes t}$ is completely reducible, $\Upsilon(e_{SS})$ projects down from $V^{\otimes t}$
onto a direct sum of irreducible $U_{q}(\mathfrak{g})$-submodules of $V^{\otimes t}$.
The matrix units $\{e_{SS} | \ (S,S) \in \omega(-q^{2n},q)_{t}\}$ are all orthogonal, 
thus all the $\Upsilon(e_{SS})$ are orthogonal. 

Let $S$ be a path of length $t$ in the Bratteli diagram for
$\widetilde{\mathscr{BW}}_{t}(-q^{2n},q)$; let $\lambda = shp(S)$.
If $\lambda$, as a Young diagram, contains no more than $n$ rows of boxes, then 
we can interpret $\Upsilon(e_{SS})$ as the projection from $V^{\otimes t}$ onto an irreducible 
$U_{q}(\mathfrak{g})$-submodule $V_{\lambda} \subseteq V^{\otimes t}$, where
we use the Young diagram $\lambda$ to label the integral dominant highest weight of
$V_{\lambda}$ as discussed in the third paragraph of 
Subsection \ref{subsect:matrixbirmanwenzlstuff}.
If $\lambda$, as a Young diagram, has more than $n$ rows of boxes (that is, $\lambda'_{1} > n$), 
 then we must consider $\Upsilon(e_{SS})$ more carefully.
In this case, $\lambda$ does not have an immediate interpretation as the 
highest weight of a finite dimensional 
irreducible representation of $U_{q}(\mathfrak{g})$.
However, there is a completely standard way of dealing with this problem.
Each such $\lambda$ satisfies $\lambda_{1}' + \lambda_{2}' \leq 2n+1$ and
we can regard $\lambda$ as labelling an irreducible representation of the Lie group
$SO(2n+1)$ as follows.  
Let $\widetilde{\lambda}$ be the following diagram:
fix $\widetilde{\lambda}_{1}' = 2n+1-\lambda_{1}'$ and $\widetilde{\lambda}_{j}' = \lambda_{j}'$
for all $j \geq 2$, then $\widetilde{\lambda}$ is a Young diagram (see the next paragraph), 
and the characters associated 
with the $SO(2n+1)$ representations labelled by
$\lambda$ and $\widetilde{\lambda}$ are the same \cite[Sec. 2]{ck}.  
Furthermore, the $OSp(1|2n)$ supercharacters of the representations labelled by $\lambda$ and
$\widetilde{\lambda}$ are the same up to a factor of $\pm 1$ \cite{farmer,ck}.

If $R$ and $S$ are paths of length $t$ in the Bratteli diagram for
$\widetilde{\mathscr{BW}}_{t}(-q^{2n},q)$ satisfying
$shp(R)=\lambda$ and $shp(S)=\widetilde{\lambda}$, then
$\Upsilon(e_{RR})$ and $\Upsilon(e_{SS})$ project down from $V^{\otimes t}$ onto 
isomorphic irreducible $U_{q}(\mathfrak{g})$-submodules of $V^{\otimes t}$.

We now show that the $\widetilde{\lambda}$ mentioned above 
is in fact a Young diagram.
Write $\lambda_{1}' = \lambda_{2}' + k$ where $k \geq 1$, then
$\widetilde{\lambda}_{1}' = 2n+1 - (\lambda_{2}' + k)$.  
Now $\widetilde{\lambda}$ is a Young diagram if
$\widetilde{\lambda}_{1}' \geq \widetilde{\lambda}_{2}'$ 
(which is just $2n+1 - k \geq 2\lambda_{2}'$) and this is  true as 
$\lambda_{1}' + \lambda_{2}' = 2\lambda_{2}' + k \leq 2n+1$.

Let $\lambda$ be a Young diagram with more than $n$ rows of boxes satisfying
$\lambda'_{1} + \lambda_{2}' \leq 2n+1$ and let 
$\widetilde{\lambda}$ be the Young diagram given by  
$\widetilde{\lambda}_{1}' = 2n+1-\lambda_{1}'$ and $\widetilde{\lambda}_{j}' = \lambda_{j}'$
for all $j \geq 2$.
We now show that there do not exist idempotent matrix units $e_{RR}$ and $e_{SS}$ in 
$\widetilde{\mathscr{BW}}_{t}(-q^{2n},q)$ where $shp(R) = \lambda$ and $shp(S)=\widetilde{\lambda}$.
We show this important result using an easy even/odd number argument.
If the number of boxes in $\lambda$ is even (resp. odd), then the number of boxes in
$\widetilde{\lambda}$ is odd (resp. even), as 
$$\widetilde{\lambda}_{1}' \bmod{2} = (2n+1-\lambda_{1}') \bmod{2} = (\lambda_{1}'+1) \bmod{2}
\hspace{5mm} \mbox{and} \hspace{5mm} \widetilde{\lambda}_{j}' = \lambda_{j}',  \ \ j \geq 2.$$
Now let $r$ be an even (resp. odd) number satisfying $0 \leq r \leq t$,
then the vertices on the $r^{th}$ level of the Bratteli diagram for
$\widetilde{\mathscr{BW}}_{t}(-q^{2n},q)$
are all the Young diagrams in  $\Gamma(-q^{2n},q)$ with  $k$ boxes where $k \leq r$ is
an even (resp. odd) number. 
Let $|\lambda|$ denote the number of boxes in the Young diagram $\lambda$, then
 $|\lambda| \bmod{2} = ( |\widetilde{\lambda}|+1 ) \bmod{2}$, and consequently
 it is not possible that
$\lambda$ and $\widetilde{\lambda}$  are vertices on the same level of 
the Bratteli diagram for $\widetilde{\mathscr{BW}}_{t}(-q^{2n},q)$.

As it is not possible that both  
$\lambda$ and $\widetilde{\lambda}$ are vertices on the $t^{th}$ level of 
Bratteli diagram for $\widetilde{\mathscr{BW}}_{t}(-q^{2n},q)$,
at most only one of $e_{RR}$ and $e_{SS}$ exists in 
$\widetilde{\mathscr{BW}}_{t}(-q^{2n},q)$
where $shp(R) = \lambda$ and $shp(S) = \widetilde{\lambda}$, and thus
no more than one of $\Upsilon(e_{RR})$ and $\Upsilon(e_{SS})$ exists 
for any $t$.

Let $S$ be a path of length $t$ in the Bratteli diagram for 
$\widetilde{\mathscr{BW}}_{t}(-q^{2n},q)$ and let $R$ be a path of length $t$ in the Bratteli
diagram for $V^{\otimes t}$.
We can directly compare the orthogonal idempotents $\Upsilon(e_{SS})$ with 
the orthogonal projectors $E_{RR} \in {\cal{C}}_{t}$ by examining the paths $S$ and $R$.
It can be seen from the definitions of the idempotents $e_{SS}$ and the map $\Upsilon$
that the idempotents $\Upsilon(e_{SS})$ act on $V^{\otimes t}$ from the left-most tensor powers as
follows.
If $R$ is a path of length $2$, then
$\Upsilon(e_{RR})$ as an element of  ${\cal{C}}_{t}$ is
$\Upsilon(e_{RR}) \otimes \mathrm{id}^{\otimes (t-2)}$.
If $S$ is a path of length $i$ where $2 \leq i \leq t$, then 
$\Upsilon(e_{SS})$ as an element of  ${\cal{C}}_{t}$ is
$\Upsilon(e_{SS}) \otimes \mathrm{id}^{\otimes (t-i)}$.

Now let $R$ be a path of length $t$ in the Bratteli diagram for 
$\widetilde{\mathscr{BW}}_{t}(-q^{2n},q)$.
Let $S$ be the same path of length $t$ as $R$ except that if a Young diagram 
$\lambda$ on the path $R$ has more
than $n$ rows of boxes, we take the transformed Young diagram 
$\widetilde{\lambda}$ to be in $S$ instead of $\lambda$.
Then $S$ is a path of length $t$ in the Bratteli diagram for $V^{\otimes t}$.
Recall that
the integral dominant highest weights of the irreducible $U_{q}(\mathfrak{g})$-submodules
of $V^{\otimes k}$, where $k \leq t$,
are the vertices on the $k^{th}$ level of the Bratteli diagram for $V^{\otimes t}$.
Then $\Upsilon(e_{RR}) \in {\cal{C}}_{t}$ and $E_{SS} \in {\cal{C}}_{t}$ project onto the same 
irreducible $U_{q}(\mathfrak{g})$-submodule of $V^{\otimes t}$, and we also have
 $Q_{shp(R)}(-q^{2n},q) = sdim_{q}(V_{shp(S)})$
from (\ref{eq:gregredinumber3}) and (\ref{eq:thegreaans(1)}).

Let $E_{SS} \in {\cal{C}}_{t}$ be a projector with the property that 
$\left(E_{SS}V^{\otimes t}\right) \cap \left(\Upsilon(e_{RR})V^{\otimes t}\right) = 0$ 
for all idempotent 
matrix units $e_{RR} \in \widetilde{\mathscr{BW}}_{t}(-q^{2n},q)$. We will show that no such
projector exists. 
Suppose that such a projector does exist, that $E_{SS}$ is a projector that is orthogonal to
$\Upsilon(e_{RR})$ for each idempotent matrix unit   
$e_{RR} \in \widetilde{\mathscr{BW}}_{t}(-q^{2n},q)$.
Then $E_{SS}$ is orthogonal to $\Upsilon(e_{RT})$ for each 
matrix unit 
$e_{RT} \in \widetilde{\mathscr{BW}}_{t}(-q^{2n},q)$ where $R \neq T$ and
$(R,T) \in \omega(-q^{2n},q)_{t}$, as
$$E_{SS} \Upsilon(e_{RT}) = E_{SS} \Upsilon(e_{RR}e_{RT}e_{TT}) = 0
 = \Upsilon(e_{RR}e_{RT}e_{TT})E_{SS} = \Upsilon(e_{RT})E_{SS},$$ 
 as $\Upsilon$ is a homomorphism.
 
As $E_{SS}$ is orthogonal to $\Upsilon(e_{RT})$ for each  matrix unit 
$e_{RT} \in \widetilde{\mathscr{BW}}_{t}(-q^{2n},q)$, $(R,T) \in \omega(-q^{2n},q)_{t}$,
it is true that $E_{SS} \in J_{t}$.
To see this, assume the contrary, then
\begin{equation}
\label{eq:john_singleton(19)}
E_{SS} = \widetilde{E}_{SS} + E_{j},
\end{equation}
 where
$\widetilde{E}_{SS} \in \Upsilon\big(\widetilde{\mathscr{BW}}_{t}(-q^{2n},q) \big)$
and $E_{j} \in J_{t}$ 
(we can write any element of ${\cal{C}}_{t}$ as a sum of elements of
$\Upsilon\big(\widetilde{\mathscr{BW}}_{t}(-q^{2n},q)$ and $J_{t}$ from (\ref{eq:milindaandmilinda})).
Then
\begin{equation}
\label{eq:illbetheretosupportyou}
\widetilde{E}_{SS} = \sum_{(R,T) \in \omega(-q^{2n},q)_{t}} c_{RT} \Upsilon(e_{RT}), 
\hspace{5mm} c_{RT} \in \mathbb{C},
\end{equation}
where $c_{RT} \neq 0$ for at least one pair $(R,T) \in \omega(-q^{2n},q)_{t}$.
Assume that $(A, B)  \in \omega(-q^{2n},q)_{t}$ is such a pair so that $c_{AB} \neq 0$, then
\begin{equation}
\label{eq:constructivist}
str_{q} \left(\Upsilon(e_{BA})\widetilde{E}_{SS} \right) = 
str_{q} \left( \sum_{T} c_{AT} \Upsilon(e_{BT}) \right) 
 =  c_{AB} str_{q} \left( \Upsilon(e_{BB}) \right) \neq  0,
\end{equation}
as $c_{AB} \neq 0$ and $str_{q}(\Upsilon(e_{BB})) \neq 0$.

We will show that (\ref{eq:constructivist}) is not true. Recall that  
$E_{SS}$ satisfies
\begin{equation}
\label{eq:caseychamberscat}
\Upsilon(e_{BA}) E_{SS} = \Upsilon(e_{BA})(\widetilde{E}_{SS} + E_{j} ) = 0,
\end{equation} 
and note that
$\big(\Upsilon(e_{BA})\widetilde{E}_{SS}\big) 
\in \Upsilon\big(\widetilde{\mathscr{BW}}_{t}(-q^{2n},q)$
and $\big(\Upsilon(e_{BA})E_{j}\big) \in J_{t}$. This last fact means that 
$\Upsilon(e_{BA})\widetilde{E}_{SS} \neq -\Upsilon(e_{BA})E_{j}$ 
if $ \Upsilon(e_{BA})\widetilde{E}_{SS} \neq 0$ and
$\Upsilon(e_{BA})E_{j} \neq 0$, and by
re-examining (\ref{eq:caseychamberscat}) it is then clear that
$\Upsilon(e_{BA})\widetilde{E}_{SS} = \Upsilon(e_{BA})E_{j} = 0$.

An implication of the result $\Upsilon(e_{BA})\widetilde{E}_{SS} = 0$
is that $str_{q} \left( \Upsilon(e_{BA})\widetilde{E}_{SS} \right) = 0$, however, this 
contradicts Eq. (\ref{eq:constructivist}).
This then implies that the assumption in (\ref{eq:illbetheretosupportyou}) that
$c_{RT} \neq 0$ for at least one pair $(R,T) \in \omega(-q^{2n},q)_{t}$ is false, thus
we have $\widetilde{E}_{SS} = 0$.

It then follows that it must be true that 
$E_{SS} = E_{j} \in J_{t}$ from (\ref{eq:john_singleton(19)}). 
However, this is not true as
$str_{q}(E_{SS}) = sdim_{q}(E_{SS}V^{\otimes t}) \neq 0$.
Thus, our original assumption that there exists 
a projector $E_{SS} \in {\cal{C}}_{t}$ with the property that 
$\left(E_{SS}V^{\otimes t}\right) \cap \left(\Upsilon(e_{RR})V^{\otimes t}\right) = 0$ 
for all idempotent 
matrix units $e_{RR} \in \widetilde{\mathscr{BW}}_{t}(-q^{2n},q)$ is not true.

\end{subsection}

\begin{subsection}{Matrix units in $\mathcal{C}_{t}$}
\label{subsection:genericcalmatrixunitsinCt}

We have not yet proved that ${\cal{L}}_{t} = {\cal{C}}_{t}$. 
We will complete the proof in this subsection by defining a complete set of 
intertwiners in ${\cal{C}}_{t}$,
which we will obtain by applying the map $\Upsilon$ to the intertwiner matrix units in
$\widetilde{\mathscr{BW}}_{t}(-q^{2n},q)$.

There is a one-to-one map between paths in the Bratteli diagram for 
$\widetilde{\mathscr{BW}}_{t}(-q^{2n},q)$
and paths in the Bratteli diagram for $V^{\otimes t}$.
Recall that each Young diagram on the $k^{th}$ level of the Bratteli diagram for 
$\widetilde{\mathscr{BW}}_{t}(-q^{2n},q)$ contains an even (resp. odd) number of boxes
if $k$ is an even (resp. odd) number.
Each vertex $\lambda$ on the $k^{th}$ level of the Bratteli diagram for 
$\widetilde{\mathscr{BW}}_{t}(-q^{2n},q)$ appears on the $k^{th}$ level of the Bratteli diagram for
$V^{\otimes t}$ unless $\lambda$ has more than $n$ rows of boxes, in which case the Young diagram
$\widetilde{\lambda}$ appears instead, where $\widetilde{\lambda}$ is the Young diagram defined in the
second paragraph after Eq. (\ref{eq:milindaandmilinda}).

Given a path $\widetilde{T}$ of length $t$ in the Bratteli diagram for $V^{\otimes t}$, 
we can write down the corresponding path $T$ of length $t$ in the Bratteli diagram for 
$\widetilde{\mathscr{BW}}_{t}(-q^{2n},q)$ as follows. 
Write $\widetilde{T} = (0, s_{1}, s_{2}, \ldots, s_{t})$ where $s_{i}$ is a Young diagram on the
$i^{th}$ level of the Bratteli diagram for $V^{\otimes t}$. 
If $i$ is an even (resp. odd) number and $s_{i}$ contains an even (resp. odd) number of boxes, then 
$s_{i}$ is also a vertex on the $i^{th}$ level of the Bratteli diagram for 
$\widetilde{\mathscr{BW}}_{t}(-q^{2n},q)$.
If, however, $i$ is an even (resp. odd) number and $s_{i}$ contains an odd (resp. even) number of
boxes, then $s_{i}=\widetilde{\lambda}$ is the vertex that is obtained by taking a vertex $\lambda$ on
the $i^{th}$ level of the Bratteli diagram for  $\widetilde{\mathscr{BW}}_{t}(-q^{2n},q)$ and
defining $\widetilde{\lambda}$ by
$\widetilde{\lambda}'_{1} = 2n+1 - \lambda'_{1}$ 
and $\widetilde{\lambda}'_{j} = \lambda_{j}'$ for $j \geq 2$.
Using this, we can define a path $T$ of length $t$ in the Bratteli diagram for
$\widetilde{\mathscr{BW}}_{t}(-q^{2n},q)$ corresponding to the path
$\widetilde{T}$ of length $t$ in the Bratteli diagram for $V^{\otimes t}$.

It is easy to define the intertwiners in ${\cal{C}}_{t}$ between the isomorphic irreducible
$U_{q}(\mathfrak{g})$-submodules of $V^{\otimes t}$ obtained by 
using the projectors $E_{RR} \in {\cal{C}}_{t}$, where
$R$ is a path of length $t$ in the Bratteli diagram for $V^{\otimes t}$.
All we need do is to check that the images in ${\cal{C}}_{t}$ of the intertwiner matrix units
are well-defined and non-zero.

We construct the intertwiners in ${\cal{C}}_{t}$ recursively.  
To do this, assume that all the matrix units in ${\cal{C}}_{t-1}$ have already been defined, 
and that they are non-zero. 
Note that the decomposition of $V \otimes V$ into irreducible $U_{q}(\mathfrak{g})$-submodules 
is multiplicity free, so no intertwiners exist in ${\cal{C}}_{2}$.

In the remainder of the subsection, let $\widetilde{M}$ and $\widetilde{P}$ 
be a pair of paths of length $t$ in the Bratteli diagram for $V^{\otimes t}$
where $shp(\widetilde{M}) = shp(\widetilde{P})$ and $\widetilde{M} \neq \widetilde{P}$. 
Let $M$ and $P$ be the corresponding paths in the Bratteli diagram for
$\widetilde{\mathscr{BW}}_{t}(-q^{2n},q)$.
The intertwiner $E_{\widetilde{M}\widetilde{P}} \in {\cal{C}}_{t}$ is precisely
$E_{\widetilde{M}\widetilde{P}} = \Upsilon(e_{MP})$.

Let us firstly deal with the situation that $shp(M)=shp(P)=\lambda$ where
$\lambda$ contains strictly fewer than $t$ boxes. 
Referring back to Subsection \ref{subsect:matrixbirmanwenzlstuff} we see that
$$E_{\widetilde{M}\widetilde{P}} =\Upsilon(e_{MP}) 
= \frac{Q_{\lambda}(-q^{2n},q)}{\sqrt{Q_{\mu}(-q^{2n},q)Q_{\widetilde{\mu}}(-q^{2n},q)}}
                 E_{\widetilde{M}'\widetilde{S}} \Upsilon(e_{t-1}) E_{\widetilde{T}\widetilde{P}'}, 
		 \hspace{5mm} e_{t-1} \in \mathscr{BW}_{t}(-q^{2n},q),$$
where $S$ and $T$ are paths of length $t-1$ in the Bratteli diagram for 
$\widetilde{\mathscr{BW}}_{t}(-q^{2n},q)$ such that 
\begin{itemize}
\item[(i)]   $shp(S)=shp(M') = \mu$, and
\item[(ii)]  $shp(T)=shp(P') = \widetilde{\mu}$, and
\item[(iii)] $S'=T'$, and
\item[(iv)]  $shp(S')= \lambda = shp(T')$.
\end{itemize}
(Recall that such paths always exist.) Note that $Q_{\mu}(-q^{2n},q) \neq 0$ and
$Q_{\widetilde{\mu}}(-q^{2n},q) \neq 0$; if it were true that
$Q_{\mu}(-q^{2n},q) = 0$ then $\mu$ would not be a vertex
in the Bratteli diagram for $\widetilde{\mathscr{BW}}_{t}(-q^{2n},q)$.
Similar remarks hold for $\widetilde{\mu}$ and $\lambda$. 
Thus $E_{\widetilde{M}\widetilde{P}}$ is well-defined.
(We will later show that $E_{\widetilde{M}\widetilde{P}}$ is non-zero.)

Now let us deal with the situation that $shp(M)=shp(P)=\lambda$ where
$\lambda$ contains exactly $t$ boxes and $shp(M') = shp(P')$.
Referring back to Subsection \ref{subsect:matrixbirmanwenzlstuff}, we see that
$$E_{\widetilde{M}\widetilde{P}} =\Upsilon(e_{MP})
 = \Upsilon((1-z_{t}) o_{MP}),$$
 where $o_{MP} = o_{M'P'}o_{PP}$ and
 $z_{t} = \sum_{S}e_{SS}$
 with the summation going over all paths $S$ of length $t$ such that $shp(S)$ contains fewer than $t$
 boxes. It is not difficult to see that each such $\Upsilon(e_{SS})$ is some projection
 $E_{\widetilde{S}\widetilde{S}} \in {\cal{C}}_{t}$.

Now let us deal with the situation that $shp(M)=shp(P)=\lambda$ where
$\lambda$ contains exactly $t$ boxes and $shp(M') \neq shp(P')$.
Choose paths $\overline{M}$ and $\overline{P}$ of length $t$ such that 
$shp(\overline{M}) = shp(M)$ and $shp(\overline{P}) = shp(P)$ and
\begin{itemize}
\item[(i)] $\overline{M}'' = \overline{P}''$, and
\item[(ii)] $shp(\overline{M}') = shp(M')$, and
\item[(iii)] $shp(\overline{P}') = shp(P')$.
\end{itemize}
Such paths can always be chosen. Then
$$E_{\widetilde{M}\widetilde{P}} =\Upsilon(e_{MP})
 = \Upsilon((1-z_{t}) o_{MP}),$$
 where 
 \begin{equation}
 \label{eq:thaiequation}
 o_{MP} = \frac{1-q^{2d}}{\sqrt{(1-q^{2(d+1)})(1-q^{2(d-1)})}} 
           o_{M'\overline{M}'} g_{t-1} o_{\overline{P}'P'}o_{PP},
	   \hspace{5mm} g_{t-1} \in \mathscr{BW}_{t}(-q^{2n},q),
\end{equation}
where $d = d(\overline{M},t-1)$ is the integer defined by (\ref{eq:brontethewhalecatslothdogmouse}).
The integer $|d(\overline{M},i)|+1$ is the number of boxes in the hook going through the boxes 
containing the numbers $i$ and $(i+1)$ \cite{tw}.

We now prove that the coefficient on the right hand side of (\ref{eq:thaiequation}) 
is well-defined and non-zero.
It is not difficult to see that the coefficient is well-defined if $|d| \neq 1$, 
and we now show that this is always true.  
As $|d| + 1$ is the length of the hook going through the boxes 
containing the numbers $(t-1)$ and $t$, it is always true that $|d| + 1 \geq 2$ as 
each such hook contains at least two boxes.
Now the only situation in which it could possibly be true that $|d| = 1$ is when the boxes 
containing the numbers $(t-1)$ and $t$ are immediately horizontally or vertically adjacent.
However this cannot occur for the following reason: from the above construction, 
the number $t$ is in the same box in $\overline{M}$ as the number $(t-1)$ is in $\overline{P}$,
and the number $t$ is in the same box in $\overline{P}$ that the number $(t-1)$ is in $\overline{M}$.
It follows  that if the numbers $(t-1)$ and $t$ are 
immediately horizontally or vertically adjacent in $\overline{M}$, 
each must be in the corresponding `swapped' box in $\overline{P}$, 
and then at least one of $\overline{M}$ or $\overline{P}$
{\emph{cannot}} be a standard tableau.  
This contradicts the assumption that both $\overline{M}$ and $\overline{P}$
are standard tableaux, thus $|d| \neq 1$  and the coefficient in
(\ref{eq:thaiequation}) is well-defined.

It remains for us to show that the coefficient in (\ref{eq:thaiequation}) is non-zero.
This follows immediately from the fact that $|d| \neq 0$.
Note that we have not yet proved that
the matrix units are all non-zero.

Let us write $E_{MP}$ to denote $E_{\widetilde{M}\widetilde{P}}$.
We note that the matrix unit $E_{MP} \in {\cal{C}}_{t}$, where $M \neq P$, is 
an intertwiner between the isomorphic irreducible $U_{q}(\mathfrak{g})$-modules
 $E_{PP} (V^{\otimes t})$ and $E_{MM} (V^{\otimes t})$:
 $$E_{MP}: E_{PP} (V^{\otimes t}) \rightarrow E_{MM} (V^{\otimes t}),$$
and that the whole collection of matrix units satisfy
$$E_{QR}E_{ST} = \delta_{RS} E_{QT}.$$
To show that each intertwiner $E_{MP}$ is non-zero, it suffices to note that
each projector $E_{PP}$ is non-zero and that $E_{PP} = E_{PM}E_{MP}$.

We then have the complete sets of projectors and intertwiners in ${\cal{C}}_{t}$.
This means that  
${\cal{L}}_{t} = {\cal{C}}_{t}$, and also that $J_{t} = 0$.
To see this last claim, 
note that the matrix units $\{E_{ST} | \ (S,T) \in \Omega^{t}\}$ are a basis for ${\cal{C}}_{t}$.
Let $X$ be an arbitrary element of ${\cal{C}}_{t}$, then
$$X = \sum_{(S,T) \in \Omega^{t}} x_{ST} E_{ST}, \hspace{10mm}
x_{ST} \in \mathbb{C},$$ 
where $x_{ST} \neq 0$ for at least one pair $(S,T)$.  Let $(A,B)$ be such a pair, then
$$str_{q}(E_{BA}X) = str_{q}(x_{AB}E_{BA}E_{AB}) = x_{AB} str_{q}(E_{BB}) \neq 0,$$
thus $X  \notin J_{t}$.  As $X$ is arbitrary, $J_{t} = 0$.

Note that we obtained  ${\cal{C}}_{t} = {\cal{L}}_{t}$ by using the fact that
$\lambda$ and $\widetilde{\lambda}$ do not appear on the same level of the Bratteli diagram for
$\widetilde{\mathscr{BW}}_{t}(-q^{2n},q)$.  
If 
$\lambda$ and $\widetilde{\lambda}$ did appear on the same level, we could only conclude 
 from our work that there is a {\emph{proper}} inclusion of ${\cal{C}}_{t}$  
in ${\cal{L}}_{t}$ rather than an equality. 
Of course, in that event, there  may actually be an equality, but
a different method would have to be used to obtain all the intertwiners.

We now present the two lemmas  used in this section.
\begin{lemma}
\label{lemma:equalityoftraces}
Let $\psi: {\cal{C}}_{t} \rightarrow \mathbb{C}$ be a map defined by
$$\psi(X) = str_{q}(X)/\big(sdim(V)\big)^{t},$$
and let $\mathrm{tr}$ be the trace functional on $\mathscr{BW}_{t}(-q^{2n},q)$ 
mentioned in (\ref{eq:tracefunctional}).
Then 
$$\psi \big( \Upsilon (a) \big) = \mathrm{tr}(a), \hspace{10mm} \forall a \in 
\mathscr{BW}_{t}(-q^{2n},q).$$
\end{lemma}
\begin{proof}
Any functional $\phi$ on $\mathscr{BW}_{\infty}(-q^{2n},q)$ satisfying
Eq. (\ref{eq:tracefunctional}) for all $t \in \mathbb{N}$ is identical to 
$\mathrm{tr}$ \cite[Lem. 3.4 (d)]{w2}, and we will show that
$\psi \circ \Upsilon$ has this property.

To show that $\psi \circ \Upsilon$ satisfies Eq. (\ref{eq:tracefunctional}), it suffices
to show that for each $t \in \mathbb{N}$, we have
\begin{equation}
 \label{eq:julietoolie(1)}
-\psi\big(\Upsilon(a) \check{R}_{t-1} \Upsilon(b)\big) = 
-\psi\big( \check{R}_{t-1}\big) \psi\big(\Upsilon(ab)\big), \hspace{10mm}
 \forall a, b \in \mathscr{BW}_{t-1}(-q^{2n},q),
 \end{equation}
as the element
$e_{t-1} \in \mathscr{BW}_{t}(-q^{2n},q)$ can be written as a function of the $g_{t-1}$'s.  
We will show that Eq. (\ref{eq:julietoolie(1)}) is true 
using Lemma \ref{lem;kilo}, which we give after this proof.

The left hand side of Eq. (\ref{eq:julietoolie(1)}) is
\begin{equation}
\label{eq:julietoolie(2)}
-str_{q}^{\otimes t} \big(A \check{R}_{t-1} B \big)/ \big(sdim_{q}(V) \big)^{t},
\end{equation}
where we write $str_{q}^{\otimes t}$ to mean 
that we take the quantum supertrace over all $t$ tensor factors,
and we also write $A = \Upsilon(a)$ and $B = \Upsilon(b)$. 
Now we can regard each $X \in {\cal{C}}_{t-1}$ as an element of ${\cal{C}}_{t}$ under the mapping
$X \mapsto X \otimes \mathrm{id}$, 
then by applying the identity to the first $t-1$ tensor powers of (\ref{eq:julietoolie(2)}) and
taking the quantum supertrace over
the $t^{th}$ tensor power of (\ref{eq:julietoolie(2)}), we obtain, 
using Lemma \ref{lem;kilo} and applying some
simple but tedious calculations,
\begin{equation}
\label{eq:julietoolie(3)}
-str_{q}^{\otimes t} \big(A \check{R}_{t-1} B \big)/ \big(sdim_{q}(V) \big)^{t}
  = \frac{-\chi_{V}(v^{\mp 1})}{sdim_{q}(V)} \ 
  \frac{str_{q}^{\otimes (t-1)} \big(A B \big)}{ \big(sdim_{q}(V) \big)^{t-1}}.
\end{equation}
 Now 
 $$\psi\big( \check{R}_{t-1}\big) = \chi_{V}(v^{\mp 1})/sdim_{q}(V),$$
 and the right hand side of Eq. (\ref{eq:julietoolie(3)}) equals
 the right hand side of Eq. (\ref{eq:julietoolie(1)}).
Now Eq. (\ref{eq:julietoolie(1)}) 
is true for all $a$ and $b$ belonging to $\mathscr{BW}_{t-1}(-q^{2n},q)$,
and it remains to show that $\psi \circ \Upsilon$ is a functional on 
$\mathscr{BW}_{\infty}(-q^{2n},q)$ satisfying Eq. (\ref{eq:tracefunctional}) for 
all natural numbers $t$.
This follows from the fact that
$\psi(A \otimes \mathrm{id}) = \psi(A)$ for all $A \in {\cal{C}}_{t}$, thus we can regard
$\psi$ as well-defined in the inductive limit  
${\cal{C}}_{2} \subset {\cal{C}}_{3} \subset {\cal{C}}_{4} \subset \cdots$.
This completes the proof.

\end{proof}

The following lemma, which we used in the proof of Lemma \ref{lemma:equalityoftraces}, appears in
\cite[Lem. 2]{lg} and is proved in \cite[Lem. 3.1]{z1}.
\begin{lemma}
\label{lem;kilo}
Let $V$ be the fundamental irreducible
$U_{q}(osp(1|2n))$-module with highest weight $\epsilon_{1}$ and
let $\pi$ be the representation of $U_{q}(osp(1|2n))$ afforded by $V$.  
Let 
$\check{\cal{R}}_{V,V} \in End_{U_{q}(osp(1|2n))}(V \otimes V)$ be as given in
Eq. (\ref{eq:bigjimmyboy(2)}). Then
$$(\mathrm{id} \otimes \mathrm{str})\big[ (\mathrm{id} \otimes \pi)
(\mathrm{id} \otimes K_{2\rho})\big]  \check{\cal{R}}_{V,V}^{\pm 1} 
= q^{\pm(\epsilon_{1}, \epsilon_{1} + 2\rho)}\mathrm{id}=\chi_{V}(v^{\mp 1})\mathrm{id}.$$
\end{lemma}

\end{subsection}

\end{section}

\end{chapter}

%% file: young100.pstex_t
\begin{picture}(0,0)%
\includegraphics{young.pstex}%
\end{picture}%
\setlength{\unitlength}{4144sp}%
\begingroup\makeatletter\ifx\SetFigFont\undefined
\def\x#1#2#3#4#5#6#7\relax{\def\x{#1#2#3#4#5#6}}%
\expandafter\x\fmtname xxxxxx\relax \def\y{splain}%
\ifx\x\y   
\gdef\SetFigFont#1#2#3{%
  \ifnum #1<17\tiny\else \ifnum #1<20\small\else
  \ifnum #1<24\normalsize\else \ifnum #1<29\large\else
  \ifnum #1<34\Large\else \ifnum #1<41\LARGE\else
     \huge\fi\fi\fi\fi\fi\fi
  \csname #3\endcsname}%
\else
\gdef\SetFigFont#1#2#3{\begingroup
  \count@#1\relax \ifnum 25<\count@\count@25\fi
  \def\x{\endgroup\@setsize\SetFigFont{#2pt}}%
  \expandafter\x
    \csname \romannumeral\the\count@ pt\expandafter\endcsname
    \csname @\romannumeral\the\count@ pt\endcsname
  \csname #3\endcsname}%
\fi
\fi\endgroup
\begin{picture}(3284,3926)(519,-3683)
\end{picture}

%% file: chap2A100.tex
\begin{chapter}{Quantum $osp(1|2n)$ at roots of unity}
\label{chap2A:titlelabel}
\pagestyle{myheadings}
\markboth{\text{Chapter \ref{chap2A:titlelabel}. Quantum $osp(1|2n)$ at roots of unity}}
{\text{ }}

In this chapter we define a $\mathbb{Z}_{2}$-graded ribbon Hopf algebra $U_{q}^{(N)}(osp(1|2n))$
that is a certain quotient of  $U_{q}(osp(1|2n))$ 
where $q = \exp{(2 \pi i/N)}$ for some integer $N \geq 3$, 
and we also study aspects of its representation
theory.
We define certain representations of $U_{q}^{(N)}(osp(1|2n))$, show that each of these representations
is self-dual, and most importantly, 
prove tensor product decomposition theorems for these representations at even $N$.
The results in this chapter are almost entirely new.

The structure of this chapter is as follows.
In Section \ref{sec:easycooljazz(1000)} we define the quotient algebra $U_{q}^{(N)}(osp(1|2n))$
and prove that it is a $\mathbb{Z}_{2}$-graded ribbon Hopf algebra.
In Section \ref{sec:projectgirl} we define a finite number of 
finite dimensional $U_{q}^{(N)}(osp(1|2n))$-modules
and prove that the dual $U_{q}^{(N)}(osp(1|2n))$-module 
to each of these modules is isomorphic to the original module.
In Section \ref{sec:fedupwiththetensors} we prove tensor product theorems for these modules 
 at even $N$.
In Section \ref{subsec:vermittin} we present the technical proof that the projections defining the 
$U_{q}^{(N)}(osp(1|2n))$-modules are all well-defined.

In this chapter we use $\mathfrak{g}$ to denote $osp(1|2n)$, and fix
$q = \exp{(2 \pi i /N)}$ where $N \geq 3$ is an integer.

\begin{section}{The $\mathbb{Z}_{2}$-graded ribbon Hopf algebra $U_{q}^{(N)}(osp(1|2n))$}
\label{sec:easycooljazz(1000)}
\markright{\text{The algebra $U_{q}^{(N)}(osp(1|2n))$}}

In this section we introduce a quotient algebra of $U_{q}(\mathfrak{g})$ and prove that it is a
$\mathbb{Z}_{2}$-graded ribbon Hopf algebra.

We firstly generalise the $q$-bracket from Chapter \ref{chap2:titlelabel} to arbitrary elements
of $U_{q}(\mathfrak{g})$.
Each element in $U_{q}(\mathfrak{g})$ is a linear combination of products of
$K_{i}^{\pm 1}, e_{i}, f_{i}, \ i=1, 2, \ldots, n$, and
every product $X$ of the generators satisfies 
\begin{equation}
\label{chap2A:eq(1)}
K_{i} X K_{i}^{-1} =  q^{(wt(X),\alpha_{i})} X, \hspace{10mm} i = 1, 2, \ldots, n,
\end{equation}
for some integral element $wt(X) \in H^{*}$.  Needless to say, $(wt(X), \alpha_{i}) \in \mathbb{Z}$ for all $i$.
Then the $q$-bracket is a bilinear map
$[ \cdot, \cdot ]_{q}: 
U_{q}(\mathfrak{g}) \times U_{q}(\mathfrak{g}) \rightarrow U_{q}(\mathfrak{g})$ defined
by 
$$\left[X,Y\right]_{q} = XY - (-1)^{[X][Y]} q^{(wt(X),wt(Y))} YX.$$
If both $X$ and $Y$ satisfy (\ref{chap2A:eq(1)}) for some $wt(X)$ and $wt(Y)$, the meaning of
$[X, Y]_{q}$ is clear.  The definition is generalised to arbitrary elements by linearity.

\begin{subsection}{Definition of $U_{q}^{(N)}(osp(1|2n))$}
\label{subsect:normalorderofourquantumospatrootsofunity}

We now define root vectors in $U_{q}(\mathfrak{g})$ 
using a particular normal ordering of the elements of $\phi$, 
the set of positive roots of the reduced root system of $\mathfrak{g}$.
Recall that 
$\phi = \{ \epsilon_{i}, \epsilon_{j} \pm \epsilon_{k} | \ 1 \leq i \leq n, \ 1 \leq j < k \leq n \}$.
We use the following notation to help in writing elements of $\phi$ in terms of the simple roots: 
\begin{eqnarray*}{}
\alpha_{i} + \cdots + \alpha_{j} & = & \sum_{k=i}^{j} \alpha_{k},    \\
\alpha_{i} + \cdots + 2\alpha_{j} & = & \sum_{k=i}^{j-1} \alpha_{k} + 2\alpha_{j},    \\
\alpha_{i} + \cdots + 2\alpha_{j} + \cdots + 2\alpha_{n} & = &  \sum_{k=i}^{j-1} \alpha_{k} + \sum_{m=j}^{n} 2\alpha_{m}.
\end{eqnarray*}
Then
$$\begin{array}{rcll}
\epsilon_{i}   & = & \alpha_{i} + \cdots + \alpha_{n}, & i=1, \ldots, n,  \\
\epsilon_{i}-\epsilon_{j} & = & \alpha_{i} + \cdots + \alpha_{j-1}, & 1 \leq i < j \leq n, \\
\epsilon_{i}+\epsilon_{j} & = & \alpha_{i} + \cdots + \alpha_{j-1}+2\alpha_{j} +  \cdots + 2\alpha_{n},  & 1 \leq i < j \leq n.
\end{array}$$
We fix the normal order ${\cal{N}}(\phi)$ that we use in this chapter to be:
$$\begin{array}{l}
\alpha_{1} \prec \alpha_{1}+\alpha_{2} \prec \alpha_{1}+\alpha_{2}+\alpha_{3} \prec
\ldots \prec \alpha_{1}+\alpha_{2}+ \cdots + \alpha_{k} \prec
\ldots \prec \alpha_{1} + \cdots + \alpha_{n} \prec \\
\hspace{5mm} \alpha_{1} + \cdots + 2\alpha_{n} \prec 
\alpha_{1} + \cdots + 2\alpha_{n-1} + 2\alpha_{n} \prec

\ldots \prec \alpha_{1}+ \cdots + 2\alpha_{k} + \cdots + 2 \alpha_{n} \prec \\
 \hspace{5mm} \alpha_{1}+ 2\alpha_{2} + \cdots + 2 \alpha_{n} \prec \\
\alpha_{2} \prec \alpha_{2}+\alpha_{3} \prec \alpha_{2}+\alpha_{3}+\alpha_{4} \prec
\ldots \prec \alpha_{2}+\alpha_{3}+ \cdots + \alpha_{k} \prec
\ldots \prec \alpha_{2} + \cdots + \alpha_{n} \prec \\
\hspace{5mm} \alpha_{2} + \cdots + 2\alpha_{n} \prec \alpha_{2} + \cdots + 2\alpha_{n-1} + 2\alpha_{n} \prec
\ldots \prec \alpha_{2}+ \cdots + 2\alpha_{k} + \cdots + 2 \alpha_{n} \prec \\
\hspace{5mm} \alpha_{2}+ 2\alpha_{3} + \cdots + 2 \alpha_{n} \prec \\

\hspace{15mm} \vdots \\

\alpha_{j} \prec \alpha_{j} + \alpha_{j+1} \prec \ldots \prec \alpha_{j}+\alpha_{j+1}+ \cdots + \alpha_{k} \prec
\ldots \prec \alpha_{j} + \cdots + \alpha_{n} \prec \\
\hspace{5mm} \alpha_{j} + \cdots + 2\alpha_{n} \prec  \alpha_{j} + \cdots + 2\alpha_{n-1}+ 2\alpha_{n} \prec
\ldots \prec \alpha_{j}+ \cdots + 2\alpha_{k} + \cdots + 2 \alpha_{n} \prec \\
\hspace{5mm} 
\alpha_{j}+ 2\alpha_{j+1} + \cdots + 2 \alpha_{n} \prec \\

\hspace{15mm} \vdots \\

\alpha_{n-1} \prec \alpha_{n-1}+\alpha_{n} \prec \alpha_{n-1}+2\alpha_{n} \prec \alpha_{n}.
\end{array}$$
We also define a second normal order $\overline{\cal{N}}(\phi)$, which we call   
the {\emph{opposite normal order}} to ${\cal{N}}(\phi)$,
by  $y \prec x$  whenever $x \prec y$ in ${\cal{N}}(\phi)$.

Using ${\cal{N}}(\phi)$, we recursively define the root vectors 
$e_{\mu}, f_{\mu} \in U_{q}(\mathfrak{g})$  for each $\mu \in \phi$ as follows:
firstly fix $e_{\alpha_{i}} = e_{i}$ and $f_{\alpha_{i}} = f_{i}$ for each $i=1, \ldots, n$, 
then recursively define
$$\begin{array}{rcll}
e_{\alpha_{i}+\alpha_{i+1}} & = & [ e_{i}, e_{i+1}]_{q}, & i=1, \ldots, n-1,  \\
e_{\alpha_{i}+\cdots +\alpha_{j}} & = & [e_{\alpha_{i}+\cdots +\alpha_{j-1}},e_{j}]_{q}, & 
   j =i+1, \ldots, n, \\
e_{\alpha_{i}+\cdots +2\alpha_{n}} & = & [e_{\alpha_{i}+\cdots +\alpha_{n}},e_{n}]_{q}, & \\
e_{\alpha_{i}+\cdots +2\alpha_{j}+\cdots + 2\alpha_{n}} & = & 
[e_{\alpha_{i}+\cdots +2\alpha_{j+1}+\cdots + 2\alpha_{n}},e_{j}]_{q}, & j=i+1, \ldots, n-1, \\
 & & & \\
f_{\alpha_{i}+\alpha_{i+1}} & = & [ f_{i+1}, f_{i}]_{q^{-1}}, & i=1, \ldots, n-1, \\
f_{\alpha_{i}+\cdots +\alpha_{j}} & = & [f_{j},f_{\alpha_{i}+\cdots +\alpha_{j-1}}]_{q^{-1}}, & 
   j =i+1, \ldots, n, \\
f_{\alpha_{i}+\cdots +2\alpha_{n}} & = & [f_{n},f_{\alpha_{i}+\cdots +\alpha_{n}}]_{q^{-1}}, & \\
f_{\alpha_{i}+\cdots +2\alpha_{j}+\cdots + 2\alpha_{n}} & = & 
[f_{j},f_{\alpha_{i}+\cdots +2\alpha_{j+1}+\cdots + 2\alpha_{n}}]_{q^{-1}}, & j=i+1, \ldots, n-1.
\end{array}$$

Using  $\overline{\cal{N}}(\phi)$, we define a further set of elements 
$\overline{e}_{\mu}, \overline{f}_{\mu} \in U_{q}(\mathfrak{g})$ for each 
$\mu \in \phi$.  Firstly fix $\overline{e}_{i}=e_{i}$ and $\overline{f}_{i}=f_{i}$
for each $i=1, \ldots, n$, then  $\overline{e}_{\mu}$ and $\overline{f}_{\mu}$ are recursively defined by
$$\begin{array}{rcll}
\overline{e}_{\alpha_{i}+\alpha_{i+1}} & = & [  e_{i+1},e_{i}]_{q}, & i=1, \ldots, n-1, \\
\overline{e}_{\alpha_{i}+\cdots +\alpha_{j}} & = & 
        [e_{j},\overline{e}_{\alpha_{i}+\cdots +\alpha_{j-1}}]_{q}, & j =i+1, \ldots, n, \\
\overline{e}_{\alpha_{i}+\cdots +2\alpha_{n}} & = & 
     [e_{n},\overline{e}_{\alpha_{i}+\cdots +\alpha_{n}}]_{q}, & \\
\overline{e}_{\alpha_{i}+\cdots +2\alpha_{j}+\cdots + 2\alpha_{n}} & = & 
  [e_{j},\overline{e}_{\alpha_{i}+\cdots +2\alpha_{j+1}+\cdots + 2\alpha_{n}}]_{q}, & j=i+1, \ldots, n-1, \\
 & & & \\
 
\overline{f}_{\alpha_{i}+\alpha_{i+1}} & = & [f_{i}, f_{i+1} ]_{q^{-1}}, & i=1, \ldots, n-1, \\
\overline{f}_{\alpha_{i}+\cdots +\alpha_{j}} & = & 
        \Big[\overline{f}_{\alpha_{i}+\cdots +\alpha_{j-1}},f_{j}\Big]_{q^{-1}}, &  j =i+1, \ldots, n, \\
\overline{f}_{\alpha_{i}+\cdots +2\alpha_{n}} & = & 
          \Big[\overline{f}_{\alpha_{i}+\cdots +\alpha_{n}},f_{n}\Big]_{q^{-1}}, & \\
\overline{f}_{\alpha_{i}+\cdots +2\alpha_{j}+\cdots + 2\alpha_{n}} & = & 
  \Big[\overline{f}_{\alpha_{i}+\cdots +2\alpha_{j+1}+\cdots + 2\alpha_{n}},f_{j} \Big]_{q^{-1}}, & j=i+1, \ldots, n-1.
\end{array}$$
With these elements of $U_{q}(\mathfrak{g})$ 
we define the quotient algebra $U_{q}^{(N)}(\mathfrak{g})$ in Theorem \ref{th:cooljazz(201)} below.
Recall that  we fix
$$N' = \left\{ 
\begin{array}{ll}
N, & \mbox{if $N$ is odd}, \\
N/2, & \mbox{if $N$ is even},
\end{array}
\right.  \  \mbox{and }
\overline{N} = \left\{ 
\begin{array}{ll}
2N, & \mbox{if $N$ is odd}, \\
N, & \mbox{if } N \equiv 0 \pmod{4}, \\
N/2, & \mbox{if } N \equiv 2 \pmod{4}.
\end{array}
\right. $$
\begin{theorem}
\label{th:cooljazz(200)}
The left ideal ${\cal{I}} \subset U_{q}(\mathfrak{g})$
generated by the elements of the set $I$ below is a two-sided Hopf ideal of $U_{q}(\mathfrak{g})$:
\begin{equation}
I = \left\{ (e_{\gamma})^{N'}, (e_{\beta})^{\overline{N}}, 
(\overline{e}_{\gamma})^{N'}, (\overline{e}_{\beta})^{\overline{N}},
(f_{\gamma})^{N'}, (f_{\beta})^{\overline{N}}, 
(\overline{f}_{\gamma})^{N'}, (\overline{f}_{\beta})^{\overline{N}}, 
(J_{i})^{\pm N}-1 | \ 1 \leq i \leq n \right\},
\end{equation}
where $J_{i} = K_{i} K_{i+1} \cdots K_{n}$ for each $i = 1, 2, \ldots, n$, and
$\gamma$ (resp. $\beta$) ranges over all the even (resp. odd) elements of 
$\phi = \left\{\epsilon_{i}, \epsilon_{j} \pm \epsilon_{k} |  
\ 1 \leq i \leq n, \ 1 \leq j < k \leq n \right\}$.
\end{theorem}
\begin{proof}
The proof of this result is technical and very lengthy.  
Thus we relegate it to Appendix \ref{chap:appendixC}.
\end{proof}
It immediately follows from this theorem that
\begin{theorem}
\label{th:cooljazz(201)}
The quotient algebra $U_{q}^{(N)}(\mathfrak{g})$ defined by
$$U_{q}^{(N)}(\mathfrak{g}) = U_{q}(\mathfrak{g})/{\cal{I}},$$
is a $\mathbb{Z}_{2}$-graded Hopf algebra.
\end{theorem}
\noindent
We denote the image of $x \in U_{q}(\mathfrak{g})$ in $U_{q}^{(N)}(\mathfrak{g})$ 
under the canonical homomorphism  by $x$.

The algebra $U_{q}^{(N)}(\mathfrak{g})$ has appeared in the literature, 
but only for the case $n=1$ and then only for $N \geq 3$ an odd integer \cite{z2}.
The representation theory of $U_{q}^{(N)}(\mathfrak{g})$ 
is unknown except in this case, and in this case it is only partially known \cite{z2,ab}.  

We now introduce a very important representation of $U_{q}^{(N)}(\mathfrak{g})$   
that we will extensively use in this thesis: the representation afforded by the
{\emph{fundamental}} irreducible $U_{q}^{(N)}(\mathfrak{g})$-module $V$.
As a matter of notation, 
we will henceforth write $V^{gen}$ to denote the fundamental irreducible module over 
$U_{q}(\mathfrak{g})$ where $q \neq 0$ is not a root of unity.
\begin{lemma}
\label{lem:manydimroot}
There exists a $(2n+1)$-dimensional irreducible $U_{q}^{(N)}(osp(1|2n))$-module $V$ 
with highest weight $\epsilon_{1} \in {\cal{P}}^{+}$.  
Let a basis of $V$ be $\{ v_{i} | \ -n \leq i \leq n \}$ where each
$v_{i}$ is a weight vector of weight $\epsilon_{i}$, where we fix $\epsilon_{-i} = -\epsilon_{i}$
and $\epsilon_{0}=0$.  
Let $v_{1}$ be the highest weight vector of $V$.
 The action of an element $x \in U_{q}^{(N)}(osp(1|2n))$ on the weight vector $v_{i}$ 
is identical to the action of the pre-image of $x$ in $U_{q}(osp(1|2n))$ 
on the weight vector $v_{i}$ of the
fundamental $U_{q}(osp(1|2n))$-module $V$ in Lemma \ref{lem:fundamentaldimensional}
if we fix $q = \exp{(2 \pi i /N)}$.

\end{lemma}
\begin{proof}
Besides what is contained in the proof of Lemma \ref{lem:fundamentaldimensional}, 
the only additional matter that we need to prove is that elements of $I$ all act by zero and
this is easily seen to be true.
\end{proof}

We always take the grading of the highest weight vector $v_{1} \in V$ to be odd.

\end{subsection}

\begin{subsection}{The universal $R$-matrix of $U_{q}^{(N)}(osp(1|2n))$}

As the ideal ${\cal{I}} \subset U_{q}(\mathfrak{g})$ is a two-sided Hopf ideal, 
we can immediately write down the universal $R$-matrix of $U_{q}^{(N)}(\mathfrak{g})$ 
following \cite{z1}:
\begin{proposition}
\label{prop:R-matrichereicome}
The universal $R$-matrix of $U_{q}^{(N)}(\mathfrak{g})$ is
$$R = \left( \prod_{a=1}^{n} \sum_{b=0}^{N-1} (J_{a})^{b} \otimes P_{a}[b] \right) \cdot
\prod_{\gamma \in \phi} R_{\gamma},$$
where
$$P_{a}[b] = 
\prod_{\stackrel{c=0}{c \neq b}}^{N-1} \frac{J_{a} - q^{c}}{q^{b}-q^{c}},$$
and where the product $\prod_{\gamma \in \phi} R_{\gamma}$  
is ordered in accordance with the normal order ${\cal{N}}(\phi)$  used to define the root vectors, 
so that $\prod_{\gamma \in \phi} R_{\gamma} = R_{\gamma_{1}} R_{\gamma_{2}} \cdots R_{\gamma_{k}}$ 
where
${\cal{N}}(\phi) = \gamma_{1} \prec \gamma_{2} \prec \cdots \prec \gamma_{k}$.
Here, $R_{\gamma}$ is 
$$R_{\gamma} = \left\{ \begin{array}{ll}
 \displaystyle{ \sum_{k=0}^{N'-1}  \frac{ (q-q^{-1})^{k} (e_{\gamma} \otimes f_{\gamma})^{k}}{[k]^{q^{-2}}!} }, 
& \mbox{if } [e_{\gamma}]=0, \\
\displaystyle{\sum_{k=0}^{\overline{N}-1} \frac{ (q^{-1}-q)^{k} (e_{\gamma} \otimes f_{\gamma})^{k}}{[k]^{-q^{-1}}!} },
& \mbox{if } [e_{\gamma}]=1.
\end{array} \right.$$
\end{proposition}
Note that $\prod_{\gamma \in \phi} R_{\gamma}$ is well-defined and each $R_{\gamma}$ 
can be thought of as a truncation of the infinite sum of the corresponding factor of $\widetilde{R}$ 
in the universal $R$-matrix of $U_{h}(\mathfrak{g})$ so that it is well-defined for
$q$ at a root of unity. 

Universal $R$-matrices have been written down 
for quotients of other quantum algebras and quantum superalgebras at roots of unity. 
This was first done for $U_{q}(sl_{2})$ at $4k^{th}$ roots of unity where $k \in \mathbb{N}$ \cite{rt}, 
then for the quantum algebras connected with the classical series of Lie algebras at 
odd roots of unity \cite{z4} (also implied in \cite{kirjnr1}) and
then for the quantum algebras connected with
the exceptional Lie algebras $G_{2}, F_{4}, E_{8}$ at odd roots of unity \cite{z6}.
Amongst quantum superalgebras, universal $R$-matrices have been written down for
quotients of $U_{q}(osp(1|2))$ \cite{z2} and $U_{q}(gl(2|1))$ both at odd roots of unity \cite{z3}. 

An immediate consequence of Proposition \ref{prop:R-matrichereicome} is the
\begin{corollary}
The quotient algebra $U_{q}^{(N)}(\mathfrak{g})$ is a $\mathbb{Z}_{2}$-graded quasitriangular Hopf
algebra.
\end{corollary}
Write the universal $R$-matrix of $U_{q}^{(N)}(\mathfrak{g})$ as
$R = \sum_{t} a_{t} \otimes b_{t}$, then the element $u = \sum_{t} S(b_{t})a_{t} (-1)^{[a_{t}]}$
satisfies
$$\epsilon(u)=1, \hspace{10mm} \Delta(u) = (u \otimes u) \left(R^{T}R\right)^{-1}, \hspace{10mm}
S^{2}(x) = u x u^{-1}, \hspace{5mm} \forall x \in U^{(N)}_{q}(\mathfrak{g}).$$
Furthermore, we have the following important theorem.
\begin{theorem}
\label{th:ribbonalgebraribbonalgebra}
The quotient algebra $U_{q}^{(N)}(\mathfrak{g})$ is a $\mathbb{Z}_{2}$-graded ribbon Hopf algebra.
\end{theorem}
\begin{proof}
Define the even element
$$v = uK_{2\rho}^{-1} \in U_{q}^{(N)}(\mathfrak{g}).$$
It suffices to prove that $v$ is central in $U^{(N)}_{q}(\mathfrak{g})$ 
and that it satisfies the following relations:
\begin{equation}
\label{chap2:boogie(1)}
\epsilon(v)=1, \hspace{10mm} v^{2} = S(u)u, \hspace{10mm} S(v)=v,
\end{equation}
\begin{equation}
\label{chap2:boogie(2)}
\Delta(v) = (v \otimes v) \left(R^{T}R\right)^{-1}.
\end{equation}

We firstly prove that $v$ is central in $U_{q}^{(N)}(\mathfrak{g})$.
The proof is standard but we repeat it here as the elements $u$ and $v$ are crucial for later applications.
Firstly, the homomorphism $S^{2}$
satisfies $S^{2}(x) = K_{2\rho} x K_{2\rho}^{-1}$ for all $x \in U^{(N)}_{q}(\mathfrak{g})$
(see Lemma \ref{lem:galliano}).
As $v$ is invertible,
$$
v x v^{-1} = u K_{2\rho}^{-1} x K_{2\rho}u^{-1} 
           = u S^{-2}(x) u^{-1} 
	   = S^{2}(S^{-2}(x)) 
	   = x, \hspace{10mm} \forall x \in U_{q}^{(N)}(\mathfrak{g}),
$$
proving that $v$ is central in $U_{q}^{(N)}(\mathfrak{g})$.
The proofs of Eq. (\ref{chap2:boogie(2)}) and the first equation in (\ref{chap2:boogie(1)}) follow 
 from the properties of $u$, and
the proof of the third equation in (\ref{chap2:boogie(1)}) is similar to the proof of
the corresponding equation in quantum algebras \cite[Prop. 5.1]{d}. 
The second equation in (\ref{chap2:boogie(1)}) follows from the third.
\end{proof}

\begin{definition}
\label{defn:bigjim(2)}
Define $\check{\cal{R}}_{V,V} \in End_{\mathbb{C}}(V \otimes V)$  by
$\check{\cal{R}}_{V,V}(v_{i} \otimes v_{j}) = 
 P \circ (\pi \otimes \pi) R (v_{i} \otimes v_{j})$ for all $v_{i}, v_{j} \in V$; 
 $\check{\cal{R}}_{V,V}$ is an element of $End_{U_{q}^{(N)}(\mathfrak{g})}(V \otimes V)$. 
Define ${\cal{C}}_{t}$ to be the subalgebra 
of $End_{U_{q}^{(N)}(\mathfrak{g})}(V^{\otimes t})$ generated by the elements
$$\left\{\check{R}_{i}^{\pm 1} \in End_{U_{q}^{(N)}(\mathfrak{g})}(V^{\otimes t}) 
\big| \ 1 \leq i \leq t-1\right\},$$ 
where
\begin{equation}
\label{eq:jackinblackbeansauce}
\check{R}_{i} = 
\mathrm{id}^{\otimes (i-1)} \otimes \check{\cal{R}}_{V,V} \otimes \mathrm{id}^{\otimes (t-(i+1))}.
\end{equation}
\end{definition}

\begin{remark}
\label{rem:ohmygodimlabellingaremark(1)}
Note that if we take the explicit expression of $\check{R}_{i}$ in (\ref{eq:tom3}) 
and set $q$ to the appropriate root of unity, we obtain the 
$\check{R}_{i}$ defined here.

\end{remark}

\end{subsection}

\end{section}

\begin{section}{$U_{q}^{(N)}(osp(1|2n))$-modules}
\label{sec:projectgirl}
\markright{\text{$U_{q}^{(N)}(osp(1|2n))$-modules}}

In this section we define certain $U_{q}^{(N)}(\mathfrak{g})$-modules for each $N \geq 3$, 
we calculate their quantum superdimensions, and we also show that each of these modules 
is self-dual, that is, each of these modules 
is isomorphic to its dual $U_{q}^{(N)}(\mathfrak{g})$-module.

\begin{subsection}{The truncated Weyl alcoves}
\label{subsec:capitallambda}

In this subsection we define the {\emph{truncated Weyl alcoves}} 
$\Lambda_{N}^{+} \subseteq \overline{\Lambda_{N}^{+}}$,  which are
proper subsets of  $X = \bigoplus_{i=1}^{n} \mathbb{Z} \epsilon_{i} \subset H^{*}$
that we will extensively use in the definition of the $U_{q}^{(N)}(\mathfrak{g})$-modules
in Subsection \ref{subsec:7zark7}.
In the usual definition, many non-integral elements of $H^{*}$ 
belong to the truncated Weyl alcoves.
However, we only ever use the term in this thesis to refer to the relevant integral elements
given below.
The truncated Weyl alcoves are defined in terms of inequalities connected 
with the formula for the quantum superdimension 
of finite dimensional irreducible $U_{q}(\mathfrak{g})$-modules for 
$q \neq 0$  not a root of unity (see Eq. (\ref{eq:qqqsuperdim})). 
In this subsection, $(\cdot, \cdot): H^{*} \times H^{*} \rightarrow \mathbb{C}$ 
is the non-degenerate, bilinear form given by  (\ref{eq:starstar(1)}) 
and $n$ is the rank of $U_{q}^{(N)}(\mathfrak{g})$.

In this work it proves convenient to introduce the following notation: we write
$\overline{\cal{P}}^{+}_{N}$ to denote $\overline{\Lambda_{N}^{+}} \cap {\cal{P}}^{+}$.
\begin{definition}
\label{def:opencat}
We define $\overline{\Lambda_{N}^{+}}\subset X$ as follows:
\begin{itemize}
\item[(i)] for $N \equiv 0, 1, 3 \pmod{4}$,
$$\overline{\Lambda_{N}^{+}} = \left\{ \lambda \in X \bigg| \ 
0 \leq \frac{2(\lambda + \rho,\alpha)}{(\alpha,\alpha)} \leq N', \ \forall \alpha \in 
\Phi^{+}_{0}   \right\},$$
\item[(ii)] for $N \equiv 2 \pmod{4}$,
$$\overline{\Lambda_{N}^{+}} = \left\{ \lambda \in X \bigg| \ 
0 \leq \frac{2(\lambda + \rho,\alpha)}{(\alpha,\alpha)} \leq N', \ \forall \alpha \in 
\overline{\Phi}_{0}^{+} \cup \Phi_{1}^{+}   \right\}.$$
\end{itemize}
\end{definition}

\begin{definition}
\label{def:opendog}
We define $\Lambda_{N}^{+} \subseteq \overline{\Lambda_{N}^{+}}$ as follows:
\begin{itemize}
\item[(i)] for $N \equiv 0, 1, 3 \pmod{4}$,
\begin{equation}
\label{eq:leah(1)egg}
\Lambda_{N}^{+} = \left\{ \lambda \in X \bigg| \
0 < \frac{2(\lambda + \rho,\alpha)}{(\alpha,\alpha)} < N', \ \forall \alpha \in 
\Phi^{+}_{0}   \right\},
\end{equation}
\item[(ii)] for $N \equiv 2 \pmod{4}$,
\begin{equation}
\label{eq:leah(2)egg}
\Lambda_{N}^{+} = \left\{ \lambda \in X \bigg| \
0 < \frac{2(\lambda + \rho,\alpha)}{(\alpha,\alpha)} < N', \ \forall \alpha \in 
\overline{\Phi}_{0}^{+} \cup \Phi_{1}^{+}   \right\}.
\end{equation}
\end{itemize}
\end{definition}

\begin{lemma}
\label{lem:subsetsofXarethesame}
For each $\lambda \in X$ let $\lambda_{i} = (\lambda,\epsilon_{i})$ for each $i=1, \ldots, n$.
There is an alternative description of  
$ \Lambda_{N}^{+}$ and $\overline{\cal{P}}^{+}_{N}=\overline{\Lambda_{N}^{+}} \cap {\cal{P}}^{+}$:
\begin{itemize}
\item[(i)] $ \Lambda_{N}^{+} =  \left\{
\begin{array}{lll}
\left\{ \lambda \in {\cal{P}}^{+}| \ \lambda_{1} + \lambda_{2} < N'-2n+2
\right\},  & \mbox{when } N \equiv 0, 1, 3 \pmod{4}, & \mbox{for } n \geq 2, \\
\left\{ \lambda \in {\cal{P}}^{+}| \ \lambda_{1}  < N'
\right\},  & \mbox{when } N \equiv 0, 1, 3 \pmod{4}, & \mbox{for } n=1, \\
\left\{ \lambda \in {\cal{P}}^{+}| \ \lambda_{1} < N/4-n+1/2
\right\},  &\mbox{when } N \equiv 2 \pmod{4}. &
\end{array}
\right. $
\item[(ii)] $\overline{\cal{P}}^{+}_{N}  = 
\left\{
\begin{array}{lll}
\left\{ \lambda \in {\cal{P}}^{+}| \ \lambda_{1} + \lambda_{2} \leq N'-2n+2
\right\}, & \mbox{when } N \equiv 0, 1, 3 \pmod{4}, & \mbox{for } n \geq 2, \\
\left\{ \lambda \in {\cal{P}}^{+}| \ \lambda_{1}  \leq N'
\right\}, & \mbox{when } N \equiv 0, 1, 3 \pmod{4}, & \mbox{for } n=1, \\
\left\{ \lambda \in {\cal{P}}^{+}| \ \lambda_{1} \leq N/4-n+1/2
\right\}, & \mbox{when } N \equiv 2 \pmod{4}. &
\end{array}
\right. $
\end{itemize}
\end{lemma}
\begin{proof}
We will only prove (i). 
The proof of (ii) is similar and will be omitted.  
We will assume throughout this proof that $\lambda \in \Lambda_{N}^{+}$.
Set $N \equiv 0, 1, 3 \pmod{4}$.  
We will firstly show that $\lambda$ must be an element of ${\cal{P}}^{+}$ if it is an element of
$\Lambda_{N}^{+}$. Fixing $\alpha = 2 \epsilon_{n}$ in
(\ref{eq:leah(1)egg}) shows that $\lambda$ satisfies
$$ 0 < \frac{2(\lambda + \rho,2\epsilon_{n})}{(2\epsilon_{n},2\epsilon_{n})} < N',$$ 
which is just $0 < \lambda_{n} + 1/2 < N'$; note that $\lambda_{n} \geq 0$ as $\lambda \in X$. 
This completes the proof for $n=1$ for $N \equiv 0, 1, 3 \pmod{4}$.  
For $n \geq 2$ we use the following arguments. 
Fixing $\alpha = \epsilon_{i} - \epsilon_{i+1}$ for $i=1, \ldots, n-1$
in (\ref{eq:leah(1)egg}) shows that $\lambda$ satisfies
$$0 < \frac{2(\lambda + \rho, \epsilon_{i} - \epsilon_{i+1})}
{(\epsilon_{i} - \epsilon_{i+1},\epsilon_{i} - \epsilon_{i+1})} < N',$$ 
which tells us that  $\lambda_{i} - \lambda_{i+1}  \geq 0$ for each $i$.  
Then $\lambda$ must be an element of ${\cal{P}}^{+}$ as $\lambda \in X$ and $\lambda_{n} \geq 0$.
Fixing $\alpha =\epsilon_{1} + \epsilon_{2}$ in  (\ref{eq:leah(1)egg}) 
shows that 
\begin{equation}
\label{eq:learhpurcell(1)}
0 < \frac{2(\lambda + \rho,\epsilon_{1} + \epsilon_{2})}
{(\epsilon_{1} + \epsilon_{2},\epsilon_{1} + \epsilon_{2})} < N',
\end{equation}
and so 
$0 \leq \lambda_{1} + \lambda_{2} < N'-2n+2$ as $\lambda \in {\cal{P}}^{+}$.
Fix $\alpha = \epsilon_{i} \pm \epsilon_{j}$ for $i<j$ where 
$\alpha \neq \epsilon_{1} + \epsilon_{2}$,   then
$$\frac{2(\lambda + \rho, \alpha)}{(\alpha,\alpha)} 
< \frac{2(\lambda + \rho, \epsilon_{1} + \epsilon_{2})}
{(\epsilon_{1} + \epsilon_{2},\epsilon_{1} + \epsilon_{2})}, \hspace{5mm}
\mbox{for each } \alpha \neq \epsilon_{1} + \epsilon_{2}.$$
Finally, note that 
 $$\frac{2(\lambda +
 \rho,2\epsilon_{i})}{(2\epsilon_{i},2\epsilon_{i})} = (\lambda + \rho,\epsilon_{i})
 < \frac{2(\lambda + \rho,\epsilon_{1} + \epsilon_{2})}
     {(\epsilon_{1} + \epsilon_{2},\epsilon_{1} + \epsilon_{2})}.$$
From this, any $\lambda \in {\cal{P}}^{+}$ satisfying (\ref{eq:learhpurcell(1)}) also belongs to
$\Lambda_{N}^{+}$, and so
 $$\Lambda_{N}^{+}  = \big\{\lambda \in {\cal{P}}^{+} | \ 
 \lambda_{1} + \lambda_{2} < N'-2n+2 \big\}, \hspace{10mm} \mbox{when } N \equiv 0, 1, 3 \pmod{4}.$$

Let us consider the case that $N \equiv 2 \pmod{4}$; again, fix $\lambda \in \Lambda_{N}^{+}$.
As above, we can show that $\lambda \in {\cal{P}}^{+}$. 
Fixing $\alpha = \epsilon_{1}$ in (\ref{eq:leah(2)egg}) shows that $\lambda$ satisfies
\begin{equation}
\label{jamestree}
0 < \frac{2(\lambda+\rho,\epsilon_{1})}{(\epsilon_{1},\epsilon_{1})} < N',
\end{equation} 
which means that $0 \leq \lambda_{1} < N/4-n+1/2$ as $\lambda \in  {\cal{P}}^{+}$.
It is also true that
$$\frac{2(\lambda+\rho,\alpha)}{(\alpha,\alpha)} < \frac{2(\lambda+\rho,\epsilon_{1})}{(\epsilon_{1},\epsilon_{1})},
\hspace{5mm} \forall \alpha = \epsilon_{i} \pm \epsilon_{j}, \epsilon_{k} 
\hspace{5mm} i < j, \ k \geq 2.$$
From this, any $\lambda \in {\cal{P}}^{+}$ satisfying (\ref{jamestree}) also belongs to
$\Lambda_{N}^{+}$, and so
 $$\Lambda_{N}^{+}  = \big\{\lambda \in {\cal{P}}^{+} | \ 
 \lambda_{1}  < N/4-n+1/2 \big\}, \hspace{5mm} \mbox{for } N \equiv 2 \pmod{4},$$
which completes the proof of part (i).
\end{proof}

\end{subsection}

\begin{subsection}{$U_{q}^{(N)}(osp(1|2n))$-modules}
\label{subsec:7zark7}

We now define the $U_{q}^{(N)}(\mathfrak{g})$-modules of interest in our work, and in doing so  
we  use the truncated Weyl alcoves.

We have already defined the fundamental irreducible $U_{q}^{(N)}(\mathfrak{g})$-module $V$ for all 
$n \geq 1$ and all $N \geq 3$.  
In this subsection we define a further set of $U_{q}^{(N)}(\mathfrak{g})$-modules 
$V_{\lambda}$ for all  
$\lambda \in \overline{\cal{P}}^{+}_{N}$ when $N \geq 4$ is even. 
For technical reasons, we only define the $U_{q}^{(N)}(\mathfrak{g})$-modules 
 for $\lambda$  an element of a proper subset of $\overline{\cal{P}}^{+}_{N}$ when $N \geq 3$ is odd.  

We define each of these new $U_{q}^{(N)}(\mathfrak{g})$-modules $V_{\lambda}$ by using a projection
operator $V^{\otimes t} \rightarrow V_{\lambda}$.  
This projection operator is an element of
$End_{U_{q}^{(N)}(\mathfrak{g})}(V^{\otimes t})$ that we obtain by setting
$q$ to the appropriate root of unity in the projection
$(V^{gen})^{\otimes t} \rightarrow V_{\lambda}^{gen}$ 
we defined in Section \ref{subsec:projectontome}.
Here, $V^{gen}$ and $V_{\lambda}^{gen}$ are finite dimensional irreducible 
$U_{q}(\mathfrak{g})$-modules where $q \neq 0$ is not a root of unity, 
and the highest weight of $V_{\lambda}^{gen}$ is $\lambda \in {\cal{P}}^{+}$.
It is crucially important that these projections are well defined and we 
deal with this problem in Section \ref{subsec:vermittin}.

The projections are all well-defined for even $N$, 
but for odd $N$ they are only well-defined if $\lambda$ 
belongs to a proper subset of $\overline{\cal{P}}^{+}_{N}$; 
this is the reason we only define these modules for $\lambda$ an element of a proper subset of 
$\overline{\cal{P}}^{+}_{N}$ when $N$ is odd.
We have not been able to resolve this ill-definedness problem, but we conjecture that
well-defined projections do exist 
for all $\lambda \in \overline{\cal{P}}^{+}_{N}$ at odd $N$.

Recall that the Bratteli diagram for $(V^{gen})^{\otimes t}$ 
defined in Chapter \ref{chap2:titlelabel} encodes the decomposition
of $(V^{gen})^{\otimes t}$ into irreducible $U_{q}(\mathfrak{g})$-submodules 
for $q \neq 0$ not a root of unity.
Recall further that the elements on the $j^{th}$ level of the Bratteli diagram for
$(V^{gen})^{\otimes t}$, for each $j \leq t$, are called
shapes, and that a sequence of $(t+1)$ elements:
$$\lambda_{i}^{t} = (0, \epsilon_{1}, s_{2},  \ldots, s_{t-1}, \lambda),$$
is called a tableau of length $t$ if 
$s_{j}$ is a shape on the $j^{th}$ level of the Bratteli diagram for each $j$.
We let ${\cal{T}}^{t}$ denote the set of all tableaux of length $t$ derived from 
the Bratteli diagram for $(V^{gen})^{\otimes t}$.

We now define two proper subsets of ${\cal{T}}^{t}$ that we will use 
in creating the projections from $V^{\otimes t}$ onto 
$U_{q}^{(N)}(\mathfrak{g})$-submodules of $V^{\otimes t}$.
\begin{definition}
Define the two proper subsets $\widetilde{\cal{T}}^{t}$ and $\widehat{\cal{T}}^{t}$ of 
${\cal{T}}^{t}$ by
$$\widetilde{\cal{T}}^{t} = 
\left\{ i^{t} =(s_{0}, s_{1}, \ldots,  s_{t}) \in {\cal{T}}^{t} | \ 
s_{j} \in \Lambda_{N}^{+} \ \mbox{for all } 0 \leq j \leq t \right\},$$
$$\widehat{\cal{T}}^{t} = \left\{
i^{t} =(s_{0}, s_{1}, \ldots,  s_{t}) \in {\cal{T}}^{t} \big| \ 
 s_{j} \in \Lambda_{N}^{+} \  \mbox{for each } 0 \leq j \leq t-1,  \mbox{and }
 s_{t} \in \overline{\cal{P}}^{+}_{N} \right\}.$$  
 \end{definition}
\begin{itemize}
\item
The set $\widetilde{\cal{T}}^{t}$ is the set of all tableaux of length $t$ 
where each shape in each tableau is an element of $\Lambda_{N}^{+}$.
\item
The set $\widehat{\cal{T}}^{t}$ is the set of all tableaux of length $t$ 
where the first $t$ shapes in each tableau are elements of $\Lambda_{N}^{+}$ 
and the last shape in each tableau is an element of $\overline{\cal{P}}^{+}_{N}$.
\item
Note that $\widetilde{\cal{T}}^{t}$ is a proper subset of $\widehat{\cal{T}}^{t}$.
Also,
it is convenient to write $\lambda_{i}^{t}$ to mean $i^{t}$ if $shp(i^{t}) = \lambda$.
\end{itemize}

Let $\lambda_{i}^{t} =(0, \epsilon_{1}, s_{2}, \ldots, s_{t-1}, \lambda) \in \widehat{\cal{T}}^{t}$ 
be a tableau and let 
$\tilde{p}_{i}^{t}[\lambda]^{gen} \in End_{U_{q}(\mathfrak{g})}  \left(V^{gen}\right)^{\otimes t}$ be
a path projection at generic $q$:
\begin{equation}
\label{eq:brontemargaretthrosby(1)}
\tilde{p}_{i}^{t}[\lambda]^{gen}: \left(V^{gen}\right)^{\otimes t} \rightarrow V_{\lambda}^{gen}
\subset \left(V^{gen}\right)^{\otimes t},
\end{equation}
then $\tilde{p}_{i}^{t}[\lambda]^{gen}$ can be written as the ratio $z_{1}/z_{2}$ where $z_{1}$ is
an ordered polynomial in the elements 
$\left\{\check{R}_{i}^{\pm 1} 
\in End_{U_{q}(\mathfrak{g})}\left(V^{gen}\right)^{\otimes t} | 
 \ i=1, \ldots, t-1 \right\}$ 
with coefficients in 
$\left\{ \pm q^{-(\mu + 2\rho, \mu)}  | \ \mu \in {\cal{P}}^{+}  \right\}$, 
and $z_{2}$ is a non-zero product of elements of the form 
$\left\{ q^{-(\mu + 2\rho, \mu)} - q^{-(\nu + 2\rho, \nu)} | \ \mu, \nu \in {\cal{P}}^{+} \right\}$.

\begin{definition}
\label{def:lessthangreaterthanwhatever}
For each tableau 
$\lambda_{i}^{t} =(0, \epsilon_{1}, s_{2}, \ldots, s_{t-1}, \lambda) \in \widehat{\cal{T}}^{t}$,
we define the path projection
\begin{equation}
\label{eq:markhopkins(yyy)}
\tilde{p}_{i}^{t}[\lambda] \in End_{U^{(N)}_{q}(\mathfrak{g})} 
\big( V^{\otimes t} \big),
\end{equation}
to be an identical expression in the $\check{R}^{\pm 1}_{i}$ and the $\pm q^{-(\mu + 2\rho, \mu)}$ 
as in the definition of $\tilde{p}_{i}^{t}[\lambda]^{gen}$ in (\ref{eq:brontemargaretthrosby(1)}) 
except that we fix $q = \exp{(2 \pi i /N)}$ and we also fix
$\check{R}^{\pm 1}_{i} \in End_{U_{q}^{(N)}(\mathfrak{g})}\left(V^{\otimes t}\right)$.
\end{definition}

Explicitly, the projection $\tilde{p}_{i}^{t}[\lambda]$ from (\ref{eq:markhopkins(yyy)}) is 
$\tilde{p}_{i}^{t}[\lambda] = 
 p^{0}_{s_{t-1}}[\lambda]p^{1}_{s_{t-2}}[s_{t-1}] \cdots p^{t-2}_{\epsilon_{1}}[s_{2}]$, where
$$p^{j}_{s_{t-(j+1)}}[s_{t-j}] = \pi^{\otimes t}
\left( \prod_{\stackrel{\nu \in {\cal{P}}^{+}_{s_{t-(j+1)}} }{ \nu \neq s_{t-j}}}
\frac{\Delta^{(t-1-j)}(v) - q^{-(\nu + 2\rho, \nu)}}{q^{-(s_{t-j} + 2\rho, s_{t-j})} - q^{-(\nu + 2\rho, \nu)}}
 \otimes  \mathrm{id}^{\otimes j}\right), \hspace{5mm} j=0, 1, \ldots, t-2,$$
 where $\pi$ is the representation of $U_{q}^{(N)}(\mathfrak{g})$
 afforded by the fundamental irreducible module $V$.

\begin{lemma}
\label{lem:projectionsarewelldefined}
For each tableau 
$\lambda_{i}^{t} =(0, \epsilon_{1}, s_{2}, \ldots, s_{t-1}, \lambda) \in \widehat{\cal{T}}^{t}$,
the path projection 
$\tilde{p}_{i}^{t}[\lambda] \in End_{U^{(N)}_{q}(\mathfrak{g})} \left( V^{\otimes t} \right)$ 
is well-defined if either one of the two following conditions is met:
\begin{itemize}
\item[(a)] $N \geq 4$ is even, or
\item[(b)] $N \geq 3$ is odd and 
\begin{itemize}
\item[(i)] $\lambda_{1} \leq (N-1)/2-n+1$, or 
\item[(ii)]  the components of $s_{t-1}= \overline{\lambda} \in \Lambda_{N}^{+}$ satisfy 
$\overline{\lambda}_{1} = (N-1)/2-n+1$ and $\overline{\lambda}_{2} = \overline{\lambda}_{1}$, and $\lambda$ is given by 
$\lambda = \overline{\lambda} + \epsilon_{1}$.
\end{itemize}
\end{itemize}
\end{lemma}

The proof of this lemma is easy but very lengthy.  To avoid interrupting the flow of thought we relegate the proof to 
Section \ref{subsec:vermittin} at the very end of this chapter.  Here we wish to make the following remarks.  
Part (a) of the lemma means that at even $N \geq 4$, 
all projections that project down by a path in $\widehat{\cal{T}}^{t}$ 
onto a $U_{q}^{(N)}(\mathfrak{g})$-module labelled by an element of
$\overline{\cal{P}}^{+}_{N}$ are well-defined. 
Part (b) means that at odd $N \geq 3$, only certain of the path projections are well defined.  
The significance of (b) is that we cannot prove tensor product theorems of the form of 
Theorem \ref{th:secondtensor} for $U_{q}^{(N)}(\mathfrak{g})$-modules at odd $N$.
As mentioned previously, we conjecture that at odd $N$ there 
exist well-defined projections onto each module labelled by an element of 
$\overline{\cal{P}}^{+}_{N}$.

In the proof of Lemma \ref{lem:likj} we use an argument that we will repeatedly use throughout
this thesis, and here we wish to make some relevant comments.  
The idea in this argument is to take an element
$p^{gen}$ belonging to $End_{\mathbb{C}}(V^{gen})^{\otimes t}$ at generic $q \neq 0$, 
and to then specialise
$q$ to a root of unity 
yielding the element $p \in End_{\mathbb{C}}(V^{\otimes t})$.
It is  essential that $p^{gen}$ is well-defined upon specialising $q$ to the desired root of unity.

Let us write $\mathbb{C}[[q,q^{-1} ]]$ to mean the ring of power series in $q$ and $q^{-1}$ with
coefficients in $\mathbb{C}$.
Let $\mathscr{V}$ be a vector subspace of $V^{\otimes t}$ over $\mathbb{C}[[q,q^{-1} ]]$ with a basis
$$\big\{v_{i_{1}} \otimes v_{i_{2}} \otimes \cdots \otimes v_{i_{t}} \left| \right. \
i_{1}, i_{2}, \ldots, i_{t} = -n, \ldots, n  \big\},$$
where each $v_{i}$ is a weight vector of the fundamental $U_{q}(\mathfrak{g})$-module $V^{gen}$.
In order that $p^{gen}$ is well-defined when $q$ is specialised to the root of unity,
we only ever apply the specialisation
argument to 
\begin{itemize}
\item[(i)] those $p^{gen}$ that we  know are well-defined if $q$ is specialised, or
\item[(ii)] those $p^{gen}$ belonging to $End_{\mathbb{C}[[q,q^{-1}]]}(\mathscr{V})$.
\end{itemize}

In the proof of Lemma \ref{lem:likj} (amongst other lemmas and propositions in this thesis), we claim
that certain elements in $End_{U_{q}^{(N)}(\mathfrak{g})}(V^{\otimes t})$ can be obtained by 
taking corresponding elements in $End_{U_{q}(\mathfrak{g})}(V^{gen})^{\otimes t}$ and
specialising $q$ to the appropriate root of unity.
In particular, we claim  this  for
$\pi^{\otimes t} (\Delta^{(t-1)}(v))$ (where $v \in U_{q}^{(N)}(\mathfrak{g})$) and for
the well-defined projections $\tilde{p}_{i}^{t}[\lambda]$ of ${\cal{C}}_{t}$
from Lemma \ref{lem:projectionsarewelldefined}.

To see this, firstly note that the matrices $\check{R}_{i}^{\pm 1}$ of
$End_{U_{q}^{(N)}(\mathfrak{g})}(V^{\otimes t})$ and $(\check{R}_{i}^{gen})^{\pm 1}$ of
$End_{U_{q}(\mathfrak{g})}(V^{gen})^{\otimes t}$
 are exactly the same relative to the basis
$$\big\{ 
v_{i_{1}} \otimes v_{i_{2}} \otimes \cdots \otimes v_{i_{t}} | \ 
i_{1}, i_{2}, \ldots, i_{t} = -n, \ldots, n \big\}$$ 
of both $V^{\otimes t}$ and $(V^{gen})^{\otimes t}$
if we consider $q$ to be an indeterminate 
(see Remark \ref{rem:ohmygodimlabellingaremark(1)}).
The same is true for
$\pi^{\otimes t}\big(\Delta^{(k)}(v) \otimes \mathrm{id}^{\otimes (t-k-1)} \big)$ for each
$k=0, 1, \ldots, t-1$ as it
is a product in the $\check{R}_{i}^{\pm 1}$ with coefficients in 
$\mathbb{C}[[q,q^{-1}]]$ in both cases.

Now the projections $\tilde{p}_{i}^{t}[\lambda]$ of ${\cal{C}}_{t}$
in Lemma \ref{lem:projectionsarewelldefined} are polynomials in products of the 
$\pi^{\otimes t}\big(\Delta^{(k)}(v) \otimes \mathrm{id}^{\otimes (t-k-1)} \big)$ 
with coefficients in $\mathbb{C}(q)$, but we know 
from Lemma \ref{lem:projectionsarewelldefined} that they are all well-defined if $q$ is specialised to 
the appropriate root of unity.

We will use specialisation arguments to prove various claims throughout this thesis. 
In each of these proofs we will rely on the comments here and we will
show that any coefficients in $\mathbb{C}(q)$ appearing in the calculations are
 well-defined if $q$ is specialised to the appropriate root of unity.

\begin{lemma}
\label{lem:likj}
Let  $\lambda_{i}^{t} =(0, \epsilon_{1}, s_{2}, \ldots, s_{t-1}, \lambda) \in \widehat{\cal{T}}^{t}$ 
be a tableau of length $t$
and $\tilde{p}_{i}^{t}[\lambda] \in {\cal{C}}_{t}$ be a well-defined projection 
referred to in Lemma \ref{lem:projectionsarewelldefined}, then
\begin{itemize}
\item[(i)] $\big(\tilde{p}_{i}^{t}[\lambda]\big)^{2} = \tilde{p}_{i}^{t}[\lambda]$,
\item[(ii)] $\tilde{p}_{i}^{t}[\lambda] \cdot \tilde{p}_{j}^{t}[\lambda] = 
\left\{  
\begin{array}{ll}
0,                          & \mbox{if } i^{t} \neq j^{t}, \\
\tilde{p}_{i}^{t}[\lambda], & \mbox{if } i^{t} = j^{t},
\end{array}
\right. $
\item[(iii)]  $\tilde{p}_{i}^{t}[\lambda] \cdot \tilde{p}_{j}^{t}[\mu]=0$ if $\lambda \neq \mu$.
\end{itemize}
\end{lemma}
\begin{proof}  

To prove this lemma we consider similar equations in 
$End_{U_{q}(\mathfrak{g})}(V^{gen})^{\otimes t}$ where $q \neq 0$ is not a root of unity,
and we apply a specialisation argument.

Consider the elements
$\tilde{p}_{i}^{t}[\lambda]^{gen} \in End_{U_{q}(\mathfrak{g})} \left(V^{gen}\right)^{\otimes t}$
that project down from $(V^{gen})^{\otimes t}$ onto 
irreducible $U_{q}(\mathfrak{g})$-submodules with highest weights in ${\cal{P}}^{+}$.
For each such $\tilde{p}_{i}^{t}[\lambda]^{gen}$, we can define
a matrix valued function $M_{\overline{q}}\left(\tilde{p}_{i}^{t}[\lambda]^{gen}\right)$
of a complex parameter $\overline{q}$, for all $\overline{q} \neq 0$ that are not roots of unity, 
such that 
$$M_{\overline{q}}\left(\tilde{p}_{i}^{t}[\lambda]^{gen}\right)\Big|_{\overline{q}=q} = 
\tilde{p}_{i}^{t}[\lambda]^{gen} \in End_{U_{q}(\mathfrak{g})} \left(V^{gen}\right)^{\otimes t},$$
for each $q \neq 0$ that is not a root of unity.
Each component of the matrix $M_{\overline{q}}\left(\tilde{p}_{i}^{t}[\lambda]^{gen}\right)$ 
is a continuous function in $\overline{q}$ and has no poles for all non-zero $\overline{q}$ that are not
roots of unity.  
Furthermore, 
$$\lim_{\overline{q} \rightarrow q} 
\big[ M_{\overline{q}}\left(\tilde{p}_{i}^{t}[\lambda]^{gen}\right) \big] = 
M_{\overline{q}}\left(\tilde{p}_{i}^{t}[\lambda]^{gen}\right)\Big|_{\overline{q}=q},$$ 
where in taking the limit we take the limit of each component of 
$M_{\overline{q}}\left(\tilde{p}_{i}^{t}[\lambda]^{gen}\right)$.

Let 
$\lambda_{i}^{t} =(0, \epsilon_{1}, s_{2}, \ldots, s_{t-1}, \lambda) \in \widehat{\cal{T}}^{t}$
be a tableau of length $t$ where the path projection
$\tilde{p}_{i}^{t}[\lambda] \in End_{U^{(N)}_{q}(\mathfrak{g})} \left( V^{\otimes t} \right)$
is well-defined from Lemma \ref{lem:projectionsarewelldefined}.
Then we can extend the domain in $\overline{q}$ of 
$M_{\overline{q}}\left(\tilde{p}_{i}^{t}[\lambda]^{gen}\right)$
to include $\overline{q} = \exp{(2 \pi i/N)}$ by defining
$$M_{e^{ 2 \pi i/N }}\left(\tilde{p}_{i}^{t}[\lambda]^{gen}\right) = 
\tilde{p}_{i}^{t}[\lambda] \in End_{U^{(N)}_{q}(\mathfrak{g})} \left(V\right)^{\otimes t}.$$
No component of $M_{e^{ 2 \pi i/N }}\left(\tilde{p}_{i}^{t}[\lambda]^{gen}\right)$  has a pole,
and furthermore,
$\tilde{p}_{i}^{t}[\lambda] = M_{e^{ 2 \pi i/N }}\left(\tilde{p}_{i}^{t}[\lambda]^{gen}\right) = 
\lim_{\overline{q} \rightarrow  e^{2 \pi i/N}}\big[ M_{\overline{q}}\left(\tilde{p}_{i}^{t}[\lambda]^{gen}\right) \big]$.

Now consider the equations in $End_{U_{q}(\mathfrak{g})} \left(V^{gen}\right)^{\otimes t}$ 
given in Lemma \ref{lem:greatgnats} (i)--(iii) corresponding to parts (i)--(iii) of this lemma.
For each non-zero $q$ that is not a root of unity, the equations in Lemma \ref{lem:greatgnats} (i)--(iii) are true.
As each component of $M_{\overline{q}}\left(\tilde{p}_{i}^{t}[\lambda]^{gen}\right)$ is 
a continuous function of $\overline{q}$ and has no poles for all  $\overline{q}$,
the product of two such matrix valued functions is again a matrix valued function 
where each component of the product is a continuous function in $\overline{q}$ 
with no poles for all $\overline{q}$.
It follows then that each equation in 
Lemma \ref{lem:greatgnats} (i)--(iii) can be obtained by taking the limit  
$\overline{q} \rightarrow q \in \mathbb{C}$ of the
products of the matrix valued functions $M_{\overline{q}}\left(\tilde{p}_{i}^{t}[\lambda]^{gen}\right)$.

The proof of Eqs. (i)--(iii) of this lemma is then given by taking the limit
$\overline{q} \rightarrow e^{2 \pi i/N}$ of the corresponding equations in 
$End_{U_{q}(\mathfrak{g})} \left(V^{gen}\right)^{\otimes t}$.
\end{proof}

We now define the specific collection of new $U_{q}^{(N)}(\mathfrak{g})$-modules of interest.
\begin{definition}
\label{lem:representationsyippe}
Let $\lambda_{i}^{t} =(0, \epsilon_{1}, s_{2}, \ldots, s_{t-1}, \lambda) \in \widehat{\cal{T}}^{t}$ 
be a tableau and let $\tilde{p}_{i}^{t}[\lambda] \in {\cal{C}}_{t}$ be a well-defined projection 
referred to in Lemma \ref{lem:projectionsarewelldefined}.
We define the finite dimensional $U_{q}^{(N)}(\mathfrak{g})$-module $V_{\lambda}$ 
depending on the tableau $\lambda_{i}^{t}$ by
\begin{equation}
\label{eq:markhopkins(zzz)}
V_{\lambda} = \tilde{p}_{i}^{t}[\lambda] \big( V^{\otimes t} \big).
\end{equation}
\end{definition}
We conjecture that each of the $U_{q}^{(N)}(\mathfrak{g})$-modules defined in 
(\ref{eq:markhopkins(zzz)}) is an irreducible $U_{q}^{(N)}(\mathfrak{g})$-module.
We now investigate properties of these $U_{q}^{(N)}(\mathfrak{g})$-modules.
Note that ${\cal{P}}^{+}_{\mu} \subseteq \overline{\Lambda_{N}^{+}}$ for each
$\mu \in \Lambda_{N}^{+}$, and we write ${\cal{P}}^{+}_{\mu} \cap \overline{\Lambda_{N}^{+}}$
below to make clear that each $\lambda \in {\cal{P}}^{+}_{\mu}$ is also an element of
$\overline{\Lambda_{N}^{+}}$.
\begin{lemma}
\label{lemma:heresyetanothermoviewlemmayouknowwhatImean}
Let $i^{t}=(0, \epsilon_{1}, s_{2}, \ldots, s_{t-1}, \mu) \in \widetilde{\cal{T}}^{t}$ 
be a tableau of length $t$  and let
$V_{\mu}$ be a $U_{q}^{(N)}(\mathfrak{g})$-module defined by 
$V_{\mu} = \tilde{p}^{t}_{i}[\mu] \big(V^{\otimes t}\big)$ as given in 
Definition \ref{lem:representationsyippe}.
Then there is a decomposition of $V_{\mu} \otimes V$ into a direct sum of
$U_{q}^{(N)}(\mathfrak{g})$-submodules:
\begin{equation}
\label{eq:moomoo}
V_{\mu} \otimes V = 
\bigoplus_{\lambda \in {\cal{P}}^{+}_{\mu} \cap \overline{\Lambda_{N}^{+}}} V_{\lambda},
\end{equation}
where each $V_{\lambda}$ is a $U_{q}^{(N)}(\mathfrak{g})$-submodule defined by
$V_{\lambda} = \tilde{p}_{j}^{(t+1)}[\lambda] \left( V^{\otimes (t+1)} \right)$, where
$\tilde{p}_{j}^{(t+1)}[\lambda] \in End_{U_{q}^{(N)}(\mathfrak{g})}(V^{\otimes (t+1)})$ 
 is a path projection of length $t+1$, and
$j^{(t+1)}=(0, \epsilon_{1}, s_{2},  \ldots, s_{t-1}, \mu, \lambda) \in \widehat{\cal{T}}^{(t+1)}$ 
is a tableau of length $t+1$.  
\end{lemma}
\begin{proof}
Firstly, note that ${\cal{P}}^{+}_{\mu} \subseteq \overline{\Lambda_{N}^{+}}$ for each
$\mu \in \Lambda_{N}^{+}$.
Lemma \ref{lem:likj} (ii) implies the only vector belonging to any pair of 
distinct summands on the right
hand side of Eq. (\ref{eq:moomoo}) is the zero vector. 
As $j^{(t+1)} \in \widehat{\cal{T}}^{(t+1)}$, the path projection
$\tilde{p}_{j}^{(t+1)}[\lambda] \in End_{ U_{q}^{(N)}(\mathfrak{g})  }(V^{\otimes (t+1)})$ is well-defined 
for each $\lambda \in {\cal{P}}^{+}_{\mu} \cap \overline{\Lambda_{N}^{+}}$ from 
Lemma \ref{lem:projectionsarewelldefined}.
To complete the proof we need only show 
that the inclusion
$\displaystyle{\bigcup_{\lambda \in {\cal{P}}^{+}_{\mu} \cap \overline{\Lambda_{N}^{+}}}
V_{\lambda} \subseteq V_{\mu} \otimes V}$ is actually an equality,
and this follows from the equation
$$ \sum_{ \lambda \in {\cal{P}}^{+}_{\mu} \cap \overline{\Lambda_{N}^{+}} } 
\tilde{p}_{j}^{(t+1)}[\lambda]
= \mathrm{id}_{V_{\mu} \otimes V},$$ 
which can be shown to be true by applying an argument similar to the one we used in the proof 
of Lemma \ref{lem:likj} to the corresponding matrix equations at generic $q$.
\end{proof}

We can prove the following lemma
by also applying a similar idea to the one we used in the proof of Lemma \ref{lem:likj}
to the corresponding matrix equations at generic $q$.
\begin{lemma}
\label{chapt2lem:anniegoingaway(1)}
Let $\lambda_{i}^{t} \in \widehat{\cal{T}}^{t}$ be a tableau of length $t$ 
and let $V_{\lambda}$ be a $U_{q}^{(N)}(\mathfrak{g})$-module defined by 
$V_{\lambda} = \tilde{p}^{t}_{i}[\lambda] \big(V^{\otimes t}\big)$.
Then
$$v \cdot w = q^{-(\lambda + 2\rho, \lambda)}w, \hspace{10mm} \forall w \in V_{\lambda},$$
where $v = u K_{2\rho}^{-1} \in U_{q}^{(N)}(\mathfrak{g})$.
\end{lemma}

\begin{lemma}
\label{lem:kilcutiecutie}
Let $\lambda \in \overline{\cal{P}}_{N}^{+}$ 
and let $V_{\lambda}$ be the $U_{q}^{(N)}(\mathfrak{g})$-module given in 
Definition \ref{lem:representationsyippe}, then
\begin{itemize}
\item[(i)]  the quantum superdimension of $V_{\lambda}$ is
\begin{equation}
\label{eq:qsuperdim}
sdim_{q}(V_{\lambda}) = (-1)^{[\lambda]}q^{-(\lambda,2\rho)}
\prod_{\alpha \in \overline{\Phi}^{+}_{0}}\left(
\frac{q^{2(\lambda+\rho,\alpha)}-1}{q^{2(\rho,\alpha)}-1}\right)
\prod_{\beta \in \Phi_{1}^{+}}\left(
\frac{q^{2(\lambda+\rho,\beta)}+1}{q^{2(\rho,\beta)}+1}\right),
\end{equation}
where $[\lambda]$ is the grading of the highest weight vector of
the irreducible $U_{q}(\mathfrak{g})$-module $V_{\lambda}^{gen}$ with highest weight
$\lambda \in {\cal{P}}^{+}$ where $q \neq 0$ is not a root of unity,
\item[(ii)] $sdim_{q}(V_{\lambda}) \neq 0$ if $\lambda \in \Lambda_{N}^{+}$,
\item[(iii)] $sdim_{q}(V_{\lambda})=0$ if 
$\lambda \in \overline{{\cal{P}}}^{+}_{N} \backslash \Lambda_{N}^{+}$.
\end{itemize}
\end{lemma}
\begin{proof}  
\begin{itemize}
\item[(i)]
The quantum superdimension of $V_{\lambda} \subseteq V^{\otimes t}$ is
$$sdim_{q}(V_{\lambda}) = 
\mathrm{str} 
\big( \tilde{p}_{i}^{t}[\lambda] \cdot \pi^{\otimes t} \big(  \Delta^{(t-1)}(K_{2 \rho}) \big) \big),$$
and by using the idea in the proof of Lemma \ref{lem:likj}, one can show that this
gives the right hand side of (\ref{eq:qsuperdim}) provided that it is well-defined.
This is indeed the case for all $\lambda \in \overline{\cal{P}}_{N}^{+}$
as the denominator of the right hand side of (\ref{eq:qsuperdim}) is non-zero.
This is very easy to see, but nevertheless we present the detailed proof.

Set $N \equiv 0, 1, 3 \pmod{4}$.  
For each $\alpha  \in \overline{\Phi}^{+}_{0}$  we have $2(\rho,\alpha) \in 2\mathbb{Z}$
and $0 < 2(\rho,\alpha) < 2N'$, which implies that  $q^{2(\rho,\alpha)}\neq 1$.  
Now if $N$ is odd, $q^{2(\rho,\beta)} \neq -1$ for all $\beta \in \Phi_{1}^{+}$. 
If $N \equiv 0 \pmod{4}$ then $q^{2(\rho,\beta)} = -1$ for some $\beta \in \Phi_{1}^{+}$ if and only if
$2(\rho,\beta) \equiv N/2 \pmod{N}$.  However, it is not possible that $2(\rho,\beta) \equiv N/2 \pmod{N}$, as
$2(\rho,\beta)$ is odd and both of $N$ and $N/2$ are even, thus $q^{2(\rho,\beta)}  \neq -1$ for all $\beta \in \Phi_{1}^{+}$.
Now if $N \equiv 2 \pmod{4}$ then $0 < 2(\rho,\alpha) < N$ for all $\alpha  \in \overline{\Phi}^{+}_{0}$, and 
$0 < 2(\rho,\beta) < N/2$ for all $\beta \in \Phi_{1}^{+}$. 

It follows that the right hand side of (\ref{eq:qsuperdim}) is well defined.

\item[(ii)]  
If $\lambda \in \Lambda_{N}^{+}$, none of the factors in the numerator of the right hand side of 
Eq. (\ref{eq:qsuperdim}) is zero.  The proof of this fact is easy thus omitted.

\item[(iii)]   Consider $\lambda \in \overline{{\cal{P}}}^{+}_{N} \backslash \Lambda_{N}^{+}$.  
From the definitions of $\overline{\Lambda_{N}^{+}}$ and $\Lambda_{N}^{+}$, we obtain
$$\begin{array}{ll}
2(\lambda + \rho, \alpha) = 2N', \hspace{5mm} \mbox{for some } \alpha \in \overline{\Phi}^{+}_{0}, & \mbox{if } N \equiv 0, 1, 3 \pmod{4},
\\
2(\lambda+\rho, \epsilon_{1}) = N/2, & \mbox{if } N \equiv 2 \pmod{4}.
\end{array} $$
Thus the right hand side of Eq. (\ref{eq:qsuperdim}) is zero.
\end{itemize}
\end{proof}

\begin{proposition}
\label{prop:moo20}
Let $\lambda \in \overline{\cal{P}}^{+}_{N}$ and let $V_{\lambda}$ 
be the $U_{q}^{(N)}(\mathfrak{g})$-module
given in Definition \ref{lem:representationsyippe}, 
then $V_{\lambda}$ is self-dual.
\end{proposition}
\begin{proof}
It suffices to show that there is a non-degenerate, $U_{q}^{(N)}(\mathfrak{g})$-invariant,
bilinear form $( \cdot, \cdot): V_{\lambda} \times V_{\lambda} \rightarrow \mathbb{C}$; and
we do this inductively.

Firstly, we will show that there exists such a form on $V_{\nu} \times V_{\nu}$
for any $U_{q}^{(N)}(\mathfrak{g})$-summand $V_{\nu}$
on the right hand side of the decomposition
\begin{equation}
\label{eq:paulnewman(00)}
V \otimes V = V_{2 \epsilon_{1}} \oplus V_{\epsilon_{1} + \epsilon_{2}} \oplus V_{0},
\end{equation}
where we write $\epsilon_{1} + \epsilon_{2}$ to mean $\epsilon_{1}$ if $n=1$.
Let $\langle \langle \ , \ \rangle \rangle: V \times V \rightarrow \mathbb{C}$ 
be the non-degenerate bilinear form
from Proposition \ref{prop:mooV}, the $U_{q}^{(N)}(\mathfrak{g})$-invariance of which is given by
$$\langle \langle   a \cdot x, y \rangle \rangle 
= (-1)^{[a][x]} \langle \langle  x,  S(a) y \rangle \rangle, 
\hspace{10mm} \forall a \in U_{q}^{(N)}(\mathfrak{g}), \hspace{5mm} x, y \in V.$$
Now define a new bilinear form
$\langle \langle \ , \ \rangle \rangle: (V \otimes V) \times (V \otimes V) \rightarrow \mathbb{C}$
by 
$$\langle \langle x_{1} \otimes y_{1}, x_{2} \otimes y_{2} \rangle \rangle
 = (-1)^{[y_{1}][x_{2}]} 
 \langle \langle x_{1},x_{2} \rangle \rangle \langle \langle y_{1}, y_{2}\rangle \rangle,
 \hspace{10mm} x_{1}, x_{2}, y_{1}, y_{2} \in V,$$
which is evidently non-degenerate.  
Elementary calculations show that this form satisfies 
$$\big\langle \big\langle   a \cdot (x_{1} \otimes y_{1}), x_{2} \otimes y_{2}  \big\rangle \big\rangle 
= (-1)^{[a]([x_{1}] + [y_{1}])} 
\big\langle \big\langle x_{1} \otimes y_{1},  
                 \Delta'(S(a))  (x_{2} \otimes y_{2}) \big\rangle \big\rangle, \hspace{5mm} 
\forall a \in U_{q}^{(N)}(\mathfrak{g}),$$
but this is not our desired $U_{q}^{(N)}(\mathfrak{g})$-invariance.
To deal with this, we introduce a new bilinear form that is non-degenerate and has 
the desired $U_{q}^{(N)}(\mathfrak{g})$-invariance.
Define the new bilinear form
$\langle \langle \ , \ \rangle \rangle^{new}: (V \otimes V) \times (V \otimes V) \rightarrow \mathbb{C}$
by
$$\big\langle \big\langle  x_{1} \otimes y_{1}, x_{2} \otimes y_{2}  \big\rangle \big\rangle^{new} =
\big\langle \big\langle  x_{1} \otimes y_{1}, 
    R \cdot (x_{2} \otimes y_{2})  \big\rangle \big\rangle, \hspace{10mm}
x_{1}, x_{2}, y_{1}, y_{2} \in V,$$
where $R$ is the universal $R$-matrix of $U_{q}^{(N)}(\mathfrak{g})$.
This new form is non-degenerate as $R$ is invertible and
$\langle \langle \ , \ \rangle \rangle$ is non-degenerate.  
We now claim that
$\langle \langle \ , \ \rangle \rangle^{new}$ is $U_{q}^{(N)}(\mathfrak{g})$-invariant.
To see this, note that for each $a \in U_{q}^{(N)}(\mathfrak{g})$ we have
\begin{eqnarray}
\lefteqn{
\langle \langle a \cdot (x_{1} \otimes y_{1}), x_{2} \otimes y_{2}  \rangle \rangle^{new} } \nonumber \\
& & = \langle \langle a \cdot (x_{1} \otimes y_{1}), R \cdot (x_{2} \otimes y_{2})  \rangle \rangle \nonumber \\
& & = \sum_{(a),t} 
\langle \langle  a_{(1)} x_{1} \otimes a_{(2)} y_{1}, \alpha_{t} x_{2} \otimes \beta_{t} y_{2} \rangle \rangle
(-1)^{[a_{2}][x_{1}] + [\beta_{t}][x_{2}]} \nonumber \\
& & = \sum_{(a),t} 
\langle \langle a_{(1)} x_{1}, \alpha_{t} x_{2} \rangle \rangle 
\langle \langle a_{(2)} y_{1}, \beta_{t} y_{2} \rangle \rangle 
(-1)^{([a_{(2)}] + [y_{1}]) ([\alpha_{t}]+[x_{2}]) + [a_{(2)}][x_{1}] + [\beta_{t}][x_{2}]}, \label{chapXX:eq(1)}
\end{eqnarray}
where we write the universal $R$-matrix as $R = \sum_{t} \alpha_{t} \otimes \beta_{t}$.
Using the $U_{q}^{(N)}(\mathfrak{g})$-invariance of 
$\langle \langle \ , \ \rangle \rangle: V \times V \rightarrow \mathbb{C}$,
we can rewrite (\ref{chapXX:eq(1)}) as
\begin{eqnarray*}
\lefteqn{
\langle \langle x_{1} \otimes y_{1}, ((S \otimes S) \Delta(a))  R \cdot (x_{2} \otimes y_{2}) \rangle \rangle
(-1)^{[a]([y_{1}]+[x_{1}])} } \\
& & =
\langle \langle x_{1} \otimes y_{1}, \Delta'(S(a))  R \cdot (x_{2} \otimes y_{2}) \rangle \rangle
(-1)^{[a]([y_{1}]+[x_{1}])} \\
& & = \langle \langle x_{1} \otimes y_{1}, R \cdot \Delta(S(a)) (x_{2} \otimes y_{2}) \rangle \rangle
(-1)^{[a]([y_{1}]+[x_{1}])} \\
& & = \langle \langle x_{1} \otimes y_{1}, \Delta(S(a))  (x_{2} \otimes y_{2}) \rangle \rangle^{new}
(-1)^{[a]([y_{1}]+[x_{1}])} \\
& & = \langle \langle x_{1} \otimes y_{1}, S(a) \cdot (x_{2} \otimes y_{2}) \rangle \rangle^{new}
(-1)^{[a]([y_{1}]+[x_{1}])},
\end{eqnarray*}
thus 
$\langle \langle \ , \ \rangle \rangle^{new}: V^{\otimes 2} \times V^{\otimes 2} \rightarrow \mathbb{C}$
is $U_{q}^{(N)}(\mathfrak{g})$-invariant.
We now show that $\langle \langle \ , \ \rangle \rangle^{new}$ 
is non-degenerate on each of the $U_{q}^{(N)}(\mathfrak{g})$-summands of $V \otimes V$ 
on the right hand side of (\ref{eq:paulnewman(00)}).
Let $\lambda, \nu \in {\cal{P}}^{+}_{\epsilon_{1}} \cap \overline{\Lambda_{N}^{+}}$ 
be non-equal and let
$x \in V_{\lambda}$ and $y \in V_{\nu}$ be arbitrary non-zero vectors.  
The fact that the projections from $V \otimes V$ onto its summands 
are well-defined comes about from the fact that
$\chi_{\lambda}(v) \neq \chi_{\nu}(v)$ for all these $\lambda$ and $\nu$ satisfying $\lambda \neq \nu$.  
Furthermore, we have $\chi_{\lambda}(v) \neq 0 \neq \chi_{\nu}(v)$, and it follows that
$$\chi_{\lambda}(v) \langle \langle x, y \rangle \rangle = 
\langle \langle v \cdot x, y \rangle \rangle = 
\langle \langle x, v \cdot y \rangle \rangle = 
\chi_{\nu}(v) \langle \langle x, y \rangle \rangle, \hspace{5mm}
\mbox{if and only if } \hspace{2mm} \langle \langle x, y \rangle \rangle =0,$$
where we have used the fact that $S(v)=v$.
The non-degeneracy of $\langle \langle \ , \ \rangle \rangle^{new}: 
V^{\otimes 2} \times V^{\otimes 2} \rightarrow \mathbb{C}$ then implies that 
$\langle \langle \ , \ \rangle \rangle^{new}$ is non-degenerate on $V_{\nu} \times V_{\nu}$ for 
each summand $V_{\nu}$ on the right hand side of (\ref{eq:paulnewman(00)}).

We have shown that 
$\langle \langle \ , \ \rangle \rangle^{new}: V^{\otimes 2} \times V^{\otimes 2} \rightarrow \mathbb{C}$
is $U_{q}^{(N)}(\mathfrak{g})$-invariant, thus the form
$\langle \langle \ , \ \rangle \rangle^{new}: V_{\nu} \otimes V_{\nu} \rightarrow \mathbb{C}$
is non-degenerate, $U_{q}^{(N)}(\mathfrak{g})$-invariant, and bilinear as desired.

Now we do the inductive step, using an almost identical argument.  
Let $V_{\mu}$, $\mu \in \Lambda_{N}^{+}$, be a $U_{q}^{(N)}(\mathfrak{g})$-module
from Definition \ref{lem:representationsyippe} and let $V_{\mu}$ be equipped with a 
non-degenerate, $U_{q}^{(N)}(\mathfrak{g})$-invariant, bilinear form
$\langle \langle \ , \ \rangle \rangle_{\mu}: V_{\mu} \times V_{\mu} \rightarrow \mathbb{C}$,
where the $U_{q}^{(N)}(\mathfrak{g})$-invariance is 
$$\langle \langle a \cdot x, y \rangle \rangle_{\mu} = (-1)^{[a][x]} \langle \langle x, S(a) y \rangle \rangle_{\mu}, 
\hspace{10mm} \forall a \in U_{q}^{(N)}(\mathfrak{g}), \ x,y \in V_{\mu}.$$

Now let $V_{\lambda}$, $\lambda \in {\cal{P}}^{+}_{\mu} \cap \overline{\Lambda_{N}^{+}}$, be the 
$U_{q}^{(N)}(\mathfrak{g})$-module defined in Definition \ref{lem:representationsyippe} by
$V_{\lambda} =  \tilde{p}_{j}^{(t+1)}[\lambda](V_{\mu} \otimes V)$.
Define a bilinear form 
\begin{equation}
\label{eq:moo5}
\langle \langle \ , \ \rangle \rangle: 
(V_{\mu} \otimes V) \times (V_{\mu} \otimes V) \rightarrow \mathbb{C}, \hspace{5mm} \mbox{by}
\end{equation}
$$\langle \langle x_{1} \otimes y_{1}, x_{2} \otimes y_{2} \rangle \rangle = 
(-1)^{[y_{1}][x_{2}]} \langle \langle x_{1}, x_{2} \rangle \rangle_{\mu}  \langle \langle y_{1}, y_{2} \rangle \rangle, 
\hspace{5mm} x_{1}, x_{2} \in V_{\mu}, \hspace{5mm} y_{1}, y_{2} \in V.$$
As
$\langle \langle \ , \ \rangle \rangle_{\mu}: V_{\mu} \times V_{\mu} \rightarrow \mathbb{C}$ 
is non-degenerate, so is the form in (\ref{eq:moo5}).
Now define a new bilinear form 
$\langle \langle \ , \ \rangle \rangle^{new}: (V_{\mu} \otimes V) \times (V_{\mu} \otimes V) \rightarrow \mathbb{C}$
by
$$\langle \langle x_{1} \otimes y_{1}, x_{2} \otimes y_{2} \rangle \rangle^{new} = 
\langle \langle x_{1} \otimes y_{1}, R \cdot (x_{2} \otimes y_{2}) \rangle \rangle,$$
where $R$ is the universal $R$-matrix, then by an almost identical argument as before,
$\langle \langle \ , \ \rangle \rangle^{new}$ is also non-degenerate, and furthermore, is 
$U_{q}^{(N)}(\mathfrak{g})$-invariant:
$$\big\langle \big\langle a \cdot (x_{1} \otimes y_{1}), 
     x_{2} \otimes y_{2} \big\rangle \big\rangle^{new} = 
 (-1)^{[a]([x_{1}] + [y_{1}])} \big\langle \big\langle x_{1} \otimes y_{1}, 
 \Delta(S(a))  (x_{2} \otimes y_{2})  \big\rangle \big\rangle^{new}, \hspace{5mm}
 \forall a \in U_{q}^{(N)}(\mathfrak{g}).$$

Recall that $V_{\mu} \otimes V$ decomposes into the following direct sum of $U_{q}^{(N)}(\mathfrak{g})$-modules:
\begin{equation}
\label{eq:moo4}
V_{\mu} \otimes V = \bigoplus_{\lambda \in {\cal{P}}^{+}_{\mu} \cap \overline{\Lambda_{N}^{+}}} V_{\lambda}.
\end{equation}
Using almost exactly the same argument as previously, we can show that
$\langle \langle \ , \ \rangle \rangle^{new}$ 
is non-degenerate on each of the summands on the right hand side of (\ref{eq:moo4}).
Now the fact that 
$\langle \langle \ , \ \rangle \rangle^{new}: 
(V_{\mu} \otimes V) \times (V_{\mu} \otimes V) \rightarrow \mathbb{C}$ 
is $U_{q}^{(N)}(\mathfrak{g})$-invariant implies that
$\langle \langle \ , \ \rangle \rangle^{new}: V_{\lambda} \times V_{\lambda} \rightarrow \mathbb{C}$
is also $U_{q}^{(N)}(\mathfrak{g})$-invariant, thus
we have our non-degenerate, $U_{q}^{(N)}(\mathfrak{g})$-invariant, bilinear form 
$\langle \langle \ , \ \rangle \rangle^{new}: V_{\lambda} \times V_{\lambda} \rightarrow \mathbb{C}$
as desired.

One easily completes the proof of the proposition using induction. 

\end{proof}

\end{subsection}

\end{section}

\begin{section}{Tensor products of $U_{q}^{(N)}(osp(1|2n))$-modules}
\label{sec:fedupwiththetensors}
\markright{\text{Tensor products of $U_{q}^{(N)}(osp(1|2n))$-modules}}

In this section we prove some of the most important results of this chapter. 
We present certain tensor product theorems
for the $U_{q}^{(N)}(\mathfrak{g})$-modules $V_{\lambda}$, 
 where $\lambda \in \Lambda_{N}^{+}$, in the following cases:
\begin{itemize}
\item[(i)]  $n =1$ and $N \geq 6 $ satisfies $N \equiv 2 \pmod{4}$, and
\item[(ii)] $n \geq 2$ and $N \geq 4$ is even.
\end{itemize}
We then give the detailed proofs of these tensor product theorems.
Note that we do not consider the case that $n=1$ and $N \geq 4$ satisfies $N \equiv 0 \pmod{4}$.
We leave this case for future study, but conjecture that the results in this case are similar to
the results for the cases that we do consider.

In obtaining this result we use many of the techniques of Sections 4 and 5 of \cite{tw},
in which a corresponding result was obtained for modules of
quantum algebras associated with the $A, B, C$ and $D$ families 
of Lie algebras at even roots of unity.

\begin{subsection}{Technical Lemmas}

Recall that ${\cal{C}}_{t}$ is the subalgebra over $\mathbb{C}$ 
of $End_{U_{q}^{(N)}(\mathfrak{g})}(V^{\otimes t})$ generated by 
$$\check{R}_{i} = 
\mathrm{id}^{\otimes (i-1)} \otimes \check{\cal{R}}_{V,V} \otimes \mathrm{id}^{\otimes (t-(i+1))},
\hspace{5mm} 1 \leq i \leq t-1.
$$
Also recall that 
\begin{itemize}
\item ${\cal{T}}^{t}$ is the set of all tableaux of length $t$ derived from the Bratteli diagram
for $(V^{gen})^{\otimes t}$, 
\item $\widetilde{\cal{T}}^{t}$ is the proper subset of ${\cal{T}}^{t}$ consisting
of all those sequences $\lambda^{t}_{i} = (0, \epsilon_{1}, s_{2},  \ldots, s_{t}) \in {\cal{T}}^{t}$
where $s_{i} \in \Lambda_{N}^{+}$ for each $1 \leq i \leq t$, 
\item
$\widehat{\cal{T}}^{t}$ is a further proper subset of ${\cal{T}}^{t}$ consisting of all those sequences 
$\lambda^{t}_{i} = (0, \epsilon_{1}, s_{2}, \ldots, s_{t}) \in {\cal{T}}^{t}$
where $s_{i} \in \Lambda_{N}^{+}$ for each $1 \leq i \leq t-1$ and $s_{t} \in \overline{\cal{P}}^{+}_{N}$.
\end{itemize}

We now define some
matrix units in ${\cal{C}}_{t}$ that will play a key role in the proof of the tensor product
theorems for certain $U_{q}^{(N)}(\mathfrak{g})$-modules later in this section.
In constructing these matrix units in ${\cal{C}}_{t}$, it is convenient to define the set 
$$\widetilde{\Omega}^{t} = \left\{(S, T) \in \Omega^{t} | \ S, T \in \widetilde{\cal{T}}^{t},
 shp(S) = shp(T) \right\}$$ 
of pairs of paths in $\widetilde{\cal{T}}^{t}$ that end at the same shape.

Now for each $(S, T) \in \widetilde{\Omega}^{t}$ we define a matrix unit 
$E_{ST} \in {\cal{C}}_{t}$
using precisely the same method we used to define the matrix unit
$E_{ST} \in End_{U_{q}(\mathfrak{g})}(V^{gen})^{\otimes t}$ 
in Subsection \ref{subsection:genericcalmatrixunitsinCt}, except that here
we fix $q = \exp{(2 \pi i/N)}$ and also fix $\check{R}_{i}$ 
to be the appropriate
 element of $End_{U_{q}^{(N)}(\mathfrak{g})}(V^{\otimes t})$.
Note that we only define a matrix unit $E_{ST}$ in ${\cal{C}}_{t}$ here 
for each pair $(S, T) \in \widetilde{\Omega}^{t}$.

We need to show that these matrix units in ${\cal{C}}_{t}$ are all well-defined and non-zero.
Firstly, each of the projectors 
$\left\{E_{SS} \in {\cal{C}}_{t} | \ (S, S) \in  \widetilde{\Omega}^{t} \right\}$ 
is defined by $E_{SS} = \tilde{p}^{t}_{i}[\lambda] \in {\cal{C}}_{t}$ where 
$S = \lambda^{t}_{i} \in \widetilde{T}^{t}$, and these are well-defined and non-zero by construction.
It takes more work to prove the corresponding result for the intertwiners 
$\left\{E_{ST} \in {\cal{C}}_{t} | \ (S, T) \in  \widetilde{\Omega}^{t}, S \neq T \right\}$, 
and we will do this inductively.  

Assume firstly that the matrix units 
$\left\{E_{ST} \in {\cal{C}}_{r} | \ (S, T) \in  \widetilde{\Omega}^{r} \right\}$, 
are all well-defined and non-zero for some positive integer $r \geq 1$.
We will show that the intertwiners
$$\left\{E_{MP} \in {\cal{C}}_{r+1} | \ (M, P) \in  \widetilde{\Omega}^{r+1}, M \neq P \right\}$$
are also all well-defined and non-zero.
To do this, we firstly recall from Subsection \ref{subsection:genericcalmatrixunitsinCt}
that each path $\widetilde{S}$ of length $t$ in the Bratteli diagram for 
$(V^{gen})^{\otimes t}$ corresponds to a path $S$ of length $t$ in the Bratteli diagram for
$\widetilde{\mathscr{BW}}_{t}(-q^{2n},q)$ and that this is a one-to-one relationship.
In the following discussion concerning intertwiner matrix units, in discussing a path 
$\widetilde{S}$ of length $t$ in $\widetilde{\cal{T}}^{t}$, we
always use this instead to refer to the corresponding path $S$ of length $t$ in the Bratteli diagram for
$\widetilde{\mathscr{BW}}_{t}(-q^{2n},q)$. 
The reason we do this is that although the intertwiner matrix units are defined using paths in the
Bratteli diagram for $\widetilde{\mathscr{BW}}_{t}(-q^{2n},q)$,
it is easier to discuss paths in $\widetilde{\cal{T}}^{t}$.

We now define the intertwiners $E_{MP} \in {\cal{C}}_{r+1}$; 
let us partition them into two sets: 
\begin{itemize}
\item[(a)] the intertwiners $E_{MP}$ for which $|shp(M)| = |shp(P)| = r+1$, 
\item[(b)] the intertwiners for which $|shp(M)| = |shp(P)| < r+1$.
\end{itemize}
We will show that the intertwiners in each set are well-defined and non-zero.
Firstly consider (a): let the paths $M$ and $P$ satisfy $|shp(M)| = |shp(P)| = r+1$ and 
$shp(M') \neq shp(P')$, then the intertwiner $E_{MP} \in {\cal{C}}_{r+1}$ is well-defined 
if the coefficient
\begin{equation}
\label{eq:johnfaulkner(1)}
\frac{1-q^{2d}}{\sqrt{(1-q^{2d+2})(1-q^{2d-2})}},
\end{equation}
is well-defined, where $d = d(\overline{M}, r)$ 
is an integer defined in Subsection \ref{subsec:rhapsodyinred}
and $\overline{M}$ is a particular path in $\widetilde{\cal{T}}^{r+1}$ also
defined in Subsection \ref{subsec:rhapsodyinred} such that $shp(\overline{M}) = shp(M)$. 
The number $|d(\overline{M}, r)|+1$ is the length of a 
hook going through the boxes containing the numbers $r$ and $r+1$ 
in the standard tableau obtained from $\overline{M}$ in the canonical way.

The question of
whether the coefficient (\ref{eq:johnfaulkner(1)}) is well-defined and non-zero 
evidently depends on the values that $d$ can take, 
and we will show that (\ref{eq:johnfaulkner(1)}) is indeed well-defined and non-zero.  
As $|d|+1$ is the length of a hook in the standard tableau obtained from $\overline{M}$, we can
calculate the values that $d$ can take by considering all possible hooks of $shp(\overline{M})$.

Clearly, the minimum length of a hook is $2$, ie $|d|+1 \geq 2$. 
We now break the problem down into a number of sub-cases:
\begin{itemize}
\item[(i)] $n \geq 2$ and $N \equiv 0 \pmod{4}$.

Here $\lambda \in \Lambda_{N}^{+}$ if and only if 
$\lambda$ is an element of ${\cal{P}}^{+}$ satisfying $0 \leq \lambda_{1} + \lambda_{2} \leq N/2-2n+1$.  
We want to find the greatest possible length of a hook over 
all allowable Young diagrams that also satisfy $0 \leq \lambda_{1} + \lambda_{2} \leq N/2-2n+1$.
Recall from Chapter \ref{chap2:titlelabel} that a Young diagram $\mu$ is 
said to be allowable if
$\mu'_{1} + \mu'_{2} \leq 2n+1$.
By considering the geometry of the relevant Young diagrams 
we can see that the greatest such length appears in a hook in any allowable
Young diagram $\lambda$  
satisfying $\lambda_{1} = N/2-2n$ and $\lambda'_{1}=2n$.
Now the greatest possible length of any hook in any such Young diagram is 
$2n+ N/2-2n-1$, thus $d$ satisfies $2 \leq |d|+1 \leq N/2-1$.

\item[(ii)] $n \geq 2$ and $N \equiv 2 \pmod{4}$.

Here $\lambda \in \Lambda_{N}^{+}$ if and only if $\lambda$ 
is an element of ${\cal{P}}^{+}$ satisfying $0 \leq \lambda_{1} \leq N/4-n-1/2$.
We want to find the greatest possible length of a hook over 
all allowable Young diagrams that also satisfy $0 \leq \lambda_{1} \leq N/4-n-1/2$.
By considering the geometry of the relevant Young diagrams 
we can easily see that the greatest such length appears in a hook in any allowable 
Young diagram $\lambda$ satisfying
$\lambda_{1} = N/4-n-1/2$ and $\lambda'_{1}=2n$.
Now the greatest possible length of a hook in any such Young diagram is 
$2n + N/4-n-3/2$, 
thus $d$ satisfies $2 \leq |d|+1 \leq N/4 +n - 3/2$.

\item[(iii)] $n=1$ and $N$ is even.

Intertwiner matrix units $E_{MP} \in {\cal{C}}_{r+1}$ do not exist in this sub-case. 
Each element of $\Lambda_{N}^{+}$ belongs to $\mathbb{Z}_{+} \epsilon_{1}$ and
after examining the geometry of the relevant Young diagrams,
the Bratteli diagram shows us that the path $M$ is
$M = (0, \epsilon_{1}, 2\epsilon_{1}, \ldots, (r+1) \epsilon_{1}) \in \widetilde{\cal{T}}^{r+1}$
as this is the only path in $\widetilde{\cal{T}}^{r+1}$ that ends on a shape with $r+1$ boxes.
As this is the only such path, there is no $P$ and no intertwiner $E_{MP}$.
\end{itemize}

This allows us to analyse the cofficient (\ref{eq:johnfaulkner(1)}) of $E_{MP}$.
Firstly, (\ref{eq:johnfaulkner(1)}) is well-defined if $q^{2d \pm 2} \neq 1$.
By examining the values of $d$ from (i)--(ii) above, 
it is not difficult to see that $q^{2d+2} \neq 1$ and $q^{2d-2} \neq 1$ 
unless $|d|=1$.  
However, in Subsection \ref{subsec:rhapsodyinred} 
we showed that $|d| \geq 2$, thus $q^{2d \pm 2} \neq 1$.  

We now show that (\ref{eq:johnfaulkner(1)}) is non-zero: 
to do this it suffices to show that $q^{2d} \neq 1$.
In sub-case (i) above we have 
$$2 \leq |2d| \leq N-4,$$
and in (ii) we have
$$2 \leq |2d| \leq N-8,$$
thus it is indeed true that $q^{2d} \neq 1$ and therefore (\ref{eq:johnfaulkner(1)}) is non-zero.
Note that we have not proved that the intertwiner $E_{MP}$ itself is non-zero, but we will show this
 momentarily.

Now we consider the remaining intertwiners in case (b): 
the intertwiners $E_{MP}$ where $|shp(M)| = |shp(P)| < r+1$.
The coefficient of $E_{MP}$ is  
\begin{equation}
\label{eq:johnfaulkner(51)}
\frac{sdim_{q}(V_{shp(M)})}{\sqrt{sdim_{q}(V_{shp(M')}) sdim_{q}(V_{shp(P')})}},
\end{equation}
which is well-defined and non-zero as $M, P \in \widetilde{\cal{T}}^{r+1}$,
each vertex on the paths $M$ and $P$ 
is an element of $\Lambda_{N}^{+}$ and
$sdim_{q}(V_{\lambda}) \neq 0$ for all $\lambda \in \Lambda_{N}^{+}$.

Now we need to prove that $E_{MP} \in {\cal{C}}_{r+1}$ is non-zero at 
the appropriate root of unity.
We showed that the coefficient of $E_{MP}$ is non-zero at roots of unity in our work above.
Now recall that 
$E_{MP} E_{PM} = E_{MM}$ at all generic $q$, then by using arguments similar to those 
we used in the proof of Lemma \ref{lem:likj}, we can show that $E_{MP}$ 
is non-zero when $q$ is specialised to
the appropriate root of unity, as $E_{PM}$ has no poles at these roots of unity.  
Furthermore,  $str_{q}(E_{PM}E_{MP}) = str_{q}(E_{PP}) \neq 0$, so $E_{MP} \neq 0$.
Thus the intertwiner $E_{MP}$ is well-defined and non-zero at roots of unity.
This completes the proof that the intertwiners in 
$\left\{ E_{MP} \in {\cal{C}}_{r+1} \left| \right. \ (M,P) \in \widetilde{\Omega}^{r+1} \right\}$ 
are all well-defined and non-zero.

We have  shown that all the 
$\left\{E_{ST} \in {\cal{C}}_{t} | \ (S, T) \in  \widetilde{\Omega}^{t} \right\}$ are well-defined
and non-zero when $q = \exp{(2 \pi i/N)}$.  
Now we wish to show that 
\begin{lemma}
\label{lem:minimalprojectors(1)}
Each projector
$E_{SS} \in {\cal{C}}_{t}$, $S \in  \widetilde{\cal{T}}^{t}$, 
is a minimal idempotent in ${\cal{C}}_{t}$, ie, 
 $E_{SS}$ cannot be written as a sum of orthogonal idempotents $E_{SS}(1)$ and $E_{SS}(2)$ 
 in ${\cal{C}}_{t}$ where both $E_{SS}(1)$ and $E_{SS}(2)$  are non-zero.
 \end{lemma}  
We will show this by drawing on work of Wenzl  \cite{w2}. 
Before proving Lemma \ref{lem:minimalprojectors(1)},
we define the annihilator ideal $J_{t} \subset {\cal{C}}_{t}$ 
with respect to the quantum supertrace, which we will use in 
the proof of Lemma \ref{lem:minimalprojectors(1)}
and in the proof of the tensor product theorems later in this chapter.
\begin{definition}
Define the  ideal $J_{t} \subset {\cal{C}}_{t}$ by
$J_{t} = \left\{y \in {\cal{C}}_{t} | \ str_{q}(xy)=0, \ \forall x \in {\cal{C}}_{t}\right\}$.
\end{definition}
\noindent
Note that for any $y \in J_{t}$, 
$str_{q}(xy) = (-1)^{[x][y]}str_{q}(yx)=0$ for all $x \in {\cal{C}}_{t}$, 
thus $J_{t}$ is a two-sided ideal.

Recall from Lemma \ref{eq:peacockinblackbeansauce} that for $\hat{q} \neq 0$ not a root of unity,
the algebra homomorphism
$$\Upsilon: g_{i} \mapsto - \check{R}_{i}^{gen}, \hspace{10mm} i=1,  \ldots, t-1,$$
furnishes a representation of  the Birman-Wenzl-Murakami algebra 
$\mathscr{BW}_{t}(-\hat{q}^{2n},\hat{q})$ in ${\cal{C}}_{t}^{gen}$.
In a similar way, we can show that 
for $q= \exp{(2 \pi i/N)}$ where $N \geq 3$ is an integer, the algebra homomorphism
$$\Upsilon: g_{i} \mapsto - \check{R}_{i}, \hspace{10mm} i=1,  \ldots, t-1,$$
furnishes a representation of $\mathscr{BW}_{t}(-q^{2n},q)$ in ${\cal{C}}_{t}$.

We now prove Lemma \ref{lem:minimalprojectors(1)}.
\begin{proof}

In this proof we shall say that an algebra $B$ is semisimple if it is isomorphic to a direct sum of
matrix algebras: $B \cong \bigoplus_{i \in I} M_{b_{i}}(\mathbb{C})$, where
$M_{b_{i}}(\mathbb{C})$ is the algebra of $b_{i} \times b_{i}$ matrices with 
complex entries \cite[Sec. 1]{w2}.

As $q = \exp{(2 \pi i/N)}$, $-q$ is a root of unity and 
 \cite[Thm. 4.4]{w2} states that there is an isomorphism
\begin{equation}
\label{eq:wenzl(1)(1)}
\mathscr{BW}_{t}\big( -q^{2n},q\big)/ \mathscr{J}_{t}\big( -q^{2n},q\big)
\cong
\mathscr{BW}_{t}\big((-q)^{2n},-q\big)/ \mathscr{J}_{t}\big((-q)^{2n},-q\big),
\end{equation}
where right hand side of Eq. (\ref{eq:wenzl(1)(1)}) is semisimple \cite[Thm. 4.4 (d), Thm. 6.4]{w2},
thus  
$$\mathscr{BW}_{t}\big( -q^{2n},q\big)/ \mathscr{J}_{t}\big( -q^{2n},q\big)
\cong \bigoplus_{i \in \Gamma (-q^{2n},q)_{t}} M_{b_{i}}(\mathbb{C}),$$
where $\Gamma (-q^{2n},q)_{t}$ is the set of Young diagrams 
containing $t-2k \geq 0$ boxes where $k \in \mathbb{Z}_{+}$  appearing in 
a certain graph $\Gamma (-q^{2n},q)$  \cite[Thm. 4.4 (d)]{w2}.  
We obtain the graph $\Gamma(-q^{2n},q)$ by following \cite[Thm. 4.4]{w2}.
To obtain $\Gamma(-q^{2n},q)$, we firstly define a subgraph 
$\widetilde{\Gamma}(-q^{2n},q)$ of the Young lattice \cite[p. 407]{w2}.
Firstly, let the Young diagram with no boxes belong to 
$\widetilde{\Gamma}(-q^{2n},q)$, then a Young diagram $\lambda$
belongs to $\widetilde{\Gamma}(-q^{2n},q)$ if
\begin{itemize}
\item[(a)] $Q_{\lambda}(-q^{2n},q) \neq 0$, and
\item[(b)] there is at least one subdiagram of $\lambda$ with $|\lambda|-1$ boxes that belongs to 
$\widetilde{\Gamma}(-q^{2n},q)$.
\end{itemize}

From $\widetilde{\Gamma}(-q^{2n},q)$, we obtain the graph $\Gamma(-q^{2n},q)$ 
using \cite[Thm. 4.4]{w2}, 
the relevant parts of which we now quote. 
In this theorem, the concept of an $N/2$ regular diagram is used.
We say that a Young diagram $\lambda$
is $N/2$ regular if its largest hook has fewer than $N/2$ boxes, ie
$\lambda_{1} + \lambda_{1}'-1 < N/2$ \cite[Eq. (2.5)]{w2}. 

Before stating the relevant parts of \cite[Thm. 4.4]{w2},
we again stress the important fact that 
$\mathscr{BW}_{t}(-q^{2n},q) / \mathscr{J}_{t}(-q^{2n},q)$ is semisimple for all
$t \in \mathbb{N}$.
\begin{theorem}{Theorem 4.4 of \cite{w2}.}
\label{thm:one_of_wenzls)theorems(1)}
Let $q^{2}$ be a primitive $(N/2)^{th}$ root of unity, then
\begin{itemize}
\item[(b)]  Assume that $\widetilde{\Gamma}(-q^{2n},q)$, as defined above, does not contain a hook
diagram with $N/2-1$ boxes.  Then $\Gamma(-q^{2n},q) = \widetilde{\Gamma}(-q^{2n},q)$ with edges
inherited from the Young lattice and it only contains $N/2$ regular diagrams.
\item[(c)] If $\widetilde{\Gamma}(-q^{2n},q)$ contains only one hook diagram $\mu$ with $N/2-1$ boxes
with, say, $\mu = [N/2-2n, 1^{2n-1}]$ and it does not contain its successor
$[N/2-2n, 2, 1^{2n-2}]$, then $\Gamma(-q^{2n},q)$ consists of all $N/2$ regular diagrams in
$\widetilde{\Gamma}(-q^{2n},q)$ and, if $Q_{\lambda}(-q^{2n},q) \neq 0$, also the diagram
$\lambda = [N/2-2n+1, 1^{2n-1}]$.  The edges of $\Gamma(-q^{2n},q)$ are exactly those inherited from the
Young lattice.
\item[(d)] These graphs completely determine the direct sum of matrix rings that is isomorphic to 
$\mathscr{BW}_{t}(-q^{2n},q)/\mathscr{J}_{t}(-q^{2n},q)$.
\end{itemize}
\end{theorem}

The first step is to construct $\widetilde{\Gamma}(-q^{2n},q)$.
Fix the Young diagram with no boxes to belong to $\widetilde{\Gamma}(-q^{2n},q)$.
For convenience, we again state
$Q_{\lambda}(-q^{2n},q)$  from  (\ref{eq:themagixQpolynomial}):
\begin{eqnarray}
Q_{\lambda}(-q^{2n},q) & = & 
\prod_{(j,j) \in \lambda} 
\frac{-q^{2n + \lambda_{j} - \lambda_{j}'}+q^{-2n-\lambda_{j} +\lambda_{j}'}+
q^{\lambda_{j} + \lambda_{j}'-2j+1}-q^{-\lambda_{j} -\lambda_{j}'+2j-1}}
{q^{h(j,j)}-q^{-h(j,j)}} \nonumber \\
& & \hspace{5mm} \times
\prod_{(i,j) \in \lambda, i \neq j} \frac{-q^{2n+d(i,j)}+q^{-2n-d(i,j)}}{q^{h(i,j)}-q^{-h(i,j)}},
\label{eq:themagixQpolynomial(2)}
\end{eqnarray}
where we recall the meanings of $d(i,j)$ and  $h(i,j)$ from Eq. (\ref{eq:themagixQpolynomial}).
As discussed after Eq. (\ref{eq:themagixQpolynomial}),
 $Q_{\lambda}(-q^{2n},q) = 0$ if and only if one (or both) of the following conditions is satisfied:
\begin{itemize}
\item[(a)] $q^{4n + 2d(i,j)}=1$ for some $(i,j) \in \lambda$ where $i \neq j$, 
\item[(b)] $q^{2n-2\lambda_{j}'+2j-1} = 1$ or $q^{2n+2\lambda_{j}-2j+1}=-1$ for some $j$.
\end{itemize}

Fix $q=\exp{(2 \pi i/N)}$ where $N \geq 4$ is an even number.
Note the overarching result
that $Q_{\lambda}(-q^{2n},q)=0$ if $\lambda_{1}' + \lambda_{2}' \geq 2n+2$
\cite[p. 422]{w2}.
We firstly determine when the numerator of $Q_{\lambda}(-q^{2n},q)$ vanishes.
Let us determine for which $\lambda$ the first equation in condition (b) above is true.
Here  $q^{2n-2\lambda_{j}'+2j-1} = 1$  if and only if 
$\lambda_{j}' = rN/2 + n + j -1/2$
for some $r \in \mathbb{Z}$, but this is not possible as $\lambda_{j}'$ must be an integer.
Similarly in relation to the second equation in condition (b) 
we have $q^{2n+2\lambda_{j}-2j+1} \neq 1$ for $N \equiv 0 \pmod{4}$.  
However, for $N \equiv 2 \pmod{4}$ we have $q^{2n+2\lambda_{j}-2j+1}=-1$ if
\begin{equation}
\label{eq:daphneandchloe}
\lambda_{j} = N/4 + rN/2 - n + j -1/2,  \hspace{10mm} \mbox{for any } r \in \mathbb{Z},
\end{equation}
where we note that the right hand side of (\ref{eq:daphneandchloe}) must be an integer.
Recall that we have fixed $N$ to satisfy the inequality  
$N/4 \geq n+1/2$ so that  $\Lambda_{N}^{+} \neq \emptyset$, 
then $N/4 - n + j -1/2 \geq j$. 
Also, we have $-N/4 -n+j-1/2 \leq -2n+j-1 \leq 0$, so the least non-negative value on
the right hand side of (\ref{eq:daphneandchloe}) is obtained by fixing $r=0$.
Now $\max{\{\lambda_{j}\}} = \lambda_{1}$, and so the numerator of 
$Q_{\lambda}(-q^{2n},q)$ is zero if
\begin{equation}
\label{eq:dukeellington(10)}
\lambda_{1} = N/4 - n + 1/2.
\end{equation}

We claim that $Q_{\lambda}(-q^{2n},q) \neq 0$ for $N \equiv 2 \pmod{4}$ if $\lambda$ satisfies
\begin{equation}
\label{eq:dukeellington(1)}
\lambda_{1} \leq N/4 - n - 1/2.
\end{equation}
It may be observed that (\ref{eq:daphneandchloe}) is never satisfied
by any $\lambda$ satisfying (\ref{eq:dukeellington(1)}).
Furthermore, any $\lambda$ satisfying (\ref{eq:dukeellington(1)}) and
$\lambda'_{1} + \lambda_{2}' \leq 2n+1$ also satisfies the hook length condition that
$h(1,1) \leq N/2-1$ as 
$$\max{(h(1,1))} = \left\{ 
\begin{array}{ll}
2n+1 \leq N/2-2, & \mbox{if } \lambda_{1} = 1, \\
N/4+n-3/2 \leq N/2-4, & \mbox{if } \lambda_{1} > 1.
\end{array} \right. $$
The statement for $\lambda_{1}=1$ comes about from the condition in (\ref{eq:dukeellington(10)})
that $1 \leq N/4-n-1/2$. 
The statement for $\lambda_{1} > 1$ comes about from the condition in (\ref{eq:dukeellington(10)}) 
that $n - 3/2 \leq N/4 - \lambda_{1} - 2 \leq N/4-4$.

Now we will show that for $N \equiv 2 \pmod{4}$, condition (a) 
mentioned just after Eq. (\ref{eq:themagixQpolynomial(2)}) is never satisfied
by any Young diagram $\lambda$ satisfying (\ref{eq:dukeellington(1)}) and
$\lambda'_{1} + \lambda_{2}' \leq 2n+1$.
Fix $N \equiv 2 \pmod{4}$.
Condition (a) is true if and only if $d(i,j) = rN/2-2n$ for some $r \in \mathbb{Z}$ where
$i \neq j$. Recall that
$$
d(i,j) = \left\{
\begin{array}{ll}
\lambda_{i} + \lambda_{j}-i-j+1, & \mbox{if } i \leq j, \\
-\lambda_{i}'-\lambda_{j}'+i+j-1, & \mbox{if } i > j, \end{array} \right. 
$$
then for all $i \neq j$, 
\begin{itemize}
\item  $\max{(d(i,j))} = d(1,2) = \lambda_{1} + \lambda_{2} -2$, and
\item $\min{(d(i,j))} = d(2,1) = -\lambda_{1}'-\lambda_{2}'+2$.
\end{itemize}
It follows that $q^{4n + 2d(i,j)}=1$ for some $i \neq j$ if and only if 
\begin{eqnarray}
\lambda_{i} + \lambda_{j} & = & rN/2-2n+i+j-1,  \hspace{10mm}  \mbox{if } i<j, \label{eq:arunram99} \\
\lambda_{i}' + \lambda_{j}' & = & rN/2+2n+i+j-1, \hspace{10mm}  \mbox{if } i>j, \label{eq:arunram100}
\end{eqnarray}
for some $r \in \mathbb{Z}$.  

Let us consider the case where $i<j$.
From the previous condition in (\ref{eq:dukeellington(1)}) that $\lambda_{1} \leq N/4-n-1/2$, we have
\begin{equation}
\label{eq:hellosailor!}
0 \leq \lambda_{i} + \lambda_{j} \leq N/2-2n-1, \hspace{5mm} \forall i<j.
\end{equation}
Note that for these $i$ and $j$, we have $i \leq n$ and $j \leq N/4-n-1/2$.
The right hand side of (\ref{eq:arunram99}) with $r=-1$ is
$-N/2-2n+i+j-1 \leq -N/4 - 2n - 3/2 <0$, and thus 
$\lambda_{i} + \lambda_{j} > -N/2-2n+i+j-1$.
The right hand side of (\ref{eq:arunram99}) with $r=1$ is
$N/2-2n+i+j-1 \geq N/2-2n+2$, but we always have
$ \lambda_{i} + \lambda_{j} < N/2-2n+2$ from (\ref{eq:hellosailor!}) and so
$ \lambda_{i} + \lambda_{j} < N/2-2n+i+j-1$.
Eq. (\ref{eq:arunram99}) with $r=0$ is
$\lambda_{i} + \lambda_{j} -i-j+1= -2n$, and as at generic $q$, this is not true.
(Recall the argument from generic $q$: we have $i \leq 2n$, 
$\lambda_{i}-j \geq 0$ and $\lambda_{j} \geq 0$, thus
$\lambda_{i} + \lambda_{j} -i-j+1 \geq -2n+1$.)
It follows that (\ref{eq:arunram99}) is not true at these roots of unity.

Let us consider the case where $i>j$ with $N$ again fixed as $N \equiv 2 \pmod{4}$.
We have already imposed the condition $\lambda_{1}' + \lambda_{2}' \leq 2n+1$,
so $\lambda_{i}' + \lambda_{j}' \leq 2n$ for all $(i,j) \neq (2,1)$.
As $\lambda_{1}' + \lambda_{2}' \leq 2n+1$,
a necessary condition on the integer $r$ for 
$\lambda_{1}' + \lambda_{2}' = rN/2+2n+2$ to be true 
is that $r \leq -1$.
Setting $r=-2$ into this equation gives $\lambda_{1}' + \lambda_{2}' = -N+2n+2$,
however, it is true that $-N+2n+2 \leq -2n$ as $N$ satisfies $N \geq 4n+2$, and thus
$\lambda_{1}' + \lambda_{2}' > -N+2n+2$. 
Fixing then $r = -1$, the equation is $\lambda_{1}' + \lambda_{2}' = -N/2+2n+2$,
and from the conditions on $N$ we have $\lambda_{1}' + \lambda_{2}' = -N/2+2n+2 \leq 1$.
Now $-N/2+2n+2$ is odd so this condition is  
$\lambda_{1}' + \lambda_{2}' = -N/2+2n+2 = 1$ which is precisely the Young diagram
$[1]$. But this Young diagram has already been ruled out for $N/2 = 2n+1$ as we have the condition
$\lambda_{1} \leq N/4-n-1/2 = 0$. 
Thus (\ref{eq:arunram100}) is not true for $(i,j) = (2,1)$ at these roots of unity.

Now consider (\ref{eq:arunram100}) for $i > j$ where $(i,j) \neq (2,1)$: 
recall that $0 \leq \lambda_{i}' + \lambda_{j}' \leq 2n$.
Here $2n+2 \geq i+j \geq 4$, so $rN/2+2n+i+j-1 \geq rN/2+2n+3$, 
and then a necessary condition on $r$ for Eq. (\ref{eq:arunram100}) to be true
is that $r \leq -1$.
For $r \leq -2$ we have
$$rN/2+2n+i+j-1 \leq -N+2n+i+j-1 \leq -N + 4n +1 \leq -1,$$
and so a necessary condition on $r$ for Eq. (\ref{eq:arunram100}) to be true
is that $r = -1$.

We now show that (\ref{eq:arunram100}) is not true with $r=-1$
for all $i > j$ where $(i,j) \neq (2,1)$ by considering all the possible relevant Young diagrams.
Consider the Young diagram $[1^{k}]$ where $k \in \{3, 4, \ldots, 2n+1 \}$.
Here $(i,j) = (i,1)$ where $i \in \{3, 4, \ldots, k \}$ and $\lambda_{i}'=0$
and $\lambda_{1}' = k$.
Eq. (\ref{eq:arunram100}) is true here only if 
\begin{equation}
\label{pink1}
\lambda_{1}'-i = -N/2+2n.
\end{equation} 
However, $\lambda_{1}'-i \geq 0$ and $N/2 \geq 2n+3$, thus the left hand
of (\ref{pink1}) is non-negative and the right hand side is strictly negative, thus
(\ref{pink1}) is not true.
Consider now a Young diagram $\lambda$ with more than one column of boxes.
Fix $i \geq n+1$, then $j=1$ (note that $i \leq 2n$).
Then (\ref{eq:arunram100}) is true here only if 
\begin{equation}
\label{pink2}
\lambda_{i}' + \lambda_{1}'-i = -N/2+2n.
\end{equation}
Now the left hand side of (\ref{pink2}) is non-negative as $\lambda_{i}' \geq 0$ 
and $\lambda_{1}'-i \geq 0$, but $N/2 \geq 2n+3$, so the right hand side is strictly negative,
thus (\ref{pink2}) is not true.
Consider again a Young diagram $\lambda$ with more than one column of boxes.
Fix $i \leq n$, then $j \leq i-1$, and (\ref{eq:arunram100}) is true only if 
\begin{equation}
\label{pink3}
\lambda_{i}' + \lambda_{j}'-i = -N/2+2n+j-1.
\end{equation}
Now $\lambda_{i}' \geq 0$ and $\lambda_{j}'-i \geq 0$ so the left hand side of (\ref{pink3})
is non-negative.
Note that $$-N/2+2n+j-1 \leq -N/4+n-3/2 \leq -3,$$
as $j \leq N/4-n-1/2$ and $N/2 \geq 2n+3$, so the right hand side of (\ref{pink3}) is strictly negative,
thus (\ref{pink3}) is not true.

It follows from these calculations that for  $N \equiv 2 \pmod{4}$,
$Q_{\lambda}(-q^{2n},q) \neq 0$ for all
the Young diagrams $\lambda$ satisfying
\begin{equation}
\label{tangerine1}
\lambda_{1} \leq N/4-n-1/2, \hspace{5mm} \mbox {and} \hspace{5mm} \lambda_{1}' + \lambda_{2} ' \leq 2n+1.
\end{equation}
We can then write down the Young diagrams that comprise the vertices of 
$\widetilde{\Gamma}(-q^{2n},q)$ at these roots of unity; these are all the Young diagrams
$\lambda$ satisfying (\ref{tangerine1}).
Note that $\widetilde{\Gamma}(-q^{2n},q)$ contains no hook diagrams with exactly
$N/2-1$ boxes.

Now we determine the vertices of $\widetilde{\Gamma}(-q^{2n},q)$ when
$N$ satisfies $N \equiv 0 \pmod{4}$; fix such an $N$. 
As previously mentioned, neither of the
equations in condition (b) is true at these roots of unity; 
we now determine whether the equation in condition
(a) is true at these roots of unity.
This condition can be expressed as the two equations (\ref{eq:arunram99})--(\ref{eq:arunram100}).
Recall the overarching result that
$Q_{\lambda}(-q^{2n},q)=0$ if $\lambda_{1}'+\lambda_{2}' \geq 2n+2$.
Let us consider Eqs. (\ref{eq:arunram99})--(\ref{eq:arunram100}) for $i<j$.

Firstly,  fix $(i, j) = (1,2)$, then Eq. (\ref{eq:arunram99}) is
\begin{equation}
\label{green1}
\lambda_{1} + \lambda_{2} = rN/2-2n+2,
\end{equation}
for some $r \in \mathbb{Z}$.
Now (\ref{green1}) can only be true for Young diagrams containing at least one box if $r \geq 1$, and
we have $Q_{\lambda}(-q^{2n},q)=0$ if $\lambda_{1} + \lambda_{2} = N/2-2n+2$.
We claim that Eqs. (\ref{eq:arunram99})--(\ref{eq:arunram100}) 
are not true for all Young diagrams $\lambda$
satisfying 
\begin{equation}
\label{tangerine2}
\lambda_{1} + \lambda_{2} \leq N/2-2n+1, \hspace{5mm} \mbox {and} \hspace{5mm} 
\lambda_{1}' + \lambda_{2}' \leq 2n+1.
\end{equation}
The proof of this claim is very similar to the proof of the corresponding result for 
$N \equiv 2 \pmod{4}$ for all Young diagrams satisfying (\ref{tangerine1}), thus we omit it.

Then for $N \equiv 0 \pmod{4}$, the vertices of the graph 
$\widetilde{\Gamma}(-q^{2n},q)$ are all the Young diagrams $\lambda$ satisfying
(\ref{tangerine2}).
We need to check that all such $\lambda$  also satisfy the hook length condition.
However this is easy to do and we omit the proof.

Now we wish to determine the number of vertices of $\widetilde{\Gamma}(-q^{2n},q)$
for $N \equiv 0 \pmod{4}$ that are hook diagrams with exactly $N/2-1$ boxes.
As $N \geq 4$, we are examining the Young diagrams with at least one box.
For all of these Young diagrams, $\lambda_{1} \geq 1$ and thus
$N/2 \geq 2n$.

Consider the case that $N/2 = 2n$.  Here 
$\lambda_{1} + \lambda_{2} \leq 1$ which means that the only non-empty Young diagram 
in $\widetilde{\Gamma}(-q^{2n},q)$ is
$[1]$. This is a hook diagram, but it contains exactly $N/2-1$ boxes only if $N=4$ and $n=1$.

Consider the case that $N/2 = 2n+2$. Here $\lambda_{1} + \lambda_{2} \leq 3$,
and the hook diagrams in $\widetilde{\Gamma}(-q^{2n},q)$ with exactly $N/2-1$ boxes are
\begin{itemize}
\item $[3]$, $[2, 1]$ and $[1^{3}]$, if $n=1$,
\item $[1^{2n+1}]$ and $[2,1^{2n-1}]$, if $n \geq 2$. 
\end{itemize}

Consider the case that $N/2 \geq 2n+4$.
The hook diagrams in $\widetilde{\Gamma}(-q^{2n},q)$ with exactly $N/2-1$ boxes are
\begin{itemize}
\item $[N/2-1]$ and $[N/2-2,1]$, if $n=1$,
\item $[N/2-2n,1^{2n-1}]$, if $n \geq 2$.
\end{itemize}

We can now work out the vertices belonging to the graph $\Gamma(-q^{2n},q)$ using 
Theorem \ref{thm:one_of_wenzls)theorems(1)}.
For $N \equiv 2 \pmod{4}$, $\widetilde{\Gamma}(-q^{2n},q)$  contains no hook 
diagrams with exactly $N/2-1$ boxes,  
and so the vertices of $\Gamma(-q^{2n},q)$ are all the Young diagrams satisfying
$$\lambda_{1} \leq N/4-n-1/2, \hspace{5mm} \mbox{and} \hspace{5mm} 
\lambda_{1}'+ \lambda_{2}' \leq 2n+1.$$

Now for $N \equiv 0 \pmod{4}$, we have the following cases.
For $N/2 = 2n$, the Young diagram $[1]$ is a hook diagram with $N/2-1$ boxes if $N=4$ and $n=1$.
However, in this case $Q_{[2]}(-q^{2n},q)=0$ and so the hook diagram $[2]$ is not a vertex in the 
graph $\Gamma(-q^{2n},q)$. Then, for $N/2 = 2n$ and $N/2 = 2n+2$, 
and also $N/2 \geq 2n+4$ where $n = 1$,
the vertices of $\Gamma(-q^{2n},q)$ are all the Young diagrams satisfying
$$\lambda_{1} + \lambda_{2} \leq N/2-2n+1, \hspace{5mm} \mbox {and} \hspace{5mm} 
\lambda_{1}' + \lambda_{2}' \leq 2n+1.$$

For $N/2 \geq 2n+4$ where $n \geq 2$, 
$\widetilde{\Gamma}(-q^{2n},q)$ contains one hook diagram 
with $N/2-1$ boxes: $[N/2-2n,1^{2n-1}]$ and it does not contain its successor
$[N/2-2n,2,1^{2n-2}]$. In addition, 
the $OSp(1|2n)$ supercharacters of $[N/2-2n+1,1^{2n-1}]$ and
$[N/2-2n+1]$ are the same up to a sign, so $Q_{[N/2-2n+1,1^{2n-1}]}(-q^{2n},q) \neq 0$. 
The vertices of
$\Gamma(-q^{2n},q)$ are then all the Young diagrams $\lambda$ in the set
$$\{ \lambda| \ \lambda_{1} + \lambda_{2} \leq N/2-2n+1, \ \lambda_{1}'+ \lambda_{2}' \leq 2n+1 \}
 \cup \{[N/2-2n+1,1^{2n-1}]\}.$$

Now we can write down the Bratteli diagram for $\mathscr{BW}_{t}(-q^{2n},q)$ as we did at generic $q$,
and then write down the matrix units in $\mathscr{BW}_{t}(-q^{2n},q) / \mathscr{J}_{t}(-q^{2n},q)$.
For each $s=0, 1, \ldots, t$, let $\Gamma(-q^{2n},q)_{s}$ be the set of Young diagrams 
belonging to $\Gamma(-q^{2n},q)$ with 
\newline
$s-2k \geq 0$ boxes, where $k$ ranges over all of $\mathbb{Z}_{+}$.
We say that $R$ is a path of length $t$ if $R$ is a sequence of $t+1$ Young diagrams:
$R=([0], [1], r_{2}, \ldots, r_{t})$ where $r_{s} \in \Gamma(-q^{2n},q)_{s}$ for each $s$
and $r_{j}$ is connected to $r_{j+1}$ for each $j=0, 1, \ldots, t-1$.
Let $\omega(-q^{2n},q)_{t}$ be the set of pairs $(R,S)$ of paths of length $t$ such that
$r_{t} = s_{t}$.

We obtain a complete set of matrix units for 
$\mathscr{BW}_{t}(-q^{2n},q) / \mathscr{J}_{t}(-q^{2n},q)$
by taking all the matrix units  $e_{RS} \in \mathscr{BW}_{t}$ where 
$(R,S) \in \omega(-q^{2n},q)_{t}$ and fixing  $q=\exp{(2 \pi i/N)}$ and $r = -q^{2n}$.
These are all well-defined and non-zero.

Let us define $\widetilde{\mathscr{BW}}_{t}(-q^{2n},q)$ to be the semisimple subalgebra 
of $\mathscr{BW}_{t}(-q^{2n},q)$ spanned by the 
matrix units $\{ e_{ST} \in \mathscr{BW}_{t}(-q^{2n},q) | \ (S,T) \in  \omega(-q^{2n},q)_{t} \}$.
Note that $\Upsilon(e_{ST}) = E_{ST} \in {\cal{C}}_{t}$ 
for each  $(S,T) \in \omega(-q^{2n},q)_{t}$, where we recall that
$\Upsilon$ is the algebra homomorphism 
$\Upsilon: g_{i} \mapsto -\check{R}_{i}$.

Now define a map $\psi: {\cal{C}}_{t} \rightarrow \mathbb{C}$ by
$$\psi(X) = str_{q}(X)/\big(sdim_{q}(V)\big)^{t},$$
then 
\begin{equation}
\label{eq:eenieweenie(1)}
\psi\big(\Upsilon(a)\big) = \mathrm{tr}(a), \hspace{10mm} \mbox{for all }
a \in \mathscr{BW}_{t}(-q^{2n},q),
\end{equation}
which we can prove in the same way that we proved Lemma \ref{lemma:equalityoftraces}.
Note here that Lemma \ref{lem;kilo} holds true if we write $U_{q}^{(N)}(\mathfrak{g})$ instead of
 $U_{q}(\mathfrak{g})$ and fix 
 $\check{\cal{R}}_{V,V} \in End_{U^{(N)}_{q}(\mathfrak{g})}(V \otimes V)$ to be as given in
Definition \ref{defn:bigjim(2)}.
Furthermore, note that $\psi(X)=0$ if and only if $str_{q}(X)=0$, thus
we can regard $J_{t}$ as the annihilator ideal of ${\cal{C}}_{t}$ with respect to $\psi$.

Note that
$\widetilde{\mathscr{BW}}_{t}(-q^{2n},q) \cap \mathscr{J}_{t}(-q^{2n},q) = 0$
for the following reasons. 
Any element $f \in \widetilde{\mathscr{BW}}_{t}(-q^{2n},q)$ 
is a linear combination of the matrix units:
$$f = \sum_{(S,T) \in \omega(-q^{2n},q)_{t}}f_{ST} e_{ST}, \hspace{5mm} f_{ST} \in \mathbb{C},$$ where  
$f_{ST} \neq 0$ for at least one pair $(S,T)$.  Fix $(A,B)$ to be such a pair, then 
Eq. (\ref{eq:eenieweenie(1)}) implies that
$$\mathrm{tr}(e_{BA}f) = str_{q}\big(f_{AB}E_{BA}E_{AB}\big)/\big(sdim_{q}(V)\big)^{t} = 
                  str_{q}\big(f_{AB}E_{BB}\big)/\big(sdim_{q}(V)\big)^{t} \neq 0,$$
as $str_{q}(E_{BB}) \neq 0$ for all $(B,B) \in \omega(-q^{2n},q)_{t}$.
Thus any non-zero $f$ belonging to $\widetilde{\mathscr{BW}}_{t}(-q^{2n},q)$ 
does not belong to $\mathscr{J}_{t}(-q^{2n},q)$, giving the direct sum decomposition 
\begin{equation}
\label{eq:directsumdecompofbirman}
\mathscr{BW}_{t}(-q^{2n},q) = \widetilde{\mathscr{BW}}_{t}(-q^{2n},q)
\oplus  \mathscr{J}_{t}(-q^{2n},q).
\end{equation}
Then for all $a \in \mathscr{BW}_{t}(-q^{2n},q)$, we have
$a = \widetilde{a} + a_{j}$, where 
$\widetilde{a} \in \widetilde{\mathscr{BW}}_{t}(-q^{2n},q)$ and
$a_{j} \in \mathscr{J}_{t}(-q^{2n},q)$.

Let us define 
$$\mathscr{P}_{t} = 
\sum_{(T,T) \in \omega(-q^{2n},q)_{t}}e_{TT} \in \mathscr{BW}_{t}(-q^{2n},q),$$
then we have $\mathscr{P}_{t} a \mathscr{P}_{t} = \widetilde{a}$ for all
$a \in \mathscr{BW}_{t}(-q^{2n},q)$, which can be seen 
by regarding $\mathscr{BW}_{t}(-q^{2n},q)$ as a matrix algebra.

Using Eq. (\ref{eq:eenieweenie(1)}), we can prove that
\begin{equation}
\label{eq:shonauy(1)}
\Upsilon\big(\mathscr{J}_{t}(-q^{2n},q)\big) = J_{t},
\end{equation} 
in the same way that we proved Eq. (\ref{eq:shonauy(99)}).
The surjectivity of $\Upsilon$ implies that 
$${\cal{C}}_{t} = \Upsilon\big(\widetilde{\mathscr{BW}}_{t}(-q^{2n},q) \big) + J_{t},$$
and we will show that this sum is direct.
To see this, assume that there exists some non-zero element $F$ of ${\cal{C}}_{t}$ belonging to
$\Upsilon\big(\widetilde{\mathscr{BW}}_{t}(-q^{2n},q) \big)$ and also to $J_{t}$, then
$str_{q}(XF) = 0$ for all $X \in {\cal{C}}_{t}$.  
However, $F$ is the image of a linear combination of matrix units: 
$F = \sum_{(S,T) \in \omega(-q^{2n},q)_{t}}f_{ST}\Upsilon(e_{ST})$, 
where $f_{ST} \in \mathbb{C}$ and $f_{ST}$ is non-zero for at least one pair $(S,T)$.  
Assume that $(A,B)$ is such a pair, then by similar reasoning as
previously, $str_{q}(\Upsilon(e_{BA})F) \neq 0$  contradicting the assumption that $F \in J_{t}$.
Thus
$\Upsilon\big(\widetilde{\mathscr{BW}}_{t}(-q^{2n},q) \big) \cap J_{t} = 0$,
and  we have
$${\cal{C}}_{t} = \Upsilon\big(\widetilde{\mathscr{BW}}_{t}(-q^{2n},q) \big) \oplus J_{t}.$$

We note the important fact that each element of
$\Upsilon\big(\widetilde{\mathscr{BW}}_{t}(-q^{2n},q) \big)$ is a linear combination
of the matrix units $\{E_{ST} \in {\cal{C}}_{t} | \ (S,T) \in \widetilde{\Omega}^{t} \}$ we defined
previously.  
Note that the two sets $\omega(-q^{2n},q)_{t}$ and $\widetilde{\Omega}^{t}$ are identical
if we apply the map $\lambda \mapsto \widetilde{\lambda}$ to the Young diagrams in the relevant paths
if appropriate, where
$\widetilde{\lambda}'_{1} = 2n+1-\lambda_{1}$ and 
$\widetilde{\lambda}'_{j} = \lambda'_{j}$ for $j \geq 2$.
As $E_{BB} = \Upsilon(e_{BB})$ for each $(B,B) \in \widetilde{\Omega}^{t}$ 
(where we think of the vertices on each path appropriately), this completes the proof.

\end{proof}

After showing that the projectors in ${\cal{C}}_{t}$ are minimal idempotents, we move back to the main
track of our argument and define $P_{t}$, which will be an extremely useful element of ${\cal{C}}_{t}$.
\begin{lemma}
\label{lem:PP}
For each $t \in \mathbb{N}$, define
$$\displaystyle{ P_{t} = \sum_{T \in \widetilde{\cal{T}}^{t}}E_{TT}} \in {\cal{C}}_{t},$$ 
then $P_{t}$ satisfies
\begin{itemize}
\item[(i)] $(P_{t})^{2} = P_{t}$,
\item[(ii)] $str_{q}(P_{t}) = sdim_{q}\left(V^{\otimes t}\right)$,
\item[(iii)] $str_{q}(1-P_{t})=0$.
\end{itemize}
\end{lemma}
Note that $P_{t}(V^{\otimes t}) = \sum_{T \in \widetilde{\cal{T}}^{t}}E_{TT}(V^{\otimes t})
\cong \bigoplus_{T \in \widetilde{\cal{T}}^{t}} V_{shp(T)}$.
\begin{proof}
The proof of (i) follows from the fact that the set of projections
$\big\{ E_{TT} \in {\cal{C}}_{t} | \ T \in \widetilde{\cal{T}}^{t} \big\}$ 
is a set of mutually orthogonal idempotents.  
We now inductively show that $str_{q}(P_{t}) = sdim_{q}\left(V^{\otimes t}\right)$.

Firstly $P_{1} = \mathrm{id}_{V}$ and $str_{q}(P_{1}) = sdim_{q}(V)$.
Secondly $V \otimes V = V_{2 \epsilon_{1}} \oplus V_{\epsilon_{1} + \epsilon_{2}} \oplus V_{0}$, where each of
$2 \epsilon_{1}, \epsilon_{1} + \epsilon_{2}, 0$ are elements of $\overline{\cal{P}}^{+}_{N}$; obviously  
$sdim_{q}(V \otimes V) = \sum_{\mu \in {\cal{P}}^{+}_{\epsilon_{1}} } sdim_{q}(V_{\mu})$
where ${\cal{P}}^{+}_{\epsilon_{1}} = \{ 2 \epsilon_{1}, \epsilon_{1} + \epsilon_{2}, 0 \}$.
(Note that for $n=1$ we write $\epsilon_{1} + \epsilon_{2}$ to mean $\epsilon_{1}$.)
Now if each element of ${\cal{P}}^{+}_{\epsilon_{1}}$ is an element of $\Lambda_{N}^{+}$, then 
$P_{2} = \mathrm{id}_{V \otimes V}$ and $str_{q}(P_{2}) = sdim_{q}(V^{\otimes 2})$.  
Now any element of ${\cal{P}}^{+}_{\epsilon_{1}}$ that is not in $\Lambda_{N}^{+}$ is an element of 
$\overline{\cal{P}}^{+}_{N} \backslash \Lambda_{N}^{+}$; denote any such element by $\mu$, then $sdim_{q}(V_{\mu}) = 0$.
Now
$$P_{2} = \sum_{\lambda \in {\cal{P}}^{+}_{\epsilon_{1}} \cap \Lambda_{N}^{+} } 
\tilde{p}^{2}_{i} [\lambda],$$
and if there is at least one element of ${\cal{P}}^{+}_{\epsilon_{1}}$ 
that is not in $\Lambda_{N}^{+}$, then $P_{2} \neq \mathrm{id}_{V \otimes V}$.  
Let ${\cal{S}}_{\epsilon_{1}}$ be the subset of 
$\{ 2 \epsilon_{1}, \epsilon_{1} + \epsilon_{2}, 0 \}$ 
consisting of those elements not in $\Lambda_{N}^{+}$, then
$str_{q}\big(\tilde{p}^{2}_{i} [\lambda]\big)=0$ for each $\lambda \in {\cal{S}}_{\epsilon_{1}}$.
As
$$
V \otimes V = V_{2 \epsilon_{1}} \oplus V_{\epsilon_{1} + \epsilon_{2}} \oplus V_{0} 
            = \left(P_{2} + \sum_{\lambda \in {\cal{S}}_{\epsilon_{1}}} \tilde{p}^{2}_{i} [\lambda]
	    \right) V \otimes V,
$$
we have
$$
str_{q}(P_{2}) = str_{q} \left( P_{2} + \sum_{\lambda \in {\cal{S}}_{\epsilon_{1}}} \tilde{p}^{2}_{i} [\lambda] \right) 
	       = sdim_{q} \big(V \otimes V \big).
$$

Now we do the inductive step.  Assume that $str_{q}(P_{j}) = sdim_{q}(V^{\otimes j})$ for some $j \geq 2$,
and let $V_{\mu} \subseteq V^{\otimes j}$, $\mu \in \Lambda_{N}^{+}$, be a $U_{q}^{(N)}(\mathfrak{g})$-module defined by
$V_{\mu} = E_{SS} \big(V^{\otimes j}\big)$ where $S \in \widetilde{\cal{T}}^{j}$ and 
$E_{SS}$ is a path projection of length $j$.
As $S \in \widetilde{\cal{T}}^{j}$, we have $sdim_{q}(V_{\mu}) \neq 0$.

Let us define a useful set: 
\begin{equation}
\label{eq:natalyakari(1)} 
{\cal{S}}_{\mu} = 
\left\{ \lambda \in {\cal{P}}^{+}_{\mu} \cap \overline{\Lambda_{N}^{+}} \big| 
\ \lambda \notin \Lambda_{N}^{+} \right\}.
\end{equation}
Recall from Lemma \ref{lemma:heresyetanothermoviewlemmayouknowwhatImean} the following result:
$V_{\mu} \otimes V =
 \bigoplus_{\lambda \in {\cal{P}}^{+}_{\mu} \cap \overline{\Lambda_{N}^{+}} } V_{\lambda}$
 where we have the important facts that  
$sdim_{q}(V_{\lambda}) \neq 0$ if $\lambda \in \Lambda_{N}^{+}$ and
$sdim_{q}(V_{\lambda}) = 0 $ if $\lambda \in \overline{\cal{P}}^{+}_{N} \backslash \Lambda_{N}^{+}$.
Then 
\begin{eqnarray*}
sdim_{q}(V_{\mu} \otimes V) & = & 
\sum_{\lambda \in {\cal{P}}^{+}_{\mu} \cap \overline{\Lambda_{N}^{+}} } 
str_{q}\big( E_{S \circ \lambda,S \circ \lambda }  \big) \\
& = & \sum_{\xi \in  {\cal{P}}^{+}_{\mu} \cap \Lambda_{N}^{+}} str_{q}\big(E_{S \circ \xi, S \circ \xi}\big)
      + \sum_{\zeta \in {\cal{S}}_{\mu} }str_{q}\big(E_{S \circ \zeta, S \circ \zeta} \big) \\
& = & \sum_{\xi \in  {\cal{P}}^{+}_{\mu} \cap \Lambda_{N}^{+}} str_{q}\big(E_{S \circ \xi, S \circ \xi}\big),
\end{eqnarray*}
as $sdim_{q}(V_{\zeta})=0$ for all $\zeta \in {\cal{S}}_{\mu}$.

Now this result is true for all $U_{q}^{(N)}(\mathfrak{g})$-modules 
$V_{\mu'} \subseteq V^{\otimes j}$ defined by
$V_{\mu'} = E_{RR} \big(V^{\otimes j}\big)$ where $R \in \widetilde{\cal{T}}^{j}$
and $E_{RR}$ is a path projection of length $j$, that is
\begin{equation}
\label{eq:bigjim(1)}
sdim_{q}\big( E_{RR} V^{\otimes j} \otimes V \big) = 
\sum_{ \xi \in {\cal{P}}^{+}_{shp(R)} \cap \Lambda_{N}^{+} } 
str_{q} \big( E_{R \circ \xi, R \circ \xi} \big).
\end{equation}
Summing over all distinct $R \in \widetilde{\cal{T}}^{j}$ 
on both sides of Eq. (\ref{eq:bigjim(1)}) gives
\begin{equation}
\label{eq:johnfaulkner(50)}
sdim_{q} \big( P_{j} V^{\otimes j} \otimes V \big) = str_{q}\big(P_{j+1}\big),
\end{equation}
as $$P_{j+1} = \sum_{Q \in \widetilde{\cal{T}}^{(j+1)} } E_{QQ} = 
\sum_{R \in \widetilde{\cal{T}}^{j}}  \left( \sum_{ \xi \in {\cal{P}}^{+}_{shp(R)} \cap \Lambda_{N}^{+}} 
E_{R \circ \xi, R \circ \xi} \right),$$
and the left hand side of Eq. (\ref{eq:johnfaulkner(50)}) is $sdim_{q} \big( V^{\otimes (j+1)} \big)$ 
by assumption.
This completes the induction and the proof of (ii), and the proof of (iii) then follows.

\end{proof}

We now consider a proposition that will be extremely useful in proving the tensor product theorems in
this chapter:
\begin{proposition}
\label{lem:concorde}
\noindent
\begin{itemize}
\item[(i)] $(1-P_{t})$ generates $J_{t}$ as a two-sided ideal in ${\cal{C}}_{t}$,
\item[(ii)] the mapping ${\cal{C}}_{t} \rightarrow P_{t} {\cal{C}}_{t}P_{t}$ defined by 
$a \mapsto P_{t}aP_{t}$ is an algebra homomorphism.
\end{itemize}
\end{proposition}

We now present two technical lemmas
 that we will use, in addition to Lemma \ref{lem;kilo}, in the proof of Proposition \ref{lem:concorde}.
\begin{lemma}
\label{lem:tomtheguineapig}
Each element in ${\cal{C}}_{t}$ can be written as a linear combination of elements 
$(a \otimes \mathrm{id})$ and $(a \otimes \mathrm{id}) \check{R}_{t-1}^{\pm 1}(b \otimes \mathrm{id})$ 
where $a, b \in {\cal{C}}_{t-1}$.
\end{lemma}
\begin{proof}
From \cite[Lem. 3.1]{bw},
each element of the Birman-Wenzl-Murakami algebra $\mathscr{BW}_{t}(z,q)$ 
can be written as a linear combination of elements $a \gamma b$ 
where $\gamma \in \{g_{t-1}, e_{t-1}, 1\}$ and $a$, $b$ are elements of $\mathscr{BW}_{t-1}(z,q)$.
We complete the proof by applying the algebra homomorphism
$\Upsilon: \mathscr{BW}_{t}(-q^{2n},q) \rightarrow {\cal{C}}_{t}$
to the appropriate equations in $\mathscr{BW}_{t}(-q^{2n},q)$.
\end{proof}
\noindent
\begin{lemma}
\label{lem:dalailama}
Let $B \in \widehat{\cal{T}}^{t} \backslash \widetilde{\cal{T}}^{t}$ and
$E_{BB} = \tilde{p}^{t}_{i}[\lambda]$, where 
$\tilde{p}^{t}_{i}[\lambda]: V^{\otimes t} \rightarrow V_{\lambda}$ is a well-defined projection.  Then
$E_{BB} \in J_{t}$.
\end{lemma}
\begin{proof}
Let $B \in \widehat{\cal{T}}^{t} \backslash \widetilde{\cal{T}}^{t}$, 
then $shp(B) = \lambda \in \overline{\cal{P}}^{+}_{N} \backslash \Lambda_{N}^{+}$ and
$sdim_{q}(V_{\lambda})=0$, which implies that $str_{q}(E_{BB})=0$.
Now let $f \in {\cal{C}}_{t}$ be arbitrary, then 
$$str_{q}(fE_{BB})=str_{q}(E_{BB}fE_{BB})=\beta str_{q}(E_{BB})=0,$$
for some complex constant $\beta$.  
\end{proof}
\noindent
We now prove Proposition \ref{lem:concorde}.
\begin{proof}
We prove (i).
We will firstly show that $P_{t} {\cal{C}}_{t} P_{t} \cap J_{t} = 0$.
Assume that this is not true, that is there exists some $f \in {\cal{C}}_{t}$ such that 
$P_{t} f P_{t} \neq 0$ and $P_{t} f P_{t} \in J_{t}$,  then 
$$P_{t} f P_{t} = \sum_{S, T \in \widetilde{\cal{T}}^{t}} E_{SS} f E_{TT} = 
\sum_{(S, T) \in \widetilde{\Omega}^{t}} f_{ST} E_{ST}, \hspace{10mm} f_{ST} \in \mathbb{C}.$$ 
The fact that
$P_{t} f P_{t} \neq 0$ implies that $f_{ST} \neq 0$ for at least one pair $(S, T) \in \widetilde{\Omega}^{t}$. 
Fix $(S, T)$ to be such a pair, then
\begin{eqnarray*}
str_{q}\big(E_{TS} P_{t} f P_{t} E_{TT}\big) & = & 
str_{q}\left(E_{TS} \sum_{(A, B) \in \widetilde{\Omega}^{t}} f_{AB} E_{AB} E_{TT} \right) \\
 & = & str_{q}\left(\sum_{B \in \widetilde{\cal{T}}^{t}} f_{SB} E_{TB} E_{TT} \right) \\
 & = & str_{q}(E_{TT}) f_{ST}  \\
 & \neq & 0,
\end{eqnarray*}
as $str_{q}(E_{TT}) \neq 0$ for all $T \in \widetilde{\cal{T}}^{t}$.
However, the fact that $P_{t} f P_{t} \in J_{t}$ implies that 
$$str_{q}\big( E_{TS} P_{t} f P_{t} E_{TT} \big) = str_{q}\big(E_{TT} E_{TS} P_{t} f P_{t} \big) = 0,$$ 
which is a contradiction, thus there does not exist any such $f \in {\cal{C}}_{t}$ and thus
$P_{t} {\cal{C}}_{t} P_{t} \cap J_{t} = 0$.

We define the inclusion $a \in {\cal{C}}_{t} \hookrightarrow {\cal{C}}_{t+1}$ by
$a \mapsto a \otimes \mathrm{id}$, and we can
regard each element of $J_{t}$ as an element of $J_{t+1}$ under this inclusion.
From Lemma \ref{lem:tomtheguineapig} each element in 
${\cal{C}}_{t+1}$ can be written as a linear combination of elements 
$(a \otimes \mathrm{id})$ and $(a \otimes \mathrm{id}) \check{R}_{t}^{\pm 1} (b\otimes \mathrm{id})$ 
where $a, b \in {\cal{C}}_{t}$ and 
$\check{R}_{t}^{\pm 1} = \mathrm{id}^{\otimes (t-1)} \otimes \check{\cal{R}}_{V,V} \in {\cal{C}}_{t+1}$.
Let $x \in J_{t}$, then we claim that $str_{q} \big( a \check{R}_{t}^{\pm 1} b x \big)=0$ for all
$a, b \in {\cal{C}}_{t}$, which we now prove using Lemma \ref{lem;kilo}, 
which we recall is still valid
for $U_{q}^{(N)}(\mathfrak{g})$ as discussed  after Eq. (\ref{eq:eenieweenie(1)}):
\begin{eqnarray*}
\lefteqn{
str_{q}^{\otimes (t+1)}\big(a \check{R}_{t}^{\pm 1}b x\big) 
= str_{q}^{\otimes t}(\mathrm{id}^{\otimes t} \otimes str_{q})\big(bxa\check{R}_{t}^{\pm 1}\big) } \\ 
& & \hspace{20mm} = \chi_{V}(v^{\mp 1}) str_{q}^{\otimes t} (bxa) 
= \chi_{V}(v^{\mp 1}) str_{q}^{\otimes t} (abx) 
= 0,
\end{eqnarray*}
where $str_{q}^{\otimes t}$ means that we take the quantum supertrace over $V^{\otimes t}$.

Let $T \in \widetilde{\cal{T}}^{t}$.  
Under the inclusion ${\cal{C}}_{t} \hookrightarrow {\cal{C}}_{t+1}$ discussed above we have
$$E_{TT} \mapsto \sum_{P \in {\cal{P}}^{+}_{shp(T)}} E_{T \circ P, T \circ P},$$
which corresponds to the decomposition of $V_{\mu} \otimes V$ 
into a direct sum of $U_{q}^{(N)}(\mathfrak{g})$-submodules
in Eq. (\ref{eq:moomoo}) where $\mu = shp(T)$ and $V_{\mu} = E_{TT} \big(V^{\otimes t}\big)$:
\begin{equation}
\label{eq:bec(10)}
\big( E_{TT} V^{\otimes t} \big) \otimes V = V_{\mu} \otimes V 
= \bigoplus_{P \in {\cal{P}}^{+}_{\mu} \cap \overline{\Lambda_{N}^{+}}} V_{P}
= \sum_{P  \in {\cal{P}}^{+}_{\mu} \cap \overline{\Lambda_{N}^{+}}} 
     E_{T \circ P, T \circ P} \left(V^{\otimes (t+1)}\right).
\end{equation}
Using the set 
${\cal{S}}_{\mu}$ defined by (\ref{eq:natalyakari(1)}), we can rewrite (\ref{eq:bec(10)}) as
\begin{eqnarray*}
 V_{\mu} \otimes V & = &
\left( \bigoplus_{P \in {\cal{P}}^{+}_{\mu} \cap \Lambda_{N}^{+}} V_{P} \right) \oplus 
\left( \bigoplus_{Q \in {\cal{S}}_{\mu}} V_{Q} \right)  \\
& = & \sum_{P \in {\cal{P}}^{+}_{\mu} \cap \Lambda_{N}^{+}} E_{T \circ P, T \circ P} V^{\otimes (t+1)} 
+ \sum_{Q \in {\cal{S}}_{\mu}} E_{T \circ Q, T \circ Q} V^{\otimes (t+1)}, 
\end{eqnarray*}
where the quantum superdimension of the module 
$E_{T \circ Q, T \circ Q} \left(V^{\otimes (t+1)} \right)$ is zero for
$Q \in {\cal{S}}_{\mu}$ as $Q \in \overline{\cal{P}}^{+}_{N} \backslash \Lambda_{N}^{+}$ and
$T \circ Q$ is a path belonging to $\widehat{\cal{T}}^{t+1} \backslash \widetilde{\cal{T}}^{t+1}$.
It is very important to note that $E_{T \circ Q, T \circ Q}$ belongs to $J_{t+1}$ from
Lemma \ref{lem:dalailama}, and  that
$\sum_{Q \in {\cal{S}}_{\mu}} E_{T \circ Q, T \circ Q}$ also  belongs to $J_{t+1}$.

We will now show that there is some integer $r \geq 2$ 
such that $(1-P_{s})$ belongs to $J_{s}$ for each integer $s \geq r$.  
To prove this we firstly note that
$\epsilon_{1} \in \Lambda_{N}^{+}$ and that $P_{1} = \mathrm{id}_{V} \notin J_{1}$.
Now $\Lambda_{N}^{+}$ is a proper subset of $\overline{\Lambda_{N}^{+}}$ and of ${\cal{P}}^{+}_{N}$,
 thus there is some integer $m \geq 2$ 
such that $\widehat{\cal{T}}^{m}$ contains at least one path that is not in 
$\widetilde{\cal{T}}^{m}$.  
Let us fix $r \geq 2$ to be the smallest integer such that 
$\widehat{\cal{T}}^{r}$ contains at least one path that is not in 
$\widetilde{\cal{T}}^{r}$, then
$$P_{(r-1)} = \mathrm{id}_{V^{\otimes (r-1)}}, \hspace{5mm} \mbox{and}
\hspace{5mm}  P_{r} \neq \mathrm{id}_{V^{\otimes r}}.$$
Note that $P_{i} = \mathrm{id}_{V^{\otimes i}}$ for all integers $i = 1, 2, \ldots, r-1$.

Recall the definition (\ref{eq:natalyakari(1)}) 
of ${\cal{S}}_{shp(T)}$ for a path $T$ of length $r-1$.
Under the inclusion ${\cal{C}}_{r-1} \hookrightarrow {\cal{C}}_{r}$ we have
\begin{eqnarray}
P_{(r-1)} & \mapsto & \sum_{ T \in \widetilde{\cal{T}}^{(r-1)} } \left(
\sum_{P \in {\cal{P}}^{+}_{shp(T)} \cap \Lambda_{N}^{+}} E_{T \circ P, T \circ P}  
+ \sum_{Q \in {\cal{S}}_{shp(T)}} E_{T \circ Q, T \circ Q}  \right), \nonumber \\
          & =   & \sum_{ T \in \widetilde{\cal{T}}^{r} } E_{TT}
+ \sum_{ T \in \widetilde{\cal{T}}^{(r-1)} } \left( \sum_{Q \in {\cal{S}}_{shp(T)}} E_{T \circ Q, T \circ Q} \right) 
\nonumber \\
	  & = & P_{r} 
+ \sum_{ T \in \widetilde{\cal{T}}^{(r-1)} } \left( \sum_{Q \in {\cal{S}}_{shp(T)}} E_{T \circ Q, T \circ Q} \right).
 \label{eq:bec(15)}
\end{eqnarray}
As $P_{(r-1)} = \mathrm{id}_{V^{\otimes (r-1)}}$ and we have
$P_{(r-1)} \mapsto P_{(r-1)} \otimes \mathrm{id}$ under the inclusion 
${\cal{C}}_{r-1} \hookrightarrow {\cal{C}}_{r}$, 
(\ref{eq:bec(15)}) equals $\mathrm{id}_{V^{\otimes r}}$.  
We reiterate once more that $P_{r} \neq \mathrm{id}_{V^{\otimes r}}$.

As $P_{(r-1)} \otimes \mathrm{id} = \mathrm{id}_{V^{\otimes r}}$,  
rewriting (\ref{eq:bec(15)}) gives us
$$(1-P_{r}) = 
\sum_{ T \in \widetilde{\cal{T}}^{(r-1)} } 
\left( \sum_{Q \in {\cal{S}}_{shp(T)}} E_{T \circ Q, T \circ Q} \right),$$
and the two summations on the right hand side of this expression can be rewritten as a single sum:
$$\sum_{ T \in \widetilde{\cal{T}}^{(r-1)} } 
\left( \sum_{Q \in {\cal{S}}_{shp(T)}} E_{T \circ Q,T \circ Q} \right) = 
\sum_{D \in \widehat{\cal{T}}^{r} \backslash \widetilde{\cal{T}}^{r}} E_{DD},$$
and therefore 
$$(1-P_{r}) = \sum_{D \in \widehat{\cal{T}}^{r} \backslash \widetilde{\cal{T}}^{r}} E_{DD},$$
where each $E_{DD}$ is an element of $J_{r}$ from Lemma \ref{lem:dalailama}.

Now for each integer $s \geq r$ we can similarly show that
$$P_{s} - P_{(s+1)} = 
\sum_{Q \in \widehat{\cal{T}}^{s+1} \backslash \widetilde{\cal{T}}^{s+1}} E_{QQ} \in J_{(s+1)},$$
and by expressing $(1-P_{(s+1)})$ as the sum:
$$(1-P_{(s+1)}) = (1-P_{r}) + (P_{r} - P_{(r+1)}) + \cdots + (P_{s} - P_{(s+1)}),$$
we have $ (1-P_{(s+1)}) \in J_{(s+1)}$ as
$(1-P_{r})$ and $(P_{i}-P_{(i+1)})$ belong to $J_{(i+1)}$ for all $i=r, \ldots, s$ under the inclusion
${\cal{C}}_{i} \hookrightarrow {\cal{C}}_{(i + 1)}$.
This proves that $(1-P_{j}) \in J_{j}$ for all $j \in \mathbb{N}$.

Note that 
\begin{eqnarray}
{\cal{C}}_{t} & = & \big(P_{t} + (1-P_{t})\big) {\cal{C}}_{t} \big(P_{t}+(1-P_{t})\big) \nonumber \\
              & = & P_{t} {\cal{C}}_{t} P_{t} + (1-P_{t}) {\cal{C}}_{t} P_{t}+
                    P_{t} {\cal{C}}_{t} (1-P_{t}) + (1-P_{t}) {\cal{C}}_{t} (1-P_{t}). 
		    \label{eq:margiesgettingmarried}
\end{eqnarray}	
Now $(1-P_{t})$ is in $J_{t}$, and each of 
$(1-P_{t}) x P_{t}$, $P_{t} x (1-P_{t})$ and
$(1-P_{t}) x (1-P_{t})$ are in $J_{t}$ for each $x \in {\cal{C}}_{t}$, thus
$$\big( (1-P_{t}) {\cal{C}}_{t} P_{t}+
 P_{t} {\cal{C}}_{t} (1-P_{t}) + (1-P_{t}) {\cal{C}}_{t} (1-P_{t}) \big) \subseteq J_{t}.$$
We previously proved that $P_{t} {\cal{C}}_{t} P_{t} \cap J_{t} = 0$, thus any element
$x \in {\cal{C}}_{t}$ belonging to $J_{t}$ must also belong to
$\big( (1-P_{t}) {\cal{C}}_{t} P_{t}+
 P_{t} {\cal{C}}_{t} (1-P_{t}) + (1-P_{t}) {\cal{C}}_{t} (1-P_{t}) \big)$
 from (\ref{eq:margiesgettingmarried}).
We thus obtain
$$            
{\cal{C}}_{t} = P_{t} {\cal{C}}_{t} P_{t} \oplus J_{t}.
$$
We wish to show that $(1 - P_{t})$ generates $J_{t}$ as a two-sided ideal.
To do this we will prove two assertions: 
\begin{itemize}
\item[(a)] $x (1-P_{t}) y$ belongs to $J_{t}$ for all $x, y \in {\cal{C}}_{t}$,  and 
\item[(b)] each element of $J_{t}$ belongs to ${\cal{C}}_{t} (1-P_{t}) {\cal{C}}_{t}$.
\end{itemize}
The proof of (a) follows easily from the fact that $(1 - P_{t}) \in J_{t}$ 
and the properties of the quantum supertrace, but the proof of (b) is more involved.
Let $z$ be an arbitrary element of $J_{t}$, then
$$z = 
P_{t} z P_{t} + (1-P_{t}) z P_{t}+P_{t} z (1-P_{t}) + (1-P_{t}) z (1-P_{t}).$$
As  $P_{t} z P_{t} \in J_{t}$, 
but $P_{t} {\cal{C}}_{t} P_{t} \cap J_{t}=0$, we have $P_{t} z P_{t} = 0$. 
Then
$$
z = (1-P_{t}) z P_{t}+ P_{t} z (1-P_{t}) + (1-P_{t}) z (1-P_{t}) 
  = z(1-P_{t}) + (1-P_{t})zP_{t},
$$ 
which belongs to ${\cal{C}}_{t} (1-P_{t}) {\cal{C}}_{t}$, proving (b).
This completes the proof that $(1 - P_{t})$ generates $J_{t}$ as a two-sided ideal.

We now prove part (ii) of Proposition \ref{lem:concorde}.
 Let $a$, $b \in {\cal{C}}_{t}$ be arbitrary, then
$$P_{t}abP_{t}=P_{t} a \big(P_{t} + (1-P_{t})\big) b P_{t} = 
P_{t}aP_{t}bP_{t} = (P_{t}aP_{t})(P_{t}bP_{t}),$$
as $P_{t} a (1-P_{t}) b P_{t} \in J_{t}$ and $P_{t} {\cal{C}}_{t}P_{t} \cap J_{t} = 0$.
\end{proof}

\end{subsection}

\begin{subsection}{The Tensor Product Theorems}
\label{subsection:turkmenistan(1)}

We now prove the tensor product theorems.

\begin{theorem}
\label{th:firsttensor}
Let $n \geq 2$ and $N \geq 4$ be even,
or let $n=1$ and $N \geq 6$ satisfy $N \equiv 2 \pmod{4}$.
Furthermore, let $N$ be sufficiently large enough so that $\epsilon_{1} \in \Lambda_{N}^{+}$.  
Then for each $t \in \mathbb{Z}_{+}$, there is a decomposition of $V^{\otimes t}$ into
a direct sum of $U_{q}^{(N)}(\mathfrak{g})$-submodules
\begin{equation}
\label{eq:bloodycold}
V^{\otimes t} = {\cal{V}} \oplus {\cal{Z}},
\end{equation}
where ${\cal{V}}$ is of the form
$${\cal{V}} = \bigoplus_{\lambda \in \Lambda_{N}^{+}} \left(V_{\lambda}\right)^{\oplus n_{t}(\lambda)},$$
with $n_{t}(\lambda) \in \mathbb{Z}_{+}$ being the 
(possibly zero) number of copies of 
the $U_{q}^{(N)}(\mathfrak{g})$-submodule
$V_{\lambda}$ in ${\cal{V}}$.  Here ${\cal{Z}}$ is a possibly vanishing 
$U_{q}^{(N)}(\mathfrak{g})$-submodule with the property that $str_{q}(f)=0$ for all 
$f \in {\cal{C}}_{t}$ satisfying $f(V^{\otimes t}) \subseteq {\cal{Z}}$.
\end{theorem}
\begin{proof}
Let $P_{t}$ be as given in Lemma \ref{lem:PP} and  
fix ${\cal{V}} = P_{t}V^{\otimes t}$ and 
${\cal{Z}} = (1-P_{t}) V^{\otimes t}$.
As $P_{t}$ and $(1-P_{t})$ are orthogonal idempotents, 
the sum on the right hand side of (\ref{eq:bloodycold}) is direct.
To prove that ${\cal{Z}}$ has the claimed property, let
$f \in {\cal{C}}_{t}$ satisfy $f(V^{\otimes t}) \subseteq {\cal{Z}}$. 
Then $f = \big(P_{t} + (1-P_{t})\big)f$. 
Now $P_{t} f (V^{\otimes t}) \subseteq P_{t} {\cal{Z}} = 
P_{t}(1-P_{t})V^{\otimes t} = 0$, thus $P_{t} f = 0$ and $f=(1-P_{t})f$.
Then $str_{q}(f) = str_{q}\big((1-P_{t})f\big)=0$ as $(1-P_{t}) \in J_{t}$, completing the proof.
\end{proof}

\begin{lemma}
Assume that the given conditions on $n$ and $N$ are the same as in Theorem \ref{th:firsttensor}.
Let $V^{\otimes t} = {\cal{V}}_{t} \oplus {\cal{Z}}_{t}$ 
be the decomposition of $V^{\otimes t}$ into $U_{q}^{(N)}(\mathfrak{g})$-submodules
for each $t \in \mathbb{N}$ given in Eq. (\ref{eq:bloodycold}).
Then ${\cal{Z}}_{t} \otimes V \subseteq {\cal{Z}}_{t+1}$.
\end{lemma}
\begin{proof}
This is similar to the proof of \cite[Thm. 5.5.2]{tw}.
From the definition of $P_{t}$ and the inclusion 
$P_{t} \hookrightarrow P_{t} \otimes \mathrm{id} \in {\cal{C}}_{t+1}$ we have
\begin{eqnarray}
P_{t} \otimes \mathrm{id} & = & \sum_{S \in \widetilde{\cal{T}}^{t}}  E_{SS} \otimes \mathrm{id} 
= \sum_{S \in \widetilde{\cal{T}}^{t} } \left(
 \sum_{ \stackrel{\mu \in \overline{\Lambda_{N}^{+}}}{shp(S) \leq \mu}} E_{S \circ \mu, S \circ \mu} 
 \right) \label{eq:guanaverras(1)} \\
& = & P_{(t+1)} + \sum_{T \in \widehat{\cal{T}}^{t+1} \backslash \widetilde{\cal{T}}^{t+1} } E_{TT}.
      \label{eq:guanaverras(2)}
\end{eqnarray}
To be certain that (\ref{eq:guanaverras(1)}) is true, note that $E_{SS}$ is a path projection 
of length $t$ projecting 
down from $V^{\otimes t}$ onto a $U_{q}^{(N)}(\mathfrak{g})$-submodule $V_{\nu}$ where 
$\nu = shp(S) \in \Lambda_{N}^{+}$, and recall from (\ref{eq:moomoo}) that
$$V_{\nu} \otimes V = 
\bigoplus_{\mu \in {\cal{P}}^{+}_{\nu} \cap \overline{\Lambda_{N}^{+}}} V_{\mu}.$$
Then (\ref{eq:guanaverras(1)}) follows from the fact that 
${\cal{P}}^{+}_{\nu} \cap \overline{\Lambda_{N}^{+}}$ is the set of all elements $\xi$ of
$\overline{\Lambda_{N}^{+}}$ where $\nu \leq \xi$, that is, 
${\cal{P}}^{+}_{\nu} \cap \overline{\Lambda_{N}^{+}}$ is the set of
all the $\xi$ connected to $\nu$ in the relevant Bratelli diagram where $\xi$ is on the 
level of the Bratteli diagram
immediately below the level containing $\nu$.

To see that (\ref{eq:guanaverras(2)}) is correct, note that each vertex 
$\nu \in \Lambda_{N}^{+}$ on the $t^{th}$ level
of the Bratteli diagram  is connected to vertices on the 
$(t+1)^{st}$ level of the Bratteli diagram where all of these latter vertices 
are elements of $\overline{\cal{P}}^{+}_{N}$.

Note that the image of  
$\displaystyle{\sum_{T \in \widehat{\cal{T}}^{t+1} \backslash \widetilde{\cal{T}}^{t+1} } E_{TT} }$ 
is contained in ${\cal{Z}}_{t+1}$.  
Hence 
$${\cal{Z}}_{t} \otimes V = 
\left( (1-P_{t+1}) - 
      \sum_{T \in \widehat{\cal{T}}^{t+1} \backslash \widetilde{\cal{T}}^{t+1} } E_{TT} \right) 
      V^{\otimes (t+1)} \subseteq {\cal{Z}}_{t+1}.$$

\end{proof}

\begin{theorem}
\label{th:secondtensor}
Assume that the given conditions on $n$ and $N$ are the same as in Theorem \ref{th:firsttensor}.
Let  $s \in \mathbb{N}$ and $\lambda_{i} \in \Lambda_{N}^{+}$ for each $1 \leq i \leq s$.  
Let $V_{\lambda_{i}}$ be a  
$U_{q}^{(N)}(\mathfrak{g})$-module defined in Definition \ref{lem:representationsyippe}:
$$\tilde{p}^{t_{i}}_{j_{i}}[\lambda_{i}]: V^{\otimes t_{i}} \rightarrow V_{\lambda_{i}}, 
\hspace{10mm}  i=1, \ldots, s.$$   
Then there is a decomposition of 
$V_{\lambda_{1}} \otimes \cdots \otimes V_{\lambda_{s}}$
into a direct sum of $U_{q}^{(N)}(\mathfrak{g})$-submodules:
\begin{equation}
\label{eq:Iamcold!}
V_{\lambda_{1}} \otimes \cdots \otimes V_{\lambda_{s}} = {\cal{V}} \oplus {\cal{Z}},
\end{equation}
where ${\cal{V}}$ is a direct sum of $q$-admissible submodules
and ${\cal{Z}}$ is a possibly vanishing submodule 
with the property that $str_{q}(f)=0$ for all $f \in {\cal{C}}_{t}$ satisfying
$f \big( V_{\lambda_{1}} \otimes \cdots \otimes V_{\lambda_{s}} \big) 
\subseteq {\cal{Z}}$.
\end{theorem}
\begin{proof}
We can write the left hand side of (\ref{eq:Iamcold!}) as
$$
V_{\lambda_{1}} \otimes \cdots \otimes V_{\lambda_{s}} = 
 e \big( V^{\otimes t_{1}} \otimes \cdots \otimes V^{\otimes t_{s}}\big) \subseteq V^{\otimes t},$$
where $t=\sum_{i=1}^{s} t_{i}$, and
$$e = \tilde{p}^{t_{1}}_{j_{1}}[\lambda_{1}]  \otimes \cdots \otimes
\tilde{p}^{t_{s}}_{j_{s}}[\lambda_{s}] \in {\cal{C}}_{t}.$$  
Then $e^{2}=e$, and we can use Proposition \ref{lem:concorde} (ii) to prove that
\begin{equation}
\label{eq:flightoftheconcordes(1)}
\big(P_{t} e P_{t}\big)^{2} = P_{t}eP_{t}eP_{t} = P_{t}eP_{t}.
\end{equation}
In the same way, we can show that $e P_{t} e P_{t} e$ is an idempotent, and thus so is
$e-e P_{t} e P_{t} e$.  Clearly these two idempotents are orthogonal to each other.
The method we use to prove the theorem is to show that
$eP_{t}eP_{t}e \big(V^{\otimes t}\big)$ is isomorphic to
$P_{t} e P_{t} (V^{\otimes t})$, which is a direct sum of $q$-admissible modules, and 
that $(e-eP_{t}eP_{t}e) V^{\otimes t}$ satisfies the given properties of ${\cal{Z}}$
in (\ref{eq:Iamcold!}).

Firstly, we will show that $(e-eP_{t}eP_{t}e) \in J_{t}$. 
Clearly, $(1-P_{t})(e-eP_{t}eP_{t}e)(1-P_{t}) \in J_{t}$, and
\begin{eqnarray*}
\lefteqn{ (1-P_{t})(e-eP_{t}eP_{t}e)(1-P_{t}) } \\
& & \hspace{15mm} = (e-eP_{t}eP_{t}e)(1-P_{t}) - (P_{t}e-P_{t}eP_{t}e)(1-P_{t}) \\
& & \hspace{15mm} = (e-eP_{t}eP_{t}e) - (eP_{t}-eP_{t}eP_{t}) + (P_{t}eP_{t}e-P_{t}e)(1-P_{t}).
\end{eqnarray*}
Now $P_{t}(eP_{t}-eP_{t}eP_{t})=0$ as $P_{t}$ and $P_{t}eP_{t}$ are idempotents, and so we can write
\begin{eqnarray*}
\lefteqn{ (1-P_{t})(e-eP_{t}eP_{t}e)(1-P_{t}) } \\
 & & \hspace{10mm} 
      = (e-eP_{t}eP_{t}e) + (P_{t}eP_{t}e-P_{t}e)(1-P_{t})-(1-P_{t})(eP_{t}-eP_{t}eP_{t}),
\end{eqnarray*}
which implies that $$(e-eP_{t}eP_{t}e) \in J_{t},$$ as both
$(P_{t}eP_{t}e-P_{t}e)(1-P_{t})$ and $(1-P_{t})(eP_{t}-eP_{t}eP_{t})$ are elements of $J_{t}$.
Clearly $str_{q}(f)=0$ for all $f \in {\cal{C}}_{t}$ satisfying
$f \big( V_{\lambda_{1}} \otimes \cdots \otimes V_{\lambda_{s}} \big) 
\subseteq (e-eP_{t}eP_{t}e)V^{\otimes t}$.

Fix $A = eP_{t}eP_{t}e$ and $B = P_{t}eP_{t}$, 
then $A$ and $B$ are $U_{q}^{(N)}(\mathfrak{g})$-linear idempotents satisfying
\begin{equation}
\label{eq:neverknowthiscouldbethelast(0)}
ABA = A, \hspace{5mm} \mbox{and} \hspace{5mm} BAB=B.
\end{equation}
Write $$V_{A} = A\big(V^{\otimes t}\big), \hspace{5mm} \mbox{and} 
                \hspace{5mm} V_{B} = B\big(V^{\otimes t}\big),$$
then $V_{A} \cong V_{B}$.  
While this is easy to show, we prove it here for clarity.
Clearly
$A \big(V_{B}\big) \subseteq V_{A}$, so
$BA \big(V_{B}\big) \subseteq B\big(V_{A}\big)$, and we can rewrite this using
(\ref{eq:neverknowthiscouldbethelast(0)}) as  
$V_{B} \subseteq B\big(V_{A}\big)$.
In addition,
$
B \big( V_{A} \big) \subseteq V_{B}$, so
$AB \big(V_{A} \big) \subseteq A \big(V_{B} \big)$, and we can rewrite this using
(\ref{eq:neverknowthiscouldbethelast(0)}) as 
$V_{A} \subseteq A \big(V_{B}\big)$.
Then 
$$A \big(V_{B}\big) = V_{A} \hspace{5mm} \mbox{and} \hspace{5mm} 
B\big(V_{A}\big) = V_{B}.$$
The idempotents $A$ and $B$ are
$U_{q}^{(N)}(\mathfrak{g})$-linear intertwiners between $V_{A}$ and $V_{B}$, so $V_{A} \cong V_{B}$.
Now $V_{B} = P_{t}eP_{t} \big( V^{\otimes t} \big)$ is a direct sum of $q$-admissible modules, so
$V_{A} = eP_{t}eP_{t}e \big(V^{\otimes t}\big)$ is also a direct sum of $q$-admissible modules.
Fixing ${\cal{V}} = eP_{t}eP_{t}e \big(V^{\otimes t}\big)$ and
${\cal{Z}} = (e-eP_{t}eP_{t}e)V^{\otimes t}$ completes the proof.

\end{proof}

We say that a $U_{q}^{(N)}(\mathfrak{g})$-module $W_{\mu}$ is {\emph{q-admissible}} 
if $W_{\mu} = p_{\mu} (V^{\otimes t})$ for some idempotent $p_{\mu} \in {\cal{C}}_{t}$ and 
if $W_{\mu}$ is isomorphic to a $U_{q}^{(N)}(\mathfrak{g})$-module 
$V_{\lambda} = E_{TT} \big(V^{\otimes t}\big)$ for some $\lambda \in \Lambda_{N}^{+}$. 
Here $E_{TT} \in End_{U_{q}^{(N)}(\mathfrak{g})}(V^{\otimes t})$ 
is the path projection associated with the path 
$T \in \widetilde{\cal{T}}^{t}$ of length $t$ with $shp(T) = \lambda$.
Each $q$-admissible $U_{q}^{(N)}(\mathfrak{g})$-module has non-vanishing quantum superdimension.

We finish off this section by making the following conjecture.  
\begin{conjecture}
Assume that the given conditions on $n$ and $N$ are the same as in Theorem \ref{th:firsttensor}.
Let  $s \in \mathbb{N}$ and $\lambda_{i} \in \Lambda_{N}^{+}$ for each $1 \leq i \leq s$.  
Let $V_{\lambda_{i}}$ be a  
$U_{q}^{(N)}(\mathfrak{g})$-module defined in Definition \ref{lem:representationsyippe}:
$$\tilde{p}^{t_{i}}_{j_{i}}[\lambda_{i}]: V^{\otimes t_{i}} \rightarrow V_{\lambda_{i}}, 
\hspace{10mm} i=1, \ldots, s.$$   
Let  
$V_{\lambda_{1}} \otimes \cdots \otimes V_{\lambda_{s}} \subseteq V^{\otimes t}$
where $t = \sum_{i=1}^{s} t_{i}$, then given any decomposition of
$V_{\lambda_{1}} \otimes \cdots \otimes V_{\lambda_{s}}$
 into a direct sum of $U_{q}^{(N)}(\mathfrak{g})$-submodules:
\begin{equation}
\label{eq:oliveoil(2)}
V_{\lambda_{1}} \otimes \cdots \otimes V_{\lambda_{s}} 
= \widetilde{\cal{V}} \oplus \widetilde{\cal{Z}},
\end{equation} 
 where
$\widetilde{\cal{V}}$ is a direct sum of $q$-admissible submodules and
$\widetilde{\cal{Z}}$ is a possibly vanishing submodule with the property that
$str_{q}(f)=0$ for all $f \in {\cal{C}}_{t}$ satisfying 
$f(V_{\lambda_{1}} \otimes \cdots \otimes V_{\lambda_{s}}) \subseteq \widetilde{\cal{Z}}$, then
$\widetilde{\cal{V}} \cong {\cal{V}}$ and $\widetilde{\cal{Z}} \cong {\cal{Z}}$
where $V_{\lambda_{1}} \otimes \cdots \otimes V_{\lambda_{s}}={\cal{V}} \oplus {\cal{Z}}$ 
is the decomposition of
$V_{\lambda_{1}} \otimes \cdots \otimes V_{\lambda_{s}}$ into 
$U_{q}^{(N)}(\mathfrak{g})$-submodules given in Eq. (\ref{eq:Iamcold!}).  
\end{conjecture}

\end{subsection}

\end{section}

\begin{section}{The well-definedness of the projection operators}
\label{subsec:vermittin}
\markright{\text{The well-definedness of the projection operators}}

In this section we present the detailed
proof of Lemma \ref{lem:projectionsarewelldefined}.
Let $\lambda_{i}^{t} = (0, \epsilon_{1}, s_{2}, s_{3}, \ldots, s_{t-1}, \lambda) \in \widehat{\cal{T}}^{t}.$
The projection $\tilde{p}_{i}^{t}[\lambda]: V^{\otimes t} \rightarrow V_{\lambda}$ 
is well defined if for each $s_{j} \in \lambda_{i}^{t}$, 
\begin{equation}
\label{eq:freedom}
\prod_{ \stackrel{\mu \in {\cal P}_{s_{j}}^{+}}{\mu \neq s_{j+1}}}
\big( \chi_{s_{j+1}}(v) - \chi_{\mu}(v) \big) \neq 0.
\end{equation}  
We now write $\lambda$ instead of $s_{j}$.  From Def. \ref{def:twopointthreepointten},
${\cal P}_{\lambda}^{+} \subseteq {\cal P}_{\lambda}^{0} =
 \big\{ \lambda, \lambda \pm \epsilon_{j} \in {\cal P}^{+} | \ 1 \leq j \leq n \big\}$.
Clearly the following equation
\begin{equation}
\label{eq:sincity3}
\prod_{\stackrel{\mu \in {\cal P}_{\lambda}^{0}} {\mu \neq s_{j+1}}} \left( \chi_{s_{j+1}}(v) - \chi_{\mu}(v) \right) \neq 0,
\end{equation}
implies (\ref{eq:freedom}).  Note that (\ref{eq:sincity3}) is true if and only if
\begin{equation}
\label{eq:abigequatliontypething}
\big( (s_{j+1} + 2\rho, s_{j+1}) - (\mu + 2\rho, \mu) \big) \not\equiv 0 \pmod{N},
\end{equation}
for all $s_{j+1}, \mu \in {\cal P}_{\lambda}^{0}$ where $\mu \neq s_{j+1}$.  
Let us write $\xi$ instead of $s_{j+1}$ for convenience.

In order to show that (\ref{eq:abigequatliontypething}) is true, 
we will find the largest and smallest elements of
$${\cal{S}}=\Big\{ \big|(\mu + 2\rho, \mu) - (\xi + 2\rho, \xi) \big| \in \mathbb{R} \ \Big|
  \ \mu, \xi \in {\cal{P}}^{0}_{\lambda}, \ \mu \neq \xi \Big\}.$$
To help with this we note the following
inequalities for  each $\lambda \in {\cal{P}}^{+}$ and each $i=1, 2, \ldots, n-1$:
$$
(\lambda + 2\rho, \lambda) \leq (\lambda + \epsilon_{i+1} + 2\rho, \lambda + \epsilon_{i+1}) \leq
(\lambda + \epsilon_{i} + 2\rho, \lambda + \epsilon_{i}),
$$
$$
(\lambda - \epsilon_{i} + 2\rho, \lambda - \epsilon_{i}) \leq (\lambda - \epsilon_{i+1} + 2\rho, \lambda - \epsilon_{i+1})
\leq (\lambda + 2\rho, \lambda).
$$
From these inequalities, the largest element of ${\cal{S}}$ is
$$(\lambda + \epsilon_{1} + 2\rho, \lambda + \epsilon_{1}) - (\lambda - \epsilon_{1} + 2\rho, \lambda - \epsilon_{1}) = 
4 \lambda_{1} + 4n - 2,$$
and the second largest element of ${\cal{S}}$ is
$$(\lambda + \epsilon_{1} + 2\rho, \lambda + \epsilon_{1}) - (\lambda - \epsilon_{2} + 2\rho, \lambda - \epsilon_{2}) =
(2\lambda + 2\rho, \epsilon_{1} + \epsilon_{2}) = 2\lambda_{1} + 2\lambda_{2} + 4n-4.$$

We now determine the smallest element of ${\cal{S}}$.
For each $\lambda \in {\cal{P}}^{+}$ and $i = 1, 2, \ldots, n-1$, we have 
$$2 \leq (\lambda + \epsilon_{i}+2\rho,\lambda+\epsilon_{i})-(\lambda+\epsilon_{i+1}+2\rho, \lambda+\epsilon_{i+1}),$$
$$2 \leq (\lambda-\epsilon_{i+1}+2\rho, \lambda-\epsilon_{i+1})-(\lambda-\epsilon_{i}+2\rho, \lambda-\epsilon_{i}),$$
and
$$2 \leq (\lambda + \epsilon_{n} + 2\rho, \lambda + \epsilon_{n})-(\lambda + 2\rho, \lambda),$$
$$2 \leq 2 \lambda_{n} \leq (\lambda + 2\rho, \lambda)-(\lambda-\epsilon_{n}+2\rho, \lambda-\epsilon_{n}),$$
as $\lambda-\epsilon_{n}$ is an integral dominant weight.

\begin{lemma}
\label{lem:spainishwooingsong}
The projections from Lemma \ref{lem:projectionsarewelldefined} (a) are well defined when 
$N \equiv 0 \pmod{4}$.
\end{lemma}
\begin{proof}
The projections are well-defined if for each $\lambda \in \Lambda_{N}^{+}$ we have
\begin{equation}
\label{eq:assetagassi}
\big| (\mu+2\rho, \mu)-(\xi + 2\rho, \xi) \big| \not\equiv 0 \pmod{N},
\end{equation}
for all $\mu,\xi \in {\cal{P}}^{0}_{\lambda}$ where $\mu \neq \xi$.
As before, the smallest value of the left hand side of (\ref{eq:assetagassi}) is $2$ 
and the largest is $4\lambda_{1}+4n-2$.  
From the definition of $\Lambda_{N}^{+}$, the components of $\lambda$ satisfy 
$0 \leq \lambda_{1} + \lambda_{2} \leq N/2-2n+1$, thus
$$0 \leq 4 \lambda_{1} + 4n-2 \leq 2N-4n+2,$$
and $2N-4n+2<2N$.  
As $\lambda_{1} \in \mathbb{Z}_{+}$, we have 
$4 \lambda_{1} + 4n-2 \in 4 \mathbb{Z}_{+} + 2$, then it follows
that $4 \lambda_{1} + 4n-2 \neq N$ as $N \equiv 0 \pmod{4}$.

To complete the proof, we note 
that the second largest value of the left hand side of (\ref{eq:assetagassi}) is
$2\lambda_{1} + 2\lambda_{2} + 4n-4$.  
Now from the inequality 
$\lambda_{1} + \lambda_{2} \leq N/2-2n+1$ 
we obtain
$2\lambda_{1} + 2\lambda_{2} + 4n-4 \leq N-2$,
completing the proof.
\end{proof}

\begin{lemma}
\label{lem:thekingofspain}
The projections from Lemma \ref{lem:projectionsarewelldefined} (a) 
are well defined for $N \equiv 2 \pmod{4}$.
\end{lemma}
\begin{proof}

For $\lambda \in \Lambda_{N}^{+}$, we have 
$0 \leq \lambda_{1} \leq N/4-n-1/2$.  As previously, the projections are well-defined if 
(\ref{eq:assetagassi}) is true.  The smallest value of the left hand side of (\ref{eq:assetagassi}) 
is $2$ and the largest value is $4\lambda_{1}+4n-2$.
Thus
$$6 \leq 4\lambda_{1} + 4n-2 \leq N-4,$$
and  the projections are well-defined.

\end{proof}

\begin{lemma}
\label{lem:oooooooodddddddd}
For each odd  $N \geq 3$, the projection $\tilde{p}_{i}^{t}[\lambda]$ is well-defined 
for each $\lambda_{i}^{t} \in \widehat{\cal{T}}^{t}$ if 
\begin{itemize}
\item[(i)] $\lambda_{1} \leq (N-1)/2-n+1$, or if
\item[(ii)] the components of $s_{t-1}= \overline{\lambda} \in \Lambda_{N}^{+}$  satisfy 
$\overline{\lambda}_{1} = (N-1)/2-n+1$ and $\overline{\lambda}_{2} = \overline{\lambda}_{1}$, and $\lambda$ is such that 
$\lambda = \overline{\lambda} + \epsilon_{1}$.
\end{itemize}
\end{lemma}
\begin{proof}
As previously, the projection is well-defined if (\ref{eq:assetagassi}) is true
for all pairs $(\mu,\xi) \in {\cal{P}}^{0}_{\lambda} \times {\cal{P}}^{0}_{\lambda}$ 
where $\mu \neq \xi$, for each $\lambda \in \Lambda_{N}^{+}$.  
We will show that this is true only when parts (i) or (ii) of the lemma are satisfied.
As before, the smallest value of the left hand side of (\ref{eq:assetagassi}) is $2$, 
the largest value is 
$4\lambda_{1}+4n-2$, and the second largest value is $2\lambda_{1} + 2\lambda_{2} + 4n-4$.

Consider the largest value.
Let $\lambda \in \Lambda_{N}^{+}$, then 
$0 \leq \lambda_{1} + \lambda_{2} \leq N-2n+1$,  thus
$$0 \leq 4 \lambda_{1} + 4n-2 < 4N.$$
Note that $4 \lambda_{1} + 4n-2 \neq N$ and that $4 \lambda_{1} + 4n-2 \neq 3N$ 
as the left hand sides of each of these is even.  However, we may have
$4 \lambda_{1} + 4n-2 = 2N$ and it transpires that this results in part (i) of the lemma.

Let us consider the second largest value of the left hand side of (\ref{eq:assetagassi}).  
From the relations $0 \leq \lambda_{1} + \lambda_{2} \leq N-2n-1$ we have
 $$4 \leq 2\lambda_{1} + 2\lambda_{2} + 4n-4 \leq 2N-2.$$
Thus if we can show that the left hand side of (\ref{eq:assetagassi}) is always even, the only possible case in which 
(\ref{eq:assetagassi}) is not satisfied is when $4\lambda_{1} + 4n-2=2N$.
We will now show that the left hand side of (\ref{eq:assetagassi}) is always even.  
To do this, we  consider
all the possible cases below, in which we let $i, j \in \{1, 2, \ldots, n\}$.
\begin{itemize}
\item[(i)]  Set $\mu=\lambda + \epsilon_{i}$ and $\xi = \lambda$.  Then
$(\lambda + \epsilon_{i} + 2\rho, \lambda + \epsilon_{i}) - (\lambda + 2\rho, \lambda) = 
(2\lambda+2\rho, \epsilon_{i})+1 \in 2 \mathbb{Z}_{+}$.
\item[(ii)]  Set $\mu=\lambda$ and $\xi=\lambda-\epsilon_{i}$.  Then
$(\lambda+2\rho,\lambda)-(\lambda-\epsilon_{i}+2\rho, \lambda-\epsilon_{i})=
(2\lambda+2\rho, \epsilon_{i})-1 \in 2\mathbb{Z}_{+}$.
\item[(iii)]  Set $\mu=\lambda+\epsilon_{i}$ and $\xi=\lambda-\epsilon_{j}$ where $i \neq j$.  Then
$(\lambda+\epsilon_{i} + 2\rho, \lambda+\epsilon_{i})-(\lambda-\epsilon_{j} + 2\rho, \lambda-\epsilon_{j}) = 
(2\lambda + 2\rho, \epsilon_{i} + \epsilon_{j}) \in 2 \mathbb{Z}_{+}$.
\item[(iv)]  Set $\mu=\lambda-\epsilon_{i}$ and $\xi=\lambda-\epsilon_{j}$ where $i \neq j$.  Then
$(\lambda-\epsilon_{i} + 2\rho, \lambda-\epsilon_{i})-(\lambda-\epsilon_{j} + 2\rho, \lambda-\epsilon_{j})
=(2\lambda+2\rho, \epsilon_{j}-\epsilon_{i}) \in 2 \mathbb{Z}$.
\item[(v)]  Set $\mu=\lambda+\epsilon_{i}$ and $\xi=\lambda+\epsilon_{j}$ where $i \neq j$.  Then
$(\lambda+\epsilon_{i} + 2\rho, \lambda+\epsilon_{i})-(\lambda+\epsilon_{j} + 2\rho, \lambda+\epsilon_{j}) = 
(2\lambda+2\rho, \epsilon_{i}-\epsilon_{j}) \in 2 \mathbb{Z}$.
\item[(vi)]  Set $\mu = \lambda+\epsilon_{i}$ and  $\xi=\lambda-\epsilon_{i}$.  Then
$(\lambda+\epsilon_{i} + 2\rho, \lambda+\epsilon_{i})-(\lambda-\epsilon_{i} + 2\rho, \lambda-\epsilon_{i}) = 
(2\lambda+2\rho, 2\epsilon_{i}) \in 2 \mathbb{Z}_{+}$ and
$$(2\lambda + 2\rho, 2 \epsilon_{i}) = 2(2\lambda_{i}+2n-2i+1) < \left\{
\begin{array}{ll}
4N, & \mbox{when } i=1, \\
2N,  & \mbox{when } i \geq 2.
\end{array} \right.$$
\end{itemize}
The only possible case in which (\ref{eq:assetagassi}) is {\emph{not}} satisfied is when
$\mu=\lambda+\epsilon_{1}$ and  $\xi=\lambda-\epsilon_{1}$, and here it may be possible that
$(\xi + 2\rho, \xi) - (\mu + 2\rho, \mu) = 2N$.  
We will show that for each odd $N \geq 3$ there always exists at least one
element $\lambda$ of $\Lambda_{N}^{+}$ with the property that 
\begin{equation}
\label{eq:alastone?}
(\lambda+\epsilon_{i} + 2\rho, \lambda+\epsilon_{i})-(\lambda-\epsilon_{i} + 2\rho, \lambda-\epsilon_{i})=2N.
\end{equation}
If (\ref{eq:alastone?}) is satisfied, the notionally existing projection
$P[\lambda+\epsilon_{1}] \in End_{U_{q}^{(N)}(\mathfrak{g})} (V_{\lambda} \otimes V)$ that would act if it was well-defined as
$$P[\lambda+\epsilon_{1}] : 
V_{\lambda} \otimes V \rightarrow V_{\lambda+\epsilon_{1}} \subset V_{\lambda} \otimes V,$$
is only well-defined if $\lambda-\epsilon_{1} \notin {\cal{P}}^{+}_{\lambda} \cap \overline{\Lambda_{N}^{+}}$.
Additionally, if (\ref{eq:alastone?}) is satisfied, the notionally existing projection
$P[\lambda-\epsilon_{1}] \in End_{U_{q}^{(N)}(\mathfrak{g})} (V_{\lambda} \otimes V)$ that would act if it was well-defined as
$$P[\lambda-\epsilon_{1}] : 
V_{\lambda} \otimes V \rightarrow V_{\lambda-\epsilon_{1}} \subset V_{\lambda} \otimes V,$$
is only well-defined if $\lambda+\epsilon_{1} \notin {\cal{P}}^{+}_{\lambda} \cap \overline{\Lambda_{N}^{+}}$.

If $\lambda-\epsilon_{1}$ and $\lambda+\epsilon_{1}$ are elements of ${\cal{P}}^{0}_{\lambda}$ for any given
$\lambda \in \Lambda_{N}^{+}$, then they are elements of 
${\cal{P}}^{+}_{\lambda} \cap \overline{\Lambda_{N}^{+}}$.
By inspection, both $\lambda-\epsilon_{1}$ and $\lambda+\epsilon_{1}$ are elements of ${\cal{P}}^{0}_{\lambda}$ if and only
if the components of $\lambda$ satisfy  
$$\lambda_{2} \leq \lambda_{1}-1, \hspace{5mm}
\mbox{and} \hspace{5mm} (\lambda_{1} + 1) + \lambda_{2} \leq N-2n+2.$$  
The first condition implies that both  
$\lambda-\epsilon_{1}$ and $\lambda+\epsilon_{1}$ are elements of ${\cal{P}}^{+}_{\lambda}$ 
and the second condition is
necessary so that $\lambda+\epsilon_{1}$ is an element of $\overline{\Lambda_{N}^{+}}$.

Let us determine the conditions on the components of $\lambda$ so that (\ref{eq:alastone?}) is true.  
This equation states that $4\lambda_{1}+4n-2=2N$ which we rewrite as $\lambda_{1} = N/2-n+1/2$.
Consider an integral dominant weight $\lambda$ where 
$\lambda_{1} = N/2-n+1/2$ and $\lambda_{2} \leq \lambda_{1}-1$.  
Then $\lambda_{1} + \lambda_{2} \leq N-2n$ and $\lambda$ is in $\Lambda_{N}^{+}$.
This means that for each $\Lambda_{N}^{+}$ there is 
{\emph{at least one}} $\lambda \in \Lambda_{N}^{+}$ where the
notional projections $P[\lambda + \epsilon_{1}]$ and $P[\lambda-\epsilon_{1}]$ are 
{\emph{not well-defined}}.
In general there are many such elements of $\Lambda_{N}^{+}$ for which this is true.  This proves (i).

We now prove (ii).
Consider an element $\overline{\lambda} \in \Lambda_{N}^{+}$ where 
$\overline{\lambda}_{1} = N/2-n+1/2$ and $\overline{\lambda}_{2}=\overline{\lambda}_{1}$.  Then
$\overline{\lambda}-\epsilon_{1} \notin {\cal{P}}^{0}_{\overline{\lambda}}$ and 
$(\overline{\lambda}_{1}+1)+\overline{\lambda}_{2} = N-2n+2$, thus we have 
$\overline{\lambda}+\epsilon_{1} \in \overline{\Lambda_{N}^{+}} \backslash \Lambda_{N}^{+}$.  
Then the projection
$P\Big[\overline{\lambda}+\epsilon_{1}\Big]: V_{\overline{\lambda}} \otimes V \rightarrow 
V_{\overline{\lambda}+\epsilon_{1}}$ is well-defined as 
$\overline{\lambda}-\epsilon_{1} \notin {\cal{P}}^{0}_{\overline{\lambda}}$, proving (ii).
Note that $sdim_{q}(V_{\overline{\lambda}+\epsilon_{1}})=0$.

\end{proof}

\end{section}

\end{chapter}

%% file: chapter4100.tex
\begin{chapter}{Topological invariants of 3-manifolds from $U_{q}^{(N)}(osp(1|2n))$}
\label{chaptitle:topologicalinvariants}
\markboth{\text{Chapter \ref{chaptitle:topologicalinvariants}. Topological invariants of 3-manifolds}}
{\text{ }}

The structure of this chapter is as follows.
In Section \ref{chapt4:sec(Knots_and_Links)} we introduce knots and links 
and the equivalence relations of ambient and regular isotopy  
generated by the Reidemeister moves on planar projections of links.  
We describe how each closed, connected, orientable $3$-manifold $M_{L}$
can be obtained by {\emph{performing surgery on the $3$-sphere $S^{3}$ along a link $L \subset S^{3}$}}, 
and we describe the equivalence relations, called the {\emph{Kirby moves}}, 
on two links $L$ and  $L'$ embedded in $S^{3}$ such that
the $3$-manifolds they give rise to upon performing surgery are homeomorphic.

In Section \ref{sec:ribbongraphs} we introduce directed ribbon tangles 
and the category of coloured directed ribbon tangles
following Reshetikhin \cite{resh} and Reshetikhin and Turaev \cite{reshtur},
as a precursor to constructing isotopy invariants of coloured directed ribbon tangles, 
and thereby regular isotopy invariants of links later in this chapter.  

In Section \ref{sec:RTfunc} we state the Reshetikhin-Turaev functor $F$, 
a covariant functor from the category of
coloured directed ribbon tangles to the category of finite dimensional representations of a $\mathbb{Z}_{2}$-graded ribbon Hopf
algebra.  
This definition follows directly from \cite{z3} and is a generalisation of the covariant functor 
from the category of coloured directed ribbon tangles to the category of 
finite dimensional representations 
of an ungraded ribbon Hopf algebra \cite{resh,reshtur,rt}.  
The functor yields isotopy invariants of coloured directed ribbon tangles 
and thereby regular isotopy invariants of links.

Sections \ref{chapter4:sectionlabel(Pseudo-Modular_Hopf_algebras)}--\ref{sec:hopefullylast?}
contain the new results in this chapter.

In Section \ref{chapter4:sectionlabel(Pseudo-Modular_Hopf_algebras)} 
we generalise the definition of a modular Hopf algebra
to the $\mathbb{Z}_{2}$-graded case.   
Drawing on this generalisation, 
we define a new algebra that we call a {\emph{pseudo-modular Hopf algebra}}.
Pseudo-modular Hopf algebras are $\mathbb{Z}_{2}$-graded ribbon Hopf algebras  
together with a finite set of finite dimensional representations 
satisfying slightly weaker conditions than those satisfied by modular Hopf algebras.  
We define these algebras as there exist quotients of quantum algebras and quantum superalgebras that
are not modular, or are not known to be modular, from which topological
invariants of $3$-manifolds can be constructed.  
Later in this chapter
we construct topological invariants of $3$-manifolds from pseudo-modular Hopf algebras.

In Section \ref{sec:rushingwater} we prove that topological invariants of 
closed, connected, orientable $3$-manifolds can be constructed using pseudo-modular Hopf algebras.  

In Section \ref{chapter4:sectionlabel(Invariants_of_3-manifolds_arising_from)} we prove that
the $\mathbb{Z}_{2}$-graded ribbon Hopf algebra $U_{q}^{(N)}(osp(1|2n))$, where $N \geq 6$ satisfies
  $N \equiv 2 \pmod{4}$, 
together with a set of non-isomorphic finite dimensional $U_{q}^{(N)}(osp(1|2n))$-modules 
$\left\{ V_{\lambda} | \ \lambda \in \Lambda_{N}^{+} \right\}$ 
defined in Chapter \ref{chap2A:titlelabel},
is a pseudo-modular Hopf algebra.
We also prove that $3$-manifold invariants cannot be constructed from $U_{q}^{(N)}(osp(1|2n))$
when $N \geq 4$ satisfies $N \equiv 0 \pmod{4}$.

In Section \ref{chapter4sectionlabel(comparing_the_)} we 
compare the $3$-manifold invariants arising from
$U_{q}^{(N)}(osp(1|2n))$ when $N \geq 6$ satisfies $N \equiv 2 \pmod{4}$,
with the invariants arising from other ribbon Hopf algebras.  
This is difficult to do, and for tractability we only compare our invariants
with those arising from $U_{q}^{(N/2)}(so(2n+1))$ at the same $N$.  
We do this as the quantised universal enveloping
algebras of $osp(1|2n)$ and $so(2n+1)$ are known to be related at generic $q$.
We show that the $3$-manifold invariants arising from 
$U_{q}^{(N)}(osp(1|2n))$ and $U_{q}^{(N/2)}(so(2n+1))$ are 
{\emph{not}} the same.

In Section \ref{sec:hopefullylast?} we give some side results.

\begin{section}{Knots and links}
\markright{\text{Knots and links}}
\label{chapt4:sec(Knots_and_Links)}

We take the following definitions from \cite{r,cp,l}.
\begin{definition}
A link $L = \bigcup_{i=1}^{m} L_{i}$ of $m$ components embedded in $S^{3}$
is a subset of $S^{3}$ consisting of $m$ disjoint smooth 
1-dimensional submanifolds of $S^{3}$.  A knot is a link with one component.
\end{definition}
\noindent
This definition of a link ensures that we do not deal with wild links in this thesis.
\begin{definition}
An oriented link $L$ is a link with an orientation assigned to each connected component of $L$.
\end{definition}
\begin{definition}
Two oriented links $L$, $L' \subset S^{3}$ are said to be equivalent if there exists
an orientation preserving diffeomorphism $f$ of $S^{3}$ such that $f(L) = L'$
and such that $f$ takes the orientation of $L$ into that of $L'$.
\end{definition}

Let $S^{3}$ denote the $3$-sphere, then $S^{3}=\mathbb{R}^{3} \cup \{\infty\}$.
Let a link $L \subset S^{3}$ be embedded in $S^{3}$ such that 
\begin{itemize}
\item[(i)] $L \cap \{\infty\} = \emptyset$, and 
\item[(ii)] given the standard projection $p: \mathbb{R}^{3} \rightarrow \mathbb{R}^{2}$, 
$p(L)$ has only transversal crossings where in some sufficiently small neighbourhood
of each crossing there are projections of at most two {\emph{branches}} of $L$.
\end{itemize}
Here we define a branch of $L$ to be a closed proper subset of a single connected component of $L$.
We define the {\emph{planar projection of a link $L \subset \mathbb{R}^{3} \subset S^{3}$}} 
to be $p(L)$ together with
information at each crossing in $p(L)$ specifying which of the branches is the
overcrossing branch.  The overcrossing branch is
represented by an unbroken line and the undercrossing branch is represented by a
broken line, the line being broken at the crossing \cite{cp}.

We detail the {\emph{Reidemeister moves}} in Figure \ref{fig:reidmoves} \cite{kau}.
Let $L \subset S^{3}$ be a link such that $L \cap \{\infty\}=\emptyset$.
Each of the Reidemeister moves replaces a configuration of arcs and crossings in the
intersection of $p(L)$ with a $2$-disc $D^{2}$, with 
another collection of arcs and crossings, such that the complement of
$p(L)$ is left unchanged \cite{l}.  
Each Reidemeister move is an equivalence relation: 
two links $L, L' \subset S^{3}$ are said to be {\emph{ambient isotopic}} if 
their planar projections are elements of the
equivalence class of planar projections of links generated by the
Reidemeister I, II and III moves.
Two links $L$, $L' \subset S^{3}$ are said to be 
\emph{regularly isotopic} if their planar projections are elements of the equivalence class generated by the
Reidemeister II and III moves.
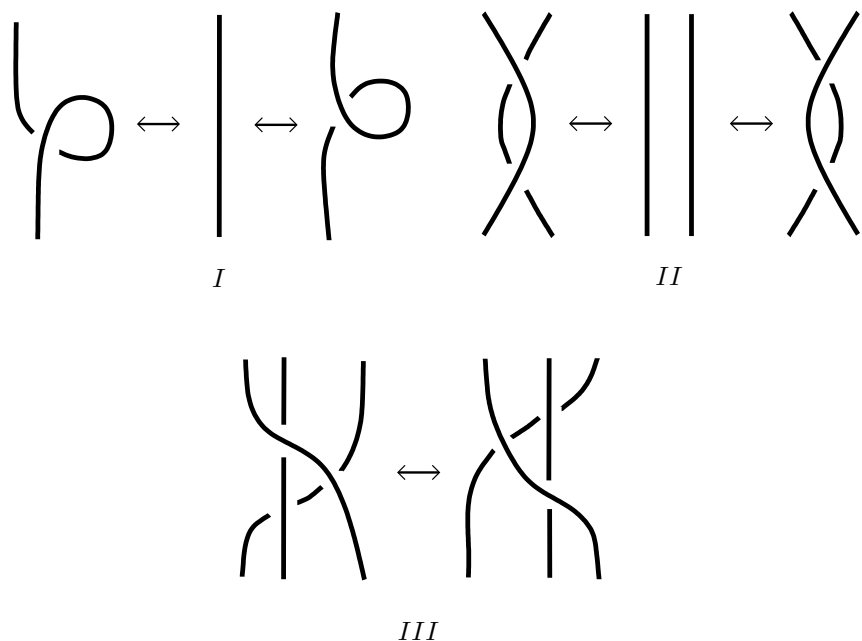
\begin{figure}[hbt]
\begin{center}
  \input{reid100.pstex_t}
\caption{The Reidemeister I, II and III moves}  \label{fig:reidmoves}
\end{center}
\end{figure}

We associate a linking number $+1$ or $-1$ with each crossing of any pair of components of an oriented link as given in 
Figure \ref{fig:linkingnumbers}.
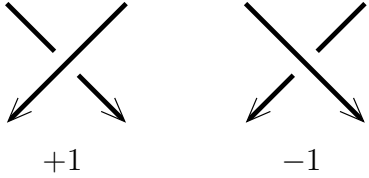
\begin{figure}[hbt]
\begin{center}
  \input{linkingnumbers100.pstex_t}
\caption{Linking numbers}  \label{fig:linkingnumbers}
\end{center}
\end{figure}
\begin{definition}
Let $L_{i}$ and $L_{j}$ be connected components of an oriented link $L$.  
The linking number $\mathrm{lk}(L_{i},L_{j})$ of $L_{i}$ and $L_{j}$ 
is half of the sum of the linking numbers
associated with each crossing of $L_{i}$ and $L_{j}$ in a planar projection of $L$.
\end{definition}
\begin{definition}
Let $L \subset S^{3}$ be an unoriented link.  
The writhing number $w(L_{i})$ of a connected component $L_{i}$ of $L$
is the sum of the linking numbers associated with each crossing of $L_{i}$ with itself in a planar projection of $L$, 
where $L$ is assigned an arbitrary orientation.
\end{definition}
The writhing number of a connected component of a link is independent of the orientation of the link, 
thus the writhing number of a connected component of an unoriented link is 
well-defined using the above definition \cite{kau}.

\begin{definition}
A framed link $L \subset S^{3}$ is a link $L$ 
together with an assignment of integers, each connected component of $L$
being assigned an integer which we call the framing number of that component.
\end{definition}
\noindent

Let $L = \bigcup_{i=1}^{m} L_{i} \subset S^{3}$ 
be a framed unoriented link with $m$ connected components where the framing number of 
the connected component $L_{i}$ is $n_{i}$.  We can equip $L_{i}$ with a normal vector field as follows:
let $U_{i} \subset S^{3}$ be a small tubular neighbourhood of $L_{i}$ such that $U_{i} \cap U_{j} = \emptyset$ if $i \neq j$.
Let $K_{i} \subset U_{i}$ be a connected link such that we can equip $L_{i}$ with a normal
vector field where the tips of the vectors trace out $K_{i}$ and such that if we assign an orientation to $L_{i}$ and a
parallel orientation to $K_{i}$, the linking number of $L_{i}$ and $K_{i}$ is $n_{i}$ \cite{kau}.
For each such normal vector field there exists an integer $n_{i}$, 
and for each integer $n_{i}$ there exists such a normal vector field.

A natural normal vector field to use gives the so-called {\emph{blackboard framing}}.  
Set the normal vector field to $L_{i}$ to lie in the planar projection of $L$.  
Then the tips of the vectors sweep out a link $K_{i}$ parallel to $L_{i}$ \cite{kau}
 and the linking number of $L_{i}$ and $K_{i}$ is $w(L_{i})$.
The normal vector field of a framed link can always be presented in the blackboard framing.

The blackboard framing allows us to define a convenient notion of `equivalence' of framed links. 
We say that two framed links are equivalent if the links are regularly isotopic when presented 
in the blackboard framing.

An important element in the study of 3-manifolds is the notion of surgery \cite[p. 252]{kau}:
\begin{definition}{Performing surgery on a 3-manifold $M^{3}$ along a link 
$L \subset M^{3}$.}

Let $M^{3}$ be a 3-manifold and $L \subset M^{3}$ an unoriented link with one
component and framing number $n_{L}$.  Let $\alpha: S^{1} \times D^{2} \rightarrow
M^{3}$ be an embedding and $\alpha(S^{1} \times 0)$ an embedding of $L$.  A
longitude of $\alpha(S^{1} \times D^{2})$ is $\alpha(S^{1} \times 1) \subset
\alpha\big(\partial(S^{1} \times D^{2})\big)$ where $\partial(S^{1} \times D^{2})$ 
denotes the boundary of $S^{1} \times D^{2}$.
Let $L' \subset \alpha\big(\partial(S^{1} \times D^{2})\big)$ be a closed twisted
longitude of $\alpha(S^{1} \times D^{2})$ such that $\mathrm{lk}(L,L') = n_{L}$ where we now fix $L$
and $L'$ to have the same orientation.  

We perform surgery on $M^{3}$ along $L$ by gluing $D^{2} \times S^{1}$ to 
$M^{3}-Int\big(\alpha(S^{1} \times D^{2})\big)$ along the boundary of 
$S^{1} \times S^{1} \subset M^{3}-Int\big(\alpha(S^{1} \times D^{2})\big)$ 
so that the meridian $m=S^{1} \times 1
\subset D^{2} \times S^{1}$ matches the closed twisted longitude $L'$.
We conduct surgery on links with more than one component by performing surgery on each
component simultaneously.
\end{definition}

The following theorem \cite{l1} (see also \cite{l}) establishes a connection between
closed, connected, orientable 3-manifolds and links embedded in $S^{3}$ that 
underpins our approach to developing topological invariants of such 3-manifolds.
\begin{theorem} 
Any closed, connected, orientable 3-manifold $M_{L}$ can be obtained by 
performing surgery on $S^{3}$ along a framed link $L \subset S^{3}$.
\end{theorem}
We refer to \cite[Chap. 9]{r} and \cite[Part I, Chap. 16]{kau} 
for detailed discussion concerning, and examples of, surgery on $S^{3}$ along framed links.

Let $L$, $L' \subset S^{3}$ be two regularly isotopic links with the same framing.
Let $M_{L}$ (resp. $M_{L'}$) be the closed, connected, orientable 3-manifold obtained by performing surgery on $S^{3}$
along $L$ (resp. $L'$), then there exists an orientation preserving homeomorphism 
$M_{L} \rightarrow M_{L'}$.

A question which immediately arises is 
whether there are any other relations between framed links such that the $3$-manifolds they give rise to 
upon performing surgery
are related by an orientation-preserving homeomorphism.  
Kirby answered that question by finding a
complete collection of such relations between framed links and proving
that closed, connected, orientable 3-manifolds $M_{L}$ and $M_{L'}$ 
are related by an orientation-preserving
homeomorphism if and only if the framed links $L$ and $L'$ are related by a 
certain set of transformations \cite{k}.

Fenn and Rourke subsequently proved a similar result but with the advantage 
that their transformations were \emph{local} transformations, 
that is, their transformations occur on some subset of the planar projection of a link \cite{fr}.
These transformations are called the {\emph{Kirby moves}} 
and are detailed in Figures \ref{fig:Kz9}--\ref{fig:Kz7}.
Let $L$, $L' \subset S^{3}$ be two links that are regularly isotopic 
outside of their intersections with a $3$-disc $D^{3}$.  
Each Kirby move in Figures \ref{fig:Kz9}--\ref{fig:Kz7} relates the intersection of 
a $2$-disc $D^{2}$ with $p(L)$ and $p(L')$, respectively.  

We say that $L$ and $L'$ are equivalent under the
Kirby moves if the intersection of $p(L)$ with $D^{2}$ is regularly isotopic to one diagram
of a move, and the intersection of $p(L')$ with $D^{2}$ 
is regularly isotopic to the other diagram of the move.
The relation {\emph{equivalent under the Kirby moves}} is an equivalence relation.
The reader is referred to \cite[Part I, Chap. 16]{kau}
and \cite[Chap. 12]{kl} for examples of links equivalent under the Kirby moves.
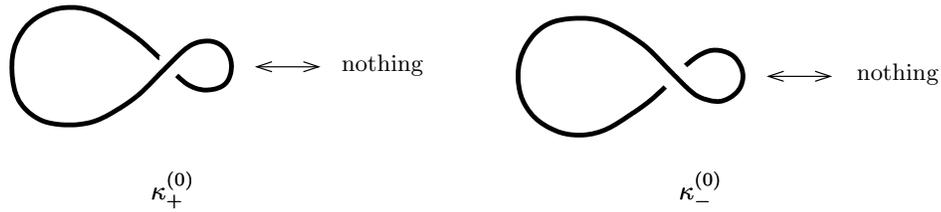
\begin{figure}[hbt]
\begin{center}
 \input{kirby1100.pstex_t}
\caption{The Kirby $\kappa_{+}^{(0)}$ and $\kappa_{-}^{(0)}$ moves}  \label{fig:Kz9}
\end{center}
\end{figure}
\begin{figure}[hbt]
\begin{center}
 \input{Kirby2100.pstex_t}
\caption{The Kirby $\kappa_{+}$ move}   \label{fig:Kz8}
\end{center}
\end{figure}
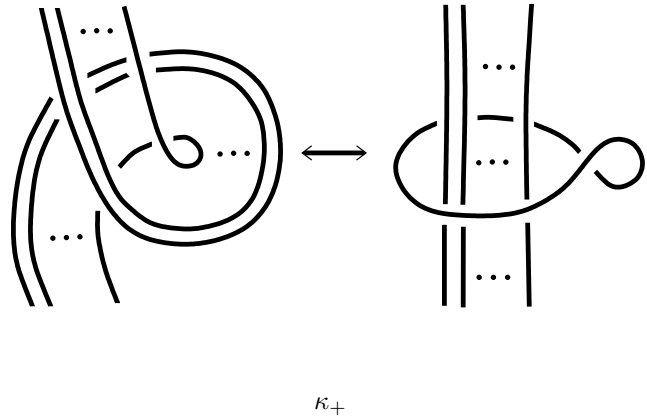
\begin{figure}[hbt]
\begin{center}
 \input{Kirby3100.pstex_t}
\caption{The Kirby $\kappa_{-}$ move}   \label{fig:Kz7}
\end{center}
\end{figure}

The $\kappa^{(0)}_{+}$ and $\kappa^{(0)}_{-}$ moves are called 
the {\emph{special Kirby (+) move}}, and the {\emph{special Kirby (-) move}}, respectively.
The $\kappa_{+}$ and $\kappa_{-}$ moves are called the 
{\emph{Kirby (+) move}}, and the {\emph{Kirby (-) move}}, respectively.

We now explain the origin of the Kirby moves \cite{kau,kl,l}.  
Let $L \subset S^{3}$ be a framed link and $L' \subset S^{3}$ an unknot with framing number $\pm 1$ 
such that $L \cup L'$ is a split link, that is, that $L \cup L'$ is the disjoint union
of links $L$, $L'$ such that $L$ and $L'$ are mutually unlinked.
Surgery on $S^{3}$ along $L'$ gives rise to a 3-manifold homeomorphic to $S^{3}$ \cite{kl}.  
Let $M_{L}$ be a closed, connected, orientable 3-manifold obtained by performing surgery on $S^{3}$ along 
the framed link $L \subset S^{3}$.  
Let $M_{L} \# S^{3}$ denote the connected sum of $M_{L}$ and $S^{3}$ and let
$\cong$ denote an orientation-preserving homeomorphism.  Then
$M_{L} \# S^{3} \cong M_{L}$, and the relations $M_{L \cup L'} \cong M_{L} \# M_{L'} \cong M_{L}$ imply the
existence of the $\kappa_{+}^{(0)}$ and $\kappa_{-}^{(0)}$ moves.

The $\kappa_{+}$ and $\kappa_{-}$ moves arise from 
an application of the {\emph{band connected sum move}} as detailed in the following theorem \cite{kau}.
\begin{theorem}
Let $K, K' \subset S^{3}$ be two links related by a band connected sum move and 
let $M_{K}, M_{K'}$ be the closed, connected, orientable 3-manifolds 
obtained by performing surgery on $S^{3}$ along $K$ and $K'$, respectively. 
Then $M_{K} \cong M_{K'}$.
\end{theorem}

We now explain the band connected sum move.
Let $L_{i}$ and $L_{j}$ be two disjoint connected components of a link $L \subset S^{3}$. 
The action of the band connected sum move on $L_{i}$ is to replace $L_{i}$ with the 
band connected sum $L_{i} \#_{b} L_{j}$, which we define as follows. 
Let $b$ be the image of a smooth embedding 
$e: [0,1] \times [0,1] \rightarrow S^{3}$
where $e\big([0,1] \times \{0\}\big)$ is glued
smoothly to $L_{i}$ and $e\big([0,1] \times \{1\}\big)$ is glued smoothly to $L_{j}$ 
such that the linking number of $e\big(\{0,1\} \times [0,1]\big)$ is zero, 
and furthermore, we require that the only part of $b$ that intersects $L_{i}$ (resp. $L_{j}$) is
$e\big([0,1] \times \{0\}\big)$ (resp. $e\big([0,1] \times \{1\}\big)$).
We define
$$L_{i} \#_{b} L_{j} = \big( L_{i} \cup b \cup L_{j} \big)  \backslash  e\big((0,1) \times [0,1] \big).$$
We now detail an example of the band connected sum move to be concrete. 
Let $L_{i}$ and $L_{j}$ be components of a split link in Figure \ref{fig:band1}, then $L_{i}
\#_{b} L_{j}$ is the left most component in Figure \ref{fig:band2}.
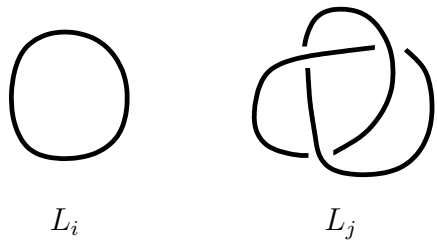
\begin{figure}[hbt]
\begin{center}
 \input{band1100.pstex_t}
\caption{Link components $L_{i}$ and $L_{j}$}   \label{fig:band1}
\end{center}
\end{figure}
\begin{figure}[hbt]
\begin{center}
 \input{band2100.pstex_t}
\caption{The links $L_{i} \#_{b} L_{j}$ and $L_{j}$}   \label{fig:band2}
\end{center}
\end{figure}
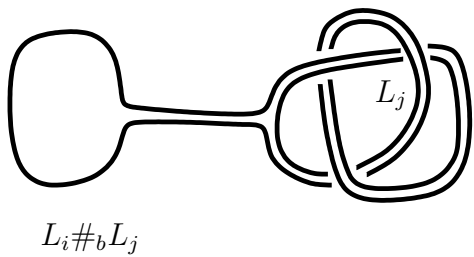

The band connected sum move gives rise to the $\kappa_{+}$ move.  
In Figure \ref{fig:ohyes} we show this in the situation that the left hand side of the
$\kappa_{+}$ move has one component.  The multicomponent case is similar.  
In Figure \ref{fig:ohyes}, $\stackrel{b}{\longleftrightarrow}$ 
indicates the band connected sum move, $\stackrel{k}{\longleftrightarrow}$
indicates the $\kappa^{(0)}_{+}$ move and the remaining move is regular isotopy.  
The generation of the $\kappa_{-}$ move is similar.
\begin{figure}[hbt]
\begin{center}
 \input{bandkir100.pstex_t}
\caption{The derivation of the $\kappa_{+}$ move}   \label{fig:ohyes}
\end{center}
\end{figure}
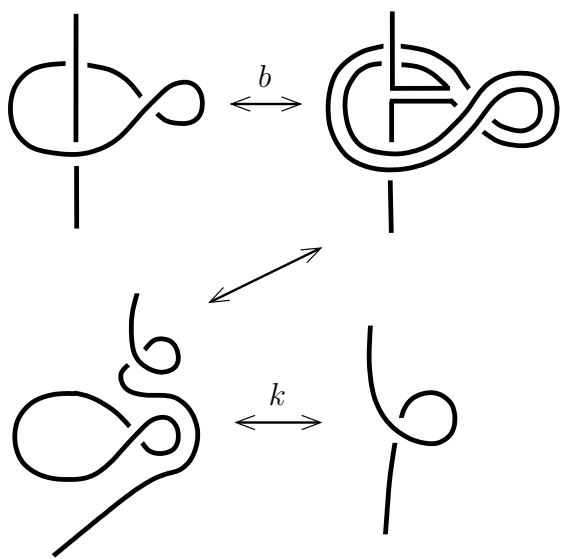

The existence of the band connected sum move suggests that there may be a large number of 
moves on links that give rise to homeomorphic 3-manifolds.  
However, a convenient result is that 
the Kirby moves are a sufficient generating set of such moves \cite[Thm. p. 1]{fr}:
\begin{theorem}
Orientation preserving homeomorphism classes of closed, connected, orientable
3-manifolds correspond bijectively to equivalence classes of framed links in
$S^{3}$ where the equivalence is generated by the Kirby moves.
\end{theorem}
It transpires that the four Kirby moves are not a minimal generating set of the
Kirby calculus: the $\kappa^{(0)}_{\pm}$ moves in addition to
either of the $\kappa_{+}$ or $\kappa_{-}$ moves generates the entire
Kirby calculus \cite{kau}.

We will use all of these facts to create topological invariants of 
closed, connected, orientable $3$-manifolds later in this thesis.
Let $M_{L}$ be a closed, connected, orientable 3-manifold
obtained by performing surgery on $S^{3}$ along a framed link $L \subset S^{3}$.
We will create a topological invariant of $M_{L}$ by taking such sums of regular isotopy 
invariants of $L$ as are unchanged after applying each of the Kirby moves to $L$.
A convenient way to study regular isotopy invariants of links is to
consider isotopy invariants of {\emph{ribbon tangles}}.

\end{section}

\begin{section}{Tangles and ribbon tangles}
\label{sec:ribbongraphs}
\markright{\text{Tangles and ribbon tangles}}

We now examine tangles, directed ribbon tangles and coloured directed ribbon tangles
before constructing isotopy invariants of directed ribbon tangles.
The work in this section is well known: it
first appeared in \cite{resh} and then in
\cite{reshtur, rt} and it often appears and is referenced in the literature.  
We state it here to be complete.

\begin{subsection}{The category of directed ribbon tangles}
\label{subsec:catribtan}

The basic element in this work is the {\emph{ribbon}}.  
A ribbon is defined to be a square $[0,1] \times [0,1]$ smoothly embedded
in $\mathbb{R}^{3}$ \cite{rt}, and  the images of the segments $[0,1] \times \{0\}$ 
and $[0,1] \times \{1\}$ under the embeddings are the {\emph{bases}} of the ribbon.  
We call the image of $\{1/2\} \times [0,1]$ the {\emph{core}} of the ribbon.  

We call the image of a cylinder $S^{1} \times [0,1]$ smoothly embedded in
$\mathbb{R}^{3}$ an {\emph{annulus}},   
and the image of $S^{1} \times \{1/2\}$ the {\emph{core}} of the annulus.

We define the {\emph{writhe}} of a ribbon to be the linking
number of the image of $\{0\} \times [0,1]$ and the image of $\{1\} \times [0,1]$ 
 where we assign the edges of the ribbon parallel orientations.

Ribbons and annuli are orientable surfaces in $\mathbb{R}^{3}$ and  
an orientation of a ribbon or annulus is equivalent to a choice of one side of 
the ribbon or annulus.  

We say that a ribbon or annulus is directed if its core is oriented, and
we orient the core of a ribbon by labelling
one of the bases of the ribbon the {\emph{initial base}} and the other base the {\emph{final base}}.  
Each ribbon and annulus can be directed in two ways and oriented in two ways.  

For $k,l \in \mathbb{Z}_{+}$ we define a $(k,l)$-ribbon tangle $\Gamma$ to be
the union of a finite number of disjoint oriented ribbons and annuli
embedded in $\mathbb{R}^{2} \times [0,1]$ such that $\Gamma$ satisfies:
$$\Gamma \cap \left(\mathbb{R}^{2} \times \{1\}\right) =
\big\{ [i-1/4,i+1/4] \times \{0\} \times \{1\} \ | \ i=1,2,\ldots,k  \big\},$$
$$\Gamma \cap \left(\mathbb{R}^{2} \times \{0\}\right) =
\big\{ [j-1/4,j+1/4] \times \{0\} \times \{0\} \ | \ j=1,2,\ldots,l  \big\},$$
and such that the same side of each ribbon faces the reader at
$\Gamma \cap(\mathbb{R}^{2} \times \{1\})$ and 
$\Gamma \cap(\mathbb{R}^{2} \times \{0\})$.

We define the $(k,l)$-tangle associated with a particular $(k,l)$-ribbon tangle $\Gamma$
to be the union of the cores of all the components of $\Gamma$.
This corresponds to our intuitive notion of a $(k,l)$-tangle.

So far we have not imposed any directions on the (ribbon) tangles.
In Reshetikhin and Turaev's isotopy invariants of ribbons, one
uses oriented and directed ribbon tangles, and so we now define
a directed $(k,l)$-ribbon tangle.

We define a {\emph{directed $(k,l)$-ribbon tangle}} $\Gamma$ to be a $(k,l)$-ribbon tangle where
all the ribbons and annuli of $\Gamma$ are directed.  
We can intuitively think of this by drawing an
arrow on each ribbon and annulus parallel to the core of the ribbon or annulus.

For each directed $(k,l)$-ribbon tangle $\Gamma$ there are two 
sequences which specify the directions of the ribbons of $\Gamma$:
\begin{eqnarray*}
\epsilon^{*}(\Gamma) & = & (\epsilon^{1},\epsilon^{2},\ldots,\epsilon^{k}), \\
\epsilon_{*}(\Gamma) & = & (\epsilon_{1},\epsilon_{2},\ldots,\epsilon_{l}).
\end{eqnarray*}
Here each $\epsilon^{i}$, $\epsilon_{i}$ is either $-1$ or $+1$. 
Intuitively, $\epsilon^{i}=+1$ (resp. $\epsilon^{i}=-1$) if the arrow on the 
$i^{th}$ ribbon from the left at the top of the ribbon tangle is pointing down (resp. up),
and $\epsilon_{i}=+1$ (resp. $\epsilon_{i}=-1$) if the arrow on the 
$i^{th}$ ribbon from the left at the bottom of the ribbon tangle is pointing down (resp. up).
Technically, we set $\epsilon^{i}=+1$ if 
$$[i-1/4,i+1/4] \times \{0\} \times \{1\}$$ 
is an initial base and  $\epsilon^{i} =-1$ if it is a final base.
Similarly, we set $\epsilon_{j} = +1$ if  
$$[j-1/4,j+1/4] \times \{0\} \times \{0\}$$ is a final 
base and $\epsilon_{j}=-1$ if it is an initial base.

We depict directed ribbon tangles as the disjoint union of ribbons 
and annuli, the directions of the ribbons and annuli 
denoted by arrows drawn on them.
On ribbons we point the arrows in the direction of the final base.
An example of this is  in Figure \ref{fig:orangesandapples}.
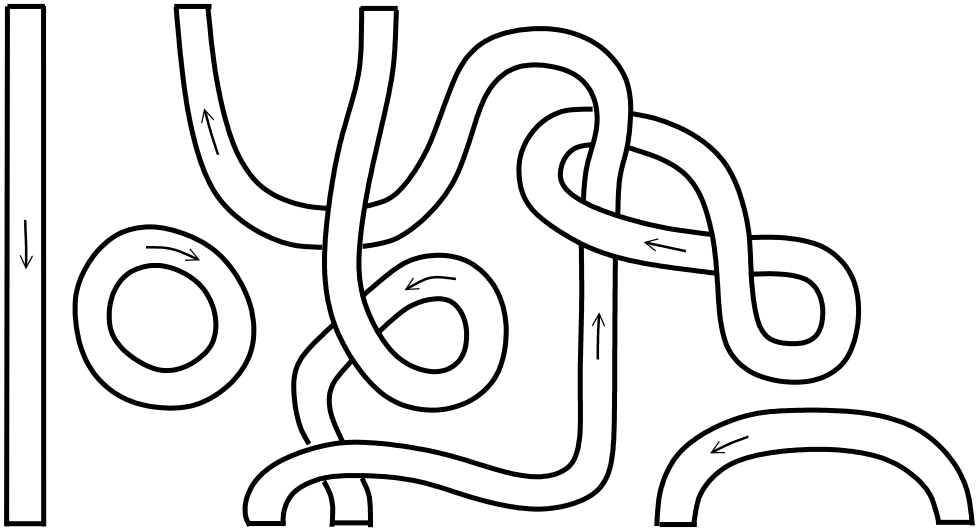
\begin{figure}[hbt]
\begin{center}
  \input{example100.pstex_t}
\caption{An example of directed ribbons and annuli}  \label{fig:orangesandapples}
\end{center}
\end{figure}

Two directed $(k,l)$-ribbon tangles $\Gamma$, $\Gamma'$ are said to be isotopic if there
exists a smooth isotopy
$$h_{t}: \mathbb{R}^{2} \times [0,1] \rightarrow \mathbb{R}^{2} \times [0,1], \hspace{10mm} t \in [0,1],$$ 
of the identity $h_{0} = \mathrm{id}$ such that each $h_{t}$ is
a  diffeomorphism of the strip $\mathbb{R}^{2} \times [0,1]$ fixing its boundary 
$\mathbb{R}^{2} \times \{0,1\}$ and $h_{1}$
transforms $\Gamma$ into $\Gamma'$ preserving the decomposition into ribbons and
annuli, the directions of cores and the orientations of the ribbons and
annuli \cite{rt}.  Isotopy is an equivalence relation.

In referring to a particular directed $(k,l)$-ribbon tangle $\Gamma$ we are
referring to any element of the equivalence class of directed $(k,l)$-ribbon tangles
containing $\Gamma$ where the equivalence is generated by isotopy.

If $\Gamma$ is a directed $(k,m)$-ribbon tangle and $\Gamma'$ is a directed
$(m,n)$-ribbon tangle then we can vertically compose them to produce a new 
directed $(k,n)$-ribbon tangle $\Gamma \circ \Gamma'$ if the directions of $\Gamma$ and $\Gamma'$ are
compatible.  
By `compatible' we mean that $\epsilon_{*}(\Gamma) = \epsilon^{*}(\Gamma')$, or
more intuitively, the `up, down' sequence of arrows at the bottom of $\Gamma$ is the same as the `up, down'
sequence of arrows at the top of $\Gamma'$.

We technically define $\Gamma \circ \Gamma'$ 
to be the directed $(k,n)$-ribbon tangle obtained as follows. 
Take the union $\Gamma \cup (\Gamma'+t)$ where $t$ is the vector $t=(0,0,-1) \in \mathbb{R}^{3}$ 
such that the bases of the ribbons are smoothly glued together.  
Then $\Gamma \circ \Gamma'$ is
the image of $\Gamma \cup (\Gamma'+t)$ under the map
$(x, y, z) \mapsto (x, y, (z+1)/2)$, for all $(x, y, z) \in \mathbb{R}^{2} \times [-1, 1]$, so that we
have $\Gamma \circ \Gamma' \subset \mathbb{R}^{2} \times [0,1]$. 
Then $\Gamma \circ \Gamma'$ is a directed $(k,n)$-ribbon tangle with two sequences encoding the
directions of the ribbons:
$$\epsilon^{*}(\Gamma \circ \Gamma') = \epsilon^{*}(\Gamma),$$
$$\epsilon_{*}(\Gamma \circ \Gamma') = \epsilon_{*}(\Gamma'),$$
as expected. 
We heuristically depict $\Gamma \circ \Gamma'$ in Figure \ref{fig:verticalcomp}.
\begin{figure}[hbt]
\begin{center}
  \input{verticcomp100.pstex_t}
\caption{The vertical composition $\Gamma \circ \Gamma'$}  \label{fig:verticalcomp}
\end{center}
\end{figure}
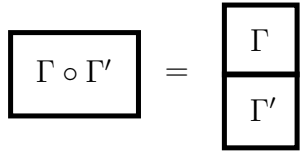

Let $\Gamma$ be a directed $(k,l)$-ribbon tangle and $\Gamma'$ a directed $(m,n)$-ribbon tangle, then 
the horizontal composition $\Gamma \otimes \Gamma'$ is always defined.
This $\Gamma \otimes \Gamma'$ is a directed $(k+m,l+n)$-ribbon tangle $\Gamma \otimes \Gamma'$ 
depicted by placing $\Gamma'$ immediately to the right of $\Gamma$
such that $\Gamma$ and $\Gamma'$ are mutually unlinked.  
We heuristically depict $\Gamma \otimes \Gamma'$ in Figure \ref{fig:horizontcomp}.
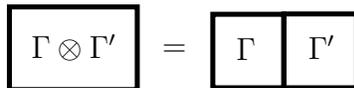
\begin{figure}[hbt]
\begin{center}
  \input{horizontcomp100.pstex_t}
\caption{The horizontal composition $\Gamma \otimes \Gamma'$}  \label{fig:horizontcomp}
\end{center}
\end{figure}

Now let $\Gamma_{1}$, $\Gamma_{1}'$, $\Gamma_{2}$ and $\Gamma_{2}'$ be directed ribbon tangles where 
$\Gamma_{1} \circ \Gamma_{1}'$ and $\Gamma_{2} \circ \Gamma_{2}'$ exist, then 
$$(\Gamma_{1} \circ \Gamma_{1}') \otimes (\Gamma_{2} \circ \Gamma_{2}') = 
(\Gamma_{1} \otimes \Gamma_{2}) \circ (\Gamma_{1}' \otimes \Gamma_{2}').$$
Note that $(\Gamma_{1} \otimes \Gamma_{2}) \circ (\Gamma_{1}' \otimes \Gamma_{2}')$ 
may exist even if $\Gamma_{1} \circ \Gamma_{1}'$ and $\Gamma_{2} \circ \Gamma_{2}'$ do not.

We now introduce the category of directed ribbon tangles {\bf{drib}}  \cite{rt,cp}. 
The objects of {\bf{drib}} are finite sequences:
$$\epsilon = (\epsilon_{1},\epsilon_{2},\ldots,\epsilon_{k}),$$ 
where $\epsilon_{i} \in \{-1,+1\}$ for each $i = 1, 2, \ldots, k$.  
We can think of an object $\epsilon$ of {\bf{drib}} as being 
$\epsilon^{*}(\Gamma)$ or $\epsilon_{*}(\Gamma)$
for some directed ribbon tangle $\Gamma$.
Isotopy equivalence classes of directed ribbon tangles are the morphisms of {\bf{drib}}.

If $\Gamma$ and $\Gamma'$ are two directed ribbon tangles where
$\Gamma \circ \Gamma'$ exists, and 
$f: \epsilon^{*}(\Gamma) \rightarrow \epsilon_{*}(\Gamma)$ and 
$f':\epsilon^{*}(\Gamma') \rightarrow \epsilon_{*}(\Gamma')$ are two morphisms of {\bf{drib}}, then  
we define 
$f \circ f': \epsilon^{*}(\Gamma \circ \Gamma') \rightarrow \epsilon_{*}(\Gamma \circ \Gamma')$ 
to be  $\Gamma \circ \Gamma'$.

Let 
$$\epsilon_{a} = (\epsilon_{a_{1}},\epsilon_{a_{2}},\ldots,
\epsilon_{a_{k}}), \hspace{10mm} \epsilon_{b} = (\epsilon_{b_{1}},\epsilon_{b_{2}},\ldots,
\epsilon_{b_{l}}),$$
be two objects of {\bf{drib}}.  
The tensor product of $\epsilon_{a}$ and $\epsilon_{b}$ 
is an object $\epsilon_{a} \otimes \epsilon_{b}$:
$$\epsilon_{a} \otimes \epsilon_{b} = 
(\epsilon_{a_{1}},\epsilon_{a_{2}},\ldots, \epsilon_{a_{k}},\epsilon_{b_{1}},\epsilon_{b_{2}},\ldots, \epsilon_{b_{l}}).$$

Now let $f:\epsilon_{a} \rightarrow \epsilon_{a}'$ and $g:\epsilon_{b} \rightarrow
\epsilon_{b}'$ be two morphisms of {\bf{drib}} that are the directed ribbon tangles 
$\Gamma$ and $\Gamma'$, respectively.  
The tensor product of $f$ and $g$ is a morphism $f \otimes g$:
$$f \otimes g: \epsilon_{a} \otimes \epsilon_{b} \rightarrow \epsilon_{a}' \otimes \epsilon_{b}',$$ 
which is just the directed ribbon tangle $\Gamma \otimes \Gamma'$.

Let $f$, $f'$, $g$ and $g'$ be morphisms where  $f \circ f'$ and 
$g \circ g'$ exist, then 
$$(f \circ f') \otimes (g \circ g') = (f \otimes g) \circ (f' \otimes g').$$
Note that $(f \otimes g) \circ (f' \otimes g')$ 
may exist even if each of $f \circ f'$ and $g \circ g'$ do not exist.

A directed ribbon tangle can be expressed as some combination of 
vertical and horizontal compositions of the directed ribbon
tangles in Figure \ref{fig:ribbontangleatoms} \cite{rt}.
For convenience we call the directed ribbon tangles in Figure \ref{fig:ribbontangleatoms} 
the {\emph{directed ribbon tangle atoms}}.
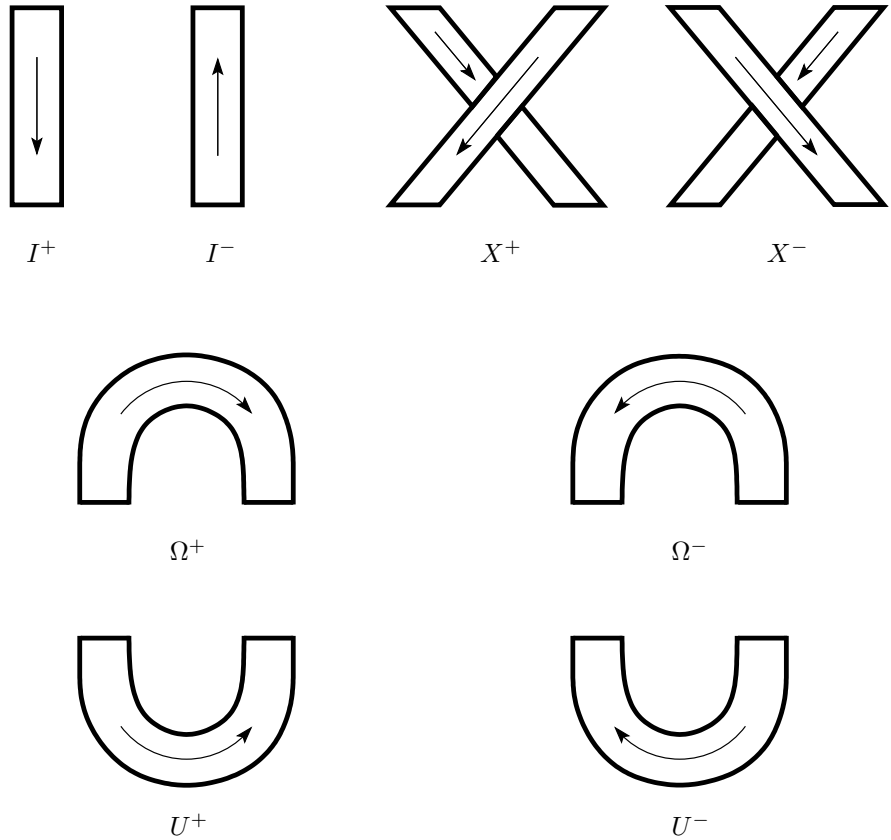
\begin{figure}[hbt]
\begin{center}
  \input{ribbon1100.pstex_t}
\caption{The directed ribbon tangle atoms}  \label{fig:ribbontangleatoms}
\end{center}
\end{figure}

\end{subsection}

\begin{subsection}{Coloured directed ribbon tangles}
\label{subsec:johnfaulkner(9991)}

We now introduce {\emph{coloured directed ribbon tangles}}, which we can intuitively think of as
directed ribbon tangles where each
component of the ribbon tangle has been `coloured' with an element of some set.

Let $S=\{s_{1},s_{2},\ldots,s_{t}\}$ be a non-empty finite set and let $\Gamma$ be a directed
$(k,l)$-ribbon tangle with $m$ disjoint directed ribbons:
$$\Gamma = \bigcup_{i=1}^{m}\Gamma_{i}.$$
We say that $\Gamma_{i}$ is coloured with $s_{a_{i}} \in S$ if we associate 
$s_{a_{i}}$  with $\Gamma_{i}$.
Let each ribbon $\Gamma_{i}$ be coloured with $s_{a_{i}}$, then we say that
$\Gamma$ is coloured with $(s_{a_{1}},s_{a_{2}},\ldots,s_{a_{m}}) \in S^{\times m}$.
We denote any involution $*$ of $S$ by $*(s_{a_{i}}) = (s_{a_{i}})^{*}$.

Two coloured directed $(k,l)$-ribbon tangles $\Gamma, \Gamma'$ are said to be isotopic if
they are isotopic as directed $(k,l)$-ribbon tangles and if the isotopy takes the colouring of
$\Gamma$ into the colouring of $\Gamma'$ \cite{rt}.

To each coloured directed $(k,l)$-ribbon tangle $\Gamma$
we associate two sequences encoding the colourings and directions of its ribbons:
\begin{eqnarray*}
(X^{*},\epsilon^{*})_{\Gamma} & = & \big( (i^{1},\epsilon^{1}),(i^{2},\epsilon^{2}),\ldots, (i^{k},\epsilon^{k}) \big),  \\
(X_{*},\epsilon_{*})_{\Gamma} & = & \big( (i_{1},\epsilon_{1}),(i_{2},\epsilon_{2}),\ldots,(i_{l},\epsilon_{l}) \big),
\end{eqnarray*}
where $i^{j}, i_{j} \in S$ and $\epsilon^{j}, \epsilon_{j} \in \{-1,+1\}$ for each $j$.  
Intuitively, the sequences 
$\epsilon^{*}(\Gamma) = (\epsilon^{1}, \epsilon^{2}, \ldots, \epsilon^{k})$ and
$\epsilon_{*}(\Gamma) = (\epsilon_{1}, \epsilon_{2}, \ldots, \epsilon_{l})$ are just the corresponding
sequences if we think of $\Gamma$ as uncoloured.
The sequences 
$$X^{*}({\Gamma}) = (i^{1}, i^{2}, \ldots, i^{k}), \hspace{10mm}
X_{*}({\Gamma}) = (i_{1}, i_{2}, \ldots, i_{l}),$$ 
encode the colourings of the ribbons, such that
$i^{j}$ (resp. $i_{j}$)
is the element of $S$ that colours the component of $\Gamma$ with non-empty intersection with
$[j-1/4, j+1/4] \times \{0\} \times \{1\}$  (resp. $[j-1/4, j+1/4] \times \{0\} \times \{0\}$).

The vertical composition $\Gamma \circ \Gamma'$ of a
coloured directed $(k,m)$-ribbon tangle $\Gamma$ and a
coloured directed $(m,n)$-ribbon tangle $\Gamma'$
exists if the directions and colourings of $\Gamma$ and 
$\Gamma'$ are compatible, that is if $(X_{*},\epsilon_{*})_{\Gamma} = (X^{*},\epsilon^{*})_{\Gamma'}$. 
We technically define $\Gamma \circ \Gamma'$ 
to be the coloured directed $(k,n)$-ribbon tangle with the underlying directed ribbon tangle to be the
vertical composition of $\Gamma$ and $\Gamma'$ and with two
sequences encoding the colouring and directions of the ribbons:
$$(X^{*},\epsilon^{*})_{\Gamma \circ \Gamma'} = (X^{*},\epsilon^{*})_{\Gamma}, \hspace{10mm}
(X_{*},\epsilon_{*})_{\Gamma \circ \Gamma'} = (X_{*},\epsilon_{*})_{\Gamma'}.$$

The horizontal composition $\Gamma \otimes \Gamma'$
of a coloured directed $(k,l)$-ribbon tangle $\Gamma$ and a
coloured directed $(m,n)$-ribbon tangle $\Gamma'$ is always defined to be
the coloured directed $(k+m,l+n)$-ribbon tangle depicted
by placing $\Gamma'$ immediately to the right of $\Gamma$ such that $\Gamma$ and
$\Gamma'$ are mutually unlinked.  

Let $\Gamma_{1}$, $\Gamma_{1}'$, $\Gamma_{2}$ and $\Gamma_{2}'$ 
be coloured directed ribbon tangles where 
$\Gamma_{1} \circ \Gamma_{1}'$ and $\Gamma_{2} \circ \Gamma_{2}'$ exist, then
$$(\Gamma_{1} \circ \Gamma_{1}') \otimes (\Gamma_{2} \circ \Gamma_{2}') = 
(\Gamma_{1} \otimes \Gamma_{2}) \circ (\Gamma_{1}' \otimes \Gamma_{2}').$$
Note that $(\Gamma_{1} \otimes \Gamma_{2}) \circ (\Gamma_{1}' \otimes \Gamma_{2}')$ may exist even if
 $\Gamma_{1} \circ \Gamma_{1}'$ and $\Gamma_{2} \circ \Gamma_{2}'$ do not.

Any coloured directed $(k,l)$-ribbon tangle, 
where the ribbon tangle is coloured with elements of a set $S$, can
be written down as some combination of vertical and horizontal compositions of
the directed ribbon tangle atoms, 
where each ribbon tangle atom is coloured with an element of $S$ \cite{rt}.
We denote a colouring of the directed ribbon tangle atoms as follows:
\begin{itemize}
\item[(i)]   $I^{\pm}_{i}$ means that $I^{\pm}$ is coloured with $i$,
\item[(ii)] $X^{+}_{i,j}$ means that $X^{+}$ is coloured so that 
the ribbon passing from the
top right hand corner to the bottom left corner is coloured with
$i$ and the other ribbon is coloured with $j$,
\item[(iii)]  $X^{-}_{i,j}$ means that $X^{-}$ is coloured in the same way that
              $X^{+}_{i,j}$ is,
\item[(iv)]  $W_{i}$ means that $W$ is coloured with $i$, where 
$W$ is any one of $\Omega^{+}$, $\Omega^{-}$, $U^{+}$ and $U^{-}$.
\end{itemize}

We now introduce the category of coloured directed ribbon tangles 
{\bf{cdrib}}(S) for a finite non-empty set $S$ \cite{resh,rt,cp}.
The objects of {\bf{cdrib}}(S) are sequences of pairs
$$(X,\epsilon)=\big((i_{1},\epsilon_{1}),(i_{2},\epsilon_{2}),\ldots, (i_{k},\epsilon_{k})  \big),$$ 
where $i_{j} \in S$ and $\epsilon_{i} \in
\{-1,+1\}$ for each $i=1, 2, \ldots, k$.  
Each object $(X,\epsilon)$ of {\bf{cdrib}}(S) can be thought of as 
$(X^{*},\epsilon^{*})_{\Gamma}$ or $(X_{*},\epsilon_{*})_{\Gamma}$
for some coloured directed ribbon tangle $\Gamma$.
The morphisms of {\bf{cdrib}}(S) 
are equivalence classes of coloured directed ribbon tangles coloured with elements of $S$.

\begin{lemma}  
All morphisms in {\bf{cdrib}}(S) can be expressed as vertical and horizontal compositions of the
coloured directed ribbon tangles $I_{i}^{+}$, $I^{-}_{i}$, $X^{+}_{i,j}$,
$X^{-}_{i,j}$, $\Omega_{i}^{+}$, $\Omega_{i}^{-}$, $U_{i}^{+}$ and $U_{i}^{-}$, where $i, j \in S$.
\end{lemma}
\begin{proof}
This follows from \cite[Prop. 1.4]{resh}.
\end{proof}

Let $\Gamma$ and $\Gamma'$ be two coloured directed ribbon tangles where $\Gamma \circ \Gamma'$ 
exists.
Let $f:(X^{*}, \epsilon^{*})_{\Gamma} \rightarrow (X_{*}, \epsilon_{*})_{\Gamma}$ and
$f':(X^{*}, \epsilon^{*})_{\Gamma'} \rightarrow (X_{*}, \epsilon_{*})_{\Gamma'}$ be two morphisms.
We define the morphism $f \circ f': (X^{*}, \epsilon^{*})_{\Gamma \circ \Gamma'} 
\rightarrow (X_{*}, \epsilon_{*})_{\Gamma \circ \Gamma'}$ to be the coloured directed ribbon tangle 
$\Gamma \circ \Gamma'$.

Let $(X_{a},\epsilon_{a})$ and $(X_{b},\epsilon_{b})$ be two objects of 
{\bf{cdrib}}(S):
\begin{eqnarray*}
(X_{a},\epsilon_{a}) & = & \big(
(i_{a_{1}},\epsilon_{a_{1}}),(i_{a_{2}},\epsilon_{a_{2}}),\ldots,
(i_{a_{k}},\epsilon_{a_{k}}) \big),  \\
(X_{b},\epsilon_{b}) & = & \big(
(i_{b_{1}},\epsilon_{b_{1}}),(i_{b_{2}},\epsilon_{b_{2}}),\ldots,
(i_{b_{m}},\epsilon_{b_{m}}) \big).
\end{eqnarray*}
The tensor product of $(X_{a},\epsilon_{a})$ and $(X_{b},\epsilon_{b})$
is an object $(X_{a},\epsilon_{a}) \otimes (X_{b},\epsilon_{b})$:
$$
(X_{a},\epsilon_{a}) \otimes (X_{b},\epsilon_{b}) 
= \big((i_{a_{1}},\epsilon_{a_{1}}),(i_{a_{2}},\epsilon_{a_{2}}),\ldots,
(i_{a_{k}},\epsilon_{a_{k}}),
(i_{b_{1}},\epsilon_{b_{1}}),(i_{b_{2}},\epsilon_{b_{2}}),\ldots,
(i_{b_{m}},\epsilon_{b_{m}})  \big).$$
Let $f$ and $g$ be morphisms of {\bf{cdrib}}(S):
$$f:(X_{a},\epsilon_{a})_{\Gamma} \rightarrow (X_{a}',\epsilon_{a}')_{\Gamma}, \hspace{10mm}
g:(X_{b},\epsilon_{b})_{\Gamma'} \rightarrow (X_{b}',\epsilon_{b}')_{\Gamma'}.$$
The tensor product of $f$ and $g$ is a morphism
$f \otimes g$ which is the coloured directed ribbon tangle $\Gamma \otimes \Gamma'$:
$$f \otimes g: (X_{a},\epsilon_{a})_{\Gamma} \otimes (X_{b},\epsilon_{b})_{\Gamma'} \rightarrow
(X_{a}',\epsilon_{a}')_{\Gamma} \otimes (X_{b}',\epsilon_{b}')_{\Gamma'}.$$
Note that if $f$, $f'$, $g$ and $g'$ are morphisms of {\bf{cdrib}}(S)
where $f \circ f'$ and $g \circ g'$ exist, then
$$(f \circ f') \otimes (g \circ g') = (f \otimes g) \circ (f' \otimes g').$$

\begin{subsubsection}{Closing a coloured directed ribbon tangle}

Let $\Gamma$ be a coloured directed $(k,k)$-ribbon tangle with
$(X_{*},\epsilon_{*})_{\Gamma} = (X^{*},\epsilon^{*})_{\Gamma}$.  
For each $j \in \{ 1,2,\ldots,k\}$ let $r_{j}$ be a ribbon 
that is the image of a smooth embedding 
$e:[0,1] \times [0,1] \rightarrow \mathbb{R}^{3}$ given as follows.
The ribbon $r_{j}$ has disjoint intersection with $\Gamma$:  
call $e\big([0,1] \times \{0\}\big)$ the {\emph{beginning}} of the ribbon and 
$e\big([0,1] \times \{1\}\big)$ the {\emph{end}} of the ribbon.  
Smoothly glue the beginning of $r_{j}$ 
to $\Gamma$ at $ [j-1/4,j+1/4] \times \{0\} \times \{1\}$ and the end of $r_{j}$ to $\Gamma$ at 
$[j-1/4,j+1/4] \times \{0\} \times \{0\}$ in such a way that 
$r_{j}$ is unlinked with each component of $\Gamma$, and also so that the
linking number of the two
edges $e\big(\{0\} \times [0,1]\big)$ and $e\big(\{1\} \times [0,1]\big)$ of $r_{j}$ is zero if they
are given parallel orientations.

Now orient and direct $r_{j}$ to be consistent with the
orientation and direction of the components of $\Gamma$ to which it is glued, and 
colour $r_{j}$ with $i_{j} \in S$.  
After attaching, orienting, directing and colouring the $k$ ribbons we obtain from 
$\Gamma$ a coloured directed $(0,0)$-ribbon tangle $\hat{\Gamma}$
heuristically depicted in Figure \ref{fig:closure} 
(here $T$ is a coloured directed $(k,k)$-ribbon tangle).
Sometimes we refer to $\hat{\Gamma}$ as the closure of $\Gamma$.
\begin{figure}[hbt]
\begin{center}
  \input{closure100.pstex_t}
\caption{A coloured directed ribbon tangle $\Gamma$ and its closure $\hat{\Gamma}$}  
          \label{fig:closure}
\end{center}
\end{figure}
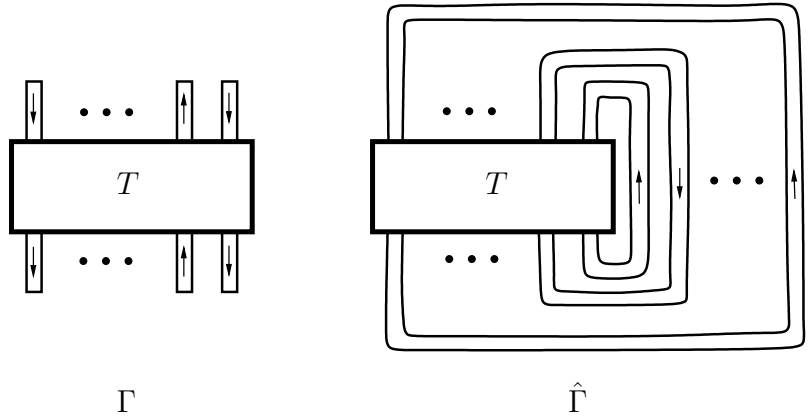

\end{subsubsection}

\end{subsection}

\begin{subsection}{Representing framed links as $(0,0)$-ribbon tangles}

Let $L = \bigcup_{i=1}^{m} L_{i} \subset S^{3}$ be a framed link with $m$ connected
components and let the framing number of the component $L_{i}$ be $n_{i}$ for each $i=1, 2, \ldots, m$.
We can associate $L$ with a $(0,0)$-ribbon tangle $\Gamma(L)$ in a natural way.
Recall that each ribbon tangle is a subset of $\mathbb{R}^{3}$ and  
that $S^{3} = \mathbb{R}^{3} \cup \{\infty\}$. 
Any link $L \subset S^{3}$ can be deformed so that $L \subset \mathbb{R}^{3}$ and ambient isotopic links in $S^{3}$ give rise to
ambient isotopic links in $\mathbb{R}^{3}$ \cite{t}.  
In the remainder of this subsection we  consider 
$L$ as a subset of $\mathbb{R}^{3}$.

We obtain the $(0,0)$-ribbon tangle $\Gamma(L)$ for each framed link $L \subset S^{3}$ as follows.  
For each connected component $L_{i}$ of $L$
fix $U_{i} \subset \mathbb{R}^{3}$ to be a small tubular neighbourhood
of $L_{i}$ such that $U_{i} \cap U_{j} = \emptyset$ if $i \neq j$.
Let $K_{i} \subset \partial(U_{i})$ be a parallel of $L_{i}$ such that $\mathrm{lk}(L_{i},K_{i})=n_{i}$, 
then we can equip $L_{i}$ with a normal vector field so that the tips of the vectors sweep out $K_{i}$.

The links $L_{i}$ and $K_{i}$ are the boundary of a $(0,0)$-ribbon tangle where 
each element of the normal vector field to $L_{i}$ coincides 
with some proper subset of the $(0,0)$-ribbon tangle.
Repeating this for each connected component of $L$ defines a $(0,0)$-ribbon tangle for $L$.  

\end{subsection}

\end{section}

\begin{section}{The Reshetikhin-Turaev functor $F$}
\label{sec:RTfunc}

\markright{\text{The Reshetikhin-Turaev functor $F$}}

The Reshetikhin-Turaev functor $F$ is a covariant functor from the category of 
coloured directed ribbon tangles to the category of finite dimensional representations 
of a $\mathbb{Z}_{2}$-graded ribbon Hopf algebra.
The functor provides the machinery allowing us to construct isotopy invariants of
ribbon tangles and thereby regular isotopy invariants of links.
The functor was first defined for ungraded ribbon Hopf algebras and their representations \cite{resh,reshtur,rt} 
and extended to $\mathbb{Z}_{2}$-graded ribbon Hopf algebras and their representations \cite{z3}.

Let $A$ be a $\mathbb{Z}_{2}$-graded ribbon Hopf algebra over $\mathbb{C}$.
Recall from Section \ref{chap1:quastriangaulrhopfalgebras} that $A$ 
admits an invertible even element, called the universal $R$-matrix,
\begin{equation}
\label{eq:tim(1)}
R =\sum_{t} a_{t} \otimes b_{t} \in A \otimes A,
\end{equation}
satisfying Eqs. (\ref{eq:asdf;lkj})--(\ref{eq:deddy}).  As usual we set 
\begin{equation}
\label{eq:tim(2)}
u=\sum_{t} S(b_{t})a_{t} (-1)^{[a_{t}]}.
\end{equation}
Recall that there exists an invertible even central element 
\begin{equation}
\label{eq:tim(3)}
v \in A,
\end{equation}
with the following properties:
$$\epsilon(v)=1, \hspace{5mm} v^{2} = uS(u), \hspace{5mm} S(v)=v, \hspace{5mm}
\Delta(v) = (v \otimes v)(R^{T}R)^{-1},$$
where $R^{T} = \sum_{t} b_{t} \otimes a_{t} (-1)^{[a_{t}]}$.

Let $Rep(A)$ be the category of finite dimensional $\mathbb{Z}_{2}$-graded left $A$-modules.
The objects of $Rep(A)$ are $\mathbb{Z}_{2}$-graded left $A$-modules over $\mathbb{C}$.  
The morphisms of $Rep(A)$ are $A$-linear homomorphisms of degree $0$. 
Let $\{V_{i}| \ i \in I\}$ be a set of objects of $Rep(A)$ for some index set $I$
such that for each $i \in I$, the dual $A$-module of $V_{i}$, which we denote by $(V_{i})^{*}$, 
is isomorphic to $V_{j}$  for some $j \in I$.   
Now let
$$\eta = \big( (i_{1},\epsilon_{1}),(i_{2},\epsilon_{2}),\ldots,(i_{k},\epsilon_{k})   \big),$$ 
be an object of ${\bf{cdrib}}(I)$, then 
for each such object $\eta$ we define the $A$-module
$$V_{\eta} = V_{i_{1}}^{\epsilon_{1}} \otimes V_{i_{2}}^{\epsilon_{2}} \otimes \cdots 
\otimes V_{i_{k}}^{\epsilon_{k}},$$ 
where we fix $V_{i}^{+1} = V_{i}$ and $V_{i}^{-1} = \left(V_{i}\right)^{*}$. 
If $\eta = \emptyset$ then we set
$V_{\emptyset} = V_{0} \cong \mathbb{C}$, which is the one-dimensional $A$-module.

We now present the theorem defining the covariant functor $F$  \cite{resh,reshtur,rt,z3}.
This theorem is a generalisation of a theorem defining 
$F$ as a covariant functor from the category of coloured directed ribbon tangles to the category of
finite dimensional representations of an ungraded ribbon Hopf algebra \cite{resh,reshtur,rt}.
The proof of Theorem \ref{the:RTfunctor} is similar to the proof of the theorem it generalises.
\begin{theorem}
\label{the:RTfunctor}
Let ${\cal{H}} = {\bf{cdrib}}(I)$.
Let $A$ be a $\mathbb{Z}_{2}$-graded ribbon Hopf algebra over $\mathbb{C}$
with universal $R$-matrix $R \in A \otimes A$ stated in (\ref{eq:tim(1)}), 
the element $u \in A$ stated in (\ref{eq:tim(2)}), 
and the even central element $v \in A$ stated in (\ref{eq:tim(3)}).  

Let $\{V_{i}| \ i \in I\}$ be a set of objects of $Rep(A)$ such that for each $i \in I$,
$(V_{i})^{*}$ is isomorphic to $V_{j}$ for some $j \in I$. 
There exists a covariant functor $F: {\cal H} \rightarrow Rep(A)$ with the following properties:
\begin{itemize}
\item[(i)] $F$ transforms any object $\eta$ of $\cal{H}$ into the object
$V_{\eta}$ of $Rep(A)$,
\item[(ii)] For any two coloured directed ribbon tangles $\Gamma, \Gamma'$,
$$F(\Gamma \otimes \Gamma')  = F(\Gamma) \otimes F(\Gamma'),$$
\item[(iii)] $F$ is defined by

\noindent
$
\begin{array}{ll}
F(I^{+}_{i}) = \mathrm{id_{i}}:V_{i} \rightarrow V_{i}; & F(I^{+}_{i})(x) = x, \\
F(I^{-}_{i}) = \mathrm{id_{i^{*}}}:\left(V_{i}\right)^{*} \rightarrow \left(V_{i}\right)^{*}; & F(I^{-}_{i})(x^{*}) = x^{*}, \\
\end{array}
$  

\noindent
$
\begin{array}{l}
F(X^{+}_{i,j})  =  P \circ R: V_{i} \otimes V_{j} \rightarrow V_{j} \otimes V_{i}; \\
F(X^{+}_{i,j})(x \otimes y)  = 
\sum_{t} b_{t} y \otimes a_{t} x \ (-1)^{[x][y]+[a_{t}]([b_{t}]+[y])}, \\
F(X^{-}_{i,j}) =  R^{-1}\circ P:V_{i} \otimes V_{j} \rightarrow V_{j} \otimes
V_{i}; \\
F(X^{-}_{i,j})(x \otimes y)  =  \sum_{t} S(a_{t})y \otimes b_{t}x \ 
(-1)^{[y]([x]+[b_{t}])},
\end{array}
$
  
\noindent
$
\begin{array}{ll}
F(\Omega_{i}^{+}): \left(V_{i}\right)^{*} \otimes V_{i} \rightarrow \mathbb{C}; &
F(\Omega_{i}^{+})(x^{*} \otimes y) = \langle x^{*}, y \rangle, \\
F(\Omega_{i}^{-}): V_{i} \otimes \left(V_{i}\right)^{*} \rightarrow \mathbb{C}; & 
F(\Omega_{i}^{-})(x \otimes y^{*}) = (-1)^{[x][y^{*}]} \langle y^{*}, (v^{-1}u) \ x\rangle, \\
F(U^{+}_{i}): \mathbb{C} \rightarrow V_{i} \otimes \left(V_{i}\right)^{*}; & 
F(U^{+}_{i})(c) = c \sum_{r}v_{r} \otimes v_{r}^{*}, \\
F(U_{i}^{-}): \mathbb{C} \rightarrow \left(V_{i}\right)^{*} \otimes V_{i}; & 
F(U_{i}^{-})(c) = c \sum_{r} (-1)^{[v_{r}]}v_{r}^{*} \otimes (vu^{-1}) \ v_{r},
\end{array}$

\noindent
where $\{v_{r}\}$ and $\{v_{r}^{*}\}$ are dual bases of $V_{i}$ and 
$\left(V_{i}\right)^{*}$, respectively, such that $\langle v_{r}^{*}, v_{s} \rangle = \delta_{rs}$ and
$[v_{r}]=[v_{r}^{*}]$.

\end{itemize}

\end{theorem}

Needless to say, if $\Gamma$ and $\Gamma'$ are
coloured directed ribbon tangles where $\Gamma \circ \Gamma'$ exists, then  we have
$$F(\Gamma \circ \Gamma') = F(\Gamma) \circ F(\Gamma').$$

Let $L$ be an oriented $(k, l)$-tangle with $m$ connected components (note that $k, l$ and $m \geq 1$ 
can all be different non-negative integers). 
Given such an oriented $(k, l)$-tangle $L$, fix $\Gamma$
to be the associated coloured directed
$(k, l)$-ribbon tangle with $m$ connected components,
where there are two associated sequences of pairs that 
encode the directions and colourings of all the ribbons of $\Gamma$
as given in Subsection \ref{subsec:johnfaulkner(9991)}:
\begin{eqnarray*}
\big(X^{*}, \epsilon^{*}\big)_{\Gamma} & = & 
\big( (\mu^{1}, \epsilon^{1}), (\mu^{2}, \epsilon^{2}), \ldots, (\mu^{k}, \epsilon^{k})
\big), \\
\big(X_{*}, \epsilon_{*}\big)_{\Gamma} & = & 
\big( (\nu_{1}, \epsilon_{1}), (\nu_{2}, \epsilon_{2}), \ldots, (\nu_{l}, \epsilon_{l})
\big),
\end{eqnarray*}
where $\mu^{i}, \nu_{j} \in I$.
The sequence $\big(X^{*}, \epsilon^{*}\big)_{\Gamma}$ uniquely specifies
the colourings and directions of the ribbon tangles intersecting the top of the ribbon tangle diagram.
In particular, the $i^{th}$ ribbon tangle from the left 
at the top of the ribbon tangle diagram is coloured with
$\mu^{i} \in I$ and is directed downwards (resp. upwards) if $\epsilon^{i} = +1$
(resp. $\epsilon^{i} = -1$).  

Similarly, the sequence $\big(X_{*}, \epsilon_{*}\big)_{\Gamma}$ uniquely specifies
the colourings and directions of the ribbon tangles intersecting the bottom 
of the ribbon tangle diagram:
the $j^{th}$ ribbon tangle from the left 
at the bottom of the ribbon tangle diagram
is coloured with $\nu_{j} \in I$ and is directed downwards (resp. upwards) if $\epsilon_{j} = +1$
(resp. $\epsilon_{j} = -1$).

For such a $\Gamma$,  $F(\Gamma)$ is a map
\begin{equation}
\label{eq:space1999}
F( \Gamma):
V_{\mu^{1}}^{\epsilon^{1}} \otimes V_{\mu^{2}}^{\epsilon^{2}} \otimes \cdots \otimes
V_{\mu^{k}}^{\epsilon^{k}} \longrightarrow
V_{\nu_{1}}^{\epsilon_{1}} \otimes V_{\nu_{2}}^{\epsilon_{2}} \otimes \cdots \otimes
V_{\nu_{l}}^{\epsilon_{l}}.
\end{equation}

As the functor $F$ maps morphisms of ${\cal{H}}$ to $A$-linear homomorphisms of degree $0$,
the map (\ref{eq:space1999}) must commute with the action of $A$, that is
$$F(\Gamma) \big( \pi_{\mu^{1}}^{\epsilon^{1}} \otimes \pi_{\mu^{2}}^{\epsilon^{2}}
\otimes \cdots \otimes \pi_{\mu^{k}}^{\epsilon^{k}} \left(\Delta^{(k-1)}(a)\right) \big) = 
\left(\pi_{\nu_{1}}^{\epsilon_{1}} \otimes \pi_{\nu_{2}}^{\epsilon_{2}} \otimes \cdots \otimes
\pi_{\nu_{l}}^{\epsilon_{l}}\big(\Delta^{(l-1)}(a)\big)\right) F(\Gamma),$$
for all $a \in A$, where $\pi_{i}^{+1}$ (resp. $\pi_{i}^{-1}$)
denotes the representation of $A$ afforded by the
$A$-module $V_{i}$ (resp. the dual $A$-module $\left(V_{i}\right)^{*}$).

To illuminate this fact, we give a direct proof of it.
Each coloured directed ribbon tangle $\Gamma$ can be expressed as some combination of
horizontal and vertical compositions of the coloured directed ribbon tangle atoms in 
Figure \ref{fig:ribbontangleatoms}, thus 
the homomorphism $F(\Gamma)$ of $A$-modules can be expressed as some 
appropriate combination of tensor products and compositions of the homomorphisms of 
$A$-modules obtained by
 applying the functor $F$ to the coloured directed ribbon tangle atoms in 
Figure \ref{fig:ribbontangleatoms}.  
To prove the theorem it suffices to prove it for each case in which $\Gamma$ is a 
coloured directed ribbon tangle atom.

The theorem is trivially true for $F(I^{+}_{i})$ and $F(I^{-}_{i})$, and it may be shown to be 
true for $F(X^{+}_{i,j})$ and $F(X^{-}_{i,j})$ by simple calculations using 
$R\Delta(x) = \Delta'(x)R, \ \forall x \in A$.  Now
\begin{eqnarray*}
F(\Omega_{i}^{+}) a (x^{*} \otimes y) 
& = & F(\Omega_{i}^{+}) \sum_{(a)} a_{(1)} x^{*} \otimes a_{(2)} y (-1)^{[x^{*}][a_{(2)}]} \\
& = & \sum_{(a)} \langle x^{*}, S(a_{(1)}) a_{(2)} y \rangle (-1)^{[x^{*}][a]} \\
& = & \epsilon(a) F(\Omega_{i}^{+}) (x^{*} \otimes y),
\end{eqnarray*}
where we have used the fact that $\epsilon(a)=0$ if $[a]=1$. 
Also,
\begin{eqnarray*}
a F(U^{+}_{i})(c) 
 & = & c \sum_{(a),r} a_{(1)} v_{r} \otimes a_{(2)} v_{r}^{*} (-1)^{[v_{r}][a_{(2)}]} \\
 & = & c \sum_{(a),r,i} a_{(1)} v_{r} \otimes 
 \langle a_{(2)} v_{r}^{*},v_{i} \rangle v_{i}^{*} (-1)^{[v_{r}][a_{(2)}]} \\
 & = & c \sum_{(a),r,i} a_{(1)}  
 \langle v_{r}^{*}, S(a_{(2)})v_{i} \rangle  v_{r} \otimes  v_{i}^{*} \\
 & = & c \sum_{(a),i} a_{(1)} S(a_{(2)})v_{i} \otimes v_{i}^{*} \\
 & = & F(U^{+}_{i})(c)  \epsilon(a).
\end{eqnarray*}
The proofs for $F(\Omega^{-}_{i})$ and $F(U_{i}^{-})$ are similar.

\begin{corollary}
\label{cor:thestoneroses}
Let $\Gamma(L, \lambda)$ be a coloured directed
$(1,1)$-ribbon tangle with $m$ components where the ribbon joining $\mathbb{R}^{2} \times \{0\}$ and
$\mathbb{R}^{2} \times \{1\}$ is coloured with $\lambda_{i} \in I$ and directed downwards at its
bases.  
Then the map $F\big(\Gamma(L, \lambda)\big): V_{\lambda_{i}} \rightarrow V_{\lambda_{i}}$ 
is an element of $End_{A}(V_{\lambda_{i}})$.
\end{corollary}

\begin{remark}
Let $\Gamma(L, \lambda)$
be a coloured directed $(0,0)$-ribbon tangle associated with the framed oriented link $L$.
Then from the definition of the functor $F$, the map
$F\big(\Gamma(L, \lambda) \big): \mathbb{C} \rightarrow \mathbb{C}$
is an invariant of isotopy of $\Gamma(L)$ and thus an invariant of regular isotopy of  $L$ \cite{resh}.
\end{remark}

\begin{subsection}{Calculations using the functor $F$}
\label{subsec:hopefullylast(2)}

We now do some calculations using the functor $F$.
These results will be needed in the later sections.

Let $A$ be a  $\mathbb{Z}_{2}$-graded ribbon Hopf algebra and
$\{V_{i} | \ i \in I\}$ a set of non-isomorphic irreducible $A$-modules 
such that for each $i \in I$, $(V_{i})^{*} \cong V_{j}$ for some $j \in I$.
Consider the two framed oriented links $L$, $L' \subset S^{3}$ in Figure \ref{fig:Kappa+move2},
the planar projections of which are equivalent with respect to the Kirby move $\kappa_{+}$.  Here
$T$ is an arbitrary oriented $(m,m)$-tangle represented by the rectangle, the orientation of which is
compatible with the orientations of $L \backslash T$ and $L' \backslash T$.  
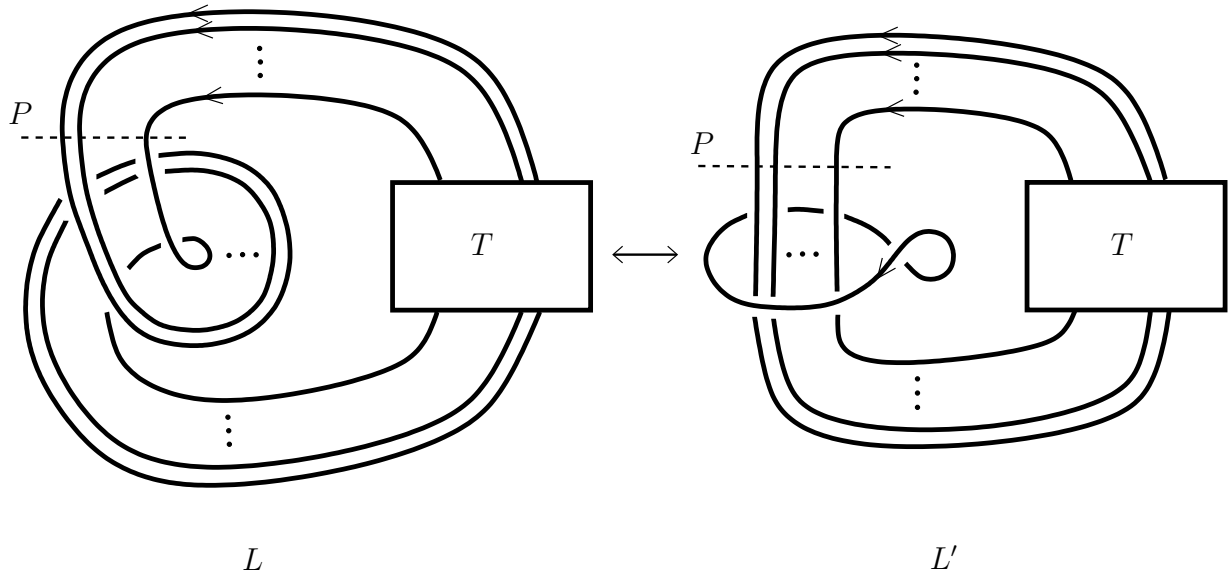
\begin{figure}[htb]
\begin{center}
 \input{Kirby+move100.pstex_t}
\caption{Two links $L$, $L'$ related by the Kirby $\kappa_{+}$ move}  \label{fig:Kappa+move2}
\end{center}
\end{figure}
Let the directed ribbon tangles associated with $L$ and $L'$ be coloured as
follows: let $\Gamma(L_{i})$ (resp. $\Gamma(L_{i}')$) denote 
the directed $(1,1)$-ribbon tangle derived from the $i^{th}$ tangle of 
$L \backslash T$ (resp. $L' \backslash T$) intersecting the line $P$ from the left.  
Colour both $\Gamma(L_{i})$ and $\Gamma(L_{i}')$ with $\lambda_{i} \in I$, and  
colour the unique component of
$\Gamma(L' \backslash T)$ isotopic to an unknotted annulus with $\mu \in I$.  
We denote these coloured directed $(m, m)$-ribbon tangles, respectively, by
$$\Gamma\big(L \backslash T,(\lambda_{1},\lambda_{2},\ldots,\lambda_{m})\big), \hspace{10mm}
\Gamma\big(L' \backslash T,(\lambda_{1},\lambda_{2},\ldots,\lambda_{m},\mu)\big).$$

\begin{lemma}
\label{lem:xylophone}
Assume that $v$ acts on $V_{\mu}$ as the multiplication by a scalar.
Fix $m = 1$, then
\begin{eqnarray}
F\big(\Gamma(L\backslash T, \lambda )\big) & = & v: V_{\lambda} \rightarrow
V_{\lambda}, \label{eq:markhopkins(aa)} \\
F\big(\Gamma(L'\backslash T,(\lambda, \mu ))\big) & = & 
\chi_{\mu}(v^{-1})C_{\mu}: V_{\lambda} \rightarrow V_{\lambda}, \label{eq:markhopkins(bb)}
\end{eqnarray}
where 
\begin{equation}
\label{eq:theCmuequation(1)}
C_{\mu} = (\mathrm{id} \otimes \mathrm{str})\left[({\mathrm{id}} \otimes \pi_{\mu})
(\mathrm{id} \otimes v^{-1}u) R^{T}R \right],
\end{equation} 
is a central element of $A$, and
$\chi_{\mu}(v^{-1})$ is the eigenvalue of $v^{-1}$ in the $A$-representation $\pi_{\mu}$.
\end{lemma}
\begin{lemma}
\label{lem:maybethelast?}
For each $m \geq 2$, the map
$$F\big(\Gamma(L\backslash T,(\lambda_{1}, \ldots, \lambda_{m}))\big):
V_{\lambda_{1}} \otimes \cdots \otimes V_{\lambda_{m}} \rightarrow V_{\lambda_{1}} 
\otimes \cdots \otimes V_{\lambda_{m}},$$
acts as
$$\Delta^{(m-1)}(v): 
V_{\lambda_{1}} \otimes \cdots \otimes V_{\lambda_{m}} \rightarrow V_{\lambda_{1}} \otimes \cdots \otimes V_{\lambda_{m}}.$$
\end{lemma}
\begin{proof}
Assume that Lemma \ref{lem:xylophone} is true (we will prove it below), then
Figure \ref{fig:bbbb} proves the lemma for $m=2$. 
Now assume that the inductive hypothesis is true for some $m \geq 2$, 
then the proof follows for $(m+1)$ by
using Figure \ref{fig:dddd} and the representation of $\Delta^{(m)}(v)$ in tensor product
representations of $A$. 
\end{proof}
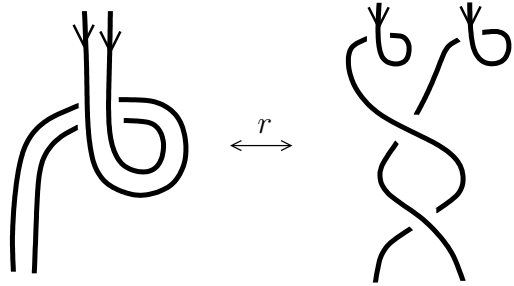
\begin{figure}[hbt]
 \begin{center}
 \input{bbbb100.pstex_t}
 \caption{Regularly isotopic $(2,2)$-tangles} \label{fig:bbbb}
  \end{center}
\end{figure}
\begin{figure}[hbt]
 \begin{center}
 \input{dddd100.pstex_t}
 \caption{Regularly isotopic $(m+1, m+1)$-tangles} \label{fig:dddd}
  \end{center}
\end{figure}
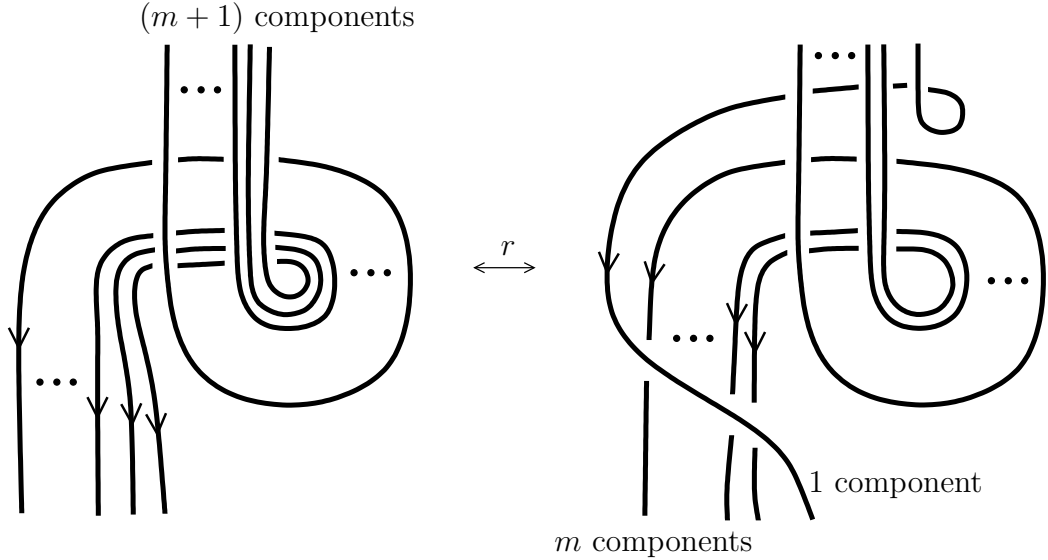
\begin{remark}
\label{rem:markhopkins(aaa)}
The results of Eq. (\ref{eq:markhopkins(aa)}) in Lemma \ref{lem:xylophone}   and
Lemma \ref{lem:maybethelast?} still hold true if the modules are not irreducible.
\end{remark}

We now prove Lemma \ref{lem:xylophone}.
\begin{proof}
We consider the first claim.  
Theorem \ref{the:RTfunctor} implies that 
$F\big(\Gamma(L\backslash T, \lambda )\big): V_{\lambda} \rightarrow V_{\lambda}$ is given by
$$
F\big(\Gamma(L\backslash T, \lambda )\big)
=
\mathrm{id}_{\lambda} \circ (\mathrm{id}_{\lambda} \otimes \Omega_{\lambda}^{-}) \circ 
(X^{-}_{\lambda, \lambda} \otimes \mathrm{id}_{\lambda}) \circ 
(\mathrm{id}_{\lambda} \otimes U_{\lambda}^{+})
\circ \mathrm{id}_{\lambda}.
$$
Let $x$ be a basis vector of $V_{\lambda}$ and let $\{v_{r}\}$, $\{v_{r}^{*}\}$ be dual bases of 
$V_{\lambda}$, $\left(V_{\lambda}\right)^{*}$, respectively, such that 
$\langle v_{r}^{*}, v_{s} \rangle = \delta_{rs}$ and $[v_{r}^{*}]=[v_{r}]$,
where $\langle \cdot, \cdot \rangle: (V_{\lambda})^{*} \times V_{\lambda} \rightarrow \mathbb{C}$ 
is the dual space pairing.
We calculate the action of $F\big(\Gamma(L\backslash T, \lambda )\big)$ on $x$ to be
\begin{eqnarray*}
x & \stackrel{(\mathrm{id}_{\lambda} \otimes U_{\lambda}^{+})}{\longmapsto}
 & \sum_{r} x \otimes  v_{r} \otimes v_{r}^{*} \\
  & \stackrel{(X^{-}_{\lambda, \lambda} \otimes \mathrm{id}_{\lambda})}{\longmapsto}
  & \sum_{t, r} S(a_{t}) v_{r} \otimes b_{t}x \otimes 
                   v_{r}^{*} (-1)^{[x][v_{r}]+[b_{t}][v_{r}]} \\
  & \stackrel{(\mathrm{id}_{\lambda} \otimes \Omega_{\lambda}^{-})}{\longmapsto} 
  & \sum_{t, r}  \langle v_{r}^{*}, v^{-1}u b_{t}x \rangle S(a_{t}) v_{r} \\
  & = & v^{-1} S(u)u x =  v x, 
\end{eqnarray*}
as $v^{2} = S(u)u$.  

We now consider the second claim.
We will calculate the action of $F\big(\Gamma(L' \backslash T,(\lambda,\mu))\big)$
on a basis vector $x \in V_{\lambda}$ where the annulus is coloured with $\mu$.  
Consider the regularly isotopic oriented $(1,1)$-tangles in Figure \ref{fig:figleaf}, 
where $\stackrel{r}{\longleftrightarrow}$ indicates regular isotopy. 
As $L'$ and $M'$ are regularly isotopic, 
$F\big(\Gamma(L',\nu)\big)=F\big(\Gamma(M',\nu)\big)$, and we will now 
calculate $F\big(\Gamma(M',\nu)\big)$.
\begin{figure}[hbt]
 \begin{center}
 \input{afigleaf100.pstex_t}
 \caption{Regularly isotopic oriented $(1,1)$-tangles} \label{fig:figleaf}
  \end{center}
\end{figure}
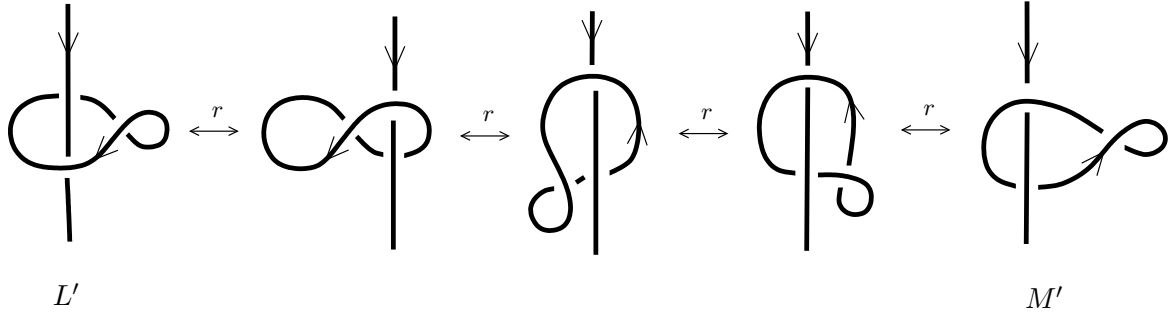
Consider 
the coloured directed $(1,1)$-ribbon tangle $\Gamma(M'', (\lambda,\mu))$ in Figure
\ref{fig:whatyoutalkingabout}: note that $\Gamma(M'', (\lambda,\mu))$ 
is identical to $\Gamma(M', (\lambda,\mu))$ `modulo' a twist.  
\begin{figure}[hbt]
 \begin{center}
 \input{zebra100.pstex_t}
 \caption{The coloured directed $(1,1)$-ribbon tangle $\Gamma(M'', ( \lambda,\mu ) )$} \label{fig:whatyoutalkingabout}
  \end{center}
\end{figure}
The map $F\big(\Gamma(M'', (\lambda,\mu))\big): V_{\lambda} \rightarrow V_{\lambda}$ is
$$
F\big(\Gamma(M'', (\lambda,\mu ) )\big) 
=
\mathrm{id}_{\lambda} \circ ( \mathrm{id}_{\lambda} \otimes \Omega_{\mu}^{-} ) \circ
( X^{+}_{\mu, \lambda} \otimes \mathrm{id}_{\mu^{*}}) \circ 
(X^{+}_{\lambda,\mu} \otimes \mathrm{id}_{\mu^{*}})
\circ (\mathrm{id}_{\lambda} \otimes U_{\mu}^{+}) \circ \mathrm{id}_{\lambda}.
$$
Let $x$ be a basis vector of $V_{\lambda}$, 
then we calculate $F\big(\Gamma(L, (\lambda,\mu))\big) (x)$ to be
\begin{eqnarray*}
x & \stackrel{(\mathrm{id}_{\lambda} \otimes U_{\mu}^{+})}{\longmapsto} & 
\sum_{r} x \otimes v_{r} \otimes v_{r}^{*} \\
  & \stackrel{(X^{+}_{\lambda,\mu} \otimes \mathrm{id}_{\mu^{*}})}{\longmapsto} & 
  \sum_{t,r} b_{t} v_{r} \otimes a_{t} x \otimes
              v_{r}^{*}  (-1)^{[x][v_{r}]+[a_{t}](1+[v_{r}])} \\
  & \stackrel{( X^{+}_{\mu, \lambda} \otimes \mathrm{id}_{\mu^{*}})}{\longmapsto} & 
  \sum_{s, t, r} b_{s}a_{t} x \otimes a_{s}b_{t}
v_{r} \otimes v_{r}^{*} (-1)^{[b_{s}]+[x]([a_{s}]+[b_{t}]) + [a_{t}][a_{s}]} \\
  & \stackrel{( \mathrm{id}_{\lambda} \otimes \Omega_{\mu}^{-} )}{\longmapsto} & 
  \sum_{s,t,r}  \langle v_{r}^{*}, (v^{-1}u) a_{s}b_{t} v_{r} \rangle b_{s}a_{t}x 
  (-1)^{[b_{s}] + [x]([b_{t}]+[a_{s}]) + [a_{t}][a_{s}] +[v_{r}^{*}]([a_{s}]+[b_{t}]+[v_{r}])} \\
  & = & (\mathrm{id} \otimes \mathrm{str}) \left[ (\mathrm{id} \otimes \pi_{\mu})(\mathrm{id} \otimes v^{-1}u)
R^{T}R \right]  x  =  C_{\mu} (x).
\end{eqnarray*}
Note that $C_{\mu}$ is a central element of $A$ \cite[Prop. 3]{zg} and 
$\chi_{0}(C_{\mu}) = sdim_{q}(V_{\mu})$.
Adding a twist with framing number $+1$ into the annulus coloured with $\mu$ gives 
 precisely $\Gamma(L' \backslash T, (\lambda, \mu))$.   
 Lemma \ref{lem:2010spaceodyssey} (ii) below implies
that  $F \big(\Gamma(L' \backslash T, (\lambda, \mu)) \big) = 
\chi_{(V_{\mu})^{*}}(v^{-1}) F\big(\Gamma(L, (\lambda, \mu))\big)$, 
where $\chi_{(V_{\mu})^{*}}(v^{-1})$ is the scalar action of $v^{-1}$ on $(V_{\mu})^{*}$.
As $v$ is even and $S(v)=v$, we have $\chi_{(V_{\mu})^{*}}(v^{-1})=\chi_{\mu}(v^{-1})$.
\end{proof}

\begin{lemma}
\label{lem:2010spaceodyssey}
{\emph{ }}
\begin{itemize}
\item[(i)]
Let $\Gamma(L, \lambda)$ be the coloured directed $(1,1)$-ribbon tangle in Figure
\ref{fig:whatyou}, then 
$F\big( \Gamma(L, \lambda) \big) = v^{-1}: V_{\lambda} \rightarrow V_{\lambda}$.
\item[(ii)]
Let $\Gamma(L',\lambda)$ be the coloured directed $(1,1)$-ribbon tangle obtained by reversing the
direction of $\Gamma(L, \lambda)$, then
$F\big(\Gamma(L', \lambda)\big) = v^{-1} : \left(V_{\lambda}\right)^{*} \rightarrow \left(V_{\lambda}\right)^{*}$.
\end{itemize}
\end{lemma}
\begin{figure}[hbt]
 \begin{center}
 \input{zebra2100.pstex_t}
 \caption{The coloured directed $(1,1)$-ribbon tangle $\Gamma(L, \lambda)$} \label{fig:whatyou}
  \end{center}
\end{figure}
\begin{proof}
We prove (i).  The map
$F\big(\Gamma(L, \lambda)\big): V_{\lambda} \rightarrow V_{\lambda}$ is 
$$
F\big(\Gamma(L, \lambda)\big)
=
\mathrm{id}_{\lambda} \circ (\mathrm{id}_{\lambda} \otimes \Omega_{\lambda}^{-} ) \circ
( X^{+}_{\lambda, \lambda} \otimes \mathrm{id}_{\lambda^{*}})
\circ (\mathrm{id}_{\lambda} \otimes U_{\lambda}^{+}) \circ \mathrm{id}_{\lambda}.
$$
Let $x$ be a basis vector of $V_{\lambda}$, 
then we calculate $F\big(\Gamma(L,\lambda)\big)(x)$ as follows:
\begin{eqnarray*}
x & \stackrel{(\mathrm{id}_{\lambda} \otimes U_{\lambda}^{+})}{\longmapsto} & 
                                           \sum_{r} x \otimes v_{r} \otimes v_{r}^{*} \\
& \stackrel{( X^{+}_{\lambda, \lambda} \otimes \mathrm{id}_{\lambda^{*}})}{\longmapsto}
& \sum_{t,r} b_{t} v_{r} \otimes a_{t} x \otimes v_{r}^{*} (-1)^{[a_{t}]+[v_{r}]([a_{t}]+[x])} \\
& \stackrel{(\mathrm{id}_{\lambda} \otimes \Omega_{\lambda}^{-} )}{\longmapsto}
& \sum_{t,r} \langle v_{r}^{*}, v^{-1} u \ a_{t} x \rangle b_{t} v_{r} (-1)^{[a_{t}]} \\
& = & \sum_{t} v^{-1} b_{t} u a_{t} \ x (-1)^{[a_{t}]} \\
& = & \sum_{t} v^{-1} b_{t} S^{2}(a_{t})u \ x (-1)^{[a_{t}]}  =  v^{-1} x.
\end{eqnarray*}

\begin{figure}[hbt]
 \begin{center}
 \input{tinytilly100.pstex_t}
 \caption{Two regularly isotopic oriented $(1,1)$-tangles} \label{fig:felinious}
  \end{center}
\end{figure}
We now prove (ii).
Consider the two regularly isotopic oriented $(1,1)$-tangles $L'$ and $L''$ in Figure \ref{fig:felinious}.
As these tangles are regularly isotopic, we can rewrite the map 
$F\big(\Gamma(L', \lambda)\big): \left(V_{\lambda}\right)^{*} \rightarrow \left(V_{\lambda}\right)^{*}$ as
\begin{eqnarray}
F\big(\Gamma(L', \lambda)\big)
& = & 
F\big(\Gamma(L'', \lambda)\big) \nonumber \\
& = & 
\mathrm{id}_{\lambda^{*}} \circ (\mathrm{id}_{\lambda^{*}} \otimes \Omega_{\lambda}^{-}) \circ
( \mathrm{id}_{\lambda^{*}} \otimes \Omega^{+}_{\lambda} \otimes \mathrm{id}_{\lambda} \otimes
\mathrm{id}_{\lambda^{*}}) \circ
(\mathrm{id}_{\lambda^{*}} \otimes \mathrm{id}_{\lambda^{*}} \otimes X^{+}_{\lambda,\lambda} \otimes
\mathrm{id}_{\lambda^{*}}) \nonumber \\
& & \hspace{5mm} \circ (\mathrm{id}_{\lambda^{*}} \otimes U^{-}_{\lambda} \otimes
\mathrm{id}_{\lambda} \otimes \mathrm{id}_{\lambda^{*}}) \circ
(U^{-}_{\lambda} \otimes \mathrm{id}_{\lambda^{*}}) \circ \mathrm{id}_{\lambda^{*}}. \label{eq:julk}
\end{eqnarray}
Let $x$ be a basis vector of $\left(V_{\lambda}\right)^{*}$, then the action of 
the right hand side of (\ref{eq:julk}) on $x$ is
\begin{eqnarray}
x & \stackrel{U_{\lambda}^{-} \otimes \mathrm{id}_{\lambda^{*}}}{\longmapsto} & 
\sum_{r} v_{r}^{*} \otimes (vu^{-1}) v_{r} \otimes x (-1)^{[v_{r}]} \nonumber \\
& \stackrel{ \mathrm{id}_{\lambda^{*}} \otimes U_{\lambda}^{-} \otimes \mathrm{id}_{\lambda} \otimes
\mathrm{id}_{\lambda^{*}}}{\longmapsto} & \sum_{r,p} v_{r}^{*} \otimes v_{p}^{*} \otimes
(v u^{-1}) v_{p} \otimes (v u^{-1}) v_{r} \otimes x (-1)^{[v_{r}]+[v_{p}]} \nonumber \\
& \stackrel{\mathrm{id}_{\lambda^{*}} \otimes \mathrm{id}_{\lambda^{*}} \otimes
X^{+}_{\lambda,\lambda} \otimes \mathrm{id}_{\lambda^{*}}}{\longmapsto} & 
\sum_{r,p,t} v_{r}^{*} \otimes v_{p}^{*} \otimes b_{t} v u^{-1} v_{r} \otimes a_{t} v u^{-1} v_{p}
\otimes x (-1)^{[v_{r}]+[v_{p}]+[a_{t}](1 + [v_{r}]) + [v_{p}][v_{r}]} \nonumber \\
& \stackrel{\mathrm{id}_{\lambda^{*}} \otimes \Omega^{+}_{\lambda} \otimes 
\mathrm{id}_{\lambda} \otimes \mathrm{id}_{\lambda^{*}}}{\longmapsto} &
\sum_{r,p,t} \langle u v^{-1} S^{-1}(b_{t}) v_{p}^{*}, v_{r} \rangle
v_{r}^{*} \otimes a_{t} v u^{-1} v_{p} \otimes x
(-1)^{[v_{p}][b_{t}] + [v_{p}] + [a_{t}] + [v_{r}](1 + [a_{t}] + [v_{p}])}
\nonumber  \\
& = & \sum_{p,t} u v^{-1} S^{-1}(b_{t}) v_{p}^{*} \otimes a_{t} vu^{-1} v_{p} \otimes x
(-1)^{[v_{p}]+[a_{t}]+[v_{p}][b_{t}]} \label{eq:tillymonster} \\
& \stackrel{\mathrm{id}_{\lambda^{*}} \otimes \Omega^{-}_{\lambda}}{\longmapsto} & 
\sum_{p,t} \langle x, v^{-1}u a_{t} v u^{-1} v_{p} \rangle uv^{-1} S^{-1}(b_{t}) v_{p}^{*}
(-1)^{[x]([a_{t}]+[v_{p}]) + [v_{p}] + [a_{t}] + [v_{p}][b_{t}]} \nonumber \\
& = & \sum_{p,t} \langle S(a_{t}) x, v_{p} \rangle u v^{-1} S^{-1}(b_{t}) v_{p}^{*}
(-1)^{[v_{p}](1+[x]+[b_{t}])+[a_{t}]} \nonumber \\
& = & \sum_{t} uv^{-1} S^{-1}(b_{t}) S(a_{t}) x (-1)^{[a_{t}]} 
= v^{-1}x. \nonumber
\end{eqnarray}
Note that (\ref{eq:tillymonster}) is obtained by using the fact that 
$\langle u v^{-1} S^{-1}(b_{t}) v_{p}^{*}, v_{r} \rangle$ 
vanishes unless $[v_{r}] \equiv ([v_{p}] + [b_{t}]) \pmod{2}$, allowing us to simplify the
sign factors.  
\end{proof}

\begin{proposition}
Let $\Gamma(L,\lambda)$ be the coloured directed
$(1,1)$-ribbon tangle in Figure \ref{fig:allanbridges} such that the map
$F\big(\Gamma(L,\lambda)\big): V_{\lambda} \rightarrow V_{\lambda}$ 
acts on $V_{\lambda}$ as the multiplication by the scalar $\zeta$.
Now let $\hat{\Gamma}(L,\lambda)$ be the closure of 
$\Gamma(L,\lambda)$ in Figure \ref{fig:allanbridges2},
then $F\big(\hat{\Gamma}(L,\lambda)\big): \mathbb{C} \rightarrow \mathbb{C}$ 
acts as the multiplication by $\zeta sdim_{q}(V_{\lambda})$.
\begin{figure}[hbt]
 \begin{center}
 \input{baby100.pstex_t}
 \caption{A coloured directed $(1,1)$-ribbon tangle $\Gamma(L,\lambda)$} 
 \label{fig:allanbridges}
  \end{center}
\end{figure}
\begin{figure}[hbt]
 \begin{center}
 \input{baby1100.pstex_t}
 \caption{The closure $\hat{\Gamma}(L,\lambda)$ 
 of the ribbon tangle $\Gamma(L,\lambda)$} \label{fig:allanbridges2}
  \end{center}
\end{figure}
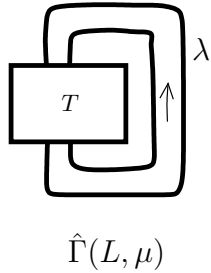
\end{proposition}
\begin{proof}
Explicitly, we have
\begin{equation}
\label{eq:roosters}
F\big(\hat{\Gamma}(L,\lambda)\big) =
\Omega^{-}_{\lambda} \circ (\zeta \mathrm{id}_{\lambda} \otimes \mathrm{id}_{\lambda^{*}}) \circ
U^{+}_{\lambda}: \mathbb{C} \rightarrow \mathbb{C},
\end{equation}
and for any $c \in \mathbb{C}$,
\begin{eqnarray*}
c & \stackrel{U^{+}_{\lambda}}{\longmapsto} & c \ \sum_{r} v_{r} \otimes v_{r}^{*} 
   \stackrel{\zeta \mathrm{id}_{\lambda} \otimes \mathrm{id}_{\lambda^{*}}}{\longmapsto}  
  c \ \zeta \sum_{r} v_{r} \otimes v_{r}^{*} \\
  & \stackrel{\Omega^{-}_{\lambda}}{\longmapsto} & 
  c \ \zeta \sum_{r} \langle v_{r}^{*}, (v^{-1}u) v_{r} \rangle (-1)^{[v_{r}]}  =  
  c \ \zeta sdim_{q}(V_{\lambda}).
\end{eqnarray*}
\end{proof}

\end{subsection}

\end{section}

\begin{section}{Pseudo-modular Hopf algebras}
\label{chapter4:sectionlabel(Pseudo-Modular_Hopf_algebras)}

\markright{\text{Pseudo-modular Hopf algebras}}

\begin{subsection}{Modular Hopf algebras}

We now introduce the notion of modular Hopf algebras that play 
a central role in the construction of topological invariants of $3$-manifolds in \cite{rt}.  
We extend the definition to the $\mathbb{Z}_{2}$-graded case in the obvious way.

\begin{definition}
\label{def:modHopfalg}
Let $A$ be a $\mathbb{Z}_{2}$-graded ribbon Hopf algebra over $\mathbb{C}$
with universal $R$-matrix $R \in A \otimes A$ stated in (\ref{eq:tim(1)}), 
the element $u \in A$ stated in (\ref{eq:tim(2)}), 
and the even central element $v \in A$ stated in (\ref{eq:tim(3)}).  

Let $I$ be a finite index set with an involution $*:I \rightarrow I$  denoted by
$*(i) = i^{*}$, and let there exist a distinguished element $0 \in I$ satisfying $*(0)=0^{*}=0$.  
Let $\{V_{i} | \ i \in I\}$ be a set of  $A$-modules, 
where $V_{0}$ is the one-dimensional $A$-module, and let
there exist a set of $A$-linear isomorphisms
$$\big\{ \omega_{i}: \left(V_{i}\right)^{*} \rightarrow V_{i^{*}}| \ i \in I  \big\},$$
where $\omega_{0} = \mathrm{id}$ and $\left(V_{i}\right)^{*}$ is the dual $A$-module to $V_{i}$.

The $\mathbb{Z}_{2}$-graded ribbon Hopf algebra 
$A$ is called a $\mathbb{Z}_{2}$-graded modular Hopf algebra if the following axioms are satisfied:
\begin{itemize}
\item[(i)] The $A$-modules $\{V_{i} | \ i \in I \}$ are irreducible, finite-dimensional, 
mutually non-isomorphic and each $V_{i}$ has non-vanishing quantum superdimension.

\item[(ii)]  For each $i \in I$, let $\{ x_{j} \}$, $\{ x_{j}^{*} \}$ and 
$\{ x_{j}^{**} \}$ be bases of $V_{i}$, $(V_{i})^{*}$ and $((V_{i})^{*})^{*}$, 
respectively, such that
$\langle \langle x_{j}^{**}, x_{k}^{*} \rangle \rangle 
= \langle x_{k}^{*}, x_{j} \rangle = \delta_{j,k}$ and 
$[x_{j}^{**}] = [x_{j}^{*}] = [x_{j}]$, where 
$\langle \langle \cdot, \cdot \rangle \rangle: 
((V_{i})^{*})^{*} \times (V_{i})^{*} \rightarrow \mathbb{C}$ is the dual space pairing.
Then the map
$$(\omega_{i})^{*} \circ (\omega_{i^{*}})^{-1}:  V_{i} \rightarrow ((V_{i})^{*})^{*},$$
is given by $x_{j} \mapsto (-1)^{[x_{j}]} u v^{-1} x_{j}^{**}, \ \forall j$.
\item[(iii)] For any $t \in \mathbb{N}$ and any sequence
$\theta = (i_{1}, i_{2}, \ldots, i_{t}) \in I^{\times t}$,
$$V_{i_{1}} \otimes V_{i_{2}} \otimes \cdots \otimes V_{i_{t}} 
= \left(\bigoplus_{i \in I}\left(V_{i}\right)^{\oplus n_{\theta}(i)} \right) 
\oplus {\cal{Z}}_{\theta},$$
where $n_{\theta}(i) \in \mathbb{Z}_{+}$ is the number of copies of the $A$-module $V_{i}$ 
in the direct sum and
${\cal{Z}}_{\theta}$ is a possibly vanishing $A$-module that satisfies axiom (iv).
\item[(iv)] For each $\theta = (i_{1}, i_{2}, \ldots, i_{t})$, $t \geq 2$, 
and all $A$-linear homomorphisms $\phi: {\cal{Z}}_{\theta} \rightarrow {\cal{Z}}_{\theta}$, $\mbox{str}_{q}(\phi) = 0$.
\item[(v)] Let $f=(f_{\lambda\mu})_{\lambda,\mu \in I}$ be the matrix with complex entries given by
$$f_{\lambda\mu} = ({\mathrm{str}} \otimes {\mathrm{str}})
\left[(\pi_{\lambda} \otimes \pi_{\mu})(v^{-1}u\otimes v^{-1}u)R^{T}R\right].$$
The matrix $f$ is invertible and there exists a 
unique set $\{d_{i} \in \mathbb{C} | \ i \in I\}$ of constants satisfying the relations
\begin{equation}
\label{eq:waggahogga}
\chi_{\lambda}(v)sdim_{q}(V_{\lambda}) = 
\sum_{\mu \in I}d_{\mu}\chi_{\mu}(v^{-1})f_{\lambda\mu}, \hspace{5mm} \mbox{for all } \lambda \in I.
\end{equation}
Here $\pi_{\mu}$ is the irreducible $A$-representation furnished by $V_{\mu}$, and
$\chi_{\mu}(v)$ denotes the eigenvalue of $v$ in $\pi_{\mu}$.
\end{itemize}
\end{definition}

Using results of Subsection \ref{subsec:hopefullylast(2)}, we have
$$f_{\lambda\mu} = sdim_{q}(V_{\lambda}) \chi_{\lambda}(C_{\mu}),$$ 
where $C_{\mu}$ is the central element of $A$ defined by Eq. (\ref{eq:theCmuequation(1)}).

An additional axiom was included in the original definition of modular Hopf algebras \cite{rt}.  
This axiom stated that the scalar
\begin{equation}
\label{eq:matklyu(1)}
z = \sum_{\lambda \in I}d_{\lambda} \chi_{\lambda}(v)sdim_{q}(V_{\lambda})
\end{equation}
should not vanish.
However, this automatically follows from the other axioms  \cite[Lem 8.20]{walk},\cite[Sec. 1]{tw}.
The proof for the $\mathbb{Z}_{2}$-graded case is exactly the same, thus we omit it here.  
However, note that the proof of Lemma \ref{eq:emotion}
for pseudo-modular Hopf algebras bears much similarity to it.

\begin{lemma}
For a $\mathbb{Z}_{2}$-graded modular Hopf algebra $A$, the complex constant $z$ defined by
\begin{equation}
\label{eq:craigemerson(2)}
z = \sum_{\lambda \in I}d_{\lambda} \chi_{\lambda}(v)sdim_{q}(V_{\lambda}),
\end{equation}
is non-zero.
\end{lemma}

Another consequence of the axioms of a modular Hopf algebra is that 
$d_{\mu^{*}} = d_{\mu}$ for all $\mu \in I$.  See \cite[Sec. 5.2]{rt} for a detailed proof.

Finally, we note for interest that Turaev and Wenzl defined a slightly more general class of ribbon Hopf
algebras that can be used in constructing $3$-manifold invariants.
They defined
{\emph{quasimodular Hopf algebras}} \cite[Sec. 2.1]{tw} 
so as to construct $3$-manifold invariants from the
quantum algebras $U_{q}(\mathfrak{g})$ associated with the $A, B, C$ and $D$ families 
of Lie algebras at even roots of unity 
without knowing whether the relevant modules were irreducible \cite{tw}.
The essential difference between quasimodular and modular Hopf algebras is that this 
irreducibility condition is relaxed for 
quasimodular Hopf algebras.

\end{subsection}

\begin{subsection}{Pseudo-modular Hopf algebras}
\label{subsec:allanbridges2}

We now introduce the notion of pseudo-modular Hopf algebras that play a central
role in this chapter.   
A {\emph{pseudo-modular Hopf algebra}}
is a $\mathbb{Z}_{2}$-graded ribbon Hopf algebra together with a collection of finite dimensional
representations satisfying slightly weaker conditions than those 
satisfied by modular Hopf algebras. 
We will prove in Section \ref{sec:rushingwater} that one can construct 
topological invariants of closed, connected, orientable $3$-manifolds from 
pseudo-modular Hopf algebras.

The essential difference between a 
pseudo-modular Hopf algebra and a modular Hopf algebra
is that relations (\ref{eq:waggahogga}) do not necessarily have a unique set of solutions in
pseudo-modular Hopf algebras, thus
the set of constants $\{d_{\nu} \in \mathbb{C} | \ \nu \in I\}$  is not necessarily unique. 
Consequently, one must independently prove that the 
$z$ in (\ref{eq:matklyu(1)}) is nonvanishing.
Examples of algebras for which the constants $\{d_{\nu} | \ \nu \in I\}$ are not unique are
$U_{q}^{(N)}(osp(1|2))$ at odd roots of unity \cite{z2} and $U_{q}^{(N)}(gl_{2})$ also at odd roots of
unity \cite{z3}.

\begin{definition}
\label{defn:pseudo}

Let $A$ be a $\mathbb{Z}_{2}$-graded ribbon Hopf algebra over $\mathbb{C}$
with universal $R$-matrix $R \in A \otimes A$ stated in (\ref{eq:tim(1)}), 
the element $u \in A$ stated in (\ref{eq:tim(2)}), 
and the even central element $v \in A$ stated in (\ref{eq:tim(3)}).  

Let $I$ be a finite index set with an involution $*:I \rightarrow I$  denoted by
$*(i) = i^{*}$, and let there exist a distinguished element $0 \in I$ satisfying $*(0)=0^{*}=0$.  
Let $\{V_{i} | \ i \in I\}$ be a set of  $A$-modules, 
where $V_{0}$ is the one-dimensional $A$-module, and let
there exist a set of $A$-linear isomorphisms
$$\big\{ \omega_{i}: \left(V_{i}\right)^{*} \rightarrow V_{i^{*}}| \ i \in I  \big\},$$
where $\omega_{0} = \mathrm{id}$ and $\left(V_{i}\right)^{*}$ is the dual $A$-module to $V_{i}$.

The $\mathbb{Z}_{2}$-graded ribbon Hopf algebra $A$ is 
said to be a pseudo-modular Hopf algebra if the following axioms are satisfied:
\begin{itemize}

\item[(I)]  The $A$-modules $\{V_{i} | \ i \in I\}$ are finite dimensional,  
mutually non-isomorphic and $sdim_{q}(V_{i}) \neq 0$ for all $i \in I$.

\item[(II)]  Let $\Gamma(L, \lambda)$
be a coloured directed $(1,1)$-ribbon tangle with $m$ components, and let
the ends of $\Gamma(L, \lambda)$ be directed
downwards and lie in a component coloured with $i \in I$.
Then the corresponding map
$F\big(\Gamma(L,\lambda)\big): V_{i} \rightarrow V_{i}$ 
acts by multiplication by a complex scalar.
Furthermore, the even central element $v \in A$ acts on each $A$-module $V_{i}, i \in I$,
as a complex scalar.

\item[(III)]
For each $i \in I$, let $\{ x_{j} \}$, $\{ x_{j}^{*} \}$ and 
$\{ x_{j}^{**} \}$ be bases of $V_{i}$, $(V_{i})^{*}$ and $((V_{i})^{*})^{*}$, 
respectively, such that
$\langle \langle x_{j}^{**}, x_{k}^{*} \rangle \rangle 
= \langle x_{k}^{*}, x_{j} \rangle = \delta_{j,k}$ and 
$[x_{j}^{**}] = [x_{j}^{*}] = [x_{j}]$, where 
$\langle \langle \cdot, \cdot \rangle \rangle: 
((V_{i})^{*})^{*} \times (V_{i})^{*} \rightarrow \mathbb{C}$ is the dual space pairing.
Then the map
$$(\omega_{i})^{*} \circ (\omega_{i^{*}})^{-1}: 
V_{i} \rightarrow ((V_{i})^{*})^{*},$$
is given by $x_{j} \mapsto (-1)^{[x_{j}]} u v^{-1} x_{j}^{**}, \ \forall j$.

\item[(IV)] For any $t \in \mathbb{N}$ and any sequence
$\theta = (i_{1}, i_{2}, \ldots, i_{t}) \in I^{\times t}$,
$$V_{i_{1}} \otimes V_{i_{2}} \otimes \cdots \otimes V_{i_{t}} 
= \left(\bigoplus_{i \in I}\left(V_{i}\right)^{\oplus n_{\theta}(i)} \right) \oplus {\cal{Z}}_{\theta},$$
where $n_{\theta}(i) \in \mathbb{Z}_{+}$ is the number of copies 
of the $A$-module $V_{i}$ in the direct sum and
${\cal{Z}}_{\theta}$ is a possibly vanishing $A$-module such that 
$str_{q}(a) = 0$ for any $A$-linear map 
$a \in End_{A}(V_{i_{1}} \otimes V_{i_{2}} \otimes \cdots \otimes V_{i_{t}})$ 
obtained by applying the
functor $F$ to a directed oriented $(m,m)$-ribbon tangle, where the module map 
$a \in End_{A}(V_{i_{1}} \otimes V_{i_{2}} \otimes \cdots \otimes V_{i_{t}})$ also satisfies
$a(V_{i_{1}} \otimes V_{i_{2}} \otimes \cdots \otimes V_{i_{t}}) \subseteq {\cal{Z}}_{\theta}$.

\item[(V)] Let $f=(f_{\lambda \nu})_{\lambda,\nu \in I}$ 
be a matrix with complex elements given by 
$$f_{\lambda \nu} = ({\mathrm{str}} \otimes {\mathrm{str}})
\left[(\pi_{\lambda} \otimes \pi_{\nu})(v^{-1}u\otimes v^{-1}u)R^{T}R\right].$$
Then there exists at least one set $\{d_{\nu} \in \mathbb{C} | \ \nu \in I\}$ 
of constants satisfying the relations
\begin{equation}
\label{eq:wowwee}
\chi_{\lambda}(v) sdim_{q}(V_{\lambda}) = 
\sum_{\nu \in I} d_{\nu} \chi_{\nu} (v^{-1})f_{\lambda \nu}, \hspace{5mm} \mbox{for all } \lambda \in I,
\end{equation}
where $d_{\nu^{*}} = d_{\nu}$ for all $\nu \in I$. 
Here $\pi_{\nu}$ is the $A$-representation furnished by $V_{\nu}$, and
$\chi_{\nu}(v) = \pi_{\nu}(v)$.

\item[(VI)]  
The scalar $\displaystyle{z = \sum_{\lambda \in I}d_{\lambda} \chi_{\lambda}(v)sdim_{q}(V_{\lambda})}$ is non-zero.  
\end{itemize}
\end{definition}

\end{subsection}

\end{section}

\begin{section}{Reshetikhin-Turaev invariant arising from pseudo-modular Hopf algebras}
\label{sec:rushingwater}
\markright{\text{Reshetikhin-Turaev invariant arising from pseudo-modular Hopf algebras}}

We will construct a topological invariant of a 
closed, connected, orientable $3$-manifold from each pseudo-modular Hopf
algebra following the general approach of Reshetikhin and Turaev \cite{rt} 
and Turaev and Wenzl \cite{tw}.
The main result of this section is Theorem \ref{th:markhopkins(a)}.

Let us fix a pseudo-modular Hopf algebra as defined in Definition \ref{defn:pseudo}.
Let $L \subset S^{3}$ be a framed oriented link with $m$ connected components.
Let $\lambda = (\lambda_{1}, \lambda_{2}, \ldots, \lambda_{m})$, $\lambda_{i} \in I$, 
and fix ${\cal{C}}(L,I)$ to be the set of all different $\lambda$.  Define
\begin{equation}
\label{eq:defsum}
\sum(L) = \sum_{\lambda \in {\cal{C}}(L,I)}  \prod_{i=1}^{m} d_{\lambda_{i}}F\big(\Gamma(L,\lambda)\big),
\end{equation}
where $\{d_{\nu} | \ \nu \in I\}$ is a set of constants satisfying Eq. (\ref{eq:wowwee}).
Note that $\sum(L)$ is a regular isotopy invariant of $L$ as each
$F\big(\Gamma(L,\lambda)\big)$ is a regular isotopy invariant of $L$.

We now prove that 
$\sum(L) = \sum(L')$ if $L$ and $L'$ are two links that are equivalent with respect 
to the $\kappa_{+}$ move.  
This $\sum(L)$ is a core part of Reshetikhin and Turaev's $3$-manifold invariants.

Consider the two framed oriented links $L$, $L' \subset S^{3}$ 
presented in the blackboard framing in Figure \ref{fig:figuremeout!}. 
The planar projections of $L$ and $L'$ are equivalent with respect to the Kirby move $\kappa_{+}$.  
In this figure, $T$ is an arbitrary oriented $(m,m)$-tangle represented by the rectangle, 
where the orientation of $T$ is
compatible with the orientations of $L \backslash T$ and $L' \backslash T$.  
\begin{figure}[htb]
\begin{center}
 \input{Kirby+move100.pstex_t}
\caption{Two links $L$, $L'$ that are equivalent with respect to the Kirby $\kappa_{+}$ move}  \label{fig:figuremeout!}
\end{center}
\end{figure}
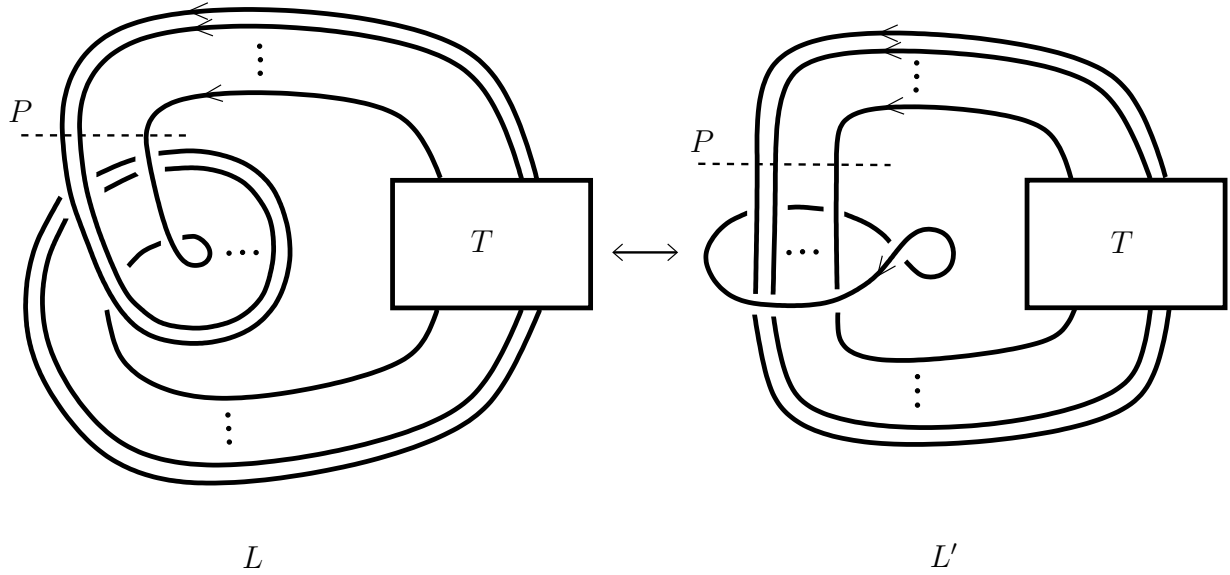
We colour the directed ribbon tangles $\Gamma(L)$ and $\Gamma(L')$
associated with $L$ and $L'$, respectively, as follows. 
Firstly, let $\Gamma(L_{i})$ (resp. $\Gamma(L_{i}')$) denote 
the directed $(1,1)$-ribbon tangle obtained from the $i^{th}$ tangle of 
$L \backslash T$ (resp. $L' \backslash T$) intersecting the line $P$ from the left
in both diagrams. 
Each of $\Gamma(L_{i})$ and $\Gamma(L_{i}')$ is coloured with $\lambda_{i} \in I$.  
Lastly, colour the unique component of
$\Gamma(L' \backslash T)$ isotopic to an unknotted annulus with $\xi \in I$.  
We denote the resulting coloured directed $(m, m)$-ribbon tangles as
$$\Gamma\big(L \backslash T,(\lambda_{1},\lambda_{2},\ldots,\lambda_{m})\big), \hspace{10mm}
\Gamma\big(L' \backslash T,(\lambda_{1},\lambda_{2},\ldots,\lambda_{m},\xi)\big),$$
respectively.
We have the following result.
\begin{proposition}
\label{prop:dragonmount}
For each $m \geq 1$, 
$$F\big(\Gamma(L\backslash T,(\lambda_{1},  \ldots, \lambda_{m}))\big) =
 \Delta^{(m-1)}(v) : V_{\lambda_{1}} \otimes \cdots \otimes V_{\lambda_{m}} 
\rightarrow V_{\lambda_{1}} \otimes \cdots \otimes V_{\lambda_{m}},$$
and
$$F\big(\Gamma(L' \backslash T,(\lambda_{1}, \ldots, \lambda_{m},\xi))\big) 
 = \chi_{\xi}(v^{-1}) \Delta^{(m-1)}(C_{\xi}): 
V_{\lambda_{1}} \otimes \cdots \otimes V_{\lambda_{m}} \rightarrow V_{\lambda_{1}} \otimes \cdots
\otimes V_{\lambda_{m}},$$
where $C_{\xi}$ is defined by (\ref{eq:theCmuequation(1)}).
\end{proposition}
\begin{proof}
The first claim follows from Remark \ref{rem:markhopkins(aaa)}.
For $m=1$, the second claim follows from Eq. (\ref{eq:markhopkins(bb)}),
which is still true as $v^{-1}$ acts as a scalar in the $A$-representation $\pi_{\xi}$, and
for $m \geq 2$  by induction.
\end{proof}

\begin{proposition}
\label{prop:craigemerson(2)}
Consider the two framed, oriented links $L$, $L'$ in Figure \ref{fig:figuremeout!} 
where $m \geq 1$ and $T$ is an arbitrary oriented (m,m)-tangle the orientation of which is 
compatible with that of $L \backslash T$ and $L' \backslash T$, then
\begin{equation}
\label{eq:majormainsexyequation}
\sum(L)=\sum(L').
\end{equation}
\end{proposition}
\begin{proof}
By definition,
\begin{eqnarray}
\sum(L') & = & \sum_{\lambda \in {\cal{C}}(L,I)} \sum_{\xi \in I} 
\prod_{i=1}^{m} d_{\lambda_{i}}
\mathrm{str}\bigg[\big(\pi_{\lambda_{1}}(v^{-1}u) \otimes \pi_{\lambda_{2}}(v^{-1}u)
\cdots \otimes \pi_{\lambda_{m}}(v^{-1}u)\big) 
\label{eq:thethunderbolt} \\
& & \hspace{10mm} \times
F\big( \Gamma (L' \backslash T, (\lambda_{1}, \lambda_{2}, \ldots, \lambda_{m}, \xi))\big) 
\circ F\big( \Gamma(T, (\lambda_{1}, \lambda_{2}, \ldots, \lambda_{m}))\big) \bigg], \nonumber
\end{eqnarray}
where in writing $\mathrm{str}$ we mean that we take the supertrace over 
$V_{\lambda_{1}} \otimes V_{\lambda_{2}} \otimes \cdots \otimes V_{\lambda_{m}}$.

Let $\theta = (\lambda_{1}, \lambda_{2}, \ldots, \lambda_{m}) \in I^{\times m}$.
From Axiom (IV) of a pseudo-modular Hopf algebra, 
$V_{\lambda_{1}} \otimes V_{\lambda_{2}} \otimes \cdots \otimes V_{\lambda_{m}} = 
\left( \bigoplus_{\xi \in I} \left(V_{\xi}\right)^{\oplus n_{\theta}(\xi)} \right) \oplus {\cal{Z}}_{\theta}$
for some non-negative constants  $n_{\theta}(\xi)$, 
which allows us to rewrite (\ref{eq:thethunderbolt}) as
\begin{eqnarray}
\lefteqn{
\sum_{\lambda \in {\cal{C}}(L,I)} \prod_{i=1}^{m} d_{\lambda_{i}}
\left( \sum_{\xi \in I} n_{\theta}(\xi)\mathrm{str}_{V_{\xi}}
 + \mathrm{str}_{{\cal{Z}}_{\theta}} \right) } \nonumber \\
& & \hspace{5mm} \times \Bigg[ \big(v^{-1}u\big) \circ 
\Bigg( \sum_{\zeta \in I} d_{\zeta} \chi_{\zeta}(v^{-1}) f_{\xi \zeta}/sdim_{q}(V_{\xi}) \Bigg)
\circ F\big( \Gamma(T, (\lambda_{1}, \lambda_{2}, \ldots, \lambda_{m}))\big) \Bigg],
\label{eq:markhopkins(g)}
\end{eqnarray}
where in writing $\mathrm{str}_{V_{\xi}}$
we mean that we take the supertrace over the submodule 
$V_{\xi} \subseteq V_{\lambda_{1}} \otimes V_{\lambda_{2}} \otimes \cdots \otimes V_{\lambda_{m}}$,
and  $\mathrm{str}_{{\cal{Z}}_{\theta}}$ has a similar meaning.
Axiom (V) allows us to rewrite the right hand side of (\ref{eq:markhopkins(g)}) as
$$
\sum_{\lambda \in {\cal{C}}(L,I)} \prod_{i=1}^{m} d_{\lambda_{i}}
\left( \sum_{\xi \in I} n_{\theta}(\xi)\mathrm{str}_{V_{\xi}}  + 
      \mathrm{str}_{{\cal{Z}}_{\theta}}\right) 
\bigg[ \big(v^{-1}u\big) \circ \left(\chi_{\xi}(v)\right) 
F\big( \Gamma(T, (\lambda_{1}, \lambda_{2}, \ldots, \lambda_{m}))\big) \bigg],
$$
which is equal to $\sum(L)$.
In the derivation of this equality, we used Eq. (\ref{eq:wowwee}) and Axiom (II).
Note that $\sum(L')=\sum(L)$ relies on the equality of the quantum supertraces 
of certain $A$-linear maps, not on an equality of the $A$-linear maps themselves.
\end{proof}
\begin{remark}

The proof of Proposition \ref{prop:craigemerson(2)} 
is similar to the proof of the corresponding theorem
for modular Hopf algebras.  See \cite[Par. 7.2]{rt} for details.
\end{remark}

\begin{proposition}
\label{prop:johnfaulkner(10000)}
Let $\Gamma(L,(\lambda_{1}, \ldots, \lambda_{m}))$ be a coloured, directed $(0,0)$-ribbon tangle, with
its $i^{th}$ component coloured by $\lambda_{i}$.
Let $\Gamma'(L,(\lambda_{1}, \ldots, \lambda_{m}))$ 
be the coloured, directed $(0,0)$-ribbon tangle obtained by changing the direction of the 
$i^{th}$ component of $\Gamma(L(\lambda_{1}, \ldots, \lambda_{m}))$.
Then 
$$F \big(\Gamma'(L,(\lambda_{1}, \ldots, \lambda_{i-1}, \lambda_{i}, \lambda_{i+1}, \ldots,\lambda_{m})) \big)
= F \big(\Gamma (L,(\lambda_{1}, \ldots, \lambda_{i-1}, (\lambda_{i})^{*}, \lambda_{i+1}, \ldots,
 \lambda_{m})) \big).$$
\end{proposition}
\begin{proof}
The proof is almost identical to the proof of \cite[Lem. 5.1]{rt}, thus omitted.
\end{proof}

We can now define the topological invariant.
Let $L = \bigcup_{i=1}^{m} L_{i} \subset S^{3}$ be an unoriented framed link with $m$ connected components,
and let $M_{L}$ be the closed, connected, orientable 3-manifold obtained by performing surgery on $S^{3}$ along $L$. 
We introduce some notation by writing
$A_{L} = (a_{ij})_{i,j=1}^{m}$ to mean the linking matrix of $L$, defined by 
\begin{itemize}
\item[(i)]  $a_{ii} = w(L_{i})$, the writhing number of $L_{i}$, for each $i$, 
\item[(ii)] $a_{ij} = \mathrm{lk}(L_{i},L_{j})$, if $i \neq j$.  
\end{itemize}
Note that $A_{L}$ is real and symmetric.  
We define $\sigma(A_{L})$ to be the number of non-positive eigenvalues of $A_{L}$. 
Following \cite{rt} we introduce the topological invariant ${\cal{F}}(M_{L})$.
\begin{theorem}
\label{th:markhopkins(a)}
Let $L$ be a framed link.
Then
$${\cal{F}}(M_{L}) = z^{-\sigma(A_{L})}\sum(L)$$
is a topological invariant of $M_{L}$, where in calculating $\sum(L)$ we assign an arbitrary orientation to $L$.
\end{theorem}
\begin{proof}
Observe that $\sum(L)$ does not depend on the orientation chosen for $L$.
Let $L$ have $m$ components and assign an orientation to each component of $L$.  Then
$$\sum(L) = \sum_{\lambda \in {\cal{C}}(L,I)} 
\prod_{j=1}^{m} d_{\lambda_{j}} F \big(\Gamma(L,(\lambda_{1}, \ldots, \lambda_{m}))\big).$$
Now let $L'$ be the link obtained by reversing the orientation of the $i^{th}$ component of $L$. 
By using Proposition \ref{prop:johnfaulkner(10000)} we have
\begin{equation}
\label{eq:biscuitbecstar(1)}
\sum(L')  =  \sum \prod_{j=1}^{m} d_{\lambda_{j}} 
F \big(\Gamma(L,(\lambda_{1}, \ldots, \lambda_{i-1}, (\lambda_{i})^{*}, 
\lambda_{i+1}, \ldots, \lambda_{m}))\big),
\end{equation}
where the right hand side of (\ref{eq:biscuitbecstar(1)})
is a $m$-fold summation over all $\lambda_{j} \in I$, $j \neq i$, and 
$(\lambda_{i})^{*} \in I$.
Since $d_{\lambda_{i}} = d_{(\lambda_{i})^{*}}$, we can rewrite $\sum(L')$ as
\begin{equation}
\label{eq:biscuitbecstarseadog(1)}
     \sum(L')   = \sum  \prod_{\stackrel{j=1}{j \neq i}}^{m} d_{\lambda_{j}} d_{(\lambda_{i})^{*}}
F \big(\Gamma(L,(\lambda_{1}, \ldots, \lambda_{i-1}, (\lambda_{i})^{*}, 
\lambda_{i+1}, \ldots, \lambda_{m}))\big).
\end{equation}
On the right hand side, we sum all the $\lambda_{k}$'s independently over $I$.  Hence
the right hand side of (\ref{eq:biscuitbecstarseadog(1)}) is equal to $\sum(L)$, that is
$$\sum(L') = \sum(L').$$

We now need to show that ${\cal{F}}(M_{L}) = {\cal{F}}(M_{\cal{L}})$ 
if the links $L$ and ${\cal{L}}$ are equivalent with respect to any of the Kirby moves.
Each Kirby move can be expressed as some composition of the 
$\kappa_{+}^{(0)}$, $\kappa_{-}^{(0)}$ and $\kappa_{+}$ Kirby moves \cite[Thm. 6.3]{rt}; 
we will show that ${\cal{F}}(M_{L}) = {\cal{F}}(M_{\cal{L}})$ 
if $L$ and ${\cal{L}}$ are equivalent with respect to any of these moves.

We firstly show that ${\cal{F}}(M_{L}) = {\cal{F}}(M_{L'})$ 
if the two links $L$ and $L'$ are equivalent with respect to the $\kappa_{+}^{(0)}$ move.  
Let ${\cal{O}}_{+1}$ denote the unknot with framing $+1$ given in Figure \ref{fig:uknots} 
and let $L' = L \cup {\cal{O}}_{+1}$ be a split link.  
\begin{figure}[hbt]
 \begin{center}
 \input{someunknots100.pstex_t}
\caption{The unknots ${\cal{O}}_{+1}$ and ${\cal{O}}_{-1}$} \label{fig:uknots}
\end{center}
\end{figure}

\noindent
It immediately follows from the definition (Eq. (\ref{eq:defsum})) 
that $\sum(L') = \sum(L)\sum({\cal{O}}_{+1})$, and
$$\sum({\cal{O}}_{+1}) = \sum_{\lambda \in I}
d_{\lambda}\chi_{\lambda}(v^{-1})sdim_{q}(V_{\lambda}) = 1,$$ 
where we have used the relation 
$f_{0 \lambda} = sdim_{q}(V_{\lambda})/sdim_{q}(V_{0}) = sdim_{q}(V_{\lambda})$.  
As $L'$ is a split link, $A_{L'}$ is the $(m+1) \times (m+1)$ matrix 
$A_{L'} = \big( (1) \oplus A_{L} \big) = \left(
\begin{array}{cccc}
1 & 0 & \cdots & 0  \\
0 &   &  &    \\
\vdots & & A_{L} &  \\
0 &   &       &
\end{array}
\right)$,
thus $\sigma(A_{L'}) = \sigma(A_{L})$ and ${\cal{F}}(M_{L}) = {\cal{F}}(M_{L'})$.

We now show that  ${\cal{F}}(M_{L}) = {\cal{F}}(M_{\widetilde{L}})$ if the two links
$L$ and $\widetilde{L}$ are equivalent with respect to the $\kappa_{-}^{(0)}$ move.  
Let ${\cal{O}}_{-1}$ denote the unknot with framing $-1$ given in Figure \ref{fig:uknots} and let
$\widetilde{L} = L \cup {\cal{O}}_{-1}$ be a split link, then 
$$\sum(\widetilde{L}) = \sum(L)\sum({\cal{O}}_{-1}) = \sum(L) \sum_{\lambda \in I}
d_{\lambda}\chi_{\lambda}(v)sdim_{q}(V_{\lambda}) = z \sum(L).$$
Now $A_{\widetilde{L}}$ is an $(m+1) \times (m+1)$ matrix given by
$A_{\widetilde{L}} = \big((-1) \oplus A_{L}\big) =
\left(
\begin{array}{cccc}
-1 & 0 & \cdots & 0  \\
0 &   &  &    \\
\vdots & & A_{L} &  \\
0 &   &       &
\end{array}
\right)$, 
 thus $\sigma(A_{\widetilde{L}}) = \sigma(A_{L})+1$
and ${\cal{F}}(M_{L}) = {\cal{F}}(M_{\widetilde{L}})$.

It remains to 
show that ${\cal{F}}(M_{L}) = {\cal{F}}(M_{L'})$ if the two links $L$ and $L'$ are
equivalent with respect to the $\kappa_{+}$ move.
Let $L$ and $L'$ be the two framed oriented links presented in the blackboard framing 
in Figure \ref{fig:figuremeout!}.  
The planar projections of $L$ and $L'$ are equivalent with respect to the $\kappa_{+}$ move. 
In Figure \ref{fig:figuremeout!}, $T$ is an arbitrary oriented $(m,m)$-tangle represented by a rectangle, 
the orientation of which is
compatible with the orientations of $L \backslash T$ and $L' \backslash T$.
From Proposition \ref{prop:craigemerson(2)}, $\sum(L')=\sum(L)$.  
We now show that $\sigma(A_{L'}) = \sigma(A_{L})$: $A_{L}$ is an $m \times m$ matrix and
$A_{L'}$ is an $(m+1) \times (m+1)$ matrix both of which we give below.  
In these matrices, $A_{T}$ is the $m \times m$ linking matrix of $T$.  We have
$$A_{L} = \left(  \begin{array}{rrrrr}
                -1 & -1 & -1 & \cdots & -1  \\
                -1 & -1 & -1 & \cdots & -1  \\
		-1 & -1 & -1 & \cdots & -1  \\
    \vdots & \vdots & \vdots &  \ddots    & \vdots \\
		-1 & -1 & -1 & \cdots & -1
\end{array}  \right) + A_{T},
$$
and $$A_{L'} = \left( \begin{array}{ccccc}
               1 & 1 & 1 & \cdots & 1 \\
	       1 & 0 & 0 & \cdots & 0 \\
	       1 & 0 & 0 & \cdots & 0 \\
    \vdots & \vdots & \vdots &  \ddots    & \vdots \\
	       1 & 0 & 0 & \cdots & 0
\end{array}
\right) + \big( (0) \oplus A_{T} \big),
$$
where 
$\big((0) \oplus A_{T}\big) = \left(
\begin{array}{cccc}
0 & 0 & \cdots & 0  \\
0 &   &  &    \\
\vdots & & A_{T} &  \\
0 &   &       &
\end{array}
\right)$ is an $(m+1) \times (m+1)$-matrix.
Let $$X = \left(  \begin{array}{rcccc}
		1  & 0 & 0 & \cdots & 0  \\
		-1 & 0 & 0 & \cdots & 0  \\
		-1 & 0 & 0 & \cdots & 0  \\
		\vdots & \vdots & \vdots & \ddots & \vdots \\
		-1 & 0 & 0 & \cdots & 0
\end{array}
\right) + \big( (0) \oplus I_{m} \big),
$$
then $$X A_{L'} X^{T} = \big( (1) \oplus A_{L} \big)  =\left(
\begin{array}{cccc}
1 & 0 & \cdots & 0  \\
0 &   &  &    \\
\vdots & & A_{L} &  \\
0 &   &       &
\end{array}
\right),  $$
where $X^{T}$ is the transpose of $X$ and $I_{m}$ is the $m \times m$ identity matrix. 
Now the matrices $A_{L'}$ and $\big((1) \oplus A_{L}\big)$ are real and symmetric.
As $X$ is invertible, $A_{L'}$ and $\big((1) \oplus A_{L}\big)$ are congruent.  
Congruent real symmetric matrices have the same numbers of 
positive, zero and negative eigenvalues \cite[Thms. 8.7, 8.9]{f}, thus  
$\sigma(A_{L'}) = \sigma(A_{L})$ and  ${\cal{F}}(M_{L}) = {\cal{F}}(M_{L'})$.

It follows that ${\cal{F}}(M_{L}) = {\cal{F}}(M_{\cal{L}})$ 
if the two links $L$ and ${\cal{L}}$ are equivalent with respect to any of the Kirby moves, and thus 
${\cal{F}}(M_{L})$ is indeed a topological invariant of $M_{L}$.
\end{proof}

Note from our definition that ${\cal{F}}(M_{L})$ is normalised to $1$ on $S^{3}$.  

Explicitly calculating ${\cal{F}}(M_{L})$ for a 
particular closed, connected, orientable $3$-manifold $M_{L}$ 
is quite a difficult problem in general, and this is true
for all the Reshetikhin-Turaev topological invariants.
While topological invariants have been calculated for various classes of $3$-manifolds 
(eg for the Lens spaces from quotients of the quantum algebras arising from the 
$A, B, C,$ and $D$ series of
Lie algebras \cite{z4} and the $G_{2}, F_{4}, E_{8}$ Lie algebras \cite{z6} at odd roots of unity),
 we are only aware of one other collection of invariants that have been
 explicitly calculated
for a substantial number of $3$-manifolds \cite[Chap. 14]{kl}.
These invariants were constructed from $U_{q}^{(N)}(sl_{2})$ where $N \equiv 0 \pmod{4}$.

\end{section}

\begin{section}{Invariants of 3-manifolds from $U_{q}^{(N)}(osp(1|2n))$}
\markright{\text{Invariants of 3-manifolds from $U_{q}^{(N)}(osp(1|2n))$ }}
\label{chapter4:sectionlabel(Invariants_of_3-manifolds_arising_from)}

The following theorem is one of the main results of this thesis.
\begin{theorem}
\label{eq:antonygreenone}
Set $q = \exp{(2 \pi i/N)}$, where $N \geq 6$ satisfies $N \equiv 2 \pmod{4}$.  
For the $\mathbb{Z}_{2}$-graded ribbon Hopf algebra 
$U_{q}^{(N)}(\mathfrak{g}) = U_{q}^{(N)}(osp(1|2n))$, fix the set of non-isomorphic  
$U_{q}^{(N)}(\mathfrak{g})$-modules $\left\{ V_{\lambda} | \ \lambda \in \Lambda_{N}^{+} \right\}$
defined in Definition \ref{lem:representationsyippe} by:
\begin{equation}
\label{eq:troybronte(a)}
V_{\lambda} = \tilde{p}_{i}^{t}[\lambda] \big(V^{\otimes t}\big),
\end{equation}
where $t \in \mathbb{N}$ is  given in Definition \ref{lem:representationsyippe} and
$\tilde{p}_{i}^{t}[\lambda] \in End_{U_{q}^{(N)}(\mathfrak{g})}(V^{\otimes t})$ is 
defined in Eq. (\ref{eq:markhopkins(yyy)}).
Fix an involution $*: \Lambda_{N}^{+} \rightarrow \Lambda_{N}^{+}$ by 
$* = id_{\Lambda_{N}^{+}}$ and
define a set of constants $\{ d_{\lambda} \in \mathbb{C} | \ \lambda \in \Lambda_{N}^{+}\}$ by
\begin{equation}
\label{eq:markhopkins(b)}
d_{\lambda} = d_{0} sdim_{q}(V_{\lambda}),
\end{equation}
where $d_{0} = \Omega Q(0)$ and 
\begin{equation}
\label{eq:markhopkins(c)}
\Omega = \frac{2^{n}t^{n} q^{n^{3}-n/2} }{\left[ (1+i) \sqrt{N}\right]^{n} }, \hspace{10mm}
t = e^{ \pi i/2N},
\end{equation}
\begin{equation}
\label{eq:markhopkins(d)}
Q(0) = \prod_{\alpha \in \overline{\Phi}_{0}^{+}  } \left( q^{(\alpha, \rho)} - q^{-(\alpha, \rho)}\right)
\prod_{\beta \in \Phi_{1}^{+}} \left( q^{(\beta, \rho)} + q^{-(\beta, \rho)}\right).
\end{equation}

Let $z=(-i)^{n} q^{2n^{3}-n} t^{2n}$. 
Let $L \subset S^{3}$ be a framed unoriented link with $m$ connected components and let $M_{L}$ be the closed, connected,
orientable $3$-manifold obtained by performing surgery on $S^{3}$ along $L$.
Let $A_{L}$ be the linking matrix of $L$ and let  
$\sigma(A_{L})$ denote the number of non-positive eigenvalues of $A_{L}$.
Then
$${\cal{F}}(M_{L}) = z^{-\sigma(A_{L})} \sum_{\lambda \in {\cal{C}}(L, \Lambda_{N}^{+})} 
\prod_{i=1}^{m} d_{\lambda_{i}} F \big(\Gamma(L, \lambda)\big),$$
is a topological invariant of $M_{L}$.
\end{theorem}
\begin{proof}
This is an immediate consequence of Theorems \ref{th:markhopkins(a)} and
\ref{theorem:bigbigbig}.
\end{proof}

\begin{theorem}
\label{theorem:bigbigbig}
Let $N \geq 6$ satisfy $N \equiv 2 \pmod{4}$, then $U_{q}^{(N)}(osp(1|2n))$ and the following data
give rise to a pseudo-modular Hopf algebra:
\begin{itemize}
\item[(i)]
the set $\{V_{\lambda} | \ \lambda \in \Lambda_{N}^{+} \}$ of non-isomorphic
$U_{q}^{(N)}(osp(1|2n))$-modules  given in (\ref{eq:troybronte(a)}),
\item[(ii)] an involution $*: \Lambda_{N}^{+} \rightarrow \Lambda_{N}^{+}$ defined by
$* = id$,
\item[(iii)] 
an isomorphism $\omega: V^{*} \rightarrow V$ where we fix $\omega^{-1}$ to be
the bijective homogeneous map of degree $0$,
$T \in End_{U_{q}^{(N)}(\mathfrak{g})}(V, V^{*})$, defined in Eq. (\ref{eq:johnfaulkner(20)}),
\item[(iv)]
a set of constants $\{ d_{\lambda} \in \mathbb{C} | \ \lambda \in \Lambda_{N}^{+}\}$ 
defined by Eq. (\ref{eq:markhopkins(b)}).
\end{itemize}
\end{theorem}

\begin{remark}
Unfortunately $U_{q}^{(N)}(\mathfrak{g})$ does not have a pseudo-modular Hopf algebra structure when 
$N \equiv 0 \pmod{4}$.  See Theorem \ref{lem:dalektron4} for details.
\end{remark}

We now prove that all the six axioms of pseudo-modular Hopf algebras are satisfied by
$U_{q}^{(N)}(osp(1|2n))$ when $N \geq 6$ satisfies $N \equiv 2 \pmod{4}$.
\begin{subsection}{Proof of Axioms (I)--(IV)}
\label{eq:brontetheevilcat(1)}
In the proof of this theorem we let $N \geq 6$ satisfy $N \equiv 2 \pmod{4}$.
\begin{itemize}
\item[(i)]  
The proof that Axiom (I) is satisfied is contained in Chapter \ref{chap2A:titlelabel}.

\item[(ii)]
To prove that Axiom (II) is satisfied, we consider all of the $U_{q}^{(N)}(\mathfrak{g})$-modules
$V_{\lambda}$, $\lambda \in \Lambda_{N}^{+}$, given in the data above, as being defined by
$V_{\lambda} = \tilde{p}^{t}_{i}[\lambda] (V^{\otimes t})$ for some $t$.

For each such module, there exists an irreducible $U_{q}(\mathfrak{g})$-module $V_{\lambda}^{gen}$
with integral dominant highest weight $\lambda$, 
for all non-zero $q$ that are not roots of unity, defined by
$V_{\lambda}^{gen} = \tilde{p}^{t}_{i}[\lambda]^{gen} (V^{gen})^{\otimes t}$,
where $\tilde{p}^{t}_{i}[\lambda]^{gen}$ is an element of the  
algebra ${\cal{C}}_{t}^{gen}$ over $\mathbb{C}$ generated by
$\left\{ \check{\cal{R}}_{i} \in End_{U_{q}(\mathfrak{g})}(V^{gen})^{\otimes t} | 
 \ i=1, \ldots, t-1 \right\}$.

Recall from the proof of Lemma \ref{lem:likj}
that we obtain the matrix $\tilde{p}^{t}_{i}[\lambda]$ 
by fixing $q$ to the appropriate root of unity in the matrix
$\tilde{p}^{t}_{i}[\lambda]^{gen}$.
We wish to show that the even central element $v \in U_{q}^{(N)}(\mathfrak{g})$ acts 
as a scalar on the
$U_{q}^{(N)}(\mathfrak{g})$-module $V_{\lambda}$.
To do this, recall from Lemma \ref{lem:markhoopkins(a)} that
there exists an invertible even element 
$v_{\lambda} \in \overline{U}^{+}_{q}(\mathfrak{g})$ 
that acts as the scalar $q^{-(\lambda + 2\rho, \lambda)}$ (here $q$ is generic)
on the finite dimensional irreducible
$U_{q}(\mathfrak{g})$-module $V_{\lambda}^{gen}$.
Furthermore,
$\Delta^{(t-1)}(v_{\lambda})$ acts as the same scalar  
on the finite dimensional irreducible
$U_{q}(\mathfrak{g})$-module $V_{\lambda}^{gen} \subseteq V^{\otimes t}$
defined by $\tilde{p}^{t}_{i}[\lambda]^{gen}: (V^{gen})^{\otimes t} \rightarrow V_{\lambda}^{gen}$.

Consider the matrix $\pi^{\otimes t} \big( \Delta^{(t-1)}(v) \big)$ where
$v \in U_{q}^{(N)}(\mathfrak{g})$.
We can obtain this by
specialising $q$ to the appropriate root of unity in  
$(\pi^{gen})^{\otimes t} \big( \Delta^{(t-1)}(v_{\lambda}) \big)$,
where $v_{\lambda} \in \overline{U}^{+}_{q}(\mathfrak{g})$.
Thus $\pi^{\otimes t} \big( \Delta^{(t-1)}(v) \big)$
also acts as a scalar on $V_{\lambda}$, and this scalar is precisely
$q^{-(\lambda+2\rho, \lambda)}$ where $q = \exp{(2 \pi i/N)}$.

Now let $\Gamma(L, \lambda)$
be a coloured directed $(1,1)$-ribbon tangle with $m$ components, and let
the ends of $\Gamma(L, \lambda)$ be directed
downwards and lie in a component coloured with $i \in I$.
We want to show that the map
$F\big(\Gamma(L,\lambda)\big): V_{i} \rightarrow V_{i}$ 
acts as the multiplication by a scalar.
Needless to say, $F\big(\Gamma(L,\lambda)\big)$ can be expressed in terms of the universal $R$-matrix
of $U_{q}^{(N)}(\mathfrak{g})$.

We embed $V_{i}$ in the $U_{q}^{(N)}(\mathfrak{g})$-module $V^{\otimes t}$.
As $\tilde{p}^{t}_{i}[\lambda]$ belongs to ${\cal{C}}_{t}$, under the embedding
$F\big(\Gamma(L,\lambda)\big)$ becomes an element of ${\cal{C}}_{t}$.

To determine the action of $F\big(\Gamma(L,\lambda)\big)$ on 
$V_{i} \subseteq V^{\otimes t}$, 
we take the element $f^{gen} \in {\cal{C}}^{gen}_{t}$ that corresponds to
$F\big(\Gamma(L,\lambda)\big)$ at all generic $q$.
Note that we obtain $F\big(\Gamma(L,\lambda)\big)$ by specialising $q$ to $\exp{(2 \pi i/N)}$ in
$f^{gen}$.  At all generic $q$, $f^{gen}$  acts as the
multiplication by a scalar on the irreducible $U_{q}(\mathfrak{g})$-module
$V_{i}^{gen} \subseteq (V^{gen})^{\otimes t}$.  
We then just take the limit of $f^{gen}$ as $q$ goes to
$\exp{(2 \pi i/N)}$.

\item[(iii)]
To prove that Axiom (III) is satisfied,  firstly note that
$V_{\lambda} \cong \tilde{p}^{t}_{i}[\lambda] \big(V^{\otimes t}\big)$ for some $t$
for each $\lambda \in \Lambda_{N}^{+}$.
Also note that $\lambda^{*} = \lambda$ for each $\lambda \in \Lambda_{N}^{+}$. 
Thus we need only show that Axiom (III) is satisfied for the map
$\omega^{*} \circ \omega^{-1}: V \rightarrow (V^{*})^{*}$, where
$\omega: V^{*} \rightarrow V$ and $\omega^{*}: V^{*} \rightarrow (V^{*})^{*}$ are isomorphisms.
Note that we fix $\lambda^{*} = \lambda$ as $(V_{\lambda})^{*} \cong V_{\lambda}$ for each 
$\lambda \in \Lambda_{N}^{+}$, which is proved in Chapter \ref{chap2A:titlelabel}.

Let $\{ v_{i} | -n \leq i \leq n \}$ be the basis of $V$ given in Lemma \ref{lem:manydimroot}
and let $\{ v_{i}^{*} \}$ and $\{ v_{i}^{**} \}$ be bases of $V^{*}$ and $(V^{*})^{*}$, 
respectively, such that
$\langle \langle v_{i}^{**}, v_{j}^{*} \rangle \rangle 
= \langle v_{j}^{*}, v_{i} \rangle = \delta_{i,j}$ and
$[v_{i}^{**}] = [v_{i}^{*}] = [v_{i}]$ for each $i$,
where 
$\langle \langle \cdot, \cdot \rangle \rangle : (V^{*})^{*} \times V^{*} \rightarrow \mathbb{C}$ is
the dual space pairing.
We want to show that the map
$\omega^{*} \circ \omega^{-1}: V \rightarrow (V^{*})^{*}$
has the following action:
\begin{equation} 
\label{eq:markhopkins(zzzzz)}
(\omega^{*} \circ \omega^{-1})(v_{i}) = (-1)^{[v_{i}]} u v^{-1} v_{i}^{**},
\hspace{10mm} -n \leq i \leq n.
\end{equation}

For $U_{q}^{(N)}(\mathfrak{g})$, $u v^{-1} = K_{2\rho}$, with $K_{2\rho}$ a product of $K_{i}$'s
satisfying $K_{2\rho} e_{i} K_{2\rho}^{-1} = q^{(2 \rho, \alpha_{i})}e_{i}$, $ \ \forall i$,
where 
$2 \rho  = \sum_{i=1}^{n} (2n-2i+1) \epsilon_{i}$.
The isomorphism $\omega^{-1}: V \rightarrow V^{*}$ is precisely the bijective map  
$T \in End_{U_{q}^{(N)}(\mathfrak{g})}(V, V^{*})$ defined in Eq. (\ref{eq:johnfaulkner(20)}).
The map $\omega^{*}: V^{*} \rightarrow (V^{*})^{*}$ is given by
$$\begin{array}{ll}
v_{i}^{*}  \mapsto (-1)^{i-1} q^{i-1} v_{-i}^{**}, \hspace{10mm} &
              v_{0}^{*}  \mapsto (-1)^{n-1} q^{n} v_{0}^{**}, \\
	      v_{-i}^{*} \mapsto (-1)^{i}   q^{2n-i} v_{i}^{**},  & 
	      1 \leq i \leq n.
	      \end{array}$$
Thus
$$(\omega^{*} \circ \omega^{-1}):
              v_{i}  \mapsto -q^{2n-2i+1} v_{i}^{**}, \hspace{10mm}
	      v_{0}  \mapsto v_{0}^{**}, \hspace{10mm}
	      v_{-i} \mapsto -q^{-(2n-2i+1)} v_{-i}^{**}, \hspace{5mm}
	      1 \leq i \leq n,$$
which shows that (\ref{eq:markhopkins(zzzzz)}) is true for all $i$.

\item[(iv)]
Axiom (IV) directly follows from the
tensor product theorems: Theorems \ref{th:firsttensor} and \ref{th:secondtensor}.

\end{itemize}
\end{subsection}

\begin{subsection}{Proof of Axiom (V)}

We now find a 
set of constants $\{d_{\lambda} \in \mathbb{C} | \ \lambda \in \Lambda_{N}^{+}\}$
satisfying the relations (\ref{eq:wowwee}).   
To do this, we need first to compute the $f_{\mu \lambda}$ for all $\mu, \lambda \in \Lambda_{N}^{+}$.

Recall Lemma \ref{lem:xylophone}.  
Because of part (ii) in Subsection \ref{eq:brontetheevilcat(1)}, $v^{-1}$ acts on $V_{\lambda}$
as the multiplication by the scalar $\chi_{\lambda}(v^{-1}) = q^{(\lambda + 2\rho, \lambda)}$
and $C_{\lambda}$ also acts on $V_{\mu}$ as the multiplication by a scalar
which we denote by $\chi_{\mu}(C_{\lambda})$.  It immediately follows that
$$f_{\mu \lambda} = \chi_{\mu} (C_{\lambda}) sdim_{q}(V_{\mu}), \hspace{10mm} 
\forall \mu, \lambda \in \Lambda_{N}^{+}.$$

Using calculations similar to those in 
Subsection \ref{subsec:hopefullylast(2)} and \cite[Lem. 2]{zg}, we obtain
$$\chi_{\mu}(C_{\lambda}) = \mathrm{str}\big( \pi_{\lambda} (K_{2 \mu + 2\rho}) \big),$$
where $\pi_{\lambda}$ is the representation of $U_{q}^{(N)}(\mathfrak{g})$ furnished by the module
$V_{\lambda}$.  Here $K_{2 \mu + 2\rho}$ is a product of $K_{i}$'s such that
$K_{2 \mu + 2\rho} e_{i} K_{2 \mu + 2\rho}^{-1} = q^{2(\mu + \rho, \alpha_{i})} e_{i}$, $\ \forall i$.

Now $V_{\lambda}$ can always be embedded in $V^{\otimes t}$ as a direct summand.
Let $\tilde{p}^{t}_{i}[\lambda]$ be the projection operator mapping $V^{\otimes t}$ onto $V_{\lambda}$.
Then
$$\chi_{\mu}(C_{\lambda}) = 
  str_{V^{\otimes t}} \big( K_{2 \mu + 2\rho} \cdot \tilde{p}^{t}_{i}[\lambda] \big).$$
As we have argued repeatedly, the right hand side can be obtained by first evaluating the
corresponding quantity at generic $q$, then specialising $q$ to a root of unity.
This way we obtain, with the help of the supercharacter formula for irreducible 
$osp(1|2n)$-representations, 
\begin{equation}
\label{chapter2:bussiethedealer(1)} 
\chi_{\mu}(C_{\lambda}) = (-1)^{[\lambda]}
\frac{\sum_{\sigma \in {\cal{W}}} \epsilon'(\sigma) q^{(2\mu + 2\rho,\sigma(\lambda + \rho))}}
{\sum_{\sigma \in {\cal{W}}}\epsilon'(\sigma)q^{(2\mu + 2\rho,\sigma(\rho))}}, \hspace{10mm}
\lambda, \mu \in \Lambda_{N}^{+},
\end{equation}
where ${\cal{W}}$ is the Weyl group of $osp(1|2n)$.  
See Appendix \ref{chap:appendixC} for further information.
It is important to observe that for all $\mu \in \Lambda_{N}^{+}$, the denominator of the above formula
is always non-zero.

For convenience, let us introduce some notation.
Recall that $X = \bigoplus_{i=1}^{n} \mathbb{Z} \epsilon_{i} \subset H^{*}$.
Let $S: X \times X \rightarrow \mathbb{C}$ and 
$Q: X \rightarrow \mathbb{C}$ be two mappings defined by
\begin{eqnarray}
(\lambda, \mu) & \mapsto &
S_{\lambda,\mu} = (-1)^{[\lambda]}\sum_{\sigma \in
{\cal{W}}}\epsilon'(\sigma)q^{2(\lambda + \rho, \sigma(\mu + \rho))}, \label{supereq:oneone} \\
\mu & \mapsto & Q(\mu) = \sum_{\sigma \in {\cal{W}}}\epsilon'(\sigma)q^{2(\mu +
\rho,\sigma(\rho))} \label{supereq:twotwo}.
\end{eqnarray}
Then we have the following result.
\begin{lemma}
\label{lem:properties}
$S_{\lambda,\mu}$ has the following properties:
\begin{itemize}
\item[(i)]  $S_{\mu,\lambda} = (-1)^{[\lambda] + [\mu]}S_{\lambda,\mu}$,
\item[(ii)] $S_{\lambda,w(\mu + \rho)-\rho} = \epsilon'(w)S_{\lambda,\mu}$, $w \in {\cal{W}}$.
\end{itemize}
\end{lemma}
\begin{proof}
\begin{itemize}
\item[(i)] 
$\displaystyle{ S_{\mu,\lambda}  =  (-1)^{[\mu]}\sum_{\sigma \in {\cal{W}}}
\epsilon'(\sigma)q^{2(\mu + \rho,\sigma(\lambda + \rho))}
 =  (-1)^{[\lambda] + [\mu]}S_{\lambda,\mu}}.$
\item[(ii)] 
$\displaystyle{ S_{\lambda,w(\mu + \rho)-\rho} = 
(-1)^{[\lambda]}\sum_{\sigma \in {\cal{W}}}\epsilon'(\sigma)
q^{2(\sigma(\lambda + \rho),w(\mu + \rho))} 
 = (-1)^{[\lambda]}\epsilon'(w)\sum_{\sigma \in {\cal{W}}}\epsilon'(\sigma)
 q^{2(\sigma(\lambda + \rho),\mu + \rho)} }.$
 \end{itemize}
\end{proof}
Restricting the domain of $S$ to $\Lambda_{N}^{+} \times \Lambda_{N}^{+}$ and the domain of $Q$ to
$\Lambda_{N}^{+}$ gives mappings with non-empty image.
We can rewrite (\ref{eq:wowwee}) as
\begin{equation}
\label{eq:sigeq1}
\chi_{\mu}(v) = \sum_{\lambda \in \Lambda_{N}^{+}} d_{\lambda} \chi_{\lambda}(v^{-1})
\left(S_{\lambda,\mu}/Q(\mu)\right), \hspace{10mm} \forall \mu \in \Lambda_{N}^{+}.
\end{equation}
The invertibility of the matrix 
$\big(S_{\lambda,\mu}/Q(\mu)\big)_{\lambda,\mu \in \Lambda_{N}^{+}}$ is unknown. 
However, $U_{q}^{(N)}(\mathfrak{g})$ may still satisfy Axiom (V) of a 
pseudo-modular Hopf algebra.  
This is the advantage of working with pseudo-modular Hopf algebras.

We shall now work towards proving Theorem \ref{theorem:bigbigbig}.  
Fix $N \equiv 2 \pmod{4}$ where $N \geq 6$, 
and also define $S_{\lambda,\mu}'=(-1)^{[\lambda]}S_{\lambda,\mu}$.  
We will solve the following set of linear equations for the constants
$\{d_{\lambda}' \in \mathbb{C} | \ \lambda \in \Lambda_{N}^{+}\}$ 
where $d_{\lambda}' = (-1)^{[\lambda]}d_{\lambda}$:
\begin{equation}
\label{eq:superQ}
Q(\mu)\chi_{\mu}(v) = 
\sum_{\lambda \in \Lambda_{N}^{+}}d_{\lambda}'\chi_{\lambda}(v^{-1}) S_{\lambda,\mu}', \hspace{10mm} 
\forall \mu \in \Lambda_{N}^{+}.
\end{equation}
Note that $Q(\mu) \neq 0$ for all $\mu \in \Lambda_{N}^{+}$ and the set
(\ref{eq:superQ}) of equations is identical to the set (\ref{eq:sigeq1}) of equations.

Define $X_{N}=X/NX$ and let
$p:X \rightarrow X_{N}$ be the canonical projection defined by 
$$p(\nu) = \sum_{i=1}^{n} \nu_{i} \epsilon_{i} + NX, \hspace{10mm}   
\nu = \sum_{i=1}^{n}\nu_{i}\epsilon_{i} \in X.$$
We distinguish between elements of $X$ and $X_{N}$ by writing $\lambda \in X$ and 
$\overline{\lambda}=p(\lambda)\in X_{N}$.
Note that $q^{(\lambda, \lambda + 2\rho)}$ makes sense for all $\lambda \in X$. 
Furthermore, if $p(\lambda) = p(\mu)$, then $q^{(\lambda, \lambda + 2\rho)}= q^{(\mu,\mu + 2\rho)}$, 
thus we can regard $q^{(\lambda,\lambda + 2\rho)}$ as a function $X_{N} \rightarrow \mathbb{C}$
defined by $p(\lambda) \mapsto q^{(\overline{\lambda},\overline{\lambda} + 2\rho)}$, where to evaluate
$q^{(\overline{\lambda}, \overline{\lambda} + 2\rho)}$ we take any representative 
$\lambda \in X$ of $\overline{\lambda}$ and calculate $q^{(\lambda, \lambda + 2\rho)}$.  
Similarly we can regard $S'_{\lambda,\mu}$ as a function $X_{N} \times X_{N} \rightarrow \mathbb{C}$.

We will encounter expressions of the form
$\sigma(\overline{\lambda})$, where $\sigma \in {\cal{W}}$, in our calculations.  
To be completely clear, 
we note that in writing $\sigma(\overline{\lambda})$ we actually mean
$p\big(\sigma(\lambda)\big)$, where $\lambda \in X$ is any representative of $\overline{\lambda}$ in $X_{N}$.

Following $\cite{z4}$, we introduce the auxiliary equations:
\begin{equation}
\label{eq:Qeq}
Q(\overline{\mu})\chi_{\overline{\mu}}(v) = \sum_{\overline{\lambda} \in
X_{N}}x_{\overline{\lambda}}\chi_{\overline{\lambda}}(v^{-1})
S_{\overline{\lambda},\overline{\mu}}',
\end{equation} 
where $\overline{\mu}$ varies over all elements of $p(\Lambda_{N}^{+})$
and we choose $\overline{\mu}$ to be the representative of $p(\mu)$ in $\Lambda_{N}^{+}$.
Note that the sum in the right hand side of (\ref{eq:Qeq}) is over all distinct 
elements of $X_{N}$.

We would like to solve (\ref{eq:Qeq}) for the $x_{\overline{\lambda}}$. 
Let us try the following solution: 
$x_{\overline{\lambda}} \equiv c q^{-(\overline{\lambda},2\rho)}$, 
where $c \in \mathbb{C}$ is some nonzero constant yet to be determined.  
Substituting this into (\ref{eq:Qeq}) gives
\begin{equation}
\label{eq:Q1eq}
Q(\overline{\mu}) q^{-(\overline{\mu},\overline{\mu}+2\rho)} = 
\sum_{\overline{\lambda} \in X_{N}}c 
q^{(\overline{\lambda},\overline{\lambda})}
\sum_{\sigma \in {\cal{W}}} \epsilon'(\sigma)q^{2\left(\overline{\lambda} + \rho, \sigma(\overline{\mu} + \rho)\right)}.
\end{equation}
We now aim to obtain an expression for $c$ independent of $\overline{\mu}$.
To help do this, we apply the following mapping to each term on the right hand side of (\ref{eq:Q1eq}) 
for each fixed pair $(\overline{\mu},\sigma) \in p(\Lambda_{N}^{+}) \times {\cal{W}}$:
$$\overline{\lambda} \mapsto \sigma(\overline{\lambda}) - \sigma(\overline{\mu}) \in X_{N}.$$
\begin{lemma}
\label{lem:proj2}
Let $(\mu,\sigma) \in \Lambda_{N}^{+} \times {\cal{W}}$ be a fixed pair.  
Then for any $\lambda$, $\lambda' \in X$,
\begin{equation}
\label{eq:joey}
p\big(\sigma(\lambda) - \sigma(\mu)\big) = p\big(\sigma(\lambda')-\sigma(\mu)\big)
 \hspace{5mm} \mbox{ if and only if }  \hspace{5mm} p(\lambda) = p(\lambda').
\end{equation}  
\end{lemma}
\begin{proof}
Write $\lambda = (\lambda_{1}, \lambda_{2}, \ldots, \lambda_{n})$ and
 $\lambda' = (\lambda'_{1}, \lambda'_{2}, \ldots, \lambda'_{n})$.
Fix $\sigma \in {\cal{W}}$ and assume that $p(\lambda) = p(\lambda')$, then
$$\lambda_{i} \equiv \lambda_{i}' \pmod{N}, \hspace{5mm} \mbox{for each } i=1, 2, \ldots, n,$$
and $\big( \sigma(\lambda) \big)_{i} \equiv \big(\sigma(\lambda')\big)_{i} \pmod{N}$ for all $i=1, 2, \ldots, n$. 
It follows that
$\big(\sigma(\lambda)-\sigma(\mu)\big)_{i} \equiv \big(\sigma(\lambda')-\sigma(\mu)\big)_{i} \pmod{N}$ 
for all $i=1, 2, \ldots, n$, and that
$$p\big(\sigma(\lambda)-\sigma(\mu)\big)=p\big(\sigma(\lambda')-\sigma(\mu)\big).$$

Now let us assume that
 $p\big(\sigma(\lambda)-\sigma(\mu)\big)=p\big(\sigma(\lambda')-\sigma(\mu)\big)$, then
clearly
$$\big(\sigma(\lambda)-\sigma(\mu)\big)_{i} \equiv \big(\sigma(\lambda')-\sigma(\mu)\big)_{i} \pmod{N}, \hspace{10mm} 
i=1, 2, \ldots, n,$$
and
\begin{equation}
\label{eq:johnedwards(a)}
\big(\sigma(\lambda)\big)_{i} \equiv \big(\sigma(\lambda')\big)_{i} \pmod{N}, \hspace{10mm} 
i=1, 2, \ldots, n.
\end{equation}
From Eq. (\ref{eq:johnedwards(a)}) we obtain
$$\lambda_{i} \equiv \lambda'_{i} \pmod{N}, \hspace{10mm} i=1, 2, \ldots, n,$$
and thus $p(\lambda) = p(\lambda')$ as desired. 
\end{proof}
It follows from Lemma \ref{lem:proj2} that for each fixed pair 
$(\overline{\mu}, \sigma) \in p(\Lambda_{N}^{+}) \times {\cal{W}}$, the 
following component of the right hand side of (\ref{eq:Q1eq})
is invariant under the mapping 
$\overline{\lambda} \mapsto \sigma(\overline{\lambda}) - \sigma(\overline{\mu}) \in X_{N}$:
\begin{equation}
\label{eq:ASDF}
c \epsilon'(\sigma)\sum_{\overline{\lambda} \in X_{N}}
q^{(\overline{\lambda},\overline{\lambda})} q^{2(\overline{\lambda}+\rho,\sigma(\overline{\mu}+\rho))}.
\end{equation}
Applying this mapping to (\ref{eq:ASDF}) gives us
\begin{eqnarray*}
\lefteqn{ c \epsilon'(\sigma) \sum_{\overline{\lambda} \in X_{N}}
q^{\left(\sigma(\overline{\lambda})-\sigma(\overline{\mu}),
\sigma(\overline{\lambda})-\sigma(\overline{\mu})  \right)}
q^{2\left(\sigma(\overline{\lambda})-\sigma(\overline{\mu})+\rho,\sigma(\overline{\mu} + \rho) \right)} }  \\
& & \hspace{15mm} = c \epsilon'(\sigma) \sum_{\overline{\lambda} \in X_{N}}
q^{(\overline{\lambda},\overline{\lambda})-2(\overline{\lambda},\overline{\mu}) +
(\overline{\mu},\overline{\mu})}
q^{2(\overline{\lambda},\overline{\mu} + \rho)-2(\overline{\mu},\overline{\mu}+\rho)}
q^{2(\rho,\sigma(\overline{\mu} + \rho))} \\
& & \hspace{15mm} = c \epsilon'(\sigma) \sum_{\overline{\lambda} \in X_{N}}
q^{(\overline{\lambda},\overline{\lambda}+2\rho)-(\overline{\mu},\overline{\mu}+2\rho)}
q^{2(\rho,\sigma(\overline{\mu} + \rho))},
\end{eqnarray*}
which allows us to rewrite (\ref{eq:Qeq}) as
$$q^{-(\overline{\mu},\overline{\mu}+2\rho)}
\sum_{\sigma \in {\cal{W}}}\epsilon'(\sigma)q^{2(\rho,\sigma(\overline{\mu}+\rho))}
=c\sum_{\overline{\lambda} \in X_{N}}
q^{(\overline{\lambda},\overline{\lambda}+2\rho)}q^{-(\overline{\mu},\overline{\mu}+2\rho)}
\sum_{\sigma \in {\cal{W}}}\epsilon'(\sigma)
q^{2(\rho,\sigma(\overline{\mu}+\rho))},$$
for each $\overline{\mu} \in p(\Lambda_{N}^{+})$.
As $Q(\overline{\mu}) \neq 0$ 
for each $\overline{\mu} \in p(\Lambda_{N}^{+})$, we easily obtain the following expression for $c^{-1}$:
\begin{equation}
\label{eq:flicky}
c^{-1} = \sum_{\overline{\lambda} \in X_{N}} q^{(\overline{\lambda},\overline{\lambda}+2\rho)},
\end{equation}
which satisfies $|c^{-1}| = \left(2N\right)^{n/2}$, as is shown presently.
We can rewrite the sum in (\ref{eq:flicky}) as
\begin{equation}
\label{eq:alisonwedding}
\sum_{\overline{\lambda} \in X_{N}}q^{(\overline{\lambda},\overline{\lambda}+2\rho)} = \prod_{k=0}^{n-1}G_{+}(N,2k+1) = 
\frac{\left[(1+i)\sqrt{N}\right]^{n}}{t^{n}}
      \left(\prod_{k=0}^{n-1}\frac{1}{q^{k(k+1)}}\right), \hspace{5mm} t=\exp{(\pi i/2N)},
\end{equation} 
where $G_{+}(N,m)$ is a Gaussian sum defined in Appendix A. 
In proving (\ref{eq:alisonwedding}), we used Lemma \ref{lem:thirdgausssum} 
and the fact that $t^{(2k+1)^{2}} = tq^{k(k+1)}$.
By inspection, the modulus of the far right hand side of (\ref{eq:alisonwedding}) is $\left(2N\right)^{n/2}$,
and we obtain the well defined expression:
\begin{equation}
\label{eq:xlambda}
x_{\overline{\lambda}} = 
\frac{q^{-(\overline{\lambda},2\rho)}t^{n}\left(\prod_{k=0}^{n-1}q^{k(k+1)}\right)}
{\left[(1+i)\sqrt{N}\right]^{n}}.
\end{equation}

Now consider the right hand side of (\ref{eq:Qeq}):
\begin{equation}
\label{eq:equanimity}
\sum_{\overline{\lambda} \in X_{N}}x_{\overline{\lambda}}
q^{(\overline{\lambda},\overline{\lambda}+2\rho)}S_{\overline{\lambda},\overline{\mu}}'.
\end{equation}
Recall that we use $N'$ to mean $N/2$.
We will rewrite (\ref{eq:equanimity}) as a sum over all the elements of $X_{N'}=X/N'X$ and
then we will use certain results in \cite{z4} to express (\ref{eq:equanimity}) as
a sum over all the elements of $\overline{\Lambda_{N}^{+}}$.  
We now consider some useful calculations.
Let $\overline{\lambda}' = \overline{\lambda}+N'\epsilon_{i}$ for some $i \in \{1,2,\ldots,n\}$, then
$$
q^{(\overline{\lambda}',\overline{\lambda}'+2\rho)} 
= q^{(\overline{\lambda} + N'\epsilon_{i},\overline{\lambda} + N' \epsilon_{i} + 2\rho)} 
= q^{(\overline{\lambda},\overline{\lambda}+2\rho)+
         N'(2 \overline{\lambda}_{i} + 2n-2i+1+N')}.
$$
As $N \equiv 2 \pmod{4}$, this equation leads to 
\begin{equation}
\label{eq:mildurabilura}
q^{(\overline{\lambda}',\overline{\lambda}' + 2\rho)}
 =q^{(\overline{\lambda},\overline{\lambda} + 2\rho)}.
\end{equation}
Thus
\begin{equation}
\label{eq:pottspointkingscross}
x_{\overline{\lambda}'} = \frac{q^{-(\overline{\lambda}+N'\epsilon_{i},2\rho)}t^{n}
\left( \prod_{k=0}^{n-1}q^{k(k+1)}\right)}{\left[(1+i)\sqrt{N}\right]^{n}} = 
\frac{-q^{-(\overline{\lambda},2\rho)}t^{n} \left( \prod_{k=0}^{n-1}q^{k(k+1)}\right)}
{\left[(1+i)\sqrt{N}\right]^{n}} = -x_{\overline{\lambda}},
\end{equation}
and furthermore
\begin{eqnarray}
S_{\overline{\lambda}',\overline{\mu}}' & = & \sum_{\sigma \in {\cal{W}}}\epsilon'(\sigma)
q^{2(\overline{\lambda}+N'\epsilon_{i} + \rho,\sigma(\overline{\mu}+\rho))} \nonumber \\
& = & \sum_{\sigma \in {\cal{W}}}\epsilon'(\sigma)
q^{2(\overline{\lambda}+\rho,\sigma(\overline{\mu}+\rho))+
N'(\epsilon_{i},\sigma(2\overline{\mu}+2\rho))} 
= -S_{\overline{\lambda},\overline{\mu}}'. \label{eq:rushcutterschaingang}
\end{eqnarray}
Using Eqs. (\ref{eq:pottspointkingscross})--(\ref{eq:rushcutterschaingang}), we obtain
\begin{equation}
\label{eq:shannongoll}
x_{\overline{\lambda}'}q^{(\overline{\lambda}',\overline{\lambda}'+2\rho)}
S_{\overline{\lambda}',\overline{\mu}}'
=x_{\overline{\lambda}}q^{(\overline{\lambda},\overline{\lambda}+2\rho)}
S_{\overline{\lambda},\overline{\mu}}'.
\end{equation}
Note that $|X_{N}| = 2^{n}|X_{N'}|$, then 
(\ref{eq:mildurabilura}) and (\ref{eq:shannongoll}) allow us to rewrite (\ref{eq:equanimity}) as
\begin{equation}
\label{eq:2ntheq}
\sum_{\overline{\lambda} \in X_{N}}
x_{\overline{\lambda}}q^{(\overline{\lambda},\overline{\lambda}+2\rho)}
S_{\overline{\lambda},\overline{\mu}}'=
2^{n}\sum_{\overline{\lambda} \in X_{N'}}x_{\overline{\lambda}}
q^{(\overline{\lambda},\overline{\lambda}+2\rho)}S_{\overline{\lambda},\overline{\mu}}'.
\end{equation}

We now aim to write the summation in the right hand side of (\ref{eq:2ntheq}) in terms of a summation 
over the elements of $\overline{\Lambda_{N}^{+}}$, and in doing this we will follow the 
general approach of \cite{z4}, 
which studied a similar problem for ordinary quantum groups at odd roots of unity.
In particular, we consider the analogous problem for $U_{q}^{(N)}(so(2n+1))$.  
Consider the affine Weyl group ${\cal{W}}_{M}^{\mathfrak{g}}$ 
of a classical Lie algebra $\mathfrak{g}$, where $M \geq 3$ is an integer, generated by
the maps 
$S_{\alpha, kM}: X \rightarrow X$, $\alpha \in \Phi_{\mathfrak{g}}^{+}$, $k \in \mathbb{Z}$, 
where $\Phi_{\mathfrak{g}}^{+}$ is the set of positive roots of $\mathfrak{g}$,
and where the action of $S_{\alpha, kM}$ on $\mu \in H^{*}$ is defined by
$$S_{\alpha, kM}: \mu \mapsto \sigma_{\alpha}(\mu+\rho)-\rho+kM\alpha.$$ 
We derive the following from \cite[Sec. 6.2]{ja}.
\begin{remark}
Set $M \geq 3$ to be odd, then
the affine Weyl group ${\cal{W}}_{M}^{\mathfrak{g}}$ of a classical Lie algebra
$\mathfrak{g}$ acts on the chambers, that is the open
connected components of
$$X-\bigcup_{\alpha \in \Phi^{+}_{\mathfrak{g}}}\bigcup_{n \in \mathbb{Z}}
\left\{\mu \in X  \Big| \ \frac{2(\mu+\rho,\alpha)}{(\alpha, \alpha)}=nM \right\},$$
transitively, with a fundamental domain 
\begin{equation}
\label{eq:dreams}
\left\{\mu \in X \Big| \ 0 \leq \frac{2(\mu+\rho,\alpha)}{(\alpha,\alpha)}
\leq M, \ \forall \alpha \in \Phi^{+}_{\mathfrak{g}}  \right\}.
\end{equation}
\end{remark}

\begin{proposition}
Set $N \equiv 2 \pmod{4}$, $N \geq 6$.  
The truncated Weyl alcove $\overline{\Lambda_{N}^{+}} \subset X$ given in Definition \ref{def:opencat}:
\begin{equation}
\label{eq:stevebracks(a)}
\overline{\Lambda_{N}^{+}} = \left\{ \lambda \in X \Big| \ 
0 \leq \frac{2(\lambda + \rho,\alpha)}{(\alpha,\alpha)} \leq N/2, \ \forall \alpha \in 
\overline{\Phi}_{0}^{+} \cup \Phi_{1}^{+}   \right\},
\end{equation}
 is identical to
the fundamental domain of $X$ under the action of the affine Weyl group ${\cal{W}}^{so(2n+1)}_{N/2}$ stated in
(\ref{eq:dreams}).
\end{proposition}
\begin{proof}
The set of positive roots $\Phi^{+}$ of $so(2n+1)$ is identical to the subset 
$\overline{\Phi}_{0}^{+} \cup \Phi_{1}^{+}$ of positive roots of $osp(1|2n)$, 
and the expression for $2\rho$ in $so(2n+1)$ is given by
$$2\rho = \sum_{i<j} \big[ (\epsilon_{i} - \epsilon_{j}) + (\epsilon_{i} + \epsilon_{j})\big] 
           + \sum_{k=1}^{n} \epsilon_{k} 
        = \sum_{i=1}^{n}(2n-2i+1)\epsilon_{i},$$ 
which is identical to the expression for $2\rho$ in $osp(1|2n)$.  The result follows.
\end{proof}

Let $\beta \in H^{*}$ be arbitrary and let $\{s_{\alpha} | \ \alpha \in \Phi^{+} \}$ be
a subset of elements of the Weyl group of $osp(1|2n)$.  
Now
$s_{\epsilon_{i}}(\beta) = s_{2\epsilon_{i}}(\beta)$ for each  $i=1,\ldots,n$,
and we identify  $s_{2\epsilon_{i}}$ with $s_{\epsilon_{i}}$, thus every element of
${\cal{W}}$ can be expressed as some ordered product of the elements of 
$\big\{ s_{\alpha} \in {\cal{W}} | \ \alpha \in \overline{\Phi}^{+}_{0} \cup \Phi^{+}_{1} \big\}$.  As
$\Phi^{+}_{so(2n+1)} = \overline{\Phi}^{+}_{0} \cup \Phi^{+}_{1}$, 
the Weyl groups of $osp(1|2n)$ and $so(2n+1)$ are  identical.

We define the action of $S_{\alpha,kM} \in {\cal{W}}_{M}^{\mathfrak{g}}$ on an element of $X_{M}$ by
$$S_{\alpha,kM}(\mu+MX) = S_{\alpha,kM}(\mu) + MX.$$
This coincides with the action of the Weyl group ${\cal{W}}^{\mathfrak{g}}$ on $X_{M}$ defined by
$$\sigma(\mu + MX) = \sigma(\mu) + MX, \hspace{10mm} \forall \sigma \in {\cal{W}}^{\mathfrak{g}}, \ \mu \in X,$$
and we can deduce from this that the image of $\overline{{\Lambda}_{N}^{+}}$ 
under the canonical projection $p_{N/2}: X \rightarrow X_{N/2}$ furnishes a fundamental domain for $X$ 
under the action of the Weyl group ${\cal{W}}$ \cite{z4}.

There is the following important result \cite{z4}:
for any $\lambda,\mu \in p_{N/2}(\Lambda_{N}^{+}) \subset  X_{N/2}$ and any $\sigma,w \in {\cal{W}}^{so(2n+1)}$,
\begin{equation}
\label{chapter2:gigiwong(2)}
\sigma(\lambda + \rho) - \rho = w(\mu + \rho) - \rho  \hspace{5mm} \mbox{iff} \hspace{5mm} 
\lambda = \mu \ {\mathrm{and }} \ \sigma = w.
\end{equation}
This result also holds for any $\sigma,w \in {\cal{W}}^{osp(1|2n)}$ as ${\cal{W}}^{osp(1|2n)} = {\cal{W}}^{so(2n+1)}$.

In order to rewrite the right hand side of (\ref{eq:2ntheq}) in terms of a
summation over the elements of $\overline{\Lambda_{N}^{+}}$ we will show that 
$S_{\lambda,\mu}'=0$ if $\lambda \in \overline{\Lambda_{N}^{+}} \backslash \Lambda_{N}^{+}$ 
or if $\mu \in \overline{\Lambda_{N}^{+}} \backslash \Lambda_{N}^{+}$.  
The corresponding result is easily proved for $U_{q}^{(N/2)}(so(2n+1))$ \cite{z4} but for $U_{q}^{(N)}(osp(1|2n))$ the
proof is more intricate, principally due to the different properties of $\epsilon(\sigma)$ and $\epsilon'(\sigma)$ where
$\sigma$ is an element of  ${\cal{W}}^{osp(1|2n)}$.
Recall that $N \equiv 2 \pmod{4}$ and that $\Lambda_{N}^{+}$ is defined by
$$\Lambda_{N}^{+} 
= \left\{ \mu \in X \Big| \ 0 <
\frac{2(\mu+\rho,\alpha)}{(\alpha,\alpha)} < N', \ \forall \alpha \in
\overline{\Phi}_{0}^{+} \cup \Phi_{1}^{+} \right\},$$
which is similar to the definition of $\overline{\Lambda_{N}^{+}}$ in (\ref{eq:stevebracks(a)}).

The following two important properties of $S_{\lambda,\mu}'$ are easily proved:
for each $\lambda, \mu \in X$,
\begin{itemize}
\item[(i)] $S_{\lambda,\mu}' = S_{\mu,\lambda}'$,
\item[(ii)] $S_{s_{\alpha}(\lambda+\rho)-\rho,\mu}' =
\epsilon'(s_{\alpha})S_{\lambda,\mu}'$ for any $s_{\alpha} \in {\cal{W}}$.
\end{itemize}

\begin{lemma}
\label{lem:guttersnap}
If $\lambda \in \overline{\Lambda_{N}^{+}} \backslash \Lambda_{N}^{+}$ or 
$\mu \in \overline{\Lambda_{N}^{+}} \backslash \Lambda_{N}^{+}$, $S_{\lambda,\mu}' = 0$.
\end{lemma}
\begin{proof}
 Define $hp_{\alpha}=\left\{ \mu \in X | \ (\mu,\alpha)=0 \right\}$ for each 
 $\alpha \in \overline{\Phi}_{0}^{+} \cup \Phi_{1}^{+}$;
$hp_{\alpha}$ is the subset of $X$ invariant under the action of $s_{\alpha} \in {\cal{W}}$.  
For an element $\lambda \in \overline{\Lambda_{N}^{+}} \backslash \Lambda_{N}^{+}$, the definitions of
$\overline{\Lambda_{N}^{+}}$ and $\Lambda_{N}^{+}$ imply that 
 $\frac{2(\lambda+\rho,\alpha)}{(\alpha,\alpha)} = kN/2$ for some 
 $k \in \mathbb{Z}_{+}$ and some $\alpha \in \overline{\Phi}_{0}^{+} \cup \Phi_{1}^{+}$.  
For $\alpha = \epsilon_{i} \pm \epsilon_{j}$, where $1 \leq i < j \leq n$, 
we have $(\lambda + \rho,\alpha) = kN/2$,
thus $(\lambda + \rho - kN\alpha/4,\alpha) = 0$, and 
$\lambda + \rho - kN\alpha/4 \in hp_{\alpha}$. 
Consequently, $s_{\alpha}(\lambda+\rho-kN\alpha/4) = \lambda+\rho-kN\alpha/4$, and
$s_{\alpha}(\lambda + \rho) - \rho + kN\alpha/2 = \lambda$, and it follows that 
$$S_{\lambda,\mu}' = S_{s_{\alpha}(\lambda +\rho)-\rho + kN'\alpha,\mu}'
=S_{s_{\alpha}(\lambda + \rho)-\rho,\mu}' = \epsilon'(s_{\alpha})S_{\lambda,\mu}',$$
which vanishes identically as $\epsilon'(s_{\alpha}) = -1$ 
for each $\alpha = \epsilon_{i} \pm \epsilon_{j} \in \overline{\Phi}^{+}_{0}$.

Now let $\lambda \in \overline{\Lambda_{N}^{+}} \backslash \Lambda_{N}^{+}$
and let $\alpha = \epsilon_{i}$, where $1 \leq i \leq n$, then
$\frac{2(\lambda+\rho,\alpha)}{(\alpha,\alpha)} = kN/2$ for some $k\in \mathbb{Z}$.  
Consequently, we have $2(\lambda + \rho,\alpha)=kN/2$ which implies that 
$(\lambda + \rho-kN\alpha/4,\alpha)=0$, and also that $\lambda + \rho - kN\alpha/4 \in hp_{\alpha}$. 
It follows that 
$s_{\alpha}(\lambda + \rho - kN\alpha/4)=\lambda + \rho - kN\alpha/4$, and that 
$$s_{\alpha}(\lambda + \rho)-\rho + kN\alpha/2 = \lambda.$$  
Consequently, we have
\begin{equation}
\label{chapter4:gigiwong(1)}
S'_{\lambda,\mu}=S'_{s_{\alpha}(\lambda+\rho)-\rho+kN'\alpha,\mu} = 
\epsilon'(s_{\alpha})S'_{\lambda+kN'\alpha,\mu} = (-1)^{k}\epsilon'(s_{\alpha})S'_{\lambda,\mu}
=(-1)^{k} S'_{\lambda,\mu},
\end{equation}
where we have used the result $\epsilon'(s_{\alpha})=1$ as $\alpha \in \Phi_{1}^{+}$,
and we have also used the following calculations:
\begin{eqnarray*}
S'_{\lambda + kN'\alpha,\mu} & = & \sum_{\sigma \in {\cal{W}}}
\epsilon'(\sigma)q^{2(\lambda+kN'\alpha + \rho,\sigma(\mu+\rho))} \\
 & = & \sum_{\sigma \in {\cal{W}}}\epsilon'(\sigma)q^{2(\lambda +
 \rho,\sigma(\mu+\rho))}q^{(kN'\alpha,\sigma(2\mu+2\rho))} 
 = \left\{ \begin{array}{ll}
 S'_{\lambda,\mu}, & {\mbox{if $k$ is even, }}   \\
 -S'_{\lambda,\mu}, & {\mbox{if $k$ is odd. }}
 \end{array} \right.
\end{eqnarray*}
Here the last equality arises from the following calculation:
$$q^{(kN'\alpha,\sigma(2\mu+2\rho))} = \left\{ \begin{array}{ll}
+1 & {\mbox{if $k$ is even, }}   \\
-1 & {\mbox{if $k$ is odd.}}
\end{array}
\right.$$  
If $k$ is odd, $S'_{\lambda,\mu}$ vanishes identically by (\ref{chapter4:gigiwong(1)}), 
and this completes the proof of the assertion
that $S'_{\lambda,\mu}=0$ if $\lambda \in \overline{\Lambda_{N}^{+}} \backslash \Lambda_{N}^{+}$.
To show that $k$ is indeed odd we note that the equation
$(\lambda+\rho,\alpha) = kN/4$ implies that $\lambda_{i}+n-i+1/2=kN/4$, 
and thus we have $k \notin 2 \mathbb{Z}$ as
$\lambda \in X$. It follows then that
$s_{\alpha}(\lambda+\rho) - \rho + kN\alpha/2 = \lambda$ for some odd $k$, and thus
$$S'_{\lambda,\mu} = -S'_{\lambda,\mu} = 0.$$

To complete the proof, we note that $S_{\lambda,\mu}' = S_{\mu,\lambda}'$, and thus
$S'_{\lambda,\mu}=0$ if $\mu \in \overline{\Lambda_{N}^{+}} \backslash \Lambda_{N}^{+}$.
\end{proof}

\begin{corollary}
\label{cor:deadstonecrow}
If
$\overline{\lambda}$ or $\overline{\mu}$ belongs to
 $p_{N/2}\big(\overline{\Lambda_{N}^{+}} \backslash \Lambda_{N}^{+}\big) \subset X_{N/2}$, then
$S'_{\overline{\lambda},\overline{\mu}} = 0$.
\end{corollary}

\begin{lemma}
The set  $\{ d_{\lambda} \in \mathbb{C} | \ \lambda \in \Lambda_{N}^{+} \}$ of constants
defined by
\begin{equation}
\label{eq:arnieee}
d_{\lambda} = d_{0}sdim_{q}(V_{\lambda}), \hspace{10mm}  d_{0} = \Omega Q(0),
\end{equation}
with
\begin{equation}
\label{eq:davideasdownideal(a)} 
\Omega = \frac{2^{n}t^{n}q^{n^{3}-n/2}}{\left[(1+i)\sqrt{N}   \right]^{n}},
\end{equation}
\begin{equation}
\label{eq:davideasdownideal(b)} 
Q(0) = \prod_{\alpha \in \overline{\Phi}^{+}_{0}} \big( q^{(\rho,\alpha)}-q^{-(\rho,\alpha)} \big)
       \prod_{\beta \in \Phi^{+}_{1}} \big( q^{(\rho,\beta)}+q^{-(\rho,\beta)} \big),
\end{equation}
satisfies the relations (\ref{eq:wowwee}).
\end{lemma}
\begin{proof}

Using (\ref{chapter2:gigiwong(2)}) and
Corollary \ref{cor:deadstonecrow}, we can rewrite the right hand side of (\ref{eq:2ntheq}) as
\begin{eqnarray*}
2^{n} \sum_{\overline{\lambda} \in X_{N'}} x_{\overline{\lambda}}
q^{(\overline{\lambda},\overline{\lambda}+2\rho)}S'_{\overline{\lambda},\overline{\mu}}
& = & 2^{n} \sum_{\lambda \in \Lambda_{N}^{+}}\sum_{\sigma \in {\cal{W}}}
x_{\sigma(\lambda + \rho)-\rho}q^{\left(\sigma(\lambda +
\rho)-\rho,\sigma(\lambda+\rho)+\rho  
\right)}S'_{\sigma(\lambda +\rho)-\rho,\overline{\mu}}\\
& = & 2^{n} \sum_{\lambda \in \Lambda_{N}^{+}}\sum_{\sigma \in {\cal{W}}}
x_{\sigma(\lambda +
\rho)-\rho}q^{(\lambda,\lambda+2\rho)}\epsilon'(\sigma)S'_{\lambda,\overline{\mu}},
\end{eqnarray*}
and we consequently obtain from Eqs. (\ref{eq:Qeq}) and (\ref{eq:superQ})
the following equation for $d_{\nu}'$:
\begin{equation}
\label{eq:tomthedeadguineapig}
\sum_{\nu \in \Lambda_{N}^{+}}d_{\nu}'q^{(\nu,\nu+2\rho)}S'_{\nu,\mu}=
2^{n} \sum_{\lambda \in \Lambda_{N}^{+}}\sum_{\sigma \in {\cal{W}}}
\epsilon'(\sigma)x_{\sigma(\lambda+\rho)-\rho}q^{(\lambda,\lambda+2\rho)} S'_{\lambda,\mu}, 
\hspace{5mm} \forall \mu \in \Lambda_{N}^{+}.
\end{equation}
Eq. (\ref{eq:tomthedeadguineapig}) is obviously satisfied by
$$d_{\lambda}' =  2^{n}\sum_{\sigma \in
{\cal{W}}}\epsilon'(\sigma)x_{\sigma(\lambda + \rho)-\rho},$$
and we now evaluate $d_{\lambda}'$:
\begin{eqnarray*}
d_{\lambda}' & = & 2^{n} \sum_{\sigma \in {\cal{W}}}
\epsilon'(\sigma)q^{-\left(\sigma(\lambda+\rho)-\rho,2\rho  \right)}/
\sum_{\overline{\nu} \in X_{N}}q^{(\overline{\nu},\overline{\nu}+2\rho)}  \\
& = & 2^{n} q^{(2\rho,\rho)}\epsilon'(w_{0})\sum_{\sigma \in {\cal{W}}}
\epsilon'(\sigma)q^{-(w_{0}\sigma(\lambda+\rho),2\rho)}/
\sum_{\overline{\nu} \in X_{N}}q^{(\overline{\nu},\overline{\nu}+2\rho)}  \\
& = & 2^{n} q^{(2\rho,\rho)}\sum_{\sigma \in {\cal{W}}}\epsilon'(\sigma)
q^{(\sigma(\lambda + \rho),2\rho)}/\sum_{\overline{\nu} \in X_{N}}
q^{(\overline{\nu},\overline{\nu}+2\rho)},   
\end{eqnarray*}
where $w_{0}= s_{\epsilon_{1}}s_{\epsilon_{2}}\cdots s_{\epsilon_{n}}$ 
is the longest element of ${\cal{W}}$;  note that $\epsilon'(w_{0}) = 1$.  
Therefore,
$$d_{\lambda}' =
(-1)^{[\lambda]}\frac{2^{n}q^{(2\rho,\rho)}sdim_{q}(V_{\lambda})Q(0)}
{\sum_{\overline{\nu} \in X_{N}}q^{(\overline{\nu},\overline{\nu}+2\rho)}},$$ 
$$d_{\lambda} = \frac{2^{n}q^{(2\rho,\rho)}sdim_{q}(V_{\lambda})Q(0)}
{\sum_{\overline{\nu} \in X_{N}}q^{(\overline{\nu},\overline{\nu}+2\rho)}}.$$
(Recall that $d_{\lambda} = (-1)^{[\lambda]}d_{\lambda}'$.)
Evaluating the Gaussian sums gives
\begin{equation}
\label{eq:dss}
d_{\lambda} = \frac{2^{n}q^{(2\rho,\rho)} sdim_{q}(V_{\lambda})Q(0)t^{n}
\left( \prod_{k=0}^{n-1}q^{k(k+1)} \right)}{\left[ (1+i)\sqrt{N}  \right]^{n}},
\end{equation}
which is non-zero for each $\lambda \in \Lambda_{N}^{+}$.  Writing 
$$d_{\lambda} = d_{0}sdim_{q}(V_{\lambda}), \hspace{15mm}  d_{0} = \Omega Q(0), \hspace{15mm} \lambda \in \Lambda_{N}^{+},$$
we then have 
\begin{eqnarray*}
\Omega & = & \frac{2^{n}q^{(2\rho,\rho)}t^{n}
\left(\prod_{k=0}^{n-1}q^{k(k+1)}\right)}{\left[(1+i)\sqrt{N}   \right]^{n}} \\
& = & \frac{2^{n} q^{(4n^{3}-n)/6} t^{n} q^{(n^{3}-n)/3} }
{\left[(1+i)\sqrt{N}   \right]^{n}} \\
& = & \frac{2^{n}t^{n}q^{n^{3}-n/2}}{\left[(1+i)\sqrt{N}   \right]^{n}}.
\end{eqnarray*}
In this calculation 
we used the result
$(2\rho,2\rho)=(4n^{3}-n)/3$.
\end{proof}

The constants $d_{\lambda}$ are proportional to 
$sdim_{q}(V_{\lambda})$ with constant of proportionality $d_{0} \neq 0$ 
for each $\lambda \in \Lambda_{N}^{+}$.
A similar phenomenon occurs for all other Reshetikhin-Turaev $3$-manifold invariants constructed from 
quotients of quantum algebras and quantum superalgebras at even and odd roots of unity 
\cite{rt,tw,z2,z3,z4,z6}.
\end{subsection}

\begin{subsection}{Proof of Axiom (VI)}

We now show that $z$ is as claimed.
\begin{lemma}
\label{eq:emotion}
Set
$\displaystyle{z=\sum(L) = \sum_{\lambda \in
\Lambda_{N}^{+}}d_{\lambda}q^{-(\lambda,\lambda+2\rho)}sdim_{q}(V_{\lambda})}$.
Then 
\begin{equation}
\label{eq:ronjames(a)}
z = (-i)^{n} q^{2n^{3}-n} t^{2n},
\end{equation}
which clearly satisfies $|z|=1$.
\end{lemma}
\begin{proof}
We calculate as follows:
\begin{eqnarray}
z & = & \sum_{\lambda \in
        \Lambda_{N}^{+}}d_{\lambda}q^{-(\lambda,\lambda+2\rho)}sdim_{q}(V_{\lambda}) \nonumber \\
  & = & \Omega Q(0) \sum_{\lambda \in \Lambda_{N}^{+}}
        q^{-(\lambda,\lambda+2\rho)}   \big(sdim_{q}(V_{\lambda})\big)^{2}  \nonumber \\
  & = & \frac{\Omega}{Q(0)} \sum_{\lambda \in \Lambda_{N}^{+}}
        q^{-(\lambda,\lambda+2\rho)}   \big(Q(\lambda)\big)^{2}, \label{eq:squarewindow}
\end{eqnarray}
where we have used (\ref{eq:arnieee}) and the relation 
$\big( sdim_{q}(V_{\lambda}) \big)^{2} = \big(S_{\lambda,0}/Q(0)\big)^{2} = \big(Q(\lambda)/Q(0)\big)^{2}$.

Let us examine $q^{-(\lambda, \lambda+2\rho)}$ and 
$\big(Q(\lambda)\big)^{2}$ under the action of the map 
$\lambda \mapsto \sigma(\lambda + \rho) -\rho$, where $\sigma \in {\cal{W}}$.
It is not difficult to show that for each $\sigma \in {\cal{W}}$ we have
$$q^{-(\sigma(\lambda + \rho)-\rho,\sigma(\lambda+\rho)-\rho+2\rho)} = q^{-(\lambda,\lambda+2\rho)},
 \hspace{5mm} \mbox{ and } \hspace{5mm} \big(Q(\sigma(\lambda+\rho)-\rho)\big)^{2} = \big(Q(\lambda)\big)^{2}.$$  
Now Eq. (\ref{eq:mildurabilura}) states that
$$q^{-(\lambda + N'\epsilon_{i},\lambda + N'\epsilon_{i} + 2\rho)} = q^{-(\lambda,\lambda + 2\rho)}.$$ 
Furthermore, by using the fact that 
$q^{(N' \epsilon_{i}, \sigma(2\rho))} = -1$ for each $i \in \{1, 2, \ldots, n\}$, we obtain
$$\big(Q(\lambda + N'\epsilon_{i})\big)^{2} = \big(Q(\lambda)\big)^{2}.$$ 

Now $p_{N/2}\big(\overline{\Lambda_{N}^{+}}\big)$ is a fundamental domain for $X$ 
under the action of the affine Weyl group ${\cal{W}}_{N'}$ and $Q(\lambda) = 0$ if 
$\lambda \in \overline{\Lambda_{N}^{+}} \backslash \Lambda_{N}^{+}$.  
The calculations in the previous paragraph imply that 
we can write the sum in (\ref{eq:squarewindow}) as a sum over the elements of $X_{N'}$, which
considerably simplifies the calculations:
\begin{equation}
\label{eq:contributionlizard}
\sum_{\lambda \in \Lambda_{N}^{+}}
  q^{-(\lambda,\lambda+2\rho)}   \big(Q(\lambda)\big)^{2} = 
  \frac{1}{|{\cal{W}|}}\sum_{\overline{\lambda} \in X_{N'}}
  q^{-(\overline{\lambda},\overline{\lambda}+2\rho)}   \big(Q(\overline{\lambda})\big)^{2},
\end{equation}
where we note that $|{\cal{W}}|=2^{n}n!$ \cite[p. 66, Table 1]{h}. 
By using  Eqs. (\ref{eq:mildurabilura}) and (\ref{eq:rushcutterschaingang})
we can further rewrite the sum on the right hand side of (\ref{eq:contributionlizard}):
\begin{equation}
\label{eq:nehru}
  \frac{1}{|{\cal{W}|}}\sum_{\overline{\lambda} \in X_{N'}}
  q^{-(\overline{\lambda},\overline{\lambda}+2\rho)}   \big(Q(\overline{\lambda})\big)^{2} 
 =  \frac{1}{2^{n}|{\cal{W}|}}\sum_{\overline{\lambda} \in X_{N}}
  q^{-(\overline{\lambda},\overline{\lambda}+2\rho)}   \big(Q(\overline{\lambda})\big)^{2},
\end{equation}
which we will now evaluate. We firstly rewrite the right hand side of (\ref{eq:nehru}):
\begin{eqnarray}
\lefteqn{\frac{1}{2^{n}|{\cal{W}|}}\sum_{\overline{\lambda} \in X_{N}}
  q^{-(\overline{\lambda},\overline{\lambda}+2\rho)}   \big(Q(\overline{\lambda})\big)^{2} } \nonumber \\
  & & = \frac{1}{2^{n}|{\cal{W}|}}\sum_{\overline{\lambda} \in X_{N}}
  q^{-(\overline{\lambda},\overline{\lambda}+2\rho)} \sum_{\sigma,w \in
  {\cal{W}}}\epsilon'(\sigma)\epsilon'(w)
q^{2(\overline{\lambda} + \rho,\sigma(\rho)+w(\rho))} \nonumber  \\
  & & = \frac{1}{2^{n}|{\cal{W}|}}\sum_{\sigma,w \in {\cal{W}}}\epsilon'(\sigma)\epsilon'(w)
  q^{2(\rho,\sigma(\rho)+w(\rho))}
  \sum_{\overline{\lambda} \in X_{N}}q^{-(\overline{\lambda},\overline{\lambda}+2\rho)}
q^{(2\overline{\lambda},\sigma(\rho) + w(\rho))}. \label{eq:eguy}
\end{eqnarray}

We wish to disentangle the summation indices in the exponents of $q$ in (\ref{eq:eguy}) 
so that we can factorise the summations into a sum over the elements of ${\cal{W}} \times {\cal{W}}$ 
and a sum over the elements of $X_{N}$, both of which we can calculate relatively easily.
To do this, for each pair $(\sigma, w) \in {\cal{W}} \times {\cal{W}}$ in (\ref{eq:eguy}) and each 
$\overline{\lambda} \in X_{N}$ we apply the following map:
\begin{equation}
\label{eq:dalekton}
\overline{\lambda} \mapsto \overline{\lambda} + \sigma(\rho) + w(\rho) \in X_{N}.
\end{equation}
Now 
$$p\big(\lambda + \sigma(\rho) + w(\rho)\big) 
= p\big(\mu + \sigma(\rho) + w(\rho)\big) \hspace{5mm} \mbox{ iff } \hspace{5mm} p(\lambda) = p(\mu), 
\hspace{5mm} \forall \lambda,\mu \in X,$$ 
thus the summation in (\ref{eq:eguy}) is unchanged under the mapping (\ref{eq:dalekton}). 
Applying this map to the right hand side of (\ref{eq:eguy})
and equating the result with the left hand side of (\ref{eq:contributionlizard}) gives
\begin{equation}
\label{eq:bronskiboi}
\sum_{\lambda \in \Lambda_{N}^{+}}
  q^{-(\lambda,\lambda+2\rho)}   (Q(\lambda))^{2}=
\frac{q^{(2\rho,\rho)}}{2^{n}|{\cal{W}|}}\sum_{\sigma,w \in {\cal{W}}}
\epsilon'(\sigma)\epsilon'(w)q^{2(\sigma(\rho),w(\rho))}
 \sum_{\overline{\lambda} \in X_{N}}q^{-(\overline{\lambda}+2\rho,\overline{\lambda})}.
\end{equation}
It is not difficult to evaluate the right hand side of (\ref{eq:bronskiboi}): note that 
$$\sum_{\sigma,w \in {\cal{W}}}
\epsilon'(\sigma)\epsilon'(w)q^{2(\sigma(\rho),w(\rho))} = |{\cal{W}}|Q(0).$$
This can be seen in the following way: firstly fix $\sigma_{1} \in {\cal{W}}$, then  
$$\sum_{w \in {\cal{W}}}\epsilon'(w)q^{2(\sigma_{1}(\rho),w(\rho))} = 
\sum_{w \in {\cal{W}}}\epsilon'(\sigma_{1})\epsilon'(w)q^{2(\rho,w(\rho))} = 
\epsilon'(\sigma_{1})Q(0),$$
thus
$$\sum_{\sigma \in {\cal{W}}}\epsilon'(\sigma)\sum_{w\in{\cal{W}}}
\epsilon'(w)q^{2(\sigma(\rho),w(\rho))} = 
\sum_{\sigma, w \in {\cal{W}}} \epsilon'(w) q^{2(\rho,w(\rho))}
= \sum_{\sigma \in {\cal{W}}}Q(0) = |{\cal{W}}|Q(0).$$
Additionally, Lemma \ref{lem:bugger} implies that
$$\sum_{\overline{\lambda} \in X_{N}}q^{-(\overline{\lambda}+2\rho,\overline{\lambda})} = 
\prod_{k=0}^{n-1}G_{-}(N,2k+1)=
t^{n}\left[(1-i)\sqrt{N}\right]^{n}\left(\prod_{k=0}^{n-1}q^{k(k+1)}\right),$$ 
where $t=\exp{(\pi i/2N)}$.  
By combining these results we obtain
\begin{eqnarray*}
z & = & \frac{\Omega}{Q(0)}\sum_{\lambda \in \Lambda_{N}^{+}}
  q^{-(\lambda,\lambda+2\rho)}   \big(Q(\lambda)\big)^{2}  \\
  & = &
  \frac{\Omega|{\cal{W}}|Q(0)q^{(2\rho,\rho)} t^{n} \left[(1-i)\sqrt{N}\right]^{n}
  \left(\prod_{k=0}^{n-1}q^{k(k+1)}\right)} {2^{n}|{\cal{W}}|Q(0)}  \\
  & = & \frac{\Omega q^{(2\rho,\rho)} t^{n} 
  \left[(1-i)\sqrt{N}\right]^{n}\left(\prod_{k=0}^{n-1}q^{k(k+1)}\right)}{2^{n}}   \\
  & = & \frac{ q^{(2\rho,2\rho)} t^{2n}
  \left[(1-i)\sqrt{N}\right]^{n}\left(\prod_{k=0}^{n-1}q^{k(k+1)}\right)^{2}} 
  {\left[(1+i)\sqrt{N}\right]^{n}}  \\
  & = & (-i)^{n} q^{2n^{3}-n} t^{2n},
\end{eqnarray*}
and $|z|=1$.
\end{proof}
This completes the proof of Theorem \ref{theorem:bigbigbig}.
\end{subsection}

\end{section}

\begin{section}{Comparing the invariants from $U_{q}^{(N)}(osp(1|2n))$ and $U_{q}^{(N/2)}(so(2n+1))$}
\markright{\text{Comparing the invariants from $U_{q}^{(N)}(osp(1|2n))$ and $U_{q}^{(N/2)}(so(2n+1))$}}
\label{chapter4sectionlabel(comparing_the_)}

Our construction of 3-manifold invariants  from $U_{q}^{(N)}(osp(1|2n))$
immediately gives rise to three important questions:
\begin{itemize}
\item[1.]  Are our topological invariants of $3$-manifolds obtained 
           from $U_{q}^{(N)}(osp(1|2n))$ {\emph{new}} invariants?
\item[2.]  Are the $3$-manifold invariants obtained 
           from $U_{q}^{(N)}(osp(1|2n))$ {\emph{complete}} invariants, 
	   that is, do they distinguish non-homeomorphic 3-manifolds
	   as well as telling us when two 3-manifolds are homeomorphic?  	   
\item[3.]  If the invariants are not complete, are they better than other invariants in 
           distinguishing non-homeomorphic 3-manifolds?
\end{itemize}

These are difficult questions to answer.  
The first requires a comparison between our invariants and all other existing
invariants, 
a positive answer to the second would solve the classification problem for closed, connected,
orientable $3$-manifolds 
and the third requires a theoretical investigation of the properties of the various invariants, 
 or the calculation of various invariants for numerous $3$-manifolds 
  and directly comparing their 
performance in distinguishing non-homeomorphic $3$-manifolds with the performance of our invariant. 
We do not know of  any theorems that would allow such a theoretical investigation, 
and the numerical work needed to compare the performance of the
various invariants is 
itself a  non-trivial exercise, similar to that done to compare how well polynomial link invariants 
 distinguish links that are not ambient isotopic (eg see \cite[Ch. 7]{ddw}).

A fact touching on question 2 is that the 
Reshetikhin-Turaev method for constructing 3-manifold invariants 
does not ensure {\emph{completeness}}: 
it does not necessarily distinguish non-homeomorphic 3-manifolds.
For example, 
the topological invariants derived from $U_{q}^{(N)}(sl_{2})$, where $N \geq 4$ satisfies
$N \equiv 0 \pmod{4}$,
do not distinguish all non-homeomorphic $3$-manifolds \cite{kb,l93,kl}.

Given the difficulty of answering these questions, 
we consider a more tractable problem.  
We will compare our invariants with the invariants derived from 
one quantum group at odd roots of unity, and ask the following question:
are the invariants from 
$U_{q}^{(N)}(osp(1|2n))$ the same as those from $U_{q}^{(N/2)}(so(2n+1))$
when $N \geq 6$ satisfies $N \equiv 2 \pmod{4}$?
By {\emph{the same}}, we mean that given a closed, connected, orientable 3-manifold $M_{L}$, we have
\begin{equation}
\label{eq:fiveone}
{\cal{F}}(M_{L})_{U_{q}^{(N)}(osp(1|2n))} = {\cal{F}}(M_{L})_{U_{q}^{(N/2)}(so(2n+1))}.
\end{equation}
It is interesting to ask this question as
the sets of integral weights in the truncated Weyl chambers $\Lambda_{N}^{+}$ 
of $U_{q}^{(N)}(osp(1|2n))$ and $U_{q}^{(N/2)}(so(2n+1))$ are the same, thus
in calculating the topological invariant for any given $3$-manifold $M_{L}$
we sum over the same module labels.
The reader may ask the obvious question
why we are not comparing the invariants from quantum $osp(1|2n)$ and 
quantum $so(2n+1)$ at the same roots of unity.  
The reason is that the sets of integral weights in the truncated Weyl alcoves of quantum $osp(1|2n)$ and 
quantum $so(2n+1)$ are different when $N \equiv 2 \pmod{4}$.

Our topological invariant has an $S^{3}$ normalisation:
$${\cal{F}}(S^{3})_{U_{q}^{(N)}(osp(1|2n))}=1={\cal{F}}(S^{3})_{U_{q}^{(N/2)}(so(2n+1))},$$
and therefore we choose a 3-manifold other than $S^{3}$
on which to compare the two families of invariants.
For calculational ease we will determine the invariants associated with $S^{2} \times S^{1}$, 
and we recall that we can obtain $S^{2} \times S^{1}$ by performing surgery on an oriented
unlink $L \subset S^{3}$ with {\emph{zero}} framing.   
We now calculate ${\cal{F}}(S^{2} \times S^{1})$:
the linking matrix of $L$ is $A_{L} = (0)$, thus $\sigma(A_{L}) = 1$ and
$$
{\cal{F}}(S^{2} \times S^{1}) = z^{-1} \sum(L) 
  = z^{-1} \sum_{\lambda \in \Lambda_{N}^{+}} d_{\lambda} sdim_{q}(V_{\lambda}),
$$
where in calculating ${\cal{F}}(S^{3})_{U_{q}^{(N/2)}(so(2n+1))}$ 
we take the quantum dimension of the
irreducible $U_{q}^{(N/2)}(so(2n+1))$-module $V_{\lambda}$ with highest weight $\lambda$ 
instead of the quantum superdimension.

We firstly calculate ${\cal{F}}(S^{2} \times S^{1})_{U_{q}^{(N)}(osp(1|2n))}$.
Let $N \geq 6$ satisfy $N \equiv 2 \pmod{4}$ and let $q = \exp{(2 \pi i/N)}$, then 
\begin{eqnarray}
\sum(L) & = & \sum_{\lambda \in \Lambda_{N}^{+}} d_{\lambda} sdim_{q}(V_{\lambda}) 
  =  \Omega Q(0) \sum_{\lambda \in \Lambda_{N}^{+}} \big(sdim_{q}(V_{\lambda})\big)^{2} \nonumber \\
& = & \frac{\Omega}{Q(0)} \sum_{\lambda \in \Lambda_{N}^{+}} \big(Q(\lambda)\big)^{2} 
  =  \frac{\Omega}{Q(0)|{\cal{W}}|} \sum_{\overline{\lambda} \in X_{N'}}
\left(Q(\overline{\lambda})\right)^{2} 
  =  \frac{\Omega}{2^{n} Q(0)|{\cal{W}}|} 
\sum_{\overline{\lambda} \in X_{N}} \left(Q(\overline{\lambda})\right)^{2}. \nonumber \\
& & \label{eq:fivetwo}
\end{eqnarray}
Using the expressions for $\sum_{\overline{\lambda} \in X_{N}} \left(Q(\overline{\lambda})\right)^{2}$ 
(Lemma \ref{lem:ballson}) and for
$\Omega$ and $z$ (Eqs. (\ref{eq:davideasdownideal(a)}) and (\ref{eq:ronjames(a)})):
$$
\Omega 
= \frac{2^{n} q^{n^{3}-n/2}t^{n}}{\left[(1+i)\sqrt{N}\right]^{n}},
\hspace{10mm} z = (-i)^{n} q^{2n^{3}-n}t^{2n},
\hspace{15mm} t=\exp{(\pi i/2N)},$$  
we have 
\begin{eqnarray*}
{\cal{F}}(S^{2} \times S^{1})_{U_{q}^{(N)}(osp(1|2n))} 
  & = &  z^{-1} \frac{  \left( N/2 \right)^{n/2} e^{-n \pi i/4} q^{n^{3}-n/4} }{Q(0)} \\
  & = & \frac{(-i)^{-n} \left( N/2 \right)^{n/2} e^{-n \pi i/4} q^{-3(\rho,\rho)} }
       { \prod_{ \alpha \in \overline{\Phi}^{+}_{0}} 
	       \left( q^{(\alpha,\rho)} - q^{-(\alpha,\rho)} \right)
	  \prod_{ \beta \in \Phi^{+}_{1}}    
	  \left( q^{(\beta,\rho)} + q^{-(\beta,\rho)} \right) }.
\end{eqnarray*}

We now calculate ${\cal{F}}(S^{2} \times S^{1})_{U_{q}^{(N/2)}(so(2n+1))}$.
With $N$ as given in the calculation of ${\cal{F}}(S^{2} \times S^{1})_{U_{q}^{(N)}(osp(1|2n))}$ above,
fix $\overline{N} = N/2$ and $\hat{q} = q^{2}$. 
From \cite[p. 635]{z4}, we immediately have
$${\cal{F}}(S^{2} \times S^{1})_{U_{q}^{(N/2)}(so(2n+1))} = 
\frac{(-1)^{\left|\Phi^{+}_{so(2n+1)}\right|}
(\hat{q})^{-(1+(\overline{N}+1)/2)(2\rho,\rho)}\big(G_{1}(\hat{q})\big)^{n}    }
{Q'_{\hat{q}}(0)},$$
where $G_{k}(\hat{q})$ is a Gaussian sum:
$\displaystyle{ G_{k}(\hat{q}) = \sum_{j=0}^{\overline{N}-1} (\hat{q})^{k j^{2}} }$,
and
$$Q'_{\hat{q}}(\mu) = 
  \sum_{\sigma \in {\cal{W}}} \epsilon(\sigma) (\hat{q})^{(\sigma(2\mu + 2\rho),\rho)}.$$
(Compare this to 
$Q(\mu) = \sum_{\sigma \in {\cal{W}}} \epsilon'(\sigma) q^{(\sigma(2\mu+2\rho),\rho)}$.)
From \cite[p. 150]{kl}, we have 
$$G_{1}(\hat{q}) = \left\{  \begin{array}{ll}
(\overline{N})^{1/2} & \mbox{if } \overline{N} \equiv 1 \pmod{4}, \\
i(\overline{N})^{1/2} & \mbox{if } \overline{N} \equiv 3 \pmod{4}.
\end{array} \right.$$
As $\Phi^{+}_{so(2n+1)} = \left\{ \epsilon_{i}, \epsilon_{j} \pm \epsilon_{k} | \ 
    1 \leq i \leq n, \ 1 \leq j < k \leq n \right\}$, we have
 $(-1)^{ \big| \Phi^{+}_{so(2n+1)} \big| } = (-1)^{n}$, 
which leads to the following expression for 
${\cal{F}}(S^{2} \times S^{1})_{U_{q}^{(N/2)}(so(2n+1))}$:
$$ {\cal{F}}(S^{2} \times S^{1})_{U_{q}^{(N/2)}(so(2n+1))} 
    =   \frac{ (-1)^{n} (-i)^{(2\rho, 2\rho)} \left(N/2\right)^{n/2} \hat{q}^{-3(\rho,\rho)} }
             { \prod_{ \alpha \in \overline{\Phi}^{+}_{0} \cup \Phi^{+}_{1} } 
	       \left( \hat{q}^{(\alpha,\rho)} - \hat{q}^{-(\alpha,\rho)}   \right)  } 
        \times \left\{ \begin{array}{ll}
1,     & \mbox{if } \overline{N} \equiv 1 \pmod{4}, \\
i^{n}, & \mbox{if } \overline{N} \equiv 3 \pmod{4}.
\end{array} \right.$$
Elementary algebra shows that
${\cal{F}}(S^{2} \times S^{1})_{U_{q}^{(N)}(osp(1|2n))}
      = {\cal{F}}(S^{2} \times S^{1})_{U_{q}^{(N/2)}(so(2n+1))}$ if and only if
\begin{equation}
\label{eq:michelledavidson(1)}
e^{n \pi i/4} 
= i^{(2\rho,2\rho)} q^{3(\rho,\rho)} 
      \prod_{ \alpha \in \overline{\Phi}^{+}_{0}} 
	       \left( q^{(\alpha,\rho)} + q^{-(\alpha,\rho)} \right)
	  \prod_{ \beta \in \Phi^{+}_{1}}    
	  \left( q^{(\beta,\rho)} - q^{-(\beta,\rho)} \right)
\times \left\{ \begin{array}{ll}     
i^{-n},     & \mbox{if } \overline{N} \equiv 1 \pmod{4}, \\
(-1)^{n},  & \mbox{if } \overline{N} \equiv 3 \pmod{4}.
\end{array} \right.
\end{equation}
Eq. (\ref{eq:michelledavidson(1)}) is never true for $n=1$, and 
for all odd $n \geq 3$ a necessary (and not sufficient) condition for
it to be true is that
$(n^{3}-n/4) \in  \mathbb{Z}(k^{2}+k+1/4)$ where $N = 2(2k+1)$.
For each odd $n \geq 3$ 
we can easily choose a sufficiently large enough
$N$ so that $(n^{3}-n/4) \notin  \mathbb{Z}(k^{2}+k+1/4)$, 
and this relation then holds true for all $N=2(2k'+1)$ where $k' > k$.  
Thus the invariants from
$U_{q}^{(N)}(osp(1|2n))$ and $U_{q}^{(N/2)}(so(2n+1))$ are not the same.

Let us now prove the following result which has been used in the derivation of
\newline
\noindent
${\cal{F}}(S^{2} \times S^{1})_{U_{q}^{(N)}(osp(1|2n))}$.
\begin{lemma}
\label{lem:ballson}
$\sum_{\overline{\lambda} \in X_{N}} \left(Q(\overline{\lambda})\right)^{2} = (2N)^{n} n!$
\end{lemma}
\begin{proof}
Recall that $Q(\overline{\lambda})$ is defined by
$Q(\overline{\lambda}) 
  = \sum_{\sigma \in {\cal{W}}} \epsilon'(\sigma) q^{2(\sigma(\rho), \overline{\lambda} + \rho)}$
for each $\overline{\lambda} \in X_{N}$
Now we claim (i) and (ii) below, where $\sigma, w \in {\cal{W}}$:  
\begin{itemize}
\item[(i)]   
$\displaystyle{\sum_{\overline{\lambda} \in X_{N}}
 q^{2 ( \sigma(\rho)+ w(\rho),\overline{\lambda} + \rho )} = 0 }$, if 
 $\sigma(\rho)+ w(\rho) \neq 0$, 
\item[(ii)] 
$\displaystyle{ \epsilon'(\sigma) \epsilon'(w) \sum_{\overline{\lambda} \in X_{N}}
q^{2 ( \sigma(\rho)+ w(\rho),\overline{\lambda} + \rho )} = N^{n} }$, 
if $\sigma(\rho)+ w(\rho) = 0$.
\end{itemize}
We will prove these results momentarily.
Using these results we calculate that:
\begin{eqnarray}
\sum_{\overline{\lambda} \in X_{N}} \left(Q(\overline{\lambda})\right)^{2}
 & = & \sum_{\overline{\lambda} \in X_{N}} \sum_{\sigma,w \in {\cal{W}}} 
        \epsilon'(\sigma) \epsilon'(w) 
	q^{2 ( \sigma(\rho)+ w(\rho),\overline{\lambda} + \rho )}  \nonumber \\
 & = & \sum_{\overline{\lambda} \in X_{N}} 
       \sum_{\stackrel{\sigma,w \in {\cal{W}}}{\sigma(\rho)+ w(\rho) \neq 0}} 
        \epsilon'(\sigma) \epsilon'(w) 
	q^{2 ( \sigma(\rho)+ w(\rho),\overline{\lambda} + \rho )} \nonumber \\
 &   & \hspace{5mm}
       + \sum_{\overline{\lambda} \in X_{N}} 
       \sum_{\stackrel{\sigma,w \in {\cal{W}}}{\sigma(\rho)+ w(\rho) = 0}} 
        \epsilon'(\sigma) \epsilon'(w) 
	q^{2 ( \sigma(\rho)+ w(\rho),\overline{\lambda} + \rho )} \nonumber \\
 & = & \sum_{\stackrel{\sigma,w \in {\cal{W}}}{w=-\sigma}} 
        \epsilon'(\sigma) \epsilon'(w)
	\sum_{\overline{\lambda} \in X_{N}}  
	q^{2 ( \sigma(\rho)+ w(\rho),\overline{\lambda} + \rho )} \label{eq:sachione} \\
 & = & |{\cal{W}}|  N^{n} 
   =  \left(2N\right)^{n} n! \label{eq:sachitwo}
\end{eqnarray}
where we obtain (\ref{eq:sachione}) from the fact that 
$\sigma(\rho)+ w(\rho) = 0$ if and only if $w = -\sigma$, 
and we obtain (\ref{eq:sachitwo}) from the fact that the order of ${\cal{W}}$ is $2^{n} n!$

Now we prove the claimed results (i) and (ii) above.
We prove (i).  Assume that $\sigma, w \in {\cal{W}}$ are such that 
$\sigma(\rho)+ w(\rho) \neq 0$, and let us write
$\sigma(\rho) + w(\rho) =  \sum_{i=1}^{n} \mu_{i} \epsilon_{i}$.
The properties of the reflections generated by the elements of ${\cal{W}}$ 
mean that $\mu_{i} \in \mathbb{Z}$ for each $i$.

By assumption, $\mu_{i} \neq 0$ for some $i=1, \ldots, n$.  Fix such an $i$, then
by considering the action of $\sigma, w \in {\cal{W}}$ on $\rho$, we have
\begin{equation}
\label{eq:muraleesix}
2 \leq |2 \mu_{i}| \leq 4n-2.
\end{equation}
Now the assumption that $\Lambda_{N}^{+}$ is non-empty means that $N \geq 4n+2$.
To see this, recall that $\Lambda_{N}^{+}$ is defined when $N \equiv 2 \pmod{4}$ by
$\Lambda_{N}^{+} = 
  \left\{ \lambda \in {\cal{P}}^{+} | \ 0 \leq \lambda_{1} \leq N/4-n-1/2 \right\}$.
Then $q^{2 \mu_{i}} \neq 1$ from Eq. (\ref{eq:muraleesix}).
To complete the proof of (i), all we need is the following trivial calculation: 
\begin{eqnarray*}
\sum_{\overline{\lambda} \in X_{N}} q^{ ( \sigma(\rho)+ w(\rho),2\overline{\lambda})} 
 & = & \sum_{\lambda_{1}, \ldots, \lambda_{i-1}, \lambda_{i+1},\ldots, \lambda_{n}=0}^{N-1}
       q^{2\mu_{1} \lambda_{1} + \cdots + 2\mu_{i-1} \lambda_{i-1} + 
        2\mu_{i+1} \lambda_{i+1} + \cdots + 2 \mu_{n} \lambda_{n}}
       \sum_{\lambda_{i}=0}^{N-1} q^{2 \mu_{i} \lambda_{i}} \\
 & = &  0.
\end{eqnarray*}

We now prove (ii).
Assume that $\sigma, w \in {\cal{W}}$ are such that
$\sigma(\rho)+w(\rho) = 0$, then 
$w = -\sigma = w_{0} \sigma$ where
$w_{0} = \sigma_{\epsilon_{1}} \sigma_{\epsilon_{2}} \cdots \sigma_{\epsilon_{n}}$ is
the longest element of ${\cal{W}}$, and we have
$\epsilon'(w) = \epsilon'(w_{0}) \epsilon'(\sigma) = \epsilon'(\sigma)$.
Finally, 
$$\epsilon'(\sigma) \epsilon'(w)\sum_{\overline{\lambda} \in X_{N}} 
q^{2 ( \sigma(\rho)+ w(\rho),\overline{\lambda} + \rho )} = 
\sum_{\overline{\lambda} \in X_{N}} 1 = N^{n}.$$

\end{proof}

\end{section}

\begin{section}{Some side results}
\label{sec:hopefullylast?}
\markright{\text{Some side results}}

\begin{subsection}{An observation}

We now discuss an observation of
Turaev and Wenzl \cite{tw}, who showed that in certain circumstances 
the Reshetikhin-Turaev $3$-manifold invariants 
could be calculated by taking a weighted sum of the
$F\big( \Gamma( L, \lambda) \big)$ where 
the link is cabled and each component of the cabled link is only ever coloured with the same module. 
Their 
observations equally apply to the $3$-manifold invariants constructed from 
pseudo-modular Hopf algebras in this thesis.
We briefly discuss this here and refer the reader to \cite{tw} for details and a proof.

Let $A$ be a modular or pseudo-modular Hopf algebra with universal $R$-matrix $R$, and let
$\{V_{\lambda} | \ \lambda \in I \}$ 
be the collection of $A$-modules used to construct the
$3$-manifold invariants.
For each $\mu \in I$,
let ${\cal{C}}_{t}(V_{\mu})$ be the algebra over $\mathbb{C}$ generated by the elements
$$\left\{ \check{\cal{R}}_{i}^{\pm 1} \in End_{A}\left(V_{\mu}\right)^{\otimes t} | 
\ 1 \leq i \leq t-1 \right\}, \hspace{10mm} \mbox{where}$$
$$
\check{\cal{R}}_{i}^{\pm 1} (v_{j_{1}} \otimes  \cdots \otimes v_{j_{t}}) 
= 
v_{j_{1}} \otimes  \cdots \otimes v_{j_{i-1}} \otimes
P \circ \big(R^{\pm 1}(v_{j_{i}} \otimes v_{j_{i+1}})\big) \otimes v_{j_{i+2}} \otimes \cdots \otimes
v_{j_{t}}.
$$
If each $V_{\lambda}$ is isomorphic to $p_{\lambda}(V_{\mu})^{\otimes t}$ for some idempotent
$p_{\lambda} \in {\cal{C}}_{t}(V_{\mu})$, for some $\mu$ in $I$,
 we say that $V_{\mu}$ is a {\emph{generating}} module.
For $U_{q}^{(N)}(osp(1|2n))$ at even roots of unity, the fundamental module $V$ is generating.
Turaev and Wenzl's observations only apply if the algebra has a generating module; 
 we only consider such algebras below and denote the generating $A$-module by $V$.

The module $V_{\lambda}$ is isomorphic to $p_{\lambda}(V^{\otimes t})$ for some idempotent
$p_{\lambda} \in {\cal{C}}_{t}$ where we can write
$p_{\lambda} = \sum_{j=1}^{r} c_{j} \mathscr{R}_{j}$, $c_{j} \in \mathbb{C}$. 
Here each $\mathscr{R}_{j}$ is an ordered product 
in the $\check{\cal{R}}_{i}^{\pm 1}$.

Now in calculating the 3-manifold invariants, one takes a weighted sum of the
$F\big(\Gamma(L,\lambda)\big)$ where the sum is over all different possible colourings $\lambda$
of $\Gamma(L)$.
The fact that $V_{\lambda}$ is isomorphic to $p_{\lambda}(V^{\otimes t})$ means that
we can express $F\big(\Gamma(L,\lambda)\big)$ differently: we can write
$F\big(\Gamma(L,\lambda)\big) = \sum_{j=1}^{r} c_{j} F\big(\Gamma(L_{j}',V)\big)$
where $L_{j}'$ is a link obtained from cabling  $L$ and  $V$ means that each component of
$\Gamma(L)$ is coloured with the generating $A$-module $V$.

To see this, we consider a framed oriented
knot $L$; the multicomponent case follows.
Regard $\Gamma(L)$ as being coloured with $\lambda \in I$.
For each $j=1, \ldots, r$, let $\overline{L}_{j}'$ be a framed oriented 
link with $t$ components obtained by cabling $L$.
By referring to $\mathscr{R}_{j}$ we construct a new link $L'_{j}$.
Of course, 
$\mathscr{R}_{j} = \check{\cal{R}}^{\epsilon_{1}}_{k_{1}}\check{\cal{R}}^{\epsilon_{2}}_{k_{2}}
 \cdots \check{\cal{R}}^{\epsilon_{m}}_{k_{m}}$, where $\epsilon_{l} \in \{-1, +1\}$ and
 $k_{l} \in \{ 1, 2, \ldots, t-1 \}$.

Now let $b_{j} = \sigma^{\epsilon_{1}}_{k_{1}} \sigma^{\epsilon_{2}}_{k_{2}}
 \cdots \sigma^{\epsilon_{m}}_{k_{m}}$ be an element of $B_{t}$, the Braid group on $t$ strings, 
  where 
$\sigma^{+1}$ and $\sigma^{-1}$ are the elements of $B_{2}$ in Figure \ref{fig:tom}, and
$\sigma^{\pm 1}_{k}  \in B_{t}$ acts as $\sigma^{\pm 1}$ on the $k^{th}$ and
$(k+1)^{st}$ components of the $(t,t)$-tangle (from the left) and leaves 
all other components unchanged.
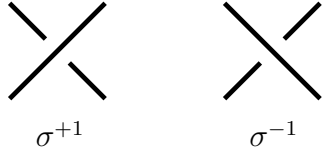
\begin{figure}[hbt]
\begin{center}
 \input{braidgen1100.pstex_t} 
\caption{The elements $\sigma^{+1}$, $\sigma^{-1}$ of $B_{2}$} \label{fig:tom}
\end{center}
\end{figure}
Now we construct the new $t$-component link $L'_{j}$ by doing the following.
Consider the intersection of $\overline{L}_{j}'$ with a $3$-disk $D^{3}$ in such a way that
all the components of $\overline{L}_{j}'$ in the intersection are parallel and oriented downwards.
We obtain $L'_{j}$ by
replacing the oriented $(t,t)$-tangle in the intersection with an 
oriented $(t,t)$-tangle given by applying $b_{j} \in B_{t}$ to the $t$ parallel components.
We do this by
firstly applying $\sigma_{k_{1}}^{\epsilon_{1}}$ to the top of the $t$ components, 
and then inductively applying
$\sigma_{k_{r+1}}^{\epsilon_{r+1}}$ below  
$\sigma_{k_{r}}^{\epsilon_{r}}$ for each $r=1, \ldots, m-1$.

After obtaining each $L'_{j}$, we have
$$
F\big(\Gamma(L,\lambda)\big) = \sum_{j=1}^{r} c_{j} F\big(\Gamma(L'_{j},V)\big),
$$
which follows from the properties of $F$.

\end{subsection}

\begin{subsection}{A further result}

\begin{theorem}
\label{lem:dalektron4}
The algebra
$U_{q}^{(N)}(osp(1|2n))$ and the set
$\{V_{\lambda} | \ \lambda \in \Lambda_{N}^{+} \}$ of modules 
is not a pseudo-modular Hopf algebra when $N \geq 4$ satisfies $N \equiv 0 \pmod{4}$.
\end{theorem}
\begin{proof}
We will show that $U_{q}^{(N)}(osp(1|2n))$ and the set of modules 
do not satisfy Axiom (V) of a pseudo-modular Hopf algebra.
To do this, we will show that the set (\ref{eq:sigeq1}) 
of equations is inconsistent if $N \equiv 0 \pmod{4}$, $N \geq 4$.

We will show that for each $\mu \in \Lambda_{N}^{+}$ there exists some 
$\mu' \in \Lambda_{N}^{+}$, where $\mu' \neq \mu$, such that
\begin{itemize}
\item[(i)] $q^{-(\mu',\mu' + 2\rho)} = -q^{-(\mu,\mu+2\rho)}$,  and
\item[(ii)] $S_{\lambda,\mu'}/Q(\mu') = S_{\lambda,\mu}/Q(\mu)$, and
\end{itemize}
we have the following pair of inconsistent equations:
\begin{equation}
\label{eq:aazz}
q^{-(\mu,\mu+2\rho)}  =  \sum_{\lambda \in \Lambda_{N}^{+}}
d_{\lambda}\chi_{\lambda}(v^{-1})S_{\lambda,\mu}/Q(\mu),
\end{equation}
\begin{equation}
\label{eq:aazz1}
q^{-(\mu',\mu'+2\rho)}  = \sum_{\lambda \in \Lambda_{N}^{+}}
d_{\lambda}\chi_{\lambda}(v^{-1})S_{\lambda,\mu'}/Q(\mu').
\end{equation}

Let $N \geq 4$ satisfy $N \equiv 0 \pmod{4}$.
Let $\mu \in \Lambda_{N}^{+}$ and fix $s_{\epsilon_{1}} \in {\cal{W}}$ 
to be the element of ${\cal{W}}$ that reflects $H^{*}$ in the subspace 
$hp_{\epsilon_{1}} = \{x \in H^{*} | \ (x,\epsilon_{1}) = 0 \}$. The element 
$s_{\epsilon_{1}}$ acts on $\nu = \sum_{i=1}^{n}\nu_{i} \epsilon_{i} \in H^{*}$ via 
$$s_{\epsilon_{1}}:\nu \mapsto -\nu_{1}\epsilon_{1} + \sum_{i=2}^{n}\nu_{i}\epsilon_{i}.$$
Define the map
$\sigma_{\epsilon_{1},N'}:X  \rightarrow X$ by
$$\sigma_{\epsilon_{1},N'}: \mu \mapsto s_{\epsilon_{1}}(\mu + \rho) - \rho + N'\epsilon_{1}.$$  
In the foregoing we will use Proposition \ref{prop:dalektron3}, to which the reader is temporarily
referred.
Set $\mu' = \sigma_{\epsilon_{1},N'}(\mu)$: the next two calculations will show that
$q^{-(\mu',\mu'+2\rho)} = -q^{-(\mu,\mu+2\rho)}$.  Now
\begin{eqnarray*}
q^{-(\mu + N'\epsilon_{1},\mu+N'\epsilon_{1}+2\rho)} & = &
q^{-\left((\mu,\mu + 2\rho)+N(\mu,\epsilon_{1})+N'(2\rho,\epsilon_{1})+(N')^{2}\right)} \\
& = & q^{-\left((\mu,\mu + 2\rho) + N'(2\rho,\epsilon_{1})\right)} \\
& = & -q^{-(\mu,\mu+2\rho)},
\end{eqnarray*}
and for any $w \in {\cal{W}}$ we have
$$
q^{-(w(\mu + \rho)-\rho,w(\mu + \rho)-\rho + 2\rho)} = 
q^{-(w(\mu + \rho),w(\mu + \rho))+(\rho,\rho)} = q^{-(\mu,\mu+2\rho)},
$$
thus $q^{-(\mu',\mu'+2\rho)} = -q^{-(\mu,\mu+2\rho)}$.

The equality $S_{\lambda,\mu'}/Q(\mu') = S_{\lambda,\mu}/Q(\mu)$ results from the following properties 
of $S_{\lambda,\mu}$, which derive from direct calculations and from Lemma \ref{lem:properties}:
\begin{eqnarray*}
S_{\lambda,s_{\epsilon_{1}}(\mu + \rho)-\rho} /Q(s_{\epsilon_{1}}(\mu + \rho)-\rho) & = & 
S_{\lambda,\mu}/Q(\mu), \\
S_{\lambda,\mu + N'\epsilon_{1}}/Q(\mu + N'\epsilon_{1})
 & = & S_{\lambda,\mu}/Q(\mu).
\end{eqnarray*}
The second relation follows from the equality
$S_{\lambda,\mu + N' \epsilon_{1}} = -S_{\lambda,\mu}$.

It follows that for each $\mu \in \Lambda_{N}^{+}$ there exists some $\mu' \in \Lambda_{N}^{+}$ 
where $\mu' \neq \mu$, such that we have 
an inconsistent pair of equations (Eqs. (\ref{eq:aazz}) and (\ref{eq:aazz1})).
As $\mu' = \sigma_{\epsilon_{1},N'}(\mu)$ and $\mu = \sigma_{\epsilon_{1},N'}(\mu')$,
there are $|\Lambda_{N}^{+}|/2$ such pairs. 
\end{proof}

\begin{proposition}
\label{prop:dalektron3}
Let $\sigma_{\epsilon_{1},N'} : X \rightarrow X$ 
be a map given in Theorem \ref{lem:dalektron4} and
let $\mu \in \Lambda_{N}^{+}$, then
$\sigma_{\epsilon_{1},N'} (\mu) = (N'-\mu_{1}-2n+1)\epsilon_{1} + \sum_{i=2}^{n} \mu_{i}\epsilon_{i}$.
Furthermore,
\begin{itemize}
\item[(i)]  $\sigma_{\epsilon_{1},N'}\big(\sigma_{\epsilon_{1},N'}(\mu)\big) = \mu$, 
\item[(ii)] $\sigma_{\epsilon_{1},N'}(\mu) \neq \mu$, 
\item[(iii)] $\sigma_{\epsilon_{1},N'}(\mu) \in \Lambda_{N}^{+}$.
\end{itemize}
\end{proposition}
\begin{proof}
The action of $\sigma_{\epsilon_{1},N'}$ on $\mu$ is shown by trivial calculations.  
The proofs of (i)--(iii) are:
\begin{itemize}
\item[(i)] $\sigma_{\epsilon_{1},N'}\big(\sigma_{\epsilon_{1},N'}(\mu)\big) =
\sigma_{\epsilon_{1},N'}\big((N' - \mu_{1} - 2n + 1)\epsilon_{1} + \sum_{i=2}^{n}\mu_{i}\epsilon_{i}\big) = \mu$.
\item[(ii)] Assume that $\sigma_{\epsilon_{1},N'}(\mu) = \mu$, then 
$(N' - \mu_{1} - 2n + 1)\epsilon_{1} = \mu_{1}\epsilon_{1}$, which implies 
that $2 \mu_{1} = N' - 2n + 1$ and $\mu_{1} \in \mathbb{Z} + 1/2$.
However, this is not possible as $\mu \in \bigoplus_{i=1}^{n} \mathbb{Z} \epsilon_{i}$ 
and therefore $\sigma_{\epsilon_{1},N'}(\mu) \neq \mu$.
\item[(iii)] 
The definition of $\Lambda_{N}^{+}$ implies that $\mu' \in \Lambda_{N}^{+}$ if and only if 
$0 \leq \mu'_{1}+\mu'_{2} < N'-2n+2$.  
Now as $\mu \in \Lambda_{N}^{+}$, the components of $\mu$ satisfy
$-1 < \mu_{1}-\mu_{2} \leq N'-2n+1$ which we can rewrite as
\begin{equation}
\label{eq:maybethelastequation}
0 \leq N'-\mu_{1}-2n+1+\mu_{2}<N'-2n+2.
\end{equation}
The statement of the proposition tells us that 
$\sigma_{\epsilon_{1},N'} (\mu) = (N'-\mu_{1}-2n+1)\epsilon_{1} + \sum_{i=2}^{n} \mu_{i}\epsilon_{i}$, 
that is $\mu'_{1} = N'-\mu_{1}-2n+1$ and $\mu'_{2}=\mu_{2}$.
This allows us to rewrite (\ref{eq:maybethelastequation}) as
$0 \leq \mu'_{1}+\mu'_{2}<N'-2n+2$, which is precisely the condition under which
$\sigma_{\epsilon_{1},N'}(\mu)$ is an element of $\Lambda_{N}^{+}$.
\end{itemize}
\end{proof}

\end{subsection}

\end{section}

\end{chapter}

%% file: reid100.pstex_t
\begin{picture}(0,0)%
\includegraphics{reid.pstex}%
\end{picture}%
\setlength{\unitlength}{3947sp}%
\begingroup\makeatletter\ifx\SetFigFont\undefined
\def\x#1#2#3#4#5#6#7\relax{\def\x{#1#2#3#4#5#6}}%
\expandafter\x\fmtname xxxxxx\relax \def\y{splain}%
\ifx\x\y   
\gdef\SetFigFont#1#2#3{%
  \ifnum #1<17\tiny\else \ifnum #1<20\small\else
  \ifnum #1<24\normalsize\else \ifnum #1<29\large\else
  \ifnum #1<34\Large\else \ifnum #1<41\LARGE\else
     \huge\fi\fi\fi\fi\fi\fi
  \csname #3\endcsname}%
\else
\gdef\SetFigFont#1#2#3{\begingroup
  \count@#1\relax \ifnum 25<\count@\count@25\fi
  \def\x{\endgroup\@setsize\SetFigFont{#2pt}}%
  \expandafter\x
    \csname \romannumeral\the\count@ pt\expandafter\endcsname
    \csname @\romannumeral\the\count@ pt\endcsname
  \csname #3\endcsname}%
\fi
\fi\endgroup
\begin{picture}(5356,4016)(3055,-5887)
\put(4322,-3630){\makebox(0,0)[lb]{\smash{\SetFigFont{10}{12.0}{rm}{\color[rgb]{0,0,0}$I$}%
}}}
\put(7107,-3614){\makebox(0,0)[lb]{\smash{\SetFigFont{10}{12.0}{rm}{\color[rgb]{0,0,0}$II$}%
}}}
\put(5488,-5847){\makebox(0,0)[lb]{\smash{\SetFigFont{10}{12.0}{rm}{\color[rgb]{0,0,0}$III$}%
}}}
\end{picture}

%% file: linkingnumbers100.pstex_t
\begin{picture}(0,0)%
\includegraphics{linkingnumbers.pstex}%
\end{picture}%
\setlength{\unitlength}{3947sp}%
\begingroup\makeatletter\ifx\SetFigFont\undefined
\def\x#1#2#3#4#5#6#7\relax{\def\x{#1#2#3#4#5#6}}%
\expandafter\x\fmtname xxxxxx\relax \def\y{splain}%
\ifx\x\y   
\gdef\SetFigFont#1#2#3{%
  \ifnum #1<17\tiny\else \ifnum #1<20\small\else
  \ifnum #1<24\normalsize\else \ifnum #1<29\large\else
  \ifnum #1<34\Large\else \ifnum #1<41\LARGE\else
     \huge\fi\fi\fi\fi\fi\fi
  \csname #3\endcsname}%
\else
\gdef\SetFigFont#1#2#3{\begingroup
  \count@#1\relax \ifnum 25<\count@\count@25\fi
  \def\x{\endgroup\@setsize\SetFigFont{#2pt}}%
  \expandafter\x
    \csname \romannumeral\the\count@ pt\expandafter\endcsname
    \csname @\romannumeral\the\count@ pt\endcsname
  \csname #3\endcsname}%
\fi
\fi\endgroup
\begin{picture}(2316,1132)(1768,-1460)
\put(2026,-1411){\makebox(0,0)[lb]{\smash{\SetFigFont{12}{14.4}{rm}{\color[rgb]{0,0,0}$+1$}%
}}}
\put(3526,-1411){\makebox(0,0)[lb]{\smash{\SetFigFont{12}{14.4}{rm}{\color[rgb]{0,0,0}$-1$}%
}}}
\end{picture}

%% file: kirby1100.pstex_t
\begin{picture}(0,0)%
\includegraphics{kirby1.pstex}%
\end{picture}%
\setlength{\unitlength}{3947sp}%
\begingroup\makeatletter\ifx\SetFigFont\undefined
\def\x#1#2#3#4#5#6#7\relax{\def\x{#1#2#3#4#5#6}}%
\expandafter\x\fmtname xxxxxx\relax \def\y{splain}%
\ifx\x\y   
\gdef\SetFigFont#1#2#3{%
  \ifnum #1<17\tiny\else \ifnum #1<20\small\else
  \ifnum #1<24\normalsize\else \ifnum #1<29\large\else
  \ifnum #1<34\Large\else \ifnum #1<41\LARGE\else
     \huge\fi\fi\fi\fi\fi\fi
  \csname #3\endcsname}%
\else
\gdef\SetFigFont#1#2#3{\begingroup
  \count@#1\relax \ifnum 25<\count@\count@25\fi
  \def\x{\endgroup\@setsize\SetFigFont{#2pt}}%
  \expandafter\x
    \csname \romannumeral\the\count@ pt\expandafter\endcsname
    \csname @\romannumeral\the\count@ pt\endcsname
  \csname #3\endcsname}%
\fi
\fi\endgroup
\begin{picture}(5340,1290)(4337,-4471)
\put(5243,-4422){\makebox(0,0)[lb]{\smash{\SetFigFont{10}{12.0}{rm}{\color[rgb]{0,0,0}$\kappa_{+}^{(0)}$}%
}}}
\put(8558,-4422){\makebox(0,0)[lb]{\smash{\SetFigFont{10}{12.0}{rm}{\color[rgb]{0,0,0}$\kappa_{-}^{(0)}$}%
}}}
\put(5243,-4422){\makebox(0,0)[lb]{\smash{\SetFigFont{10}{12.0}{rm}{\color[rgb]{0,0,0}$\kappa_{+}^{(0)}$}%
}}}
\put(8558,-4422){\makebox(0,0)[lb]{\smash{\SetFigFont{10}{12.0}{rm}{\color[rgb]{0,0,0}$\kappa_{-}^{(0)}$}%
}}}
\put(6441,-3612){\makebox(0,0)[lb]{\smash{\SetFigFont{9}{10.8}{rm}{\color[rgb]{0,0,0}nothing}%
}}}
\put(9677,-3670){\makebox(0,0)[lb]{\smash{\SetFigFont{9}{10.8}{rm}{\color[rgb]{0,0,0}nothing}%
}}}
\end{picture}

%% file: Kirby2100.pstex_t
\begin{picture}(0,0)%
\includegraphics{Kirby2.pstex}%
\end{picture}%
\setlength{\unitlength}{3947sp}%
\begingroup\makeatletter\ifx\SetFigFont\undefined
\def\x#1#2#3#4#5#6#7\relax{\def\x{#1#2#3#4#5#6}}%
\expandafter\x\fmtname xxxxxx\relax \def\y{splain}%
\ifx\x\y   
\gdef\SetFigFont#1#2#3{%
  \ifnum #1<17\tiny\else \ifnum #1<20\small\else
  \ifnum #1<24\normalsize\else \ifnum #1<29\large\else
  \ifnum #1<34\Large\else \ifnum #1<41\LARGE\else
     \huge\fi\fi\fi\fi\fi\fi
  \csname #3\endcsname}%
\else
\gdef\SetFigFont#1#2#3{\begingroup
  \count@#1\relax \ifnum 25<\count@\count@25\fi
  \def\x{\endgroup\@setsize\SetFigFont{#2pt}}%
  \expandafter\x
    \csname \romannumeral\the\count@ pt\expandafter\endcsname
    \csname @\romannumeral\the\count@ pt\endcsname
  \csname #3\endcsname}%
\fi
\fi\endgroup
\begin{picture}(4017,2628)(2068,-5734)
\put(3990,-5685){\makebox(0,0)[lb]{\smash{\SetFigFont{10}{12.0}{rm}{\color[rgb]{0,0,0}$\kappa_{+}$}%
}}}
\end{picture}

%% file: Kirby3100.pstex_t
\begin{picture}(0,0)%
\includegraphics{Kirby3.pstex}%
\end{picture}%
\setlength{\unitlength}{3947sp}%
\begingroup\makeatletter\ifx\SetFigFont\undefined
\def\x#1#2#3#4#5#6#7\relax{\def\x{#1#2#3#4#5#6}}%
\expandafter\x\fmtname xxxxxx\relax \def\y{splain}%
\ifx\x\y   
\gdef\SetFigFont#1#2#3{%
  \ifnum #1<17\tiny\else \ifnum #1<20\small\else
  \ifnum #1<24\normalsize\else \ifnum #1<29\large\else
  \ifnum #1<34\Large\else \ifnum #1<41\LARGE\else
     \huge\fi\fi\fi\fi\fi\fi
  \csname #3\endcsname}%
\else
\gdef\SetFigFont#1#2#3{\begingroup
  \count@#1\relax \ifnum 25<\count@\count@25\fi
  \def\x{\endgroup\@setsize\SetFigFont{#2pt}}%
  \expandafter\x
    \csname \romannumeral\the\count@ pt\expandafter\endcsname
    \csname @\romannumeral\the\count@ pt\endcsname
  \csname #3\endcsname}%
\fi
\fi\endgroup
\begin{picture}(3826,2710)(3184,-4858)
\put(4942,-4809){\makebox(0,0)[lb]{\smash{\SetFigFont{10}{12.0}{rm}{\color[rgb]{0,0,0}$\kappa_{-}$}%
}}}
\end{picture}

%% file: band1100.pstex_t
\begin{picture}(0,0)%
\includegraphics{band1.pstex}%
\end{picture}%
\setlength{\unitlength}{3947sp}%
\begingroup\makeatletter\ifx\SetFigFont\undefined
\def\x#1#2#3#4#5#6#7\relax{\def\x{#1#2#3#4#5#6}}%
\expandafter\x\fmtname xxxxxx\relax \def\y{splain}%
\ifx\x\y   
\gdef\SetFigFont#1#2#3{%
  \ifnum #1<17\tiny\else \ifnum #1<20\small\else
  \ifnum #1<24\normalsize\else \ifnum #1<29\large\else
  \ifnum #1<34\Large\else \ifnum #1<41\LARGE\else
     \huge\fi\fi\fi\fi\fi\fi
  \csname #3\endcsname}%
\else
\gdef\SetFigFont#1#2#3{\begingroup
  \count@#1\relax \ifnum 25<\count@\count@25\fi
  \def\x{\endgroup\@setsize\SetFigFont{#2pt}}%
  \expandafter\x
    \csname \romannumeral\the\count@ pt\expandafter\endcsname
    \csname @\romannumeral\the\count@ pt\endcsname
  \csname #3\endcsname}%
\fi
\fi\endgroup
\begin{picture}(2708,1481)(1302,-1694)
\put(3301,-1636){\makebox(0,0)[lb]{\smash{\SetFigFont{12}{14.4}{rm}{\color[rgb]{0,0,0}$L_{j}$}%
}}}
\put(1576,-1636){\makebox(0,0)[lb]{\smash{\SetFigFont{12}{14.4}{rm}{\color[rgb]{0,0,0}$L_{i}$}%
}}}
\end{picture}

%% file: band2100.pstex_t
\begin{picture}(0,0)%
\includegraphics{band2.pstex}%
\end{picture}%
\setlength{\unitlength}{3947sp}%
\begingroup\makeatletter\ifx\SetFigFont\undefined
\def\x#1#2#3#4#5#6#7\relax{\def\x{#1#2#3#4#5#6}}%
\expandafter\x\fmtname xxxxxx\relax \def\y{splain}%
\ifx\x\y   
\gdef\SetFigFont#1#2#3{%
  \ifnum #1<17\tiny\else \ifnum #1<20\small\else
  \ifnum #1<24\normalsize\else \ifnum #1<29\large\else
  \ifnum #1<34\Large\else \ifnum #1<41\LARGE\else
     \huge\fi\fi\fi\fi\fi\fi
  \csname #3\endcsname}%
\else
\gdef\SetFigFont#1#2#3{\begingroup
  \count@#1\relax \ifnum 25<\count@\count@25\fi
  \def\x{\endgroup\@setsize\SetFigFont{#2pt}}%
  \expandafter\x
    \csname \romannumeral\the\count@ pt\expandafter\endcsname
    \csname @\romannumeral\the\count@ pt\endcsname
  \csname #3\endcsname}%
\fi
\fi\endgroup
\begin{picture}(2954,1564)(1127,-1694)
\put(3451,-736){\makebox(0,0)[lb]{\smash{\SetFigFont{12}{14.4}{rm}{\color[rgb]{0,0,0}$L_{j}$}%
}}}
\put(1351,-1636){\makebox(0,0)[lb]{\smash{\SetFigFont{12}{14.4}{rm}{\color[rgb]{0,0,0}$L_{i} \#_{b} L_{j}$}%
}}}
\end{picture}

%% file: bandkir100.pstex_t
\begin{picture}(0,0)%
\includegraphics{bandkir.pstex}%
\end{picture}%
\setlength{\unitlength}{3947sp}%
\begingroup\makeatletter\ifx\SetFigFont\undefined
\def\x#1#2#3#4#5#6#7\relax{\def\x{#1#2#3#4#5#6}}%
\expandafter\x\fmtname xxxxxx\relax \def\y{splain}%
\ifx\x\y   
\gdef\SetFigFont#1#2#3{%
  \ifnum #1<17\tiny\else \ifnum #1<20\small\else
  \ifnum #1<24\normalsize\else \ifnum #1<29\large\else
  \ifnum #1<34\Large\else \ifnum #1<41\LARGE\else
     \huge\fi\fi\fi\fi\fi\fi
  \csname #3\endcsname}%
\else
\gdef\SetFigFont#1#2#3{\begingroup
  \count@#1\relax \ifnum 25<\count@\count@25\fi
  \def\x{\endgroup\@setsize\SetFigFont{#2pt}}%
  \expandafter\x
    \csname \romannumeral\the\count@ pt\expandafter\endcsname
    \csname @\romannumeral\the\count@ pt\endcsname
  \csname #3\endcsname}%
\fi
\fi\endgroup
\begin{picture}(3499,3488)(213,-2814)
\put(1800,164){\makebox(0,0)[lb]{\smash{\SetFigFont{12}{14.4}{rm}{\color[rgb]{0,0,0}$b$}%
}}}
\put(1869,-1831){\makebox(0,0)[lb]{\smash{\SetFigFont{12}{14.4}{rm}{\color[rgb]{0,0,0}$k$}%
}}}
\end{picture}

%% file: example100.pstex_t
\begin{picture}(0,0)%
\includegraphics{example.pstex}%
\end{picture}%
\setlength{\unitlength}{4144sp}%
\begingroup\makeatletter\ifx\SetFigFont\undefined
\def\x#1#2#3#4#5#6#7\relax{\def\x{#1#2#3#4#5#6}}%
\expandafter\x\fmtname xxxxxx\relax \def\y{splain}%
\ifx\x\y   
\gdef\SetFigFont#1#2#3{%
  \ifnum #1<17\tiny\else \ifnum #1<20\small\else
  \ifnum #1<24\normalsize\else \ifnum #1<29\large\else
  \ifnum #1<34\Large\else \ifnum #1<41\LARGE\else
     \huge\fi\fi\fi\fi\fi\fi
  \csname #3\endcsname}%
\else
\gdef\SetFigFont#1#2#3{\begingroup
  \count@#1\relax \ifnum 25<\count@\count@25\fi
  \def\x{\endgroup\@setsize\SetFigFont{#2pt}}%
  \expandafter\x
    \csname \romannumeral\the\count@ pt\expandafter\endcsname
    \csname @\romannumeral\the\count@ pt\endcsname
  \csname #3\endcsname}%
\fi
\fi\endgroup
\begin{picture}(5848,3187)(519,-2413)
\end{picture}

%% file: verticcomp100.pstex_t
\begin{picture}(0,0)%
\includegraphics{verticcomp.pstex}%
\end{picture}%
\setlength{\unitlength}{4144sp}%
\begingroup\makeatletter\ifx\SetFigFont\undefined
\def\x#1#2#3#4#5#6#7\relax{\def\x{#1#2#3#4#5#6}}%
\expandafter\x\fmtname xxxxxx\relax \def\y{splain}%
\ifx\x\y   
\gdef\SetFigFont#1#2#3{%
  \ifnum #1<17\tiny\else \ifnum #1<20\small\else
  \ifnum #1<24\normalsize\else \ifnum #1<29\large\else
  \ifnum #1<34\Large\else \ifnum #1<41\LARGE\else
     \huge\fi\fi\fi\fi\fi\fi
  \csname #3\endcsname}%
\else
\gdef\SetFigFont#1#2#3{\begingroup
  \count@#1\relax \ifnum 25<\count@\count@25\fi
  \def\x{\endgroup\@setsize\SetFigFont{#2pt}}%
  \expandafter\x
    \csname \romannumeral\the\count@ pt\expandafter\endcsname
    \csname @\romannumeral\the\count@ pt\endcsname
  \csname #3\endcsname}%
\fi
\fi\endgroup
\begin{picture}(1780,921)(628,-1732)
\put(2094,-1119){\makebox(0,0)[lb]{\smash{\SetFigFont{12}{14.4}{rm}{\color[rgb]{0,0,0}$\Gamma$}%
}}}
\put(2092,-1536){\makebox(0,0)[lb]{\smash{\SetFigFont{12}{14.4}{rm}{\color[rgb]{0,0,0}$\Gamma'$}%
}}}
\put(813,-1305){\makebox(0,0)[lb]{\smash{\SetFigFont{12}{14.4}{rm}{\color[rgb]{0,0,0}$\Gamma \circ \Gamma'$}%
}}}
\put(1585,-1292){\makebox(0,0)[lb]{\smash{\SetFigFont{12}{14.4}{rm}{\color[rgb]{0,0,0}$=$}%
}}}
\end{picture}

%% file: horizontcomp100.pstex_t
\begin{picture}(0,0)%
\includegraphics{horizontcomp.pstex}%
\end{picture}%
\setlength{\unitlength}{4144sp}%
\begingroup\makeatletter\ifx\SetFigFont\undefined
\def\x#1#2#3#4#5#6#7\relax{\def\x{#1#2#3#4#5#6}}%
\expandafter\x\fmtname xxxxxx\relax \def\y{splain}%
\ifx\x\y   
\gdef\SetFigFont#1#2#3{%
  \ifnum #1<17\tiny\else \ifnum #1<20\small\else
  \ifnum #1<24\normalsize\else \ifnum #1<29\large\else
  \ifnum #1<34\Large\else \ifnum #1<41\LARGE\else
     \huge\fi\fi\fi\fi\fi\fi
  \csname #3\endcsname}%
\else
\gdef\SetFigFont#1#2#3{\begingroup
  \count@#1\relax \ifnum 25<\count@\count@25\fi
  \def\x{\endgroup\@setsize\SetFigFont{#2pt}}%
  \expandafter\x
    \csname \romannumeral\the\count@ pt\expandafter\endcsname
    \csname @\romannumeral\the\count@ pt\endcsname
  \csname #3\endcsname}%
\fi
\fi\endgroup
\begin{picture}(2143,561)(668,-1529)
\put(2060,-1315){\makebox(0,0)[lb]{\smash{\SetFigFont{12}{14.4}{rm}{\color[rgb]{0,0,0}$\Gamma$}%
}}}
\put(2490,-1310){\makebox(0,0)[lb]{\smash{\SetFigFont{12}{14.4}{rm}{\color[rgb]{0,0,0}$\Gamma'$}%
}}}
\put(833,-1305){\makebox(0,0)[lb]{\smash{\SetFigFont{12}{14.4}{rm}{\color[rgb]{0,0,0}$\Gamma \otimes \Gamma'$}%
}}}
\put(1618,-1297){\makebox(0,0)[lb]{\smash{\SetFigFont{12}{14.4}{rm}{\color[rgb]{0,0,0}$=$}%
}}}
\end{picture}

%% file: ribbon1100.pstex_t
\begin{picture}(0,0)%
\includegraphics{ribbon1.pstex}%
\end{picture}%
\setlength{\unitlength}{3947sp}%
\begingroup\makeatletter\ifx\SetFigFont\undefined
\def\x#1#2#3#4#5#6#7\relax{\def\x{#1#2#3#4#5#6}}%
\expandafter\x\fmtname xxxxxx\relax \def\y{splain}%
\ifx\x\y   
\gdef\SetFigFont#1#2#3{%
  \ifnum #1<17\tiny\else \ifnum #1<20\small\else
  \ifnum #1<24\normalsize\else \ifnum #1<29\large\else
  \ifnum #1<34\Large\else \ifnum #1<41\LARGE\else
     \huge\fi\fi\fi\fi\fi\fi
  \csname #3\endcsname}%
\else
\gdef\SetFigFont#1#2#3{\begingroup
  \count@#1\relax \ifnum 25<\count@\count@25\fi
  \def\x{\endgroup\@setsize\SetFigFont{#2pt}}%
  \expandafter\x
    \csname \romannumeral\the\count@ pt\expandafter\endcsname
    \csname @\romannumeral\the\count@ pt\endcsname
  \csname #3\endcsname}%
\fi
\fi\endgroup
\begin{picture}(5540,5278)(2203,-7006)
\put(3226,-5236){\makebox(0,0)[lb]{\smash{\SetFigFont{10}{12.0}{rm}{\color[rgb]{0,0,0}$\Omega^{+}$}%
}}}
\put(6976,-3361){\makebox(0,0)[lb]{\smash{\SetFigFont{10}{12.0}{rm}{\color[rgb]{0,0,0}$X^{-}$}%
}}}
\put(6376,-5236){\makebox(0,0)[lb]{\smash{\SetFigFont{10}{12.0}{rm}{\color[rgb]{0,0,0}$\Omega^{-}$}%
}}}
\put(6376,-6961){\makebox(0,0)[lb]{\smash{\SetFigFont{10}{12.0}{rm}{\color[rgb]{0,0,0}$U^{-}$}%
}}}
\put(3226,-6961){\makebox(0,0)[lb]{\smash{\SetFigFont{10}{12.0}{rm}{\color[rgb]{0,0,0}$U^{+}$}%
}}}
\put(3451,-3361){\makebox(0,0)[lb]{\smash{\SetFigFont{10}{12.0}{rm}{\color[rgb]{0,0,0}$I^{-}$}%
}}}
\put(2326,-3361){\makebox(0,0)[lb]{\smash{\SetFigFont{10}{12.0}{rm}{\color[rgb]{0,0,0}$I^{+}$}%
}}}
\put(5176,-3361){\makebox(0,0)[lb]{\smash{\SetFigFont{10}{12.0}{rm}{\color[rgb]{0,0,0}$X^{+}$}%
}}}
\end{picture}

%% file: closure100.pstex_t
\begin{picture}(0,0)%
\includegraphics{closure.pstex}%
\end{picture}%
\setlength{\unitlength}{4144sp}%
\begingroup\makeatletter\ifx\SetFigFont\undefined
\def\x#1#2#3#4#5#6#7\relax{\def\x{#1#2#3#4#5#6}}%
\expandafter\x\fmtname xxxxxx\relax \def\y{splain}%
\ifx\x\y   
\gdef\SetFigFont#1#2#3{%
  \ifnum #1<17\tiny\else \ifnum #1<20\small\else
  \ifnum #1<24\normalsize\else \ifnum #1<29\large\else
  \ifnum #1<34\Large\else \ifnum #1<41\LARGE\else
     \huge\fi\fi\fi\fi\fi\fi
  \csname #3\endcsname}%
\else
\gdef\SetFigFont#1#2#3{\begingroup
  \count@#1\relax \ifnum 25<\count@\count@25\fi
  \def\x{\endgroup\@setsize\SetFigFont{#2pt}}%
  \expandafter\x
    \csname \romannumeral\the\count@ pt\expandafter\endcsname
    \csname @\romannumeral\the\count@ pt\endcsname
  \csname #3\endcsname}%
\fi
\fi\endgroup
\begin{picture}(4796,2524)(508,-1919)
\put(3871,-1861){\makebox(0,0)[lb]{\smash{\SetFigFont{12}{14.4}{rm}{\color[rgb]{0,0,0}$\hat{\Gamma}$}%
}}}
\put(1171,-1861){\makebox(0,0)[lb]{\smash{\SetFigFont{12}{14.4}{rm}{\color[rgb]{0,0,0}$\Gamma$}%
}}}
\put(3376,-556){\makebox(0,0)[lb]{\smash{\SetFigFont{12}{14.4}{rm}{\color[rgb]{0,0,0}$T$}%
}}}
\put(1171,-556){\makebox(0,0)[lb]{\smash{\SetFigFont{12}{14.4}{rm}{\color[rgb]{0,0,0}$T$}%
}}}
\end{picture}

%% file: Kirby+move100.pstex_t
\begin{picture}(0,0)%
\includegraphics{Kirby+move.pstex}%
\end{picture}%
\setlength{\unitlength}{3947sp}%
\begingroup\makeatletter\ifx\SetFigFont\undefined
\def\x#1#2#3#4#5#6#7\relax{\def\x{#1#2#3#4#5#6}}%
\expandafter\x\fmtname xxxxxx\relax \def\y{splain}%
\ifx\x\y   
\gdef\SetFigFont#1#2#3{%
  \ifnum #1<17\tiny\else \ifnum #1<20\small\else
  \ifnum #1<24\normalsize\else \ifnum #1<29\large\else
  \ifnum #1<34\Large\else \ifnum #1<41\LARGE\else
     \huge\fi\fi\fi\fi\fi\fi
  \csname #3\endcsname}%
\else
\gdef\SetFigFont#1#2#3{\begingroup
  \count@#1\relax \ifnum 25<\count@\count@25\fi
  \def\x{\endgroup\@setsize\SetFigFont{#2pt}}%
  \expandafter\x
    \csname \romannumeral\the\count@ pt\expandafter\endcsname
    \csname @\romannumeral\the\count@ pt\endcsname
  \csname #3\endcsname}%
\fi
\fi\endgroup
\begin{picture}(7677,3614)(2011,-6110)
\put(8936,-4080){\makebox(0,0)[lb]{\smash{\SetFigFont{12}{14.4}{rm}{\color[rgb]{0,0,0}$T$}%
}}}
\put(2011,-3256){\makebox(0,0)[lb]{\smash{\SetFigFont{12}{14.4}{rm}{\color[rgb]{0,0,0}$P$}%
}}}
\put(3480,-6061){\makebox(0,0)[lb]{\smash{\SetFigFont{12}{14.4}{rm}{\color[rgb]{0,0,0}$L$}%
}}}
\put(7801,-6053){\makebox(0,0)[lb]{\smash{\SetFigFont{12}{14.4}{rm}{\color[rgb]{0,0,0}$L'$}%
}}}
\put(4917,-4077){\makebox(0,0)[lb]{\smash{\SetFigFont{12}{14.4}{rm}{\color[rgb]{0,0,0}$T$}%
}}}
\put(6292,-3448){\makebox(0,0)[lb]{\smash{\SetFigFont{12}{14.4}{rm}{\color[rgb]{0,0,0}$P$}%
}}}
\end{picture}

%% file: bbbb100.pstex_t
\begin{picture}(0,0)%
\includegraphics{bbbb.pstex}%
\end{picture}%
\setlength{\unitlength}{4144sp}%
\begingroup\makeatletter\ifx\SetFigFont\undefined
\def\x#1#2#3#4#5#6#7\relax{\def\x{#1#2#3#4#5#6}}%
\expandafter\x\fmtname xxxxxx\relax \def\y{splain}%
\ifx\x\y   
\gdef\SetFigFont#1#2#3{%
  \ifnum #1<17\tiny\else \ifnum #1<20\small\else
  \ifnum #1<24\normalsize\else \ifnum #1<29\large\else
  \ifnum #1<34\Large\else \ifnum #1<41\LARGE\else
     \huge\fi\fi\fi\fi\fi\fi
  \csname #3\endcsname}%
\else
\gdef\SetFigFont#1#2#3{\begingroup
  \count@#1\relax \ifnum 25<\count@\count@25\fi
  \def\x{\endgroup\@setsize\SetFigFont{#2pt}}%
  \expandafter\x
    \csname \romannumeral\the\count@ pt\expandafter\endcsname
    \csname @\romannumeral\the\count@ pt\endcsname
  \csname #3\endcsname}%
\fi
\fi\endgroup
\begin{picture}(3036,1743)(302,-1409)
\put(1801,-481){\makebox(0,0)[lb]{\smash{\SetFigFont{12}{14.4}{rm}{\color[rgb]{0,0,0}$r$}%
}}}
\end{picture}

%% file: dddd100.pstex_t
\begin{picture}(0,0)%
\includegraphics{dddd.pstex}%
\end{picture}%
\setlength{\unitlength}{4144sp}%
\begingroup\makeatletter\ifx\SetFigFont\undefined
\def\x#1#2#3#4#5#6#7\relax{\def\x{#1#2#3#4#5#6}}%
\expandafter\x\fmtname xxxxxx\relax \def\y{splain}%
\ifx\x\y   
\gdef\SetFigFont#1#2#3{%
  \ifnum #1<17\tiny\else \ifnum #1<20\small\else
  \ifnum #1<24\normalsize\else \ifnum #1<29\large\else
  \ifnum #1<34\Large\else \ifnum #1<41\LARGE\else
     \huge\fi\fi\fi\fi\fi\fi
  \csname #3\endcsname}%
\else
\gdef\SetFigFont#1#2#3{\begingroup
  \count@#1\relax \ifnum 25<\count@\count@25\fi
  \def\x{\endgroup\@setsize\SetFigFont{#2pt}}%
  \expandafter\x
    \csname \romannumeral\the\count@ pt\expandafter\endcsname
    \csname @\romannumeral\the\count@ pt\endcsname
  \csname #3\endcsname}%
\fi
\fi\endgroup
\begin{picture}(6234,3354)(-610,-2269)
\put(2344,-443){\makebox(0,0)[lb]{\smash{\SetFigFont{12}{14.4}{rm}{\color[rgb]{0,0,0}$r$}%
}}}
\put(4186,-1861){\makebox(0,0)[lb]{\smash{\SetFigFont{12}{14.4}{rm}{\color[rgb]{0,0,0}$1$ component}%
}}}
\put(181,929){\makebox(0,0)[lb]{\smash{\SetFigFont{12}{14.4}{rm}{\color[rgb]{0,0,0}$(m+1)$ components}%
}}}
\put(2666,-2211){\makebox(0,0)[lb]{\smash{\SetFigFont{12}{14.4}{rm}{\color[rgb]{0,0,0}$m$ components}%
}}}
\end{picture}

%% file: afigleaf100.pstex_t
\begin{picture}(0,0)%
\includegraphics{afigleaf.pstex}%
\end{picture}%
\setlength{\unitlength}{3947sp}%
\begingroup\makeatletter\ifx\SetFigFont\undefined
\def\x#1#2#3#4#5#6#7\relax{\def\x{#1#2#3#4#5#6}}%
\expandafter\x\fmtname xxxxxx\relax \def\y{splain}%
\ifx\x\y   
\gdef\SetFigFont#1#2#3{%
  \ifnum #1<17\tiny\else \ifnum #1<20\small\else
  \ifnum #1<24\normalsize\else \ifnum #1<29\large\else
  \ifnum #1<34\Large\else \ifnum #1<41\LARGE\else
     \huge\fi\fi\fi\fi\fi\fi
  \csname #3\endcsname}%
\else
\gdef\SetFigFont#1#2#3{\begingroup
  \count@#1\relax \ifnum 25<\count@\count@25\fi
  \def\x{\endgroup\@setsize\SetFigFont{#2pt}}%
  \expandafter\x
    \csname \romannumeral\the\count@ pt\expandafter\endcsname
    \csname @\romannumeral\the\count@ pt\endcsname
  \csname #3\endcsname}%
\fi
\fi\endgroup
\begin{picture}(7320,2002)(-500,-2732)
\put(5906,-2692){\makebox(0,0)[lb]{\smash{\SetFigFont{11}{13.2}{rm}{\color[rgb]{0,0,0}$M'$}%
}}}
\put(-201,-2672){\makebox(0,0)[lb]{\smash{\SetFigFont{11}{13.2}{rm}{\color[rgb]{0,0,0}$L'$}%
}}}
\put(5268,-1479){\makebox(0,0)[lb]{\smash{\SetFigFont{9}{10.8}{rm}{\color[rgb]{0,0,0}$r$}%
}}}
\put(797,-1487){\makebox(0,0)[lb]{\smash{\SetFigFont{9}{10.8}{rm}{\color[rgb]{0,0,0}$r$}%
}}}
\put(2499,-1517){\makebox(0,0)[lb]{\smash{\SetFigFont{9}{10.8}{rm}{\color[rgb]{0,0,0}$r$}%
}}}
\put(3872,-1501){\makebox(0,0)[lb]{\smash{\SetFigFont{9}{10.8}{rm}{\color[rgb]{0,0,0}$r$}%
}}}
\end{picture}

%% file: zebra100.pstex_t
\begin{picture}(0,0)%
\includegraphics{zebra.pstex}%
\end{picture}%
\setlength{\unitlength}{3947sp}%
\begingroup\makeatletter\ifx\SetFigFont\undefined
\def\x#1#2#3#4#5#6#7\relax{\def\x{#1#2#3#4#5#6}}%
\expandafter\x\fmtname xxxxxx\relax \def\y{splain}%
\ifx\x\y   
\gdef\SetFigFont#1#2#3{%
  \ifnum #1<17\tiny\else \ifnum #1<20\small\else
  \ifnum #1<24\normalsize\else \ifnum #1<29\large\else
  \ifnum #1<34\Large\else \ifnum #1<41\LARGE\else
     \huge\fi\fi\fi\fi\fi\fi
  \csname #3\endcsname}%
\else
\gdef\SetFigFont#1#2#3{\begingroup
  \count@#1\relax \ifnum 25<\count@\count@25\fi
  \def\x{\endgroup\@setsize\SetFigFont{#2pt}}%
  \expandafter\x
    \csname \romannumeral\the\count@ pt\expandafter\endcsname
    \csname @\romannumeral\the\count@ pt\endcsname
  \csname #3\endcsname}%
\fi
\fi\endgroup
\begin{picture}(1033,1314)(3480,-2517)
\put(4081,-1335){\makebox(0,0)[lb]{\smash{\SetFigFont{12}{14.4}{rm}{\color[rgb]{0,0,0}$\lambda$}%
}}}
\put(4492,-1791){\makebox(0,0)[lb]{\smash{\SetFigFont{12}{14.4}{rm}{\color[rgb]{0,0,0}$\mu$}%
}}}
\end{picture}

%% file: zebra2100.pstex_t
\begin{picture}(0,0)%
\includegraphics{zebra2.pstex}%
\end{picture}%
\setlength{\unitlength}{3947sp}%
\begingroup\makeatletter\ifx\SetFigFont\undefined
\def\x#1#2#3#4#5#6#7\relax{\def\x{#1#2#3#4#5#6}}%
\expandafter\x\fmtname xxxxxx\relax \def\y{splain}%
\ifx\x\y   
\gdef\SetFigFont#1#2#3{%
  \ifnum #1<17\tiny\else \ifnum #1<20\small\else
  \ifnum #1<24\normalsize\else \ifnum #1<29\large\else
  \ifnum #1<34\Large\else \ifnum #1<41\LARGE\else
     \huge\fi\fi\fi\fi\fi\fi
  \csname #3\endcsname}%
\else
\gdef\SetFigFont#1#2#3{\begingroup
  \count@#1\relax \ifnum 25<\count@\count@25\fi
  \def\x{\endgroup\@setsize\SetFigFont{#2pt}}%
  \expandafter\x
    \csname \romannumeral\the\count@ pt\expandafter\endcsname
    \csname @\romannumeral\the\count@ pt\endcsname
  \csname #3\endcsname}%
\fi
\fi\endgroup
\begin{picture}(802,1060)(2738,-3887)
\put(3013,-2959){\makebox(0,0)[lb]{\smash{\SetFigFont{12}{14.4}{rm}{\color[rgb]{0,0,0}$\lambda$}%
}}}
\end{picture}

%% file: tinytilly100.pstex_t
\begin{picture}(0,0)%
\includegraphics{tinytilly.pstex}%
\end{picture}%
\setlength{\unitlength}{3947sp}%
\begingroup\makeatletter\ifx\SetFigFont\undefined
\def\x#1#2#3#4#5#6#7\relax{\def\x{#1#2#3#4#5#6}}%
\expandafter\x\fmtname xxxxxx\relax \def\y{splain}%
\ifx\x\y   
\gdef\SetFigFont#1#2#3{%
  \ifnum #1<17\tiny\else \ifnum #1<20\small\else
  \ifnum #1<24\normalsize\else \ifnum #1<29\large\else
  \ifnum #1<34\Large\else \ifnum #1<41\LARGE\else
     \huge\fi\fi\fi\fi\fi\fi
  \csname #3\endcsname}%
\else
\gdef\SetFigFont#1#2#3{\begingroup
  \count@#1\relax \ifnum 25<\count@\count@25\fi
  \def\x{\endgroup\@setsize\SetFigFont{#2pt}}%
  \expandafter\x
    \csname \romannumeral\the\count@ pt\expandafter\endcsname
    \csname @\romannumeral\the\count@ pt\endcsname
  \csname #3\endcsname}%
\fi
\fi\endgroup
\begin{picture}(2613,1717)(454,-2050)
\put(601,-2001){\makebox(0,0)[lb]{\smash{\SetFigFont{12}{14.4}{rm}{\color[rgb]{0,0,0}$L'$}%
}}}
\put(2251,-1996){\makebox(0,0)[lb]{\smash{\SetFigFont{12}{14.4}{rm}{\color[rgb]{0,0,0}$L''$}%
}}}
\end{picture}

%% file: baby100.pstex_t
\begin{picture}(0,0)%
\includegraphics{baby.pstex}%
\end{picture}%
\setlength{\unitlength}{3947sp}%
\begingroup\makeatletter\ifx\SetFigFont\undefined
\def\x#1#2#3#4#5#6#7\relax{\def\x{#1#2#3#4#5#6}}%
\expandafter\x\fmtname xxxxxx\relax \def\y{splain}%
\ifx\x\y   
\gdef\SetFigFont#1#2#3{%
  \ifnum #1<17\tiny\else \ifnum #1<20\small\else
  \ifnum #1<24\normalsize\else \ifnum #1<29\large\else
  \ifnum #1<34\Large\else \ifnum #1<41\LARGE\else
     \huge\fi\fi\fi\fi\fi\fi
  \csname #3\endcsname}%
\else
\gdef\SetFigFont#1#2#3{\begingroup
  \count@#1\relax \ifnum 25<\count@\count@25\fi
  \def\x{\endgroup\@setsize\SetFigFont{#2pt}}%
  \expandafter\x
    \csname \romannumeral\the\count@ pt\expandafter\endcsname
    \csname @\romannumeral\the\count@ pt\endcsname
  \csname #3\endcsname}%
\fi
\fi\endgroup
\begin{picture}(732,1986)(1534,-2614)
\put(1838,-1483){\makebox(0,0)[lb]{\smash{\SetFigFont{9}{10.8}{rm}{\color[rgb]{0,0,0}$T$}%
}}}
\put(2101,-2161){\makebox(0,0)[lb]{\smash{\SetFigFont{12}{14.4}{rm}{\color[rgb]{0,0,0}$\lambda$}%
}}}
\put(2101,-811){\makebox(0,0)[lb]{\smash{\SetFigFont{12}{14.4}{rm}{\color[rgb]{0,0,0}$\lambda$}%
}}}
\put(1631,-2556){\makebox(0,0)[lb]{\smash{\SetFigFont{12}{14.4}{rm}{\color[rgb]{0,0,0}$\Gamma(L,\mu)$}%
}}}
\end{picture}

%% file: baby1100.pstex_t
\begin{picture}(0,0)%
\includegraphics{baby1.pstex}%
\end{picture}%
\setlength{\unitlength}{3947sp}%
\begingroup\makeatletter\ifx\SetFigFont\undefined
\def\x#1#2#3#4#5#6#7\relax{\def\x{#1#2#3#4#5#6}}%
\expandafter\x\fmtname xxxxxx\relax \def\y{splain}%
\ifx\x\y   
\gdef\SetFigFont#1#2#3{%
  \ifnum #1<17\tiny\else \ifnum #1<20\small\else
  \ifnum #1<24\normalsize\else \ifnum #1<29\large\else
  \ifnum #1<34\Large\else \ifnum #1<41\LARGE\else
     \huge\fi\fi\fi\fi\fi\fi
  \csname #3\endcsname}%
\else
\gdef\SetFigFont#1#2#3{\begingroup
  \count@#1\relax \ifnum 25<\count@\count@25\fi
  \def\x{\endgroup\@setsize\SetFigFont{#2pt}}%
  \expandafter\x
    \csname \romannumeral\the\count@ pt\expandafter\endcsname
    \csname @\romannumeral\the\count@ pt\endcsname
  \csname #3\endcsname}%
\fi
\fi\endgroup
\begin{picture}(1170,1711)(1651,-2534)
\put(2002,-1521){\makebox(0,0)[lb]{\smash{\SetFigFont{9}{10.8}{rm}{\color[rgb]{0,0,0}$T$}%
}}}
\put(2821,-1206){\makebox(0,0)[lb]{\smash{\SetFigFont{12}{14.4}{rm}{\color[rgb]{0,0,0}$\lambda$}%
}}}
\put(2041,-2476){\makebox(0,0)[lb]{\smash{\SetFigFont{12}{14.4}{rm}{\color[rgb]{0,0,0}$\hat{\Gamma}(L,\mu)$}%
}}}
\end{picture}

%% file: someunknots100.pstex_t
\begin{picture}(0,0)%
\includegraphics{someunknots.pstex}%
\end{picture}%
\setlength{\unitlength}{4144sp}%
\begingroup\makeatletter\ifx\SetFigFont\undefined
\def\x#1#2#3#4#5#6#7\relax{\def\x{#1#2#3#4#5#6}}%
\expandafter\x\fmtname xxxxxx\relax \def\y{splain}%
\ifx\x\y   
\gdef\SetFigFont#1#2#3{%
  \ifnum #1<17\tiny\else \ifnum #1<20\small\else
  \ifnum #1<24\normalsize\else \ifnum #1<29\large\else
  \ifnum #1<34\Large\else \ifnum #1<41\LARGE\else
     \huge\fi\fi\fi\fi\fi\fi
  \csname #3\endcsname}%
\else
\gdef\SetFigFont#1#2#3{\begingroup
  \count@#1\relax \ifnum 25<\count@\count@25\fi
  \def\x{\endgroup\@setsize\SetFigFont{#2pt}}%
  \expandafter\x
    \csname \romannumeral\the\count@ pt\expandafter\endcsname
    \csname @\romannumeral\the\count@ pt\endcsname
  \csname #3\endcsname}%
\fi
\fi\endgroup
\begin{picture}(3409,1046)(590,-1154)
\put(991,-1096){\makebox(0,0)[lb]{\smash{\SetFigFont{12}{14.4}{rm}{\color[rgb]{0,0,0}${\cal{O}}_{+1}$}%
}}}
\put(3151,-1096){\makebox(0,0)[lb]{\smash{\SetFigFont{12}{14.4}{rm}{\color[rgb]{0,0,0}${\cal{O}}_{-1}$}%
}}}
\end{picture}

%% file: braidgen1100.pstex_t
\begin{picture}(0,0)%
\includegraphics{braidgen1.pstex}%
\end{picture}%
\setlength{\unitlength}{3947sp}%
\begingroup\makeatletter\ifx\SetFigFont\undefined
\def\x#1#2#3#4#5#6#7\relax{\def\x{#1#2#3#4#5#6}}%
\expandafter\x\fmtname xxxxxx\relax \def\y{splain}%
\ifx\x\y   
\gdef\SetFigFont#1#2#3{%
  \ifnum #1<17\tiny\else \ifnum #1<20\small\else
  \ifnum #1<24\normalsize\else \ifnum #1<29\large\else
  \ifnum #1<34\Large\else \ifnum #1<41\LARGE\else
     \huge\fi\fi\fi\fi\fi\fi
  \csname #3\endcsname}%
\else
\gdef\SetFigFont#1#2#3{\begingroup
  \count@#1\relax \ifnum 25<\count@\count@25\fi
  \def\x{\endgroup\@setsize\SetFigFont{#2pt}}%
  \expandafter\x
    \csname \romannumeral\the\count@ pt\expandafter\endcsname
    \csname @\romannumeral\the\count@ pt\endcsname
  \csname #3\endcsname}%
\fi
\fi\endgroup
\begin{picture}(2016,991)(1618,-2069)
\put(1816,-2011){\makebox(0,0)[lb]{\smash{\SetFigFont{12}{14.4}{rm}{\color[rgb]{0,0,0}$\sigma^{+1}$}%
}}}
\put(3158,-2011){\makebox(0,0)[lb]{\smash{\SetFigFont{12}{14.4}{rm}{\color[rgb]{0,0,0}$\sigma^{-1}$}%
}}}
\end{picture}

%% file: chapterappendixA100.tex
\begin{chapter}{Gaussian binomial identities and Gaussian sums}
\label{appen:gausssss}

\markboth{\text{Appendix \ref{appen:gausssss}}}
{\text{ }}

\markright{\text{Gaussian binomial coefficient identities and Gaussian sums}}

In this Appendix we give certain identities 
involving the Gaussian binomial coefficients and we also
investigate Gaussian sums.
The pseudo-Gaussian binomial coefficients are closely related to the Gaussian binomial coefficients by
 $(n)_{q} = [n]^{q^{-1}}$, 
but we consider them separately to aid comprehension.

\begin{section}{Gaussian binomial identities}

In this section we give certain identities 
involving the Gaussian and pseudo-Gaussian binomial coefficients.
Lemmas \ref{lem:abcdefghijkl} and \ref{lem:paperlessoffice} are given without proof; 
they are easily proved by
induction.

\begin{lemma}
\label{lem:abcdefghijkl}
The Gaussian binomial coefficients
satisfy the following relations, 
where $i, n \in \mathbb{Z}_{+}$ and $i \leq n$:
\begin{itemize}
\item[(i)]
$\left[ \begin{array}{c}
n+1 \\
i
\end{array} \right]^{q} = \left[ \begin{array}{c}
n \\
i
\end{array} \right]^{q} + q^{n+1-i} \left[ \begin{array}{c}
n \\
i-1
\end{array} \right]^{q}$,
\item[(ii)]
$\left[ \begin{array}{c}
n+1 \\
i
\end{array} \right]^{q} = \left[ \begin{array}{c}
n \\
i-1
\end{array} \right]^{q} + q^{i} \left[ \begin{array}{c}
n \\
i
\end{array} \right]^{q}$.
\end{itemize}
\end{lemma}

\begin{lemma}
\label{lem:paperlessoffice}
The pseudo-Gaussian binomial coefficients 
satisfy the following relations,
where $i, n \in \mathbb{Z}_{+}$ and $i \leq n$:
\begin{itemize}
\item[(i)]
$\left( \begin{array}{c}
n+1 \\
i
\end{array} \right)_{q} = \left( \begin{array}{c}
n \\
i
\end{array} \right)_{q} + (-q)^{n+1-i} \left( \begin{array}{c}
n \\
i-1
\end{array}\right)_{q}$,
\item[(ii)]
$
\left(
\begin{array}{c}
n+1 \\
i
\end{array}
\right)_{q}  =  \left(
\begin{array}{c}
n \\
i-1
\end{array}
\right)_{q} + (-q)^{i} \left(
\begin{array}{c}
n \\
i
\end{array}
\right)_{q}$.
\end{itemize}
\end{lemma}

\end{section}

\begin{section}{Gaussian sums}

In this section we consider some properties of Gaussian sums.
Here we fix $q = \exp{(2 \pi i/N)}$ where $N \geq 3$ is some integer, and we fix 
$G_{+}(N,m)$ and $G_{-}(N,m)$ to mean the following Gaussian sums:
$$G_{+}(N,m) = \sum_{n=0}^{N-1}q^{n(n+m)}, \hspace{10mm} G_{-}(N,m) = \sum_{n=0}^{N-1}q^{-n(n+m)}.$$
\begin{lemma}
\label{lem:firstgausssum}
Let $N \equiv 2 \pmod{4}$ where $N \in \mathbb{Z}_{+}$, then $G_{+}(N,0) = 0$.
\end{lemma}
\begin{proof}
As $q^{(j+N/2)^{2}} = q^{j^{2}+N^{2}/4} = -q^{j^{2}}$, for each $j=0, 1, \ldots, N/2-1$, 
we have $G_{+}(N,0) = \sum_{k=0}^{N-1} q^{k^{2}}=0$.
\end{proof}
\begin{lemma}
\label{lem:secondgausssum}
Let $N \equiv 0 \pmod{4}$ where $N \in \mathbb{Z}_{+}$, then $G_{+}(N,0) = (1+i)\sqrt{N}$.
\end{lemma}
\begin{proof}
See \cite[Sect. 12.8]{kl}.
\end{proof}
\begin{lemma}
\label{lem:thirdgausssum}
Let $m \geq 1$ be an odd integer, $N \equiv 2 \pmod{4}$ where $N \in \mathbb{Z}_{+}$, and
let $t = \exp{(\pi i/2N)}$, then
$$G_{+}(N,m) = \frac{(1+i)\sqrt{N}}{t^{m^{2}}}.$$
\end{lemma}
\begin{proof}
Observe that 
$$\sum_{n=0}^{N-1}q^{n(n+m)} = \frac{1}{2}\sum_{n=0}^{2N-1}q^{n(n+m)},$$ 
and that $q^{n(n+m)} = t^{(2n+m)^{2}}/t^{m^{2}}$, then
\begin{eqnarray*}
G_{+}(N,m) & = &  \frac{1}{2t^{m^{2}}}\sum_{n=0}^{2N-1}t^{(2n+m)^{2}} 
  = \frac{1}{2t^{m^{2}}}\left(t^{m^{2}} + t^{(2+m)^{2}} + t^{(4+m)^{2}} + \cdots +
t^{(4N-2+m)^{2}}   \right) \\
 & = & \frac{1}{2t^{m^{2}}} \left(t^{1^{2}} + t^{3^{2}} + t^{5^{2}} + \cdots +
 t^{(4N-1)^{2}}   \right)  
   =  \frac{1}{2t^{m^{2}}}\big(G_{+}(4N,0) - 2G_{+}(N,0)   \big) \\
 & = & \frac{(1+i)\sqrt{N}}{t^{m^{2}}},
\end{eqnarray*}
which follows from the observation that $G_{+}(N,0)=0$ as $N \equiv 2 \pmod{4}$.
\end{proof}

Lemmas \ref{lemjulie:00} and \ref{lem:bugger} are proved by noting that 
$G_{-}(N,m)$ is the complex conjugate of $G_{+}(N,m)$.
\begin{lemma}
\label{lemjulie:00}
\emph{ }
\begin{itemize}
\item[(i)]  Let $N \equiv 2 \pmod{4}$ where $N \in \mathbb{Z}_{+}$, then $G_{-}(N,0) = 0$,
\item[(ii)] Let $N \equiv 0 \pmod{4}$ where $N \in \mathbb{Z}_{+}$, then $G_{-}(N,0) = (1-i)\sqrt{N}$.
\end{itemize}
\end{lemma}
\begin{lemma}
\label{lem:bugger}
Let $m \geq 1$ be an odd integer, let $N \equiv 2 \pmod{4}$ where $N \in \mathbb{Z}_{+}$, and
let $t = \exp{(\pi i/2N)}$, then $G_{-}(N,m) = t^{m^{2}}(1-i)\sqrt{N}$.
\end{lemma}

\end{section}

\end{chapter}

%% file: chapterappendixB100.tex
\begin{chapter}{The $q$-binomial theorem and generalisations}
\label{chap:appendixB}

\markboth{\text{Appendix \ref{chap:appendixB}}}
{\text{ }}

\markright{\text{The $q$-binomial theorem and generalisations}}

The $q$-binomial theorem is given in 
Lemmas \ref{lem:redrovercomeover} and \ref{lem:noncommute} without proof;
each of these lemmas is easily proved.  
The lemmas are identical in the sense that applying the map $q \mapsto -q$ to 
Lemma \ref{lem:redrovercomeover} gives  Lemma \ref{lem:noncommute}.
We state both lemmas as their properties differ 
when $q$ is a primitive root of unity.
\begin{lemma}
\label{lem:redrovercomeover}
Let $a$ and $b$ be elements of an associative algebra over $\mathbb{C}$  
satisfying the relation $ba=qab$, where $0 \neq q \in \mathbb{C}$.  Then
$\displaystyle{ (a+b)^{n} = \sum_{i=0}^{n} \left[ \begin{array}{c}
n \\
i
\end{array} \right]^{q} a^{i} b^{n-i} }$ for all $n \in \mathbb{N}$.
\end{lemma}

\begin{lemma}
\label{lem:noncommute}
Let $a$ and $b$ be elements of an associative algebra over $\mathbb{C}$  
satisfying the relation $ba=-qab$ where $0 \neq q \in \mathbb{C}$.  Then
$\displaystyle{ (a+b)^{n} = \sum_{i=0}^{n} \left(
\begin{array}{c}
 n \\
 i
\end{array} \right)_{q} a^{i} b^{n-i}}$ for all $n \in \mathbb{N}$.
\end{lemma}

Now we investigate these lemmas when $q$ is a primitive root of unity.
Set $q = \exp{(2 \pi i /N)}$ for some integer $N \geq 3$ and
let $a$ and $b$ be elements of an associative algebra over $\mathbb{C}$ satisfying  
$ab=q^{2}ba$, then we have
$$(a+b)^{N'} = \sum_{i=0}^{N'} \left[  \begin{array}{c}
N' \\
i
\end{array} \right]^{q^{2}} b^{i} a^{N'-i} =  a^{N'} + b^{N'}, $$
as $$\left[  \begin{array}{c}
N' \\
i
\end{array} \right]^{q^{2}} = \left\{ \begin{array} {ll}
1, & i \in \{ 0, N' \}, \\
0, & \mbox{otherwise},
\end{array} \right. $$ 
which follows from the fact that $[i]^{q^{2}}=0$ if and only if $i =kN'$ for some $k \in \mathbb{Z}$.

Now let $a$ and $b$ be elements of an associative algebra over $\mathbb{C}$ satisfying 
$ab=-q ba$, then
$$(a+b)^{\overline{N}} = \sum_{i=0}^{\overline{N}} \left(  \begin{array}{c}
\overline{N} \\
i
\end{array} \right)_{q} b^{i} a^{\overline{N}-i} =  a^{\overline{N}} + b^{\overline{N}}, $$
as 
$$\left(  \begin{array}{c}
\overline{N} \\
i
\end{array} \right)_{q} = \left\{ \begin{array} {ll}
1, & i \in \{ 0, \overline{N} \}, \\
0, & \mbox{otherwise},
\end{array} \right. $$ 
which follows from the fact that $(i)_{q}=0$ if and only if
$i=k \overline{N}$ for some $k \in \mathbb{Z}$.

We give two further generalisations of Pascal's binomial theorem below, 
each of which is a generalisation of the $q$-binomial theorem.  
We are  unaware of these generalisations appearing 
in the literature, so we present them with the relevant proofs.
Note that we obtain the $q$-binomial theorem 
in Lemma \ref{appendixB:generalisationnewtonsbinomial1} if we fix $c=0$, and
we obtain the $q$-multinomial theorem if we artificially fix $\xi=0$
in Lemma \ref{appendixB:generalisationnewtonsbinomial2}. 
\begin{lemma}
\label{appendixB:generalisationnewtonsbinomial1}
Let $a, b$ and $c$ be elements of an associative algebra over $\mathbb{C}$ satisfying 
$$ab = -q ba + c, \hspace{10mm} ac = q^{2} ca, \hspace{10mm} cb = q^{2} bc,$$
where $0 \neq q \in \mathbb{C}$ and $q^{2} \neq 1$, then
$$(a + b)^{n} = \sum_{\stackrel{\alpha,\beta,\gamma \in \mathbb{Z}_{+}}{\alpha + 2\beta + \gamma = n}}  
\frac{ (n)_{q}!  }{ (\alpha)_{q}! (\gamma)_{q}! (2)_{q} (4)_{q} \cdots (2\beta)_{q} }   \ 
b^{\alpha}c^{\beta} a^{\gamma}, \hspace{10mm} n \in \mathbb{N}.$$
\end{lemma}
\begin{proof}
By using the algebra relations we can inductively prove that
$$a^{n}b = (-q)^{n}  b a^{n} + (-q)^{n-1}(n)_{q}  c a^{n-1}, \hspace{10mm} n \in \mathbb{N},$$
which we can use to obtain the following relations, 
where we use $\cdot$ to denote the algebra multiplication 
and  let $\alpha, \beta, \gamma$ be non-negative integers:
\begin{eqnarray*}
b^{\alpha}c^{\beta} a^{\gamma} \cdot  a & = & b^{\alpha}c^{\beta} a^{\gamma+1}, \\
b^{\alpha}c^{\beta} a^{\gamma} \cdot  b & = & 
(-q)^{\gamma + 2\beta} \ b^{\alpha + 1} c^{\beta} a^{\gamma} + 
(-q)^{\gamma-1} (\gamma)_{q} \ b^{\alpha} c^{\beta + 1} a^{\gamma-1}.
\end{eqnarray*}
We can prove that $\alpha + 2\beta + \gamma = n$
if $b^{\alpha}c^{\beta} a^{\gamma}$ is a component in $(a+b)^{n}$, thus we have
\begin{equation}
\label{appendixB:simoneyoung}
(a + b)^{n} = \sum_{\stackrel{\alpha, \beta, \gamma \in \mathbb{Z}_{+}}{\alpha + 2\beta + \gamma = n}}  
\theta(\alpha,\beta,\gamma) \ b^{\alpha}c^{\beta} a^{\gamma}, \hspace{10mm} n \in \mathbb{N}, 
\end{equation}
for some set of coefficients 
$\left\{\theta(\alpha,\beta,\gamma) \in \mathbb{C} | \ \alpha, \beta, \gamma \in \mathbb{Z}_{+} \right\}$.

From (\ref{appendixB:simoneyoung}) and the algebra relations, 
the coefficients $\theta(\alpha,\beta,\gamma)$ satisfy the recursion relation
\begin{equation}
\label{eq:appenB:binom1}
\theta(\alpha,\beta,\gamma)
= \theta(\alpha,\beta,\gamma-1) + (-q)^{\gamma + 2\beta} \theta(\alpha-1,\beta,\gamma)
+ (-q)^{\gamma}(\gamma+1)_{q} \theta(\alpha,\beta-1,\gamma+1)  
\end{equation}
and the boundary conditions $\theta(1,0,0)=\theta(0,0,1)=1$.
In (\ref{eq:appenB:binom1}) we fix $\theta(\alpha,\beta,\gamma) = 0$ 
if any of $\alpha, \beta$ or $\gamma$ are negative.
To complete the proof we just need to show that the
\begin{equation}
\label{appendixB:simoneyoung2}
\theta(\alpha, \beta, \gamma) = 
\frac{ (\alpha + 2\beta + \gamma)_{q}!  }{ (\alpha)_{q}! (\gamma)_{q}! (2)_{q} (4)_{q} \cdots (2\beta)_{q} },
\end{equation}
furnish a solution to the recurrence relation that also satisfy the boundary conditions.  
It is easy to see that the $\theta(\alpha, \beta, \gamma)$ satisfy 
the boundary conditions, and 
substituting them into the right hand side of (\ref{eq:appenB:binom1}) gives
\begin{eqnarray*}
\lefteqn{
\frac{ (\alpha + 2\beta + \gamma - 1)_{q}! }
{ (\alpha)_{q}! (\gamma-1)_{q}! (2)_{q} (4)_{q} \cdots (2\beta)_{q} } 
+ (-q)^{\gamma + 2\beta}  \frac{  (\alpha + 2\beta + \gamma - 1)_{q}! }
{ (\alpha-1)_{q}! (\gamma)_{q}! (2)_{q} (4)_{q} \cdots (2\beta)_{q} } } \\
& & +  (-q)^{\gamma}(\gamma+1)_{q} \frac{  (\alpha + 2\beta + \gamma - 1)_{q}! }
{ (\alpha)_{q}! (\gamma+1)_{q}! (2)_{q} (4)_{q} \cdots (2\beta - 2)_{q} } \\
& = & \frac{ (\alpha + 2\beta + \gamma - 1)_{q}! }
{ (\alpha)_{q}! (\gamma)_{q}! (2)_{q} (4)_{q} \cdots (2\beta)_{q} } 
\left( (\gamma)_{q} + (-q)^{\gamma + 2\beta}(\alpha)_{q} + (-q)^{\gamma}(2 \beta)_{q}  \right)  \\
& = & \frac{ (\alpha + 2\beta + \gamma)_{q}!  }{ (\alpha)_{q}! (\gamma)_{q}! (2)_{q} (4)_{q} \cdots (2\beta)_{q} },
\end{eqnarray*}
as required.
\end{proof}

\begin{lemma}
\label{appendixB:generalisationnewtonsbinomial2}
Let $a, b$ and $c$ be elements of an associative algebra over $\mathbb{C}$ satisfying 
$$ac = q^{2} ca + \xi b^{2}, \hspace{10mm} ab = q^{2} ba, \hspace{10mm} bc = q^{2} cb,$$
where $0 \neq q \in \mathbb{C}$, $q^{2} \neq 1$ and $\xi =  -(1+q)^{2}/(q-q^{-1})$, then
$$(a + b + c)^{n} = \sum_{\stackrel{\alpha,\beta,\gamma \in \mathbb{Z}_{+}}{\alpha + \beta + \gamma = n}} 
\frac{ [n]^{q^{2}}! \ \phi_{\beta} }{ [\alpha]^{q^{2}}! [\beta]^{q^{2}}! [\gamma]^{q^{2}}! } \ 
c^{\alpha} b^{\beta} a^{\gamma}, \hspace{10mm}  n \in \mathbb{N},$$
where $\phi_{\beta} \in \mathbb{C}$ is recursively defined by
$$\phi_{0}=1, \hspace{10mm}  \phi_{1} = 1,  
\hspace{10mm} \phi_{\beta} = \phi_{\beta-1} + \xi [\beta-1]^{q^{2}} \phi_{\beta-2},
\hspace{10mm} \beta \in \mathbb{N} \backslash \{1\}.$$

\end{lemma}
\begin{proof}
By using the algebra relations we can inductively prove that
$$a^{n}c = q^{2n}ca^{n} + \xi q^{2(n-1)}[n]^{q^{2}} b^{2} a^{n-1}, \hspace{10mm} n \in \mathbb{N},$$
and by using this we obtain the following relations,
where we use $\cdot$ to denote the algebra multiplication, 
and let $\alpha, \beta, \gamma$ be non-negative integers:
\begin{eqnarray*}
c^{\alpha} b^{\beta} a^{\gamma} \cdot a & = & c^{\alpha} b^{\beta} a^{\gamma+1} \\
c^{\alpha} b^{\beta} a^{\gamma} \cdot b & = & q^{2\gamma} c^{\alpha} b^{\beta+1} a^{\gamma} \\
c^{\alpha} b^{\beta} a^{\gamma} \cdot c & = & q^{2\gamma + 2\beta} c^{\alpha+1} b^{\beta} a^{\gamma}
+ \xi q^{2(\gamma-1)}[\gamma]^{q^{2}} c^{\alpha} b^{\beta+2} a^{\gamma-1}.
\end{eqnarray*}
We can prove that $\alpha + \beta + \gamma = n$ if $c^{\alpha} b^{\beta} a^{\gamma}$ is a component in 
$(a+b+c)^{n}$,  thus 
\begin{equation}
\label{appendixB:simoneyoung3}
(a+b+c)^{n} = \sum_{\stackrel{\alpha, \beta, \gamma \in \mathbb{Z}_{+}}{\alpha + \beta + \gamma = n}} 
\theta(\alpha,\beta,\gamma) \ c^{\alpha} b^{\beta} a^{\gamma}, \ 
\hspace{10mm} n \in \mathbb{N},
\end{equation}
for some collection of coefficients 
$\{ \theta(\alpha,\beta,\gamma) \in \mathbb{C} | \ \alpha,\beta,\gamma \in \mathbb{Z}_{+} \}$.

From (\ref{appendixB:simoneyoung3}) and the algebra relations,
the coefficients $\theta(\alpha,\beta,\gamma)$ satisfy the recursion relation
\begin{eqnarray}
\theta(\alpha,\beta,\gamma) & = & \theta(\alpha,\beta,\gamma-1) + q^{2\gamma} \theta(\alpha,\beta-1,\gamma)
+ q^{2\gamma + 2\beta} \theta(\alpha-1,\beta,\gamma)  \nonumber  \\ 
& & + \xi q^{2\gamma} [\gamma+1]^{q^{2}} \theta(\alpha,\beta-2,\gamma+1) \label{eq:appenB:binom2}
\end{eqnarray}
and the boundary conditions $\theta(1,0,0)= \theta(0,1,0) = \theta(0,0,1) = 1$.
Here we fix $\theta(\alpha,\beta,\gamma)=0$ if any of $\alpha, \beta, \gamma$ are negative.
To complete the proof all we need do is show that 
\begin{equation}
\label{appendixB:simoneyoung4}
\theta(\alpha,\beta,\gamma) = 
\frac{ [\alpha + \beta + \gamma]^{q^{2}}! \ \phi_{\beta} }
{ [\alpha]^{q^{2}}! [\beta]^{q^{2}}! [\gamma]^{q^{2}}! },
\end{equation}
solves the recurrence relation and satisfies the boundary conditions,
where $\phi_{\beta}$ is itself recursively defined (as stated in the lemma).
Clearly, (\ref{appendixB:simoneyoung4}) satisfies the boundary conditions, 
and substituting (\ref{appendixB:simoneyoung4}) into the right hand side of
(\ref{eq:appenB:binom2}) gives 
\begin{eqnarray}
\lefteqn{  
\frac{ [\alpha+\beta+\gamma-1]^{q^{2}}! \ \phi_{\beta} }{ [\alpha]^{q^{2}}! [\beta]^{q^{2}}! [\gamma-1]^{q^{2}}! }
+ 
q^{2\gamma}\frac{ [\alpha+\beta+\gamma-1]^{q^{2}}! \ \phi_{\beta-1}}{[\alpha]^{q^{2}}![\beta-1]^{q^{2}}![\gamma]^{q^{2}}! }
 } \nonumber \\
 & & + q^{2\gamma + 2\beta} 
\frac{ [\alpha+\beta+\gamma-1]^{q^{2}}! \ \phi_{\beta} }{ [\alpha-1]^{q^{2}}! [\beta]^{q^{2}}! [\gamma]^{q^{2}}! }
 + \xi q^{2\gamma}[\gamma+1]^{q^{2}}
 \frac{ [\alpha+\beta+\gamma-1]^{q^{2}}! \ \phi_{\beta-2} }{ [\alpha]^{q^{2}}! [\beta-2]^{q^{2}}!
 [\gamma+1]^{q^{2}}! } \nonumber \\
 & = & \frac{ [\alpha+\beta+\gamma-1]^{q^{2}}!  }{ [\alpha]^{q^{2}}! [\beta]^{q^{2}}! [\gamma]^{q^{2}}! }
 \left( [\gamma]^{q^{2}} \phi_{\beta} + q^{2\gamma} [\beta]^{q^{2}} \phi_{\beta-1}
 + q^{2\gamma + 2\beta}[\alpha]^{q^{2}} \phi_{\beta} + \xi q^{2\gamma} [\beta]^{q^{2}} [\beta-1]^{q^{2}}
 \phi_{\beta-2}\right) \nonumber \\
 & = & \frac{ [\alpha+\beta+\gamma-1]^{q^{2}}!  }{ [\alpha]^{q^{2}}! [\beta]^{q^{2}}! [\gamma]^{q^{2}}! }
 \left( [\gamma]^{q^{2}} \phi_{\beta} + q^{2\gamma + 2\beta}[\alpha]^{q^{2}}\phi_{\beta}
 + q^{2\gamma} [\beta]^{q^{2}} \left[\phi_{\beta-1} + \xi[\beta-1]^{q^{2}}\phi_{\beta-2}\right]\right).
 \label{appendixB:simoneyoung5}
\end{eqnarray}
By writing $\phi_{\beta} = \phi_{\beta-1} + \xi[\beta-1]^{q^{2}} \phi_{\beta-2}$
for each $\beta \in \mathbb{N} \backslash \{1\}$,
 we can rewrite
(\ref{appendixB:simoneyoung5}) as
$$\frac{ [\alpha+\beta+\gamma-1]^{q^{2}}!  }{ [\alpha]^{q^{2}}! [\beta]^{q^{2}}! [\gamma]^{q^{2}}! }
 \left( [\gamma]^{q^{2}} \phi_{\beta} + q^{2\gamma + 2\beta}[\alpha]^{q^{2}}\phi_{\beta}
 + q^{2\gamma} [\beta]^{q^{2}} \phi_{\beta} \right)
 =\frac{ [\alpha+\beta+\gamma]^{q^{2}}! \ \phi_{\beta}}{ [\alpha]^{q^{2}}! [\beta]^{q^{2}}! [\gamma]^{q^{2}}!},$$
which proves the lemma.
\end{proof}

We obtain an explicit expression for the $\phi_{\beta}$ appearing
in Lemma \ref{appendixB:generalisationnewtonsbinomial2} below.
\begin{lemma}
\label{appendixB:kermitthefrog(1)}
Let $0 \neq q \in \mathbb{C}$ satisfy $q^{2} \neq 1$ and let 
$\phi_{\beta} \in \mathbb{C}$ be recursively defined by
$$\phi_{0}=1, \hspace{10mm}  \phi_{1} = 1,  \hspace{10mm} 
\phi_{\beta} = \phi_{\beta-1} + \xi [\beta-1]^{q^{2}} \phi_{\beta-2},
\hspace{10mm} \beta \in \mathbb{N} \backslash \{1\},$$
where $\xi = -(1+q)^{2}/(q-q^{-1})$.  Then $\phi_{\beta}$ is given by
$$\phi_{0} = 1, \hspace{10mm} \phi_{1}=1, \hspace{10mm}
\phi_{2i}   = (1-q)^{-i} \Psi_{2i},  \hspace{10mm} \phi_{2i+1} = [2i+1]^{q} \phi_{2i},$$
for each $i \in \mathbb{N}$, where
$$\Psi_{2i} = \frac{[4]^{q}}{[2]^{q}} [3]^{q} \frac{[8]^{q}}{[4]^{q}} [5]^{q} \frac{[12]^{q}}{[6]^{q}} [7]^{q}
\cdots [2i-1]^{q} \frac{[4i]^{q}}{[2i]^{q}}.$$
\end{lemma}
\begin{proof}
We firstly calculate $\phi_{2}$:
$$
\phi_{2}= 1 + \xi [1]^{q^{2}}  = (1 + q^{2})/(1-q) = (1-q)^{-1} [4]^{q}/[2]^{q},
$$
thus the claimed solution for $\phi_{\beta}$ is true for $\beta=0,1,2$.
Assume that $\phi_{2i}$ is as given in the lemma for some $i \in \mathbb{N}$, 
then we calculate that
\begin{eqnarray*}
\phi_{2i+1} 
& = & \phi_{2i} - \frac{(1+q)^{2}}{q-q^{-1}}[2i]^{q^{2}}\phi_{2i-1}  \\
& = & [2i-1]^{q}\frac{1+q}{q-q^{-1}} 
             \frac{[4i]^{q}}{[2i]^{q}}\left(-q^{-1} - [2i]^{q} \right)\phi_{2i-2}  \\
& = & (1-q)^{-1} [2i-1]^{q} \frac{[4i]^{q}}{[2i]^{q}} [2i+1]^{q} \phi_{2i-2}  \\
& = & [2i+1]^{q} \phi_{2i},
\end{eqnarray*}
as required, and we have
\begin{eqnarray}
\phi_{2i+2} & = & [2i+1]^{q}\phi_{2i}  - \frac{(1+q)^{2}}{q-q^{-1}}[2i+1]^{q^{2}} \phi_{2i} \nonumber \\
            & = & \phi_{2i} \left(   [2i+1]^{q} - \frac{(1+q)}{q-q^{-1}}[4i+2]^{q}  \right) 
	    \label{appendixB:simoneyoung7} \\
	    & = & \phi_{2i} (q-q^{-1})^{-1} \left( \frac{-q^{-1}(1+q) [2i+1]^{q}[4i+4]^{q}}{  [2i+2]^{q} } 
	    \right) \nonumber \\
	    & = & (1-q)^{-1}  \frac{[2i+1]^{q}[4i+4]^{q}}{[2i+2]^{q}} \phi_{2i}.  \nonumber
\end{eqnarray}
Here we used $(1+q)\left[2i+1\right]^{q^{2}} = \left[4i+2\right]^{q}$ to obtain (\ref{appendixB:simoneyoung7}).

\end{proof}

We now examine the two generalisations of the binomial theorem when $q$ is a primitive root of unity.
Set $q=\exp(2 \pi i/N)$ where $N \geq 3$ is an integer, and
let $a, b$ and $c$ be elements of an associative algebra over $\mathbb{C}$ satisfying 
$$ab = -q ba + c, \hspace{10mm} ac = q^{2} ca, \hspace{10mm} cb = q^{2} bc.$$
Then Lemma \ref{appendixB:generalisationnewtonsbinomial1} implies that
$$(a+b)^{\overline{N}} = \left\{ \begin{array}{ll}
a^{2N} + b^{2N} + (1)_{q}(3)_{q}(5)_{q} \cdots (2N-1)_{q}c^{N}, & N \equiv 1,3 \pmod{4}, \\
a^{N} + b^{N} + (1)_{q}(3)_{q}(5)_{q} \cdots (N-1)_{q}c^{N/2}, & N \equiv 0 \pmod{4}, \\
a^{N/2} + b^{N/2}, & N \equiv 2 \pmod{4}.
\end{array} \right.
$$

Now redefine $a, b$ and $c$ to be 
elements of an associative algebra over $\mathbb{C}$ satisfying 
$$ac = q^{2} ca + \xi b^{2}, \hspace{10mm} ab = q^{2} ba, \hspace{10mm} bc = q^{2} cb,$$
where  $\xi =  -(1+q)^{2}/(q-q^{-1})$. 
Then Lemmas \ref{appendixB:generalisationnewtonsbinomial2}--\ref{appendixB:kermitthefrog(1)} imply that
$$(a+b+c)^{N'} = a^{N'} + \phi_{N'} b^{N'} + c^{N'},$$
where
$$\phi_{N'} = \left\{ \begin{array}{lcll}
 (1-q)^{-(N-1)/2} [N]^{q} \Psi_{N-1}         & =     & 0,       & \mbox{if } N \equiv 1,3 \pmod{4}, \\
(1-q)^{-N/4} \Psi_{N/2}                      & =     & 0,       & \mbox{if } N \equiv 0 \pmod{4}, \\
 (1-q)^{-(N/2-1)/2} [N/2]^{q}  \Psi_{N/2-1}  & \neq  & 0,       & \mbox{if } N \equiv 2 \pmod{4}.
\end{array}  \right. $$
Using this result for $\phi_{N'}$ we have
$$(a+b+c)^{N'} = \left\{ \begin{array}{ll}
a^{N} + c^{N},                              & \mbox{if } N \equiv 1,3 \pmod{4}, \\
a^{N/2} + c^{N/2},                          & \mbox{if } N \equiv 0   \pmod{4}, \\
a^{N/2} + \phi_{N/2}b^{N/2} + c^{N/2},      & \mbox{if } N \equiv 2   \pmod{4}.
\end{array} \right. $$
			       
\end{chapter}

%% file: chapterappendixD100.tex
\begin{chapter}{The Weyl supercharacter formula}
\label{appendixD:title}

\markboth{\text{Appendix \ref{appendixD:title}}}
{\text{ }}

\markright{\text{The Weyl supercharacter formula}}

In this appendix we recall
the Weyl supercharacter of a finite dimensional irreducible 
$U(osp(1|2n))$-module with integral dominant highest weight \cite{k0,k1}.

For each $\Lambda \in H^{*}$ let $D(\Lambda) \subset H^{*}$ be defined by
$$D(\Lambda) = \left\{\Lambda - \sum_{\alpha \in \Phi^{+}} n_{\alpha}\alpha
\bigg| \ n_{\alpha} \in \mathbb{Z}_{+} \right\}.$$  Let $E$ be the algebra of functions on
$H^{*}$ that vanish outside the union of finitely many sets of the form $D(\Lambda)$.
The convolution of two elements $f, g \in E$ is defined by
$f \cdot g (\lambda) = \sum_{\mu \in H^{*}} f(\lambda - \mu)g(\mu)$.
This sum is well-defined as only a finite number of terms in the sum are non-zero.
The algebra $E$ is a commutative algebra with respect to convolution.

Let $e^{\lambda} \in E$ be a function defined by
$$e^{\lambda}(\nu) = \delta_{\lambda \nu}.$$  The convolution 
of two such functions $e^{\lambda}$, $e^{\mu} \in E$ is
$$
e^{\lambda} \cdot e^{\mu}(\nu) 
   = \sum_{\gamma \in H^{*}}e^{\lambda}(\nu - \gamma) e^{\mu}(\gamma)  
   = \sum_{\gamma \in H^{*}} \delta_{\lambda,\nu - \gamma} \delta_{\mu,\gamma} 
   = \sum_{\mu \in H^{*}} \delta_{\lambda, \nu-\mu} 
   = e^{\lambda + \mu}(\nu).
$$
The element $e^{0}$ is the unit of $E$.  
Any function $f \in E$ can be expressed as a sum
$f = \sum_{\lambda \in H^{*}} f(\lambda) e^{\lambda}$.

Let ${\cal{W}}^{\mathfrak{g}}$ be the Weyl group of the Lie (super)algebra $\mathfrak{g}$ and let
$\Phi_{\mathfrak{g}}$ be the set of roots of $\mathfrak{g}$.
The Weyl group ${\cal{W}}^{\mathfrak{g}}$ is generated by the elements 
$\left\{s_{\alpha} | \ \alpha \in \Phi_{\mathfrak{g}} \right\}$, 
where $s_{\alpha}$ is the map $s_{\alpha}: H^{*} \rightarrow H^{*}$ 
defined by 
$$s_{\alpha}(\lambda) = \lambda - \frac{2(\alpha,\lambda)\alpha}{(\alpha,\alpha)},$$
where $( \cdot, \cdot): H^{*} \times H^{*} \rightarrow \mathbb{C}$ 
is the non-degenerate bilinear form defined by $(\epsilon_{i}, \epsilon_{j}) = \delta_{ij}$.

Let $E'$ be the set of rational expressions in the elements of $E$.
By definition, the Weyl supercharacter of a finite dimensional irreducible 
$U(osp(1|2n))$-module 
$V_{\Lambda}$ with integral dominant highest weight $\Lambda \in H^{*}$ is 
$$sch_{\Lambda} = \sum_{\lambda} (-1)^{[\lambda]} m(\lambda) e^{\lambda},$$
where the sum is over all weight spaces of $V_{\Lambda}$, $[\lambda]$ is the grading of the
vectors of $V_{\Lambda}$ in the weight space $\lambda$, 
and $m(\lambda)$ is the multiplicity of the weight space $\lambda$. 
Define a homomorphism $\epsilon':{\cal{W}} \rightarrow \{-1,+1\}$ by: 
$$\epsilon'(\sigma) = \left\{ \begin{array}{rl}
-1, & \mbox{if the
number of reflections in the expression of $\sigma$ with respect to} \\
 & \mbox{the elements of $\overline{\Phi}^{+}_{0}$ is odd}, \\ 
+1, &  \mbox{otherwise}.
 \end{array}  \right. $$
 (Recall that $\overline{\Phi}^{+}_{0}$ is the set of roots
 $\{ \epsilon_{i} \pm \epsilon_{j} | \ 1 \leq i < j \leq n \}$.)
 Then the supercharacter of $V_{\Lambda}$ is \cite{k1}:  
$$sch_{\Lambda} = (-1)^{[\Lambda]}(L')^{-1}
\sum_{\sigma \in {\cal{W}}} \epsilon'(\sigma) e^{\sigma(\Lambda + \rho)},$$ 
where we write $[\Lambda]$ to mean the grading of the highest weight vector of
$V_{\Lambda}$, and
$$L' = \frac{\prod_{\alpha \in \Phi^{+}_{0}}\left(e^{\alpha/2}-e^{-\alpha/2}\right)}
{\prod_{\beta \in \Phi^{+}_{1}}\left(e^{\beta/2}-e^{-\beta/2}\right)}.$$

The expression for $sch_{\Lambda}$ dramatically simplifies for $U(osp(1|2n))$.  
From the root system of $osp(1|2n)$, we have
$$L' = \prod_{\alpha \in \overline{\Phi}^{+}_{0}}
\left(e^{\alpha/2} - e^{-\alpha/2}\right) 
\prod_{\beta \in \Phi^{+}_{1}}\left(e^{\beta/2}+e^{-\beta/2}\right).$$
Now the supercharacter of the trivial (one-dimensional) $U(osp(1|2n))$-module $V_{0}$ is 
$e^{0}$ (as the grading of the highest weight vector of $V_{0}$ is even), thus
$$e^{0} = \frac{\sum_{\sigma \in{\cal{W}}}\epsilon'(\sigma)e^{\sigma(\rho)}}
{\prod_{\alpha \in \overline{\Phi}^{+}_{0}} (e^{\alpha/2} - e^{-\alpha/2}) 
\prod_{\beta \in \Phi^{+}_{1}}(e^{\beta/2}+e^{-\beta/2})},$$
which yields a variant of Weyl's denominator formula for $U(osp(1|2n))$:
$$\sum_{\sigma \in{\cal{W}}}\epsilon'(\sigma)e^{\sigma(\rho)} = 
\prod_{\alpha \in \overline{\Phi}^{+}_{0}} (e^{\alpha/2} - e^{-\alpha/2}) 
\prod_{\beta \in \Phi^{+}_{1}}(e^{\beta/2}+e^{-\beta/2}).$$
We thus obtain the following expression for $sch_{\Lambda}$:
$$sch_{\Lambda} = (-1)^{[\Lambda]}
\frac{\sum_{\sigma \in {\cal{W}}} \epsilon'(\sigma) e^{\sigma(\Lambda + \rho)}}
{\sum_{\sigma \in{\cal{W}}}\epsilon'(\sigma)e^{\sigma(\rho)}}.$$

\end{chapter}

%% file: chapterappendixC100.tex
\begin{chapter}{Hopf ideal of $U_{q}(osp(1|2n))$ at roots of unity}
\label{chap:appendixC}

\markboth{\text{Appendix \ref{chap:appendixC}}}
{\text{ }}

\markright{\text{Hopf ideal of $U_{q}(osp(1|2n))$ at roots of unity}}

In this appendix we prove that the left ideal ${\cal{I}} \subset U_{q}(osp(1|2n))$ 
given in Chapter \ref{chap2A:titlelabel},
where $q = \exp{(2 \pi i /N)}$ and $N \geq 3$ is an integer, is a two-sided Hopf ideal.  
We have not seen the results in this appendix appearing in the literature, and we present them in full,
together with all relevant proofs, for completeness.

The calculations in this appendix are often quite involved.
In particular, the calculations showing 
that ${\cal{I}}$ is a two-sided co-ideal in which we obtain expressions for 
powers of the co-multiplication of each root vector in $U_{q}(osp(1|2n))$, are very intricate. 
We do these particular calculations 
by using the generalisations of the binomial theorem in Appendix \ref{chap:appendixB}.
Note that in this appendix we always fix $q = \exp{(2 \pi i /N)}$ where $N \geq 3$ is an integer,
and we use $\mathfrak{g}$ to denote $osp(1|2n)$.

\begin{section}{Preliminaries}

Recall that the $q$-bracket $[ \cdot, \cdot]_{q}$
is defined for homogeneous $x, y \in  U_{q}(\mathfrak{g})$ 
with weights $wt(x)$, $wt(y)$ respectively, by
$$[x, y]_{q} = xy-(-1)^{[x][y]} q^{(wt(x),wt(y))} yx.$$
We obtain the graded commutator if we formally fix $q=1$ in the $q$-bracket.
The $q$-bracket satisfies the following useful identities:
\begin{eqnarray}
\left[x,yz\right]_{q} & = & 
\left[x,y\right]_{q} z + (-1)^{[x][y]} q^{(wt(x),wt(y))} y\left[x,z\right]_{q}, \label{eq:acer(11)} \\
\left[xy,z\right]_{q} & = & 
x \left[y,z\right]_{q} + (-1)^{[y][z]} q^{(wt(y),wt(z))} \left[x,z\right]_{q}y, \label{eq:acer(22)} \\
\left[x,y^{n}\right]_{q} & = & 
\sum_{i=0}^{n-1} (-1)^{i [x][y]} q^{i(wt(x),wt(y))} y^{i}\left[x,y\right]_{q} y^{n-1-i},
\hspace{5mm} n \in \mathbb{N}, \label{eq:acer(1)} \\
\left[x^{n},y\right]_{q} & = & \sum_{i=0}^{n-1} (-1)^{i [x][y]} q^{i(wt(x),wt(y))} x^{n-1-i} \left[x,y\right]_{q}x^{i},
\hspace{5mm} n \in \mathbb{N}. \label{eq:acer(2)}
\end{eqnarray}
Note that we can obtain (\ref{eq:acer(22)}) (resp. (\ref{eq:acer(2)})) from 
(\ref{eq:acer(11)}) (resp. (\ref{eq:acer(1)}))  
by using the obvious symmetry properties of the $q$-bracket.
We will extensively use Eqs. (\ref{eq:acer(11)})--(\ref{eq:acer(2)}) in this appendix.

\end{section}

\begin{section}{Root vectors in $U_{q}(osp(1|2n))$}

We will extensively use the following theorem due to 
Khoroshkin and Tolstoy in this appendix \cite[Prop. 3.3]{kt}, which we
will refer to as \emph{Khoroshkin and Tolstoy's proposition}.
\begin{theorem}
Let ${\cal{N}}(\phi)$ be a normal order of the elements of $\phi$ and
let the root vectors $e_{\gamma} \in U_{q}(\mathfrak{g})$, $\gamma \in \phi$, be constructed with
respect to ${\cal{N}}(\phi)$ following Subsection \ref{subsect:Rmatricesfrorepresofquantumosp}.
Let $\alpha, \beta \in \phi$ satisfy $\alpha \prec \beta$ with respect to ${\cal{N}}(\phi)$, then
$$\left[e_{\alpha}, e_{\beta}\right]_{q} =  \sum_{\alpha \prec \gamma_{1} \prec \cdots \prec \gamma_{t} \prec \beta }
		C_{(\gamma_{1},k_{1},\ldots,\gamma_{t},k_{t})}
		 (e_{\gamma_{1}})^{k_{1}} (e_{\gamma_{2}})^{k_{2}} \cdots (e_{\gamma_{t}})^{k_{t}},$$
where $\gamma_{1}, \ldots, \gamma_{t} \in \phi$,
$\sum_{i=1}^{t} k_{i} \gamma_{i} = \alpha + \beta$ and the coefficients 
$C_{(\gamma_{1},k_{1},\ldots,\gamma_{t},k_{t})}$ are complex constants.
\end{theorem}
\noindent
An important consequence of this theorem is that 
$\left[e_{\alpha}, e_{\beta}\right]_{q} = 0$ if there does not exist any set of
elements $\gamma_{1}, \ldots, \gamma_{t} \in \phi$ satisfying
$\alpha \prec \gamma_{1} \prec \ldots \prec \gamma_{t} \prec \beta$ and 
$\sum_{i=1}^{t} k_{i} \gamma_{i} = \alpha + \beta$ 
for any set of constants $k_{i} \in \mathbb{N}$. 
A further useful result is that
 $\left[\overline{e}_{\beta}, \overline{e}_{\gamma}\right]_{q}=0$ 
 if $\left[e_{\gamma}, e_{\beta}\right]_{q}=0$.

Henceforth in this Appendix we fix the normal order ${\cal{N}}(\phi)$ and the root vectors in
$U_{q}(\mathfrak{g})$ to be as we defined in 
Subsection \ref{subsect:normalorderofourquantumospatrootsofunity}.

We now prove some very useful identities.
\begin{proposition}
{\emph{ }}
\begin{itemize}
\item[(i)]
For all $1 \leq i < j \leq n$,
$\left[ e_{i}, e_{\alpha_{i+1}+\cdots +\alpha_{j}}\right]_{q} = e_{\alpha_{i}+\cdots +\alpha_{j}}$.
\item[(ii)]
For each $i=1, \ldots, n-2$, 
$\left[e_{i}, e_{\alpha_{i+1}+\cdots +2\alpha_{n}}\right]_{q} = e_{\alpha_{i}+\cdots +2\alpha_{n}}$.
\item[(iii)]
For each $i=1, \ldots, n-2$ and each $j=i+2, \ldots, n-1$,  
$$\left[e_{i}, e_{\alpha_{i+1}+\cdots +2\alpha_{j}+ \cdots+2\alpha_{n}}\right]_{q} 
 = e_{\alpha_{i}+\cdots +2\alpha_{j}+ \cdots+2\alpha_{n}}.$$
\end{itemize}
\end{proposition}
\begin{proof}  
We prove (i).  Firstly, for each $i=1, \ldots, n-2$,
\begin{eqnarray*}
\left[e_{i}, e_{\alpha_{i+1}+\alpha_{i+2}}\right]_{q} 
& = & \left[e_{i}, e_{i+1}\right]_{q} e_{i+2} + q^{-1} e_{i+1}\left[e_{i}, e_{i+2}\right]_{q}
      -q^{-1} \left[e_{i},e_{i+2}\right]_{q}e_{i+1} -q^{-1}e_{i+2}\left[e_{i},e_{i+1}\right]_{q} \\
& = & e_{\alpha_{i}+\alpha_{i+1}}e_{i+2}-q^{-1}e_{i+2}e_{\alpha_{i}+\alpha_{i+1}} 
  =  e_{\alpha_{i}+\alpha_{i+1}+\alpha_{i+2}}.
\end{eqnarray*}
Keeping $i$ fixed, 
assume that $\left[ e_{i}, e_{\alpha_{i+1}+\cdots +\alpha_{j}}\right]_{q} = e_{\alpha_{i}+\cdots +\alpha_{j}}$ for some 
$j = i+2, \ldots, n-1$, then
\begin{eqnarray*}
\left[ e_{i}, e_{\alpha_{i+1}+\cdots +\alpha_{j+1}}\right]_{q}
& = & \left[ e_{i}, e_{\alpha_{i+1}+\cdots +\alpha_{j}}\right]_{q}e_{j+1}
      + q^{-1} e_{\alpha_{i+1}+\cdots +\alpha_{j}}\left[ e_{i},e_{j+1}\right]_{q} \\
& & -q^{-1} \left[ e_{i},e_{j+1}\right]_{q} e_{\alpha_{i+1}+\cdots +\alpha_{j}}
-q^{-1} e_{j+1}\left[ e_{i}, e_{\alpha_{i+1}+\cdots +\alpha_{j}}\right]_{q} \\
& = & e_{\alpha_{i}+\cdots +\alpha_{j}}e_{j+1}-q^{-1} e_{j+1}e_{\alpha_{i}+\cdots +\alpha_{j}} 
  = e_{\alpha_{i}+\cdots +\alpha_{j+1}},
\end{eqnarray*}
as $\left[ e_{i},e_{j+1}\right]_{q}=0$.

We now prove (ii).  A simple calculation shows that 
$\left[e_{n-2}, e_{\alpha_{n-1}+2\alpha_{n}}\right]_{q} = e_{\alpha_{n-2}+ \alpha_{n-1}+2\alpha_{n}}$.
Assume that 
$\left[e_{i}, e_{\alpha_{i+1}+\cdots +2\alpha_{n}}\right]_{q} = e_{\alpha_{i}+\cdots +2\alpha_{n}}$ for some
$i=2, \ldots, n-2$, then
\begin{eqnarray*}
\lefteqn{
\left[e_{i-1}, e_{\alpha_{i}+\cdots +2\alpha_{n}}\right]_{q} } \\
& = & \left[e_{i-1},e_{\alpha_{i}+\cdots +\alpha_{n}}\right]_{q} e_{n} 
+q^{-1}e_{\alpha_{i}+\cdots +\alpha_{n}}\left[e_{i-1},e_{n}\right]_{q} 
 +\left[e_{i-1},e_{n}\right]_{q}e_{\alpha_{i}+\cdots +\alpha_{n}} 
+ e_{n}\left[e_{i-1},e_{\alpha_{i}+\cdots +\alpha_{n}} \right]_{q} \\
& = & e_{\alpha_{i-1}+\cdots +\alpha_{n}}e_{n}+ e_{n}e_{\alpha_{i-1}+\cdots +\alpha_{n}} 
  = e_{\alpha_{i-1}+\cdots +2\alpha_{n}},
\end{eqnarray*}
as $\left[e_{i-1},e_{n}\right]_{q}=0$.

We now prove (iii).  A simple calculation shows that
$\left[e_{i}, e_{\alpha_{i+1}+\cdots +2\alpha_{n-1}+2\alpha_{n}}\right]_{q} = e_{\alpha_{i}+\cdots
+2\alpha_{n-1}+2\alpha_{n}}$
for each $i=1, \ldots, n-3$.
Assume that $\left[e_{i}, e_{\alpha_{i+1}+\cdots +2\alpha_{j}+ \cdots+2\alpha_{n}}\right]_{q} 
= e_{\alpha_{i}+\cdots +2\alpha_{j}+ \cdots+2\alpha_{n}}$ for some $i=1, \ldots, n-3$ and some
$j=i+3, \ldots, n-1$, then
\begin{eqnarray*}
\lefteqn{
\left[e_{i}, e_{\alpha_{i+1}+\cdots +2\alpha_{j-1}+ \cdots+2\alpha_{n}}\right]_{q} } \\
& = & \left[ e_{i},e_{\alpha_{i+1}+\cdots +2\alpha_{j}+ \cdots+2\alpha_{n}}\right]_{q}e_{j-1}
     +q^{-1}e_{\alpha_{i+1}+\cdots +2\alpha_{j}+ \cdots+2\alpha_{n}}\left[ e_{i},e_{j-1}\right]_{q} \\
& & -q^{-1}\left[e_{i}, e_{j-1}\right]_{q}e_{\alpha_{i+1}+\cdots +2\alpha_{j}+ \cdots+2\alpha_{n}}
    -q^{-1} e_{j-1}\left[e_{i},e_{\alpha_{i+1}+\cdots +2\alpha_{j}+ \cdots+2\alpha_{n}}\right]_{q} \\
& = &  e_{\alpha_{i}+\cdots +2\alpha_{j}+ \cdots+2\alpha_{n}}e_{j-1}
        -q^{-1} e_{j-1}e_{\alpha_{i}+\cdots +2\alpha_{j}+ \cdots+2\alpha_{n}} 
  =      e_{\alpha_{i}+\cdots +2\alpha_{j-1}+ \cdots+2\alpha_{n}},
\end{eqnarray*}
as $\left[ e_{i},e_{j-1}\right]_{q}=0$.
\end{proof}

\end{section}

\begin{section}{The left ideal ${\cal{I}} \subset U_{q}(osp(1|2n))$}

From Chapter \ref{chap2A:titlelabel},
${\cal{I}} \subset U_{q}(osp(1|2n))$ is the left ideal generated by the elements of the set
\begin{equation}
\label{appendixC:theidealII}
I = \left\{ (e_{\gamma})^{N'}, (e_{\beta})^{\overline{N}}, 
(\overline{e}_{\gamma})^{N'}, (\overline{e}_{\beta})^{\overline{N}},
(f_{\gamma})^{N'}, (f_{\beta})^{\overline{N}}, 
(\overline{f}_{\gamma})^{N'}, (\overline{f}_{\beta})^{\overline{N}}, 
(J_{i})^{\pm N}-1 | \ 1 \leq i \leq n \right\},
\end{equation}
where $\gamma$ (resp. $\beta$) ranges over all the even (resp. odd) elements of $\phi$.
Recall that the even (resp. odd) elements of $\phi$ are 
$\left\{ \epsilon_{j} \pm \epsilon_{k} | \ 1 \leq j < k \leq n \right\}$ 
(resp. $\left\{ \epsilon_{i} | \ 1 \leq i \leq n \right\}$).
It is convenient to introduce a convention in this section and in
Sections \ref{appendixB:arnie(1)}--\ref{appendixC:arnietermincarrotator(1)}
 that
$\gamma$ (resp. $\beta$) means an even (resp. odd) element of $\phi$, and    
 $\eta$ means any element of $\phi$.

We now define a graded antiautomorphism $\omega: U_{q}(\mathfrak{g}) \rightarrow U_{q}(\mathfrak{g})$ 
to map between  $e_{\eta}$ and $f_{\eta}$.
Slightly generalising the definition of a similar map in \cite{z1}, we define $\omega$ by
$$\omega:  e_{i} \mapsto f_{i}, \hspace{10mm} f_{i} \mapsto e_{i}, 
\hspace{10mm} K_{i}^{\pm 1} \mapsto K_{i}^{\mp 1}, \hspace{10mm} c \mapsto \overline{c},$$
where $\overline{c}$ is the complex conjugate of $c \in \mathbb{C}$.
By definition,
\begin{equation}
\label{appendixB:omegaantiauto}
\omega(xy) = (-1)^{[x][y]} \omega(y) \omega(x), \hspace{10mm} \forall x, y \in U_{q}(\mathfrak{g}),
\end{equation}
and it is easy to see that $[\omega(x)] = [x]$ for each $x \in U_{q}(\mathfrak{g})$.
\begin{proposition}
The graded antiautomorphism $\omega$ is an involution.
\end{proposition}
\begin{proof}
Firstly $\omega^{2}(x)=x$ for each generator $x \in U_{q}(\mathfrak{g})$.  
Now assume that
$\omega^{2}(y)=y$ and $\omega^{2}(z)=z$ for two elements $y$ and $z$ of $U_{q}(\mathfrak{g})$, then
$$\omega^{2}(yz) = (-1)^{[y][z]}\omega \big( \omega(z) \omega(y)\big) 
 = \omega^{2}(y)\omega^{2}(z) = yz, \hspace{10mm} \mbox{and we also have}$$
$$\omega^{2}\left( c_{1} y + c_{2}z\right) = \omega \big(\overline{c}_{1}\omega(y) + \overline{c}_{2} \omega(z) \big) = 
c_{1} y + c_{2}z, \hspace{10mm} \forall c_{1}, c_{2} \in \mathbb{C}.$$
\end{proof}
\begin{proposition}
The graded antiautomorphism $\omega$ commutes with the antipode, ie
$$S \circ \omega = \omega \circ S.$$
\end{proposition}
\begin{proof}
Recall that the action of the antipode on the generators of $U_{q}(\mathfrak{g})$ is
$$S: e_{i} \mapsto -e_{i} K_{i}^{-1}, \hspace{10mm} f_{i} \mapsto -K_{i}f_{i},
\hspace{10mm} K_{i}^{\pm 1} \mapsto K_{i}^{\mp 1},$$
and that $[S(x)] = [x]$ for each $x \in U_{q}(\mathfrak{g})$.
A direct calculation shows that $S \big(\omega(x)\big) = \omega \big(S(x)\big)$ 
for each generator $x \in U_{q}(\mathfrak{g})$.
Now assume that there exist elements $y$ and $z$ in $U_{q}(\mathfrak{g})$ such that
$S \big(\omega(y)\big) = \omega \big(S(y)\big)$ and 
$S \big(\omega(z)\big) = \omega \big(S(z)\big)$, then
an elementary calculation shows that $S \big(\omega(yz)\big) = \omega \big( S(yz) \big)$
and that
$S \big( \omega (c_{1} y + c_{2} z) \big) = \omega \big( S (c_{1} y + c_{2} z) \big)$ 
for all complex constants $c_{1}$ and $c_{2}$.

\end{proof}
\begin{proposition}
The graded antiautomorphism $\omega$ satisfies the relation
\begin{equation}
\label{appendB:propomegadelta}
(\omega \otimes \omega) \circ  \Delta' = \Delta \circ \omega.
\end{equation}
\end{proposition}
\begin{proof}
A direct calculation
shows that (\ref{appendB:propomegadelta}) is true for each generator of $U_{q}(\mathfrak{g})$, and 
the result then follows by a straightforward calculation.
\end{proof}

Recall that if $e_{\mu}$ is defined by $e_{\mu}=\left[e_{\eta}, e_{i}\right]_{q}$, then
$f_{\mu}$ is defined by $f_{\mu}=\left[f_{i},f_{\eta}\right]_{q^{-1}}$.  
The graded antiautomorphism $\omega$ maps between $e_{\mu}$ and $f_{\mu}$ as follows.
\begin{proposition}
\label{sppendixB:steveieg(1)}
For all $1 \leq i < j \leq n$, we have
\begin{itemize}
\item[(i)]
$\omega( e_{\alpha_{i} + \cdots + \alpha_{j}} ) = f_{\alpha_{i} + \cdots + \alpha_{j}}$,
\item[(ii)] $\omega( e_{\alpha_{i} + \cdots + 2\alpha_{j} + \cdots + 2\alpha_{n}} ) 
		= -f_{\alpha_{i} + \cdots + 2\alpha_{j} + \cdots + 2\alpha_{n}}$,
\item[(iii)] $\omega( \overline{e}_{\alpha_{i} + \cdots + \alpha_{j}}) = \overline{f}_{\alpha_{i} + \cdots + \alpha_{j}}$,
\item[(iv)] $\omega( \overline{e}_{\alpha_{i} + \cdots + 2\alpha_{j} + \cdots + 2\alpha_{n}} ) 
		= -\overline{f}_{\alpha_{i} + \cdots + 2\alpha_{j} + \cdots + 2\alpha_{n}}$.
\end{itemize}
\end{proposition}
\begin{proof}
By definition, $\omega(e_{i}) = f_{i}$ for each simple root $\alpha_{i}$.  We now prove (i):
assume that $\omega(e_{\alpha_{i} + \cdots + \alpha_{j}}) = f_{\alpha_{i} + \cdots + \alpha_{j}}$ 
for some $j = i, \ldots, n-1$, then
$$
\omega(e_{\alpha_{i} + \cdots + \alpha_{j+1}})
 =  f_{j+1} f_{\alpha_{i} + \cdots + \alpha_{j}} -q f_{\alpha_{i} + \cdots + \alpha_{j}} f_{j+1} 
  =  f_{\alpha_{i} + \cdots + \alpha_{j+1}}.
$$
We now prove (ii): for each $i=1, \ldots, n-1$, 
$$
\omega(e_{\alpha_{i} + \cdots + 2\alpha_{n}})
  = -\left( f_{n} f_{\alpha_{i} + \cdots + \alpha_{n}} + f_{\alpha_{i} + \cdots + \alpha_{n}} f_{n} \right) 
  = -f_{\alpha_{i} + \cdots + 2\alpha_{n}}.
$$
Now assume that $\omega(e_{\alpha_{i} + \cdots + 2\alpha_{j} + \cdots + 2\alpha_{n}}) = 
-f_{\alpha_{i} + \cdots + 2\alpha_{j} + \cdots + 2\alpha_{n}}$ for some $j=i+2, \ldots, n$, then
\begin{eqnarray*}
\omega(e_{\alpha_{i} + \cdots + 2\alpha_{j-1} + \cdots + 2\alpha_{n}}) 
 & = & -\left(f_{j-1} f_{\alpha_{i} + \cdots + 2\alpha_{j} + \cdots + 2\alpha_{n}} 
           - q f_{\alpha_{i} + \cdots + 2\alpha_{j} + \cdots + 2\alpha_{n}}f_{j-1} \right) \\
 & = & -f_{\alpha_{i} + \cdots + 2\alpha_{j-1} + \cdots + 2\alpha_{n}}.
\end{eqnarray*}
The proof of (iii) (resp. (iv)) is almost identical to the proof of (i) (resp. (ii)).
\end{proof}

\begin{proposition}
\label{sppendixB:steveieg(2)}
The action of the antipode on $e_{\eta}$, for each $\eta \in \phi$, is
$S(e_{\eta}) = c_{\eta} \overline{e}_{\eta} K_{\eta}^{-1}$
for some scalar $c_{\eta} \neq 0$.
\end{proposition}
\begin{proof}
Fixing $\overline{e}_{i} = e_{i}$, 
the proposition is trivially true for all simple roots $\alpha_{i}$ with  
$c_{\alpha_{i}}=-1$.  
Assume now that the proposition is true for some $\eta \in \phi$ and 
define $e_{\mu} = [e_{\eta},e_{i}]_{q}$ where $\mu = \eta + \alpha_{i} \in \phi$, then
$\overline{e}_{\mu} = [e_{i}, \overline{e}_{\eta}]_{q}$, and 
\begin{eqnarray*}
S(e_{\mu}) & = & S(e_{\eta}e_{i}) - (-1)^{[e_{\eta}][e_{i}]}q^{(\eta,\alpha_{i})} S(e_{i}e_{\eta}) \\
           & = & -(-1)^{[e_{\eta}][e_{i}]} q^{-(\eta,\alpha_{i})} c_{\eta}
	   \left( e_{i} \overline{e}_{\eta} - (-1)^{[e_{\eta}][e_{i}]} q^{(\eta,\alpha_{i})} \overline{e}_{\eta} e_{i}
	   \right)K_{\mu}^{-1} \\
	   & = & c_{\mu} \overline{e}_{\mu} K_{\mu}^{-1},
\end{eqnarray*}
where $c_{\mu} = -(-1)^{[e_{\eta}][e_{i}]} q^{-(\eta,\alpha_{i})} c_{\eta}$.
This formula determines $c_{\mu}$ recursively.
\end{proof}
\begin{proposition}
\label{sppendixB:steveieg(3)}
The action of the antipode on $\overline{e}_{\eta}$, for each $\eta \in \phi$, is
$S(\overline{e}_{\eta}) = d_{\eta} e_{\eta} K_{\eta}^{-1}$
for some scalar $d_{\eta} \neq 0$.
\end{proposition}
\begin{proof}
The proposition is trivially true for all simple roots $\alpha_{i}$
with $d_{\alpha_{i}}=-1$.
Assume now that the proposition is true for some $\eta \in \phi$ and 
define $\overline{e}_{\mu} = [e_{i},\overline{e}_{\eta}]_{q}$ where $\mu = \eta + \alpha_{i} \in \phi$,  
then $e_{\mu} = \left[e_{\eta}, e_{i} \right]_{q}$, and
\begin{eqnarray*}
S(\overline{e}_{\mu}) & = & S(e_{i} \overline{e}_{\eta}) 
                            - (-1)^{[e_{\eta}][e_{i}]}q^{(\eta,\alpha_{i})}S(\overline{e}_{\eta} e_{i})  \\
		      & = & -(-1)^{[e_{i}][\overline{e}_{\eta}]} q^{-(\eta, \alpha_{i})} d_{\eta} 
		            \left( e_{\eta} e_{i} - (-1)^{[e_{i}][e_{\eta}]} q^{(\eta, \alpha_{i})}e_{i}e_{\eta}\right)
			    K_{\mu}^{-1} \\
		      & = & d_{\mu} e_{\mu} K_{\mu}^{-1},
\end{eqnarray*}
where $d_{\mu} = -(-1)^{[e_{\eta}][e_{i}]} q^{-(\eta,\alpha_{i})} d_{\eta}$.
This formula determines $d_{\mu}$ recursively.

\end{proof}

\begin{proposition}
\label{sppendixB:steveieg(4)}
The action of the antipode on $f_{\eta}$, for each $\eta \in \phi$, is
$S(f_{\eta}) = \overline{c}_{\eta} K_{\eta} \overline{f}_{\eta}$
for some scalar $\overline{c}_{\eta} \neq 0$.
\end{proposition}
\begin{proof}  
As the antipode commutes with the graded antiautomorphism $\omega$,
for each $1 \leq i < j \leq n$ we have
$$
S(f_{\alpha_{i} + \cdots + \alpha_{j}}) = S\big(\omega(e_{\alpha_{i} + \cdots + \alpha_{j}})\big) 
    = \omega \big( c_{\alpha_{i} + \cdots + \alpha_{j}} \overline{e}_{\alpha_{i} + \cdots + \alpha_{j}} 
            K_{\alpha_{i} + \cdots + \alpha_{j}}^{-1} \big) 
    = \overline{c}_{\alpha_{i} + \cdots + \alpha_{j}} 
  K_{\alpha_{i} + \cdots + \alpha_{j}}\overline{f}_{\alpha_{i} + \cdots + \alpha_{j}},
$$
where $\overline{c}_{\alpha_{i} + \cdots + \alpha_{j}} \neq 0$ is a scalar,
and we have used Proposition \ref{sppendixB:steveieg(2)}.
An almost identical calculation proves the corresponding result when 
$\eta = \alpha_{i} + \cdots + \alpha_{j-1} + 2\alpha_{j} + \cdots + 2\alpha_{n}$. 
\end{proof}

\begin{proposition}
\label{sppendixB:steveieg(5)}
The action of the antipode on $\overline{f}_{\eta}$, for each $\eta \in \phi$, is
$S(\overline{f}_{\eta}) = \overline{d}_{\eta} K_{\eta} f_{\eta}$
for some scalar $\overline{d}_{\eta} \neq 0$.
\end{proposition}
\begin{proof}
The fact that the antipode commutes with the graded antiautomorphism $\omega$
means that for each $1 \leq i < j \leq n$ we have
$$
S(\overline{f}_{\alpha_{i} + \cdots + \alpha_{j}}) 
  = S\big(\omega(\overline{e}_{\alpha_{i} + \cdots + \alpha_{j}})\big) 
  = \omega \big(d_{\alpha_{i} + \cdots + \alpha_{j}} e_{\alpha_{i} + \cdots + \alpha_{j}} 
      K_{\alpha_{i} + \cdots + \alpha_{j}}^{-1} \big) 
  = \overline{d}_{\alpha_{i} + \cdots + \alpha_{j}}K_{\alpha_{i} + \cdots + \alpha_{j}}f_{\alpha_{i} + \cdots + \alpha_{j}},
$$
where $\overline{d}_{\alpha_{i} + \cdots + \alpha_{j}} \neq 0$ is a scalar, 
and we have used Proposition \ref{sppendixB:steveieg(3)}.
An almost identical calculation proves the corresponding result when 
$\eta = \alpha_{i} + \cdots + \alpha_{j-1} + 2\alpha_{j} + \cdots + 2\alpha_{n}$. 
\end{proof}

\end{section}

\begin{section}{The left ideal ${\cal{I}}$ is a two-sided ideal}
\label{appendixB:arnie(1)}

In this section we prove that ${\cal{I}}$ is a two-sided ideal of $U_{q}(osp(1|2n))$.
To do this, we show that each element of the set $I$ from (\ref{appendixC:theidealII})
commutes or anticommutes with each generator of
$U_{q}(\mathfrak{g})$.
We firstly show that $(e_{\gamma})^{N'}$  and $(e_{\beta})^{\overline{N}}$ have this property for each
$\gamma, \beta \in \phi$, then it is not difficult to show that 
$(\overline{e}_{\gamma})^{N'}, (\overline{e}_{\beta})^{\overline{N}}, (f_{\gamma})^{N'}, (f_{\beta})^{\overline{N}},
(\overline{f}_{\gamma})^{N'}$ and $(\overline{f}_{\beta})^{\overline{N}}$ 
also have this property using 
the antipode and  $\omega$.

Trivially,
$(K_{i})^{\pm N}$ is central in $U_{q}(\mathfrak{g})$ and thus so is $(J_{i})^{\pm N}-1$.
We now show that $(e_{i})^{N'}$ and $(e_{n})^{\overline{N}}$  commute or anticommute with each generator 
of $U_{q}(\mathfrak{g})$ for each $i= 1, \ldots, n-1$.

\begin{subsection}{$(e_{i})^{N'}$ and $(e_{n})^{\overline{N}}$}
Set $1 \leq i \leq n-1$ and $1 \leq j \leq n$, then 
trivially $(e_{i})^{N'}$ and $K_{j}^{\pm 1}$ (anti)commute.
We now show that $(e_{i})^{N'}f_{j} = f_{j}(e_{i})^{N'}$.
The $U_{q}(\mathfrak{g})$ 
relations state that $e_{j}$ and $f_{i}$ commute for all $j \neq i$, thus we need only consider the relations
between $e_{i}$ and $f_{i}$.
The quadruple $\left\{e_{i},f_{i},K_{i}^{\pm 1}\right\}$ generates a 
$U_{q}(sl_{2})$ subalgebra of $U_{q}(osp(1|2n))$, and for each $t \in \mathbb{N}$,
\begin{eqnarray*}
(e_{i})^{t}f_{i} & = & f_{i}(e_{i})^{t} \\
& & + \frac{1}{q-q^{-1}}
\bigg[\big(1+q^{-2}+\cdots+q^{-2(t-1)}\big)K_{i}   
- \big(1+q^{2}+\cdots+q^{2(t-1)}\big)K_{i}^{-1} \bigg] (e_{i})^{t-1} \\
& & = f_{i}(e_{i})^{t} + [t]^{q^{2}} \left(\frac{q^{-2t+2}K_{i}-K_{i}^{-1}}{q-q^{-1}}\right) (e_{i})^{t-1}.
\end{eqnarray*}
As $[N']^{q^{2}}=0$, $(e_{i})^{N'}f_{i} = f_{i}(e_{i})^{N'}$.

We now compute the relations between $(e_{i})^{N'}$ and $e_{j}$ for all $i \neq j$.
If $|i-j|>1$, the relevant Serre relation tells us that $e_{i}$ and
$e_{j}$ commute, so assume that $j=i \pm 1$.  The Serre relation is
\begin{equation}
\label{eq:sweet}
(e_{i})^{2}e_{i \pm 1}-(q+q^{-1})e_{i}e_{i \pm1}e_{i} + e_{i \pm 1}(e_{i})^{2}=0.
\end{equation}
By repeatedly using (\ref{eq:sweet}), we have
\begin{eqnarray*}
(e_{i})^{t} e_{i \pm 1} & = & (q+q^{-1})(e_{i})^{t-1}e_{i \pm 1}e_{i}-(e_{i})^{t-2}e_{i \pm 1}(e_{i})^{2} \\
                        & = & \left[(q+q^{-1})^{2}-1\right](e_{i})^{t-2}e_{i \pm 1}(e_{i})^{2}-
			      (q+q^{-1})(e_{i})^{t-3}e_{i \pm 1}(e_{i})^{3}  \\
                        & = & a_{r}(e_{i})^{t-1-r}e_{i \pm 1}(e_{i})^{1+r} + b_{r}(e_{i})^{t-2-r}e_{i \pm 1}(e_{i})^{2+r},
\end{eqnarray*}
for all $0 \leq r \leq t-2$, where $a_{r}, b_{r} \in \mathbb{C}$ satisfy the recurrence relations
\begin{eqnarray*}
a_{r+1} & = & (q+q^{-1})a_{r}+b_{r}, \\
b_{r+1} & = & -a_{r},
\end{eqnarray*}
where $a_{0} = q + q^{-1}$ and $b_{0} = -1$.
The solution for $a_{r}$ and $b_{r}$ is
$$\left( \begin{array}{c}
a_{r} \\
b_{r} \end{array} \right)
= 
\frac{q^{-1-r}}{q^{-1}-q}\left(
\begin{array}{c}
q^{-1} \\
-1
\end{array}
\right) + \frac{q^{1+r}}{q-q^{-1}} \left(
\begin{array}{c}
q \\
-1
\end{array}
\right),$$
and by setting $r=t-2$, we obtain
$$
(e_{i})^{t}e_{i \pm 1} = \left( \frac{q^{t}-q^{-t}}{q-q^{-1}} \right) e_{i}e_{i \pm 1}(e_{i})^{t-1}
+ \left( \frac{-q^{t-1}+q^{1-t}}{q-q^{-1}} \right) e_{i \pm 1}(e_{i})^{t}.
$$
Setting $t=N'$, we have
$$\begin{array}{rcll}
(e_{i})^{N}e_{i \pm 1}   & = & e_{i \pm 1}(e_{i})^{N},    & \mbox{if $N$ is odd}, \\
(e_{i})^{N/2}e_{i \pm 1} & = & -e_{i \pm 1}(e_{i})^{N/2}, & \mbox{if $N$ is even}.
\end{array} $$

We now consider $(e_{n})^{\overline{N}}$: fix $1 \leq i \leq n$.  
Trivially, $(e_{n})^{\overline{N}}$ and $K_{i}^{\pm 1}$ (anti)commute.
The $U_{q}(\mathfrak{g})$ 
relations state that $e_{n}$ and $f_{i}$ commute for all $i < n$, thus we consider the relations between
$e_{n}$ and $f_{n}$.  The quadruple $\left\{e_{n},f_{n},K_{n}^{\pm 1}\right\}$
generates a $U_{q}(osp(1|2))$ subalgebra of $U_{q}(osp(1|2n))$, and for each $t \in \mathbb{N}$,
\begin{eqnarray*}
(e_{n})^{t}f_{n} & = & (-1)^{t} f (e_{n})^{t} \\
& & +    (e_{n})^{t-1} \frac{1}{q-q^{-1}}
		    \Big[ (1-q+\cdots + (-q)^{n-1})K_{n}  \\
& & \hspace{35mm} -(1-q^{-1}+ \cdots + (-q)^{-(n-1)}) K_{n}^{-1}\Big] \\
& = & (-1)^{t} f (e_{n})^{t} + (e_{n})^{t-1} (n)_{q} \left( \frac{K_{n} - (-q)^{1-n}K_{n}^{-1}}{q-q^{-1}} \right).
\end{eqnarray*}
As $(\overline{N})_{q}=0$, $(e_{n})^{\overline{N}}f_{n} = (-1)^{\overline{N}} f_{n} (e_{n})^{\overline{N}}$.

We now compute the relations between $(e_{n})^{\overline{N}}$ and $e_{i}$ for all $i < n$.
If $n-i>1$, the relevant Serre relation tells us that $e_{n}$ and $e_{i}$ commute, so
fix $i=n-1$.  The Serre relation is
\begin{equation}
\label{eq:adjoint3}
(e_{n})^{3}e_{n-1}-(q-1+q^{-1})(e_{n})^{2}e_{n-1}e_{n} -(q-1+q^{-1})e_{n}e_{n-1}(e_{n})^{2}+e_{n-1}(e_{n})^{3}=0.
\end{equation}
By repeatedly using (\ref{eq:adjoint3}), we have
\begin{equation}
\label{eq:pigsarmpit}
(e_{n})^{t} e_{n-1}=
a_{r}(e_{n})^{t-1-r}e_{n-1}(e_{n})^{1+r} + 
b_{r}(e_{n})^{t-2-r}e_{n-1}(e_{n})^{2+r} + 
c_{r}(e_{n})^{t-3-r}e_{n-1}(e_{n})^{3+r}, 
\end{equation}
for each $0 \leq r \leq t-3$, where $a_{r}, b_{r}, c_{r} \in \mathbb{C}$ satisfy the
recurrence relations
\begin{eqnarray*}
a_{r} & = & (q-1+q^{-1})a_{r-1} + b_{r-1}, \\
b_{r} & = & (q-1+q^{-1})a_{r-1} + c_{r-1}, \\
c_{r} & = & -a_{r-1},
\end{eqnarray*} 
where $a_{0} = b_{0} = (q-1+q^{-1})$ and $c_{0} = -1$. 
The solution for $a_{r}$, $b_{r}$, $c_{r}$ is
$$\left( \begin{array}{c}
a_{r} \\
b_{r} \\
c_{r}
\end{array}
\right) = 
 d_{1} q^{-r} \left(
\begin{array}{c}
-q^{-1} \\
1-q^{-1} \\
1 \\
\end{array}
\right) + 
d_{2} (-1)^{r} \left(
\begin{array}{c}
1 \\
-q-q^{-1} \\
1 \\
\end{array}
\right) + 
d_{3} q^{r}\left(
\begin{array}{c}
q \\
q-1 \\
-1 \\
\end{array}
\right),$$
where $$d_{1} = \frac{1}{(q^{2}-1)(q+1)}, \  \ \ d_{2} = \frac{-q}{(q+1)^{2}}, \ \ \ 
d_{3} = \frac{q^{3}}{(q^{2}-1)(q+1)}.$$
By setting $r=t-3$, we obtain
\begin{eqnarray*}
a_{t-3} & = & \frac{-q^{2-t}(1+q) + (-1)^{t}q(q^{2}-1) + q^{t+1}(q+1)}{(q^{2}-1)(1+q)^{2}}, \\
b_{t-3} & = & \frac{ q^{3-t}(1-q^{-1})(1+q) + (-1)^{t-1}q(q+q^{-1})(q^{2}-1) + q^{t}(q-1)(q+1) }{(q^{2}-1)(1+q)^{2}}, \\
c_{t-3} & = & \frac{q^{3-t}(1+q) + (-1)^{t}q(q^{2}-1) - q^{t}(1+q)}{(q^{2}-1)(1+q)^{2}}.
\end{eqnarray*}
Fixing $t=\overline{N}$, we have $a_{\overline{N}-3}=0$, $b_{\overline{N}-3}=0$ and 
$c_{\overline{N}-3}= (-1)^{\overline{N}}$,  thus
$$(e_{n})^{\overline{N}} e_{n-1} = (-1)^{\overline{N}} e_{n-1}(e_{n})^{\overline{N}}.$$

\end{subsection}

\begin{subsection}{Relations between $(e_{\gamma})^{N'}$, $(e_{\beta})^{\overline{N}}$ and $e_{i}$}

We now prove that $(e_{\gamma})^{N'}$ and $(e_{\beta})^{\overline{N}}$ 
(anti)commute with each generator of 
$U_{q}(\mathfrak{g})$ for each non-simple root $\gamma, \beta \in \phi$.  
These calculations are simplified by noting that if 
we can show that $[e_{\gamma}, e_{i}]=0$ or that
$[e_{\gamma}, e_{i}]_{q}=0$, then $(e_{\gamma})^{N'}$ automatically (anti)commutes with $e_{i}$, 
and if we can show that $[e_{\beta}, e_{i}]=0$ or that
$[e_{\beta}, e_{i}]_{q}=0$, then $(e_{\beta})^{\overline{N}}$ 
automatically (anti)commutes with $e_{i}$.
In the following we write $\mu$ to mean any element of $\phi$.

We will determine the relations between $e_{i}$ and $(e_{\gamma})^{N'}$, $(e_{\beta})^{\overline{N}}$,
 for a fixed $i$, by breaking the problem into 4 sub-problems:
\begin{itemize}
\item[(i)] $\alpha_{i+2} \prec \mu$, 
\item[(ii)] $\alpha_{i} \prec \mu \prec \alpha_{i+1}$,
\item[(iii)] $\alpha_{i+1} \prec \mu \prec \alpha_{i+2}$, 
\item[(iv)] $\mu \prec \alpha_{i}$.
\end{itemize}

\begin{subsubsection}{Case 1: $\alpha_{i+2} \prec \mu$.}
Here $\mu =\alpha_{i+2} + \cdots + \alpha_{j}$ or  
$\mu = \alpha_{i+2} + \cdots + \alpha_{j-1} + 2\alpha_{j} + \cdots + 2\alpha_{n}$ for some $j=i+3, \ldots, n$.
From Khoroshkin and Tolstoy's proposition,
$\left[e_{i}, e_{\mu}\right]_{q}=0$,
which also follows from the following Serre relation: $e_{i}e_{j}=e_{j}e_{i}$ if $|i-j|>1$.
\end{subsubsection}

\begin{subsubsection}{Case 2:  $\alpha_{i} \prec \mu \prec \alpha_{i+1}$.}  
Here $\mu = \alpha_{i} + \cdots + \alpha_{j}$ or  
$\mu = \alpha_{i} + \cdots +  \alpha_{j-1} + 2\alpha_{j} + \cdots + 2\alpha_{n}$ for some $j=i+1, \ldots, n$.
From Khoroshkin and Tolstoy's proposition, we have
\begin{eqnarray*}
\left[ e_{i}, e_{\alpha_{i} + \cdots + \alpha_{j}} \right]_{q} & = & 0, \hspace{10mm} i<j \leq n, 
                      \\
\left[ e_{i}, e_{\alpha_{i} + \cdots + 2\alpha_{j} + \cdots + 2\alpha_{n}} \right]_{q} & = & 0,
                           \hspace{10mm} i+2 \leq j \leq n, \hspace{10mm} \mbox{and}
\end{eqnarray*}
\begin{eqnarray*}
\lefteqn{
\left[ e_{i}, e_{\alpha_{i} + 2\alpha_{i+1} + \cdots + 2\alpha_{n}} \right]_{q} } \\
& & =
\sum_{k=1}^{n-i-1} C_{k }e_{\alpha_{i} + \cdots + \alpha_{i+k}} 
                   e_{\alpha_{i} + \cdots +  2\alpha_{i+k+1} + \cdots + 2\alpha_{n}}
		   + C_{n-i} \left( e_{\alpha_{i} + \cdots + \alpha_{n}}\right)^{2}, \ C_{k} \in \mathbb{C}.
\end{eqnarray*}
The first two identities dispose of much of this case.  
To deal with the remaining problem, we claim that the right hand 
side of Eq. (\ref{appendixB:julietassone(3)}) below vanishes identically:
\begin{eqnarray}
\lefteqn{
\left[ e_{i}, \left(e_{\alpha_{i} + 2\alpha_{i+1} + \cdots + 2\alpha_{n}}\right)^{N'} \right]_{q} } \nonumber \\
& = & 
\sum_{m=0}^{N'-1} \left(e_{\alpha_{i} + 2\alpha_{i+1} + \cdots + 2\alpha_{n}}\right)^{m}
\left[ e_{i}, e_{\alpha_{i} + 2\alpha_{i+1} + \cdots + 2\alpha_{n}} \right]_{q}
\left(e_{\alpha_{i} + 2\alpha_{i+1} + \cdots + 2\alpha_{n}}\right)^{N'-1-m}. \label{appendixB:julietassone(3)}
\end{eqnarray}
To show this, we use the identities:
$$\begin{array}{rcll}
\left[e_{\alpha_{i} + \cdots + \alpha_{k}}, e_{\alpha_{i} + 2\alpha_{i+1} + \cdots + 2\alpha_{n}} \right]_{q} & = & 0, &
     i+1 \leq k \leq n, \\
\left[e_{\alpha_{i} +  \cdots + 2\alpha_{k} + \cdots + 2\alpha_{n}}, 
e_{\alpha_{i} + 2\alpha_{i+1} + \cdots + 2\alpha_{n}} \right]_{q} & = & 0, &
   i+2 \leq k \leq n,
   \end{array} $$
which can be rewritten, respectively, as
\begin{eqnarray*}
e_{\alpha_{i} + \cdots + \alpha_{k}}e_{\alpha_{i} + 2\alpha_{i+1} + \cdots + 2\alpha_{n}}
 & = & q e_{\alpha_{i} + 2\alpha_{i+1} + \cdots + 2\alpha_{n}}e_{\alpha_{i} + \cdots + \alpha_{k}}, \\
e_{\alpha_{i} +  \cdots + 2\alpha_{k} + \cdots + 2\alpha_{n}}e_{\alpha_{i} + 2\alpha_{i+1} + \cdots + 2\alpha_{n}}
& = & q e_{\alpha_{i} + 2\alpha_{i+1} + \cdots + 2\alpha_{n}}e_{\alpha_{i} +  \cdots + 2\alpha_{k} +
\cdots + 2\alpha_{n}}.
\end{eqnarray*}
Using these, we rewrite the right hand side of (\ref{appendixB:julietassone(3)}) as
$$\sum_{m=0}^{N'-1} q^{-2m} \left[ e_{i}, e_{\alpha_{i} + 2\alpha_{i+1} + \cdots + 2\alpha_{n}} \right]_{q}
\left(e_{\alpha_{i} + 2\alpha_{i+1} + \cdots + 2\alpha_{n}}\right)^{N'-1}=0,$$
 as $\sum_{m=0}^{N'-1} q^{-2m}=0$. Thus 
$\left(e_{\alpha_{i} + 2\alpha_{i+1} + \cdots + 2\alpha_{n}}\right)^{N'}$ (anti)commutes with $e_{i}$.
\end{subsubsection}

\begin{subsubsection}{Case 3: $\alpha_{i+1} \prec \mu \prec \alpha_{i+2}$.}
This is the most difficult of the four cases. The identity:
$\left[e_{\alpha_{i}+\cdots+\alpha_{j}}, e_{\alpha_{i+1}+\cdots+\alpha_{j}}\right]_{q} = 0$, for 
$j < n$,
will be useful here.  Note that we can rewrite this identity as
$e_{\alpha_{i}+\cdots+\alpha_{j}}e_{\alpha_{i+1}+\cdots+\alpha_{j}}
   =q  e_{\alpha_{i+1}+\cdots+\alpha_{j}}e_{\alpha_{i}+\cdots+\alpha_{j}}$.
To prove this case we break it into a number of sub-cases, 
each of which we consider in the following proposition.
\begin{proposition}  We have
\begin{itemize}
\item[(i)] $\left[  e_{i}, \left(e_{\alpha_{i+1}+\cdots+\alpha_{j}} \right)^{N'} \right]_{q} = 0$, 
             for each $j=i+2, \ldots, n-1$,
\item[(ii)] $\left[  e_{i}, \left(e_{\alpha_{i+1}+\cdots+\alpha_{n}} \right)^{\overline{N}} \ \right]_{q} = 0$,
\item[(iii)] $\left[  e_{i}, \left(e_{\alpha_{i+1}+\cdots+2\alpha_{j}+\cdots+2\alpha_{n}}\right)^{N'} \right]_{q} = 0$,
              for each $j=i+1, \ldots, n$.
\end{itemize}
\end{proposition}
\begin{proof}
We use Eqs. (\ref{eq:acer(1)})--(\ref{eq:acer(2)}) in this proof.
We firstly prove (i):
$$
\left[  e_{i}, \left(e_{\alpha_{i+1}+\cdots+\alpha_{j}} \right)^{N'} \right]_{q} 
= \sum_{k=0}^{N'-1} q^{-2k} e_{\alpha_{i}+\cdots+\alpha_{j}}
      \left(e_{\alpha_{i+1}+\cdots+\alpha_{j}} \right)^{N'-1} 
  = 0.
$$
We prove (ii): 
\begin{eqnarray}
\lefteqn{
\left[  e_{i}, \left(e_{\alpha_{i+1}+\cdots+\alpha_{n}} \right)^{\overline{N}} \ \right]_{q} } \nonumber \\
& = & \sum_{k=0}^{\overline{N}-1} q^{-k} \left(e_{\alpha_{i+1}+\cdots+\alpha_{n}} \right)^{k}
 e_{\alpha_{i}+\cdots+\alpha_{n}}\left(e_{\alpha_{i+1}+\cdots+\alpha_{n}} \right)^{\overline{N}-1-k} \nonumber \\
& = & e_{\alpha_{i}+\cdots+\alpha_{n}}\left(e_{\alpha_{i+1}+\cdots+\alpha_{n}} \right)^{\overline{N}-1} 
  + \sum_{k=1}^{\overline{N}-1} q^{-k} \left(e_{\alpha_{i+1}+\cdots+\alpha_{n}} \right)^{k}
 e_{\alpha_{i}+\cdots+\alpha_{n}}\left(e_{\alpha_{i+1}+\cdots+\alpha_{n}} \right)^{\overline{N}-1-k} \nonumber \\
& = &  e_{\alpha_{i}+\cdots+\alpha_{n}}\left(e_{\alpha_{i+1}+\cdots+\alpha_{n}} \right)^{\overline{N}-1} 
 + \sum_{k=1}^{\overline{N}-1} (-q)^{-k} 
e_{\alpha_{i}+\cdots+\alpha_{n}}\left(e_{\alpha_{i+1}+\cdots+\alpha_{n}} \right)^{\overline{N}-1} \nonumber \\
& & + \sum_{k=1}^{\overline{N}-1} q^{1-2k} (k)_{q} 
\left[  e_{\alpha_{i}+\cdots+\alpha_{n}}, e_{\alpha_{i+1}+\cdots+\alpha_{n}}  \right]_{q}
\left(e_{\alpha_{i+1}+\cdots+\alpha_{n}} \right)^{\overline{N}-2}, \label{appendixB:julietassone(10)}
\end{eqnarray}
where we have used the following calculation:
\begin{eqnarray*}
\lefteqn{
e_{\alpha_{i}+\cdots+\alpha_{n}} \left(e_{\alpha_{i+1}+\cdots+\alpha_{n}} \right)^{k}
-(-1)^{k}\left(e_{\alpha_{i+1}+\cdots+\alpha_{n}} \right)^{k}e_{\alpha_{i}+\cdots+\alpha_{n}} } \\
& = & \left[e_{\alpha_{i}+\cdots+\alpha_{n}}, \left(e_{\alpha_{i+1}+\cdots+\alpha_{n}} \right)^{k}\right]_{q} \\
& = & 
\sum_{m=0}^{k-1} (-1)^{m} \left(e_{\alpha_{i+1}+\cdots+\alpha_{n}} \right)^{m}
\left[e_{\alpha_{i}+\cdots+\alpha_{n}}, e_{\alpha_{i+1}+\cdots+\alpha_{n}}\right]_{q}
\left(e_{\alpha_{i+1}+\cdots+\alpha_{n}} \right)^{k-1-m} \\
& = & \sum_{m=0}^{k-1} (-1)^{m} q^{-m} \left[e_{\alpha_{i}+\cdots+\alpha_{n}}, e_{\alpha_{i+1}+\cdots+\alpha_{n}}\right]_{q}
      \left(e_{\alpha_{i+1}+\cdots+\alpha_{n}} \right)^{k-1} \\
& = & (-q)^{-k+1}(k)_{q}\left[e_{\alpha_{i}+\cdots+\alpha_{n}}, e_{\alpha_{i+1}+\cdots+\alpha_{n}}\right]_{q}
      \left(e_{\alpha_{i+1}+\cdots+\alpha_{n}} \right)^{k-1}.
\end{eqnarray*}
Here we used the following calculation:  from Khoroshkin and Tolstoy's proposition, we have
\begin{eqnarray}
\lefteqn{
\left[e_{\alpha_{i}+\cdots+\alpha_{n}}, e_{\alpha_{i+1}+\cdots+\alpha_{n}}\right]_{q} } \nonumber \\
& &  = C_{i+1} e_{\alpha_{i}+2\alpha_{i+1}+ \cdots + 2\alpha_{n}} +  
 \sum_{p=i+2}^{n} C_{p} e_{\alpha_{i}+\cdots+2\alpha_{p}+ \cdots + 2\alpha_{n}}e_{\alpha_{i+1}+\cdots+\alpha_{p-1}}, 
 \ C_{p} \in \mathbb{C}, \label{appendixB:MarkLatham(1)}
\end{eqnarray}
and we also have
$$
\begin{array}{rcll}
\left[e_{\alpha_{i+1}+\cdots+\alpha_{p-1}},e_{\alpha_{i+1}+\cdots+\alpha_{n}}\right]_{q} & = & 0, & p=i+2, \ldots, n, \\
\left[e_{\alpha_{i}+\cdots+2\alpha_{p}+ \cdots + 2\alpha_{n}},e_{\alpha_{i+1}+\cdots+\alpha_{n}}\right]_{q} & = & 0, & 
p=i+2, \ldots, n, \\
\left[ e_{\alpha_{i}+2\alpha_{i+1}+ \cdots + 2\alpha_{n}}, e_{\alpha_{i+1}+\cdots+\alpha_{n}}\right]_{q} & = & 0. & 
\end{array}$$
We can now re-write (\ref{appendixB:julietassone(10)}) as
\begin{eqnarray}
\lefteqn{
(-q)^{1-\overline{N}}(\overline{N})_{q}
e_{\alpha_{i}+\cdots+\alpha_{n}}\left(e_{\alpha_{i+1}+\cdots+\alpha_{n}} \right)^{\overline{N}-1} } \nonumber \\
& &  + \sum_{k=1}^{\overline{N}-1} q^{1-2k} (k)_{q} 
\left[  e_{\alpha_{i}+\cdots+\alpha_{n}}, e_{\alpha_{i+1}+\cdots+\alpha_{n}}  \right]_{q}
\left(e_{\alpha_{i+1}+\cdots+\alpha_{n}} \right)^{\overline{N}-2} \label{eq:acer(99)} \\
& = & \left\{  \begin{array}{ll}
 q^{3} (2N-1)_{q}\left[N\right]^{q^{2}}\left[  e_{\alpha_{i}+\cdots+\alpha_{n}}, e_{\alpha_{i+1}+\cdots+\alpha_{n}}  \right]_{q}
\left(e_{\alpha_{i+1}+\cdots+\alpha_{n}} \right)^{\overline{N}-2},   & N \equiv 1, 3 \pmod{4}, \nonumber \\
q^{3} (N-1)_{q}\left[N/2\right]^{q^{2}}\left[  e_{\alpha_{i}+\cdots+\alpha_{n}}, e_{\alpha_{i+1}+\cdots+\alpha_{n}}  \right]_{q}
\left(e_{\alpha_{i+1}+\cdots+\alpha_{n}} \right)^{\overline{N}-2},  & N \equiv 0 \pmod{4}, \nonumber \\
q^{3}(N/2)_{q}\left[(N-2)/4\right]^{q^{2}}
\left[  e_{\alpha_{i}+\cdots+\alpha_{n}}, e_{\alpha_{i+1}+\cdots+\alpha_{n}}  \right]_{q}
\left(e_{\alpha_{i+1}+\cdots+\alpha_{n}} \right)^{\overline{N}-2}, & N \equiv 2 \pmod{4},
\end{array} \right.  \nonumber \\
& = & 0, \label{appendixB:MarkLatham(2)}
\end{eqnarray}
as  the first term in (\ref{eq:acer(99)}) vanishes from $(\overline{N})_{q}=0$,
and the second term in (\ref{eq:acer(99)}) vanishes from Proposition \ref{appendixB:bronte(1)}
as $\left[N\right]^{q^{2}}=0$ for an odd integer $N$, $\left[N/2\right]^{q^{2}}=0$ if $N \equiv 0 \pmod{4}$, and 
$(N/2)_{q}=0$ if $N \equiv 2 \pmod{4}$.

We prove (iii): for each $j=i+1, \ldots, n$, 
\begin{eqnarray*}
\lefteqn{
\left[  e_{i}, \left(e_{\alpha_{i+1}+\cdots+2\alpha_{j}+\cdots+2\alpha_{n}}\right)^{N'} \right]_{q} } \\
& = & 
\sum_{k=0}^{N'-1} q^{-k} \left(e_{\alpha_{i+1}+\cdots+2\alpha_{j}+\cdots+2\alpha_{n}}\right)^{k}
e_{\alpha_{i}+\cdots+2\alpha_{j}+\cdots+2\alpha_{n}}
\left(e_{\alpha_{i+1}+\cdots+2\alpha_{j}+\cdots+2\alpha_{n}}\right)^{N'-1-k} \\
& = & \sum_{k=0}^{N'-1} q^{-2k}e_{\alpha_{i}+\cdots+2\alpha_{j}+\cdots+2\alpha_{n}}
      \left(e_{\alpha_{i+1}+\cdots+2\alpha_{j}+\cdots+2\alpha_{n}}\right)^{N'-1} 
  = 0,
\end{eqnarray*}
where we have used the identity
$\left[e_{\alpha_{i}+\cdots+2\alpha_{j}+\cdots+2\alpha_{n}},
e_{\alpha_{i+1}+\cdots+2\alpha_{j}+\cdots+2\alpha_{n}}\right]_{q}=0$.
\end{proof}

\end{subsubsection}

\begin{subsubsection}{Case 4: $\mu \prec \alpha_{i}$.}
We deal with this case by breaking it into a number of subcases in the following proposition.
\begin{proposition}  We have
\begin{itemize}
\item[(i)]  $\left[ e_{\alpha_{j}+\cdots+\alpha_{k}}, e_{i}\right]_{q}=0$, for all $1 \leq j < k \leq i-2$,
\item[(ii)] $\left[ \left(e_{\alpha_{j}+\cdots+\alpha_{i-1}}\right)^{N'}, e_{i}\right]_{q}=0$, for each $j=1, \ldots, i-2$,
\item[(iii)]  $\left[ e_{\alpha_{j}+\cdots+\alpha_{i}}, e_{i}\right]_{q}=0$, for all $1 \leq j < i \leq n-1$,
\item[(iv)]  $\left[ \left(e_{\alpha_{j}+\cdots+\alpha_{n}}\right)^{\overline{N}}, e_{n}\right]_{q}=0$,
for each $j=1, \ldots, n-1$,
\item[(v)] $\left[ e_{\alpha_{j}+\cdots+\alpha_{k}}, e_{i}\right]_{q}=0$, for all $1 \leq j < i < k \leq n$,
\item[(vi)] $\left[ e_{\alpha_{j}+\cdots+2\alpha_{k}+ \cdots + 2\alpha_{n}}, e_{i}\right]_{q}=0$,
for all $1 \leq j < k \leq i \leq n$,
\item[(vii)] $\left[ \left(e_{\alpha_{j}+\cdots+2\alpha_{i+1}+ \cdots + 2\alpha_{n}}\right)^{N'}, e_{i}\right]_{q}=0$,
for each $j=1, \ldots, i-1$,
\item[(viii)] $\left[ e_{\alpha_{j}+\cdots+2\alpha_{k}+ \cdots + 2\alpha_{n}}, e_{i}\right]_{q}=0$,
for all $1 \leq j < i \leq k-2$, where $k \leq n$.
\end{itemize}
\end{proposition}
\begin{proof}
The proof of (i) follows from the Serre relation stating that
$e_{r}$ and $e_{t}$ commute if $|r-t|>1$.  We prove (ii):
$$
\left[ \left(e_{\alpha_{j}+\cdots+\alpha_{i-1}}\right)^{N'}, e_{i}\right]_{q} 
 =  \sum_{k=0}^{N'-1} q^{N'-1-2k} e_{\alpha_{j}+\cdots+\alpha_{i}}
    \left(e_{\alpha_{j}+\cdots+\alpha_{i-1}}\right)^{N'-1} 
 =  0.
$$
A trivial calculation proves (iii), and we now prove (iv):
$$
\left[ \left(e_{\alpha_{j}+\cdots+\alpha_{n}}\right)^{\overline{N}}, e_{n}\right]_{q} 
 = \sum_{k=0}^{\overline{N}-1} (-1)^{k} q^{\overline{N}-1-k} e_{\alpha_{j}+\cdots+2\alpha_{n}}
\left(e_{\alpha_{j}+\cdots+\alpha_{n}}\right)^{\overline{N}-1} 
 =  0.
$$
The proofs of (v), (vi) and (viii) are trivial, and to complete the proof we prove (vii):
$$
\left[ \left(e_{\alpha_{j}+\cdots+2\alpha_{i+1}+ \cdots + 2\alpha_{n}}\right)^{N'}, e_{i}\right]_{q} 
 = \sum_{k=0}^{N'-1} q^{N'-1-2k}e_{\alpha_{j}+\cdots+2\alpha_{i}+ \cdots + 2\alpha_{n}}
      \left(e_{\alpha_{j}+\cdots+2\alpha_{i+1}+ \cdots + 2\alpha_{n}}\right)^{N'-1} 
 =  0.
$$
\end{proof}

\end{subsubsection}

\end{subsection}

\begin{subsection}{The relations between $(e_{\gamma})^{N'}$, $(e_{\beta})^{\overline{N}}$ and $f_{i}$}

To complete the proof that $(e_{\gamma})^{N'}$ and $(e_{\beta})^{\overline{N}}$ 
(anti)commute with each generator of $U_{q}(\mathfrak{g})$, we now prove that
 $(e_{\gamma})^{N'}$ and $(e_{\beta})^{\overline{N}}$ each 
 (anti)commute with $f_{i}$ for each $i=1, \ldots, n$.
By writing $\mu = \sum_{i=1}^{n} \mu_{i} \alpha_{i}$ and examining the definition of $f_{\mu}$
 and the $U_{q}(\mathfrak{g})$ relations, it is apparent that 
$e_{\mu}$ commutes with $f_{j}$ if $\mu_{j}=0$, which gives a partial solution to the problem.
To complete the consideration of this problem, 
we prove Propositions \ref{appendixB:ralphnader(1)} and \ref{appendixB:ralphnader(2)}.
\begin{proposition}
\label{appendixB:ralphnader(1)}
We have
\begin{itemize}
\item[(i)] $\left[ e_{\alpha_{i}+\cdots+\alpha_{j}}, f_{i} \right] = -q^{-1} K_{i}^{-1} e_{\alpha_{i+1}+\cdots+\alpha_{j}}$,
            for all $1 \leq i<j \leq n$, 
\item[(ii)] $\left[ e_{\alpha_{i}+\cdots+\alpha_{j}}, f_{j} \right] = K_{j} e_{\alpha_{i}+\cdots+\alpha_{j-1}}$,
              for all  $1 \leq i<j \leq n$, 
\item[(iii)] $\left[ e_{\alpha_{i}+\cdots+\alpha_{j}}, f_{k} \right] = 0$, for all $1 \leq i<k<j \leq n$, 
\item[(iv)]  $\left[ e_{\alpha_{i}+\cdots+2\alpha_{k} + \cdots + 2\alpha_{n}},f_{j} \right]=0$, for all 
              $1 \leq i<j<k \leq n$, 
\item[(v)]  $\left[e_{\alpha_{n-1}+2\alpha_{n}},f_{n-1}\right] = -q^{-1}(1+q^{-1}) K_{n-1}^{-1} (e_{n})^{2}$,
\item[(vi)]   $\left[ e_{\alpha_{i}+\cdots + 2\alpha_{n}},f_{n} \right]=-K_{n}e_{\alpha_{i}+\cdots + \alpha_{n}}$, 
              for each $i = 1, \ldots, n-2$,
\item[(vii)]  $\left[ e_{\alpha_{i}+\cdots + 2\alpha_{j}+\cdots+2\alpha_{n}},f_{j} \right] 
            = K_{j}e_{\alpha_{i}+\cdots + 2\alpha_{j+1}+\cdots + 2\alpha_{n}}$, for each $j=i+1, \ldots, n-1$,
\item[(viii)] $\left[ e_{\alpha_{i}+\cdots + 2\alpha_{j}+\cdots+2\alpha_{n}},f_{n} \right]=0$, for each $j=i+1, \ldots, n-1$,
\item[(ix)] $\left[ e_{\alpha_{i}+\cdots + 2\alpha_{j}+\cdots+2\alpha_{n}},f_{k} \right] = 0$, 
               for each $k = i+2, \ldots, n-1$, and $j=i+1, \ldots, k-1$,
\item[(x)]  $\left[e_{\alpha_{i}+\cdots + 2\alpha_{j}+\cdots+2\alpha_{n}},f_{i}\right] 
               = -q^{-1}K_{i}^{-1}e_{\alpha_{i+1}+\cdots + 2\alpha_{j}+\cdots+2\alpha_{n}}$, for each $j=i+2, \ldots, n$,
\item[(xi)]   $\left[e_{\alpha_{i} + 2\alpha_{i+1}+\cdots+2\alpha_{n}},f_{i}\right] = 
		-q^{-1}K_{i}^{-1}e_{\alpha_{i+1}+ 2\alpha_{i+2}+\cdots+2\alpha_{n}}e_{i+1}
     		+q^{-3}K_{i}^{-1}e_{i+1} e_{\alpha_{i+1}+ 2\alpha_{i+2}+\cdots+2\alpha_{n}}$, for each 
		$i=1, \ldots, n-2$,
\end{itemize}
\end{proposition}
\begin{proof}
We prove (i): firstly, for each $i=1, \ldots, n-1$,
\begin{eqnarray*}
\left[ e_{\alpha_{i}+\alpha_{i+1}}, f_{i} \right] & = & \left[ e_{i} e_{i+1} -q^{-1} e_{i+1}e_{i},f_{i}\right] \\
& = & e_{i} \left[e_{i+1},f_{i}\right] + \left[e_{i},f_{i}\right]e_{i+1}
     -q^{-1} e_{i+1}\left[e_{i},f_{i}\right] -q^{-1} \left[e_{i+1},f_{i}\right] e_{i} 
 = -q^{-1}  K_{i}^{-1} e_{i+1}.
\end{eqnarray*}
Now for each $j=i+2, \ldots, n$, 
\begin{eqnarray*}
\left[ e_{\alpha_{i}+\cdots+\alpha_{j}}, f_{i} \right] & = & 
\left[ e_{i} e_{\alpha_{i+1}+\cdots+\alpha_{j}} -q^{-1}e_{\alpha_{i+1}+\cdots+\alpha_{j}}e_{i}, f_{i} \right] \\
& = & e_{i} \left[e_{\alpha_{i+1}+\cdots+\alpha_{j}},f_{i} \right]
+ \left[e_{i},f_{i} \right] e_{\alpha_{i+1}+\cdots+\alpha_{j}} \\
& & -q^{-1} e_{\alpha_{i+1}+\cdots+\alpha_{j}} \left[e_{i}, f_{i} \right]
-q^{-1} \left[ e_{\alpha_{i+1}+\cdots+\alpha_{j}}, f_{i} \right]e_{i} 
= -q^{-1} K_{i}^{-1} e_{\alpha_{i+1}+\cdots+\alpha_{j}},
\end{eqnarray*}
as $\left[e_{\alpha_{i+1}+\cdots+\alpha_{j}},f_{i} \right]=0$.

We prove (ii) using a similar approach to the proof of (i).  Firstly, for each $j=2, \ldots, n$, 
\begin{eqnarray*}
\left[ e_{\alpha_{j-1}+\alpha_{j}},f_{j}\right] & = &  \left[ e_{j-1}e_{j}-q^{-1}e_{j}e_{j-1},f_{j}\right] \\
& = & e_{j-1} \left[e_{j}, f_{j}\right] + (-1)^{[f_{j}]} \left[e_{j-1},f_{j}\right]e_{j}
       -q^{-1} e_{j} \left[e_{j-1},f_{j}\right]-q^{-1}\left[e_{j},f_{j}\right]e_{j-1} \\
& = & K_{j} e_{j-1}.
\end{eqnarray*}
For each $i=1, \ldots, j-2$, 
\begin{eqnarray*}
\left[ e_{\alpha_{i}+\cdots+\alpha_{j}}, f_{j} \right] & = & 
\left[e_{\alpha_{i}+\cdots+\alpha_{j-1}}e_{j}-q^{-1}e_{j}e_{\alpha_{i}+\cdots+\alpha_{j-1}}, f_{j} \right] \\
& = & e_{\alpha_{i}+\cdots+\alpha_{j-1}}\left[e_{j}, f_{j} \right]
+(-1)^{[f_{j}]}\left[ e_{\alpha_{i}+\cdots+\alpha_{j-1}}, f_{j} \right]e_{j} \\
& & -q^{-1} e_{j} \left[e_{\alpha_{i}+\cdots+\alpha_{j-1}}, f_{j} \right]
    -q^{-1} \left[e_{j}, f_{j} \right] e_{\alpha_{i}+\cdots+\alpha_{j-1}} 
= K_{j}e_{\alpha_{i}+\cdots+\alpha_{j-1}},
\end{eqnarray*}
as $\left[e_{\alpha_{i}+\cdots+\alpha_{j-1}}, f_{j} \right]=0$.
We now prove (iii) using (ii) and induction.  For each $k=2, \ldots,n-1$, 
\begin{eqnarray*}
\left[e_{\alpha_{k-1}+\alpha_{k}+\alpha_{k+1}},f_{k}\right] & = & 
\left[ e_{\alpha_{k-1}+\alpha_{k}} e_{k+1}-q^{-1}e_{k+1}e_{\alpha_{k-1}+\alpha_{k}}, f_{k}\right] \\
& = & e_{\alpha_{k-1}+\alpha_{k}}\left[ e_{k+1}, f_{k}\right] + 
\left[ e_{\alpha_{k-1}+\alpha_{k}}, f_{k}\right] e_{k+1} \\
& & -q^{-1}e_{k+1} \left[ e_{\alpha_{k-1}+\alpha_{k}}, f_{k}\right]
    -q^{-1} \left[e_{k+1}, f_{k}\right]e_{\alpha_{k-1}+\alpha_{k}} \\
& = & K_{k} e_{k-1} e_{k+1} - K_{k}e_{k+1}e_{k-1} 
  = 0.
\end{eqnarray*}
Keeping $k$ fixed, assume that
$\left[ e_{\alpha_{i}+\cdots+\alpha_{k+1}}, f_{k} \right] = 0$ for some $i=2, \ldots, k-1$, then
\begin{eqnarray*}
\left[ e_{\alpha_{i-1}+\cdots+\alpha_{k+1}}, f_{k} \right] & = & 
\left[ e_{i-1} e_{\alpha_{i}+\cdots+\alpha_{k+1}} -q^{-1}e_{\alpha_{i}+\cdots+\alpha_{k+1}}e_{i-1}, f_{k} \right] \\
& = & e_{i-1} \left[ e_{\alpha_{i}+\cdots+\alpha_{k+1}}, f_{k} \right] 
      + \left[e_{i-1}, f_{k} \right]e_{\alpha_{i}+\cdots+\alpha_{k+1}} \\
& & -q^{-1}e_{\alpha_{i}+\cdots+\alpha_{k+1}} \left[e_{i-1}, f_{k} \right]
    -q^{-1} \left[e_{\alpha_{i}+\cdots+\alpha_{k+1}}, f_{k} \right]e_{i-1} 
  = 0,
\end{eqnarray*}
as $\left[ e_{\alpha_{i}+\cdots+\alpha_{k+1}}, f_{k} \right] = [e_{i-1},f_{k}]=0$.  Now assume that
$\left[ e_{\alpha_{i}+\cdots+\alpha_{j}}, f_{k} \right] = 0$ for some $i = 1, \ldots, k-1$,
and some $j=k+1, \ldots, n-1$, then
$$
\left[ e_{\alpha_{i}+\cdots+\alpha_{j+1}}, f_{k} \right] =  
\left[ e_{\alpha_{i}+\cdots+\alpha_{j}}e_{j+1}-q^{-1}e_{j+1}e_{\alpha_{i}+\cdots+\alpha_{j}}, f_{k} \right] = 0, $$
as $\left[e_{\alpha_{i}+\cdots+\alpha_{j}}, f_{k} \right]=\left[e_{j+1}, f_{k}\right]=0$, 
which completes the proof of (iii).

We now prove (iv).  For each $j=i+1, \ldots, n-1$, 
$$\left[e_{\alpha_{i}+ \cdots + 2\alpha_{n}},f_{j} \right] =
\left[e_{\alpha_{i}+ \cdots + \alpha_{n}}e_{n}+e_{n}e_{\alpha_{i}+ \cdots + \alpha_{n}},f_{j}\right]
 =  0,$$
as $\left[e_{\alpha_{i}+ \cdots + \alpha_{n}},f_{j} \right] = 0$ from (iii) and $[e_{n},f_{j}]=0$.  
Keeping $j$ fixed,  assume that
$\left[e_{\alpha_{i}+ \cdots + 2\alpha_{k}+\cdots+2\alpha_{n}},f_{j} \right]=0$ for some
$k=j+2, \ldots, n$, then
$$
\left[e_{\alpha_{i}+ \cdots + 2\alpha_{k-1}+\cdots+2\alpha_{n}},f_{j} \right] = 
\left[e_{\alpha_{i}+ \cdots + 2\alpha_{k}+\cdots+2\alpha_{n}}e_{k-1}
   -q^{-1}e_{k-1}e_{\alpha_{i}+ \cdots + 2\alpha_{k}+\cdots+2\alpha_{n}},f_{j}\right] 
 = 0,
$$
as $\left[e_{\alpha_{i}+ \cdots + 2\alpha_{k}+\cdots+2\alpha_{n}},f_{j} \right]=
\left[e_{k-1},f_{j}\right]=0$.

We prove (v) using (i):
\begin{eqnarray*}
\left[e_{\alpha_{n-1}+2\alpha_{n}},f_{n-1}\right] & = & 
\left[ e_{\alpha_{n-1}+\alpha_{n}} e_{n} + e_{n}e_{\alpha_{n-1}+\alpha_{n}},f_{n-1}\right] \\
& = & e_{\alpha_{n-1}+\alpha_{n}}\left[e_{n},f_{n-1}\right]
     + \left[e_{\alpha_{n-1}+\alpha_{n}},f_{n-1}\right]e_{n} \\
& &  + e_{n}\left[e_{\alpha_{n-1}+\alpha_{n}},f_{n-1}\right]
     + \left[e_{n},f_{n-1}\right]e_{\alpha_{n-1}+\alpha_{n}} \\
& = & -q^{-1}(1+q^{-1}) K_{n-1}^{-1} (e_{n})^{2}.
\end{eqnarray*}

We prove (vi) using (ii):
\begin{eqnarray*}
\left[ e_{\alpha_{i}+\cdots + 2\alpha_{n}},f_{n} \right] & = & 
\left[ e_{\alpha_{i}+\cdots + \alpha_{n}} e_{n} + e_{n}e_{\alpha_{i}+\cdots + \alpha_{n}}, f_{n} \right] \\
& = & e_{\alpha_{i}+\cdots + \alpha_{n}}\left[e_{n}, f_{n} \right]
      -\left[e_{\alpha_{i}+\cdots + \alpha_{n}}, f_{n} \right]e_{n} \\
&   & + e_{n} \left[e_{\alpha_{i}+\cdots + \alpha_{n}}, f_{n} \right]    
     - \left[e_{n}, f_{n} \right] e_{\alpha_{i}+\cdots + \alpha_{n}} \\
& = &  -K_{n}e_{\alpha_{i}+\cdots + \alpha_{n-1}}e_{n}
      + q^{-1} K_{n}e_{n}e_{\alpha_{i}+\cdots + \alpha_{n-1}} 
  = -K_{n}e_{\alpha_{i}+\cdots + \alpha_{n}},
\end{eqnarray*}
as $\left[e_{n}, f_{n} \right]$ commutes with $e_{\alpha_{i}+\cdots+\alpha_{n}}$.

We prove (vii) using (iv). For each $j=i+1, \ldots, n-1$, 
\begin{eqnarray*}
\left[ e_{\alpha_{i}+\cdots + 2\alpha_{j}+\cdots+2\alpha_{n}},f_{j} \right] 
& = & 
\left[ e_{\alpha_{i}+\cdots + 2\alpha_{j+1}+\cdots + 2\alpha_{n}}e_{j}
     -q^{-1}e_{j}e_{\alpha_{i}+\cdots + 2\alpha_{j+1}+\cdots + 2\alpha_{n}},f_{j} \right] \\
& = & e_{\alpha_{i}+\cdots + 2\alpha_{j+1}+\cdots + 2\alpha_{n}}\left[e_{j},f_{j} \right]
      +\left[e_{\alpha_{i}+\cdots + 2\alpha_{j+1}+\cdots + 2\alpha_{n}},f_{j} \right]e_{j} \\
& &      -q^{-1}e_{j}\left[e_{\alpha_{i}+\cdots + 2\alpha_{j+1}+\cdots + 2\alpha_{n}} ,f_{j} \right]
      -q^{-1} \left[e_{j},f_{j} \right]e_{\alpha_{i}+\cdots + 2\alpha_{j+1}+\cdots + 2\alpha_{n}} \\
& = & K_{j}e_{\alpha_{i}+\cdots + 2\alpha_{j+1}+\cdots + 2\alpha_{n}}.
\end{eqnarray*}

We now prove (viii) using (vi).  Firstly, for each $i=1, \ldots, n-2$, 
\begin{eqnarray*}
\left[ e_{\alpha_{i}+\cdots + 2\alpha_{n-1}+2\alpha_{n}},f_{n} \right] & = & 
\left[ e_{\alpha_{i}+\cdots + 2\alpha_{n}}e_{n-1}-q^{-1}e_{n-1}e_{\alpha_{i}+\cdots + 2\alpha_{n}},f_{n} \right] \\
& = & e_{\alpha_{i}+\cdots + 2\alpha_{n}}\left[e_{n-1},f_{n} \right]
      +\left[e_{\alpha_{i}+\cdots + 2\alpha_{n}},f_{n} \right]e_{n-1} \\
&  &  -q^{-1} e_{n-1}\left[e_{\alpha_{i}+\cdots + 2\alpha_{n}},f_{n} \right]
      -q^{-1} \left[e_{n-1},f_{n} \right]e_{\alpha_{i}+\cdots + 2\alpha_{n}} \\
& = & K_{n}(e_{n-1}e_{\alpha_{i}+\cdots + \alpha_{n}}-e_{\alpha_{i}+\cdots + \alpha_{n}}e_{n-1})
  =  0,
\end{eqnarray*}
as $e_{n-1}$ commutes with $e_{\alpha_{i}+\cdots + \alpha_{n}}$.  Assume that 
$\left[ e_{\alpha_{i}+\cdots + 2\alpha_{j}+\cdots+2\alpha_{n}},f_{n} \right]=0$ for some 
$j=i+2, \ldots, n-1$, then
$$
\left[ e_{\alpha_{i}+\cdots + 2\alpha_{j-1}+\cdots+2\alpha_{n}},f_{n} \right] = 
\left[ e_{\alpha_{i}+\cdots + 2\alpha_{j}+\cdots+2\alpha_{n}}e_{j-1}
       -q^{-1}e_{j-1}e_{\alpha_{i}+\cdots + 2\alpha_{j}+\cdots+2\alpha_{n}},f_{n} \right]
= 0,
$$
as $\left[e_{\alpha_{i}+\cdots + 2\alpha_{j}+\cdots+2\alpha_{n}},f_{n} \right]=
\left[e_{j-1},f_{n} \right]=0$, which proves (viii).

We now prove (ix).  Recall from (vii) that
$$\left[ e_{\alpha_{i}+\cdots + 2\alpha_{k}+\cdots+2\alpha_{n}},f_{k} \right] 
            = K_{k}e_{\alpha_{i}+\cdots + 2\alpha_{k+1}+\cdots + 2\alpha_{n}},$$ 
for each $k = i+1, \ldots, n-1$.  Now for each $k=i+2, \ldots, n-1$, 
\begin{eqnarray*}
\left[ e_{\alpha_{i}+\cdots + 2\alpha_{k-1}+\cdots+2\alpha_{n}},f_{k} \right] & = &
\left[ e_{\alpha_{i}+\cdots + 2\alpha_{k}+\cdots+2\alpha_{n}}e_{k-1} 
-q^{-1}e_{k-1}e_{\alpha_{i}+\cdots + 2\alpha_{k}+\cdots+2\alpha_{n}},f_{k} \right] \\
 & = & e_{\alpha_{i}+\cdots + 2\alpha_{k}+\cdots+2\alpha_{n}}\left[e_{k-1} ,f_{k} \right]
       +\left[e_{\alpha_{i}+\cdots + 2\alpha_{k}+\cdots+2\alpha_{n}},f_{k} \right]e_{k-1} \\
 &  & -q^{-1} e_{k-1}\left[e_{\alpha_{i}+\cdots + 2\alpha_{k}+\cdots+2\alpha_{n}},f_{k} \right]
      -q^{-1} \left[e_{k-1},f_{k} \right]e_{\alpha_{i}+\cdots + 2\alpha_{k}+\cdots+2\alpha_{n}} \\
 & = & K_{k}e_{\alpha_{i}+\cdots + 2\alpha_{k+1}+\cdots + 2\alpha_{n}}e_{k-1}
       -K_{k}e_{k-1}e_{\alpha_{i}+\cdots + 2\alpha_{k+1}+\cdots + 2\alpha_{n}} 
   = 0,
\end{eqnarray*}
as $e_{k-1}$ commutes with $e_{\alpha_{i}+\cdots + 2\alpha_{k+1}+\cdots + 2\alpha_{n}}$.
We now do the inductive step: assume that
$$\left[e_{\alpha_{i}+\cdots + 2\alpha_{j}+\cdots+2\alpha_{n}},f_{k} \right]=0,$$
for some $j=i+2, \ldots, k-1$, then
$$
\left[e_{\alpha_{i}+\cdots + 2\alpha_{j-1}+\cdots+2\alpha_{n}},f_{k} \right] = 
\left[e_{\alpha_{i}+\cdots + 2\alpha_{j}+\cdots+2\alpha_{n}}e_{j-1} 
      -q^{-1}e_{j-1}e_{\alpha_{i}+\cdots + 2\alpha_{j}+\cdots+2\alpha_{n}},f_{k} \right] 
    = 0,
$$
as $\left[e_{\alpha_{i}+\cdots + 2\alpha_{j}+\cdots+2\alpha_{n}},f_{k} \right]=
\left[e_{j-1},f_{k} \right]=0$, which proves (ix).

We now prove (x).  Firstly, for each $i=1, \ldots, n-2$  
(the $i=n-1$ case is dealt with in (v)), we have
\begin{eqnarray*}
\left[e_{\alpha_{i}+\cdots + 2\alpha_{n}},f_{i}\right] & = & 
\left[e_{\alpha_{i}+\cdots + \alpha_{n}}e_{n}+e_{n}e_{\alpha_{i}+\cdots + \alpha_{n}},f_{i}\right] \\
& = & e_{\alpha_{i}+\cdots + \alpha_{n}}\left[e_{n},f_{i}\right]
      + \left[e_{\alpha_{i}+\cdots + \alpha_{n}},f_{i}\right]e_{n} 
        + e_{n}\left[e_{\alpha_{i}+\cdots + \alpha_{n}},f_{i}\right] 
     + \left[ e_{n},f_{i}\right]e_{\alpha_{i}+\cdots + \alpha_{n}} \\
& = & -q^{-1}K_{i}^{-1} \left( e_{\alpha_{i+1}+\cdots + \alpha_{n}}e_{n}
       + e_{n}e_{\alpha_{i+1}+\cdots + \alpha_{n}} \right) 
  = -q^{-1}K_{i}^{-1}e_{\alpha_{i+1}+\cdots + 2\alpha_{n}}.
\end{eqnarray*}
If $i<n-2$, 
\begin{eqnarray*}
\lefteqn{
\left[e_{\alpha_{i}+\cdots + 2\alpha_{n-1}+2\alpha_{n}},f_{i}\right] } \\
& & = 
\left[e_{\alpha_{i}+\cdots + 2\alpha_{n}}e_{n-1}-q^{-1}e_{n-1}e_{\alpha_{i}+\cdots + 2\alpha_{n}},f_{i}\right] \\
& & = e_{\alpha_{i}+\cdots + 2\alpha_{n}}\left[e_{n-1},f_{i}\right]
+\left[e_{\alpha_{i}+\cdots + 2\alpha_{n}},f_{i}\right]e_{n-1} 
  -q^{-1}e_{n-1}\left[e_{\alpha_{i}+\cdots + 2\alpha_{n}},f_{i}\right]
-q^{-1}\left[e_{n-1},f_{i}\right]e_{\alpha_{i}+\cdots + 2\alpha_{n}} \\
& & = -q^{-1}K_{i}^{-1} \left( e_{\alpha_{i+1}+\cdots + 2\alpha_{n}}e_{n-1}
      -q^{-1} e_{n-1}e_{\alpha_{i+1}+\cdots + 2\alpha_{n}} \right) 
 = -q^{-1}K_{i}^{-1}e_{\alpha_{i+1}+\cdots + 2\alpha_{n-1}+2\alpha_{n}}.
\end{eqnarray*}
Let us assume that
$\left[e_{\alpha_{i}+\cdots + 2\alpha_{j}+\cdots+2\alpha_{n}},f_{i}\right] 
 = -q^{-1}K_{i}^{-1}e_{\alpha_{i+1}+\cdots + 2\alpha_{j}+\cdots+2\alpha_{n}}$
for some $j=i+3, \ldots, n-1$, then
\begin{eqnarray*}
\left[e_{\alpha_{i}+\cdots + 2\alpha_{j-1}+\cdots+2\alpha_{n}},f_{i}\right] 
& = & 
\left[e_{\alpha_{i}+\cdots + 2\alpha_{j}+\cdots+2\alpha_{n}}e_{j-1}
   -q^{-1} e_{j-1}e_{\alpha_{i}+\cdots + 2\alpha_{j}+\cdots+2\alpha_{n}},f_{i}\right] \\
& = & e_{\alpha_{i}+\cdots + 2\alpha_{j}+\cdots+2\alpha_{n}}\left[e_{j-1},f_{i}\right]
 + \left[e_{\alpha_{i}+\cdots + 2\alpha_{j}+\cdots+2\alpha_{n}},f_{i}\right]e_{j-1} \\
& & -q^{-1} e_{j-1}\left[e_{\alpha_{i}+\cdots + 2\alpha_{j}+\cdots+2\alpha_{n}},f_{i}\right] 
   -q^{-1}\left[e_{j-1},f_{i}\right]e_{\alpha_{i}+\cdots + 2\alpha_{j}+\cdots+2\alpha_{n}} \\
& = & -q^{-1}K_{i}^{-1} \left( e_{\alpha_{i+1}+\cdots + 2\alpha_{j}+\cdots+2\alpha_{n}}e_{j-1}
     -q^{-1} e_{j-1}e_{\alpha_{i+1}+\cdots + 2\alpha_{j}+\cdots+2\alpha_{n}} \right) \\
& = & -q^{-1}K_{i}^{-1}e_{\alpha_{i+1}+\cdots + 2\alpha_{j-1}+\cdots+2\alpha_{n}},
\end{eqnarray*}
which proves (x).

We now prove (xi).  For each $i=1, \ldots, n-2$, 
\begin{eqnarray*}
\left[e_{\alpha_{i}+2\alpha_{i+1}+ \cdots+2\alpha_{n}},f_{i}\right] & = & 
\left[e_{\alpha_{i}+\alpha_{i+1}+ 2\alpha_{i+2}+\cdots+2\alpha_{n}}e_{i+1}
      -q^{-1}e_{i+1}e_{\alpha_{i}+\alpha_{i+1}+ 2\alpha_{i+2}+\cdots+2\alpha_{n}},f_{i}\right] \\
& = & e_{\alpha_{i}+\alpha_{i+1}+ 2\alpha_{i+2}+\cdots+2\alpha_{n}}\left[e_{i+1},f_{i}\right]
  + \left[e_{\alpha_{i}+\alpha_{i+1}+ 2\alpha_{i+2}+\cdots+2\alpha_{n}},f_{i}\right]e_{i+1} \\
&  & -q^{-1} e_{i+1}\left[e_{\alpha_{i}+\alpha_{i+1}+ 2\alpha_{i+2}+\cdots+2\alpha_{n}},f_{i}\right]
    -q^{-1} \left[e_{i+1},f_{i}\right]e_{\alpha_{i}+\alpha_{i+1}+ 2\alpha_{i+2}+\cdots+2\alpha_{n}} \\
& = & -q^{-1}K_{i}^{-1}e_{\alpha_{i+1}+ 2\alpha_{i+2}+\cdots+2\alpha_{n}}e_{i+1}
     +q^{-3}K_{i}^{-1}e_{i+1} e_{\alpha_{i+1}+ 2\alpha_{i+2}+\cdots+2\alpha_{n}}.
\end{eqnarray*}

\end{proof}

The following proposition completes the proof that 
$(e_{\gamma})^{N'}$ and $(e_{\beta})^{\overline{N}}$ (anti)commute with $f_{i}$ for all $\gamma, \beta \in \phi$.
\begin{proposition}
\label{appendixB:ralphnader(2)}
There are the identities:
\begin{itemize}
\item[(i)]   $\left[ \left(e_{\alpha_{i}+\cdots+\alpha_{j}}\right)^{N'}, f_{i} \right] = 0$,
             for all $1 \leq i<j \leq n-1$,
	     
\item[(ii)]  $\left[ \left(e_{\alpha_{i}+\cdots+\alpha_{n}}\right)^{\overline{N}}, f_{i} \right]  = 0 $,
             for each $i=1, \ldots, n-1$,
	     
\item[(iii)]  $\left[ \left(e_{\alpha_{i}+\cdots+\alpha_{j}}\right)^{N'}, f_{j} \right] = 0$, 
             for all $1 \leq i<j \leq n-1$,
\item[(iv)]  $\left[ \left(e_{\alpha_{i}+\cdots+\alpha_{n}}\right)^{\overline{N}}, f_{n} \right] = 0$, 
              for each $i=1, \ldots, n-1$,

\item[(v)]  $\left[\left(e_{\alpha_{n-1}+2\alpha_{n}}\right)^{N'},f_{n-1}\right] = 0$,
\item[(vi)]   $\left[ \left(e_{\alpha_{i}+\cdots + 2\alpha_{n}}\right)^{N'},f_{n} \right]=0$, 
              for each $i = 1, \ldots, n-2$,
\item[(vii)]  $\left[ \left(e_{\alpha_{i}+\cdots + 2\alpha_{j}+\cdots+2\alpha_{n}}\right)^{N'},f_{j} \right] 
               = 0$, for each $j=i+1, \ldots, n-1$,
\item[(viii)]  $\left[ \left( e_{\alpha_{i}+\cdots + 2\alpha_{j}+\cdots+2\alpha_{n}} \right)^{N'} ,f_{i}\right] 
               = 0$, for each $j=i+2, \ldots, n$,
\item[(ix)]   $\left[\left(e_{\alpha_{i} + 2\alpha_{i+1}+\cdots+2\alpha_{n}}\right)^{N'},f_{i}\right] = 0$, for each 
		$i=1, \ldots, n-2$.
\end{itemize}
\end{proposition}
\begin{proof}
We prove (i):  for all $1 \leq i < j \leq n-1$, 
$$
\left[ \left(e_{\alpha_{i}+\cdots+\alpha_{j}}\right)^{N'}, f_{i} \right] 
  =  -q^{-1} K_{i}^{-1} \sum_{t=0}^{N'-1} q^{2N'-2-2t}   e_{\alpha_{i+1}+\cdots+\alpha_{j}}
              \left(e_{\alpha_{i}+\cdots+\alpha_{j}}\right)^{N'-1} 
  =  0.
$$

We prove (ii).  For each $i=1, \ldots, n-1$, 
\begin{eqnarray*}
\left[\left(e_{\alpha_{i}+\cdots+\alpha_{n}}\right)^{\overline{N}}, f_{i} \right] 
& = & -q^{\overline{N}-2} K_{i}^{-1}   \left(e_{\alpha_{i}+\cdots+\alpha_{n}}\right)^{\overline{N}-1}
      e_{\alpha_{i+1}+\cdots+\alpha_{n}} \\
& &   -q^{\overline{N}-2} K_{i}^{-1} \sum_{t=1}^{\overline{N}-1} q^{-t}
        \left(e_{\alpha_{i}+\cdots+\alpha_{n}}\right)^{\overline{N}-1-t}
        e_{\alpha_{i+1}+\cdots+\alpha_{n}} \left(e_{\alpha_{i}+\cdots+\alpha_{n}}\right)^{t} \\ 
& = & -q^{\overline{N}-2} K_{i}^{-1} \sum_{t=0}^{\overline{N}-1} 
       (-q)^{-t} \left(e_{\alpha_{i}+\cdots+\alpha_{n}}\right)^{\overline{N}-1}
                          e_{\alpha_{i+1}+\cdots+\alpha_{n}} \\
& &   -q^{\overline{N}-2} K_{i}^{-1} \sum_{t=1}^{\overline{N}-1} 
         q^{1-2t} (t)_{q} \left(e_{\alpha_{i}+\cdots+\alpha_{n}}\right)^{\overline{N}-2}
	               e_{\alpha_{i+1}+\cdots+\alpha_{n}} 
 =  0,
\end{eqnarray*}
which vanishes for the same reasons that (\ref{appendixB:MarkLatham(2)}) vanishes.
Here we used the following result:
\begin{eqnarray*}
\lefteqn{
\left(e_{\alpha_{i}+\cdots+\alpha_{n}} \right)^{t} e_{\alpha_{i+1}+\cdots+\alpha_{n}}  
-(-1)^{t}e_{\alpha_{i+1}+\cdots+\alpha_{n}} \left(e_{\alpha_{i}+\cdots+\alpha_{n}} \right)^{t} } \\
& & \hspace{5mm} = \left[ \left(e_{\alpha_{i}+\cdots+\alpha_{n}} \right)^{t},e_{\alpha_{i+1}+\cdots+\alpha_{n}} \right]_{q} \\
& & \hspace{5mm}= \sum_{s=0}^{t-1} (-1)^{s} q^{-s}  \left(e_{\alpha_{i}+\cdots+\alpha_{n}} \right)^{t-1}
  \left[e_{\alpha_{i}+\cdots+\alpha_{n}},e_{\alpha_{i+1}+\cdots+\alpha_{n}} \right]_{q}
    \\
& & \hspace{5mm}= (-q)^{1-t} (t)_{q}  \left(e_{\alpha_{i}+\cdots+\alpha_{n}} \right)^{t-1}
 \left[e_{\alpha_{i}+\cdots+\alpha_{n}},e_{\alpha_{i+1}+\cdots+\alpha_{n}} \right]_{q},
\end{eqnarray*}
in which we have used Eq. (\ref{appendixB:MarkLatham(1)}) and the following results:
$$\begin{array}{rcll}
\left[e_{\alpha_{i}+\cdots+\alpha_{n}}, e_{\alpha_{i}+2\alpha_{i+1}+\cdots+2\alpha_{n}} \right]_{q} & = & 0, & \\
\left[e_{\alpha_{i}+\cdots+\alpha_{n}}, e_{\alpha_{i}+\cdots + 2\alpha_{p}+\cdots+2\alpha_{n}} \right]_{q} & = & 0,
& p=i+2, \ldots, n, \\
\left[e_{\alpha_{i}+\cdots+\alpha_{n}}, e_{\alpha_{i+1}+\cdots+\alpha_{p-1}} \right]_{q} & = & 0,
& p=i+2, \ldots, n.
\end{array} $$

We prove (iii).  For all $1 \leq i < j \leq n-1$,
$$
\left[ \left(e_{\alpha_{i}+\cdots+\alpha_{j}}\right)^{N'}, f_{j} \right] 
= K_{j} \sum_{t=0}^{N'-1} q^{-N'+1+2t} 
       \left(e_{\alpha_{i}+\cdots+\alpha_{j}}\right)^{N'-1}e_{\alpha_{i}+\cdots+\alpha_{j-1}} 
 = 0.
$$

We prove (iv): for each $i=1, \ldots, n-1$,
$$
\left[ \left(e_{\alpha_{i}+\cdots+\alpha_{n}}\right)^{\overline{N}}, f_{n} \right] 
 =  K_{n}  \sum_{t=0}^{\overline{N}-1} (-q)^{t} \left(e_{\alpha_{i}+\cdots+\alpha_{n}}\right)^{\overline{N}-1}
      e_{\alpha_{i}+\cdots+\alpha_{n-1}} 
 = 0.
$$

We prove (v):
$$
\left[\left(e_{\alpha_{n-1}+2\alpha_{n}}\right)^{N'},f_{n-1}\right] 
 = -q^{-1}(1+q^{-1}) K_{n-1}^{-1}  \sum_{t=0}^{N'-1} q^{-2t} \left(e_{\alpha_{n-1}+2\alpha_{n}}\right)^{N'-1}
 (e_{n})^{2} 
 = 0. 
$$

We prove (vi): for each $i = 1, \ldots, n-2$,
$$
\left[ \left(e_{\alpha_{i}+\cdots + 2\alpha_{n}}\right)^{N'},f_{n} \right] 
 = - K_{n} \sum_{t=0}^{N'-1} q^{-N'+1+2t} \left(e_{\alpha_{i}+\cdots + 2\alpha_{n}}\right)^{N'-1}
      e_{\alpha_{i}+\cdots + \alpha_{n}} 
 =  0.
$$

We prove (vii): for each $j=i+1, \ldots, n-1$,
\begin{eqnarray*}
\lefteqn{
\left[ \left(e_{\alpha_{i}+\cdots + 2\alpha_{j}+\cdots+2\alpha_{n}}\right)^{N'},f_{j} \right] } \\
& = & K_{j}  \sum_{t=0}^{N'-1} q^{-N'+1+2t}
           \left(e_{\alpha_{i}+\cdots + 2\alpha_{j}+\cdots+2\alpha_{n}}\right)^{N'-1}
	   e_{\alpha_{i}+\cdots + 2\alpha_{j+1}+\cdots + 2\alpha_{n}} 
  =  0.
\end{eqnarray*}

We prove (viii): for each $j=i+2, \ldots, n$,
\begin{eqnarray*}
\lefteqn{
\left[ \left( e_{\alpha_{i}+\cdots + 2\alpha_{j}+\cdots+2\alpha_{n}} \right)^{N'} ,f_{i}\right] } \\
& = &  -q^{-1} K_{i}^{-1} \sum_{t=0}^{N'-1} q^{N'-1-2t} 
      \left( e_{\alpha_{i}+\cdots + 2\alpha_{j}+\cdots+2\alpha_{n}} \right)^{N'-1}
      e_{\alpha_{i+1}+\cdots + 2\alpha_{j}+\cdots+2\alpha_{n}} 
  =  0.
\end{eqnarray*}

We prove (ix): for each $i=1, \ldots, n-2$,
\begin{eqnarray*}
\lefteqn{
\left[ \left( e_{\alpha_{i} + 2\alpha_{i+1}+\cdots+2\alpha_{n}}\right)^{N'},f_{i}\right] } \\
& = & -q^{-1} K_{i}^{-1} \sum_{t=0}^{N'-1} q^{-2t}
   \left( e_{\alpha_{i} + 2\alpha_{i+1}+\cdots+2\alpha_{n}}\right)^{N'-1}
   e_{\alpha_{i+1}+ 2\alpha_{i+2}+\cdots+2\alpha_{n}} e_{i+1} \\
& & + q^{-3} K_{i}^{-1} \sum_{t=0}^{N'-1}  q^{-2t}
    \left( e_{\alpha_{i} + 2\alpha_{i+1}+\cdots+2\alpha_{n}}\right)^{N'-1}
     e_{i+1} e_{\alpha_{i+1}+ 2\alpha_{i+2}+\cdots+2\alpha_{n}}   
 = 0,
\end{eqnarray*}
where we have used 
$\left[ e_{\alpha_{i} + 2\alpha_{i+1}+\cdots+2\alpha_{n}}, e_{i+1} \right]_{q}=
\left[ e_{\alpha_{i} + 2\alpha_{i+1}+\cdots+2\alpha_{n}},
e_{\alpha_{i+1}+ 2\alpha_{i+2}+\cdots+2\alpha_{n}} \right]_{q}=0$.

\end{proof}

\end{subsection}

\begin{subsection}{The remaining elements of $I$}

We have proved that $(e_{\gamma})^{N'}$ and $(e_{\beta})^{\overline{N}}$ (anti)commute 
with each generator of $U_{q}(\mathfrak{g})$ for each positive root $\gamma, \beta \in \phi$.  
In this subsection we use the antipode $S$ and 
the graded antiautomorphism $\omega$ to prove that each of
$(f_{\gamma})^{N'}$, $(f_{\beta})^{\overline{N}}$, $(\overline{e}_{\gamma})^{N'}$, $(\overline{e}_{\beta})^{\overline{N}}$,
$(\overline{f}_{\gamma})^{N'}$ and $(\overline{f}_{\beta})^{\overline{N}}$ also 
(anti)commute 
with each generator of $U_{q}(\mathfrak{g})$ for each positive root $\gamma, \beta \in \phi$.

As $\omega$ is a graded antiautomorphism,
$\omega(x^{m}) = (-1)^{m(m-1)[x]/2}\big(\omega(x)\big)^{m}$ for each $m \in \mathbb{N}$.
Now for each generator $x \in U_{q}(\mathfrak{g})$, we have
$$\left(e_{\alpha_{i}+\cdots+\alpha_{j}}\right)^{N'} x = 
\pm  x \left(e_{\alpha_{i}+\cdots+\alpha_{j}}\right)^{N'},$$
and applying $\omega$ to this equation shows that
$\left(f_{\alpha_{i}+\cdots+\alpha_{j}}\right)^{N'}$ also (anti)commutes with $x$. 
Using almost identical arguments, we can show that each of 
$\left(f_{\alpha_{i}+\cdots+\alpha_{n}}\right)^{\overline{N}}$
and $\left(f_{\alpha_{i}+\cdots+2\alpha_{j}+\cdots + 2\alpha_{n}}\right)^{N'}$ 
also (anti)commute with each generator of $U_{q}(\mathfrak{g})$ for each $j=i+1, \ldots, n$.

Recall that $S$ is a graded antiautomorphism, that
$S(e_{\eta}) = c_{\eta} \overline{e}_{\eta} K_{\eta}^{-1}$ 
where $0 \neq c_{\eta} \in \mathbb{C}$, and that
$S(f_{\eta}) = c_{\eta} K_{\eta} \overline{f}_{\eta}$ where $0 \neq c_{\eta} \in \mathbb{C}$, 
for each $\eta \in \phi$.
It is then almost trivial to prove that
the remaining elements of $I$ (anti)commute with each generator of $U_{q}(\mathfrak{g})$.
This completes the proof that the left ideal ${\cal{I}} \subset U_{q}(\mathfrak{g})$ 
is a two-sided ideal.

\end{subsection}

\end{section}

\begin{section}{The co-multiplication of $e_{\mu}$, $\mu \in \phi$}
\label{appendixC:arnietermincarrotator(1)}

We now prove that the two-sided ideal ${\cal{I}}$ is a two-sided co-ideal, that is, that
$$\Delta(x) \in {\cal{I}} \otimes U_{q}(\mathfrak{g}) + U_{q}(\mathfrak{g}) \otimes {\cal{I}}, \hspace{20mm}  
\forall x \in {\cal{I}}.$$
To prove that ${\cal{I}}$ is a two-sided co-ideal, it suffices to show that 
$$\Delta(x) \in {\cal{I}} \otimes U_{q}(\mathfrak{g}) + U_{q}(\mathfrak{g}) \otimes {\cal{I}}, \hspace{20mm} 
\forall x \in I,$$
which we will show in a straightforward way.  We will firstly prove that
\begin{equation}
\label{appendixC:Ineedlunchnow(1)} 
\Delta(e_{\gamma})^{N'} \in {\cal{I}} \otimes U_{q}(\mathfrak{g}) + U_{q}(\mathfrak{g}) \otimes {\cal{I}},
\end{equation}  
\begin{equation}
\label{appendixC:Ineedlunchnow(2)}
\Delta(e_{\beta})^{\overline{N}} \in {\cal{I}} \otimes U_{q}(\mathfrak{g}) + U_{q}(\mathfrak{g}) \otimes {\cal{I}},
\end{equation}
for each $\gamma, \beta \in \phi$, and then by using the antipode $S$, the graded antiautomorphism $\omega$, the relations
$$\Delta \circ S = (S \otimes S) \circ \Delta', \hspace{10mm}
\Delta \circ \omega = (\omega \otimes \omega) \circ \Delta',$$
and Propositions \ref{sppendixB:steveieg(1)}--\ref{sppendixB:steveieg(5)}, we will prove the appropriate results for all other
elements of $I$ except $(J_{i}^{\pm N}-1)$:
trivially,
$\Delta(J_{i}^{\pm N}-1) \in 
{\cal{I}} \otimes U_{q}(\mathfrak{g}) + U_{q}(\mathfrak{g}) \otimes {\cal{I}}$.

We now explicitly determine the components of $\Delta(e_{\mu})$
before calculating powers of $\Delta(e_{\mu})$ for each $\mu \in \phi$.
In succeeding sections, we will calculate the commutation
relations between the components of $\Delta(e_{\mu})$, 
and then by using these commutation relations, 
the $q$-binomial theorem and the two generalisations of the binomial theorem in Appendix B, 
we will prove relations
(\ref{appendixC:Ineedlunchnow(1)})--(\ref{appendixC:Ineedlunchnow(2)})
for all $\gamma$ and $\beta$ in $\phi$.

The co-multiplication of $e_{\alpha_{i}}$ for each simple root $\alpha_{i}$ is
given in the definition of $U_{q}(\mathfrak{g})$ and will not be further considered.

It is not a trivial matter to calculate the co-multiplication of $e_{\mu}$ for each non-simple root 
$\mu \in \phi$, as the method for constructing $e_{\mu}$ does not `commute' with the co-multiplication.  
(A similar situation occurs for quantum algebras \cite{cp}.)
We will directly calculate $\Delta(e_{\mu})$:
it is not difficult to calculate  $\Delta(e_{\alpha_{i} + \cdots + \alpha_{j}})$ 
for all $1 \leq i < j \leq n$, nor
is it difficult to calculate $\Delta(e_{\alpha_{i} + \cdots + 2\alpha_{n}})$.  
However, it is much more difficult to calculate
$\Delta(e_{\alpha_{i} + \cdots + 2\alpha_{j} + \cdots + 2\alpha_{n}})$ 
for each $j=i+1, \ldots, n-1$.

After calculating $\Delta(e_{\mu})$, we can calculate
 $\Delta(\overline{e}_{\mu})$, $\Delta(f_{\mu})$ and $\Delta(\overline{f}_{\mu})$
by applying the antipode $S$ and the graded antiautomorphism $\omega$ to $\Delta(e_{\mu})$ 
as in the following equations.
In these equations the proportionality sign means that the left hand side is proportional to
the right hand side with a non-zero scalar constant of proportionality:
$$\begin{array}{rcccl}
\Delta(f_{\mu}) & \propto & \Delta \big(\omega(e_{\mu})\big) & \propto 
& (\omega \otimes \omega) \circ \Delta'(e_{\mu}), \\
\Delta(\overline{e}_{\mu}) & \propto & \Delta \big(S(e_{\mu})\big)\Delta(K_{\mu}) & \propto 
& \left[ (S \otimes S ) \circ \Delta'(e_{\mu}) \right] \Delta(K_{\mu}), \\
\Delta(\overline{f}_{\mu}) & \propto & \Delta \big( \omega(\overline{e}_{\mu})\big)
 & \propto & (\omega \otimes \omega) \circ \Delta'(\overline{e}_{\mu}).
\end{array}$$

\begin{subsection}{$\Delta(e_{\alpha_{i} + \cdots + \alpha_{j}})$, $1 \leq i < j \leq n$}

\begin{lemma}
\label{appendixB:simoneyounglem1}
  For each $j = i+1, \ldots, n$, the co-multiplication of $e_{\alpha_{i} + \cdots + \alpha_{j}}$ is
\begin{eqnarray}
\Delta(e_{\alpha_{i} + \cdots + \alpha_{j}}) & = & 
e_{\alpha_{i} + \cdots + \alpha_{j}} \otimes K_{\alpha_{i} + \cdots + \alpha_{j}} + 
1 \otimes e_{\alpha_{i} + \cdots + \alpha_{j}}  \nonumber \\
& & + (q-q^{-1}) \sum_{k=i+1}^{j} e_{\alpha_{k} + \cdots + \alpha_{j}} \otimes K_{\alpha_{k} + \cdots + \alpha_{j}}
 e_{\alpha_{i} + \cdots + \alpha_{k-1}}.  \label{appendB:coidealA}
\end{eqnarray}
\end{lemma}
\begin{proof}
Obviously,
$\Delta(e_{\alpha_{i} + \cdots + \alpha_{j}}) = 
\Delta(e_{\alpha_{i} + \cdots + \alpha_{j-1}})\Delta( e_{j}) 
- q^{-1} \Delta(e_{j})\Delta(e_{\alpha_{i} + \cdots + \alpha_{j-1}})$.
We will prove the lemma inductively.  
Firstly, we calculate $\Delta(e_{\alpha_{i} + \alpha_{i+1}})$ for each $i=1, \ldots, n-1$:
\begin{eqnarray*}
\Delta(e_{\alpha_{i} + \alpha_{i+1}}) & = & 
(e_{i} \otimes K_{i} + 1 \otimes e_{i})(e_{i+1} \otimes K_{i+1} + 1 \otimes e_{i+1}) \\
& &   -q^{-1}(e_{i+1} \otimes K_{i+1} + 1 \otimes e_{i+1})(e_{i} \otimes K_{i} + 1 \otimes e_{i}) \\
  & = & e_{i} e_{i+1} \otimes K_{i}K_{i+1} -q^{-1}e_{i+1}e_{i} \otimes K_{i+1}K_{i} 
       + 1 \otimes e_{i}e_{i+1} -q^{-1} \left(1 \otimes e_{i+1}e_{i}\right) \\
& & + e_{i} \otimes K_{i}e_{i+1} -q^{-1}\left(e_{i+1} \otimes K_{i+1}e_{i} \right)  
   + e_{i+1} \otimes e_{i}K_{i+1} - q^{-1}\left( e_{i} \otimes e_{i+1}K_{i} \right) \\
& = & e_{\alpha_{i} + \alpha_{i+1}} \otimes K_{\alpha_{i} + \alpha_{i+1}}
     + 1 \otimes e_{\alpha_{i} + \alpha_{i+1}} 
     + (q-q^{-1}) \left(e_{i+1} \otimes K_{i+1}e_{i}\right).
\end{eqnarray*}
Now assume that (\ref{appendB:coidealA}) is true for some $j \geq i+1$, 
then we need to calculate
$$\Delta(e_{\alpha_{i} + \cdots + \alpha_{j+1}}) = 
\Delta(e_{\alpha_{i} + \cdots + \alpha_{j}}) \Delta(e_{j+1}) 
- q^{-1} \Delta(e_{j+1})\Delta(e_{\alpha_{i} + \cdots + \alpha_{j}}).$$
Firstly,
\begin{eqnarray*}
\lefteqn{
\left(e_{\alpha_{i} + \cdots + \alpha_{j}} \otimes K_{\alpha_{i} + \cdots + \alpha_{j}} + 
1 \otimes e_{\alpha_{i} + \cdots + \alpha_{j}}\right)
\left(e_{j+1} \otimes K_{j+1} + 1 \otimes e_{j+1}\right)  } \\
& & -q^{-1}\left(e_{j+1} \otimes K_{j+1} + 1 \otimes e_{j+1}\right) 
\left(e_{\alpha_{i} + \cdots + \alpha_{j}} \otimes K_{\alpha_{i} + \cdots + \alpha_{j}} + 
1 \otimes e_{\alpha_{i} + \cdots + \alpha_{j}}\right) \\
& = & e_{\alpha_{i} + \cdots + \alpha_{j+1}} \otimes K_{\alpha_{i} + \cdots + \alpha_{j+1}} + 
1 \otimes e_{\alpha_{i} + \cdots + \alpha_{j+1}} 
+ (q-q^{-1}) e_{j+1} \otimes K_{j+1} e_{\alpha_{i} + \cdots \alpha_{j}},
\end{eqnarray*}
and to determine the remaining components of $\Delta(e_{\alpha_{i} + \cdots + \alpha_{j+1}})$  
we calculate that
\begin{eqnarray*}
\lefteqn{ \left(e_{\alpha_{k} + \cdots + \alpha_{j}} \otimes K_{\alpha_{k} + \cdots + \alpha_{j}}
 e_{\alpha_{i} + \cdots + \alpha_{k-1}} \right)
 \left( e_{j+1} \otimes K_{j+1} + 1 \otimes e_{j+1} \right) } \\
 & & - q^{-1} \left( e_{j+1} \otimes K_{j+1} + 1 \otimes e_{j+1} \right)
 \left(e_{\alpha_{k} + \cdots + \alpha_{j}} \otimes K_{\alpha_{k} + \cdots + \alpha_{j}}
 e_{\alpha_{i} + \cdots + \alpha_{k-1}}\right) \\
 & = & e_{\alpha_{k} + \cdots + \alpha_{j+1}} \otimes K_{\alpha_{k} + \cdots + \alpha_{j+1}}
 e_{\alpha_{i} + \cdots + \alpha_{k-1}}.
\end{eqnarray*}
\end{proof}

\end{subsection}

\begin{subsection}{$\Delta(e_{\alpha_{i} + \cdots + 2\alpha_{n}})$, $i=1, \ldots, n-1$}

\begin{lemma}
For each $i = 1, 2, \ldots, n-1$, the co-multiplication of $e_{\alpha_{i} + \cdots + 2\alpha_{n}}$ is
\begin{eqnarray*}
\Delta(e_{\alpha_{i} + \cdots + 2\alpha_{n}})
& = &
e_{\alpha_{i} + \cdots + 2\alpha_{n}} \otimes K_{\alpha_{i} + \cdots + 2\alpha_{n}} 
+ 1 \otimes e_{\alpha_{i} + \cdots + 2\alpha_{n}}  \\
& &  + (q-q^{-1})  
 \sum_{k=i+1}^{n-1} e_{\alpha_{k} + \cdots + 2\alpha_{n}} \otimes 
 K_{\alpha_{k} + \cdots + 2\alpha_{n}} e_{\alpha_{i} + \cdots + \alpha_{k-1}} \\
& & + (q-q^{-1}) \Big( (1+q)(e_{n})^{2} \otimes (K_{n})^{2} e_{\alpha_{i} + \cdots + \alpha_{n-1}}
 + e_{n} \otimes K_{n} e_{\alpha_{i} + \cdots + \alpha_{n}} \Big).
\end{eqnarray*}
\end{lemma}
\begin{proof}
Recall that 
$e_{\alpha_{i} + \cdots + 2 \alpha_{n}} = 
e_{\alpha_{i} + \cdots + \alpha_{n}} e_{n} + e_{n} e_{\alpha_{i} + \cdots + \alpha_{n}}$. 
Now from Lemma \ref{appendixB:simoneyounglem1},
\begin{eqnarray}
\Delta(e_{\alpha_{i} + \cdots + \alpha_{n}}) & = & 
e_{\alpha_{i} + \cdots + \alpha_{n}} \otimes K_{\alpha_{i} + \cdots + \alpha_{n}} + 
1 \otimes e_{\alpha_{i} + \cdots + \alpha_{n}}  \nonumber \\
& & + (q-q^{-1}) \sum_{k=i+1}^{n} e_{\alpha_{k} + \cdots + \alpha_{n}} \otimes K_{\alpha_{k} + \cdots + \alpha_{n}}
 e_{\alpha_{i} + \cdots + \alpha_{k-1}}, \label{appendixB:simoneyoungeq1}
\end{eqnarray}
and we determine $\Delta(e_{\alpha_{i} + \cdots + 2\alpha_{n}})$ from the following calculations.
Firstly,
\begin{eqnarray*}
\lefteqn{
\big( e_{\alpha_{i} + \cdots + \alpha_{n}} \otimes K_{\alpha_{i} + \cdots + \alpha_{n}} + 
1 \otimes e_{\alpha_{i} + \cdots + \alpha_{n}} \big) \Delta(e_{n})  
 + \Delta(e_{n})\big( e_{\alpha_{i} + \cdots + \alpha_{n}} \otimes K_{\alpha_{i} + \cdots + \alpha_{n}} + 
1 \otimes e_{\alpha_{i} + \cdots + \alpha_{n}} \big) }  \\
& = & e_{\alpha_{i} + \cdots + \alpha_{n}}e_{n} \otimes K_{\alpha_{i} + \cdots + 2\alpha_{n}}
      + e_{n} e_{\alpha_{i} + \cdots + \alpha_{n}} \otimes K_{\alpha_{i} + \cdots + 2\alpha_{n}} 
+ 1 \otimes e_{\alpha_{i} + \cdots + \alpha_{n}}e_{n} + 1 \otimes e_{n}e_{\alpha_{i} + \cdots + \alpha_{n}} \\
& = &  e_{\alpha_{i} + \cdots + 2\alpha_{n}} \otimes K_{\alpha_{i} + \cdots + 2\alpha_{n}} + 
1 \otimes e_{\alpha_{i} + \cdots + 2\alpha_{n}},
\end{eqnarray*}
and furthermore,
\begin{eqnarray*}
\lefteqn{
\sum_{k=i+1}^{n} \Big[
e_{\alpha_{k} + \cdots + \alpha_{n}} \otimes K_{\alpha_{k} + \cdots + \alpha_{n}}e_{\alpha_{i} + \cdots + \alpha_{k-1}}
\Delta(e_{n}) + 
\Delta(e_{n})e_{\alpha_{k} + \cdots + \alpha_{n}} \otimes K_{\alpha_{k} + \cdots + \alpha_{n}}e_{\alpha_{i} + \cdots + \alpha_{k-1}}
 \Big] } \\
 & & =  q (e_{n})^{2} \otimes (K_{n})^{2} e_{\alpha_{i} + \cdots + \alpha_{n-1}}
      + (e_{n})^{2} \otimes (K_{n})^{2} e_{\alpha_{i} + \cdots + \alpha_{n-1}} \\      
& & \hspace{5mm}
      + e_{n} \otimes K_{n} e_{\alpha_{i} + \cdots + \alpha_{n-1}} e_{n}
      - q^{-1} e_{n} \otimes K_{n} e_{n} e_{\alpha_{i} + \cdots + \alpha_{n-1}} \\
& &   \hspace{5mm} + \sum_{k=i+1}^{n-1} \Big[
     e_{\alpha_{k} + \cdots + \alpha_{n}}e_{n} \otimes K_{\alpha_{k} + \cdots + 2\alpha_{n}}
     e_{\alpha_{i} + \cdots + \alpha_{k-1}}
     +e_{n} e_{\alpha_{k} + \cdots + \alpha_{n}} \otimes K_{\alpha_{k} + \cdots + 2\alpha_{n}}
     e_{\alpha_{i} + \cdots + \alpha_{k-1}}  \\
& & \hspace{23mm}
     + e_{\alpha_{k} + \cdots + \alpha_{n}} \otimes K_{\alpha_{k} + \cdots + \alpha_{n}}e_{\alpha_{i} + \cdots + \alpha_{k-1}}e_{n}
     - e_{\alpha_{k} + \cdots + \alpha_{n}} \otimes K_{\alpha_{k} + \cdots + \alpha_{n}}e_{\alpha_{i} + \cdots +
 \alpha_{k-1}}e_{n} \Big] \\
& & = \sum_{k=i+1}^{n-1} 
     e_{\alpha_{k} + \cdots + 2\alpha_{n}}e_{n} \otimes K_{\alpha_{k} + \cdots + 2\alpha_{n}}
     e_{\alpha_{i} + \cdots + \alpha_{k-1}} \\
& & \hspace{10mm}  + (1+q)(e_{n})^{2} \otimes (K_{n})^{2} e_{\alpha_{i} + \cdots + \alpha_{n-1}}
     + e_{n} \otimes K_{n} e_{\alpha_{i} + \cdots + \alpha_{n}}.
\end{eqnarray*}
\end{proof}

\end{subsection}

\begin{subsection}{$\Delta(e_{\alpha_{i} + \cdots + 2\alpha_{j} + \cdots + 2 \alpha_{n}})$, $1 \leq i < j < n$}

\begin{lemma}
\label{appenB:lem10}
For all $1 \leq i < j < n$, 
the co-multiplication of $e_{\alpha_{i} + \cdots + 2\alpha_{j} + \cdots + 2 \alpha_{n}}$ is
$$\Delta( e_{\alpha_{i} + \cdots + 2\alpha_{j} + \cdots + 2 \alpha_{n}} ) =
D_{0} + (q-q^{-1})\sum_{k=i+1}^{n} D_{k} + (q-q^{-1})\overline{D} + (q-q^{-1})\sum_{l=j}^{n-1} F_{l} + F_{j-1},$$
where
$$
\begin{array}{rcll}
D_{0} & = & e_{\alpha_{i} + \cdots + 2\alpha_{j} + \cdots + 2 \alpha_{n}} \otimes 
K_{\alpha_{i} + \cdots + 2\alpha_{j} + \cdots + 2 \alpha_{n}}, & \\
D_{k} & = & e_{\alpha_{k} + \cdots + 2\alpha_{j} + \cdots + 2 \alpha_{n}}
\otimes K_{\alpha_{k} + \cdots + 2\alpha_{j} + \cdots + 2 \alpha_{n}}e_{\alpha_{i} + \cdots + \alpha_{k-1}}, &
  k = i+1, \ldots, j-1, \\
\overline{D} & = &  \left( 
q e_{\alpha_{j} + 2\alpha_{j+1} + \cdots + 2 \alpha_{n}} e_{j} - q^{-1} e_{j}e_{\alpha_{j} + 2\alpha_{j+1} + \cdots + 2 \alpha_{n}}
\right)  & \\
& & \hspace{50mm} \otimes (K_{j} \cdots K_{n})^{2} e_{\alpha_{i} + \cdots + \alpha_{j-1}}, & \\
D_{k} & = &  \overline{e}_{\alpha_{j} + \cdots + 2\alpha_{k+1} + \cdots + 2\alpha_{n}}
\otimes K_{\alpha_{j} + \cdots + 2\alpha_{k+1} + \cdots + 2\alpha_{n}}e_{\alpha_{i} + \cdots + \alpha_{k}},
 &  k = j, \ldots, n-1, \\
D_{n} & = &  \overline{e}_{\alpha_{j} + \cdots + \alpha_{n}} 
\otimes K_{\alpha_{j} + \cdots + \alpha_{n}} e_{\alpha_{i} + \cdots + \alpha_{n}}, & \\
F_{p} & = &  \overline{e}_{\alpha_{j} + \cdots + \alpha_{p}} \otimes K_{\alpha_{j} + \cdots + \alpha_{p}}
e_{\alpha_{i} + \cdots + 2\alpha_{p+1} + \cdots + 2\alpha_{n}}, &  p = j, \ldots, n-1,  \\
F_{j-1} & = & 1 \otimes e_{\alpha_{i} + \cdots + 2\alpha_{j} + \cdots + 2\alpha_{n}}. &
\end{array} $$
\end{lemma}
\begin{proof}
Firstly, recall that
 $$e_{\alpha_{i} + \cdots + 2\alpha_{j} + \cdots + 2 \alpha_{n}} = 
e_{\alpha_{i} + \cdots + 2\alpha_{j+1} + \cdots + 2 \alpha_{n}}e_{j} 
-q^{-1}e_{j}e_{\alpha_{i} + \cdots + 2\alpha_{j+1} + \cdots + 2 \alpha_{n}},
\hspace{7mm} j=i+1, \ldots, n-1.$$
In proving this lemma, we firstly
determine $\Delta(e_{\alpha_{i} + \cdots + 2\alpha_{n-1} + 2\alpha_{n}})$, 
and then prove the remaining cases by induction.

We now obtain two very useful results.  
Set $k=i+1, \ldots, n-2$ and $k < j < n$, and in writing
$e_{\alpha_{i} + \cdots + \alpha_{k-1}}$ in the calculations below we will always consider this to be 
identically equal to $1$ if $k-1 < i$.  Furthermore, if $k-1=i$, then in writing 
$\alpha_{i} + \cdots + \alpha_{k-1}$ we will always take this to mean $\alpha_{i}$.
The first useful result is
\begin{eqnarray}\lefteqn{
\left(e_{\alpha_{k} + \cdots + 2\alpha_{j+1} + \cdots + 2\alpha_{n}} \otimes K_{\alpha_{k} + \cdots + 2\alpha_{j+1} + \cdots + 2\alpha_{n}}
e_{\alpha_{i} + \cdots + \alpha_{k-1}}\right) \Delta(e_{j}) } \nonumber \\
& &
-q^{-1} \Delta(e_{j})
\left(e_{\alpha_{k} + \cdots + 2\alpha_{j+1} + \cdots + 2\alpha_{n}} \otimes K_{\alpha_{k} + \cdots + 2\alpha_{j+1} + \cdots + 2\alpha_{n}}
e_{\alpha_{i} + \cdots + \alpha_{k-1}}\right) \nonumber \\
& = & e_{\alpha_{k} + \cdots + 2\alpha_{j+1} + \cdots + 2\alpha_{n}} e_{j} \otimes 
K_{\alpha_{k} + \cdots + 2\alpha_{j} + \cdots + 2\alpha_{n}}e_{\alpha_{i} + \cdots + \alpha_{k-1}} \nonumber \\
& &  -q^{-1} e_{j} e_{\alpha_{k} + \cdots + 2\alpha_{j+1} + \cdots + 2\alpha_{n}} \otimes
K_{\alpha_{k} + \cdots + 2\alpha_{j} + \cdots + 2\alpha_{n}}e_{\alpha_{i} + \cdots + \alpha_{k-1}} \nonumber \\
& & + e_{\alpha_{k} + \cdots + 2\alpha_{j+1} + \cdots + 2\alpha_{n}} \otimes 
K_{\alpha_{k} + \cdots + 2\alpha_{j+1} + \cdots + 2\alpha_{n}}e_{\alpha_{i} + \cdots + \alpha_{k-1}} e_{j}
\nonumber \\
& & - e_{\alpha_{k} + \cdots + 2\alpha_{j+1} + \cdots + 2\alpha_{n}} \otimes 
K_{\alpha_{k} + \cdots + 2\alpha_{j+1} + \cdots + 2\alpha_{n}}e_{j} e_{\alpha_{i} + \cdots + \alpha_{k-1}}
\nonumber \\
& = & e_{\alpha_{k} + \cdots + 2\alpha_{j} + \cdots + 2\alpha_{n}} \otimes 
K_{\alpha_{k} + \cdots + 2\alpha_{j} + \cdots + 2\alpha_{n}} e_{\alpha_{i} + \cdots + \alpha_{k-1}},
\label{appendixB:rel(1)}
\end{eqnarray}
as $e_{j}$ commutes with $e_{\alpha_{i} + \cdots + \alpha_{k-1}}$.
The second useful result is as follows:  set $i<j<n$, then
\begin{eqnarray}
\lefteqn{
\left(1 \otimes e_{\alpha_{i} + \cdots + 2\alpha_{j+1} + \cdots + 2\alpha_{n}}\right) \Delta(e_{j})
 -q^{-1}\Delta(e_{j}) \left(1 \otimes e_{\alpha_{i} + \cdots + 2\alpha_{j+1} + \cdots + 2\alpha_{n}}\right)  }
 \nonumber \\
 & & = (q-q^{-1}) \left(e_{j} \otimes K_{j} e_{\alpha_{i} + \cdots + 2\alpha_{j+1} + \cdots + 2\alpha_{n}}\right)
       + 1 \otimes e_{\alpha_{i} + \cdots + 2\alpha_{j} + \cdots + 2\alpha_{n}}. \label{appendixB:rel(2)}
\end{eqnarray}
These results imply that
\begin{eqnarray*}
\lefteqn{
\Big[
e_{\alpha_{i} + \cdots + 2\alpha_{n}} \otimes K_{\alpha_{i} + \cdots + 2\alpha_{n}} 
+ 1 \otimes e_{\alpha_{i} + \cdots + 2\alpha_{n}} }  \\
& &  + (q-q^{-1})  
 \sum_{k=i+1}^{n-2} e_{\alpha_{k} + \cdots + 2\alpha_{n}} \otimes 
 K_{\alpha_{k} + \cdots + 2\alpha_{n}} e_{\alpha_{i} + \cdots + \alpha_{k-1}} \Big] \Delta(e_{n-1}) \\
\lefteqn{
 -q^{-1}\Delta(e_{n-1}) \Big[
e_{\alpha_{i} + \cdots + 2\alpha_{n}} \otimes K_{\alpha_{i} + \cdots + 2\alpha_{n}} 
+ 1 \otimes e_{\alpha_{i} + \cdots + 2\alpha_{n}} } \\
& &  \hspace{20mm} + (q-q^{-1})  
 \sum_{k=i+1}^{n-2} e_{\alpha_{k} + \cdots + 2\alpha_{n}} \otimes 
 K_{\alpha_{k} + \cdots + 2\alpha_{n}} e_{\alpha_{i} + \cdots + \alpha_{k-1}} \Big] \\
& = & e_{\alpha_{i} + \cdots + 2\alpha_{n-2}+2\alpha_{n}} \otimes K_{\alpha_{i} + \cdots + 2\alpha_{n-2}+2\alpha_{n}}
+ 1 \otimes e_{\alpha_{i} + \cdots + 2\alpha_{n-2}+2\alpha_{n}} \\
& & + (q-q^{-1}) \sum_{k=i+1}^{n-2} e_{\alpha_{k} + \cdots + 2\alpha_{n-1}+2\alpha_{n}} \otimes 
 K_{\alpha_{k} + \cdots + 2\alpha_{n-1}+2\alpha_{n}} e_{\alpha_{i} + \cdots + \alpha_{k-1}} \\
& &  +(q-q^{-1}) \left(e_{n-1} \otimes K_{n-1} e_{\alpha_{i} + \cdots + 2\alpha_{n}}\right).
\end{eqnarray*}
To complete the proof, we need only simplify
\begin{eqnarray}
\lefteqn{
\Big[ e_{\alpha_{n-1}+ 2\alpha_{n}} \otimes K_{\alpha_{n-1}+ 2\alpha_{n}}e_{\alpha_{i}+ \cdots + \alpha_{n-2}} } \nonumber \\
& &  \hspace{15mm} +(1+q) (e_{n})^{2} \otimes (K_{n})^{2} e_{\alpha_{i}+ \cdots + \alpha_{n-1}}  
  +e_{n} \otimes K_{n}e_{\alpha_{i}+ \cdots + \alpha_{n}}  \Big] \Delta(e_{n-1}) \nonumber  \\
& & -q^{-1} \Delta(e_{n-1}) \Big[  e_{\alpha_{n-1}+ 2\alpha_{n}} \otimes K_{\alpha_{n-1}+ 2\alpha_{n}}e_{\alpha_{i}+ \cdots + \alpha_{n-2}}
 \nonumber \\
& & \hspace{30mm} +(1+q) (e_{n})^{2} \otimes (K_{n})^{2} e_{\alpha_{i}+ \cdots + \alpha_{n-1}} 
 +e_{n} \otimes K_{n}e_{\alpha_{i}+ \cdots + \alpha_{n}}  \Big].  \label{appendixBcomul1}
\end{eqnarray}
To simplify (\ref{appendixBcomul1}) we perform the following calculations: 
\begin{eqnarray*}
\lefteqn{
\left(e_{\alpha_{n-1}+ 2\alpha_{n}} \otimes K_{\alpha_{n-1}+ 2\alpha_{n}}e_{\alpha_{i}+ \cdots + \alpha_{n-2}}\right) \Delta(e_{n-1}) } \\
& & 
\hspace{20mm} -q^{-1} \Delta(e_{n-1}) \left(
e_{\alpha_{n-1}+ 2\alpha_{n}} \otimes K_{\alpha_{n-1}+ 2\alpha_{n}}e_{\alpha_{i}+ \cdots + \alpha_{n-2}} \right) \\
& & \hspace{5mm} =  \left( q e_{\alpha_{n-1}+ 2\alpha_{n}} e_{n-1} - q^{-1} e_{n-1} e_{\alpha_{n-1}+ 2\alpha_{n}} \right)
\otimes (K_{\alpha_{n-1}+ \alpha_{n}})^{2} e_{\alpha_{i}+ \cdots + \alpha_{n-2}} \\
& &   \hspace{10mm} + e_{\alpha_{n-1}+ 2\alpha_{n}} \otimes K_{\alpha_{n-1}+ 2\alpha_{n}} e_{\alpha_{i}+ \cdots + \alpha_{n-1}},
\end{eqnarray*}
and
\begin{eqnarray*}
\lefteqn{
\Big((e_{n})^{2} \otimes (K_{n})^{2} e_{\alpha_{i} + \cdots + \alpha_{n-1}}\Big) \Delta(e_{n-1})
-q^{-1} \Delta(e_{n-1})\Big((e_{n})^{2} \otimes (K_{n})^{2} e_{\alpha_{i} + \cdots + \alpha_{n-1}}\Big) } \\
& = & q^{-1} (e_{n})^{2} e_{n-1} \otimes K_{n-1} (K_{n})^{2} e_{\alpha_{i} + \cdots + \alpha_{n-1}}
-q^{-1} e_{n-1}(e_{n})^{2} \otimes K_{n-1}(K_{n})^{2} e_{\alpha_{i} + \cdots + \alpha_{n-1}} \\
& & + (e_{n})^{2} \otimes (K_{n})^{2} e_{\alpha_{i} + \cdots + \alpha_{n-1}} e_{n-1}
-q (e_{n})^{2} \otimes (K_{n})^{2} e_{n-1} e_{\alpha_{i} + \cdots + \alpha_{n-1}} \\
& = & q^{-1} \Big( (e_{n})^{2} e_{n-1} - e_{n-1}(e_{n})^{2} \Big) 
              \otimes K_{n-1}(K_{n})^{2}e_{\alpha_{i} + \cdots + \alpha_{n-1}},
\end{eqnarray*}
where we have used the fact that $[e_{\alpha_{i} + \cdots + \alpha_{n-1}},e_{n-1}]_{q}=0$.
Furthermore, 
\begin{eqnarray*}
\lefteqn{
\left(e_{n} \otimes K_{n} e_{\alpha_{i} + \cdots + \alpha_{n}}\right) \Delta(e_{n-1}) 
-q^{-1} \Delta(e_{n-1})\left(e_{n} \otimes K_{n} e_{\alpha_{i} + \cdots + \alpha_{n}}\right)  } \\
& = & e_{n}e_{n-1} \otimes K_{n-1}K_{n} e_{\alpha_{i} + \cdots + \alpha_{n}}
      -q^{-1} e_{n-1} e_{n} \otimes K_{n-1}K_{n} e_{\alpha_{i} + \cdots + \alpha_{n}}  \\
 & & + e_{n} \otimes K_{n} e_{\alpha_{i} + \cdots + \alpha_{n}} e_{n-1}
     -q^{-1} e_{n} \otimes e_{n-1} K_{n} e_{\alpha_{i} + \cdots + \alpha_{n}} \\
& = & \overline{e}_{\alpha_{n-1} + \alpha_{n}} \otimes K_{n-1}K_{n} e_{\alpha_{i} + \cdots + \alpha_{n}},
\end{eqnarray*}
where we have used the fact that $[e_{\alpha_{i} + \cdots + \alpha_{n}}, e_{n-1}]_{q}=0$.

Combining these results, we can rewrite (\ref{appendixBcomul1}) as
\begin{eqnarray*}
\lefteqn{
\left( q e_{\alpha_{n-1}+ 2\alpha_{n}} e_{n-1} - q^{-1} e_{n-1} e_{\alpha_{n-1}+ 2\alpha_{n}} \right)
\otimes (K_{n-1})^{2}(K_{n})^{2} e_{\alpha_{i}+ \cdots + \alpha_{n-2}} } \\
& &    + \Big[ e_{\alpha_{n-1}+ 2\alpha_{n}} + 
(1+q^{-1})\left( (e_{n})^{2} e_{n-1} - e_{n-1}(e_{n})^{2} \right) \Big] 
\otimes K_{\alpha_{n-1}+ 2\alpha_{n}} e_{\alpha_{i}+ \cdots + \alpha_{n-1}} \\
& & \hspace{25mm} + \overline{e}_{\alpha_{n-1} + \alpha_{n}} \otimes K_{n-1}K_{n} e_{\alpha_{i} + \cdots + \alpha_{n}} \\
& = & \left( q e_{\alpha_{n-1}+ 2\alpha_{n}} e_{n-1} - q^{-1} e_{n-1} e_{\alpha_{n-1}+ 2\alpha_{n}} \right)
\otimes (K_{n-1})^{2}(K_{n})^{2} e_{\alpha_{i}+ \cdots + \alpha_{n-2}}  \\
& &    + \overline{e}_{\alpha_{n-1}+ 2\alpha_{n}} \otimes K_{\alpha_{n-1}+ 2\alpha_{n}} e_{\alpha_{i}+ \cdots + \alpha_{n-1}} 
   + \overline{e}_{\alpha_{n-1} + \alpha_{n}} \otimes K_{n-1}K_{n} e_{\alpha_{i} + \cdots + \alpha_{n}},
\end{eqnarray*}
where we have used the elementary result that
$$\overline{e}_{\alpha_{n-1}+ 2\alpha_{n}}=
e_{\alpha_{n-1}+ 2\alpha_{n}} + (1+q^{-1})\left[ (e_{n})^{2} e_{n-1} - e_{n-1}(e_{n})^{2} \right].$$

Combining these calculations gives us
\begin{eqnarray*}
\lefteqn{
\Delta(e_{\alpha_{i} + \cdots + 2\alpha_{n-1} + 2\alpha_{n}}) } \\
& = &  e_{\alpha_{i} + \cdots + 2\alpha_{n-1} + 2\alpha_{n}} \otimes K_{\alpha_{i} + \cdots + 2\alpha_{n-1} + 2\alpha_{n}}  \\
      & & + (q-q^{-1})\left[\sum_{k=i+1}^{n-2} e_{\alpha_{k} + \cdots + 2\alpha_{n-1} + 2\alpha_{n}}
           \otimes K_{\alpha_{k} + \cdots + 2\alpha_{n-1} + 2\alpha_{n}} e_{\alpha_{i} + \cdots + \alpha_{k-1}}
	   \right] \\
     & & + (q-q^{-1}) \Big[ \left( q e_{\alpha_{n-1} + 2 \alpha_{n}} e_{n-1} 
           -q^{-1} e_{n-1} e_{\alpha_{n-1} + 2 \alpha_{n}} \right)
            \otimes (K_{n-1}K_{n})^{2} e_{\alpha_{i} + \cdots + \alpha_{n-2}} \Big] \\
     & & + (q-q^{-1}) \Big[\overline{e}_{\alpha_{n-1} + 2\alpha_{n}} \otimes 
           K_{\alpha_{n-1} + 2\alpha_{n}}e_{\alpha_{i} + \cdots + \alpha_{n-1}} \Big] \\
      & & + (q-q^{-1})\Big[
	   \overline{e}_{\alpha_{n-1} + \alpha_{n}} \otimes 
           K_{\alpha_{n-1} + \alpha_{n}}e_{\alpha_{i} + \cdots + \alpha_{n}} \Big] \\
     & & + (q-q^{-1}) \Big[e_{n-1} \otimes K_{n-1} e_{\alpha_{i} + \cdots + 2\alpha_{n}} \Big] \\
     & & + 1 \otimes e_{\alpha_{i} + \cdots + 2\alpha_{n-1} + 2\alpha_{n}},
\end{eqnarray*}
which proves the lemma for $\Delta(e_{\alpha_{i} + \cdots + 2\alpha_{n-1} + 2\alpha_{n}})$.

Now assume that the lemma is true for 
$\Delta( e_{\alpha_{i} + \cdots + 2\alpha_{j} + \cdots + 2 \alpha_{n}})$ for some 
$j= i+2, \ldots, n-1$, then 
Eqs. (\ref{appendixB:rel(1)})--(\ref{appendixB:rel(2)}) allow us to write
\begin{eqnarray*}
\lefteqn{
\left[ D_{0} + (q-q^{-1})\sum_{k=i+1}^{j-2} D_{k} + F_{j-1}\right] \Delta(e_{j-1}) 
- q^{-1} \Delta(e_{j-1}) \left[ D_{0} + (q-q^{-1})\sum_{k=i+1}^{j-2} D_{k} + F_{j-1}\right] }  \\
& = & e_{\alpha_{i} + \cdots + 2\alpha_{j-1} + \cdots + 2\alpha_{n}} \otimes 
K_{\alpha_{i} + \cdots + 2\alpha_{j-1} + \cdots + 2\alpha_{n}} \\
& & + (q-q^{-1})\left( \sum_{k=i+1}^{j-2}e_{\alpha_{k} + \cdots + 2\alpha_{j-1} + \cdots + 2\alpha_{n}} \otimes 
K_{\alpha_{k} + \cdots + 2\alpha_{j-1} + \cdots + 2\alpha_{n}}e_{\alpha_{i}+\cdots+\alpha_{k-1}} \right) \\
& & + (q-q^{-1})\left( e_{j-1} \otimes K_{j-1} e_{\alpha_{i}+\cdots+2\alpha_{j} + \cdots +2\alpha_{n}} \right) \\
& & + 1 \otimes e_{\alpha_{i}+\cdots+2\alpha_{j-1} + \cdots +2\alpha_{n}}.
\end{eqnarray*}

We now perform another calculation that will greatly assist the proof of this lemma.
Let $\varphi$ be any element of $\phi$ satisfying $\alpha_{j} \preceq \varphi \prec \alpha_{j+1}$ with respect to
${\cal{N}}(\phi)$, then
\begin{eqnarray*}
\lefteqn{
\left(\overline{e}_{\varphi} \otimes K_{\varphi} e_{(\alpha_{i}+\cdots+2\alpha_{j}+\cdots+2\alpha_{n}-\varphi)}
\right) \Delta(e_{j-1}) 
-q^{-1}\Delta(e_{j-1})
\left(\overline{e}_{\varphi} \otimes K_{\varphi} e_{(\alpha_{i}+\cdots+2\alpha_{j}+\cdots+2\alpha_{n}-\varphi)}
\right) } \\
& = & \overline{e}_{\varphi} e_{j-1} \otimes 
K_{\varphi}e_{(\alpha_{i}+\cdots+2\alpha_{j}+\cdots+2\alpha_{n}-\varphi)}K_{j-1}
   + \overline{e}_{\varphi} \otimes
   K_{\varphi}e_{(\alpha_{i}+\cdots+2\alpha_{j}+\cdots+2\alpha_{n}-\varphi)}e_{j-1} \\
 & & -q^{-1} e_{j-1} \overline{e}_{\varphi} \otimes K_{j-1} K_{\varphi}
      e_{(\alpha_{i}+\cdots+2\alpha_{j}+\cdots+2\alpha_{n}-\varphi)}e_{j-1}
      -q^{-1} \overline{e}_{\varphi} \otimes e_{j-1}K_{\varphi}
      e_{(\alpha_{i}+\cdots+2\alpha_{j}+\cdots+2\alpha_{n}-\varphi)} \\
 & = & \left[ \overline{e}_{\varphi}, e_{j-1} \right]_{q} \otimes K_{j-1}K_{\varphi}
        e_{(\alpha_{i}+\cdots+2\alpha_{j}+\cdots+2\alpha_{n}-\varphi)} \\
  & & + \overline{e}_{\varphi} \otimes K_{\varphi}
  e_{(\alpha_{i}+\cdots+2\alpha_{j}+\cdots+2\alpha_{n}-\varphi)} e_{j-1}
  -\overline{e}_{\varphi} \otimes K_{\varphi} e_{j-1}
  e_{(\alpha_{i}+\cdots+2\alpha_{j}+\cdots+2\alpha_{n}-\varphi)} \\
 & = & \overline{e}_{\alpha_{j-1} + \varphi} \otimes K_{j-1} K_{\varphi}
     e_{(\alpha_{i}+\cdots+2\alpha_{j}+\cdots+2\alpha_{n}-\varphi)},
\end{eqnarray*}
where we have used the fact that
$\left[e_{(\alpha_{i}+\cdots+2\alpha_{j}+\cdots+2\alpha_{n}-\varphi)},e_{j-1}\right]_{q}=0$
(note that this can be rewritten as 
$e_{(\alpha_{i}+\cdots+2\alpha_{j}+\cdots+2\alpha_{n}-\varphi)} e_{j-1}
 = e_{j-1}e_{(\alpha_{i}+\cdots+2\alpha_{j}+\cdots+2\alpha_{n}-\varphi)}$).
This allows us to write  
\begin{eqnarray*}
\lefteqn{
\left[ (q-q^{-1})\sum_{k=j}^{n} D_{k} + (q-q^{-1})\sum_{l=j}^{n-1}F_{l} \right] \Delta(e_{j-1}) } \\
& &  -q^{-1}\Delta(e_{j-1})\left[ (q-q^{-1})\sum_{k=j}^{n} D_{k} + (q-q^{-1})\sum_{l=j}^{n-1}F_{l} \right]   \\
& = & (q-q^{-1}) \left( 
\sum_{k=j}^{n-1} \overline{e}_{\alpha_{j-1}+\cdots+2\alpha_{k+1}+\cdots+2\alpha_{n}} \otimes
                   K_{\alpha_{j-1}+\cdots+2\alpha_{k+1}+\cdots+2\alpha_{n}} 
		   e_{\alpha_{i}+\cdots+\alpha_{k}} \right) \\
 & & + (q-q^{-1}) \left( \overline{e}_{\alpha_{j-1}+\cdots+2\alpha_{n}} \otimes K_{\alpha_{j-1}+\cdots+\alpha_{n}}
       e_{\alpha_{i}+\cdots+2\alpha_{n}}\right) \\
 & & + (q-q^{-1}) \sum_{l=j}^{n-1}\left( \overline{e}_{\alpha_{j-1}+\cdots+\alpha_{l}} 
       \otimes K_{\alpha_{j-1}+\cdots+\alpha_{l}} e_{\alpha_{i}+\cdots+2\alpha_{l+1}+\cdots+2\alpha_{n}} \right),     
\end{eqnarray*}
and to complete the proof of the lemma we need only calculate
$$\left(\overline{D}+D_{j-1}\right) \Delta(e_{j-1})-q^{-1}\Delta(e_{j-1})\left(\overline{D}+D_{j-1}\right).$$
Now 
\begin{eqnarray}
\lefteqn{
D_{j-1}\Delta(e_{j-1})-q^{-1}\Delta(e_{j-1})D_{j-1} } \nonumber \\
& = & \left(e_{\alpha_{j-1} + 2\alpha_{j} + \cdots + 2 \alpha_{n}}
\otimes K_{\alpha_{j-1} + 2\alpha_{j} + \cdots + 2 \alpha_{n}}e_{\alpha_{i} + \cdots + \alpha_{j-2}}\right)
\Delta(e_{j-1}) \nonumber \\
& & - q^{-1}\Delta(e_{j-1})
\left(e_{\alpha_{j-1} + 2\alpha_{j} + \cdots + 2 \alpha_{n}}
\otimes K_{\alpha_{j-1} + 2\alpha_{j} + \cdots + 2 \alpha_{n}}e_{\alpha_{i} + \cdots + \alpha_{j-2}}\right) \\
& = & q e_{\alpha_{j-1} + 2\alpha_{j} + \cdots + 2 \alpha_{n}} e_{j-1} \otimes 
(K_{\alpha_{j-1}+\cdots+\alpha_{n}})^{2} e_{\alpha_{i}+\cdots+\alpha_{j-2}} \nonumber \\
 & &  -q^{-1} e_{j-1} e_{\alpha_{j-1} + 2\alpha_{j} + \cdots + 2 \alpha_{n}} e_{j-1} \otimes
     (K_{\alpha_{j-1}+\cdots+\alpha_{n}})^{2}e_{\alpha_{i}+\cdots+\alpha_{j-2}} \nonumber \\
 & & +  e_{\alpha_{j-1} + 2\alpha_{j} + \cdots + 2 \alpha_{n}} \otimes 
        K_{\alpha_{j-1} + 2\alpha_{j} + \cdots + 2 \alpha_{n}}  e_{\alpha_{i}+\cdots+\alpha_{j-2}}e_{j-1} \nonumber \\
 & &  -q^{-1} e_{\alpha_{j-1} + 2\alpha_{j} + \cdots + 2 \alpha_{n}}  
 \otimes K_{\alpha_{j-1} + 2\alpha_{j} + \cdots + 2 \alpha_{n}} e_{j-1}e_{\alpha_{i}+\cdots+\alpha_{j-2}} \nonumber \\
 & = & \left( q e_{\alpha_{j-1} + 2\alpha_{j} + \cdots + 2 \alpha_{n}} e_{j-1}
              -q^{-1} e_{j-1} e_{\alpha_{j-1} + 2\alpha_{j} + \cdots + 2 \alpha_{n}}\right)
	      \otimes (K_{\alpha_{j-1}+\cdots + \alpha_{n}})^{2} e_{\alpha_{i}+\cdots+\alpha_{j-2}} \nonumber \\
  & & + e_{\alpha_{j-1} + 2\alpha_{j} + \cdots + 2 \alpha_{n}} \otimes 
  K_{\alpha_{j-1} + 2\alpha_{j} + \cdots + 2 \alpha_{n}} e_{\alpha_{i}+\cdots+\alpha_{j-1}}, \label{appendixB:simoneyoungeq99}
\end{eqnarray}
and
\begin{eqnarray}
\lefteqn{
\overline{D}\Delta(e_{j-1})-q^{-1}\Delta(e_{j-1})\overline{D} } \nonumber \\
& = & \bigg[ e_{\alpha_{j} + 2\alpha_{j+1} + \cdots + 2 \alpha_{n}} e_{j} e_{j-1} 
      -q^{-2}e_{j} e_{\alpha_{j} + 2\alpha_{j+1} + \cdots + 2\alpha_{n}} e_{j-1} \nonumber \\
  & & \hspace{5mm} -e_{j-1}e_{\alpha_{j} + 2\alpha_{j+1} + \cdots + 2\alpha_{n}} e_{j}
      +q^{-2} e_{j-1}e_{j}e_{\alpha_{j} + 2\alpha_{j+1} + \cdots + 2\alpha_{n}}\bigg] \nonumber \\
  & & \hspace{60mm} \otimes K_{\alpha_{j-1} + 2\alpha_{j} + \cdots + 2\alpha_{n}}e_{\alpha_{i}+\cdots+\alpha_{j-1}}.
  \label{appendixB:simoneyoungeq100}
\end{eqnarray}
We combine (\ref{appendixB:simoneyoungeq100}) and the second line of (\ref{appendixB:simoneyoungeq99}) 
using the following calculation:
\begin{eqnarray*}
\lefteqn{
e_{\alpha_{j-1} + 2\alpha_{j} + \cdots + 2\alpha_{n}} + e_{\alpha_{j}+2\alpha_{j+1}+\cdots+2\alpha_{n}}e_{j}e_{j-1}
-e_{j-1}e_{\alpha_{j}+2\alpha_{j+1}+\cdots+2\alpha_{n}}e_{j} } \\
& & -q^{-2}e_{j}e_{\alpha_{j}+2\alpha_{j+1}+\cdots+2\alpha_{n}}e_{j-1}
+q^{-2}e_{j-1}e_{j}e_{\alpha_{j}+2\alpha_{j+1}+\cdots+2\alpha_{n}}  \\
& = & -q^{-1}e_{\alpha_{j}+2\alpha_{j+1}+\cdots+2\alpha_{n}}e_{j-1}e_{j}
      -q^{-1} e_{j}e_{j-1}e_{\alpha_{j}+2\alpha_{j+1}+\cdots+2\alpha_{n}} \\
& & + e_{\alpha_{j}+2\alpha_{j+1}+\cdots+2\alpha_{n}} e_{j} e_{j-1} 
    + q^{-2} e_{j-1} e_{j} e_{\alpha_{j}+2\alpha_{j+1}+\cdots+2\alpha_{n}} \\
& = & e_{\alpha_{j}+2\alpha_{j+1}+\cdots+2\alpha_{n}} \overline{e}_{\alpha_{j-1}+\alpha_{j}}
      -q^{-1}\overline{e}_{\alpha_{j-1}+\alpha_{j}}e_{\alpha_{j}+2\alpha_{j+1}+\cdots+2\alpha_{n}} \\
& = & [e_{\alpha_{j}+2\alpha_{j+1}+\cdots+2\alpha_{n}},\overline{e}_{\alpha_{j-1}+\alpha_{j}}]_{q} \\
& = & \overline{e}_{\alpha_{j-1} + 2\alpha_{j} + \cdots + 2\alpha_{n}},
\end{eqnarray*}
where the last equality arises from Proposition \ref{appendixB:lem20}.  It follows that
\begin{eqnarray*}
\lefteqn{
\left(\overline{D}+D_{j-1}\right) \Delta(e_{j-1})-q^{-1}\Delta(e_{j-1})\left(\overline{D}+D_{j-1}\right) } \\
& = & \left( q e_{\alpha_{j-1} + 2\alpha_{j} + \cdots + 2 \alpha_{n}} e_{j-1}
              -q^{-1} e_{j-1} e_{\alpha_{j-1} + 2\alpha_{j} + \cdots + 2 \alpha_{n}}\right)
	      \otimes (K_{\alpha_{j-1}+\cdots + \alpha_{n}})^{2} e_{\alpha_{i}+\cdots+\alpha_{j-2}} \\
 & & + \overline{e}_{\alpha_{j-1} + 2\alpha_{j} + \cdots + 2\alpha_{n}} \otimes 
         K_{\alpha_{j-1} + 2\alpha_{j} + \cdots + 2\alpha_{n}}e_{\alpha_{i}+\cdots+\alpha_{j-1}}.
\end{eqnarray*}
\end{proof}

\end{subsection}

\end{section}

\begin{section}{The commutation relations between the components of $\Delta(e_{\mu})$}

The commutation relations between the components of $\Delta(e_{\mu})$ are, in general, 
quite complicated, and
the algebraic structures they give rise to come in three families depending on the nature of $\mu$.  

In this section we prove the following results.
The algebra underlying the components of $\Delta(e_{i})$, where $i=1, \ldots, n$, 
is the $q$-binomial theorem.
The algebra underlying the components of $\Delta(e_{\alpha_{i}+ \cdots + \alpha_{j}})$, where 
$i<j<n$, is the $q$-multinomial theorem.
The algebra underlying the components of $\Delta(e_{\alpha_{i}+ \cdots + \alpha_{n}})$
 is more complicated:
it is a combination of the $q$-binomial theorem and the 
associative algebra given in Lemma \ref{appendixB:generalisationnewtonsbinomial1}.
The algebra underlying the components of $\Delta(e_{\alpha_{i}+\cdots+2\alpha_{j}+\cdots+\alpha_{n}})$ 
is a combination of the $q$-binomial theorem and the associative algebra
given in Lemma \ref{appendixB:generalisationnewtonsbinomial2}.

The $q$-binomial theorem and the two generalisations of the binomial theorem given in Appendix B allow us to compute
any desired power of $\Delta(e_{\mu})$.

\begin{subsection}{$\Delta(e_{i})$, $i=1, \ldots, n$}

Here $\Delta(e_{i}) = e_{i} \otimes K_{i} + 1 \otimes e_{i}$,
and elementary calculations show that 
$$(e_{i} \otimes K_{i})( 1 \otimes e_{i})  =  q^{2}( 1 \otimes e_{i})(e_{i} \otimes K_{i}), \hspace{5mm} i<n,$$
$$(e_{n} \otimes K_{n})( 1 \otimes e_{n}) =  -q( 1 \otimes e_{n})(e_{n} \otimes K_{n}).$$

\end{subsection}

\begin{subsection}{$\Delta(e_{\alpha_{i} + \cdots + \alpha_{j}})$, $1 \leq i<j<n$}

For each $j = i+1, \ldots, n-1$, we write 
$$\Delta(e_{\alpha_{i} + \cdots + \alpha_{j}}) = D_{i} + (q-q^{-1}) \sum_{k=i+1}^{j}D_{k} + D_{\infty},
\hspace{10mm} \mbox{where}$$
\begin{eqnarray*}
D_{i} & = & e_{\alpha_{i} + \cdots + \alpha_{j}} \otimes K_{\alpha_{i} + \cdots + \alpha_{j}}, \\
D_{k} & = & e_{\alpha_{k} + \cdots + \alpha_{j}} \otimes K_{\alpha_{k} + \cdots + \alpha_{j}}
 e_{\alpha_{i} + \cdots + \alpha_{k-1}}, \\
D_{\infty} & = & 1 \otimes e_{\alpha_{i} + \cdots + \alpha_{j}}.
\end{eqnarray*}
We will show that the commutation relations between the components of 
$\Delta(e_{\alpha_{i} + \cdots + \alpha_{j}})$ are
$$D_{r} D_{s} = q^{2} D_{s} D_{r}, \hspace{10mm} \forall r < s.$$
For each $k=i+1, \ldots, j$, the relations between $D_{i}$ and $D_{k}$ are
\begin{eqnarray*}
D_{i}D_{k} 
& = & q  e_{\alpha_{k} + \cdots + \alpha_{j}} e_{\alpha_{i} + \cdots + \alpha_{j}} \otimes 
      K_{\alpha_{i} + \cdots + \alpha_{j}} K_{\alpha_{k} + \cdots + \alpha_{j}} e_{\alpha_{i} + \cdots + \alpha_{k-1}}   \\
& = & q^{2} e_{\alpha_{k} + \cdots + \alpha_{j}} e_{\alpha_{i} + \cdots + \alpha_{j}} \otimes 
      K_{\alpha_{k} + \cdots + \alpha_{j}} e_{\alpha_{i} + \cdots + \alpha_{k-1}}K_{\alpha_{i} + \cdots + \alpha_{j}}   
 = q^{2} D_{k}D_{i},
\end{eqnarray*}
where we have used 
$\left[e_{\alpha_{i} + \cdots + \alpha_{j}},e_{\alpha_{k} + \cdots + \alpha_{j}}\right]_{q}=0$. 
The relation between $D_{i}$ and $D_{\infty}$ is
$$
D_{i}D_{\infty} 
= q^{2} e_{\alpha_{i} + \cdots + \alpha_{j}} \otimes 
      e_{\alpha_{i} + \cdots + \alpha_{j}}K_{\alpha_{i} + \cdots + \alpha_{j}}   
 =  q^{2} D_{\infty}D_{i}.
$$
Now for each $k = i+1, \ldots, j-1$, and each $l$ satisfying $k< l \leq j$, 
the relation between $D_{k}$ and $D_{l}$ is
\begin{eqnarray*}
D_{k} D_{l} 
& = & q e_{\alpha_{l} + \cdots + \alpha_{j}} e_{\alpha_{k} + \cdots + \alpha_{j}} \otimes
K_{\alpha_{k} + \cdots + \alpha_{j}} e_{\alpha_{i} + \cdots + \alpha_{k-1}}
K_{\alpha_{l} + \cdots + \alpha_{j}} e_{\alpha_{i} + \cdots + \alpha_{l-1}} \\ 
& = & q^{2} e_{\alpha_{l} + \cdots + \alpha_{j}} e_{\alpha_{k} + \cdots + \alpha_{j}} \otimes
K_{\alpha_{l} + \cdots + \alpha_{j}}K_{\alpha_{k} + \cdots + \alpha_{j}}e_{\alpha_{i} + \cdots + \alpha_{l-1}}
e_{\alpha_{i} + \cdots + \alpha_{k-1}} \\
& = & q^{2} D_{l} D_{k}.
\end{eqnarray*}
Finally, the relation between $D_{k}$ and $D_{\infty}$, for each $k=i+1,\ldots, j$, is
\begin{eqnarray*}
D_{k}D_{\infty} 
 & = &  e_{\alpha_{k} + \cdots + \alpha_{j}} \otimes K_{\alpha_{k} + \cdots + \alpha_{j}}
 e_{\alpha_{i} + \cdots + \alpha_{k-1}}e_{\alpha_{i} + \cdots + \alpha_{j}} \\
 & = & q^{2} e_{\alpha_{k} + \cdots + \alpha_{j}} \otimes e_{\alpha_{i} + \cdots + \alpha_{j}}
 K_{\alpha_{k} + \cdots + \alpha_{j}}e_{\alpha_{i} + \cdots + \alpha_{k-1}} 
 = q^{2} D_{\infty}D_{k}.
\end{eqnarray*}
These calculations give the flavour of many of the calculations in this section.

\end{subsection}

\begin{subsection}{$\Delta\left( e_{\alpha_{i} + \cdots + \alpha_{n}}\right)$, $i=1, \ldots, n-1$}

This case is more difficult than the previous one.
We write 
$$\Delta ( e_{\alpha_{i} + \cdots + \alpha_{n}}) = D_{i} + \sum_{k=i+1}^{n} D_{k} + D_{\infty},
\hspace{10mm} \mbox{where}$$
\begin{eqnarray*}
D_{i}      & = & e_{\alpha_{i} + \cdots + \alpha_{n}} \otimes K_{\alpha_{i} + \cdots + \alpha_{n}}, \\
D_{k}      & = & (q-q^{-1})\left(e_{\alpha_{k} + \cdots + \alpha_{n}} \otimes K_{\alpha_{k} + \cdots + \alpha_{n}}
       		 e_{\alpha_{i} + \cdots + \alpha_{k-1}}\right),  \hspace{10mm} k=i+1, \ldots, n, \\
D_{\infty} & = &  1 \otimes e_{\alpha_{i} + \cdots + \alpha_{n}}.  
\end{eqnarray*}

We will show that the commutation relations between the components of 
$\Delta (e_{\alpha_{i} + \cdots + \alpha_{n}} )$ are  
(\ref{appendixB:abstralgrel(1)})--(\ref{appendixB:abstralgrel(666)}) below, 
where we fix $i \leq p, r, s, t \leq n$ throughout these equations,  
$r<t$ in (\ref{appendixB:abstralgrel(1)})--(\ref{appendixB:abstralgrel(6)}), 
$r<s<t$ in (\ref{appendixB:abstralgrel(4)}), 
and $p<r<t$ in (\ref{appendixB:abstralgrel(5)})--(\ref{appendixB:abstralgrel(6)}):
\begin{eqnarray}
D_{r} D_{t}    & = & -q D_{t} D_{r} + E_{r,t},         \label{appendixB:abstralgrel(1)} \\
D_{r} E_{r,t}  & = & q^{2} E_{r,t} D_{r},              \label{appendixB:abstralgrel(2)} \\
E_{r,t} D_{t}  & = & q^{2} D_{t} E_{r,t},              \label{appendixB:abstralgrel(3)} \\
D_{s} E_{r,t}  & = & E_{r,t} D_{s},                    \label{appendixB:abstralgrel(4)} \\
D_{p} E_{r,t}  & = & q^{2} E_{r,t} D_{p} + F_{p,r,t},  \label{appendixB:abstralgrel(5)} \\
0 & = & F_{p,r,t} + (1-q^{2})E_{p,t}D_{r},             \label{appendixB:abstralgrel(6)} \\
D_{r} D_{\infty} & = & -q D_{\infty}D_{r}.             \label{appendixB:abstralgrel(666)}
\end{eqnarray}
Note that relations (\ref{appendixB:abstralgrel(1)})--(\ref{appendixB:abstralgrel(3)}) 
imply the following relations
\begin{equation}
\label{appendixB:abstralgrel(7)}
D_{p} E_{r,t} - q^{2}E_{r,t}D_{p} = E_{p,r}D_{t}-q^{2}D_{t}E_{p,r}, \hspace{10mm}  p < r < t,
\end{equation}
\begin{equation}
\label{appendixB:abstralgrel(8)}
E_{p,r} E_{p,t} = q^{2} E_{p,t} E_{p,r}, \hspace{10mm} E_{p,t}E_{r,t} = q^{2} E_{r,t} E_{p,t}, 
\hspace{10mm}  p < r < t.
\end{equation}

The following proposition shows that if the components of 
$\Delta( e_{\alpha_{i} + \cdots + \alpha_{n}})$ satisfy 
(\ref{appendixB:abstralgrel(1)})--(\ref{appendixB:abstralgrel(666)}), then we can use
Lemma \ref{appendixB:generalisationnewtonsbinomial1}
to give an expansion of $\left( D_{i} + D_{i+1} + \cdots + D_{n}\right)^{m}$ for each $m \in \mathbb{N}$.
Note that relations (\ref{appendixB:abstralgrel(9)})--(\ref{appendixB:abstralgrel(11)}) in 
Proposition \ref{appendixC:prop:gigwong(1)}
are precisely the relations between the elements $a,b,c$, respectively, 
in the associative algebra given in
Lemma \ref{appendixB:generalisationnewtonsbinomial1}.

\begin{proposition}  
\label{appendixC:prop:gigwong(1)}
For each $k=i+1, \ldots, n$, and each $j$ satisfying $i \leq j < k$, we have
\begin{eqnarray}
\lefteqn{
 (D_{j} + D_{j+1} + \cdots + D_{k-1})D_{k} } \nonumber \\
 &  & = -q D_{k}(D_{j} + D_{j+1} + \cdots + D_{k-1}) 
                       + (E_{j,k} + E_{j+1,k} + \cdots + E_{k-1,k}), \label{appendixB:abstralgrel(9)}
\end{eqnarray}
where 
\begin{equation}
\label{appendixB:abstralgrel(10)}
\left( \sum_{b=j}^{k-1} E_{b,k} \right) D_{k} = q^{2} D_{k}\left( \sum_{b=j}^{k-1} E_{b,k} \right),
\end{equation}
and
\begin{equation}
\label{appendixB:abstralgrel(11)}
\left( \sum_{a=j}^{k-1} D_{a} \right) \left( \sum_{b=j}^{k-1} E_{b,k} \right)
   = q^{2}  \left( \sum_{b=j}^{k-1} E_{b,k} \right) \left( \sum_{a=j}^{k-1} D_{a} \right).
\end{equation}
\end{proposition}
\begin{proof}
Relation (\ref{appendixB:abstralgrel(9)}) follows from (\ref{appendixB:abstralgrel(1)}) and
(\ref{appendixB:abstralgrel(10)}) follows from (\ref{appendixB:abstralgrel(3)}).
The proof of (\ref{appendixB:abstralgrel(11)}) is only slightly more difficult:
firstly, the elements $D_{k-1}$ and $D_{k}$ satisfy the relations
$$D_{k-1} D_{k} = -q D_{k} D_{k-1} + E_{k-1,k}, \hspace{5mm} D_{k-1} E_{k-1,k} = q^{2} E_{k-1,k} D_{k-1},
\hspace{5mm} E_{k-1,k} D_{k} = q^{2} D_{k} E_{k-1,k},$$
which proves (\ref{appendixB:abstralgrel(11)}) for $j=k-1$.
Now assume that (\ref{appendixB:abstralgrel(11)}) is true for some $j \leq k-1$, we will prove that 
(\ref{appendixB:abstralgrel(11)}) is true for $j-1$ if $j-1 \geq i$:
\begin{eqnarray}
\lefteqn{
\left( \sum_{a=j-1}^{k-1} D_{a} \right) \left( \sum_{b=j-1}^{k-1} E_{b,k} \right) } \nonumber \\
& = & \left[ D_{j-1} + \left( \sum_{a=j}^{k-1} D_{a} \right) \right] 
\left[ E_{j-1,k} + \left( \sum_{b=j}^{k-1} E_{b,k} \right) \right] \nonumber \\
& = & D_{j-1} E_{j-1,k} + \left( \sum_{a=j}^{k-1} D_{a} \right)\left( \sum_{b=j}^{k-1} E_{b,k} \right)
      + D_{j-1} \left( \sum_{c=j}^{k-1} E_{c,k} \right) + \left( \sum_{d=j}^{k-1} D_{d} \right) E_{j-1,k} \nonumber \\
& = & q^{2} E_{j-1,k} D_{j-1} + q^{2} \left( \sum_{b=j}^{k-1} E_{b,k} \right)\left( \sum_{a=j}^{k-1} D_{a} \right)
+ \sum_{c=j}^{k-1} \left( q^{2} E_{c,k}D_{j-1} + F_{j-1,c,k} \right) \nonumber \\
& & + q^{2} E_{j-1,k} \left( \sum_{d=j}^{k-1} D_{d} \right)
    + (1-q^{2}) E_{j-1,k} \left( \sum_{e=j}^{k-1} D_{e} \right) \nonumber \\
& = & q^{2} \left(\sum_{c=j-1}^{k-1} E_{c,k}\right) D_{j-1} + q^{2} E_{j-1,k} \left( \sum_{a=j}^{k-1} D_{a} \right)
      + q^{2} \left( \sum_{b=j}^{k-1} E_{b,k} \right)\left( \sum_{a=j}^{k-1} D_{a} \right) \label{appendixC:simoneyoung(10)} \\
& = & q^{2} \left(\sum_{b=j-1}^{k-1} E_{b,k}\right) \left( \sum_{a=j-1}^{k-1} D_{a} \right), \nonumber 
\end{eqnarray}
where  (\ref{appendixC:simoneyoung(10)})  follows from (\ref{appendixB:abstralgrel(6)}), completing the induction.
\end{proof}

We now show that the components of $\Delta(e_{\alpha_{i} + \cdots + \alpha_{n}})$ satisfy the claimed
commutation relations.  We firstly prove relation (\ref{appendixB:abstralgrel(666)}): 
\begin{eqnarray*}
D_{i} D_{\infty} 
		 & = & e_{\alpha_{i} + \cdots +\alpha_{n}} \otimes
		 K_{\alpha_{i} + \cdots +\alpha_{n}}e_{\alpha_{i} + \cdots +\alpha_{n}} \\
		 & = & q e_{\alpha_{i} + \cdots +\alpha_{n}} \otimes
		 e_{\alpha_{i} + \cdots +\alpha_{n}}K_{\alpha_{i} + \cdots +\alpha_{n}} 
		 = -q D_{\infty} D_{i},
\end{eqnarray*}
and for each $k=i+1, \ldots, n$, the relation between $D_{k}$ and $D_{\infty}$ is
\begin{eqnarray*}
D_{k} D_{\infty} 
& = &   (q-q^{-1}) \left( e_{\alpha_{k} + \cdots +\alpha_{n}}   
\otimes K_{\alpha_{k} + \cdots +\alpha_{n}} e_{\alpha_{i} + \cdots +\alpha_{k-1}}e_{\alpha_{i} + \cdots +\alpha_{n}} \right)
\\
& = & q (q-q^{-1})\left( e_{\alpha_{k} + \cdots +\alpha_{n}}   
\otimes e_{\alpha_{i} + \cdots +\alpha_{n}} 
K_{\alpha_{k} + \cdots +\alpha_{n}} e_{\alpha_{i} + \cdots +\alpha_{k-1}} \right) 
= -q D_{\infty}D_{k},
\end{eqnarray*}
proving relation (\ref{appendixB:abstralgrel(666)}).  

We now prove relation (\ref{appendixB:abstralgrel(1)}).
For each $k = i+1, \ldots, n$, the relation between $D_{i}$ and $D_{k}$ is
\begin{eqnarray}
D_{i} D_{k} 
	    & = & (q-q^{-1})\left(e_{\alpha_{i} + \cdots + \alpha_{n}} e_{\alpha_{k} + \cdots + \alpha_{n}} \otimes
	           K_{\alpha_{i} + \cdots + \alpha_{n}}K_{\alpha_{k} + \cdots + \alpha_{n}}
		   e_{\alpha_{i} + \cdots + \alpha_{k-1}}\right)  \nonumber \\
	    & = & (q-q^{-1})\left( -e_{\alpha_{k} + \cdots + \alpha_{n}}e_{\alpha_{i} + \cdots + \alpha_{n}} 
	         + \left[ e_{\alpha_{i} + \cdots + \alpha_{n}},e_{\alpha_{k} + \cdots + \alpha_{n}}\right]_{q} \right) \nonumber \\
	    &  & \hspace{50mm} 
		   \otimes K_{\alpha_{i} + \cdots + \alpha_{n}}K_{\alpha_{k} + \cdots + \alpha_{n}}
		   e_{\alpha_{i} + \cdots + \alpha_{k-1}}  \nonumber \\
	    & = & (q-q^{-1})
	           \Big(-qe_{\alpha_{k} + \cdots + \alpha_{n}}e_{\alpha_{i} + \cdots + \alpha_{n}}
	          \otimes K_{\alpha_{k} + \cdots + \alpha_{n}} e_{\alpha_{i} + \cdots + \alpha_{k-1}} 
		  K_{\alpha_{i} + \cdots + \alpha_{n}} \label{appendixB:normanlamont(1)} \\
	    &   &  \hspace{20mm} + \left[ e_{\alpha_{i} + \cdots + \alpha_{n}},e_{\alpha_{k} + \cdots + \alpha_{n}}\right]_{q}
	           \otimes K_{\alpha_{i} + \cdots + \alpha_{n}}K_{\alpha_{k} + \cdots + \alpha_{n}}
		   e_{\alpha_{i} + \cdots + \alpha_{k-1}} \Big)  \nonumber \\
	& = & -q D_{k} D_{i} + E_{i,k}, \nonumber 
\end{eqnarray}
where we have used the relation 
$\left[e_{\alpha_{i} + \cdots + \alpha_{n}},e_{\alpha_{k} + \cdots + \alpha_{n}}\right]_{q}=0$ 
to obtain (\ref{appendixB:normanlamont(1)}), and where we set 
$$E_{i,k} = (q-q^{-1})\left(\left[ e_{\alpha_{i} + \cdots + \alpha_{n}},e_{\alpha_{k} + \cdots + \alpha_{n}}\right]_{q}
	           \otimes K_{\alpha_{i} + \cdots + \alpha_{n}}K_{\alpha_{k} + \cdots + \alpha_{n}}
		   e_{\alpha_{i} + \cdots + \alpha_{k-1}}\right).$$
For all $k, p$ satisfying $i+1 \leq k < p \leq n$, the relation between $D_{k}$ and $D_{p}$ is
\begin{eqnarray}
D_{k} D_{p} 
            & = & (q-q^{-1})^{2}\left(e_{\alpha_{k} + \cdots + \alpha_{n}}e_{\alpha_{p} + \cdots + \alpha_{n}} \otimes
	          K_{\alpha_{k} + \cdots + \alpha_{n}}e_{\alpha_{i} + \cdots + \alpha_{k-1}}
		  K_{\alpha_{p} + \cdots + \alpha_{n}}e_{\alpha_{i} + \cdots + \alpha_{p-1}}\right) \nonumber \\
	    & = & (q-q^{-1})^{2}\left(-e_{\alpha_{p} + \cdots + \alpha_{n}}e_{\alpha_{k} + \cdots + \alpha_{n}}
	          + \left[ e_{\alpha_{k} + \cdots + \alpha_{n}},e_{\alpha_{p} + \cdots + \alpha_{n}}\right]_{q}\right) \nonumber \\
            & & \hspace{40mm}
		  \otimes K_{\alpha_{k} + \cdots + \alpha_{n}}K_{\alpha_{p} + \cdots + \alpha_{n}}
		  e_{\alpha_{i} + \cdots + \alpha_{k-1}}e_{\alpha_{i} + \cdots + \alpha_{p-1}} \nonumber \\
            & = & -q (q-q^{-1})^{2} \left(
	        e_{\alpha_{p} + \cdots + \alpha_{n}}e_{\alpha_{k} + \cdots + \alpha_{n}}
	         \otimes K_{\alpha_{k} + \cdots + \alpha_{n}}K_{\alpha_{p} + \cdots + \alpha_{n}}
		  e_{\alpha_{i} + \cdots + \alpha_{p-1}}e_{\alpha_{i} + \cdots + \alpha_{k-1}} \right)
		  \nonumber
		   \\
	    &   &  + (q-q^{-1})^{2}
		   \left( \left[ e_{\alpha_{k} + \cdots + \alpha_{n}},e_{\alpha_{p} + \cdots + \alpha_{n}}\right]_{q}
	          \otimes K_{\alpha_{k} + \cdots + \alpha_{n}}K_{\alpha_{p} + \cdots + \alpha_{n}}
		  e_{\alpha_{i} + \cdots + \alpha_{k-1}}e_{\alpha_{i} + \cdots + \alpha_{p-1}} \right) \nonumber \\
		  & & \label{appendixB:charleswooley(1)} \\
	    & = & -q D_{p} D_{k} + E_{k,p}, \nonumber
\end{eqnarray}
where we have used the relation
$\left[e_{\alpha_{i} + \cdots + \alpha_{k-1}},e_{\alpha_{i} + \cdots + \alpha_{p-1}}\right]_{q}=0$ 
to obtain (\ref{appendixB:charleswooley(1)}), and we set
$$E_{k,p} = (q-q^{-1})^{2}
            \left( \left[ e_{\alpha_{k} + \cdots + \alpha_{n}},e_{\alpha_{p} + \cdots + \alpha_{n}}\right]_{q}
	          \otimes K_{\alpha_{k} + \cdots + \alpha_{n}}K_{\alpha_{p} + \cdots + \alpha_{n}}
		  e_{\alpha_{i} + \cdots + \alpha_{k-1}}e_{\alpha_{i} + \cdots + \alpha_{p-1}}\right).$$
It follows that the components of 
$\Delta( e_{\alpha_{i} + \cdots + \alpha_{n}})$ satisfy relation (\ref{appendixB:abstralgrel(1)}).

We now consider the relations between the $D_{k}$ and $E_{k,p}$: 
we will show that $D_{k}$ and $E_{k,p}$ satisfy relation (\ref{appendixB:abstralgrel(2)}) for all $k<p$.
For calculational ease we will write $E_{i,p}$ and $E_{k,p}$ as the following sums: 
$$E_{i,p} = \sum_{j=p}^{n} E_{i,p}^{j}, \hspace{20mm} 
  E_{k,p} = \sum_{j=p}^{n} E_{k,p}^{j}, \hspace{10mm} k>i,$$
where the components of these sums are given in Proposition \ref{appendixB:germany(1)}:
\begin{eqnarray*}
E_{i,p}^{p} & = & (q-q^{-1}) (-q)^{n-p} 
              \left( e_{\alpha_{i}+\cdots+2\alpha_{p}+\cdots+2\alpha_{n}} \otimes 
              K_{\alpha_{i} + \cdots + \alpha_{n}}K_{\alpha_{p} + \cdots + \alpha_{n}} 
	      e_{\alpha_{i} + \cdots + \alpha_{p-1}}\right), \\
E_{i,p}^{j} & = &  (q-q^{-1})^{2}(-q)^{n-j} 
               \big( e_{\alpha_{i}+\cdots+2\alpha_{j}+\cdots+2\alpha_{n}} e_{\alpha_{p}+\cdots+\alpha_{j-1}} \\
& &  \hspace{50mm}        \otimes 
              K_{\alpha_{i} + \cdots + \alpha_{n}}K_{\alpha_{p} + \cdots + \alpha_{n}} 
	      e_{\alpha_{i} + \cdots + \alpha_{p-1}}\big), \hspace{20mm} j>p, \\
E_{k,p}^{p} & = & (q-q^{-1})^{2} (-q)^{n-p} 
		   \left( e_{\alpha_{k}+\cdots+2\alpha_{p}+\cdots+2\alpha_{n}} \otimes
                     K_{\alpha_{k} + \cdots + \alpha_{n}}K_{\alpha_{p} + \cdots + \alpha_{n}} 
	             e_{\alpha_{i} + \cdots + \alpha_{k-1}}e_{\alpha_{i} + \cdots + \alpha_{p-1}} \right), \\
E_{k,p}^{j} & = & (q-q^{-1})^{3} (-q)^{n-j}
                   \big( e_{\alpha_{k}+\cdots+2\alpha_{j}+\cdots+2\alpha_{n}}e_{\alpha_{p}+\cdots+\alpha_{j-1}} \\
 & & \hspace{40mm} \otimes
                     K_{\alpha_{k} + \cdots + \alpha_{n}}K_{\alpha_{p} + \cdots + \alpha_{n}}
		     e_{\alpha_{i} + \cdots + \alpha_{k-1}}e_{\alpha_{i} + \cdots + \alpha_{p-1}} \big). 
		     \hspace{12mm} j>p.
\end{eqnarray*}	
The following identities allow us to determine the relations between $D_{i}$ and $E_{i,p}$, 
and between $D_{k}$ and $E_{k,p}$, quite easily:
\begin{eqnarray}
\left[e_{\alpha_{i} + \cdots + \alpha_{n}}, e_{\alpha_{i}+\cdots+2\alpha_{j}+\cdots+2\alpha_{n}} \right]_{q} & = & 0, 
\hspace{10mm}   1 \leq i \leq n-1, \ \ i+1 \leq j \leq n, \label{appendixB:identities1000} \\
\left[ e_{\alpha_{i} + \cdots + \alpha_{n}},  e_{\alpha_{p}+\cdots+\alpha_{j-1}} \right]_{q} & = & 0,
\hspace{10mm}   1 \leq i \leq n-2, \ \ i<p<j-1<n. \label{appendixB:identities1001}
\end{eqnarray}
We can rewrite (\ref{appendixB:identities1000}) and (\ref{appendixB:identities1001}), respectively, as 
\begin{eqnarray*}
e_{\alpha_{i} + \cdots + \alpha_{n}}e_{\alpha_{i}+\cdots+2\alpha_{j}+\cdots+2\alpha_{n}}
& = & q e_{\alpha_{i}+\cdots+2\alpha_{j}+\cdots+2\alpha_{n}}e_{\alpha_{i} + \cdots + \alpha_{n}}, \\
e_{\alpha_{i} + \cdots + \alpha_{n}} e_{\alpha_{p}+\cdots+\alpha_{j-1}}
&= & e_{\alpha_{p}+\cdots+\alpha_{j-1}}e_{\alpha_{i} + \cdots + \alpha_{n}}.
\end{eqnarray*}
Using these identities we compute the relations between $D_{i}$ and $E_{i,p}$ 
where $i=1, \ldots, n-1$ and $p > i$:
\begin{eqnarray*}
\lefteqn{
D_{i}E_{i,p}^{p} } \\
& = & (q-q^{-1}) (-q)^{n-p} \\
& & \hspace{5mm} \times \left( e_{\alpha_{i} + \cdots + \alpha_{n}}e_{\alpha_{i}+\cdots+2\alpha_{p}+\cdots+2\alpha_{n}} \otimes
	      K_{\alpha_{i} + \cdots + \alpha_{n}} K_{\alpha_{i} + \cdots + \alpha_{n}}K_{\alpha_{p} + \cdots + \alpha_{n}} 
	      e_{\alpha_{i} + \cdots + \alpha_{p-1}} \right) \\
 & = & q^{2} (q-q^{-1}) (-q)^{n-p} \\
 & & \hspace{5mm} \times 
 \left( e_{\alpha_{i}+\cdots+2\alpha_{p}+\cdots+2\alpha_{n}}e_{\alpha_{i} + \cdots + \alpha_{n}}\otimes
             K_{\alpha_{i} + \cdots + \alpha_{n}}K_{\alpha_{p} + \cdots + \alpha_{n}}e_{\alpha_{i} + \cdots + \alpha_{p-1}}
	     K_{\alpha_{i} + \cdots + \alpha_{n}} \right) \\
 & = & q^{2} E_{i,p}^{p} D_{i},
\end{eqnarray*}
and for $p>i$ and $j > p$ we have
\begin{eqnarray*}
\lefteqn{
D_{i}E_{i,p}^{j} } \\ 
 & = & (q-q^{-1})^{2}(-q)^{n-j} \\
 & & \hspace{5mm} \times 
 \left(
 e_{\alpha_{i} + \cdots + \alpha_{n}}e_{\alpha_{i}+\cdots+2\alpha_{j}+\cdots+2\alpha_{n}} e_{\alpha_{p}+\cdots+\alpha_{j-1}}
  \otimes K_{\alpha_{i} + \cdots + \alpha_{n}}K_{\alpha_{i} + \cdots + \alpha_{n}}K_{\alpha_{p} + \cdots + \alpha_{n}} 
	      e_{\alpha_{i} + \cdots + \alpha_{p-1}} \right) \\
 & = & q^{2} (q-q^{-1})^{2}(-q)^{n-j} \\
 & & \hspace{5mm} \times  
 \left(
 e_{\alpha_{i}+\cdots+2\alpha_{j}+\cdots+2\alpha_{n}} e_{\alpha_{p}+\cdots+\alpha_{j-1}} e_{\alpha_{i} + \cdots + \alpha_{n}}
 \otimes K_{\alpha_{i} + \cdots + \alpha_{n}}K_{\alpha_{p} + \cdots + \alpha_{n}} e_{\alpha_{i} + \cdots + \alpha_{p-1}}
 K_{\alpha_{i} + \cdots + \alpha_{n}} \right) \\
 & = & q^{2} E_{i,p}^{j}D_{i}.
\end{eqnarray*}
These calculations show that $D_{i} E_{i,p} = q^{2} E_{i,p} D_{i}$ for $p = i+1, \ldots, n$.
Now we determine the relations between $D_{k}$ and the components of $E_{k,p}$.
For $k = i+1, \ldots, n-1$ and $p>k$ we have
\begin{eqnarray*}
\lefteqn{
D_{k}E_{k,p}^{p} } \\
 & = & (q-q^{-1})^{3}(-q)^{n-p}   \\
 & & \hspace{5mm} \times 
\big( e_{\alpha_{k} + \cdots + \alpha_{n}} e_{\alpha_{k}+\cdots+2\alpha_{p}+\cdots+2\alpha_{n}} \\
 & & \hspace{15mm} \otimes  K_{\alpha_{k} + \cdots + \alpha_{n}} e_{\alpha_{i} + \cdots + \alpha_{k-1}}
 K_{\alpha_{k} + \cdots + \alpha_{n}}K_{\alpha_{p} + \cdots + \alpha_{n}} 
	             e_{\alpha_{i} + \cdots + \alpha_{k-1}}e_{\alpha_{i} + \cdots + \alpha_{p-1}} \big) \\
& = & q^{2}(q-q^{-1})^{3}(-q)^{n-p}   \\
& & \hspace{5mm} \times \big( 
e_{\alpha_{k}+\cdots+2\alpha_{p}+\cdots+2\alpha_{n}} e_{\alpha_{k} + \cdots + \alpha_{n}} \\
& & \hspace{15mm} \otimes  K_{\alpha_{k} + \cdots + \alpha_{n}}K_{\alpha_{p} + \cdots + \alpha_{n}} 
	             e_{\alpha_{i} + \cdots + \alpha_{k-1}}e_{\alpha_{i} + \cdots + \alpha_{p-1}}
		      K_{\alpha_{k} + \cdots + \alpha_{n}} e_{\alpha_{i} + \cdots + \alpha_{k-1}} \big) \\
& = & q^{2}E_{k,p}^{p} D_{k},
\end{eqnarray*}
and for $j > p$, 
\begin{eqnarray*}
\lefteqn{
D_{k}E_{k,p}^{j} } \\
& = & (q-q^{-1})^{4} (-q)^{n-j}  \\
& & \hspace{5mm} \times  
\big( e_{\alpha_{k} + \cdots + \alpha_{n}}e_{\alpha_{k}+\cdots+2\alpha_{j}+\cdots+2\alpha_{n}}e_{\alpha_{p}+\cdots+\alpha_{j-1}} \\
& & \hspace{15mm} \otimes K_{\alpha_{k} + \cdots + \alpha_{n}}e_{\alpha_{i} + \cdots + \alpha_{k-1}}
           K_{\alpha_{k} + \cdots + \alpha_{n}}K_{\alpha_{p} + \cdots + \alpha_{n}}
		     e_{\alpha_{i} + \cdots + \alpha_{k-1}}e_{\alpha_{i} + \cdots + \alpha_{p-1}} \big) \\
& = & q^{2}(q-q^{-1})^{4} (-q)^{n-j}  \\
& & \hspace{5mm} \times
\big( 
e_{\alpha_{k}+\cdots+2\alpha_{j}+\cdots+2\alpha_{n}}e_{\alpha_{p}+\cdots+\alpha_{j-1}}e_{\alpha_{k} + \cdots + \alpha_{n}} \\
& & \hspace{15mm} \otimes K_{\alpha_{k} + \cdots + \alpha_{n}}K_{\alpha_{p} + \cdots + \alpha_{n}}
		     e_{\alpha_{i} + \cdots + \alpha_{k-1}}e_{\alpha_{i} + \cdots + \alpha_{p-1}}
		     K_{\alpha_{k} + \cdots + \alpha_{n}}e_{\alpha_{i} + \cdots + \alpha_{k-1}} \big) \\
& = & q^{2}E_{k,p}^{j}D_{k}.
\end{eqnarray*}
This shows that $D_{k} E_{k,p} = q^{2} E_{k,p} D_{k}$ for $p = k+1, \ldots, n$.
Together with $D_{i} E_{i,p} = q^{2} E_{i,p}D_{i}$, this proves that
 $D_{k}$ and $E_{k,p}$  satisfy relation (\ref{appendixB:abstralgrel(2)}).

We  now show that $D_{k}$ and $E_{k,p}$ satisfy relation (\ref{appendixB:abstralgrel(3)}).
The identity
$$\left[e_{\alpha_{i}+\cdots+2\alpha_{p}+\cdots+2\alpha_{n}},e_{\alpha_{p} + \cdots + \alpha_{n}}\right]_{q}=0,
\hspace{5mm}  i+1 \leq p \leq n,$$
which we can rewrite as
$e_{\alpha_{i}+\cdots+2\alpha_{p}+\cdots+2\alpha_{n}} e_{\alpha_{p} + \cdots + \alpha_{n}}
 = q e_{\alpha_{p} + \cdots + \alpha_{n}} e_{\alpha_{i}+\cdots+2\alpha_{p}+\cdots+2\alpha_{n}}$,
 will be useful.
The following calculations are similar to those immediately above.
For each $k=i+1, \ldots, n-1$ and $p>k$, 
\begin{eqnarray*}
E_{i,p}^{p} D_{p} 
 & = & (q-q^{-1})^{2} (-q)^{n-p} \\
 & & \hspace{5mm} \times \big(
  e_{\alpha_{i}+\cdots+2\alpha_{p}+\cdots+2\alpha_{n}}e_{\alpha_{p} + \cdots + \alpha_{n}}  \\
  & & \hspace{15mm} \otimes
     K_{\alpha_{i} + \cdots + \alpha_{n}}K_{\alpha_{p} + \cdots + \alpha_{n}} 
	      e_{\alpha_{i} + \cdots + \alpha_{p-1}}K_{\alpha_{p} + \cdots + \alpha_{n}}
       		 e_{\alpha_{i} + \cdots + \alpha_{p-1}} \big) \\
 & = & q (q-q^{-1})^{2} (-q)^{n-p} \\
  & & \hspace{5mm} \times \big(
  e_{\alpha_{p} + \cdots + \alpha_{n}}e_{\alpha_{i}+\cdots+2\alpha_{p}+\cdots+2\alpha_{n}}  \\
  & & \hspace{15mm} \otimes
     K_{\alpha_{i} + \cdots + \alpha_{n}}K_{\alpha_{p} + \cdots + \alpha_{n}} 
	      e_{\alpha_{i} + \cdots + \alpha_{p-1}}K_{\alpha_{p} + \cdots + \alpha_{n}}
       		 e_{\alpha_{i} + \cdots + \alpha_{p-1}} \big) \\
 & = & q^{2} (q-q^{-1})^{2} (-q)^{n-p} \\
   & & \hspace{5mm} \times \big(
  e_{\alpha_{p} + \cdots + \alpha_{n}}e_{\alpha_{i}+\cdots+2\alpha_{p}+\cdots+2\alpha_{n}}  \\
  & & \hspace{15mm} \otimes
     K_{\alpha_{p} + \cdots + \alpha_{n}} e_{\alpha_{i} + \cdots + \alpha_{p-1}}
     K_{\alpha_{i} + \cdots + \alpha_{n}}K_{\alpha_{p} + \cdots + \alpha_{n}}
       		 e_{\alpha_{i} + \cdots + \alpha_{p-1}} \big) \\
& = & q^{2} D_{p} E_{i,p}^{p}.
\end{eqnarray*}
For each $j=p+1,\ldots, n$, we have
$\left[e_{\alpha_{i}+\cdots+2\alpha_{j}+\cdots+2\alpha_{n}},
e_{\alpha_{p} + \cdots + \alpha_{n}}\right]_{q}=0$, and 
\begin{eqnarray*}
E_{i,p}^{j} D_{p} 
& = & (q-q^{-1})^{3}(-q)^{n-j}  \\
& & \hspace{5mm} \times e_{\alpha_{i}+\cdots+2\alpha_{j}+\cdots+2\alpha_{n}} e_{\alpha_{p}+\cdots+\alpha_{j-1}}
                     e_{\alpha_{p} + \cdots + \alpha_{n}} \\
& & \hspace{15mm}     \otimes K_{\alpha_{i} + \cdots + \alpha_{n}}K_{\alpha_{p} + \cdots + \alpha_{n}} 
	      e_{\alpha_{i} + \cdots + \alpha_{p-1}}K_{\alpha_{p} + \cdots + \alpha_{n}}
       		 e_{\alpha_{i} + \cdots + \alpha_{p-1}} \\
& = & q (q-q^{-1})^{3}(-q)^{n-j}  \\
& & \hspace{5mm} \times e_{\alpha_{p} + \cdots + \alpha_{n}}
                  e_{\alpha_{i}+\cdots+2\alpha_{j}+\cdots+2\alpha_{n}} e_{\alpha_{p}+\cdots+\alpha_{j-1}} \\
& & \hspace{15mm}     \otimes K_{\alpha_{i} + \cdots + \alpha_{n}}K_{\alpha_{p} + \cdots + \alpha_{n}} 
	      e_{\alpha_{i} + \cdots + \alpha_{p-1}}K_{\alpha_{p} + \cdots + \alpha_{n}}
       		 e_{\alpha_{i} + \cdots + \alpha_{p-1}} \\
& = & q^{2} (q-q^{-1})^{3}(-q)^{n-j}  \\
& & \hspace{5mm} \times e_{\alpha_{p} + \cdots + \alpha_{n}}
                  e_{\alpha_{i}+\cdots+2\alpha_{j}+\cdots+2\alpha_{n}} e_{\alpha_{p}+\cdots+\alpha_{j-1}} \\
& & \hspace{15mm}     \otimes K_{\alpha_{p} + \cdots + \alpha_{n}} e_{\alpha_{i} + \cdots + \alpha_{p-1}}
	      K_{\alpha_{i} + \cdots + \alpha_{n}}K_{\alpha_{p} + \cdots + \alpha_{n}}
       		 e_{\alpha_{i} + \cdots + \alpha_{p-1}} \\
& = & q^{2} D_{p} E_{i,p}^{j}.
\end{eqnarray*}
For each $k = i+1, \ldots, n-1$ and $p>k$,
\begin{eqnarray*}
E_{k,p}^{p} D_{p} 
	& = & (q-q^{-1})^{3} (-q)^{n-p} \\
	& & \hspace{5mm} \times
	     e_{\alpha_{k}+\cdots+2\alpha_{p}+\cdots+2\alpha_{n}}e_{\alpha_{p} + \cdots + \alpha_{n}}
	      \\
	 & & \hspace{15mm} \otimes K_{\alpha_{k} + \cdots + \alpha_{n}}K_{\alpha_{p} + \cdots + \alpha_{n}} 
	             e_{\alpha_{i} + \cdots + \alpha_{k-1}}e_{\alpha_{i} + \cdots + \alpha_{p-1}}
		     K_{\alpha_{p} + \cdots + \alpha_{n}}
       		 e_{\alpha_{i} + \cdots + \alpha_{p-1}} \\
	& = & q (q-q^{-1})^{3} (-q)^{n-p} \\
	& & \hspace{5mm} \times
	e_{\alpha_{p} + \cdots + \alpha_{n}} e_{\alpha_{k}+\cdots+2\alpha_{p}+\cdots+2\alpha_{n}} \\
	& & \hspace{15mm} \otimes K_{\alpha_{k} + \cdots + \alpha_{n}}K_{\alpha_{p} + \cdots + \alpha_{n}} 
	             e_{\alpha_{i} + \cdots + \alpha_{k-1}}e_{\alpha_{i} + \cdots + \alpha_{p-1}}
		     K_{\alpha_{p} + \cdots + \alpha_{n}}
       		 e_{\alpha_{i} + \cdots + \alpha_{p-1}}  \\
	& = & q^{2} (q-q^{-1})^{3} (-q)^{n-p} \\
	& & \hspace{5mm} \times
	e_{\alpha_{p} + \cdots + \alpha_{n}} e_{\alpha_{k}+\cdots+2\alpha_{p}+\cdots+2\alpha_{n}} \\
	& & \hspace{15mm} \otimes K_{\alpha_{p} + \cdots + \alpha_{n}}e_{\alpha_{i} + \cdots + \alpha_{p-1}} 
		     K_{\alpha_{k} + \cdots + \alpha_{n}}K_{\alpha_{p} + \cdots + \alpha_{n}}
		     e_{\alpha_{i} + \cdots + \alpha_{k-1}} e_{\alpha_{i} + \cdots + \alpha_{p-1}}  \\
	& = & q^{2} D_{p} E_{k,p}^{p},
\end{eqnarray*}
and for $j=p+1, \ldots, n$,
\begin{eqnarray*}
E_{k,p}^{j} D_{p} 
& = & (q-q^{-1})^{4} (-q)^{n-j} \\
& & \hspace{5mm} \times 
   e_{\alpha_{k}+\cdots+2\alpha_{j}+\cdots+2\alpha_{n}}e_{\alpha_{p}+\cdots+\alpha_{j-1}}
   e_{\alpha_{p} + \cdots + \alpha_{n}} \\
& & \hspace{15mm} \otimes 
    K_{\alpha_{k} + \cdots + \alpha_{n}}K_{\alpha_{p} + \cdots + \alpha_{n}}
    e_{\alpha_{i} + \cdots + \alpha_{k-1}}e_{\alpha_{i} + \cdots + \alpha_{p-1}}
    K_{\alpha_{p} + \cdots + \alpha_{n}}e_{\alpha_{i} + \cdots + \alpha_{p-1}} \\
 & = & q(q-q^{-1})^{4} (-q)^{n-j} \\
& & \hspace{5mm} \times 
   e_{\alpha_{p} + \cdots + \alpha_{n}} 
   e_{\alpha_{k}+\cdots+2\alpha_{j}+\cdots+2\alpha_{n}}e_{\alpha_{p}+\cdots+\alpha_{j-1}} \\
& & \hspace{15mm} \otimes 
    K_{\alpha_{k} + \cdots + \alpha_{n}}K_{\alpha_{p} + \cdots + \alpha_{n}}
    e_{\alpha_{i} + \cdots + \alpha_{k-1}}e_{\alpha_{i} + \cdots + \alpha_{p-1}}
    K_{\alpha_{p} + \cdots + \alpha_{n}}e_{\alpha_{i} + \cdots + \alpha_{p-1}} \\
 & = & q^{2}(q-q^{-1})^{4} (-q)^{n-j} \\
& & \hspace{5mm} \times 
   e_{\alpha_{p} + \cdots + \alpha_{n}} 
   e_{\alpha_{k}+\cdots+2\alpha_{j}+\cdots+2\alpha_{n}}e_{\alpha_{p}+\cdots+\alpha_{j-1}} \\
& & \hspace{15mm} \otimes 
    K_{\alpha_{p} + \cdots + \alpha_{n}} e_{\alpha_{i} + \cdots + \alpha_{p-1}} 
    K_{\alpha_{k} + \cdots + \alpha_{n}} K_{\alpha_{p} + \cdots + \alpha_{n}}
    e_{\alpha_{i} + \cdots + \alpha_{k-1}}e_{\alpha_{i} + \cdots + \alpha_{p-1}} \\
& = & q^{2} D_{p} E_{k,p}^{j}.
\end{eqnarray*}
These calculations show that $E_{k,p}D_{p} = q^{2} D_{p}E_{k,p}$ 
for each $k=i, \ldots, n-1$, and  $p>k$, proving (\ref{appendixB:abstralgrel(3)}).

We now show that $D_{s}E_{k,p} = E_{k,p}D_{s}$ for all $i \leq k < s < p \leq n$; 
and in doing so we will 
use the following identity: for all $i \leq k < s < p \leq n$ we have
$$\left[e_{\alpha_{k}+\cdots+2\alpha_{p}+\cdots+2\alpha_{n}},e_{\alpha_{s} + \cdots + \alpha_{n}}\right]_{q}=0,$$
which we can rewrite as
$e_{\alpha_{k}+\cdots+2\alpha_{p}+\cdots+2\alpha_{n}}e_{\alpha_{s} + \cdots + \alpha_{n}}
   =e_{\alpha_{s} + \cdots + \alpha_{n}}e_{\alpha_{k}+\cdots+2\alpha_{p}+\cdots+2\alpha_{n}}$.
The relation between $D_{s}$ and $E_{i,p}^{p}$ is
\begin{eqnarray*}
D_{s} E_{i,p}^{p} 
& = & (q-q^{-1})^{2} (-q)^{n-p} 
             e_{\alpha_{s} + \cdots + \alpha_{n}}e_{\alpha_{i}+\cdots+2\alpha_{p}+\cdots+2\alpha_{n}} \\
& & \hspace{5mm} \otimes  K_{\alpha_{s} + \cdots + \alpha_{n}}e_{\alpha_{i} + \cdots + \alpha_{s-1}}
                          K_{\alpha_{i} + \cdots + \alpha_{n}}K_{\alpha_{p} + \cdots + \alpha_{n}} 
	                  e_{\alpha_{i} + \cdots + \alpha_{p-1}} \\
& = & (q-q^{-1})^{2} (-q)^{n-p} 
             e_{\alpha_{i}+\cdots+2\alpha_{p}+\cdots+2\alpha_{n}}e_{\alpha_{s} + \cdots + \alpha_{n}} \\
& & \hspace{5mm} \otimes  K_{\alpha_{i} + \cdots + \alpha_{n}}K_{\alpha_{p} + \cdots + \alpha_{n}} 
	                  e_{\alpha_{i} + \cdots + \alpha_{p-1}}
			  K_{\alpha_{s} + \cdots + \alpha_{n}}e_{\alpha_{i} + \cdots + \alpha_{s-1}} \\
& = & E_{i,p}^{p} D_{s}.
\end{eqnarray*}
For each $j=p+1, \ldots, n$,
\begin{eqnarray*}
D_{s} E_{i,p}^{j} 
& = & (q-q^{-1})^{3}(-q)^{n-j} \\
& & \hspace{5mm} \times   
e_{\alpha_{s} + \cdots + \alpha_{n}}e_{\alpha_{i}+\cdots+2\alpha_{j}+\cdots+2\alpha_{n}} e_{\alpha_{p}+\cdots+\alpha_{j-1}} \\
& & \hspace{15mm} \otimes K_{\alpha_{s} + \cdots + \alpha_{n}} e_{\alpha_{i} + \cdots + \alpha_{s-1}}
              K_{\alpha_{i} + \cdots + \alpha_{n}}K_{\alpha_{p} + \cdots + \alpha_{n}} 
	      e_{\alpha_{i} + \cdots + \alpha_{p-1}} \\
	      & = & (q-q^{-1})^{3}(-q)^{n-j} \\
& & \hspace{5mm} \times   
e_{\alpha_{i}+\cdots+2\alpha_{j}+\cdots+2\alpha_{n}} e_{\alpha_{p}+\cdots+\alpha_{j-1}}e_{\alpha_{s} + \cdots + \alpha_{n}} \\
& & \hspace{15mm} \otimes 
              K_{\alpha_{i} + \cdots + \alpha_{n}}K_{\alpha_{p} + \cdots + \alpha_{n}} 
	      e_{\alpha_{i} + \cdots + \alpha_{p-1}} 
	      K_{\alpha_{s} + \cdots + \alpha_{n}} e_{\alpha_{i} + \cdots + \alpha_{s-1}} \\
& & \hspace{15mm} \times \left(e_{\alpha_{s} + \cdots + \alpha_{n}} \otimes K_{\alpha_{s} + \cdots + \alpha_{n}}
       		 e_{\alpha_{i} + \cdots + \alpha_{s-1}}\right) \\
& = & E_{i,p}^{j}D_{s}.
\end{eqnarray*}
For $k=i+1, \ldots, n-1$ and $k < s < p$, 
\begin{eqnarray*}
D_{s} E_{k,p}^{p} 
& = & (q-q^{-1})^{3} (-q)^{n-p} \\
& & \hspace{5mm} \times
               e_{\alpha_{s} + \cdots + \alpha_{n}}e_{\alpha_{k}+\cdots+2\alpha_{p}+\cdots+2\alpha_{n}} \\
& & \hspace{15mm}	    \otimes K_{\alpha_{s} + \cdots + \alpha_{n}}e_{\alpha_{i} + \cdots + \alpha_{s-1}}
	       K_{\alpha_{k} + \cdots + \alpha_{n}}K_{\alpha_{p} + \cdots + \alpha_{n}} 
	             e_{\alpha_{i} + \cdots + \alpha_{k-1}}e_{\alpha_{i} + \cdots + \alpha_{p-1}} \\
& = & (q-q^{-1})^{3} (-q)^{n-p} \\
& & \hspace{5mm} \times
               e_{\alpha_{k}+\cdots+2\alpha_{p}+\cdots+2\alpha_{n}}e_{\alpha_{s} + \cdots + \alpha_{n}} \\
& & \hspace{15mm}	    \otimes 
	       K_{\alpha_{k} + \cdots + \alpha_{n}}K_{\alpha_{p} + \cdots + \alpha_{n}} 
	             e_{\alpha_{i} + \cdots + \alpha_{k-1}}e_{\alpha_{i} + \cdots + \alpha_{p-1}}
		     K_{\alpha_{s} + \cdots + \alpha_{n}}e_{\alpha_{i} + \cdots + \alpha_{s-1}} \\
& = & E_{k,p}^{p}D_{s},
\end{eqnarray*}
and for $j>p$,
\begin{eqnarray*}
D_{s} E_{k,p}^{j} 
& = & (q-q^{-1})^{4} (-q)^{n-j} \\
& & \hspace{5mm} \times
  e_{\alpha_{s} + \cdots + \alpha_{n}}e_{\alpha_{k}+\cdots+2\alpha_{j}+\cdots+2\alpha_{n}}e_{\alpha_{p}+\cdots+\alpha_{j-1}} \\
 & & \hspace{15mm} \otimes 
    K_{\alpha_{s} + \cdots + \alpha_{n}} e_{\alpha_{i} + \cdots + \alpha_{s-1}}
    K_{\alpha_{k} + \cdots + \alpha_{n}}K_{\alpha_{p} + \cdots + \alpha_{n}}
		     e_{\alpha_{i} + \cdots + \alpha_{k-1}}e_{\alpha_{i} + \cdots + \alpha_{p-1}} \\
& = & (q-q^{-1})^{4} (-q)^{n-j} \\
& & \hspace{5mm} \times
e_{\alpha_{k}+\cdots+2\alpha_{j}+\cdots+2\alpha_{n}}e_{\alpha_{p}+\cdots+\alpha_{j-1}}
e_{\alpha_{s} + \cdots + \alpha_{n}} \\
 & & \hspace{15mm} \otimes  K_{\alpha_{k} + \cdots + \alpha_{n}}K_{\alpha_{p} + \cdots + \alpha_{n}}
		     e_{\alpha_{i} + \cdots + \alpha_{k-1}}e_{\alpha_{i} + \cdots + \alpha_{p-1}}
		     K_{\alpha_{s} + \cdots + \alpha_{n}} e_{\alpha_{i} + \cdots + \alpha_{s-1}} \\
& = & E_{k,p}^{j}D_{s}.
\end{eqnarray*}
These calculations show that $D_{s}E_{k,p} = E_{k,p}D_{s}$ 
for all $i \leq k < s < p \leq n$, proving
relation (\ref{appendixB:abstralgrel(4)}).

For all $i \leq p < r < t \leq n$, $F_{p,r,t}$ is defined by
$F_{p,r,t}=D_{p}E_{r,t}-q^{2}E_{r,t}D_{p}$, and we calculate it from its alternative definition
in relation (\ref{appendixB:abstralgrel(7)}):
$F_{p,r,t}=E_{p,r} D_{t} - q^{2} D_{t} E_{p,r}$.
We firstly calculate $F_{i,r,t}$: 
\begin{eqnarray*}
F_{i,r,t} 
& = &  E_{i,r} D_{t} - q^{2} D_{t} E_{i,r}   \\
& = & (q-q^{-1})^{3} 
\Bigg[ \sum_{j=r+1}^{n} (-q)^{n-j} e_{\alpha_{i} + \cdots + 2\alpha_{j} + \cdots + 2\alpha_{n}} 
e_{\alpha_{r} + \cdots + \alpha_{j-1}}e_{\alpha_{t} + \cdots + \alpha_{n}} \\
& & \hspace{30mm} - \sum_{j=r+1}^{t-1} (-q)^{n-j} e_{\alpha_{t} + \cdots + \alpha_{n}}
e_{\alpha_{i} + \cdots + 2\alpha_{j} + \cdots +2\alpha_{n}}e_{\alpha_{r} + \cdots + \alpha_{j-1}} \\
& & \hspace{30mm} - (-q)^{n-t} e_{\alpha_{t} + \cdots + \alpha_{n}}
e_{\alpha_{i} + \cdots + 2\alpha_{t} + \cdots +2\alpha_{n}}e_{\alpha_{r} + \cdots + \alpha_{t-1}} \\
& & \hspace{30mm} - \sum_{j=t+1}^{n} (-q)^{n-j} e_{\alpha_{t} + \cdots + \alpha_{n}}
e_{\alpha_{i} + \cdots + 2\alpha_{j} + \cdots + 2\alpha_{n}}e_{\alpha_{r} + \cdots + \alpha_{j-1}} \Bigg]  \\
& & \hspace{45mm} \otimes K_{\alpha_{i} + \cdots + \alpha_{n}}K_{\alpha_{r} + \cdots + \alpha_{n}}
                          K_{\alpha_{t} + \cdots + \alpha_{n}}
                          e_{\alpha_{i} + \cdots + \alpha_{r-1}}e_{\alpha_{i} + \cdots + \alpha_{t-1}} \\
& = & (q-q^{-1})^{3} \Bigg[
          \sum_{j=t}^{n} (-q)^{n-j} e_{\alpha_{i} + \cdots + 2\alpha_{j} + \cdots + 2\alpha_{n}}
	  e_{\alpha_{r} + \cdots + \alpha_{j-1}}e_{\alpha_{t} + \cdots + \alpha_{n}} \\
& & \hspace{30mm} -(-q)^{n-t}q^{-1} e_{\alpha_{i} + \cdots + 2\alpha_{t} + \cdots + 2\alpha_{n}}
	  e_{\alpha_{t} + \cdots + \alpha_{n}}e_{\alpha_{r} + \cdots + \alpha_{t-1}} \\
& & \hspace{30mm} - \sum_{j=t+1}^{n} (-q)^{n-j}	e_{\alpha_{i} + \cdots + 2\alpha_{j} + \cdots + 2\alpha_{n}} 
	  e_{\alpha_{t} + \cdots + \alpha_{n}}e_{\alpha_{r} + \cdots + \alpha_{j-1}} \Bigg] \\
& & \hspace{45mm} \otimes K_{\alpha_{i} + \cdots + \alpha_{n}}K_{\alpha_{r} + \cdots + \alpha_{n}}
                          K_{\alpha_{t} + \cdots + \alpha_{n}}
                          e_{\alpha_{i} + \cdots + \alpha_{r-1}}e_{\alpha_{i} + \cdots + \alpha_{t-1}} 	\\
& = &  (q-q^{-1})^{3} \Bigg[
     (-q)^{n-t} e_{\alpha_{i} + \cdots + 2\alpha_{t} + \cdots + 2\alpha_{n}}e_{\alpha_{r} + \cdots + \alpha_{n}} \\
& & \hspace{30mm} + \sum_{j=t+1}^{n} (q-q^{-1}) (-q)^{n-j}
      e_{\alpha_{i} + \cdots + 2\alpha_{j} + \cdots + 2\alpha_{n}}e_{\alpha_{r} + \cdots + \alpha_{n}}
      e_{\alpha_{t} + \cdots + \alpha_{j-1}} \Bigg] \\
& & \hspace{45mm} \otimes K_{\alpha_{i} + \cdots + \alpha_{n}}K_{\alpha_{r} + \cdots + \alpha_{n}}
                          K_{\alpha_{t} + \cdots + \alpha_{n}}
                          e_{\alpha_{i} + \cdots + \alpha_{r-1}}e_{\alpha_{i} + \cdots + \alpha_{t-1}}.
\end{eqnarray*}
In obtaining the last equality we used the two identities
$$\begin{array}{rcll}
q^{-1} e_{\alpha_{t} + \cdots + \alpha_{n}}e_{\alpha_{r} + \cdots + \alpha_{t-1}}
 & = & -e_{\alpha_{r} + \cdots + \alpha_{n}}+e_{\alpha_{r} + \cdots + \alpha_{t-1}}e_{\alpha_{t} + \cdots + \alpha_{n}}, & \\
e_{\alpha_{r} + \cdots + \alpha_{j-1}} e_{\alpha_{t} + \cdots + \alpha_{n}}
  - e_{\alpha_{t} + \cdots + \alpha_{n}} e_{\alpha_{r} + \cdots + \alpha_{j-1}}
 & = & (q-q^{-1}) e_{\alpha_{r} + \cdots + \alpha_{n}}e_{\alpha_{t} + \cdots + \alpha_{j-1}}, & t \leq j-1,
 \end{array} $$
which arise from Propositions \ref{appendixC:incrediblehulk(1)} and \ref{appendixC:incrediblehulk(2)} 
 respectively.

Now we will show that $F_{i,r,t} + (1-q^{2}) E_{i,t}D_{r} = 0$.
Note that
\begin{eqnarray*}
(1-q^{2}) E_{i,t}D_{r} 
& = & (1-q^{2}) (q-q^{-1})^{2}  \\
& &  \times
\Bigg[ (-q)^{n-t} e_{\alpha_{i} + \cdots + 2\alpha_{t} + \cdots + 2\alpha_{n}}e_{\alpha_{r} + \cdots + \alpha_{n}} \\
& & \hspace{7mm} + (q-q^{-1})\sum_{j=t+1}^{n} (-q)^{n-j} e_{\alpha_{i} + \cdots + 2\alpha_{j} + \cdots + 2\alpha_{n}}
e_{\alpha_{t} + \cdots + \alpha_{j-1}}e_{\alpha_{r} + \cdots + \alpha_{n}} \Bigg] \\
& & \hspace{25mm} \otimes K_{\alpha_{i} + \cdots + \alpha_{n}}K_{\alpha_{t} + \cdots + \alpha_{n}}
   e_{\alpha_{i} + \cdots + \alpha_{t-1}}K_{\alpha_{r} + \cdots + \alpha_{n}}e_{\alpha_{i} + \cdots + \alpha_{r-1}}.
\end{eqnarray*}
An elementary calculation shows that $F_{i,r,t} + (1-q^{2}) E_{i,t}D_{r} = 0$
after manipulating the expansion of $E_{i,t}D_{r}$ using the following identities:
\begin{eqnarray*}
e_{\alpha_{t} + \cdots + \alpha_{j-1}}e_{\alpha_{r} + \cdots + \alpha_{n}} 
& = & e_{\alpha_{r} + \cdots + \alpha_{n}}e_{\alpha_{t} + \cdots + \alpha_{j-1}}, \\
e_{\alpha_{i} + \cdots + \alpha_{t-1}} K_{\alpha_{r} + \cdots + \alpha_{n}}
& = & K_{\alpha_{r} + \cdots + \alpha_{n}}e_{\alpha_{i} + \cdots + \alpha_{t-1}}, \\
e_{\alpha_{i} + \cdots + \alpha_{t-1}}e_{\alpha_{i} + \cdots + \alpha_{r-1}}
& = & q^{-1}e_{\alpha_{i} + \cdots + \alpha_{r-1}}e_{\alpha_{i} + \cdots + \alpha_{t-1}}.
\end{eqnarray*}
Set $i<p<r<t \leq n$.  Using similar calculations as those above gives
\begin{eqnarray*}
F_{p,r,t} 
 & = & E_{p,r} D_{t} - q^{2} D_{t} E_{p,r} \\
          & = & (q-q^{-1})^{4} \\
	  & & \times \Bigg[
	        (-q)^{n-t} e_{\alpha_{p} + \cdots + 2\alpha_{t}+ \cdots + 2\alpha_{n}}e_{\alpha_{r} + \cdots +\alpha_{n}} \\
	& & \hspace{7mm} + (q-q^{-1})\sum_{j=t+1}^{n} (-q)^{n-j}
	   e_{\alpha_{p} + \cdots + 2\alpha_{j}+ \cdots + 2\alpha_{n}}e_{\alpha_{r} + \cdots +\alpha_{n}}
	   e_{\alpha_{t} + \cdots +\alpha_{j-1}} \Bigg] \\
	& & \hspace{25mm} \otimes K_{\alpha_{t} + \cdots +\alpha_{n}}K_{\alpha_{p} + \cdots +\alpha_{n}}
	   K_{\alpha_{r} + \cdots +\alpha_{n}}e_{\alpha_{i} + \cdots +\alpha_{p-1}}
	   e_{\alpha_{i} + \cdots +\alpha_{r-1}}e_{\alpha_{i} + \cdots +\alpha_{t-1}},
\end{eqnarray*}
and by using the following identities
\begin{eqnarray*}
e_{\alpha_{t} + \cdots +\alpha_{j-1}}e_{\alpha_{r} + \cdots +\alpha_{n}}
& = & e_{\alpha_{r} + \cdots +\alpha_{n}}e_{\alpha_{t} + \cdots +\alpha_{j-1}}, \\
e_{\alpha_{i} + \cdots +\alpha_{t-1}}e_{\alpha_{i} + \cdots +\alpha_{r-1}}
& = & q^{-1}e_{\alpha_{i} + \cdots +\alpha_{r-1}}e_{\alpha_{i} + \cdots +\alpha_{t-1}}, \\
e_{\alpha_{i} + \cdots +\alpha_{p-1}}e_{\alpha_{i} + \cdots +\alpha_{t-1}}K_{\alpha_{r} + \cdots +\alpha_{n}}
& = & K_{\alpha_{r} + \cdots +\alpha_{n}}e_{\alpha_{i} + \cdots +\alpha_{p-1}}e_{\alpha_{i} + \cdots +\alpha_{t-1}},
\end{eqnarray*}
one can easily show that 
$$F_{p,r,t} + (1-q^{2}) E_{p,t}D_{r} = 0,$$
proving relation (\ref{appendixB:abstralgrel(6)}).

We have shown that the components $D_{i}, D_{i+1}, \ldots, D_{n}, D_{\infty}$ of
$\Delta(e_{\alpha_{i} + \cdots +\alpha_{n}})$ satisfy all the claimed relations, 
thus we can calculate
$\big(\Delta(e_{\alpha_{i} + \cdots +\alpha_{n}})\big)^{m}$ for each $m \in\mathbb{N}$.

\end{subsection}

\begin{subsection}{$\Delta(e_{\alpha_{i} + \cdots + 2\alpha_{n}})$, $i=1, \ldots, n-1$}

We now determine the commutation relations between the components of 
$\Delta(e_{\alpha_{i} + \cdots + 2\alpha_{n}})$.
We write $\Delta(e_{\alpha_{i} + \cdots + 2\alpha_{n}})$ as
$$\Delta(e_{\alpha_{i} + \cdots + 2\alpha_{n}}) = 
\sum_{k=i}^{n-1} D_{k} + D_{n} + D_{0} + D_{\infty}, \hspace{10mm} \mbox{where}$$
\begin{eqnarray*}
D_{i} & = & e_{\alpha_{i} + \cdots + 2\alpha_{n}} \otimes K_{\alpha_{i} + \cdots + 2\alpha_{n}}, \\
D_{k} & = & (q-q^{-1}) \left( e_{\alpha_{k} + \cdots + 2\alpha_{n}} \otimes K_{\alpha_{k} + \cdots + 2\alpha_{n}}
                                  e_{\alpha_{i} + \cdots + \alpha_{k-1}}\right), \hspace{10mm} k=i+1, \ldots, n-1, \\
D_{n} & = & (q-q^{-1})(1+q) (e_{n})^{2} \otimes (K_{n})^{2} e_{\alpha_{i} + \cdots + \alpha_{n-1}} \\
D_{0} & = & (q-q^{-1}) \left(e_{n} \otimes K_{n} e_{\alpha_{i} + \cdots + \alpha_{n}}\right), \\
D_{\infty} & = & 1 \otimes e_{\alpha_{i} + \cdots + 2\alpha_{n}}.
\end{eqnarray*}
We claim that the commutation relations between these components are
$$\begin{array}{rcll}
D_{i}D_{k} & = & q^{2} D_{k} D_{i}, & k=0,i+1, \ldots, n, \infty,  \\
D_{j}D_{k} & = & q^{2} D_{k}D_{j}, & i+1 \leq j<k \leq n-1, \\
D_{j}D_{m} & = & q^{2} D_{m}D_{j}, & i+1 \leq j \leq n-1, \ \ m=0, n, \infty,  \\
D_{n}D_{0} & = & q^{2} D_{0} D_{n}, & \\
D_{n} D_{\infty} & = & q^{2} D_{\infty}D_{n} + \xi(D_{0})^{2}, & \xi = -(1+q)^{2}/(q-q^{-1}), \\
D_{0}D_{\infty} & = & q^{2} D_{\infty}D_{0}. &
\end{array}$$
The following easily proved identities will assist in proving these commutation relations:
$$\begin{array}{rcll}
\left[e_{\alpha_{i} + \cdots + 2\alpha_{n}}, e_{\alpha_{k} + \cdots + 2\alpha_{n}}\right]_{q} & = & 0, 
& 1 \leq i < k \leq n-1, \\
\left[e_{\alpha_{i} + \cdots + 2\alpha_{n}}, e_{n}\right]_{q} & = & 0, & 1 \leq i \leq n-1, \\
\left[e_{\alpha_{i} + \cdots + \alpha_{k}},e_{\alpha_{i} + \cdots + 2\alpha_{n}}\right]_{q} & = & 0,
& 1 \leq i < k \leq n-2, \\
\left[e_{\alpha_{i} + \cdots + \alpha_{j}}, e_{\alpha_{i} + \cdots + \alpha_{k}}\right]_{q} & = & 0,
& 1 \leq i<j<k \leq n, \\
\left[e_{\alpha_{i} + \cdots + \alpha_{n-1}},e_{\alpha_{i} + \cdots + 2\alpha_{n}}\right]_{q} & = & 
(1+q)(e_{\alpha_{i} + \cdots + \alpha_{n}})^{2}. &
\end{array} $$
We now prove that the components of $\Delta(e_{\alpha_{i} + \cdots + 2\alpha_{n}})$ 
do satisfy the claimed commutation relations.
For each $k=i+1, \ldots, n-1$ the relation between $D_{i}$ and $D_{k}$ is
\begin{eqnarray*}
D_{i} D_{k} 
 & = & (q-q^{-1}) \left( e_{\alpha_{i} + \cdots + 2\alpha_{n}}e_{\alpha_{k} + \cdots + 2\alpha_{n}}
     \otimes K_{\alpha_{i} + \cdots + 2\alpha_{n}}K_{\alpha_{k} + \cdots + 2\alpha_{n}}
                                  e_{\alpha_{i} + \cdots + \alpha_{k-1}} \right) \\
& = & q^{2} (q-q^{-1}) \left(e_{\alpha_{k} + \cdots + 2\alpha_{n}}e_{\alpha_{i} + \cdots + 2\alpha_{n}}
\otimes K_{\alpha_{k} + \cdots + 2\alpha_{n}}e_{\alpha_{i} + \cdots + \alpha_{k-1}} 
K_{\alpha_{i} + \cdots + 2\alpha_{n}} \right) 
 =  q^{2}  D_{k}D_{i}.
\end{eqnarray*}
Similar (although not identical) calculations show the following:
\begin{itemize}
\item $D_{i} D_{n}= q^{2} D_{n} D_{i}$,
\item $D_{i} D_{0}= q^{2} D_{0} D_{i}$,
\item $D_{i} D_{\infty} = q^{2} D_{\infty} D_{i}$,
\item $D_{k} D_{j} = q^{2} D_{j} D_{k}$ for each $k=i+1, \ldots, n-2$ and each $j = k+1, \ldots, n-1$,
\item $D_{k} D_{n} = q^{2} D_{n} D_{k}$ for each $k=i+1, \ldots, n-1$,
\item $D_{k} D_{0} =  q^{2} D_{0} D_{k}$ for each $k=i+1, \ldots, n-1$,
\item $D_{k} D_{\infty} =  q^{2} D_{\infty}D_{k}$ for each $k=i+1, \ldots, n-1$,
\item $D_{n}D_{0} =  q^{2} D_{0}D_{n}$,
\item $D_{0} D_{\infty} =  q^{2} D_{\infty} D_{0}$,
\end{itemize}
and to complete the proof, the relation between $D_{n}$ and $D_{\infty}$ is
\begin{eqnarray*}
D_{n}D_{\infty} 
& = & (q-q^{-1})(1+q) \left( (e_{n})^{2} \otimes
	(K_{n})^{2} e_{\alpha_{i} + \cdots + \alpha_{n-1}}e_{\alpha_{i} + \cdots + 2\alpha_{n}}\right) \\
 & = & (q-q^{-1})(1+q) \left( (e_{n})^{2} \otimes
   (K_{n})^{2} \Big[ e_{\alpha_{i} + \cdots + 2\alpha_{n}} e_{\alpha_{i} + \cdots + \alpha_{n-1}}
                     + (1+q)(e_{\alpha_{i} + \cdots + \alpha_{n}})^{2} \Big] \right) \\
 & = & q^{2} (q-q^{-1})(1+q) \left( (e_{n})^{2} \otimes
    e_{\alpha_{i} + \cdots + 2\alpha_{n}} (K_{n})^{2} e_{\alpha_{i} + \cdots + \alpha_{n-1}} \right) \\
   & & + (q-q^{-1})(1+q)^{2} \left( (e_{n})^{2} \otimes
          (K_{n})^{2} (e_{\alpha_{i} + \cdots + \alpha_{n}})^{2} \right) 
   =  q^{2} D_{\infty} D_{n} + \xi (D_{0})^{2},
\end{eqnarray*}
where $\xi = -(1+q)^{2}/(q-q^{-1})$.

\end{subsection}

\begin{subsection}{$\Delta(e_{\alpha_{i} + \cdots + 2\alpha_{j} + \cdots + 2 \alpha_{n}})$, 
                   $1 \leq i < j \leq n-1$}

We write the co-multiplication of $e_{\alpha_{i} + \cdots + 2\alpha_{j} + \cdots + 2 \alpha_{n}}$ 
for all $1 \leq i < j \leq n-1$ as
$$\Delta( e_{\alpha_{i} + \cdots + 2\alpha_{j} + \cdots + 2 \alpha_{n}} ) =
D_{0} + \sum_{k=i+1}^{j-1} D_{k} + \overline{D} +  \sum_{k=j}^{n} D_{k} + \sum_{p=j}^{n-1} F_{p} + F_{j-1},$$
where
$$
\begin{array}{rcll}
D_{0} & = & e_{\alpha_{i} + \cdots + 2\alpha_{j} + \cdots + 2 \alpha_{n}} \otimes 
		K_{\alpha_{i} + \cdots + 2\alpha_{j} + \cdots + 2 \alpha_{n}}, & \\
D_{k} & = &  (q-q^{-1})
		\left( e_{\alpha_{k} + \cdots + 2\alpha_{j} + \cdots + 2 \alpha_{n}}
		\otimes K_{\alpha_{k} + \cdots + 2\alpha_{j} + \cdots + 2 \alpha_{n}}e_{\alpha_{i} + \cdots + \alpha_{k-1}}
		\right),
 		 &  k = i+1, \ldots, j-1, \\
\overline{D} & = & (q-q^{-1}) \left( q e_{\alpha_{j} + 2\alpha_{j+1} + \cdots + 2 \alpha_{n}} e_{j} 
   		-q^{-1} e_{j}e_{\alpha_{j} + 2\alpha_{j+1} + \cdots + 2 \alpha_{n}}\right) & \\
   		& & \hspace{60mm} \otimes (K_{\alpha_{j} + \cdots + \alpha_{n}})^{2} e_{\alpha_{i} + \cdots + \alpha_{j-1}}, & \\
D_{k} & = &  (q-q^{-1}) \left( \overline{e}_{\alpha_{j} + \cdots + 2\alpha_{k+1} + \cdots + 2\alpha_{n}}
		\otimes K_{\alpha_{j} + \cdots + 2\alpha_{k+1} + \cdots + 2\alpha_{n}}e_{\alpha_{i} + \cdots + \alpha_{k}}
		\right),
		&  k = j, \ldots, n-1, \\
D_{n} & = &  (q-q^{-1})
             \left( \overline{e}_{\alpha_{j} + \cdots + \alpha_{n}} 
		\otimes K_{\alpha_{j} + \cdots + \alpha_{n}} e_{\alpha_{i} + \cdots + \alpha_{n}}\right), & \\
F_{p} & = &  (q-q^{-1}) 
             \left( \overline{e}_{\alpha_{j} + \cdots + \alpha_{p}} \otimes K_{\alpha_{j} + \cdots + \alpha_{p}}
		e_{\alpha_{i} + \cdots + 2\alpha_{p+1} + \cdots + 2\alpha_{n}} \right), &  p = j, \ldots, n-1, \\
F_{j-1} & = & 1 \otimes e_{\alpha_{i} + \cdots + 2\alpha_{j} + \cdots + 2\alpha_{n}}, &
\end{array} $$
noting the slightly different normalisations compared to those used in the previous section.
We claim that the commutation relations between the components of 
$\Delta( e_{\alpha_{i} + \cdots + 2\alpha_{j} + \cdots + 2 \alpha_{n}})$ are
\begin{eqnarray}
\lefteqn{
(\overline{D} + D_{j} + D_{j+1} + \cdots + D_{n-1})(F_{n-1}+F_{n-2} + \cdots + F_{j-1}) } \nonumber \\
& & 
=q^{2} (F_{n-1}+F_{n-2} + \cdots + F_{j-1})(\overline{D} + D_{j} + D_{j+1} + \cdots + D_{n-1})
 + \xi (D_{n})^{2}, \label{appendixB:marklathbabes(1)}
\end{eqnarray}
where $\xi = -(1+q)^{2}/(q-q^{-1})$, and
$$\begin{array}{rcll}
D_{0} D_{k} & = & q^{2} D_{k} D_{0}, & k=i+1, \ldots, j-1, \\
D_{0} \overline{D} & = & q^{2} \overline{D}D_{0}, & \\
D_{0} D_{k} & = & q^{2} D_{k} D_{0}, & k=j, \ldots, n, \\
D_{0} F_{p} & = & q^{2} F_{p} D_{0}, & p=j-1, \ldots, n-1, \\
D_{k} D_{p} & = & q^{2} D_{p} D_{k}, & i+1 \leq k < p \leq j-1, \\
D_{k} \overline{D} & = & q^{2} \overline{D} D_{k}, & k=i+1, \ldots, j-1, \\
D_{k} D_{p} & = & q^{2} D_{p} D_{k}, & k= i+1, \ldots, j-1, \ p= j, \ldots, n, \\
D_{k} F_{p} & = & q^{2} F_{p} D_{k}, & k=i+1, \ldots, j-1, \ p=j-1, \ldots, n-1, \\
\overline{D} D_{k} & = & q^{2} D_{k} \overline{D}, & k=j, \ldots, n, \\
\overline{D} F_{p} & = & q^{2} F_{p} \overline{D}, & p=j, \ldots, n-1, \\
D_{k} D_{p} & = & q^{2} D_{p} D_{k}, & j \leq k < p \leq n, \\
D_{k} F_{p} & = & q^{2} F_{p} D_{k}, & k=j, \ldots, n, \ p=j-1, \ldots, n-1, \ k \neq p, \\
F_{k} F_{p} & = & q^{2} F_{p} F_{k}, & n-1 \geq k > l \geq j-1.
\end{array}$$
If the components of $\Delta( e_{\alpha_{i} + \cdots + 2\alpha_{j} + \cdots + 2 \alpha_{n}})$
do satisfy these commutation relations then the following relations
\begin{eqnarray*}
\lefteqn{
(\overline{D} + D_{j} + D_{j+1} + \cdots + D_{n-1})(F_{n-1}+F_{n-2} + \cdots + F_{j-1}) } \\
& & =q^{2} (F_{n-1}+F_{n-2} + \cdots + F_{j-1})(\overline{D} + D_{j} + D_{j+1} + \cdots + D_{n-1})
 + \xi (D_{n})^{2}, \\
\lefteqn{
(\overline{D} + D_{j} + D_{j+1} + \cdots + D_{n-1})D_{n}
=q^{2} D_{n}(\overline{D} + D_{j} + D_{j+1} + \cdots + D_{n-1}), } \\
\lefteqn{
D_{n}(F_{n-1}+F_{n-2} + \cdots + F_{j-1}) = q^{2}(F_{n-1}+F_{n-2} + \cdots + F_{j-1})D_{n}, } 
\end{eqnarray*}
show that
$(\overline{D} + D_{j} + D_{j+1} + \cdots + D_{n-1})$, $D_{n}$ and $(F_{n-1}+F_{n-2} + \cdots + F_{j-1})$
satisfy the same relations as do the elements $a, b$ and $c$, respectively,
 in Lemma \ref{appendixB:generalisationnewtonsbinomial2}.
 This means that we can immediately use the $q$-binomial theorem and 
 Lemma \ref{appendixB:generalisationnewtonsbinomial2} 
 (one of the generalisations of the Binomial theorem) 
 to obtain an expression for
$\big(\Delta( e_{\alpha_{i} + \cdots + 2\alpha_{j} + \cdots + 2 \alpha_{n}} ) \big)^{m}$ 
for each $m \in \mathbb{N}$.

We now prove that the components of 
$\Delta( e_{\alpha_{i} + \cdots + 2\alpha_{j} + \cdots + 2 \alpha_{n}})$ do
satisfy the claimed commutation relations.  
We firstly consider the easier relations, and then we consider 
the more complicated calculations needed to prove (\ref{appendixB:marklathbabes(1)}).  
In showing this last relation we must consider the most complicated commutation relations 
in this problem, namely those between 
 $\overline{D}$ and $F_{j-1}$, and between $D_{k}$ and $F_{k}$ for each $k=j, \ldots, n-1$.

We firstly consider some relations that we will extensively use.
For each $p=j, \ldots, n$,
\begin{equation}
\label{appendixC:eqpig}
\left[e_{\alpha_{k} + \cdots + 2\alpha_{j} + \cdots + 2 \alpha_{n}},e_{p} \right]_{q}=0,
\end{equation}
which we can rewrite as
$$\begin{array}{rcll}
e_{\alpha_{k} + \cdots + 2\alpha_{j} + \cdots + 2 \alpha_{n}}e_{j} & = & 
 q e_{j} e_{\alpha_{k} + \cdots + 2\alpha_{j} + \cdots + 2 \alpha_{n}}, & \\
e_{\alpha_{k} + \cdots + 2\alpha_{j} + \cdots + 2 \alpha_{n}}e_{m} & = & 
  e_{m} e_{\alpha_{k} + \cdots + 2\alpha_{j} + \cdots + 2 \alpha_{n}}, & \mbox{if } m=j+1, \ldots, n.
  \end{array} $$
Equation (\ref{appendixC:eqpig}) implies the following result:
for each $\gamma \in \phi$ satisfying $ \alpha_{j} \preceq \gamma \prec \alpha_{j+1}$ we have
$$\left[ e_{\alpha_{k} + \cdots + 2\alpha_{j} + \cdots + 2\alpha_{n}}, \overline{e}_{\gamma}\right]_{q}=0,$$
which we can rewrite as
$ e_{\alpha_{k} + \cdots + 2\alpha_{j} + \cdots + 2\alpha_{n}} \overline{e}_{\gamma}
 = q \overline{e}_{\gamma} e_{\alpha_{k} + \cdots + 2\alpha_{j} + \cdots + 2\alpha_{n}}$.
 For all $i<k$ we have
$$\left[ e_{\alpha_{i} + \cdots + 2\alpha_{j} + \cdots + 2\alpha_{n}},
          e_{\alpha_{k} + \cdots + 2\alpha_{j} + \cdots + 2\alpha_{n}}\right]_{q}=0,$$
which we can rewrite as
$$e_{\alpha_{i} + \cdots + 2\alpha_{j} + \cdots + 2\alpha_{n}}
          e_{\alpha_{k} + \cdots + 2\alpha_{j} + \cdots + 2\alpha_{n}}
	  =q e_{\alpha_{k} + \cdots + 2\alpha_{j} + \cdots + 2\alpha_{n}}
	  e_{\alpha_{i} + \cdots + 2\alpha_{j} + \cdots + 2\alpha_{n}}.$$
A further useful identity is
$\left[ e_{\alpha_{i} + \cdots + 2\alpha_{j} + \cdots + 2\alpha_{n}},
         e_{\alpha_{j} + 2\alpha_{j+1} + \cdots + 2\alpha_{n}} \right]_{q}=0$,
which we can rewrite as
$$e_{\alpha_{i} + \cdots + 2\alpha_{j} + \cdots + 2\alpha_{n}}e_{\alpha_{j} + 2\alpha_{j+1} + \cdots + 2\alpha_{n}}=
 q e_{\alpha_{j} + 2\alpha_{j+1} + \cdots + 2\alpha_{n}}e_{\alpha_{i} + \cdots + 2\alpha_{j} + \cdots + 2\alpha_{n}}.$$

We now prove that the components of $\Delta(e_{\alpha_{i} + \cdots + 2\alpha_{j} + \cdots + 2 \alpha_{n}})$ do satisfy the
claimed commutation relations.
We will firstly prove that $D_{0} D_{k} = q^{2} D_{k} D_{0}$ for each $k=i+1, \ldots, j-1$.  
For each such $k$,
\begin{eqnarray*}
D_{0} D_{k} 
& = & (q-q^{-1}) e_{\alpha_{i} + \cdots + 2\alpha_{j} + \cdots + 2 \alpha_{n}}
                e_{\alpha_{k} + \cdots + 2\alpha_{j} + \cdots + 2 \alpha_{n}} \\
& & \hspace{25mm} \otimes K_{\alpha_{i} + \cdots + 2\alpha_{j} + \cdots + 2 \alpha_{n}}
                K_{\alpha_{k} + \cdots + 2\alpha_{j} + \cdots + 2 \alpha_{n}}e_{\alpha_{i} + \cdots + \alpha_{k-1}} \\
& = & q^{2} (q-q^{-1})	e_{\alpha_{k} + \cdots + 2\alpha_{j} + \cdots + 2 \alpha_{n}}
                        e_{\alpha_{i} + \cdots + 2\alpha_{j} + \cdots + 2 \alpha_{n}} \\
& & \hspace{25mm} \otimes 
                K_{\alpha_{k} + \cdots + 2\alpha_{j} + \cdots + 2 \alpha_{n}}e_{\alpha_{i} + \cdots + \alpha_{k-1}}
		K_{\alpha_{i} + \cdots + 2\alpha_{j} + \cdots + 2 \alpha_{n}} \\	
& = & q^{2} D_{k} D_{0}.
\end{eqnarray*}
We now prove that $D_{0} \overline{D} = q^{2} \overline{D}D_{0}$:
\begin{eqnarray*}
D_{0} \overline{D} & = & (q-q^{-1})\left(e_{\alpha_{i} + \cdots + 2\alpha_{j} + \cdots + 2 \alpha_{n}} \otimes 
		K_{\alpha_{i} + \cdots + 2\alpha_{j} + \cdots + 2 \alpha_{n}}\right) \\
& & \hspace{5mm} \times \left( q e_{\alpha_{j} + 2\alpha_{j+1} + \cdots + 2 \alpha_{n}} e_{j} 
   		-q^{-1} e_{j}e_{\alpha_{j} + 2\alpha_{j+1} + \cdots + 2 \alpha_{n}}\right) 
                 \otimes (K_{\alpha_{j}+\cdots+\alpha_{n}})^{2} e_{\alpha_{i} + \cdots + \alpha_{j-1}} \\
& = & 	q^{2}(q-q^{-1}) \left( q e_{\alpha_{j} + 2\alpha_{j+1} + \cdots + 2 \alpha_{n}} e_{j} 
   		-q^{-1} e_{j}e_{\alpha_{j} + 2\alpha_{j+1} + \cdots + 2 \alpha_{n}}\right) 
		e_{\alpha_{i} + \cdots + 2\alpha_{j} + \cdots + 2 \alpha_{n}} \\
& & \hspace{5mm} \otimes (K_{\alpha_{j}+\cdots+\alpha_{n}})^{2} e_{\alpha_{i} + \cdots + \alpha_{j-1}}
		K_{\alpha_{i} + \cdots + 2\alpha_{j} + \cdots + 2 \alpha_{n}}  
  =  q^{2}\overline{D}	D_{0}.
\end{eqnarray*}
In a similar way, we can prove the following:
\begin{itemize}
\item $D_{0}D_{k} = q^{2} D_{k} D_{0}$ for each $k=j, \ldots, n-1$,
\item $D_{0} D_{n} = q^{2} D_{n} D_{0}$,
\item $D_{0} F_{k} = q^{2} F_{k} D_{0}$ for each $k=n-1, \ldots, j$,
\item $D_{0} F_{j-1} = q^{2}    F_{j-1} D_{0}$.
\end{itemize}

We now prove that $D_{k} D_{p} = q^{2} D_{p} D_{k}$ for each $k=i+1, \ldots, j-2$ 
and each $p=k+1, \ldots, j-1$:
\begin{eqnarray*}
D_{k} D_{p} 
& = & 	(q-q^{-1})^{2}	e_{\alpha_{k} + \cdots + 2\alpha_{j} + \cdots + 2 \alpha_{n}}
                      e_{\alpha_{p} + \cdots + 2\alpha_{j} + \cdots + 2 \alpha_{n}} \\
& &  \hspace{25mm}    \otimes K_{\alpha_{k} + \cdots + 2\alpha_{j} + \cdots + 2 \alpha_{n}}
e_{\alpha_{i} + \cdots + \alpha_{k-1}}
     K_{\alpha_{p} + \cdots + 2\alpha_{j} + \cdots + 2 \alpha_{n}}e_{\alpha_{i} + \cdots + \alpha_{p-1}} \\
& = & q	(q-q^{-1})^{2} e_{\alpha_{p} + \cdots + 2\alpha_{j} + \cdots + 2 \alpha_{n}}
		      e_{\alpha_{k} + \cdots + 2\alpha_{j} + \cdots + 2 \alpha_{n}} \\ 
& &  \hspace{25mm}    \otimes     K_{\alpha_{k} + \cdots + 2\alpha_{j} + \cdots + 2 \alpha_{n}}
     K_{\alpha_{p} + \cdots + 2\alpha_{j} + \cdots + 2 \alpha_{n}}e_{\alpha_{i} + \cdots + \alpha_{k-1}}
     e_{\alpha_{i} + \cdots + \alpha_{p-1}} \\
 & = & q^{2}     (q-q^{-1})^{2} e_{\alpha_{p} + \cdots + 2\alpha_{j} + \cdots + 2 \alpha_{n}}
		      e_{\alpha_{k} + \cdots + 2\alpha_{j} + \cdots + 2 \alpha_{n}} \\
 & &  \hspace{25mm}    \otimes K_{\alpha_{p} + \cdots + 2\alpha_{j} + \cdots + 2 \alpha_{n}}
     e_{\alpha_{i} + \cdots + \alpha_{p-1}} 
      K_{\alpha_{k} + \cdots + 2\alpha_{j} + \cdots + 2 \alpha_{n}}
     e_{\alpha_{i} + \cdots + \alpha_{k-1}} \\
  & = & q^{2} D_{p} D_{k}.   
\end{eqnarray*}
We now prove that $D_{k} \overline{D} = q^{2} \overline{D} D_{k}$ for each $k=i+1, \ldots, j-1$:
\begin{eqnarray*}
D_{k} \overline{D} 
& = & q^{2} (q-q^{-1})^{2}	
        \left( q e_{\alpha_{j} + 2\alpha_{j+1} + \cdots + 2 \alpha_{n}} e_{j} 
   		-q^{-1} e_{j}e_{\alpha_{j} + 2\alpha_{j+1} + \cdots + 2 \alpha_{n}}\right)  
		  e_{\alpha_{k} + \cdots + 2\alpha_{j} + \cdots + 2 \alpha_{n}} \\
& & \hspace{15mm} \otimes 	
           K_{\alpha_{k} + \cdots + 2\alpha_{j} + \cdots + 2 \alpha_{n}}
	   (K_{\alpha_{j} + \cdots + \alpha_{n}})^{2}
	   e_{\alpha_{i} + \cdots + \alpha_{k-1}}	 
	    e_{\alpha_{i} + \cdots + \alpha_{j-1}} \\
& = & q^{3} (q-q^{-1})^{2}	
        \left( q e_{\alpha_{j} + 2\alpha_{j+1} + \cdots + 2 \alpha_{n}} e_{j} 
   		-q^{-1} e_{j}e_{\alpha_{j} + 2\alpha_{j+1} + \cdots + 2 \alpha_{n}}\right)  
		  e_{\alpha_{k} + \cdots + 2\alpha_{j} + \cdots + 2 \alpha_{n}} \\
& & \hspace{15mm} \otimes 	
           K_{\alpha_{k} + \cdots + 2\alpha_{j} + \cdots + 2 \alpha_{n}}
	   (K_{\alpha_{j} + \cdots + \alpha_{n}})^{2}
	    e_{\alpha_{i} + \cdots + \alpha_{j-1}} e_{\alpha_{i} + \cdots + \alpha_{k-1}} \\
& = & q^{2} (q-q^{-1})^{2}	
        \left( q e_{\alpha_{j} + 2\alpha_{j+1} + \cdots + 2 \alpha_{n}} e_{j} 
   		-q^{-1} e_{j}e_{\alpha_{j} + 2\alpha_{j+1} + \cdots + 2 \alpha_{n}}\right)  
		  e_{\alpha_{k} + \cdots + 2\alpha_{j} + \cdots + 2 \alpha_{n}} \\
& & \hspace{15mm} \otimes 	
	   (K_{\alpha_{j} + \cdots + \alpha_{n}})^{2}
	    e_{\alpha_{i} + \cdots + \alpha_{j-1}} 
	    K_{\alpha_{k} + \cdots + 2\alpha_{j} + \cdots + 2 \alpha_{n}}
	    e_{\alpha_{i} + \cdots + \alpha_{k-1}} \\
& = & q^{2} \overline{D} D_{k}.   	    	    
\end{eqnarray*}
In a similar way, we can show the following:
\begin{itemize}
\item $D_{k} D_{p} = q^{2} D_{p} D_{k}$ for each $k=i+1, \ldots, j-1$ 
and each $p=j, \ldots, n-1$,
\item $D_{k} D_{n} = q^{2} D_{n} D_{k}$ for each $k=i+1, \ldots, j-1$,
\item $D_{k} F_{p} = q^{2} F_{p} D_{k}$ for each $k=i+1, \ldots, j-1$ 
and each $p=j, \ldots, n-1$,
\item $D_{k} F_{j-1} = q^{2} F_{j-1} D_{k}$ for each $k=i+1, \ldots, j-1$.
\end{itemize}
We now prove that  
$ \overline{D} F_{k} = q^{2} F_{k} \overline{D}$ for all $k=j, \ldots, n-1$:
\begin{eqnarray*}
\lefteqn{ \overline{D} F_{k} } \\
& = & (q-q^{-1})^{2} 
    \overline{e}_{\alpha_{j} + \cdots + \alpha_{k}}\left( q e_{\alpha_{j} + 2\alpha_{j+1} + \cdots + 2 \alpha_{n}} e_{j} 
   		-q^{-1} e_{j}e_{\alpha_{j} + 2\alpha_{j+1} + \cdots + 2 \alpha_{n}}\right) \\
& & \hspace{25mm} \otimes (K_{\alpha_{j} + \cdots +\alpha_{n}})^{2} e_{\alpha_{i} + \cdots + \alpha_{j-1}}
                K_{\alpha_{j} + \cdots + \alpha_{k}}e_{\alpha_{i} + \cdots + 2\alpha_{k+1} + \cdots + 2\alpha_{n}} \\
& = & q (q-q^{-1})^{2} 
    \overline{e}_{\alpha_{j} + \cdots + \alpha_{k}}\left( q e_{\alpha_{j} + 2\alpha_{j+1} + \cdots + 2 \alpha_{n}} e_{j} 
   		-q^{-1} e_{j}e_{\alpha_{j} + 2\alpha_{j+1} + \cdots + 2 \alpha_{n}}\right) \\
& & \hspace{25mm} \otimes K_{\alpha_{j} + \cdots + \alpha_{k}}
 (K_{\alpha_{j} + \cdots +\alpha_{n}})^{2} e_{\alpha_{i} + \cdots + \alpha_{j-1}}
 e_{\alpha_{i} + \cdots + 2\alpha_{k+1} + \cdots + 2\alpha_{n}} \\
& = & q^{2} (q-q^{-1})^{2} 
    \overline{e}_{\alpha_{j} + \cdots + \alpha_{k}}\left( q e_{\alpha_{j} + 2\alpha_{j+1} + \cdots + 2 \alpha_{n}} e_{j} 
   		-q^{-1} e_{j}e_{\alpha_{j} + 2\alpha_{j+1} + \cdots + 2 \alpha_{n}}\right) \\ 
& & \hspace{25mm} \otimes K_{\alpha_{j} + \cdots + \alpha_{k}} (K_{\alpha_{j} + \cdots +\alpha_{n}})^{2}
      e_{\alpha_{i} + \cdots + 2\alpha_{k+1} + \cdots + 2\alpha_{n}}e_{\alpha_{i} + \cdots + \alpha_{j-1}} \\
& = & q^{2} (q-q^{-1})^{2} 
    \overline{e}_{\alpha_{j} + \cdots + \alpha_{k}}\left( q e_{\alpha_{j} + 2\alpha_{j+1} + \cdots + 2 \alpha_{n}} e_{j} 
   		-q^{-1} e_{j}e_{\alpha_{j} + 2\alpha_{j+1} + \cdots + 2 \alpha_{n}}\right) \\
& & \hspace{25mm} \otimes K_{\alpha_{j} + \cdots + \alpha_{k}}e_{\alpha_{i} + \cdots + 2\alpha_{k+1} + \cdots + 2\alpha_{n}}
(K_{\alpha_{j} + \cdots +\alpha_{n}})^{2}e_{\alpha_{i} + \cdots + \alpha_{j-1}} 
  =  q^{2} F_{k}\overline{D}. 
\end{eqnarray*}
In this calculation we used the fact that the element
$$\left( q e_{\alpha_{j} + 2\alpha_{j+1} + \cdots + 2 \alpha_{n}} e_{j} 
   		-q^{-1} e_{j}e_{\alpha_{j} + 2\alpha_{j+1} + \cdots + 2 \alpha_{n}}\right)$$ 
commutes with $e_{k}$ for each $k=j, \ldots, n$ (see Proposition \ref{appendixB:iceage(1)}).
We can also show the following:
\begin{itemize}
\item $D_{k}D_{p} = q^{2} D_{p} D_{k}$ for all $j \leq k < p \leq n-1$,
\item $D_{k} D_{n} = q^{2} D_{n} D_{k}$ for each $k=j, \ldots, n-1$.
\end{itemize}

We now prove that $D_{k} F_{p} = q^{2}  F_{p} D_{k}$ 
for all $k=j, \ldots, n-1$ and all $p = j, \ldots, n-1$ satisfying $k \neq p$:
\begin{eqnarray*}
D_{k} F_{p} 
& = & (q-q^{-1})^{2} \big( \overline{e}_{\alpha_{j} + \cdots + 2\alpha_{k+1} + \cdots + 2\alpha_{n}}
      \overline{e}_{\alpha_{j} + \cdots + \alpha_{p}} \\
 &  & \hspace{25mm} \otimes 
     K_{\alpha_{j} + \cdots + 2\alpha_{k+1} + \cdots + 2\alpha_{n}} e_{\alpha_{i} + \cdots + \alpha_{k}}
       K_{\alpha_{j} + \cdots + \alpha_{p}}
		e_{\alpha_{i} + \cdots + 2\alpha_{p+1} + \cdots + 2\alpha_{n}} \big) \\
& = & q (q-q^{-1})^{2} \overline{e}_{\alpha_{j} + \cdots + \alpha_{p}}
      \overline{e}_{\alpha_{j} + \cdots + 2\alpha_{k+1} + \cdots + 2\alpha_{n}} \\
& & \hspace{25mm} \otimes
     K_{\alpha_{j} + \cdots + \alpha_{p}} K_{\alpha_{j} + \cdots + 2\alpha_{k+1} + \cdots + 2\alpha_{n}} 
       e_{\alpha_{i} + \cdots + \alpha_{k}}e_{\alpha_{i} + \cdots + 2\alpha_{p+1} + \cdots + 2\alpha_{n}} \\
& = & q^{2} (q-q^{-1})^{2} \overline{e}_{\alpha_{j} + \cdots + \alpha_{p}}
      \overline{e}_{\alpha_{j} + \cdots + 2\alpha_{k+1} + \cdots + 2\alpha_{n}} \\
& & \hspace{25mm} \otimes
 K_{\alpha_{j} + \cdots + \alpha_{p}}e_{\alpha_{i} + \cdots + 2\alpha_{p+1} + \cdots + 2\alpha_{n}}
 K_{\alpha_{j} + \cdots + 2\alpha_{k+1} + \cdots + 2\alpha_{n}}e_{\alpha_{i} + \cdots + \alpha_{k}} \\
& = &  q^{2}  F_{p} D_{k},
\end{eqnarray*}
where we used the relations
$$\left[ e_{\alpha_{j} + \cdots + \alpha_{p}}, e_{\alpha_{j} + \cdots + 2\alpha_{k+1} + \cdots + 2\alpha_{n}}\right]_{q}=0,
\hspace{5mm} \left[ \overline{e}_{\alpha_{j} + \cdots + 2\alpha_{k+1} + \cdots + 2\alpha_{n}},
  \overline{e}_{\alpha_{j} + \cdots + \alpha_{p}}\right]_{q}=0,
\hspace{10mm} k \neq p.$$

We can also show the following are true:
\begin{itemize}
\item $\overline{D} D_{k} = q^{2} D_{k} \overline{D}$ for each $k=j, \ldots, n-1$,
\item $\overline{D} D_{n} =  q^{2} D_{n} \overline{D}$,
\item $D_{k} F_{j-1} = q^{2} F_{j-1} D_{k}$ for each $k=j, \ldots, n-1$.
\end{itemize}
We now prove that $D_{n} F_{k} = q^{2} F_{k}D_{n}$ for each $k=n-1, \ldots, j$:
\begin{eqnarray*}
D_{n} F_{k} 
& = & (q-q^{-1})^{2}
     \overline{e}_{\alpha_{j} + \cdots + \alpha_{n}}\overline{e}_{\alpha_{j} + \cdots + \alpha_{k}}
     \otimes K_{\alpha_{j} + \cdots + \alpha_{n}} e_{\alpha_{i} + \cdots + \alpha_{n}}
     K_{\alpha_{j} + \cdots + \alpha_{k}}
		e_{\alpha_{i} + \cdots + 2\alpha_{k+1} + \cdots + 2\alpha_{n}} \\
& = & q (q-q^{-1})^{2}
      \overline{e}_{\alpha_{j} + \cdots + \alpha_{k}}\overline{e}_{\alpha_{j} + \cdots + \alpha_{n}}
      \otimes K_{\alpha_{j} + \cdots + \alpha_{k}}K_{\alpha_{j} + \cdots + \alpha_{n}}
      e_{\alpha_{i} + \cdots + \alpha_{n}} e_{\alpha_{i} + \cdots + 2\alpha_{k+1} + \cdots + 2\alpha_{n}} \\
& = & q^{2}   (q-q^{-1})^{2} 
      \overline{e}_{\alpha_{j} + \cdots + \alpha_{k}}\overline{e}_{\alpha_{j} + \cdots + \alpha_{n}}
       \otimes  K_{\alpha_{j} + \cdots + \alpha_{k}} e_{\alpha_{i} + \cdots + 2\alpha_{k+1} + \cdots + 2\alpha_{n}} 
         K_{\alpha_{j} + \cdots + \alpha_{n}}e_{\alpha_{i} + \cdots + \alpha_{n}} \\
& = & q^{2} F_{k}D_{n},
\end{eqnarray*}
as $\left[e_{\alpha_{i}+\cdots+\alpha_{n}},e_{\alpha_{i}+\cdots+2\alpha_{j}+\cdots+2\alpha_{n}}\right]_{q}=0$, 
and the relation between $D_{n}$ and $F_{j-1}$ is
$$
D_{n} F_{j-1} 
 =  q^{2} (q-q^{-1}) \overline{e}_{\alpha_{j} + \cdots + \alpha_{n}} 
      \otimes e_{\alpha_{i} + \cdots + 2\alpha_{j} + \cdots + 2\alpha_{n}}
              K_{\alpha_{j} + \cdots + \alpha_{n}} e_{\alpha_{i} + \cdots + \alpha_{n}} 
  =  q^{2}  F_{j-1} D_{n}.
$$
We now prove that $F_{k} F_{p} = q^{2} F_{p} F_{k}$ for each $k=n-1, \ldots, j+1$ 
and each $p=k-1, \ldots, j$:
\begin{eqnarray*}
F_{k} F_{p} & = & (q-q^{-1})^{2}
	  \overline{e}_{\alpha_{j} + \cdots + \alpha_{k}}\overline{e}_{\alpha_{j} + \cdots + \alpha_{p}} \\
& & \hspace{25mm}	 \otimes K_{\alpha_{j} + \cdots + \alpha_{k}}
		e_{\alpha_{i} + \cdots + 2\alpha_{k+1} + \cdots + 2\alpha_{n}}
		K_{\alpha_{j} + \cdots + \alpha_{p}}
		e_{\alpha_{i} + \cdots + 2\alpha_{p+1} + \cdots + 2\alpha_{n}} \\
 & = & q^{2} (q-q^{-1})^{2}
     \overline{e}_{\alpha_{j} + \cdots + \alpha_{p}} \overline{e}_{\alpha_{j} + \cdots + \alpha_{k}} \\
 & & \hspace{25mm}  \otimes  K_{\alpha_{j} + \cdots + \alpha_{p}}  K_{\alpha_{j} + \cdots + \alpha_{k}}
              e_{\alpha_{i} + \cdots + 2\alpha_{p+1} + \cdots + 2\alpha_{n}}
	      e_{\alpha_{i} + \cdots + 2\alpha_{k+1} + \cdots + 2\alpha_{n}} \\
& = & q^{2} F_{p} F_{k},
\end{eqnarray*}
where we have used the identity $\left[e_{\alpha_{i} + \cdots + 2\alpha_{k+1} + \cdots + 2\alpha_{n}},
               e_{\alpha_{i} + \cdots + 2\alpha_{p+1} + \cdots + 2\alpha_{n}}\right]_{q}=0$.
In addition, for each $k=n-1, \ldots, j$, 
$$
F_{k} F_{j-1} 
 =  q^{2} (q-q^{-1}) \overline{e}_{\alpha_{j} + \cdots + \alpha_{k}}
  \otimes e_{\alpha_{i} + \cdots + 2\alpha_{j} + \cdots + 2\alpha_{n}}
  K_{\alpha_{j} + \cdots + \alpha_{k}}
		 e_{\alpha_{i} + \cdots + 2\alpha_{k+1} + \cdots + 2\alpha_{n}} 
 =  q^{2} F_{j-1} F_{k}.
$$
This completes the proof of the easier commutation relations.  
We now prove the more complicated relations.
The relation between $\overline{D}$ and $F_{j-1}$ is:
\begin{eqnarray*}
\lefteqn{
\overline{D} F_{j-1} } \\
 & = & (q-q^{-1}) \left( 
q e_{\alpha_{j} + 2\alpha_{j+1} + \cdots + 2 \alpha_{n}} e_{j} 
-q^{-1} e_{j}e_{\alpha_{j} + 2\alpha_{j+1} + \cdots + 2 \alpha_{n}} \right) \\
& &  \hspace{10mm} \otimes (K_{\alpha_{j} + \cdots + \alpha_{n}})^{2} e_{\alpha_{i} + \cdots + \alpha_{j-1}}
e_{\alpha_{i} + \cdots + 2\alpha_{j} + \cdots + 2\alpha_{n}} \\
& = & (q-q^{-1}) \left(q e_{\alpha_{j} + 2\alpha_{j+1} + \cdots + 2 \alpha_{n}} e_{j} 
-q^{-1} e_{j}e_{\alpha_{j} + 2\alpha_{j+1} + \cdots + 2 \alpha_{n}} \right) \\
& & \hspace{10mm}  \otimes (K_{\alpha_{j} + \cdots + \alpha_{n}})^{2} 
\left(  e_{\alpha_{i} + \cdots + 2\alpha_{j} + \cdots + 2\alpha_{n}} e_{\alpha_{i} + \cdots + \alpha_{j-1}}
  + \left[e_{\alpha_{i} + \cdots + \alpha_{j-1}},e_{\alpha_{i} + \cdots + 2\alpha_{j} + \cdots + 2\alpha_{n}} \right]_{q}
      \right) \\
& = & q^{2} (q-q^{-1}) \left(q e_{\alpha_{j} + 2\alpha_{j+1} + \cdots + 2 \alpha_{n}} e_{j} 
-q^{-1} e_{j}e_{\alpha_{j} + 2\alpha_{j+1} + \cdots + 2 \alpha_{n}} \right) \\
& & \hspace{10mm}  \otimes e_{\alpha_{i} + \cdots + 2\alpha_{j} + \cdots + 2\alpha_{n}}
(K_{\alpha_{j} + \cdots + \alpha_{n}})^{2}  e_{\alpha_{i} + \cdots + \alpha_{j-1}} \\
& & + (q-q^{-1}) \left(q e_{\alpha_{j} + 2\alpha_{j+1} + \cdots + 2 \alpha_{n}} e_{j} 
-q^{-1} e_{j}e_{\alpha_{j} + 2\alpha_{j+1} + \cdots + 2 \alpha_{n}} \right) \\
& & \hspace{10mm}  \otimes (K_{\alpha_{j} + \cdots + \alpha_{n}})^{2} 
\left[e_{\alpha_{i} + \cdots + \alpha_{j-1}},e_{\alpha_{i} + \cdots + 2\alpha_{j} + \cdots + 2\alpha_{n}} \right]_{q} \\
& = & q^{2}F_{j-1}\overline{D} 
 + (q-q^{-1}) \left(q e_{\alpha_{j} + 2\alpha_{j+1} + \cdots + 2 \alpha_{n}} e_{j} 
-q^{-1} e_{j}e_{\alpha_{j} + 2\alpha_{j+1} + \cdots + 2 \alpha_{n}} \right) \\
& & \hspace{50mm}  \otimes (K_{\alpha_{j} + \cdots + \alpha_{n}})^{2} 
\left[e_{\alpha_{i} + \cdots + \alpha_{j-1}},e_{\alpha_{i} + \cdots + 2\alpha_{j} + \cdots + 2\alpha_{n}} \right]_{q}. 
\end{eqnarray*}

\noindent
The relation between $D_{k}$ and $F_{k}$ for each $k=j, \ldots, n-1$ is:
\begin{eqnarray*}
D_{k} F_{k} & = & 
 (q-q^{-1})^{2} \overline{e}_{\alpha_{j} + \cdots + 2\alpha_{k+1} + \cdots + 2\alpha_{n}}
  \overline{e}_{\alpha_{j} + \cdots + \alpha_{k}} \\
  & & \hspace{25mm} \otimes 
  K_{\alpha_{j} + \cdots + 2\alpha_{k+1} + \cdots + 2\alpha_{n}}e_{\alpha_{i} + \cdots + \alpha_{k}}
  K_{\alpha_{j} + \cdots + \alpha_{k}}e_{\alpha_{i} + \cdots + 2\alpha_{k+1} + \cdots + 2\alpha_{n}} \\
& = & (q-q^{-1})^{2} \overline{e}_{\alpha_{j} + \cdots + \alpha_{k}}\overline{e}_{\alpha_{j} + \cdots + 2\alpha_{k+1} + \cdots + 2\alpha_{n}}
\\
& & \hspace{25mm} \otimes 
  K_{\alpha_{j} + \cdots + 2\alpha_{k+1} + \cdots + 2\alpha_{n}}e_{\alpha_{i} + \cdots + \alpha_{k}}
  K_{\alpha_{j} + \cdots + \alpha_{k}}e_{\alpha_{i} + \cdots + 2\alpha_{k+1} + \cdots + 2\alpha_{n}} \\
& & + (q-q^{-1})^{2} \left[\overline{e}_{\alpha_{j} + \cdots + 2\alpha_{k+1} + \cdots + 2\alpha_{n}}
  \overline{e}_{\alpha_{j} + \cdots + \alpha_{k}}\right]_{q} \\
  & & \hspace{25mm} \otimes 
  K_{\alpha_{j} + \cdots + 2\alpha_{k+1} + \cdots + 2\alpha_{n}}e_{\alpha_{i} + \cdots + \alpha_{k}}
  K_{\alpha_{j} + \cdots + \alpha_{k}}e_{\alpha_{i} + \cdots + 2\alpha_{k+1} + \cdots + 2\alpha_{n}} \\
& = & q^{-1} (q-q^{-1})^{2} 
\overline{e}_{\alpha_{j} + \cdots + \alpha_{k}}\overline{e}_{\alpha_{j} + \cdots + 2\alpha_{k+1} + \cdots + 2\alpha_{n}}
\\
& & \hspace{25mm} \otimes  (K_{\alpha_{j} + \cdots + \alpha_{n}})^{2}
  e_{\alpha_{i} + \cdots + \alpha_{k}}
  e_{\alpha_{i} + \cdots + 2\alpha_{k+1} + \cdots + 2\alpha_{n}} \\
& & + q^{-1} (q-q^{-1})^{2}  \left[\overline{e}_{\alpha_{j} + \cdots + 2\alpha_{k+1} + \cdots + 2\alpha_{n}}
  \overline{e}_{\alpha_{j} + \cdots + \alpha_{k}}\right]_{q} \\
  & & \hspace{25mm} \otimes 
  (K_{\alpha_{j} + \cdots + \alpha_{n}})^{2}e_{\alpha_{i} + \cdots + \alpha_{k}}
  e_{\alpha_{i} + \cdots + 2\alpha_{k+1} + \cdots + 2\alpha_{n}} \\
& = & (q-q^{-1})^{2} \overline{e}_{\alpha_{j} + \cdots + \alpha_{k}}
      \overline{e}_{\alpha_{j} + \cdots + 2\alpha_{k+1} + \cdots + 2\alpha_{n}} \\
& & \hspace{25mm} \otimes  
  K_{\alpha_{j} + \cdots + \alpha_{k}}
  e_{\alpha_{i} + \cdots + 2\alpha_{k+1} + \cdots + 2\alpha_{n}}
  K_{\alpha_{j} + \cdots + 2\alpha_{k+1} + \cdots + 2\alpha_{n}}
  e_{\alpha_{i} + \cdots + \alpha_{k}}  \\
& & +q^{-1} (q-q^{-1})^{2}  \overline{e}_{\alpha_{j} + \cdots + \alpha_{k}}
      \overline{e}_{\alpha_{j} + \cdots + 2\alpha_{k+1} + \cdots + 2\alpha_{n}} \\ 
& & \hspace{25mm} \otimes
    (K_{\alpha_{j} + \cdots + \alpha_{n}})^{2} 
    \left[e_{\alpha_{i} + \cdots + \alpha_{k}}, e_{\alpha_{i} + \cdots + 2\alpha_{k+1} + \cdots + 2\alpha_{n}} \right]_{q} \\
& & +q^{-1} (q-q^{-1})^{2} 
\left[ \overline{e}_{\alpha_{j} + \cdots + 2\alpha_{k+1} + \cdots + 2\alpha_{n}},
     \overline{e}_{\alpha_{j} + \cdots + \alpha_{k}} \right]_{q} \\
 & & \hspace{25mm} \otimes (K_{\alpha_{j} + \cdots + \alpha_{n}})^{2} 
    e_{\alpha_{i} + \cdots + \alpha_{k}}e_{\alpha_{i} + \cdots + 2\alpha_{k+1} + \cdots + 2\alpha_{n}} \\
& = & F_{k} D_{k} 
 +q^{-1} (q-q^{-1})^{2}  \overline{e}_{\alpha_{j} + \cdots + \alpha_{k}}
      \overline{e}_{\alpha_{j} + \cdots + 2\alpha_{k+1} + \cdots + 2\alpha_{n}} \\ 
& & \hspace{45mm} \otimes
    (K_{\alpha_{j} + \cdots + \alpha_{n}})^{2} 
    \left[e_{\alpha_{i} + \cdots + \alpha_{k}}, e_{\alpha_{i} + \cdots + 2\alpha_{k+1} + \cdots + 2\alpha_{n}} \right]_{q} \\
& & \hspace{11mm} +q^{-1} (q-q^{-1})^{2} 
\left[ \overline{e}_{\alpha_{j} + \cdots + 2\alpha_{k+1} + \cdots + 2\alpha_{n}},
     \overline{e}_{\alpha_{j} + \cdots + \alpha_{k}} \right]_{q} \\
 & & \hspace{45mm} \otimes (K_{\alpha_{j} + \cdots + \alpha_{n}})^{2} 
    e_{\alpha_{i} + \cdots + \alpha_{k}}e_{\alpha_{i} + \cdots + 2\alpha_{k+1} + \cdots + 2\alpha_{n}}.
\end{eqnarray*}

\noindent
We can rewrite this relation as
\begin{eqnarray*}
D_{k} F_{k} 
& = & q^{2} F_{k} D_{k} +q^{-1}(1-q^{2})(q-q^{-1})^{2} \overline{e}_{\alpha_{j} + \cdots + \alpha_{k}}
      \overline{e}_{\alpha_{j} + \cdots + 2\alpha_{k+1} + \cdots + 2\alpha_{n}} \\
& & \hspace{25mm} \otimes  
  (K_{\alpha_{j} + \cdots + \alpha_{n}})^{2}
  e_{\alpha_{i} + \cdots + 2\alpha_{k+1} + \cdots + 2\alpha_{n}}e_{\alpha_{i} + \cdots + \alpha_{k}}  \\
& & +q^{-1} (q-q^{-1})^{2}  \overline{e}_{\alpha_{j} + \cdots + \alpha_{k}}
      \overline{e}_{\alpha_{j} + \cdots + 2\alpha_{k+1} + \cdots + 2\alpha_{n}} \\ 
& & \hspace{25mm} \otimes
    (K_{\alpha_{j} + \cdots + \alpha_{n}})^{2} 
    \left[e_{\alpha_{i} + \cdots + \alpha_{k}}, e_{\alpha_{i} + \cdots + 2\alpha_{k+1} + \cdots + 2\alpha_{n}} \right]_{q} \\
& & +q^{-1} (q-q^{-1})^{2} 
\left[ \overline{e}_{\alpha_{j} + \cdots + 2\alpha_{k+1} + \cdots + 2\alpha_{n}},
     \overline{e}_{\alpha_{j} + \cdots + \alpha_{k}} \right]_{q} \\
 & & \hspace{25mm} \otimes (K_{\alpha_{j} + \cdots + \alpha_{n}})^{2} 
    e_{\alpha_{i} + \cdots + \alpha_{k}}e_{\alpha_{i} + \cdots + 2\alpha_{k+1} + \cdots + 2\alpha_{n}} \\
& = & q^{2} F_{k} D_{k}  + (q-q^{-1})^{2} \Big[ (q^{-1}-q)
          \overline{e}_{\alpha_{j} + \cdots + \alpha_{k}}
	  \overline{e}_{\alpha_{j} + \cdots + 2\alpha_{k+1} + \cdots + 2\alpha_{n}} \\
& & \hspace{45mm}	 +q^{-1} \left[ \overline{e}_{\alpha_{j} + \cdots + 2\alpha_{k+1} + \cdots + 2\alpha_{n}},
	             \overline{e}_{\alpha_{j} + \cdots + \alpha_{k}} \right]_{q} \Big] \\
& & \hspace{60mm}  \otimes (K_{\alpha_{j} + \cdots + \alpha_{n}})^{2} 
	e_{\alpha_{i} + \cdots + \alpha_{k}}e_{\alpha_{i} + \cdots + 2\alpha_{k+1} + \cdots + 2\alpha_{n}} \\   
& & + q(q-q^{-1})^{2} \overline{e}_{\alpha_{j} + \cdots + \alpha_{k}}
                      \overline{e}_{\alpha_{j} + \cdots + 2\alpha_{k+1} + \cdots + 2\alpha_{n}} \\
 & & \hspace{25mm} \otimes (K_{\alpha_{j} + \cdots + \alpha_{n}})^{2}
   \left[ e_{\alpha_{i} + \cdots + \alpha_{k}}, 
     e_{\alpha_{i} + \cdots + 2\alpha_{k+1} + \cdots + 2\alpha_{n}} \right]_{q}.
\end{eqnarray*}

From these calculations we have the following relation:
\begin{eqnarray}
\lefteqn{
\big(\overline{D} + D_{j} + D_{j+1} + \cdots + D_{n-1}\big)
\big(F_{n-1}+F_{n-2} + \cdots + F_{j-1}\big) } \nonumber \\
& = & q^{2} \big(F_{n-1}+F_{n-2} + \cdots + F_{j-1}\big)
            \big(\overline{D} + D_{j} + D_{j+1} + \cdots + D_{n-1}\big) \nonumber \\
&   & + (q-q^{-1}) \left( q e_{\alpha_{j}+2\alpha_{j+1}+\cdots+2\alpha_{n}} e_{j}
                         -q^{-1} e_{j}e_{\alpha_{j}+2\alpha_{j+1}+\cdots+2\alpha_{n}} \right) \nonumber \\
& & \hspace{30mm} \otimes (K_{\alpha_{j} + \cdots + \alpha_{n}})^{2}
      \left[e_{\alpha_{i}+\cdots+\alpha_{j-1}},e_{\alpha_{i}+\cdots+2\alpha_{j}+\cdots+2\alpha_{n}} \right]_{q} 
      \label{appendixB:amos(1)} \\
& &  +(q-q^{-1})^{2} \\
& & \hspace{5mm} \times \sum_{k=j}^{n-1} \Bigg(
	\Big[ (q^{-1}-q) \overline{e}_{\alpha_{j}+\cdots+\alpha_{k}}
	                 \overline{e}_{\alpha_{j}+\cdots+2\alpha_{k+1}+\cdots+2\alpha_{n}} \nonumber \\
& & \hspace{25mm}
	    +q^{-1} \left[ \overline{e}_{\alpha_{j}+\cdots+2\alpha_{k+1}+\cdots+2\alpha_{n}},
	                    \overline{e}_{\alpha_{j}+\cdots+\alpha_{k}} \right]_{q}  \Big] \nonumber \\
& & \hspace{50mm} \otimes (K_{\alpha_{j}+\cdots+\alpha_{n}})^{2}
                    e_{\alpha_{i}+\cdots+\alpha_{k}}e_{\alpha_{i}+\cdots+2\alpha_{k+1}+\cdots+2\alpha_{n}} 
		    \label{appendixB:amos(2)} \\
& & \hspace{21mm} + q \overline{e}_{\alpha_{j}+\cdots+\alpha_{k}}
\overline{e}_{\alpha_{j}+\cdots+2\alpha_{k+1}+\cdots+2\alpha_{n}} \nonumber \\
& & \hspace{50mm} \otimes (K_{\alpha_{j}+\cdots+\alpha_{n}})^{2}
  \left[  e_{\alpha_{i}+\cdots+\alpha_{k}},e_{\alpha_{i}+\cdots+2\alpha_{k+1}+\cdots+2\alpha_{n}}\right]_{q} \Bigg).
   \nonumber \\
& &   \label{appendixB:amos(3)}
\end{eqnarray}
To simplify this expression we will expand out the terms in the second tensor power 
and sum the terms accordingly.
The easiest part of this calculation is as follows:
the component in the expansion for which the second tensor power contains a term of the form
$e_{\alpha_{i}+\cdots+\alpha_{j}}e_{\alpha_{i}+\cdots+2\alpha_{j+1}+\cdots+2\alpha_{n}}$, is
\begin{eqnarray}
\lefteqn{
(q-q^{-1})^{2} \Big[ q e_{\alpha_{j}+2\alpha_{j+1}+\cdots+2\alpha_{n}} e_{j}
                         -q^{-1} e_{j}e_{\alpha_{j}+2\alpha_{j+1}+\cdots+2\alpha_{n}} } \nonumber \\
& & \hspace{20mm}   +q^{-1} \overline{e}_{\alpha_{j}+2\alpha_{j+1}+\cdots+2\alpha_{n}} e_{j}
                          -q e_{j} \overline{e}_{\alpha_{j}+2\alpha_{j+1}+\cdots+2\alpha_{n}} \Big] \nonumber \\
& & \hspace{40mm} \otimes (K_{\alpha_{j} + \cdots + \alpha_{n}})^{2}
e_{\alpha_{i}+\cdots+\alpha_{j}}e_{\alpha_{i}+\cdots+2\alpha_{j+1}+\cdots+2\alpha_{n}} \label{appendixB:troythebattleeq} \\
& = & 0, \nonumber
\end{eqnarray}
where  (\ref{appendixB:troythebattleeq}) 
arises from (\ref{appendixB:amos(1)}) and the $k=j$ term in (\ref{appendixB:amos(2)}), and
Proposition \ref{appendixB:troythebattleprop} then implies that 
(\ref{appendixB:troythebattleeq}) vanishes.

To simplify the remaining terms in the above expression we need to know the expansions of
$$\left[e_{\alpha_{i}+\cdots+\alpha_{j-1}},
e_{\alpha_{i}+\cdots+2\alpha_{j}+\cdots+2\alpha_{n}} \right]_{q}
\hspace{5mm} \mbox{and} \hspace{5mm}
\left[ \overline{e}_{\alpha_{j}+\cdots+2\alpha_{k+1}+\cdots+2\alpha_{n}},
	                    \overline{e}_{\alpha_{j}+\cdots+\alpha_{k}} \right]_{q}$$
which we have detailed in Propositions \ref{appendixB:propamos(1)} and \ref{appendixB:propamos(2)} respectively.
Using these propositions, we can determine
  the component in the expansion of the above expression for which the second tensor power contains a term of the form
$e_{\alpha_{i}+\cdots+\alpha_{j+p}}e_{\alpha_{i}+\cdots+2\alpha_{j+p+1}+\cdots+2\alpha_{n}}$
for each $p=1, \ldots, n-j-1$.  This component is
\begin{eqnarray*}
\lefteqn{ \Big[
(-q)^{p}(q-q^{-1})^{2} 
\left( q e_{\alpha_{j}+2\alpha_{j+1}+\cdots+2\alpha_{n}} e_{j} -q^{-1} e_{j} e_{\alpha_{j}+2\alpha_{j+1}+\cdots+2\alpha_{n}}
\right) } \\
& &  + \sum_{k=j}^{j+p-1} (-q)^{j+p-k-1} q (q-q^{-1})^{3} 
 \overline{e}_{\alpha_{j}+\cdots+\alpha_{k}}
                      \overline{e}_{\alpha_{j}+\cdots+2\alpha_{k+1}+\cdots+2\alpha_{n}} \\
& & + (q-q^{-1})^{2} \left( -q \overline{e}_{\alpha_{j}+\cdots+\alpha_{j+p}}
            \overline{e}_{\alpha_{j}+\cdots+2\alpha_{j+p+1}+\cdots+2\alpha_{n}} 
	    +q^{-1} \overline{e}_{\alpha_{j}+\cdots+2\alpha_{j+p+1}+\cdots+2\alpha_{n}} 
	    \overline{e}_{\alpha_{j}+\cdots+\alpha_{j+p}} \right) \Big] \\
& & \hspace{50mm}  
    \otimes (K_{\alpha_{j}+\cdots+\alpha_{n}})^{2} e_{\alpha_{i}+\cdots+\alpha_{j+p}}
         e_{\alpha_{i}+\cdots+2\alpha_{j+p+1}+\cdots+2\alpha_{n}}.
\end{eqnarray*}
The element of $U_{q}(\mathfrak{g})$ in the first tensor power of this component is
\begin{eqnarray*}
\lefteqn{
(q-q^{-1})^{2}  \Big[ (-q)^{p} q e_{j} \overline{e}_{\alpha_{j}+2\alpha_{j+1}+\cdots+2\alpha_{n}}
               -q^{-1} (-q)^{p} \overline{e}_{\alpha_{j}+2\alpha_{j+1}+\cdots+2\alpha_{n}} e_{j} } \\
& & \hspace{15mm} + (-q)^{p-1} q (q-q^{-1}) e_{j} \overline{e}_{\alpha_{j}+2\alpha_{j+1}+\cdots+2\alpha_{n}} \\
& & \hspace{15mm} + \sum_{k=j+1}^{j+p-1} (-q)^{j+p-k-1} q (q-q^{-1}) 
   \overline{e}_{\alpha_{j}+\cdots+\alpha_{k}}\overline{e}_{\alpha_{j}+\cdots+2\alpha_{k+1}+\cdots+2\alpha_{n}} \\
& & \hspace{15mm} -q \overline{e}_{\alpha_{j}+\cdots+\alpha_{j+p}}\overline{e}_{\alpha_{j}+\cdots+2\alpha_{j+p+1}+\cdots+2\alpha_{n}}
    +q^{-1}\overline{e}_{\alpha_{j}+\cdots+2\alpha_{j+p+1}+\cdots+2\alpha_{n}} \overline{e}_{\alpha_{j}+\cdots+\alpha_{j+p}}
    \Big],
\end{eqnarray*}
which we can rewrite as
\begin{eqnarray}
\lefteqn{
(q-q^{-1})^{2}  
  \Big[ (-1)^{p} q^{p-1} e_{j} \overline{e}_{\alpha_{j}+2\alpha_{j+1}+\cdots+2\alpha_{n}}
        -(-1)^{p} q^{p-1} \overline{e}_{\alpha_{j}+2\alpha_{j+1}+\cdots+2\alpha_{n}} e_{j} } \label{appendixC:margaratT(1)} \\	       
& & \hspace{15mm} + \sum_{k=j+1}^{j+p-1} (-q)^{j+p-k-1} q (q-q^{-1}) 
   \overline{e}_{\alpha_{j}+\cdots+\alpha_{k}}\overline{e}_{\alpha_{j}+\cdots+2\alpha_{k+1}+\cdots+2\alpha_{n}} \nonumber  \\
& & \hspace{15mm} -q \overline{e}_{\alpha_{j}+\cdots+\alpha_{j+p}}\overline{e}_{\alpha_{j}+\cdots+2\alpha_{j+p+1}+\cdots+2\alpha_{n}}
    +q^{-1}\overline{e}_{\alpha_{j}+\cdots+2\alpha_{j+p+1}+\cdots+2\alpha_{n}} \overline{e}_{\alpha_{j}+\cdots+\alpha_{j+p}}
    \Big]. \nonumber
\end{eqnarray}

Now we calculate that
\begin{eqnarray*}
\lefteqn{ 
(-1)^{p} q^{p-1} \left( e_{j} \overline{e}_{\alpha_{j}+2\alpha_{j+1}+\cdots+2\alpha_{n}}
        - \overline{e}_{\alpha_{j}+2\alpha_{j+1}+\cdots+2\alpha_{n}} e_{j}\right) } \\
& = & (-q)^{p-1} \left[ \overline{e}_{\alpha_{j}+2\alpha_{j+1}+\cdots+2\alpha_{n}}, e_{j} \right]_{q} \\
& = & (-q)^{p-1} \left[ e_{j+1} \overline{e}_{\alpha_{j}+\alpha_{j+1}+2\alpha_{j+2}+\cdots+2\alpha_{n}}
                -q^{-1} \overline{e}_{\alpha_{j}+\alpha_{j+1}+2\alpha_{j+2}+\cdots+2\alpha_{n}} e_{j+1}, e_{j} \right]_{q} \\
& = & (-q)^{p-1}  
    \big( q \left[e_{j+1},e_{j}\right]_{q} \overline{e}_{\alpha_{j}+\alpha_{j+1}+2\alpha_{j+2}+\cdots+2\alpha_{n}}
     -q^{-1} \overline{e}_{\alpha_{j}+\alpha_{j+1}+2\alpha_{j+2}+\cdots+2\alpha_{n}}\left[e_{j+1},e_{j}\right]_{q} \big),
\end{eqnarray*}
and substituting this into (\ref{appendixC:margaratT(1)}) gives 
\begin{eqnarray}
\lefteqn{
(q-q^{-1})^{2}  } \nonumber \\
&  & \times \Big[ (-q)^{p-1} 
\left( q \overline{e}_{\alpha_{j}+\alpha_{j+1}} \overline{e}_{\alpha_{j}+\alpha_{j+1}+2\alpha_{j+2}+\cdots+2\alpha_{n}}
-q^{-1} \overline{e}_{\alpha_{j}+\alpha_{j+1}+2\alpha_{j+2}+\cdots+2\alpha_{n}}\overline{e}_{\alpha_{j}+\alpha_{j+1}} \right) 
\nonumber  \\
& & \hspace{10mm} + (-q)^{j+p-j-2} q (q-q^{-1}) \overline{e}_{\alpha_{j}+\alpha_{j+1}}
                                  \overline{e}_{\alpha_{j}+\alpha_{j+1}+2\alpha_{j+2}+\cdots+2\alpha_{n}} \nonumber \\
& & \hspace{10mm}  + \sum_{k=j+2}^{j+p-1} (-q)^{j+p-k-1} (q^{2}-1) \overline{e}_{\alpha_{j}+\cdots+\alpha_{k}}
                                  \overline{e}_{\alpha_{j}+\cdots+2\alpha_{k+1}+\cdots+2\alpha_{n}} \nonumber \\
& & \hspace{10mm} -q \overline{e}_{\alpha_{j}+\cdots+\alpha_{j+p}}\overline{e}_{\alpha_{j}+\cdots+2\alpha_{j+p+1}+\cdots+2\alpha_{n}}
    +q^{-1}\overline{e}_{\alpha_{j}+\cdots+2\alpha_{j+p+1}+\cdots+2\alpha_{n}} \overline{e}_{\alpha_{j}+\cdots+\alpha_{j+p}}
    \Big] \nonumber  \\
& = & (q-q^{-1})^{2} \label{eq:franceseveleigh(1)} \\
& & \times \Big[ (-q)^{p-2}  
\left[\overline{e}_{\alpha_{j}+\alpha_{j+1}+2\alpha_{j+2}+\cdots+2\alpha_{n}},\overline{e}_{\alpha_{j}+\alpha_{j+1}}\right]_{q}
 \label{appendixB:charles(1)} \\
& & \hspace{10mm}  + \sum_{k=j+2}^{j+p-1} (-q)^{j+p-k-1} (q^{2}-1) \overline{e}_{\alpha_{j}+\cdots+\alpha_{k}}
                                  \overline{e}_{\alpha_{j}+\cdots+2\alpha_{k+1}+\cdots+2\alpha_{n}}
				  \label{appendixB:charles(2)} \\
& & \hspace{10mm} -q \overline{e}_{\alpha_{j}+\cdots+\alpha_{j+p}}\overline{e}_{\alpha_{j}+\cdots+2\alpha_{j+p+1}+\cdots+2\alpha_{n}}
    +q^{-1}\overline{e}_{\alpha_{j}+\cdots+2\alpha_{j+p+1}+\cdots+2\alpha_{n}} \overline{e}_{\alpha_{j}+\cdots+\alpha_{j+p}}
    \Big]. \nonumber \\
 & &    \label{eq:franceseveleigh(2)} 
\end{eqnarray}
We can dramatically simplify this expression by using the following calculation:
\begin{eqnarray}
\lefteqn{
(-q)^{p-m}\left[ \overline{e}_{\alpha_{j}+\cdots+2\alpha_{j+m}+\cdots+2\alpha_{n}},
                                   \overline{e}_{\alpha_{j}+\cdots+\alpha_{j+m-1}} \right]_{q} } \nonumber \\
& &  + (-q)^{p-m-1} (q^{2}-1) \overline{e}_{\alpha_{j}+\cdots+\alpha_{j+m}}
                         \overline{e}_{\alpha_{j}+\cdots+2\alpha_{j+m+1}+\cdots+2\alpha_{n}}  \nonumber \\
& = & (-q)^{p-m} \left[ e_{j+m} \overline{e}_{\alpha_{j}+\cdots+2\alpha_{j+m+1}+\cdots+2\alpha_{n}} 
              -q^{-1} \overline{e}_{\alpha_{j}+\cdots+2\alpha_{j+m+1}+\cdots+2\alpha_{n}}e_{j+m},
	      \overline{e}_{\alpha_{j}+\cdots+\alpha_{j+m-1}}\right]_{q} \nonumber \\
& &  + (-q)^{p-m-1} (q^{2}-1) \overline{e}_{\alpha_{j}+\cdots+\alpha_{j+m}}
                         \overline{e}_{\alpha_{j}+\cdots+2\alpha_{j+m+1}+\cdots+2\alpha_{n}}  \nonumber \\
& = & (-q)^{p-m}  \left( q \overline{e}_{\alpha_{j}+\cdots+\alpha_{j+m}}
                     \overline{e}_{\alpha_{j}+\cdots+2\alpha_{j+m+1}+\cdots+2\alpha_{n}} 
		   -q^{-1}     \overline{e}_{\alpha_{j}+\cdots+2\alpha_{j+m+1}+\cdots+2\alpha_{n}}
		               \overline{e}_{\alpha_{j}+\cdots+\alpha_{j+m}} \right) \nonumber \\
& & + (-q)^{p-m-1} (q^{2}-1) \overline{e}_{\alpha_{j}+\cdots+\alpha_{j+m}}
                             \overline{e}_{\alpha_{j}+\cdots+2\alpha_{j+m+1}+\cdots+2\alpha_{n}} \nonumber \\
& = & (-q)^{p-m-1} \left[ \overline{e}_{\alpha_{j}+\cdots+2\alpha_{j+m+1}+\cdots+2\alpha_{n}},
                                   \overline{e}_{\alpha_{j}+\cdots+\alpha_{j+m}} \right]_{q}. \label{appendixC:julianardas(3)}
\end{eqnarray}
By repeatedly using this identity in (\ref{appendixB:charles(1)})--(\ref{appendixB:charles(2)}),
we can rewrite (\ref{eq:franceseveleigh(1)})--(\ref{eq:franceseveleigh(2)}) as
the following expression for each $t=2,\ldots,p-1$:
\begin{eqnarray}
\lefteqn{
(q-q^{-1})^{2}  } \nonumber \\
& & \times \Big[ (-q)^{p-t} 
\left[\overline{e}_{\alpha_{j}+\cdots+2\alpha_{j+t}+\cdots+2\alpha_{n}},
       \overline{e}_{\alpha_{j}+\cdots + \alpha_{j+t-1}}\right]_{q} \nonumber \\
& & \hspace{5mm}  + \sum_{k=j+t}^{j+p-1} (-q)^{j+p-k-1} (q^{2}-1) 
      \overline{e}_{\alpha_{j}+\cdots+\alpha_{k}}
      \overline{e}_{\alpha_{j}+\cdots+2\alpha_{k+1}+\cdots+2\alpha_{n}}  \nonumber \\
& & \hspace{5mm} -q \overline{e}_{\alpha_{j}+\cdots+\alpha_{j+p}}\overline{e}_{\alpha_{j}+\cdots+2\alpha_{j+p+1}+\cdots+2\alpha_{n}}
    +q^{-1}\overline{e}_{\alpha_{j}+\cdots+2\alpha_{j+p+1}+\cdots+2\alpha_{n}} \overline{e}_{\alpha_{j}+\cdots+\alpha_{j+p}}
    \Big]. \nonumber \\
& & \label{eq:franceseveleigh(3)}
\end{eqnarray}
Substituting $t=p-1$ into (\ref{eq:franceseveleigh(3)}) and using (\ref{appendixC:julianardas(3)}), 
we can rewrite
(\ref{eq:franceseveleigh(1)})--(\ref{eq:franceseveleigh(2)}) as
\begin{eqnarray*}
\lefteqn{
 (q-q^{-1})^{2}  \Big[ \left[\overline{e}_{\alpha_{j}+\cdots+2\alpha_{j+p}+\cdots+2\alpha_{n}},
       \overline{e}_{\alpha_{j}+\cdots +\alpha_{j+p-1}}\right]_{q} }  \\
& & \hspace{5mm}  
    -q \overline{e}_{\alpha_{j}+\cdots+\alpha_{j+p}}\overline{e}_{\alpha_{j}+\cdots+2\alpha_{j+p+1}+\cdots+2\alpha_{n}}
    +q^{-1}\overline{e}_{\alpha_{j}+\cdots+2\alpha_{j+p+1}+\cdots+2\alpha_{n}} \overline{e}_{\alpha_{j}+\cdots+\alpha_{j+p}}
    \Big]  \\
& = & 0,  
\end{eqnarray*}
where we have used the following result: 
\begin{eqnarray*}
\lefteqn{
\left[\overline{e}_{\alpha_{j}+\cdots+2\alpha_{j+p}+\cdots+2\alpha_{n}},
       \overline{e}_{\alpha_{j}+\cdots +\alpha_{j+p-1}}\right]_{q} } \\
 & = & 
  \left[ e_{j+p} \overline{e}_{\alpha_{j}+\cdots+2\alpha_{j+p+1}+\cdots+2\alpha_{n}}
         -q^{-1} \overline{e}_{\alpha_{j}+\cdots+2\alpha_{j+p+1}+\cdots+2\alpha_{n}}e_{j+p}, 
	          \overline{e}_{\alpha_{j}+\cdots+\alpha_{j+p-1}} \right]_{q} \\
& = & q \overline{e}_{\alpha_{j}+\cdots+\alpha_{j+p}}\overline{e}_{\alpha_{j}+\cdots+2\alpha_{j+p+1}+\cdots+2\alpha_{n}}
     -q^{-1} \overline{e}_{\alpha_{j}+\cdots+2\alpha_{j+p+1}+\cdots+2\alpha_{n}}\overline{e}_{\alpha_{j}+\cdots+\alpha_{j+p}}.
\end{eqnarray*}
This substantially simplifies the problem.  
All we now have to do is to determine the component for which the second tensor power contains a term of the form 
$(e_{\alpha_{j}+\cdots+\alpha_{n}})^{2}$.  This is not difficult to do:
the first tensor power of this component is 
\begin{eqnarray}
\lefteqn{
(-q)^{n-j} (1+q)(q-q^{-1}) \left(q e_{\alpha_{j}+2\alpha_{j+1}+\cdots+2\alpha_{n}} e_{j}
                            -q^{-1} e_{j} e_{\alpha_{j}+2\alpha_{j+1}+\cdots+2\alpha_{n}} \right) } \nonumber \\
& & \hspace{10mm} \sum_{k=j}^{n-1} (-q)^{n-k-1} q(1+q)(q-q^{-1})^{2}
     \overline{e}_{\alpha_{j}+\cdots+\alpha_{k}} \overline{e}_{\alpha_{j}+\cdots+2\alpha_{k+1}+\cdots+2\alpha_{n}}. 
     \label{appendixB:spock(1)}
\end{eqnarray}
Now by repeatedly using (\ref{appendixC:julianardas(3)}), we have 
\begin{eqnarray*}
\lefteqn{
(-q)^{n-j} \left( q e_{\alpha_{j}+2\alpha_{j+1}+\cdots+2\alpha_{n}} e_{j} 
          -q^{-1}  e_{j} e_{\alpha_{j}+2\alpha_{j+1}+\cdots+2\alpha_{n}} \right) } \\
 & & \hspace{10mm} + \sum_{k=j}^{n-1} (-q)^{n-k-1} (q^{2}-1)
         \overline{e}_{\alpha_{j}+\cdots+\alpha_{k}} \overline{e}_{\alpha_{j}+\cdots+2\alpha_{k+1}+\cdots+2\alpha_{n}} \\
& = & \left[\overline{e}_{\alpha_{j}+\cdots+2\alpha_{n}},\overline{e}_{\alpha_{j}+\cdots+\alpha_{n-1}}\right]_{q} \\
& = & (1+q) \left( \overline{e}_{\alpha_{j}+\cdots+\alpha_{n}} \right)^{2},
\end{eqnarray*}
and thus (\ref{appendixB:spock(1)}) is 
$$(1+q)^{2}(q-q^{-1})\left( \overline{e}_{\alpha_{j}+\cdots+\alpha_{n}} \right)^{2}
 = \xi (D_{n})^{2}, \hspace{10mm} \xi = -(1+q)^{2}/(q-q^{-1}),$$ 
 as required.

\end{subsection}

\end{section}

\begin{section}{${\cal{I}}$ is a two-sided co-ideal}

In this section we prove that ${\cal{I}}$ is a two-sided co-ideal.
We firstly note the almost trivial result that $\epsilon(x)=0$ for each $x \in I$ and thus
$\epsilon(x)=0$ for all $x \in {\cal{I}}$.
We now deal with the more substantial problem: we will prove that
$$\Delta(x) \in {\cal{I}} \otimes U_{q}(\mathfrak{g}) + U_{q}(\mathfrak{g}) \otimes {\cal{I}},$$
for each $x \in I$.  
Firstly, $$\Delta(J_{i}^{\pm N}-1) = J_{i}^{\pm N} \otimes J_{i}^{\pm N} - 1 \otimes 1 = 
J_{i}^{\pm N} \otimes (J_{i}^{\pm N}-1) + (J_{i}^{\pm N}-1) \otimes 1,$$
which is an element of ${\cal{I}} \otimes U_{q}(\mathfrak{g}) + U_{q}(\mathfrak{g}) \otimes {\cal{I}}$.
The remaining problems are harder and we break them down into a number of subcases.
Initially we will show that $(e_{\gamma})^{N'}$ and $(e_{\beta})^{\overline{N}}$ are elements of 
${\cal{I}} \otimes U_{q}(\mathfrak{g}) + U_{q}(\mathfrak{g}) \otimes {\cal{I}}$; 
then we will use the antipode $S$ and the
graded antiautomorphism $\omega$ to prove the remaining cases.

In writing down the components of $\Delta(e_{\mu})$ we will use the notation used in 
Section \ref{appendixC:arnietermincarrotator(1)}, 
although we may use slightly different normalisations. 
Any alternative normalisations will not significantly affect the calculations.

\begin{subsubsection}{Case 1. $(e_{i})^{N'}$,  $(e_{n})^{\overline{N}}$, $1 \leq i<n$ }

The components of $\Delta(e_{i})$ satisfy the $q$-binomial theorem for each $i$, thus 
we have the following results from Appendix B:
\begin{eqnarray*}
\big(\Delta(e_{i})\big)^{N'} & = & (e_{i} \otimes K_{i})^{N'} + (1 \otimes e_{i})^{N'} 
= (e_{i})^{N'} \otimes (K_{i})^{N'} + 1 \otimes (e_{i})^{N'}, \hspace{5mm} i<n, \\
\big(\Delta(e_{n})\big)^{\overline{N}} 
& = & (e_{n} \otimes K_{n})^{\overline{N}} + (1 \otimes e_{n})^{\overline{N}} 
= (e_{n})^{\overline{N}} \otimes (K_{n})^{\overline{N}} + 1 \otimes (e_{n})^{\overline{N}}.
\end{eqnarray*}
As $(e_{i})^{N'} \in {\cal{I}}$ and $(e_{n})^{\overline{N}} \in {\cal{I}}$, 
$\big(\Delta(e_{i})\big)^{N'}$ and $\big(\Delta(e_{n})\big)^{\overline{N}}$ are elements of 
${\cal{I}} \otimes U_{q}(\mathfrak{g}) + U_{q}(\mathfrak{g}) \otimes {\cal{I}}$.

\end{subsubsection}

\begin{subsubsection}{Case 2. $(e_{\alpha_{i} + \cdots + \alpha_{j}})^{N'}$, $1 \leq i<j<n$ }

By writing 
$\Delta(e_{\alpha_{i} + \cdots + \alpha_{j}}) = \sum_{k=i}^{j} D_{k} + D_{\infty}$ 
and noting the  relations
$D_{r} D_{s} = q^{2} D_{s} D_{r}$ for all $r < s$
we can use the $q$-multinomial theorem to immediately obtain
$$\big(\Delta(e_{\alpha_{i} + \cdots + \alpha_{j}})\big)^{N'} = 
\sum_{k=i}^{j} (D_{k})^{N'} + (D_{\infty})^{N'}.$$
As $(e_{\alpha_{k} + \cdots + \alpha_{j}})^{N'} \in I$ 
for all $1 \leq k \leq j < n$, we have that
$(D_{k})^{N'}$ and $(D_{\infty})^{N'}$ belong to 
${\cal{I}} \otimes U_{q}(\mathfrak{g}) + U_{q}(\mathfrak{g}) \otimes {\cal{I}}$ for all $i \leq k \leq j$, and thus
$$\big(\Delta(e_{\alpha_{i} + \cdots + \alpha_{j}})\big)^{N'} \in
{\cal{I}} \otimes U_{q}(\mathfrak{g}) + U_{q}(\mathfrak{g}) \otimes {\cal{I}}.$$

\end{subsubsection}

\begin{subsubsection}{Case 3. $(e_{\alpha_{i} + \cdots + \alpha_{n}})^{\overline{N}}$, $i=1, \ldots, n-1$}

Recall that 
$\Delta(e_{\alpha_{i} + \cdots + \alpha_{n}}) = \sum_{k=i}^{n} D_{k} + D_{\infty}$.
We claim that
\begin{eqnarray*}
\lefteqn{
\big(\Delta\left(e_{\alpha_{i} + \cdots + \alpha_{n}}\right)\big)^{\overline{N}} 
   = (D_{\infty})^{\overline{N}} } \\
& &  + \left\{
\begin{array}{ll}
\displaystyle{
\sum_{k=i}^{n} (D_{k})^{2N} + (1)_{q} (3)_{q} \cdots (2N-1)_{q} \left[\sum_{i \leq r<t \leq n} (E_{r,t})^{N} \right] },
 & \mbox{if } N \equiv 1,3 \pmod{4}, \\
\displaystyle{ 
\sum_{k=i}^{n} (D_{k})^{N} + (1)_{q} (3)_{q} \cdots (N-1)_{q} \left[\sum_{i \leq r<t \leq n} (E_{r,t})^{N/2} \right] },
 & \mbox{if } N \equiv 0 \pmod{4}, \\
\displaystyle{ 
\sum_{k=i}^{n} (D_{k})^{N/2} }, & \mbox{if } N \equiv 2 \pmod{4}.
\end{array}  \right.  
\end{eqnarray*}
Firstly, from the relation 
$D_{k} D_{\infty} = -q D_{\infty}D_{k}$ for all $i \leq k < \infty$, we have  
$$\left(\sum_{k=i}^{n} D_{k} + D_{\infty}\right)^{\overline{N}}=
(D_{\infty})^{\overline{N}}+\left(\sum_{k=i}^{n} D_{k}\right)^{\overline{N}},$$
and we now use induction to obtain an expression for $\big(\sum_{k=i}^{n} D_{k}\big)^{\overline{N}}$.
The elements $D_{i}$ and $D_{i+1}$ satisfy the relations
$$D_{i} D_{i+1} = -q D_{i+1}D_{i} + E_{i,i+1}, \hspace{5mm}  D_{i} E_{i,i+1} = q^{2} E_{i,i+1}D_{i},  \hspace{5mm} 
E_{i,i+1}D_{i+1} = q^{2} D_{i+1}E_{i,i+1},$$ 
and therefore 
\begin{eqnarray*}
\lefteqn{
\big(D_{i} + D_{i+1}\big)^{\overline{N}} } \\
& = & 
\left\{
\begin{array}{ll}
\displaystyle{
 (D_{i})^{2N} + (D_{i+1})^{2N} + (1)_{q} (3)_{q} \cdots (2N-1)_{q} (E_{i,i+1})^{N} },
 & \mbox{if } N \equiv 1,3 \pmod{4}, \\
\displaystyle{ 
(D_{i})^{N} + (D_{i+1})^{N} +(1)_{q} (3)_{q} \cdots (N-1)_{q} (E_{i,i+1})^{N/2} },
 & \mbox{if } N \equiv 0 \pmod{4}, \\
\displaystyle{ 
(D_{i})^{N/2} + (D_{i+1})^{N/2} }, & \mbox{if } N \equiv 2 \pmod{4}.
\end{array}  \right.  
\end{eqnarray*}
Using induction, it is quite simple to prove the following expansion for each $j=i+1, \ldots, n$:
\begin{eqnarray*}
\lefteqn{
\left(\sum_{k=i}^{j}D_{k}\right)^{\overline{N}} = } \\
& = & 
\left\{
\begin{array}{ll}
\displaystyle{
\sum_{k=i}^{j} (D_{k})^{2N} + (1)_{q} (3)_{q} \cdots (2N-1)_{q}\left[\sum_{i \leq r<t \leq j} (E_{r,t})^{N} \right] },
 & \mbox{if } N \equiv 1,3 \pmod{4}, \\
\displaystyle{ 
\sum_{k=i}^{j} (D_{k})^{N} +(1)_{q} (3)_{q} \cdots (N-1)_{q}\left[\sum_{i \leq r<t \leq j} (E_{r,t})^{N/2} \right] },
 & \mbox{if } N \equiv 0 \pmod{4}, \\
\displaystyle{ 
\sum_{k=i}^{j} (D_{k})^{N/2} }, & \mbox{if } N \equiv 2 \pmod{4},
\end{array}  \right.  
\end{eqnarray*}
where for the inductive step we use the relations
\begin{eqnarray*}
\lefteqn{
(D_{i} + D_{i+1} + \cdots + D_{j})D_{j+1} } \\
& & =  -q D_{j+1}(D_{i} + D_{i+1} + \cdots + D_{j}) + (E_{i,j+1} + E_{i+1,j+1} + \cdots + E_{j,j+1}), \\
\lefteqn{
(D_{i} + D_{i+1} + \cdots + D_{j})(E_{i,j+1} + E_{i+1,j+1} + \cdots + E_{j,j+1}) } \\
&  & =  q^{2}(E_{i,j+1} + E_{i+1,j+1} + \cdots + E_{j,j+1})(D_{i} + D_{i+1} + \cdots + D_{j}), \\
\lefteqn{
(E_{i,j+1} + E_{i+1,j+1} + \cdots + E_{j,j+1})D_{j+1} } \\
& & =    q^{2}D_{j+1}(E_{i,j+1} + E_{i+1,j+1} + \cdots + E_{j,j+1}).
\end{eqnarray*}
This proves the claimed expression for 
$\big(\Delta (e_{\alpha_{i} + \cdots + \alpha_{n}})\big)^{\overline{N}}$.

Now  $\left(e_{\alpha_{k} + \cdots + \alpha_{n}}\right)^{\overline{N}} \in {\cal{I}}$
for each $k=i, \ldots, n$ and thus
$$(D_{k})^{\overline{N}} \in {\cal{I}} \otimes U_{q}(\mathfrak{g}) + U_{q}(\mathfrak{g}) \otimes {\cal{I}},
 \hspace{10mm}  i \leq k \leq \infty.$$
It remains to show that 
$(E_{r,t})^{N'} \in {\cal{I}} \otimes U_{q}(\mathfrak{g}) + U_{q}(\mathfrak{g}) \otimes {\cal{I}}$
when $N \equiv 0, 1, 3 \pmod{4}$ for all $i \leq r < t \leq n$.
To show this, we note that the second tensor power of $E_{i,p}$ is
$$K_{\alpha_{i} + \cdots + \alpha_{n}}K_{\alpha_{p} + \cdots + \alpha_{n}}
   e_{\alpha_{i} + \cdots + \alpha_{p-1}}, \hspace{10mm} p=i+1, \ldots, n,$$ 
   and the fact that
 $\left(e_{\alpha_{i} + \cdots + \alpha_{p-1}}\right)^{N'} \in {\cal{I}}$ means that
$$(E_{i,p})^{N'} \in {\cal{I}} \otimes U_{q}(\mathfrak{g}) + U_{q}(\mathfrak{g}) \otimes {\cal{I}}, 
\hspace{10mm} p=i+1, \ldots, n.$$
The second tensor power of $E_{k,p}$ for each $k=i+1, \ldots, n-1$ and each $p=k+1, \ldots n$ is
$$K_{\alpha_{k} + \cdots + \alpha_{n}}K_{\alpha_{p} + \cdots + \alpha_{n}}
  e_{\alpha_{i} + \cdots + \alpha_{k-1}}e_{\alpha_{i} + \cdots + \alpha_{p-1}}.$$
Using the fact that $\left[e_{\alpha_{i} + \cdots + \alpha_{k-1}},e_{\alpha_{i} + \cdots + \alpha_{p-1}}\right]_{q}=0$
and that $\left(e_{\alpha_{i} + \cdots + \alpha_{p-1}}\right)^{N'} \in {\cal{I}}$, it is quite easy to see that
$$(E_{k,p})^{N'} \in {\cal{I}} \otimes U_{q}(\mathfrak{g}) + U_{q}(\mathfrak{g}) \otimes {\cal{I}},$$
which completes the proof of this case.

\end{subsubsection}

\begin{subsubsection}{Case 4. $(e_{\alpha_{i} + \cdots + 2\alpha_{j}+\cdots + 2\alpha_{n}})^{N'}$, $1 \leq i<j \leq n$}
We firstly consider 
$\big(\Delta (e_{\alpha_{i} + \cdots  + 2\alpha_{n}})\big)^{N'}$.
Writing $\Delta(e_{\alpha_{i} + \cdots  + 2\alpha_{n}}) = \sum_{k=i}^{n-1}D_{k} + D_{n} + D_{0} + D_{\infty}$,
the $q$-binomial theorem and one of the generalisations of the binomial theorem immediately gives 
\begin{eqnarray*}
\lefteqn{
\big(\Delta (e_{\alpha_{i} + \cdots  + 2\alpha_{n}} )\big)^{N'} } \\
 & = & 
\sum_{k=i}^{n-1} \left(D_{k}\right)^{N'} + \left( D_{n} + D_{0} + D_{\infty}\right)^{N'} \\
& = & \sum_{k=i}^{n-1} \left(D_{k}\right)^{N'} + \left\{ \begin{array}{ll}
(D_{n})^{N'} + (D_{\infty})^{N'}, & \mbox{if } N \equiv 0, 1, 3 \pmod{4}, \\
(D_{n})^{N'} + (D_{\infty})^{N'} + \phi_{N/2} (D_{0})^{N'} & \mbox{if } N \equiv 2 \pmod{4},
\end{array} \right.
\end{eqnarray*}
where $\phi_{N/2} \neq 0$.  
Now $$(e_{\alpha_{k} + \cdots  + 2\alpha_{n}})^{N'} \in {\cal{I}}, \hspace{10mm} k=i, \ldots, n-1,$$ 
and the fact that $(e_{n})^{\overline{N}} \in {\cal{I}}$ means that $(e_{n})^{2N'} \in {\cal{I}}$.  
By using these facts and examining some simple calculations it follows that 
$$(D_{k})^{N'}, (D_{\infty})^{N'} \in
{\cal{I}} \otimes U_{q}(\mathfrak{g}) + U_{q}(\mathfrak{g}) \otimes {\cal{I}},  \hspace{10mm} k=i, \ldots, n.$$

It remains to show that $(D_{0})^{N'} \in {\cal{I}} \otimes U_{q}(\mathfrak{g}) + U_{q}(\mathfrak{g}) \otimes {\cal{I}}$
when $N \equiv 2 \pmod{4}$; this follows from  $(e_{n})^{\overline{N}} \in {\cal{I}}$.

We now consider the more general problem: write 
$$\Delta(e_{\alpha_{i} + \cdots + 2\alpha_{j}+\cdots + 2\alpha_{n}}) = 
D_{0} + \sum_{k=i+1}^{n} D_{k} + \overline{D} + \sum_{l=j-1}^{n-1} F_{l}.$$
Using the $q$-binomial theorem we immediately obtain
\begin{eqnarray}
\big(\Delta (e_{\alpha_{i} + \cdots + 2\alpha_{j}+\cdots + 2\alpha_{n}} )\big)^{N'} & = & 
(D_{0})^{N'} + \sum_{k=i+1}^{j-1} (D_{k})^{N'}  \label{appendixC:blinco(1)} \\
& & + \big(\overline{D} + D_{j} + \cdots + D_{n-1} + D_{n} + F_{n-1} + F_{n-2} + \cdots + F_{j-1}\big)^{N'}. \nonumber
\end{eqnarray}
We can expand out the last component of the right hand side of this expression by using 
the following relations:
\begin{eqnarray*}
\lefteqn{
(\overline{D} + D_{j} + \cdots + D_{n-1})(F_{n-1} + F_{n-2} + \cdots + F_{j-1}) } \\
& & = q^{2}(F_{n-1} + F_{n-2} + \cdots + F_{j-1})(\overline{D} + D_{j} + \cdots + D_{n-1})
 + \xi(D_{n})^{2},
 \end{eqnarray*}
$$(\overline{D} + D_{j} + \cdots + D_{n-1}) D_{n} = q^{2} D_{n}(\overline{D} + D_{j} + \cdots + D_{n-1}),$$
$$D_{n}(F_{n-1} + F_{n-2} + \cdots + F_{j-1}) = q^{2} (F_{n-1} + F_{n-2} + \cdots + F_{j-1})D_{n},$$
 where $\xi=-(1+q)^{2}/(q-q^{-1})$.
Consequently,
\begin{eqnarray*}
\lefteqn{
\big(  \overline{D} + D_{j} + \cdots + D_{n-1} + D_{n} + F_{n-1} + F_{n-2} + \cdots + 
    F_{j-1}\big)^{N'} } \\
& = & \left\{  \begin{array}{ll}
(\overline{D} + D_{j} + \cdots + D_{n-1})^{N'} + (F_{n-1} + F_{n-2} + \cdots + F_{j-1})^{N'}, 
  & \mbox{if } N \equiv 0, 1, 3 \pmod{4}, \\
(\overline{D} + D_{j} + \cdots + D_{n-1})^{N'} + (F_{n-1} + F_{n-2} + \cdots + F_{j-1})^{N'} & \\
 \hspace{78mm} + \phi_{N/2}(D_{n})^{N'},  & \mbox{if } N \equiv 2 \pmod{4},
\end{array} \right.
\end{eqnarray*} 
where $\phi_{N/2} \neq 0$.  The $q$-binomial theorem then implies that
\begin{equation}
\label{appendixC:blinco(2)}
(\overline{D} + D_{j} + \cdots + D_{n-1})^{N'} = (\overline{D})^{N'} + \sum_{k=j}^{n-1} (D_{k})^{N'},
\end{equation}
\begin{equation}
\label{appendixC:blinco(3)}
(F_{n-1} + F_{n-2} + \cdots + F_{j-1})^{N'} = \sum_{k=j-1}^{n-1} (F_{k})^{N'}.
\end{equation}
Now the facts that 
$(e_{\alpha_{k} + \cdots + 2\alpha_{j}+\cdots + 2\alpha_{n}})^{N'} \in {\cal{I}}$ for each 
$k=i, \ldots, n-1$ and each $j=k+1, \ldots, n$, and that 
$(e_{\alpha_{i} + \cdots + \alpha_{k}})^{N'} \in {\cal{I}}$ for each $k=i+1, \ldots, n-1$, 
mean that each component of the right hand sides of 
(\ref{appendixC:blinco(1)})--(\ref{appendixC:blinco(3)})
belongs to 
${\cal{I}} \otimes U_{q}(\mathfrak{g}) + U_{q}(\mathfrak{g}) \otimes {\cal{I}}$.  

It remains to show that the same is true for
$(D_{n})^{N'}$ when $N \equiv 2 \pmod{4}$.
To see this, note that 
$(\overline{e}_{\alpha_{j} +\cdots + \alpha_{n}})^{\overline{N}} \in {\cal{I}}$ and that 
$\overline{N} = N/2$ for $N \equiv 2 \pmod{4}$, thus 
$(D_{n})^{N'} \in {\cal{I}} \otimes U_{q}(\mathfrak{g}) + U_{q}(\mathfrak{g}) \otimes {\cal{I}}$ when $N \equiv 2 \pmod{4}$.
This completes the proof of this case.

\end{subsubsection}

\begin{subsubsection}{The remaining elements of $I$}

We have shown that $\Delta\left(e_{\gamma}\right)^{N'}$ and 
$\Delta\left(e_{\beta}\right)^{\overline{N}}$
are elements of ${\cal{I}} \otimes U_{q}(\mathfrak{g}) + U_{q}(\mathfrak{g}) \otimes {\cal{I}}$.
We now prove that the same is true for the remaining elements of ${\cal{I}}$ 
in the following calculations, in which
we write the proportionality sign to mean that the left hand side is proportional 
to the right hand side with a non-zero scalar constant of
proportionality.  The cases we do not consider here are almost identically proved:
$$\Delta (f_{\gamma})^{N'} \propto \Delta\left(\omega \big(e_{\gamma}^{N'}\big) \right) 
\propto (\omega \otimes \omega) \circ \Delta' (e_{\gamma} )^{N'},$$
$$\Delta ( \overline{e}_{\gamma} )^{N'} 
\propto \Delta \left(S\big(e_{\gamma}^{N'}\big)\right) \Delta\left(K_{\gamma}\right)^{N'}
  \propto \left[ (S \otimes S) \circ  \Delta' (e_{\gamma})^{N'} \right] \Delta(K_{\gamma})^{N'},$$
$$\Delta \big(\overline{f}_{\gamma} \big)^{N'} 
\propto \Delta \left(\omega \big(\overline{e}_{\gamma}^{N'}\big) \right)
    \propto (\omega \otimes \omega) \circ \Delta'(\overline{e}_{\gamma})^{N'}.$$
Each of these expressions is an element of
 ${\cal{I}} \otimes U_{q}(\mathfrak{g}) + U_{q}(\mathfrak{g}) \otimes {\cal{I}}$
 from Proposition \ref{appendixC:carnegiehall(1)}, thus ${\cal{I}}$ is a two-sided co-ideal.

\end{subsubsection}

\end{section}

\begin{section}{${\cal{I}}$ is a two-sided Hopf ideal}
\label{appendixC:andrewblinco(4)}

We have proved that ${\cal{I}}$ is a two-sided ideal and a two-sided co-ideal
of $U_{q}(osp(1|2n))$.
To prove that ${\cal{I}}$ is a two-sided Hopf ideal all we need do is prove that 
$S(x) \in {\cal{I}}$ for each $x \in {\cal{I}}$, 
and to show this it suffices to show that $S(x) \in {\cal{I}}$ for each $x \in I$.
\begin{proposition}
\label{appendixC:carnegiehall(1)}
For each $x \in I$, $\omega(x) \in {\cal{I}}$ and $S(x) \in {\cal{I}}$.
\end{proposition}
\begin{proof}
Firstly note that $\omega(J_{i}^{\pm N}-1) = (J_{i}^{\mp N}-1) \in {\cal{I}}$ 
and $S(J_{i}^{\pm N}-1) = (J_{i}^{\mp N}-1) \in {\cal{I}}$. The rest of the proof follows from
Propositions \ref{sppendixB:steveieg(1)}--\ref{sppendixB:steveieg(5)} and the facts that $\omega$ is an involution and that
$\omega$ and $S$ are graded antiautomorphisms.
\end{proof}

\end{section}

\begin{section}{Technical results}

In this section we prove technical results used previously in this appendix.
\begin{proposition}
\label{appendixB:bronte(1)}
For each $j \in \mathbb{Z}_{+}$,
$$\sum_{k=1}^{2j} q^{-2k} (k)_{q} = q^{-4j} (2j+1)_{q} [j]^{q^{2}}, \hspace{10mm}
\sum_{k=1}^{2j+1} q^{-2k}(k)_{q} = q^{-2(2j+1)} (2j+1)_{q} [j+1]^{q^{2}}.$$
\end{proposition}
\begin{proof}
By direct calculation we have
\begin{itemize}
\item $q^{-2} (1)_{q} = q^{-2}$,
\item $\sum_{k=1}^{2} q^{-2k} (k)_{q} = q^{-4} (3)_{q}$,
\item $\sum_{k=1}^{3} q^{-2k} (k)_{q} = q^{-6} (3)_{q}[2]^{q^{2}}$,
\item $\sum_{k=1}^{4} q^{-2k} (k)_{q} = q^{-8} (5)_{q}[2]^{q^{2}}$.
\end{itemize}
Assume that the proposition is true for $\sum_{k=1}^{2j} q^{-2k}(k)_{q}$, then
\begin{eqnarray*}
\sum_{k=1}^{2j+1} q^{-2k}(k)_{q} & = & q^{-4j} (2j+1)_{q} [j]^{q^{2}} + q^{-4j-2} (2j+1)_{q} \\
& = & q^{-4j-2} (2j+1)_{q} ( q^{2} + q^{4} + \cdots + q^{2j}) + q^{-4j-2} (2j+1)_{q} \\
& = & q^{-4j-2} (2j+1)_{q} [j+1]^{q^{2}},
\end{eqnarray*}
as required.
Now assume that the proposition is true for $\sum_{k=1}^{2j+1} q^{-2k}(k)_{q}$, then
\begin{eqnarray*}
\sum_{k=1}^{2j+2} q^{-2k}(k)_{q} 
& = & q^{-4j-2} (2j+1)_{q} [j+1]^{q^{2}} + q^{-4j-4}(2j+2)_{q} \\
& = & q^{-4j-4} \big((2j+1)_{q} (q^{2} + q^{4} +\cdots + q^{2j+2})+ (2j+2)_{q} \big) \\
& = & q^{-4j-4} \big( (2j+1)_{q} \left[ j+1 \right]^{q^{2}} + (-q)^{2j+1} + (-q)^{2j+2} + \cdots + (-q)^{4j+2} \big) \\
& = & q^{-4j-4} \big( (2j+1)_{q} \left[ j+1 \right]^{q^{2}}+\left((-q)^{2j+1}+(-q)^{2j+2}\right)\left[ j+1 \right]^{q^{2}}
                       \big) \\
& = & q^{-4j-4}(2j+3)_{q}\left[ j+1 \right]^{q^{2}},
\end{eqnarray*}
completing the proof.
\end{proof}

\begin{proposition}
\label{appendixB:germany(1)}
For all $k, p$ satisfying $i \leq k<p \leq n$, 
\begin{eqnarray}
\left[ e_{\alpha_{k} + \cdots + \alpha_{n}}, e_{\alpha_{p} + \cdots + \alpha_{n}} \right]_{q} 
& = & (-q)^{n-p} e_{\alpha_{k} + \cdots + 2\alpha_{p} + \cdots + 2\alpha_{n}} \label{appendixB:margarat(2)} \\
&   &  + (q-q^{-1}) \left[ \sum_{j=1}^{n-p} (-q)^{j-1} 
e_{\alpha_{k} + \cdots + 2\alpha_{n-j+1} + \cdots + 2\alpha_{n}} e_{\alpha_{p} + \cdots + \alpha_{n-j}} \right].
\nonumber
\end{eqnarray}

\end{proposition}
\begin{proof}
The relation 
$\left[ e_{\alpha_{k} + \cdots + \alpha_{n}}, e_{n} \right]_{q} = e_{\alpha_{k} + \cdots + 2\alpha_{n}}$ proves
the proposition for $p=n$ and all $k$.  Assume that $k<n-1$, then
\begin{eqnarray*}
\left[ e_{\alpha_{k} + \cdots + \alpha_{n}}, e_{\alpha_{n-1}+\alpha_{n}} \right]_{q} & = &
       \left[ e_{\alpha_{k} + \cdots + \alpha_{n}}, e_{n-1}e_{n}-q^{-1}e_{n}e_{n-1} \right]_{q} \\
       & = & \left[ e_{\alpha_{k} + \cdots + \alpha_{n}}, e_{n-1}\right]_{q}e_{n} 
             + e_{n-1} \left[ e_{\alpha_{k} + \cdots + \alpha_{n}},e_{n}\right]_{q} \\
	& & -q^{-1} \left[ e_{\alpha_{k} + \cdots + \alpha_{n}},e_{n}\right]_{q} e_{n-1}
	    +q^{-1} e_{n} \left[ e_{\alpha_{k} + \cdots + \alpha_{n}},e_{n-1}\right]_{q} \\
       & = & e_{n-1} e_{\alpha_{k} + \cdots + 2\alpha_{n}} -q^{-1} e_{\alpha_{k} + \cdots + 2\alpha_{n}}e_{n-1},
\end{eqnarray*}
as $\left[ e_{\alpha_{k} + \cdots + \alpha_{n}}, e_{n-1}\right]_{q}=0$.
By re-writing the relation $\left[e_{\alpha_{k} + \cdots + 2\alpha_{n}},e_{n-1}\right]_{q} 
= e_{\alpha_{k} + \cdots + 2\alpha_{n-1}+ 2\alpha_{n}}$ as
$e_{n-1}e_{\alpha_{k} + \cdots + 2\alpha_{n}} = q e_{\alpha_{k} + \cdots + 2\alpha_{n}} e_{n-1}
  -q e_{\alpha_{k} + \cdots + 2\alpha_{n-1}+ 2\alpha_{n}}$, 
we obtain
$$\left[ e_{\alpha_{k} + \cdots + \alpha_{n}}, e_{\alpha_{n-1}+\alpha_{n}} \right]_{q}
 = -q e_{\alpha_{k} + \cdots + 2\alpha_{n-1}+ 2\alpha_{n}} 
 + (q-q^{-1})e_{\alpha_{k} + \cdots + 2\alpha_{n}}e_{n-1},$$
 proving the proposition for $p=n-1$.  
 Now assume that the proposition is true for some $p=k+2, \ldots, n-1$, then
\begin{eqnarray}
\lefteqn{
\left[ e_{\alpha_{k} + \cdots + \alpha_{n}}, e_{\alpha_{p-1} + \cdots + \alpha_{n}} \right]_{q} } \nonumber \\
 & = &
\left[ e_{\alpha_{k} + \cdots + \alpha_{n}}, e_{p-1}e_{\alpha_{p} + \cdots + \alpha_{n}}
                                             -q^{-1}e_{\alpha_{p} + \cdots + \alpha_{n}}e_{p-1} \right]_{q} \nonumber \\
         & = & \left[e_{\alpha_{k} + \cdots + \alpha_{n}}, e_{p-1} \right]_{q}e_{\alpha_{p} + \cdots + \alpha_{n}}
	       + e_{p-1} \left[e_{\alpha_{k} + \cdots + \alpha_{n}},e_{\alpha_{p} + \cdots + \alpha_{n}} \right]_{q} \nonumber \\
	  & & -q^{-1} \left[e_{\alpha_{k} + \cdots + \alpha_{n}},e_{\alpha_{p} + \cdots + \alpha_{n}} \right]_{q} e_{p-1} 
	     +q^{-1}e_{\alpha_{p} + \cdots + \alpha_{n}}\left[e_{\alpha_{k} + \cdots + \alpha_{n}},e_{p-1}\right]_{q}  \nonumber \\
	 & = &  e_{p-1} \left[e_{\alpha_{k} + \cdots + \alpha_{n}},e_{\alpha_{p} + \cdots + \alpha_{n}} \right]_{q}
	        -q^{-1} \left[e_{\alpha_{k} + \cdots + \alpha_{n}},e_{\alpha_{p} + \cdots + \alpha_{n}} \right]_{q} e_{p-1},
		\label{appendixB:margarat(1)}
\end{eqnarray}
as  $\left[e_{\alpha_{k} + \cdots + \alpha_{n}}, e_{p-1} \right]_{q}=0$.
Substituting (\ref{appendixB:margarat(2)}) into (\ref{appendixB:margarat(1)}) gives
\begin{eqnarray*}
& & 
(-q)^{n-p} e_{p-1} e_{\alpha_{k} + \cdots + 2\alpha_{p} + \cdots + 2\alpha_{n}}  \\
& & + (q-q^{-1}) \left[ \sum_{j=1}^{n-p} (-q)^{j-1} e_{p-1}
e_{\alpha_{k} + \cdots + 2\alpha_{n-j+1} + \cdots + 2\alpha_{n}} e_{\alpha_{p} + \cdots + \alpha_{n-j}} \right]  \\
& & -q^{-1}(-q)^{n-p} e_{\alpha_{k} + \cdots + 2\alpha_{p} + \cdots + 2\alpha_{n}} e_{p-1}  \\
& &  -q^{-1}(q-q^{-1}) \left[ \sum_{j=1}^{n-p} (-q)^{j-1} 
e_{\alpha_{k} + \cdots + 2\alpha_{n-j+1} + \cdots + 2\alpha_{n}} e_{\alpha_{p} + \cdots + \alpha_{n-j}}e_{p-1} \right]  \\
& = & -q(-q)^{n-p} e_{\alpha_{k} + \cdots + 2\alpha_{p-1} + \cdots + 2\alpha_{n}}
     + (q-q^{-1})(-q)^{n-p} e_{\alpha_{k} + \cdots + 2\alpha_{p} + \cdots + 2\alpha_{n}} e_{p-1} \\
& & + (q-q^{-1}) \Bigg[ \sum_{j=1}^{n-p} (-q)^{j-1}
e_{\alpha_{k} + \cdots + 2\alpha_{n-j+1} + \cdots + 2\alpha_{n}}e_{p-1}e_{\alpha_{p}+\cdots+\alpha_{n-j}}  \\
& & \hspace{30mm} -q^{-1} (-q)^{j-1}
e_{\alpha_{k} + \cdots + 2\alpha_{n-j+1} + \cdots + 2\alpha_{n}}e_{\alpha_{p}+\cdots+\alpha_{n-j}}e_{p-1} \Bigg] \\
& = & (-q)^{n-p+1}e_{\alpha_{k} + \cdots + 2\alpha_{p-1} + \cdots + 2\alpha_{n}} \\
& & + (q-q^{-1}) \left[ \sum_{j=1}^{n-p+1} (-q)^{j-1} 
e_{\alpha_{k} + \cdots + 2\alpha_{n-j+1} + \cdots + 2\alpha_{n}}e_{\alpha_{p-1}+\cdots+\alpha_{n-j}} \right],
\end{eqnarray*}
completing the induction.

\end{proof}

\begin{proposition} 
\label{appendixC:incrediblehulk(1)}
For all $1 \leq i < t \leq n$, 
$$\left[ e_{\alpha_{i} + \cdots + \alpha_{t-1}}, e_{\alpha_{t} + \cdots + \alpha_{n}}\right]_{q} = 
    e_{\alpha_{i} + \cdots + \alpha_{n}}.$$
\end{proposition}
\begin{proof}
Firstly $\left[e_{\alpha_{i} + \cdots + \alpha_{n-1}}, e_{n} \right]_{q} = e_{\alpha_{i} + \cdots + \alpha_{n}}$.
Now assume that 
$\left[ e_{\alpha_{i} + \cdots + \alpha_{t-1}}, e_{\alpha_{t} + \cdots + \alpha_{n}}\right]_{q} = 
    e_{\alpha_{i} + \cdots + \alpha_{n}}$ for some $i+2 \leq t \leq n$, then
\begin{eqnarray*}
\left[ e_{\alpha_{i} + \cdots + \alpha_{t-2}}, e_{\alpha_{t-1} + \cdots + \alpha_{n}}\right]_{q} & = & 
\left[ e_{\alpha_{i} + \cdots + \alpha_{t-2}}, 
               e_{t-1} e_{\alpha_{t} + \cdots + \alpha_{n}}-q^{-1}e_{\alpha_{t} + \cdots + \alpha_{n}}e_{t-1} \right]_{q} \\
	       & = & \left[ e_{\alpha_{i} + \cdots + \alpha_{t-2}},e_{t-1} \right]_{q}e_{\alpha_{t} + \cdots + \alpha_{n}}
	    +q^{-1} e_{t-1}\left[ e_{\alpha_{i} + \cdots + \alpha_{t-2}},e_{\alpha_{t} + \cdots + \alpha_{n}}\right]_{q} \\
& & -q^{-1} \left[ e_{\alpha_{i} + \cdots + \alpha_{t-2}},e_{\alpha_{t} + \cdots + \alpha_{n}}\right]_{q}e_{t-1}
-q^{-1}e_{\alpha_{t} + \cdots + \alpha_{n}} \left[ e_{\alpha_{i} + \cdots + \alpha_{t-2}},e_{t-1}\right]_{q} \\
& = & e_{\alpha_{i} + \cdots + \alpha_{t-1}}e_{\alpha_{t} + \cdots + \alpha_{n}}
-q^{-1}e_{\alpha_{t} + \cdots + \alpha_{n}}e_{\alpha_{i} + \cdots + \alpha_{t-1}} \\
& = & \left[ e_{\alpha_{i} + \cdots + \alpha_{t-1}},e_{\alpha_{t} + \cdots + \alpha_{n}}\right]_{q} 
=  e_{\alpha_{i} + \cdots + \alpha_{n}},
\end{eqnarray*}
as $\left[ e_{\alpha_{i} + \cdots + \alpha_{t-2}},e_{\alpha_{t} + \cdots + \alpha_{n}}\right]_{q}=0$.

\end{proof}

\begin{proposition}
\label{appendixC:incrediblehulk(2)}
For all $1 \leq i < t \leq j-1 < n$, 
$$\left[ e_{\alpha_{i} + \cdots + \alpha_{j-1}},e_{\alpha_{t} + \cdots + \alpha_{n}}\right]_{q}
    = (q-q^{-1}) e_{\alpha_{i} + \cdots + \alpha_{n}}e_{\alpha_{t} + \cdots + \alpha_{j-1}}.$$
\end{proposition}
\begin{proof}
Firstly $\left[ e_{\alpha_{i} + \cdots + \alpha_{j-1}},e_{\alpha_{j} + \cdots + \alpha_{n}}\right]_{q}
 = e_{\alpha_{i} + \cdots + \alpha_{n}}$, and
\begin{eqnarray*}
\left[ e_{\alpha_{i} + \cdots + \alpha_{j-1}},e_{\alpha_{j-1} + \cdots + \alpha_{n}}\right]_{q} & = & 
\left[ e_{\alpha_{i} + \cdots + \alpha_{j-1}},
e_{j-1} e_{\alpha_{j} + \cdots + \alpha_{n}} - q^{-1} e_{\alpha_{j} + \cdots + \alpha_{n}}e_{j-1} \right]_{q} \\
& = & \left[ e_{\alpha_{i} + \cdots + \alpha_{j-1}},e_{j-1} \right]_{q}e_{\alpha_{j} + \cdots + \alpha_{n}}
      +q e_{j-1}\left[ e_{\alpha_{i} + \cdots + \alpha_{j-1}},e_{\alpha_{j} + \cdots + \alpha_{n}}\right]_{q} \\
& &   -q^{-1} \left[ e_{\alpha_{i} + \cdots + \alpha_{j-1}},e_{\alpha_{j} + \cdots + \alpha_{n}}\right]_{q}e_{j-1}
      -q^{-2} e_{\alpha_{j} + \cdots + \alpha_{n}}\left[ e_{\alpha_{i} + \cdots + \alpha_{j-1}},e_{j-1}\right]_{q} \\
& = & q e_{j-1} e_{\alpha_{i} + \cdots + \alpha_{n}}-q^{-1}e_{\alpha_{i} + \cdots + \alpha_{n}}e_{j-1} 
  = (q-q^{-1})  e_{\alpha_{i} + \cdots + \alpha_{n}}e_{j-1},
\end{eqnarray*}
as $e_{j-1} e_{\alpha_{i} + \cdots + \alpha_{n}}=e_{\alpha_{i} + \cdots + \alpha_{n}}e_{j-1}$.
Now assume that 
$$\left[ e_{\alpha_{i} + \cdots + \alpha_{j-1}},e_{\alpha_{t} + \cdots + \alpha_{n}}\right]_{q}
    = (q-q^{-1}) e_{\alpha_{i} + \cdots + \alpha_{n}}e_{\alpha_{t} + \cdots + \alpha_{j-1}},$$ 
for some $i+2 \leq t \leq j-1$, then
\begin{eqnarray*}
\lefteqn{
\left[ e_{\alpha_{i} + \cdots + \alpha_{j-1}},e_{\alpha_{t-1} + \cdots + \alpha_{n}}\right]_{q} } \\
& = & 
\left[ e_{\alpha_{i} + \cdots + \alpha_{j-1}}, 
e_{t-1}e_{\alpha_{t} + \cdots + \alpha_{n}}-q^{-1}e_{\alpha_{t} + \cdots + \alpha_{n}}e_{t-1}\right]_{q} \\
& = & \left[ e_{\alpha_{i} + \cdots + \alpha_{j-1}},e_{t-1} \right]_{q}e_{\alpha_{t} + \cdots + \alpha_{n}}
     + e_{t-1} \left[ e_{\alpha_{i} + \cdots + \alpha_{j-1}},e_{\alpha_{t} + \cdots + \alpha_{n}}\right]_{q} \\
& & -q^{-1}\left[ e_{\alpha_{i} + \cdots + \alpha_{j-1}},e_{\alpha_{t} + \cdots + \alpha_{n}}\right]_{q}e_{t-1}
    -q^{-1} e_{\alpha_{t} + \cdots + \alpha_{n}}\left[ e_{\alpha_{i} + \cdots + \alpha_{j-1}},e_{t-1}\right]_{q} \\
& = & e_{t-1} \left[ e_{\alpha_{i} + \cdots + \alpha_{j-1}},e_{\alpha_{t} + \cdots + \alpha_{n}}\right]_{q}
      -q^{-1}\left[ e_{\alpha_{i} + \cdots + \alpha_{j-1}},e_{\alpha_{t} + \cdots + \alpha_{n}}\right]_{q}e_{t-1} \\
& = & (q-q^{-1}) e_{t-1}e_{\alpha_{i} + \cdots + \alpha_{n}}e_{\alpha_{t} + \cdots + \alpha_{j-1}}
-q^{-1}(q-q^{-1})e_{\alpha_{i} + \cdots + \alpha_{n}}e_{\alpha_{t} + \cdots + \alpha_{j-1}}e_{t-1} \\
& = & (q-q^{-1})e_{\alpha_{i} + \cdots + \alpha_{n}}e_{\alpha_{t-1} + \cdots + \alpha_{j-1}},
\end{eqnarray*}
as $\left[ e_{\alpha_{i} + \cdots + \alpha_{j-1}},e_{t-1} \right]_{q}=0$ and
$e_{t-1}$ commutes with $e_{\alpha_{i} + \cdots + \alpha_{n}}$.
\end{proof}

We can easily prove Propositions \ref{appendixB:propamos(1)} and
\ref{appendixB:propamos(2)} inductively using elementary calculations.
\begin{proposition}  
\label{appendixB:propamos(1)}
For all $1 \leq i < j < n$, 
\begin{eqnarray*}
\left[ e_{\alpha_{i} + \cdots + \alpha_{j}}, e_{\alpha_{i} + \cdots + 2\alpha_{j+1} + \cdots + 2 \alpha_{n}} \right]_{q}
& = & (-q)^{n-j-1} (1+q) (e_{\alpha_{i} + \cdots + \alpha_{n}})^{2} \\
& & + (q-q^{-1}) \sum_{p=1}^{n-j-1} (-q)^{p-1} 
          e_{\alpha_{i} + \cdots + \alpha_{j+p}}e_{\alpha_{i} + \cdots + 2\alpha_{j+p+1} + \cdots + 2\alpha_{n}}.
\end{eqnarray*}

\end{proposition}
\begin{proposition}  
\label{appendixB:propamos(2)}
For all $1 \leq i < j < n$, 
\begin{eqnarray*}
\left[\overline{e}_{\alpha_{i} + \cdots + 2\alpha_{j+1} + \cdots + 2\alpha_{n}}, 
     \overline{e}_{\alpha_{i} + \cdots + \alpha_{j}}, \right]_{q}
& = & (-q)^{n-j-1} (1+q) (\overline{e}_{\alpha_{i} + \cdots + \alpha_{n}})^{2} \\
& & + (q-q^{-1}) \sum_{p=1}^{n-j-1} (-q)^{p-1} 
      \overline{e}_{\alpha_{i} + \cdots + 2\alpha_{j+p+1} + \cdots + 2\alpha_{n}}
      \overline{e}_{\alpha_{i} + \cdots + \alpha_{j+p}}.
\end{eqnarray*}

\end{proposition}

\begin{proposition}
\label{appendixB:iceage(1)} 
The element
 $\left( q e_{\alpha_{j} + 2\alpha_{j+1} + \cdots + 2 \alpha_{n}} e_{j} 
   		-q^{-1} e_{j}e_{\alpha_{j} + 2\alpha_{j+1} + \cdots + 2 \alpha_{n}}\right)$
commutes with $e_{k}$ for each  $k=1, \ldots, n$.
\end{proposition}
\begin{proof}
Firstly, the Serre relations state that $e_{k}$ commutes with $e_{i}$ if $|k-i|>1$.
For each $k = j+2, \ldots, n$, we have the relation
$$\left[e_{\alpha_{j} + 2\alpha_{j+1} + \cdots + 2 \alpha_{n}}, e_{k} \right]_{q}=0,$$ 
which states that each such $e_{k}$ commutes with $e_{\alpha_{j} + 2\alpha_{j+1} + \cdots + 2 \alpha_{n}}$.

Now we will show that $e_{j+1}$ commutes with
$\left( q e_{\alpha_{j} + 2\alpha_{j+1} + \cdots + 2 \alpha_{n}} e_{j} 
   		-q^{-1} e_{j}e_{\alpha_{j} + 2\alpha_{j+1} + \cdots + 2 \alpha_{n}}\right)$. We calculate that
\begin{eqnarray*}
\lefteqn{
\left( q e_{\alpha_{j} + 2\alpha_{j+1} + \cdots + 2 \alpha_{n}} e_{j} 
   		-q^{-1} e_{j}e_{\alpha_{j} + 2\alpha_{j+1} + \cdots + 2 \alpha_{n}}\right) e_{j+1} } \\
& & \hspace{5mm} = q e_{\alpha_{j} + 2\alpha_{j+1} + \cdots + 2\alpha_{n}} e_{j} e_{j+1}
      -e_{j} e_{j+1} e_{\alpha_{j} + 2\alpha_{j+1} + \cdots + 2 \alpha_{n}} \\
& & \hspace{5mm}= q e_{\alpha_{j} + 2\alpha_{j+1} + \cdots + 2\alpha_{n}} e_{\alpha_{j} + \alpha_{j+1}}
      + e_{\alpha_{j} + 2\alpha_{j+1} + \cdots + 2\alpha_{n}} e_{j+1} e_{j} \\
&   & \hspace{15mm} -e_{\alpha_{j} + \alpha_{j+1}} e_{\alpha_{j} + 2\alpha_{j+1} + \cdots + 2\alpha_{n}}
      -q^{-1} e_{j+1} e_{j}  e_{\alpha_{j} + 2\alpha_{j+1} + \cdots + 2\alpha_{n}} \\
& & \hspace{5mm}= e_{j+1}\left(q  e_{\alpha_{j} + 2\alpha_{j+1} + \cdots + 2\alpha_{n}} e_{j}
      -q^{-1} e_{j} e_{\alpha_{j} + 2\alpha_{j+1} + \cdots + 2\alpha_{n}}\right),
\end{eqnarray*}
where we have used the relations
$$\left[e_{\alpha_{j} + 2\alpha_{j+1} + \cdots + 2 \alpha_{n}},e_{j+1} \right]_{q}=0, \hspace{10mm}
\left[e_{\alpha_{j}+\alpha_{j+1}}, e_{\alpha_{j} + 2\alpha_{j+1} + \cdots + 2\alpha_{n}}\right]_{q}=0.$$

To complete the proof we will show that 
\begin{eqnarray*}
\lefteqn{
e_{j} \left( q e_{\alpha_{j} + 2\alpha_{j+1} + \cdots + 2 \alpha_{n}} e_{j} 
   		-q^{-1} e_{j}e_{\alpha_{j} + 2\alpha_{j+1} + \cdots + 2 \alpha_{n}}\right) } \\
& & \hspace{10mm} =	\left( q e_{\alpha_{j} + 2\alpha_{j+1} + \cdots + 2 \alpha_{n}} e_{j} 
   		-q^{-1} e_{j}e_{\alpha_{j} + 2\alpha_{j+1} + \cdots + 2 \alpha_{n}}\right) e_{j},
		\end{eqnarray*}	
and to prove this we note that
\begin{eqnarray*}
\left[ e_{j}, e_{\alpha_{j} + 2\alpha_{j+1} + \cdots + 2 \alpha_{n}} \right]_{q} 
& = & e_{j} e_{\alpha_{j} + 2\alpha_{j+1} + \cdots + 2 \alpha_{n}} 
      - e_{\alpha_{j} + 2\alpha_{j+1} + \cdots + 2 \alpha_{n}} e_{j} \\
& = & \sum_{k=j}^{n-1} C_{k} e_{\alpha_{j} + \cdots+ \alpha_{k}}e_{\alpha_{j}+\cdots + 2\alpha_{k+1} + \cdots + 2 \alpha_{n}}
   + C_{n} (e_{\alpha_{j} + \cdots+ \alpha_{n}})^{2}, \ C_{k} \in \mathbb{C}, 
\end{eqnarray*}
and that $\left[ e_{j},  e_{\alpha_{j} + \cdots+ \alpha_{k}}\right]_{q}=0$ and
$\left[ e_{j}, e_{\alpha_{j}+\cdots + 2\alpha_{k+1} + \cdots + 2 \alpha_{n}} \right]_{q}=0$ for each
$k=j+1, \ldots, n-1$,
and thus $$e_{j} e_{\alpha_{j} + \cdots+ \alpha_{k}}=q e_{\alpha_{j} + \cdots+ \alpha_{k}} e_{j},$$
$$e_{j} e_{\alpha_{j}+\cdots + 2\alpha_{k+1} + \cdots + 2 \alpha_{n}}
      = q e_{\alpha_{j}+\cdots + 2\alpha_{k+1} + \cdots + 2 \alpha_{n}} e_{j}.$$ 
Consequently,
\begin{eqnarray*}
\lefteqn{
e_{j} \left( q e_{\alpha_{j} + 2\alpha_{j+1} + \cdots + 2 \alpha_{n}} e_{j} 
   		-q^{-1} e_{j}e_{\alpha_{j} + 2\alpha_{j+1} + \cdots + 2 \alpha_{n}}\right)  } \\
& = & q \left( e_{\alpha_{j} + 2\alpha_{j+1} + \cdots + 2 \alpha_{n}} (e_{j})^{2}
 +\sum_{k=j}^{n-1} C_{k} e_{\alpha_{j} + \cdots+ \alpha_{k}}e_{\alpha_{j}+\cdots + 2\alpha_{k+1} + \cdots + 2 \alpha_{n}}e_{j}
 + C_{n} (e_{\alpha_{j} + \cdots+ \alpha_{n}})^{2} e_{j} \right) \\
& & -q^{-1} e_{j} \left(
e_{\alpha_{j} + 2\alpha_{j+1} + \cdots + 2 \alpha_{n}} e_{j}
+\sum_{k=j}^{n-1} C_{k} e_{\alpha_{j} + \cdots+ \alpha_{k}}e_{\alpha_{j}+\cdots + 2\alpha_{k+1} + \cdots + 2 \alpha_{n}}
 + C_{n} (e_{\alpha_{j} + \cdots+ \alpha_{n}})^{2}  \right) \\
& = & \left( q e_{\alpha_{j} + 2\alpha_{j+1} + \cdots + 2 \alpha_{n}} e_{j}
             -q^{-1} e_{j}e_{\alpha_{j} + 2\alpha_{j+1} + \cdots + 2 \alpha_{n}}\right) e_{j}.
\end{eqnarray*}

\end{proof}

\begin{proposition}
\label{appendixB:lem20}
For each $j=2, 3, \ldots, n-1$,
$$[e_{\alpha_{j} + 2\alpha_{j+1} + \cdots + 2\alpha_{n}}, \overline{e}_{\alpha_{j-1}+\alpha_{j}}]_{q}
               = \overline{e}_{\alpha_{j-1} + 2\alpha_{j} + \cdots + 2\alpha_{n}}.$$
\end{proposition}
\begin{proof}
The proof of this proposition is lengthy, and we prove in in a number of stages.
Firstly consider 
$\left[e_{\alpha_{j} + \alpha_{j+1}}, \overline{e}_{\alpha_{j-1} + \alpha_{j}}\right]_{q}$ for each
$j=2, 3, \ldots, n-1$:
\begin{eqnarray*}
\left[e_{\alpha_{j} + \alpha_{j+1}}, \overline{e}_{\alpha_{j-1} + \alpha_{j}}\right]_{q} & = &
                \left[e_{j} e_{j+1} - q^{-1} e_{j+1} e_{j} , \overline{e}_{\alpha_{j-1} + \alpha_{j}}\right]_{q} \\
        & = &    e_{j} \overline{e}_{\alpha_{j-1} + \alpha_{j} + \alpha_{j+1}}
		 + q^{-1} \left[ e_{j}, [e_{j}, e_{j-1}]_{q} \right]_{q} e_{j+1} \\
& & 	\hspace{5mm}	 -q^{-1} e_{j+1} \left[ e_{j}, [e_{j}, e_{j-1}]_{q} \right]_{q}
		 - \overline{e}_{\alpha_{j-1} + \alpha_{j} + \alpha_{j+1}}e_{j} \\
		 & = & e_{j} \overline{e}_{\alpha_{j-1} + \alpha_{j} + \alpha_{j+1}}
		       - \overline{e}_{\alpha_{j-1} + \alpha_{j} + \alpha_{j+1}}e_{j} 
		   = 0,
\end{eqnarray*}
as $\left[ e_{j}, [e_{j}, e_{j-1}]_{q} \right]_{q}=0$ is just a restatement of the Serre relation 
$(ad_{q} e_{j})^{2}(e_{j-1})=0$ for $j<n$,   
and $\left[e_{\alpha_{j-1} + \alpha_{j} + \alpha_{j+1}},e_{j}\right]_{q}=0$ implies
$\left[ e_{j}, \overline{e}_{\alpha_{j-1} + \alpha_{j} + \alpha_{j+1}}\right]_{q}=0$.

Now we will show that 
$$\left[e_{\alpha_{j} + \cdots + \alpha_{k}}, \overline{e}_{\alpha_{j-1} + \alpha_{j}} \right]_{q}=0,$$
for each $k = j+1, \ldots, n$.  The calculation immediately above states that 
this is true for $k=j+1$, 
now assume that it is true for some $k \geq j+1$, then
\begin{eqnarray*}
\left[e_{\alpha_{j} + \cdots + \alpha_{k+1}}, \overline{e}_{\alpha_{j-1} + \alpha_{j}}\right]_{q} & = & 
\left[e_{\alpha_{j} + \cdots + \alpha_{k}} e_{k+1} - q^{-1} e_{k+1}e_{\alpha_{j} + \cdots + \alpha_{k}},
\overline{e}_{\alpha_{j-1} + \alpha_{j}}\right]_{q} \\
& = &  e_{\alpha_{j} + \cdots + \alpha_{k}} \left[e_{k+1}, \overline{e}_{\alpha_{j-1} + \alpha_{j}}\right]_{q}
   + \left[e_{\alpha_{j} + \cdots + \alpha_{k}}, \overline{e}_{\alpha_{j-1} + \alpha_{j}}\right]_{q} e_{k+1} \\
   & & -q^{-1} e_{k+1} \left[e_{\alpha_{j} + \cdots + \alpha_{k}}, \overline{e}_{\alpha_{j-1} + \alpha_{j}}\right]_{q}
       -q^{-1} \left[e_{k+1}, \overline{e}_{\alpha_{j-1} + \alpha_{j}}\right]_{q} e_{\alpha_{j} + \cdots + \alpha_{k}} \\
    & = & 0,
\end{eqnarray*}
as $\left[e_{k+1}, \overline{e}_{\alpha_{j-1} + \alpha_{j}}\right]_{q}=0$ and
$\left[e_{\alpha_{j} + \cdots + \alpha_{k}}, \overline{e}_{\alpha_{j-1} + \alpha_{j}} \right]_{q}=0$ by assumption.

Now for each $j < n-1$, we claim that 
$$\left[e_{\alpha_{j}+\cdots+2\alpha_{k}+\cdots+2\alpha_{n}}, \overline{e}_{\alpha_{j-1} + \alpha_{j}} \right]_{q}=0,$$
for each $k=j+2, \ldots, n$.
We firstly show that this is true for $k=n$:
\begin{eqnarray*}
\left[e_{\alpha_{j} + \cdots + 2 \alpha_{n}}, \overline{e}_{\alpha_{j-1} + \alpha_{j}} \right]_{q} & = &
\left[e_{\alpha_{j} + \cdots + \alpha_{n}}e_{n} + e_{n}e_{\alpha_{j} + \cdots + \alpha_{n}}, 
\overline{e}_{\alpha_{j-1} + \alpha_{j}} \right]_{q}  \\
& = & e_{\alpha_{j} + \cdots + \alpha_{n}} \left[ e_{n}, \overline{e}_{\alpha_{j-1} + \alpha_{j}} \right]_{q}
+ \left[e_{\alpha_{j} + \cdots + \alpha_{n}}, \overline{e}_{\alpha_{j-1} + \alpha_{j}} \right]_{q} e_{n} \\
& & + e_{n} \left[e_{\alpha_{j} + \cdots + \alpha_{n}}, \overline{e}_{\alpha_{j-1} + \alpha_{j}} \right]_{q}
    + \left[ e_{n}, \overline{e}_{\alpha_{j-1} + \alpha_{j}} \right]_{q}e_{\alpha_{j} + \cdots + \alpha_{n}} \\
& = & 0,
\end{eqnarray*}
as the preceding calculation implies that 
$\left[e_{\alpha_{j} + \cdots + \alpha_{n}}, \overline{e}_{\alpha_{j-1} + \alpha_{j}}\right]_{q}=0$ 
and $\left[ e_{n}, \overline{e}_{\alpha_{j-1} + \alpha_{j}} \right]_{q}=0$ from the Serre relations.  
Now assume that 
$$\left[e_{\alpha_{j} + \cdots + 2\alpha_{k+1} + \cdots + 2 \alpha_{n}}, \overline{e}_{\alpha_{j-1} + \alpha_{j}} \right]_{q}=0,$$
for some $k+1 = j+3, \ldots, n$, then
\begin{eqnarray*}
\lefteqn{
\left[e_{\alpha_{j} + \cdots + 2\alpha_{k} + \cdots + 2 \alpha_{n}}, \overline{e}_{\alpha_{j-1} + \alpha_{j}} \right]_{q} } \\
 & = & \left[e_{\alpha_{j} + \cdots + 2\alpha_{k+1} + \cdots + 2 \alpha_{n}}e_{k}
-q^{-1} e_{k}e_{\alpha_{j} + \cdots + 2\alpha_{k+1} + \cdots + 2 \alpha_{n}}, 
 \overline{e}_{\alpha_{j-1} + \alpha_{j}} \right]_{q}  \\
& = & e_{\alpha_{j} + \cdots + 2\alpha_{k+1} + \cdots + 2 \alpha_{n}} \left[ e_{k}, \overline{e}_{\alpha_{j-1} + \alpha_{j}}\right]_{q}
+ \left[e_{\alpha_{j} + \cdots + 2\alpha_{k+1} + \cdots + 2 \alpha_{n}}, \overline{e}_{\alpha_{j-1} + \alpha_{j}} \right]_{q}e_{k} \\
& & 
-q^{-1}e_{k}\left[e_{\alpha_{j} + \cdots + 2\alpha_{k+1} + \cdots + 2 \alpha_{n}}, \overline{e}_{\alpha_{j-1} + \alpha_{j}} \right]_{q}
-q^{-1} \left[ e_{k}, \overline{e}_{\alpha_{j-1} + \alpha_{j}}\right]_{q}e_{\alpha_{j} + \cdots + 2\alpha_{k+1} + \cdots + 2 \alpha_{n}}
\\
& = & 0,
\end{eqnarray*}
as $\left[ e_{k}, \overline{e}_{\alpha_{j-1} + \alpha_{j}}\right]_{q}=0$ and
$\left[e_{\alpha_{j} + \cdots + 2\alpha_{k+1} + \cdots + 2 \alpha_{n}}, \overline{e}_{\alpha_{j-1} + \alpha_{j}} \right]_{q}=0$
by assumption.

Now consider $\left[e_{\alpha_{j} + 2\alpha_{j+1} + \cdots + 2 \alpha_{n}}, \overline{e}_{\alpha_{j-1} + \alpha_{j}}
\right]_{q}$:
\begin{eqnarray*}
\lefteqn{
\left[e_{\alpha_{j} + 2\alpha_{j+1} + \cdots + 2 \alpha_{n}}, \overline{e}_{\alpha_{j-1} + \alpha_{j}} \right]_{q} } \\
& = & \left[e_{\alpha_{j} + \alpha_{j+1} + 2\alpha_{j+2}+\cdots + 2 \alpha_{n}}e_{j+1} 
-q^{-1}e_{j+1}e_{\alpha_{j} + \alpha_{j+1} + 2\alpha_{j+2}+\cdots + 2 \alpha_{n}}, 
\overline{e}_{\alpha_{j-1} + \alpha_{j}} \right]_{q} \\
& = & e_{\alpha_{j} + \alpha_{j+1} + 2\alpha_{j+2} + \cdots + 2 \alpha_{n}} \left[ e_{j+1}, \overline{e}_{\alpha_{j-1} + \alpha_{j}}\right]_{q}
+ q^{-1}\left[e_{\alpha_{j} + \alpha_{j+1} + 2\alpha_{j+2} + \cdots + 2 \alpha_{n}}, \overline{e}_{\alpha_{j-1} + \alpha_{j}} \right]_{q}e_{j+1} \\
& & 
-q^{-1}e_{j+1}\left[e_{\alpha_{j}+\alpha_{j+1}+2\alpha_{j+2}+\cdots+2\alpha_{n}}, \overline{e}_{\alpha_{j-1} + \alpha_{j}} \right]_{q}
-q^{-1} \left[ e_{j+1}, \overline{e}_{\alpha_{j-1}+\alpha_{j}}\right]_{q}
e_{\alpha_{j}+\alpha_{j+1}+2\alpha_{j+2}+\cdots+2\alpha_{n}} \\
& = & e_{\alpha_{j} + \alpha_{j+1} + 2\alpha_{j+2} + \cdots + 2 \alpha_{n}}\overline{e}_{\alpha_{j-1} + \alpha_{j}+\alpha_{j+1}}
-q^{-1} \overline{e}_{\alpha_{j-1}+\alpha_{j}+\alpha_{j+1}} e_{\alpha_{j}+\alpha_{j+1}+2\alpha_{j+2}+\cdots+2\alpha_{n}} \\
& = & \left[ e_{\alpha_{j} + \alpha_{j+1} + 2\alpha_{j+2} + \cdots + 2 \alpha_{n}},
\overline{e}_{\alpha_{j-1} + \alpha_{j}+\alpha_{j+1}} \right]_{q},
\end{eqnarray*}
as
$\left[e_{\alpha_{j} + \alpha_{j+1} + 2\alpha_{j+2} + \cdots + 2 \alpha_{n}}, \overline{e}_{\alpha_{j-1} + \alpha_{j}} \right]_{q}=0$
from the preceding calculation.

Now we claim that
$$\left[e_{\alpha_{j} + \cdots + 2\alpha_{k} + \cdots + 2\alpha_{n}},\overline{e}_{\alpha_{j-1}+ \cdots +\alpha_{k-1}} \right]_{q}
 = \left[e_{\alpha_{j} + \cdots + 2\alpha_{k+1} + \cdots + 2\alpha_{n}},\overline{e}_{\alpha_{j-1}+ \cdots +\alpha_{k}}
 \right]_{q},$$
 for each $k=j+1, \ldots, n$.
This is true for $k=j+1$ from the preceding calculation, now assume that it is true for some $k \geq j+1$, then
\begin{eqnarray*}
\lefteqn{
\left[e_{\alpha_{j} + \cdots + 2\alpha_{k} + \cdots + 2\alpha_{n}},\overline{e}_{\alpha_{j-1}+ \cdots +\alpha_{k-1}} \right]_{q} } \\
& = & \left[e_{\alpha_{j} + \cdots + 2\alpha_{k+1} + \cdots + 2\alpha_{n}}e_{k}
-q^{-1}e_{k}e_{\alpha_{j} + \cdots + 2\alpha_{k+1} + \cdots + 2\alpha_{n}} ,\overline{e}_{\alpha_{j-1}+ \cdots +\alpha_{k-1}}
\right]_{q} \\
& = &
e_{\alpha_{j} + \cdots + 2\alpha_{k+1} + \cdots + 2\alpha_{n}}\left[e_{k},\overline{e}_{\alpha_{j-1}+ \cdots +\alpha_{k-1}} \right]_{q}
+q^{-1}\left[e_{\alpha_{j} + \cdots + 2\alpha_{k+1}+\cdots+2\alpha_{n}},\overline{e}_{\alpha_{j-1}+\cdots+\alpha_{k-1}}\right]_{q}
e_{k} \\
& & -q^{-1} e_{k}\left[e_{\alpha_{j}+\cdots+2\alpha_{k+1}+\cdots+2\alpha_{n}},\overline{e}_{\alpha_{j-1}+\cdots+\alpha_{k-1}}\right]_{q} 
 -q^{-1} \left[e_{k},\overline{e}_{\alpha_{j-1}+ \cdots +\alpha_{k-1}} \right]_{q}e_{\alpha_{j}+\cdots+2\alpha_{k+1}+\cdots+2\alpha_{n}}
 \\
& = & e_{\alpha_{j} + \cdots + 2\alpha_{k+1} + \cdots + 2\alpha_{n}}\left[e_{k},\overline{e}_{\alpha_{j-1}+ \cdots +\alpha_{k-1}} \right]_{q}
-q^{-1} \left[e_{k},\overline{e}_{\alpha_{j-1}+ \cdots +\alpha_{k-1}}
\right]_{q}e_{\alpha_{j}+\cdots+2\alpha_{k+1}+\cdots+2\alpha_{n}} \\
& = & \left[e_{\alpha_{j}+\cdots+2\alpha_{k+1}+\cdots+2\alpha_{n}},\overline{e}_{\alpha_{j-1}+\cdots+\alpha_{k}} \right]_{q},
\end{eqnarray*}
where we have used the result  
$\left[e_{\alpha_{j}+\cdots+2\alpha_{k+1}+\cdots+2\alpha_{n}},\overline{e}_{\alpha_{j-1}+\cdots+\alpha_{k-1}} \right]_{q}=0$
which we will now prove.  To prove this last result recall that 
$\left[e_{\alpha_{j}+\cdots+2\alpha_{k+1}+\cdots+2\alpha_{n}},\overline{e}_{\alpha_{j-1}+\alpha_{j}} \right]_{q}=0$. 
Now assume that 
$\left[e_{\alpha_{j}+\cdots+2\alpha_{k+1}+\cdots+2\alpha_{n}},\overline{e}_{\alpha_{j-1}+\cdots + \alpha_{m}} \right]_{q}=0$
for some $m$ satisfying $j \leq m \leq k-2$, then 
\begin{eqnarray*}
\lefteqn{
\left[e_{\alpha_{j}+\cdots+2\alpha_{k+1}+\cdots+2\alpha_{n}},\overline{e}_{\alpha_{j-1}+\cdots + \alpha_{m+1}} \right]_{q} } \\
& = & \left[e_{\alpha_{j}+\cdots+2\alpha_{k+1}+\cdots+2\alpha_{n}},
e_{m+1}\overline{e}_{\alpha_{j-1}+\cdots +\alpha_{m}}-q^{-1}\overline{e}_{\alpha_{j-1}+\cdots +\alpha_{m}}e_{m+1} \right]_{q} \\
& = & \left[e_{\alpha_{j}+\cdots+2\alpha_{k+1}+\cdots+2\alpha_{n}},e_{m+1}\right]_{q}\overline{e}_{\alpha_{j-1}+\cdots +\alpha_{m}}
+ e_{m+1}\left[e_{\alpha_{j}+\cdots+2\alpha_{k+1}+\cdots+2\alpha_{n}},\overline{e}_{\alpha_{j-1}+\cdots +\alpha_{m}}\right]_{q}
\\
& & -q^{-1}\left[e_{\alpha_{j}+\cdots+2\alpha_{k+1}+\cdots+2\alpha_{n}},\overline{e}_{\alpha_{j-1}+\cdots +\alpha_{m}}
 \right]_{q}e_{m+1} - q^{-1}\overline{e}_{\alpha_{j-1}+\cdots +\alpha_{m}}
 \left[e_{\alpha_{j}+\cdots+2\alpha_{k+1}+\cdots+2\alpha_{n}},e_{m+1}\right]_{q} \\
& = & 0,
\end{eqnarray*}
as $\left[e_{\alpha_{j}+\cdots+2\alpha_{k+1}+\cdots+2\alpha_{n}},\overline{e}_{\alpha_{j-1}+\cdots + \alpha_{m}} \right]_{q}=0$
by assumption and $\left[e_{\alpha_{j}+\cdots+2\alpha_{k+1}+\cdots+2\alpha_{n}},e_{m+1}\right]_{q}=0$.
We have thus shown that
$$[e_{\alpha_{j} + 2\alpha_{j+1} + \cdots + 2\alpha_{n}}, \overline{e}_{\alpha_{j-1}+\alpha_{j}}]_{q} = 
[e_{\alpha_{j} + \cdots + 2\alpha_{n}}, \overline{e}_{\alpha_{j-1}+\cdots+\alpha_{n-1}}]_{q}.$$
Now
\begin{eqnarray*}
\lefteqn{
[e_{\alpha_{j} + \cdots + 2\alpha_{n}}, \overline{e}_{\alpha_{j-1}+\cdots +\alpha_{n-1}}]_{q} } \\
& & \hspace{5mm} = [e_{\alpha_{j} + \cdots + \alpha_{n}}e_{n} 
+ e_{n}e_{\alpha_{j} + \cdots + \alpha_{n}}, \overline{e}_{\alpha_{j-1}+\cdots +\alpha_{n-1}}]_{q} \\
& & \hspace{5mm}= e_{\alpha_{j} + \cdots + \alpha_{n}} [e_{n}, \overline{e}_{\alpha_{j-1}+\cdots +\alpha_{n-1}}]_{q}
     + q^{-1} [e_{\alpha_{j} + \cdots + \alpha_{n}}, \overline{e}_{\alpha_{j-1}+\cdots +\alpha_{n-1}}]_{q}e_{n} \\
 & & \hspace{10mm}+ e_{n} [e_{\alpha_{j} + \cdots + \alpha_{n}}, \overline{e}_{\alpha_{j-1}+\cdots +\alpha_{n-1}}]_{q}
     + [e_{n}, \overline{e}_{\alpha_{j-1}+\cdots +\alpha_{n-1}}]_{q}e_{\alpha_{j} + \cdots + \alpha_{n}} \\
& &\hspace{5mm} = e_{\alpha_{j} + \cdots + \alpha_{n}} [e_{n}, \overline{e}_{\alpha_{j-1}+\cdots +\alpha_{n-1}}]_{q}
      + [e_{n}, \overline{e}_{\alpha_{j-1}+\cdots +\alpha_{n-1}}]_{q}e_{\alpha_{j} + \cdots + \alpha_{n}} \\
& & \hspace{5mm}= [e_{\alpha_{j} + \cdots + \alpha_{n}}, \overline{e}_{\alpha_{j-1}+\cdots +\alpha_{n}}]_{q},
\end{eqnarray*}
where we have used the result
$[e_{\alpha_{j} + \cdots + \alpha_{n}}, \overline{e}_{\alpha_{j-1}+\cdots +\alpha_{n-1}}]_{q}=0$ which we now prove.
Recall that
$\left[e_{\alpha_{j} + \cdots + \alpha_{n}}, \overline{e}_{\alpha_{j-1} + \alpha_{j}} \right]_{q}=0$
for each $j=2, \ldots, n-1$.  Assume that
$\left[e_{\alpha_{j} + \cdots + \alpha_{n}}, \overline{e}_{\alpha_{j-1} + \cdots + \alpha_{k}} \right]_{q}=0$
for some $k=j, \ldots, n-2$, then
\begin{eqnarray*}
\lefteqn{
\left[e_{\alpha_{j} + \cdots + \alpha_{n}}, \overline{e}_{\alpha_{j-1} + \cdots + \alpha_{k+1}} \right]_{q} } \\
& = &  \left[e_{\alpha_{j} + \cdots + \alpha_{n}}, e_{k+1}\overline{e}_{\alpha_{j-1} + \cdots + \alpha_{k}}
-q^{-1}\overline{e}_{\alpha_{j-1} + \cdots + \alpha_{k}}e_{k+1} \right]_{q} \\
& = & \left[e_{\alpha_{j} + \cdots + \alpha_{n}}, e_{k+1} \right]_{q}\overline{e}_{\alpha_{j-1} + \cdots + \alpha_{k}} 
+ e_{k+1}\left[e_{\alpha_{j} + \cdots + \alpha_{n}},\overline{e}_{\alpha_{j-1} + \cdots + \alpha_{k}}\right]_{q} \\
& & -q^{-1} \left[e_{\alpha_{j} + \cdots + \alpha_{n}},\overline{e}_{\alpha_{j-1} + \cdots + \alpha_{k}}\right]_{q}e_{k+1}
    -q^{-1}\overline{e}_{\alpha_{j-1} + \cdots + \alpha_{k}}\left[e_{\alpha_{j} + \cdots + \alpha_{n}},e_{k+1}\right]_{q} \\
& = & 0,
\end{eqnarray*}
as $\left[e_{\alpha_{j} + \cdots + \alpha_{n}}, e_{k+1} \right]_{q}=0$ and 
$\left[e_{\alpha_{j} + \cdots + \alpha_{n}}, \overline{e}_{\alpha_{j-1} + \cdots + \alpha_{k}} \right]_{q}=0$ by assumption.
An implication of this is
 $[e_{\alpha_{j} + \cdots + \alpha_{n}}, \overline{e}_{\alpha_{j-1}+\cdots +\alpha_{n-1}}]_{q}=0$.
It follows that
$$[e_{\alpha_{j} + \cdots + 2\alpha_{n}}, \overline{e}_{\alpha_{j-1}+\cdots +\alpha_{n-1}}]_{q} = 
[e_{\alpha_{j} + \cdots + \alpha_{n}}, \overline{e}_{\alpha_{j-1}+\cdots +\alpha_{n}}]_{q}.$$

Now
\begin{eqnarray*}
\lefteqn{
[e_{\alpha_{j} + \cdots + \alpha_{n}}, \overline{e}_{\alpha_{j-1}+\cdots +\alpha_{n}}]_{q} } \\
& & \hspace{5mm}= [e_{\alpha_{j} + \cdots + \alpha_{n-1}}e_{n} - q^{-1} e_{n}e_{\alpha_{j} + \cdots + \alpha_{n-1}},
	\overline{e}_{\alpha_{j-1}+\cdots +\alpha_{n}}]_{q} \\
& & \hspace{5mm}= e_{\alpha_{j} + \cdots + \alpha_{n-1}}[e_{n},\overline{e}_{\alpha_{j-1}+\cdots +\alpha_{n}}]_{q}
      -[e_{\alpha_{j} + \cdots + \alpha_{n-1}},\overline{e}_{\alpha_{j-1}+\cdots +\alpha_{n}}]_{q}e_{n} \\
& & \hspace{10mm} -q^{-1} e_{n} [e_{\alpha_{j} + \cdots + \alpha_{n-1}},\overline{e}_{\alpha_{j-1}+\cdots +\alpha_{n}}]_{q}
    -q^{-1}[e_{n},\overline{e}_{\alpha_{j-1}+\cdots +\alpha_{n}}]_{q}e_{\alpha_{j} + \cdots + \alpha_{n-1}} \\
& & \hspace{5mm}= e_{\alpha_{j} + \cdots + \alpha_{n-1}}\overline{e}_{\alpha_{j-1}+\cdots +2\alpha_{n}}
      -q^{-1} \overline{e}_{\alpha_{j-1}+\cdots +2\alpha_{n}}e_{\alpha_{j} + \cdots + \alpha_{n-1}} \\
& & \hspace{5mm}= [e_{\alpha_{j} + \cdots + \alpha_{n-1}},\overline{e}_{\alpha_{j-1}+\cdots +2\alpha_{n}}]_{q},
\end{eqnarray*}
as $[e_{\alpha_{j} + \cdots + \alpha_{n-1}}, \overline{e}_{\alpha_{j-1}+\cdots +\alpha_{n}}]_{q}=0$, and
\begin{eqnarray*}
\lefteqn{
[e_{\alpha_{j} + \cdots + \alpha_{n-1}},\overline{e}_{\alpha_{j-1}+\cdots +2\alpha_{n}}]_{q} } \\
& = & [e_{\alpha_{j} + \cdots + \alpha_{n-2}}e_{n-1}-q^{-1}e_{n-1}e_{\alpha_{j} + \cdots + \alpha_{n-2}},
\overline{e}_{\alpha_{j-1}+\cdots +2\alpha_{n}}]_{q} \\
& = & e_{\alpha_{j} + \cdots + \alpha_{n-2}}[e_{n-1},\overline{e}_{\alpha_{j-1}+\cdots +2\alpha_{n}}]_{q}
+q^{-1}[e_{\alpha_{j} + \cdots + \alpha_{n-2}},\overline{e}_{\alpha_{j-1}+\cdots +2\alpha_{n}}]_{q}e_{n-1} \\
& & -q^{-1}e_{n-1}[e_{\alpha_{j} + \cdots + \alpha_{n-2}},\overline{e}_{\alpha_{j-1}+\cdots +2\alpha_{n}}]_{q}
    -q^{-1}[e_{n-1},\overline{e}_{\alpha_{j-1}+\cdots +2\alpha_{n}}]_{q}e_{\alpha_{j} + \cdots + \alpha_{n-2}} \\
& = & e_{\alpha_{j} + \cdots + \alpha_{n-2}}[e_{n-1},\overline{e}_{\alpha_{j-1}+\cdots +2\alpha_{n}}]_{q}
-q^{-1}[e_{n-1},\overline{e}_{\alpha_{j-1}+\cdots +2\alpha_{n}}]_{q}e_{\alpha_{j} + \cdots + \alpha_{n-2}} \\
& = & [e_{\alpha_{j} + \cdots + \alpha_{n-2}},\overline{e}_{\alpha_{j-1}+\cdots +2\alpha_{n-1}+2\alpha_{n}}]_{q},
\end{eqnarray*}
as $[e_{\alpha_{j} + \cdots + \alpha_{n-2}},\overline{e}_{\alpha_{j-1}+\cdots +2\alpha_{n}}]_{q}=0$.

We now claim that 
$$[e_{\alpha_{j} + \cdots + \alpha_{k+1}}, \overline{e}_{\alpha_{j-1}+\cdots +2\alpha_{k+2}+\cdots+2\alpha_{n}}]_{q}
= [e_{\alpha_{j} + \cdots + \alpha_{k}}, \overline{e}_{\alpha_{j-1}+\cdots +2\alpha_{k+1}+\cdots+2\alpha_{n}}]_{q},$$
for each $k+1=j+2,\ldots,n-1$.  This is true for $k+1=n-1$, and assume that it is true for
some $k+1 = j+3, \ldots, n-1$, then 
\begin{eqnarray*}
\lefteqn{[e_{\alpha_{j} + \cdots + \alpha_{k}}, \overline{e}_{\alpha_{j-1}+\cdots +2\alpha_{k+1}+\cdots+2\alpha_{n}}]_{q} } \\
& = & [e_{\alpha_{j} + \cdots + \alpha_{k-1}}e_{k}-q^{-1}e_{k}e_{\alpha_{j} + \cdots + \alpha_{k-1}},
\overline{e}_{\alpha_{j-1}+\cdots +2\alpha_{k+1}+\cdots+2\alpha_{n}}]_{q} \\
& = & e_{\alpha_{j} + \cdots + \alpha_{k-1}}[e_{k},\overline{e}_{\alpha_{j-1}+\cdots +2\alpha_{k+1}+\cdots+2\alpha_{n}}]_{q}
+q^{-1}[e_{\alpha_{j} + \cdots + \alpha_{k-1}},\overline{e}_{\alpha_{j-1}+\cdots +2\alpha_{k+1}+\cdots+2\alpha_{n}}]_{q}e_{k} \\
& &-q^{-1}e_{k}[e_{\alpha_{j} + \cdots + \alpha_{k-1}},\overline{e}_{\alpha_{j-1}+\cdots +2\alpha_{k+1}+\cdots+2\alpha_{n}}]_{q}
-q^{-1}[e_{k},\overline{e}_{\alpha_{j-1}+\cdots +2\alpha_{k+1}+\cdots+2\alpha_{n}}]_{q}e_{\alpha_{j} + \cdots + \alpha_{k-1}} \\
& = & e_{\alpha_{j} + \cdots + \alpha_{k-1}}[e_{k},\overline{e}_{\alpha_{j-1}+\cdots +2\alpha_{k+1}+\cdots+2\alpha_{n}}]_{q}
-q^{-1}[e_{k},\overline{e}_{\alpha_{j-1}+\cdots +2\alpha_{k+1}+\cdots+2\alpha_{n}}]_{q}e_{\alpha_{j} + \cdots + \alpha_{k-1}} \\
& = & [e_{\alpha_{j} + \cdots + \alpha_{k-1}},\overline{e}_{\alpha_{j-1}+\cdots +2\alpha_{k}+\cdots+2\alpha_{n}}]_{q},
\end{eqnarray*}
as $[e_{\alpha_{j} + \cdots + \alpha_{k-1}},\overline{e}_{\alpha_{j-1}+\cdots +2\alpha_{k+1}+\cdots+2\alpha_{n}}]_{q}=0$.

To complete the proof it is a simple matter to show
the following results using calculations almost identical to those immediately above:
$$[e_{\alpha_{j} + \alpha_{j+1}}, \overline{e}_{\alpha_{j-1}+\alpha_{j}+\alpha_{j+1}+2\alpha_{j+2}+\cdots+2\alpha_{n}}]_{q}
=[e_{\alpha_{j}}, \overline{e}_{\alpha_{j-1}+\alpha_{j}+2\alpha_{j+1}+\cdots+2\alpha_{n}}]_{q},$$
$$[e_{\alpha_{j}}, \overline{e}_{\alpha_{j-1}+\alpha_{j}+2\alpha_{j+1}+\cdots+2\alpha_{n}}]_{q} = 
\overline{e}_{\alpha_{j-1}+2\alpha_{j}+\cdots+2\alpha_{n}}.$$
\end{proof}

\begin{proposition}
\label{appendixB:troythebattleprop}
\begin{eqnarray}
\lefteqn{
q e_{j-1} \overline{e}_{\alpha_{j-1} + 2\alpha_{j} + \cdots + 2\alpha_{n}} 
  -q^{-1}\overline{e}_{\alpha_{j-1} + 2\alpha_{j} + \cdots + 2\alpha_{n}} e_{j-1} } \nonumber \\
  &  & \hspace{5mm}= q e_{\alpha_{j-1}+2\alpha_{j}+\cdots+2\alpha_{n}} e_{j-1}
       -q^{-1} e_{j-1} e_{\alpha_{j-1}+2\alpha_{j}+\cdots+2\alpha_{n}}. \label{appendixB:sudan(3)}
\end{eqnarray}
\end{proposition}
\begin{proof}
Proposition \ref{appendixB:lem20} implies the following result for each $j=2, \ldots, n-1$:
$$\overline{e}_{\alpha_{j-1} + 2\alpha_{j} + \cdots + 2\alpha_{n}} = 
\left[ e_{\alpha_{j}+2\alpha_{j+1}+\cdots+2\alpha_{n}}, \overline{e}_{\alpha_{j-1}+\alpha_{j}} \right]_{q},$$
which we can use to rewrite the left hand side of (\ref{appendixB:sudan(3)}):
\begin{eqnarray}
\lefteqn{
q e_{j-1} \overline{e}_{\alpha_{j-1} + 2\alpha_{j} + \cdots + 2\alpha_{n}} 
  -q^{-1}\overline{e}_{\alpha_{j-1} + 2\alpha_{j} + \cdots + 2\alpha_{n}} e_{j-1} } \nonumber \\
 & = & q e_{j-1} e_{\alpha_{j}+2\alpha_{j+1}+\cdots+2\alpha_{n}}\overline{e}_{\alpha_{j-1}+\alpha_{j}}
       -e_{j-1} \overline{e}_{\alpha_{j-1}+\alpha_{j}}e_{\alpha_{j}+2\alpha_{j+1}+\cdots+2\alpha_{n}} \nonumber \\
 & & -q^{-1} e_{\alpha_{j}+2\alpha_{j+1}+\cdots+2\alpha_{n}} \overline{e}_{\alpha_{j-1}+\alpha_{j}} e_{j-1}
     +q^{-2} \overline{e}_{\alpha_{j-1}+\alpha_{j}}e_{\alpha_{j}+2\alpha_{j+1}+\cdots+2\alpha_{n}} e_{j-1} \nonumber \\
 & = & q e_{j-1} e_{\alpha_{j}+2\alpha_{j+1}+\cdots+2\alpha_{n}}\overline{e}_{\alpha_{j-1}+\alpha_{j}}
      -e_{\alpha_{j}+2\alpha_{j+1}+\cdots+2\alpha_{n}}e_{j-1}\overline{e}_{\alpha_{j-1}+\alpha_{j}} \nonumber \\
 & & -q^{-1} \overline{e}_{\alpha_{j-1}+\alpha_{j}}e_{j-1}e_{\alpha_{j}+2\alpha_{j+1}+\cdots+2\alpha_{n}}
     +q^{-2} \overline{e}_{\alpha_{j-1}+\alpha_{j}}e_{\alpha_{j}+2\alpha_{j+1}+\cdots+2\alpha_{n}} e_{j-1}
     \label{appendixB:sudan} \\
& = &  q e_{\alpha_{j-1}+\alpha_{j}+2\alpha_{j+1}+\cdots+2\alpha_{n}}\overline{e}_{\alpha_{j-1}+\alpha_{j}}
      -q^{-1} \overline{e}_{\alpha_{j-1}+\alpha_{j}}e_{\alpha_{j-1}+\alpha_{j}+2\alpha_{j+1}+\cdots+2\alpha_{n}} \nonumber \\
& = &  q e_{\alpha_{j-1}+\alpha_{j}+2\alpha_{j+1}+\cdots+2\alpha_{n}} e_{j} e_{j-1}
       -e_{\alpha_{j-1}+\alpha_{j}+2\alpha_{j+1}+\cdots+2\alpha_{n}} e_{j-1}e_{j} \nonumber \\
& & -q^{-1} e_{j} e_{j-1} e_{\alpha_{j-1}+\alpha_{j}+2\alpha_{j+1}+\cdots+2\alpha_{n} }
    +q^{-2}e_{j-1} e_{j} e_{\alpha_{j-1}+\alpha_{j}+2\alpha_{j+1}+\cdots+2\alpha_{n} }  \nonumber \\
& = & q e_{\alpha_{j-1}+\alpha_{j}+2\alpha_{j+1}+\cdots+2\alpha_{n}} e_{j} e_{j-1}
     -q^{-1} e_{j-1} e_{\alpha_{j-1}+\alpha_{j}+2\alpha_{j+1}+\cdots+2\alpha_{n}} e_{j} \label{appendixB:sudan(2)} \\
& & -e_{j}e_{\alpha_{j-1}+\alpha_{j}+2\alpha_{j+1}+\cdots+2\alpha_{n}}  e_{j-1}
    +q^{-2} e_{j-1} e_{j} e_{\alpha_{j-1}+\alpha_{j}+2\alpha_{j+1}+\cdots+2\alpha_{n}} \nonumber \\
& = & q e_{\alpha_{j-1}+2\alpha_{j}+\cdots+2\alpha_{n}} e_{j-1}
       -q^{-1} e_{j-1} e_{\alpha_{j-1}+2\alpha_{j}+\cdots+2\alpha_{n}}. \nonumber 
\end{eqnarray}
 We used the relation 
$\overline{e}_{\alpha_{j-1}+\alpha_{j}} e_{j-1} = q e_{j-1}\overline{e}_{\alpha_{j-1}+\alpha_{j}}$
to obtain (\ref{appendixB:sudan}) and  the relation
$e_{j-1} e_{\alpha_{j-1}+\alpha_{j}+2\alpha_{j+1}+\cdots+2\alpha_{n}}
 = q e_{\alpha_{j-1}+\alpha_{j}+2\alpha_{j+1}+\cdots+2\alpha_{n}}e_{j-1}$
  to obtain (\ref{appendixB:sudan(2)}).

\end{proof}

\end{section}

\end{chapter}